\documentclass[leqno]{article} 

\usepackage[frenchb,english]{babel}
\usepackage[T1]{fontenc}

\usepackage{amsmath}
\usepackage{amssymb}
\usepackage{amsfonts}
\usepackage{enumerate}
\usepackage{vmargin}
\usepackage[all]{xy}
\usepackage{mathrsfs}
\usepackage{mathtools}
\usepackage{lmodern}
\usepackage[colorlinks=true,pagebackref=true,urlcolor=blue]{hyperref}
\usepackage{ulem}
\usepackage[all]{xy}
\usepackage{comment}
\usepackage{centernot}

\setmarginsrb{3cm}{3cm}{3.5cm}{3cm}{0cm}{0cm}{1.5cm}{3cm}

\newcommand{\A}{\ensuremath{\mathcal{A}}}
\newcommand{\B}{\mathrm{B}}
\newcommand{\C}{\ensuremath{\mathbb{C}}}
\newcommand{\D}{\ensuremath{\mathcal{D}}}
\newcommand{\E}{\ensuremath{\mathbb{E}}}
\newcommand{\F}{\ensuremath{\mathbb{F}}}
\newcommand{\Fc}{\ensuremath{\mathcal{F}}}
\newcommand{\G}{\mathrm{Gap}} 

\let\H\relax 
\newcommand{\H}{\mathrm{H}}
\newcommand{\HS}{\mathrm{HS}} 
\newcommand{\HI}{\mathrm{H}^\infty} 

\newcommand{\I}{\mathrm{I}}
\let\L\relax
\newcommand{\L}{\mathrm{L}}
\newcommand{\M}{\mathrm{M}}
\newcommand{\Mat}{\ensuremath{\mathbb{M}}}
\newcommand{\N}{\ensuremath{\mathbb{N}}}
\let\O\relax
\newcommand{\O}{\mathrm{O}} 
\let\P\relax 
\newcommand{\P}{\mathcal{P}} 
\newcommand{\Q}{\ensuremath{\mathbb{Q}}}
\newcommand{\R}{\ensuremath{\mathbb{R}}}
\let\S\relax 
\newcommand{\S}{\mathrm{S}} 
\newcommand{\T}{\ensuremath{\mathbb{T}}}
\newcommand{\Z}{\ensuremath{\mathbb{Z}}}

\newcommand{\W}{\mathrm{W}}
\newcommand{\vect}{\ensuremath{\mathop{\rm span\,}\nolimits}}

\newcommand{\scr}{\mathscr} 
\renewcommand{\leq}{\ensuremath{\leqslant}}
\renewcommand{\geq}{\ensuremath{\geqslant}}
\newcommand{\qed}{\hfill \vrule height6pt  width6pt depth0pt}

\newcommand{\bnorm}[1]{ \big\| #1  \big\|}
\newcommand{\Bnorm}[1]{ \Big\| #1  \Big\|}
\newcommand{\bgnorm}[1]{ \bigg\| #1  \bigg\|}
\newcommand{\Bgnorm}[1]{ \Bigg\| #1  \Bigg\|}
\newcommand{\norm}[1]{\left\Vert#1\right\Vert}
\newcommand{\xra}{\xrightarrow}

\newcommand{\ot}{\otimes}
\newcommand{\epsi}{\varepsilon}
\newcommand{\ovl}{\overline}
\newcommand{\otvn}{\ovl\ot}
\newcommand{\ul}{\mathcal{U}}
\newcommand{\dsp}{\displaystyle}
\newcommand{\co}{\colon}
\renewcommand{\d}{\mathop{}\mathopen{}\mathrm{d}} 
\let\i\relax 
\newcommand{\i}{\mathrm{i}} 
\newcommand{\ov}{\overset}

\newcommand{\rad}{\mathrm{rad}}

\newcommand{\Aut}{\mathrm{Aut}}

\newcommand{\QWEP}{\mathrm{QWEP}}
\newcommand{\sign}{\mathrm{sign}}

\newcommand{\disc}{\mathrm{disc}}
\newcommand{\dist}{\mathrm{dist}} 
\newcommand{\fin}{\mathrm{fin}} 
\newcommand{\Id}{\mathrm{Id}} 
\newcommand{\VN}{\mathrm{VN}} 

\newcommand{\CBAP}{\mathrm{CBAP}}
\newcommand{\CCAP}{\mathrm{CCAP}}
\newcommand{\UMD}{\mathrm{UMD}}
\newcommand{\e}{\mathrm{e}} 
\let\ker\relax 
\DeclareMathOperator{\ker}{Ker} 
\DeclareMathOperator{\Ran}{Ran} 

\DeclareMathOperator{\weakstar}{w*-}
\DeclareMathOperator{\supp}{supp} 
\DeclareMathOperator{\card}{card} 
\DeclareMathOperator{\Span}{span} 
\DeclareMathOperator{\dom}{dom} 
\DeclareMathOperator{\pv}{p.v.} 
\DeclareMathOperator{\tr}{\mathrm{Tr}} 
\DeclareMathOperator{\sgn}{\mathrm{sgn}} 
\DeclareMathOperator{\vol}{\mathrm{vol}} 
\DeclareMathOperator{\Lip}{\mathrm{Lip}} 
\let\Re\relax 
\DeclareMathOperator{\Re}{Re} 
\let\Im\relax 
\DeclareMathOperator{\Im}{Im} 
\newcommand{\cb}{\mathrm{cb}} 
\newcommand{\Gauss}{\mathrm{Gauss}} 
\newcommand{\sa}{\mathrm{sa}} 
\newcommand{\mk}{\mathrm{mk}} 
\newcommand{\Tron}{\mathcal{T}}
\newcommand{\VNGfin}{\Span\left\{ \lambda_s : \: s \in G \right\}}

\newcommand{\SpI}{\ovl{\Ran A_p}}
\newcommand{\StwoI}{\ovl{\Ran A_2}}
\newcommand{\SinftyI}{\ovl{\Ran A_\infty}}
\newcommand{\MIfinzero}{\M_{I,\fin} \cap \Ran A}

\allowdisplaybreaks

\newtheorem{thm}{Theorem}[section]
\newtheorem{defi}[thm]{Definition}

\newtheorem{prop}[thm]{Proposition}

\newtheorem{cor}[thm]{Corollary}
\newtheorem{lemma}[thm]{Lemma}

\newtheorem{remark}[thm]{Remark}

\newenvironment{proof}[1][]{\noindent {\it Proof #1} : }{\hbox{~}\qed
\smallskip
}

\numberwithin{equation}{section}
\usepackage[nottoc,notlot,notlof]{tocbibind}

\let\OLDthebibliography\thebibliography
\renewcommand\thebibliography[1]{
  \OLDthebibliography{#1}
  \setlength{\parskip}{0pt}
  \setlength{\itemsep}{0pt plus 0.3ex}
}

\begin{document}
\selectlanguage{english}
\title{\bfseries{Riesz transforms, Hodge-Dirac operators and functional calculus for multipliers I}}
\date{}
\author{\bfseries{C\'edric Arhancet - Christoph Kriegler}}

\maketitle


\begin{abstract}
In this work, we solve the problem explicitly stated at the end of a paper of Junge, Mei and Parcet [JEMS2018, Problem C.5] for a large class of groups including all amenable groups and free groups. More precisely, we prove that the Hodge-Dirac operator of the canonical ``hidden'' noncommutative geometry associated with a Markov semigroup $(T_t)_{t \geq 0}$ of  Fourier multipliers is bisectorial and admits a bounded $\H^\infty$ functional calculus on a bisector which implies a positive answer to the quoted problem. Our result can be seen as a strengthening of the dimension free estimates of Riesz transforms of the above authors and also allows us to provide Hodge decompositions. A part of our proof relies on a new transference argument between multipliers which is of independent interest. Our results are even new for the Poisson semigroup on $\T^n$. We also provide a similar result for Markov semigroups of Schur multipliers and dimension free estimates for noncommutative Riesz transforms associated with these semigroups. Along the way, we also obtain new Khintchine type equivalences for $q$-Gaussians in $\L^p$-spaces associated to crossed products. Our results allow us to introduce new spectral triples (i.e. noncommutative manifolds) and new quantum (locally) compact metric spaces, in connection with the carr\'e du champ, which summarize the underlying geometry of our setting. Finally, our examples lead us to introduce a Banach space variant of the notion of spectral triple suitable for our context. 
\end{abstract}


\makeatletter
 \renewcommand{\@makefntext}[1]{#1}
 \makeatother

{
  \hypersetup{linkcolor=blue}
 \tableofcontents
}

\newpage

\section{Introduction}
\label{sec:Introduction}

The continuity of the Hilbert transform on $\L^p(\R)$ by Riesz \cite{Ri1} is known as one of the greatest discoveries in analysis of the twentieth century. This transformation is at the heart of many areas : complex analysis, harmonic analysis, Banach space geometry, martingale theory and signal processing. We refer to the thick books \cite{HvNVW1}, \cite{HvNVW2}, \cite{Kin1} and \cite{Kin2} and references therein for more information. Directional Riesz transforms $R_j$ are higher-dimensional generalizations of the Hilbert transform defined by the formula
\begin{equation}
\label{}
R_j
\ov{\mathrm{def}}{=}  \partial_j \circ (-\Delta)^{-\frac{1}{2}}, \quad j=1,\ldots,n
\end{equation}
where $\Delta$ is the laplacian on $\R^n$. Generalizing Riesz's result, Calder\'on and Zygmund proved in \cite{CaZ1} that these operators are bounded on $\L^p(\R^n)$ if $1<p<\infty$. It is known that the $\L^p$-norms do not depend on the dimension $n$ 
by \cite{IwG1} and \cite{BaW1}. In \cite{Ste3} (see also \cite{IwG1} and \cite{BaW1}), Stein showed that the vectorial Riesz transform $\partial (-\Delta)^{-\frac{1}{2}}$ satisfies
\begin{equation}
\label{ine-intro-6}
\bnorm{\partial  (-\Delta)^{-\frac{1}{2}}f}_{\L^p(\R^n,\ell^2)}
\lesssim_p \norm{f}_{\L^p(\R^n)} 
\end{equation}
with free dimensional bound where $\partial f \ov{\mathrm{def}}{=} (\partial_1 f,\ldots,\partial_n f)$ is the gradient of a function $f$ belonging to some suitable subspace of $\L^p(\R^n)$. Furthermore, by duality we have by e.g. \cite[Proposition 2.1]{CoD1} an equivalence of the form
\begin{equation}
\label{Equiv-Riesz-intro}
\bnorm{(-\Delta)^{\frac12}f}_{\L^p(\R^n)} 
\approx_{p} \norm{\partial_{}f}_{\L^p(\R^n,\ell^2_n)}. 
\end{equation}
Note that this equivalence can be seen as a variant of the famous Kato square root problem solved in \cite{AHLMT1} and in \cite{AKM1}, see also \cite{Tch1} and \cite{HLM1}.

The study of these operators and variants in many contexts has a long history and was the source of many fundamental developments such as the Calder\'on-Zygmund theory \cite{CaZ1} or the classical work of Stein \cite{Ste1} on Littlewood-Paley Theory. Nowadays, Riesz transforms associated to various geometric structures is a recurrent theme in analysis and geometry, see for example the papers \cite{Bak1}, \cite{CCH1} (riemannian manifolds), \cite{Ste1} (compact Lie groups), \cite{Mey6}, \cite{Pis14} (Gaussian spaces), \cite{Lus1} (fermions algebras), \cite{Lus3} (deformed Gaussians algebras), \cite{Lus2}, \cite{DoP1} (abelian groups), \cite{Lus4} (generalized Heisenberg groups), \cite{BaR1} (graphs), \cite{AsO1} (Schr\"odinger operators) and  \cite{Aus1} (elliptic operators). Indeed, Riesz transforms have become a cornerstone of analysis and the literature is quite huge and it would be impossible to give complete references here.

An important generalization was given by Meyer \cite{Mey6}. It consists in replacing the Laplacian $\Delta$ by the $\L^p$-realization $A_p$ of the negative infinitesimal generator $A$ of a Markov semigroup $(T_t)_{t \geq 0}$ of operators acting on the $\L^p$-spaces of a measure space $\Omega$ and to replace the gradient $\partial$ by the ``carr\'e du champ'' $\Gamma$ introduced by Roth \cite{Rot1} (see also \cite{Hir1}) defined\footnote{\thefootnote. Here, the domain of $A$ must contain a suitable involutive algebra.} by
\begin{equation}
\label{carre-du-champ}
\Gamma(f,g)
\ov{\mathrm{def}}{=}  \frac{1}{2}\big[ A(\ovl{f}) g + \ovl{f} A(g) - A(\ovl{f} g)\big].
\end{equation}
In the case of the Heat semigroup $(\e^{-t \Delta})_{t \geq 0}$ with generator $\Delta$, we recover the gradient form $\langle \partial f, \partial g \rangle_{\ell^2_n}$. Meyer was interested by the  equivalence
\begin{equation}
\label{Meyer-estimates}
\bnorm{A_p^{\frac{1}{2}}(f)}_{\L^p(\Omega)}
\approx_{p}  \bnorm{\Gamma(f,f)^{\frac{1}{2}}}_{\L^p(\Omega)} 
\end{equation}
on some suitable subspace (ideally $\dom A_p^{\frac{1}{2}}$) of $\L^p(\Omega)$. Meyer proved such equivalence for the Ornstein-Uhlenbeck semigroup. Nevertheless, with sharp contrast, if $1<p<2$ these estimates are surprisingly false for the Poisson semigroup on $\L^p(\R^n)$ which is a Markov semigroup of Fourier multipliers, see \cite[Appendix D]{JMP2}. Actually, other examples of semigroups illustrating this phenomenon are already present in the papers of Lust-Piquard \cite[Proposition 2.9]{Lus2} and \cite[p.~283]{Lus1} relying on an observation of Lamberton. 

Of course, when something goes wrong with a mathematical problem it is rather natural to change slightly the formulation of the problem in order to obtain a natural positive statement. By introducing some gradients with values in a \textit{noncommutative} space, Junge, Mei and Parcet obtained in \cite{JMP2} dimension free estimates for Riesz transforms associated with arbitrary Markov semigroups $(T_t)_{t \geq 0}$ of Fourier multipliers acting on classical $\L^p$-spaces $\L^p(\hat{G})$ where $G$ is for example an abelian discrete group with (compact) dual group $\hat{G}$ (and more generally on the noncommutative $\L^p$-spaces $\L^p(\VN(G))$ associated with a nonabelian group $G$). We denote by $\psi \co G \to \C$ the symbol of the (negative) infinitesimal generator $A$ of the semigroup. In the spirit of \eqref{Equiv-Riesz-intro}, the above authors proved estimates of the form 
\begin{equation} 
\label{Riesz-Parcet-commutatif}
\bnorm{A_p^{\frac12}(f)}_{\L^p(\hat{G})} 
\approx_{p} \norm{\partial_{\psi,1,p}(f)}_{\L^p(\L^\infty(\Omega) \rtimes_\alpha G)}
\end{equation}
where $\partial_{\psi,1,p}$ is some kind of gradient defined on a dense subspace of the classical $\L^p$-space $\L^p(\hat{G})$.
It takes values in a closed subpace $\Omega_{\psi,1,p} $ of a \textit{noncommutative} $\L^p$-space $\L^p(\L^\infty(\Omega) \rtimes_\alpha G)$ associated with some crossed product $\L^\infty(\Omega) \rtimes_\alpha G$ where $\Omega$ is a probability space and where $\alpha \co G \to \Aut(\L^\infty(\Omega))$ is an action of $G$ on $\L^\infty(\Omega)$ determined by the semigroup. Let us explain the simplest case, i.e. the case where $\alpha$ is trivial. In this \textit{non-crossed and very particular situation}, we have an identification of $\L^p(\L^\infty(\Omega) \rtimes_\alpha G)$ with the classical $\L^p$-space $\L^p(\Omega \ot \hat{G})$ and the map $\partial_{\psi,1,p}$ is defined on the span of characters $\langle s, \cdot \rangle_{G,\hat{G}}$ in $\L^p(\hat{G})$ with values in $\L^p(\Omega \ot \hat{G})$. It is defined by
\begin{equation}
\label{def-partial-psi-very-particular}
\partial_{\psi,1,p}\big(\langle s, \cdot \rangle_{G,\hat{G}}\big)
\ov{\mathrm{def}}{=} \W(b_\psi(s)) \ot \langle s, \cdot \rangle_{G,\hat{G}}.
\end{equation}
where $\W \co H \to \L^0(\Omega)$ is an $H$-isonormal \textit{Gaussian process}\footnote{\thefootnote. In particular, for any $h \in H$ the random variable $\W(h)$ is a centred real Gaussian.} for some real Hilbert space $H$ and where $b_\psi \co G \to H$ is a specific function satisfying
\begin{equation}
\label{liens-psi-bpsi}
\psi(s)
=\norm{b_\psi(s)}_H^2,\quad s \in G.
\end{equation}
We refer to Subsection \ref{Sec-infos-on-Gamma-Fourier} for the (crossed) general situation where the action $\alpha$ is obtained by second quantization from an orthogonal representation $\pi \co G \to \B(H)$ associated to the semigroup, see \eqref{equ-markovian-automorphism-group}. 

The approach by Junge, Mei and Parcet highlights an intrinsic noncommutativity since $\Omega_{\psi,1,p} $ is in general a highly noncommutative object \textit{although} the group $G$ may be abelian. It is fair to say that this need of noncommutativity was first noticed and explicitly written by Lust-Piquard in \cite{Lus1} and \cite{Lus2} in some particular cases under a somewhat different but essentially equivalent form of \eqref{Riesz-Parcet-commutatif}. Moreover, it is remarkable that the estimates of \cite{Lus2} were exploited in a decisive way by Naor \cite{Nao1} to understand subtle geometric phenomena. Finally, note that the existence of gradients suitable for arbitrary Markov semigroups of linear operators appears already in the work of Sauvageot and Cipriani, see \cite{Sau1}, \cite{CiS1} and the survey \cite{Cip1}. Finally, we refer to \cite{ArK1}, \cite{Kos1}, \cite{JMX} and \cite{PiX} and references therein for more information on noncommutative $\L^p$-spaces.
 
In the context of Riesz transforms, the authors of the classical and remarkable paper \cite{AKM1} were the first to introduce suitable Hodge-Dirac operators. The $\L^p$-boundedness of the $\HI$ calculus of this unbounded operator allows everyone to obtain immediately the $\L^p$-boundedness of Riesz transforms. The authors of \cite{JMP2} introduced a similar operator in the context of Markov semigroups $(T_t)_{t \geq 0}$ of Fourier multipliers acting on classical $\L^p$-spaces and more generally on noncommutative $\L^p$-spaces $\L^p(\VN(G))$ associated with group von Neumann algebras $\VN(G)$ where $1<p<\infty$ and where $G$ is a discrete group. We refer to the papers \cite{CGIS1}, \cite[Definition 10.4]{Cip1}, \cite{HMP1}, \cite{HMP2}, \cite{MaN1}, \cite{McM1}and \cite{NeV1} for Hodge-Dirac operators in related contexts.

Recall that if $G$ is a discrete group then the von Neumann algebra $\VN(G)$, whose elements are bounded operators acting on the Hilbert space $\ell^2_G$, is generated by the left translation unitaries $\lambda_s \co \ell^2_G \to \ell^2_G$, $\delta_r \mapsto \delta_{sr}$ where $r,s \in G$. If $G$ is abelian, then $\VN(G)$ is $*$-isomorphic to the algebra $\L^\infty(\hat{G})$ of essentially bounded functions on the dual group $\hat{G}$ of $G$. Recall that we can see the functions of $\L^\infty(\hat{G})$ as multiplication operators on $\L^2(\hat{G})$. As basic models of quantum groups, these von Neumann algebras play a fundamental role in operator algebras. Moreover, we can equip $\VN(G)$ with a normalized trace (=noncommutative integral) and if $1 \leq p \leq \infty$ we have a canonical identification
$$
\L^p(\VN(G))
=\L^p(\hat{G}).
$$
A Markov semigroup $(T_t)_{t \geq 0}$ of Fourier multipliers on $\VN(G)$ is characterized by a conditionally negative length $\psi \co G \to \C$ such that the symbol of each operator $T_t$ of the semigroup is $\e^{-t\psi}$. Moreover, the symbol of the (negative) infinitesimal generator $A_p$ on the noncommutative $\L^p$-space $\L^p(\VN(G))$ of the semigroup is $\psi$. 

Introducing the Banach space $\L^p(\VN(G)) \oplus_p \Omega_{\psi,1,p}$, the authors of \cite{JMP2} define the Hodge-Dirac operator
\begin{equation}
\label{Def-Dirac-operator-Fourier-intro}
D_{\psi,1,p} 
\ov{\mathrm{def}}{=}\begin{bmatrix} 
0 & (\partial_{\psi,1,p})^* \\ 
\partial_{\psi,1,p}& 0 
\end{bmatrix}
\end{equation}
which is an unbounded operator defined on a dense subspace. In \cite[Problem C.5]{JMP2}, the authors ask for dimension free estimates for the operator $\sgn D_{\psi,1,p} \ov{\mathrm{def}}{=} D_{\psi,1,p}|D_{\psi,1,p}|^{-1}$. We affirmatively answer this question for a large class of groups including all amenable discrete groups and free groups by showing the following result in the spirit of \cite{AKM1} (see also \cite{AAM1}).

\begin{thm}[see Theorem \ref{Th-functional-calculus-bisector-Fourier}, Theorem \ref{prop-Hodge-Fourier-HI-on-Omega-psi} and Remark \ref{rem-dimension-free-Hodge-Fourier}]
\label{Th-functional-calculus-bisector-Fourier-intro}
Suppose $1<p<\infty$. Let $G$ be a weakly amenable discrete group such that the crossed product $\L^\infty(\Omega) \rtimes_\alpha G$ is $\QWEP$. The Hodge-Dirac operator $D_{\psi,1,p}$ is bisectorial on $\L^p(\VN(G)) \oplus_p \Omega_{\psi,1,p}$ and admits a bounded $\HI(\Sigma^\pm_\omega)$ functional calculus on a bisector $\Sigma^\pm_\omega$. Moreover, the norm of the functional calculus is bounded by a constant $K_\omega$ which depends neither on $G$ nor on the semigroup\footnote{\thefootnote. In particular it is independent of the dimension of the Hilbet space $H$ associated to the 1-cocycle by Proposition \ref{prop-Schoenberg}.}. 
\end{thm}

We refer to \cite{HvNVW2} for more information on bisectorial operators. Our result can be seen as a strengthening of the dimension free estimates \eqref{Riesz-Parcet-commutatif} of Riesz transforms of the above authors since it is almost immediate that this bounded functional calculus implies the equivalence \eqref{Riesz-Parcet-commutatif} see Remark \ref{remark-sgn-Fourier}. Our argument relies in part on a new transference argument between Fourier multipliers on crossed products and classical Fourier multipliers (see Proposition \ref{prop-Fourier-mult-crossed-product}) which is of independent interest. With the help of an extension of this result (Theorem \ref{Thm-full-operator-bisectorial-Fourier}), we obtain a Hodge decomposition (see Theorem \ref{Th-Hodge-decomposition-Fourier}).
In particular, we are able to deduce functional calculus for an extension of $D_{\psi,1,p}$ on the whole space $\L^p(\VN(G)) \oplus_p \L^p(\L^\infty(\Omega) \rtimes_\alpha G)$.
Note that the QWEP assumption is satisfied for amenable groups by \cite[Proposition 4.1]{Oza} and for free groups $\F_n$ by the same reasoning used in the proof of \cite[Proposition 4.8]{Arh1}.

Below, we describe concrete semigroups in which Theorem \ref{Th-functional-calculus-bisector-Fourier-intro} applies.
	
\paragraph{Semigroups on abelian groups}
Recall that a particular case of \cite[Corollary 18.20]{BeF1} says that a function $\psi \co G \to \R$ on a discrete abelian group $G$ is a conditionally negative length if and only if there exists a quadratic form\footnote{\thefootnote. That means that $2q(s)+2q(t)=q(s+t)+q(s-t)$ for any $s,t \in G$.} $q \co G \to \R^+$ and a symmetric positive measure $\mu$ on $\hat{G}-\{0\}$ such that $\int_{\hat{G}-\{0\}} \big(1-\Re \chi(s)\big) \d\mu(\chi)<\infty$ for any $s \in G$ satisfying 
\begin{equation}
\label{Levy-Khint}
\psi(s)
=q(s)+\int_{\hat{G}-\{0\}} \big(1-\Re \chi(s)\big) \d\mu(\chi),\quad s \in G.
\end{equation}
In this case, $\mu$ is the so called L\'evy measure of $\psi$ and $q$ is determined by the formula $q(s)=\lim_{n \to +\infty} \frac{\psi(ns)}{n^2}$. This is the L\'evy-Khinchin representation of $\psi$ as a continuous sum of elementary conditionally negative lengths\footnote{\thefootnote. Recall that a quadratic form $q \co G \to \R^+$ and a function $G \to \R$, $s \mapsto 1-\Re \chi(s)$ are conditionally negative lengths by \cite[Proposition 7.19]{BeF1}, \cite[Proposition 7.4 (ii)]{BeF1} and \cite[Corollary 7.7]{BeF1}.}. 

\begin{enumerate}
\item[(a)] If $G=\Z^n$, we recover the Markov semigroups on the $\L^p$-spaces $\L^p(\VN(\Z^n))=\L^p(\T^n)$ of the torus $\T^n$. For example, taking $\mu=0$ and $\psi(k_1,\ldots,k_n)=q(k_1,\ldots,k_n)=k_1^2+\cdots+ k_n^2$, we obtain the function defining the heat semigroup $(\e^{t \Delta})_{t \geq 0}$. By choosing $q=0$ and the right measure $\mu$, we can obtain the Poisson semigroup $(\e^{-t(-\Delta)^{\frac{1}{2}}})_{t \geq 0}$ or more generally the semigroups associated to the fractional laplacians $(-\Delta)^\alpha$ with $0 \leq \alpha <2$ \cite{RoS1}, \cite{RoS2}.
Note that in the particular case of the Poisson semigroup on the $\L^p$-space $\L^p(\T)$ associated with the torus $\T$, the L\'evy measure is given by 
\[ 
\d\mu(\e^{\i x}) 
= \Re \left[ \frac{-2 \e^{\i x}}{(1-\e^{\i x})^2} \right] \d \mu_\T(\e^{\i x}),
\]
where $\mu_\T$ denotes the normalized Haar measure on the torus.



\item[(b)] 	
Fix some integer $n \geq 1$. We consider the group $G$ of Walsh functions defined by 
$w_A(\epsi)=\prod_{i \in A} \epsi_i$ where $A$ is a subset of $\{1,\ldots,n\}$ and where $\epsi=(\epsi_1,\ldots,\epsi_n)$ belongs to the discrete abelian group $\hat{G} = \{-1,1\}^n$. For any $1 \leq j \leq n$, we let $\sigma_j \ov{\mathrm{def}}{=} (1,\ldots,1,-1,1,\ldots,1)$. If we consider the atomic measure $\mu=\frac{1}{2}\sum_{j=1}^{n} \delta_{\sigma_j}$ on $\hat{G}-\{(1,\ldots,1)\}$ and $q=0$, we obtain\footnote{\thefootnote. For any $1 \leq j \leq n$, note that $w_A(\sigma_j)=\prod_{i \in A} (\sigma_j)_i$ is equal to $1$ if $j \not\in A$ and to $-1$ if $j \in A$. So we have
\begin{align*}
\MoveEqLeft
\psi(w_A)            
\ov{\eqref{Levy-Khint}}{=} \frac{1}{2}\sum_{j=1}^{m}  \big(1-w_A(\sigma_j)\big) 
=\frac{1}{2} \bigg(m-\sum_{j=1}^{m} w_A(\sigma_j)\bigg)
=\card A.
\end{align*}} 
$\psi(w_A)= \card A$. So we recover the discrete Heat semigroup of \cite[p.~19]{HvNVW2}. It is also related to \cite{EfL1}, \cite{Lus1}, \cite{Lus2} and \cite{Nao1}.

\end{enumerate}

\paragraph{Semigroups on finitely generated groups} Let $G$ be a finitely generated group and $S$ be a generating set for $G$ such that $S^{-1}=S$ and $e \not\in S$. Any element $s$ admits a decomposition
\begin{equation}
\label{rep-reduced}
s
=s_1 s_2\cdots s_{n}
\end{equation}
where $s_1, \ldots, s_n$ are elements of $S$. The word length $|s|$ of $s$ with respect to the generating set $S$ is defined to be the minimal integer $n$ of such a decomposition and is a basic notion in geometric group theory. As a special case, the neutral element $e$ has length zero.

\begin{enumerate}
	\item[(a)]\textbf{Coxeter groups} Here, we refer to \cite{BGM1} and references therein for more information. Recall that a group $G=W$ is called a Coxeter group if $W$ admits the following presentation:
$$
W
=\left\langle S \, \big|\, (s_1 s_2)^{m(s_1,s_2)}=1:s_1,s_2\in S,m(s_1,s_2)\ne\infty\right\rangle
$$
where $m \co S \times S \to \{1,2,3,\ldots,\infty\}$ is a function such that $m(s_1,s_2)=m(s_2,s_1)$ for any $s_1,s_2\in S$ and $m(s_1,s_2)=1$ if and only if $s_1=s_2$. The pair $(W,S)$ is called a Coxeter system. In particular, every generator $s \in S$ has order two. By \cite[Theorem 7.3.3]{Boz2}, the word length $|s|$ is a conditionally negative length and our results can be used with the semigroup generated by this function. Recall that dihedral groups $\mathrm{D}_n=\langle r,s \mid r^n=s^2=(sr)^2=1\rangle$, $\Z_2 \times \cdots \times \Z_2$, symmetric groups $\mathrm{S}_{n}$ with $S=\{(n,n+1) : n \in \N\}$ and the group $\mathrm{S}_\infty$ of all finite permutations of the set $\N$, called the infinite symmetric group, with $S=\{(n,n+1) : n \in \N\}$ are examples of Coxeter groups. In the case of symmetric groups, the length $|\sigma|$ of $\sigma$ is the number of crossings in the diagram which represents the permutation $\sigma$. 

Consider a Coxeter group $W$. If $s \in W$, note that the sequence $s_1,\ldots,s_n$ in \eqref{rep-reduced} chosen in such a way that $n$ is minimal is not unique in general. However, the set of involved generators is unique, i.e. if $s=s_1 s_2\cdots s_n=s_1' s_2'\cdots s_n'$ are minimal words of $s \in W$ then $\{s_1,s_2,\ldots,s_n\}=\{s_1',s_2',\ldots,s_n'\}$. This subset $\{s_1,s_2,\ldots,s_n\}$ of $S$ is denoted $S_s$ and is called the colour of $s$, following \cite[p.~585]{BGM1}. We define the colour-length of $s$ putting $\norm{s} \ov{\mathrm{def}}{=}\card S_s$. We always have $\norm{w} \leq |w|$. By \cite[Theorem 4.3 and Corollary 5.4]{BGM1}, if $0 \leq \alpha \leq 1$ the functions $|\cdot|^\alpha$ and $\norm{\cdot}$ are conditionally negative lengths on $\mathrm{S}_\infty$. Finally, see also \cite[p.~1971]{JMP1} for other examples for $\mathrm{S}_n$ for $n<\infty$.


\item[(b)] \textbf{Free groups} Our results can be used with the noncommutative Poisson semigroup \cite{Haag1}, \cite[Definition 10.1]{JMX} on free groups $\mathbb{F}_n$ ($1 \leq n \leq \infty$) whose generator is the length $|\cdot|$. Moreover, we can use the characterization \cite[Theorem 1.2]{HaK1} of radial functions $\psi \co \F_n \to \C$ with $\psi(e)=0$ which are conditionally negative definite. If $\varphi_z \co \F_n \to \C$ denotes the spherical function\footnote{\thefootnote. In the case $n=\infty$, we have $\varphi_z(x)=z^{|x|}$.} of parameter $z \in \C$, these functions can be written $\psi(s) = \int_{-1}^1 \psi_z(s) \d\nu(z)$, for $s \in \F_n$ and for some finite positive Borel measure $\nu$ on $[-1,1]$ where $\psi_z(s) \ov{\mathrm{def}}{=} \frac{1 - \varphi_z(s)}{1-z}$, $z \in \C \setminus\{1\}$ and $\psi_1(x) \ov{\mathrm{def}}{=} \lim_{z \to 1} \frac{1-\varphi_z(s)}{1-z}$. Finally, \cite{JPPP1} contains (but without proof) examples of weighted forms of the word length which are conditionally negative.

We can write $\F_n= \ast_{i=1}^n G_i$ with $G_i=\Z$. Every element $s$ of $\F_n - \{ e \}$ has a unique representation $s = s_{i_1}s_{i_2} \cdots s_{i_m}$, where $s_{i_k} \in G_{i_k}$ are distinct from the corresponding neutral elements and $i_1 \not= i_2 \not= \cdots\not= i_k$. The number $m$ is called the block length of $s$ and denoted $\norm{s}$. By \cite[Example 6.14]{JuZ3}, this function is a conditionally negative length. 


\item[(c)]\textbf{Cyclic groups} The word length on $\Z_n$ is given by $|k|=\min\{n,n-k\}$. It is known that this length is negative definite, see for example \cite[Appendix B]{JPPP1}, \cite[Example 5.9]{JuZ1}, \cite[Section 5.3]{JuZ3} and \cite[p.~553]{JMP2} for more information. By \cite[p.~925]{JRZ1}, the function $\psi_n$ defined on $\Z_n$ by $\psi_n(k)=\frac{n^2}{2\pi^2}\big(1-\cos(\frac{2\pi k}{n})\big)$ is another example of conditionally negative length on $\Z_n$.
\end{enumerate}

\paragraph{Semigroups on the discrete Heisenberg group}
Let $\mathbb{H}=\Z^{2n+1}$ be the discrete Heisenberg group with group operations
\begin{equation}
\label{Operations-Heisenberg}
(a,b,t)\cdot(a',b',t') 
=(a+a',b+b',t+t'+a b')
\quad \text{and} \quad
(a,b,t)^{-1}=(-a,-b,-t+a b)	
\end{equation}
where $a,b,a',b' \in \Z^n$ and $t,t' \in \Z$. By \cite[Proposition 5.13]{JuZ3} and \cite[p.~261]{JuZ1}, the map $\psi \co \mathbb{H} \to \R$, $(a,b,t) \mapsto |a|+|b|$ is a conditionally negative length.

\vspace{0.2cm}

We also prove in this paper an analogue of the equivalences \eqref{Riesz-Parcet-commutatif} for markovian semigroups $(T_t)_{t \geq 0}$ of Schur multipliers acting on Schatten spaces $S^p_I \ov{\mathrm{def}}{=}S^p(\ell^2_I)$ for $1 < p < \infty$ where $I$ is an index set. In this case, by \cite[Proposition 5.4]{Arh1}, the Schur multiplier symbol $[a_{ij}]$ of the negative generator $A$ of $(T_t)_{t \geq 0}$ is given by $a_{ij} = \norm{ \alpha_i  - \alpha_j }^2_H$ for some family $\alpha = (\alpha_i)_{i \in I}$ of vectors of a real Hilbert space $H$. We define a gradient operator $\partial_{\alpha,1,p}$ as the closure of the unbounded linear operator $\M_{I,\fin} \to \L^p(\Omega,S^p_I)$, $e_{ij} \mapsto \W(\alpha_i-\alpha_j) \ot e_{ij}$ where $\W \co H \to \L^0(\Omega)$ is a $H$-isonormal Gaussian process, $\Omega$ is the associated probability space and where $\M_{I,\fin}$ is the subspace of $S^p_I$ of matrices with a finite number of non null entries. Then the result reads as follows.

\begin{thm}[see Theorem \ref{Thm-Riesz-equivalence-Schur} and \eqref{Equivalence-square-root-domaine-Schur}]
\label{Thm-Riesz-equivalence-Schur-intro}
Let $I$ be an index set and $A$ be the generator of a markovian semigroup $(T_t)_{t \geq 0}$ of Schur multipliers on $\B(\ell^2_I)$. Suppose $1<p<\infty$. For any $x \in \M_{I,\fin}$, we have
\begin{equation} 
\label{OneLineDimFree-intro}
\bnorm{A_p^{\frac12}(x)}_{S^p_I} 
\approx_{p} \norm{\partial_{\alpha,1,p}(x)}_{\L^p(\Omega,S^p_I)}. 
\end{equation}
\end{thm}

\noindent We also obtain an analogue of Theorem \ref{Th-functional-calculus-bisector-Fourier-intro}. With this result, we are equally able to obtain a Hodge decomposition, see Theorem \ref{Th-Hodge-decomposition}.

\begin{thm}[see Theorems \ref{Thm-full-operator-bisectorial} and \ref{prop-explicitly-dimension-free}]
\label{Thm-full-operator-bisectorial-intro}
Suppose $1 < p < \infty$. The unbounded operator $\D_{\alpha,1,p} \ov{\mathrm{def}}{=} \begin{bmatrix} 
0 & (\partial_{\alpha,1,p^*})^* \\ 
\partial_{\alpha,1,p} & 0 
\end{bmatrix}$ on the Banach space $S^p_I \oplus \L^p(\Omega,S^p_I)$ is bisectorial and has a bounded $\HI(\Sigma_\omega^\pm)$ functional calculus on a bisector $\Sigma_\omega^\pm$. Moreover, the norm of the functional calculus is bounded by a constant $K_\omega$ which depends neither on $I$ nor on the semigroup. In particular it is independent of the dimension of $H$.
\end{thm}

Moreover, we also relate the equivalences \eqref{Riesz-Parcet-commutatif} and \eqref{OneLineDimFree-intro} with the ones of Meyer's formulation \eqref{Meyer-estimates}. To achieve this, we define and study in the spirit of \eqref{carre-du-champ} a carr\'e du champ $\Gamma$ (see \eqref{Def-Gamma-sgrp} and \eqref{Def-Gamma}) and its closed extension in the sense of Definition \ref{defi-closure} and we connect this notion to some approximation properties of groups. It leads us to obtain alternative formulations of \eqref{Riesz-Parcet-commutatif} and \eqref{OneLineDimFree-intro}. Note that some carr\'e du champ were studied in the papers \cite[Section 9]{CiS1}, \cite{Cip1}, \cite{JM1}, \cite{Sau1} and \cite{JuZ3} mainly in the $\sigma$-finite case and for $\L^2$-spaces (see \cite{DaL1} for related things) but unfortunately their approach does not suffice for our work on $\L^p$-spaces. By the way, it is rather surprising that even in the \textit{commutative} setting, no one has examined the carr\'e du champ on $\L^p$-spaces with $p \not=2$. The following is an example of result that we have achieved, and which can be compared with \eqref{Meyer-estimates}. 

\begin{thm}[see Theorem \ref{cor-Riesz-equivalence-Schur}]
\label{cor-Riesz-equivalence-Schur-MIfin-intro}
Suppose $2 \leq p < \infty$. Let $A$ be the generator of a markovian semigroup of Schur multipliers on $\B(\ell^2_I)$. For any $x \in \dom A_p^{\frac12}$, we have
\begin{equation}
\label{max-Gamma-2-intro}
\bnorm{A_p^{\frac12}(x)}_{S^p_I} 
\approx_p 
\max \left\{ \bnorm{\Gamma(x,x)^{\frac12}}_{S^p_I}, \bnorm{\Gamma(x^*,x^*)^{\frac12}}_{S^p_I} \right\}.
\end{equation}		
\end{thm}
The maximum is natural in noncommutative analysis due to the use of noncommutative Khintchine inequalities.

It is remarkable that the point of view of Hodge-Dirac operators fits perfectly in the setting of noncommutative geometry if $p=2$. If $G$ is a discrete group, the Hilbert space $H \ov{\mathrm{def}}{=} \L^2(\VN(G)) \oplus_2 \Omega_{\psi,1,2}$, the $*$-algebra $A \ov{\mathrm{def}}{=} \VNGfin$ of trigonometric polynomials and the Hodge-Dirac operator $D_{\psi,1,2}$ on $\L^2(\VN(G)) \oplus_2 \Omega_{\psi,1,2}$ define a triple $(A,H,D_{\psi,1,2})$ in the spirit of noncommutative geometry \cite{Con3}, \cite{GVF1}, \cite{Var1}. Recall that the notion of spectral triple $(A,H,D)$ (= noncommutative manifold) \`a la Connes covers a huge variety of different geometries such as Riemmannian manifolds, fractals, quantum groups or even non-Hausdorff spaces. We refer to \cite{CoM1} for an extensive list of examples and to \cite{CCHKM}, \cite{Con5}, \cite{Haj1}, \cite{Lan1} and \cite{Sit1} for some surveys. From here, it is apparent that we can see Markov semigroups of Fourier multipliers as geometric objects. The same observation is true for Markov semigroups of Schur multipliers. Nevertheless, the Hilbert space setting of the noncommutative geometry is too narrow to encompass our setting on $\L^p$-spaces. So, we develop in Subsection \ref{Sec-Banach-spectral-triples-first} a natural Banach space variant $(A,X,D)$ of a spectral triple where the selfadjoint operator $D$ acting on the Hilbert space $H$ is replaced by a bisectorial operator $D$ acting on a (reflexive) Banach space $X$, allowing us to use (noncommutative) $\L^p$-spaces ($1<p<\infty$). 

It is well-known that Gaussian variables are not bounded, i.e. do not belong to $\L^\infty(\Omega)$. From the perspective of noncommutative geometry, it is problematic under technical aspects as the boundedness of the commutators $[D,\pi(a)]$ of a spectral triple $(A,H,D)$. Indeed, the noncommutative gradients \eqref{def-partial-psi-very-particular} appear naturally in the commutators of our spectral triples and these gradients are defined with Gaussian variables. Fortunately, the noncommutative setting is very flexible and allows us to introduce a continuum of gradients $\partial_{\psi,q,p}$ and $\partial_{\alpha,q,p}$ indexed by a new parameter $-1 \leq q \leq 1$ replacing Gaussian variables by \textit{bounded} noncommutative $q$-deformed Gaussian variables and $\L^\infty(\Omega)$ by the von Neumann algebra $\Gamma_q(H)$ of \cite{BoS} and \cite{BKS}. Note that $\Gamma_{-1}(H)$ is the fermion algebra and that $\L^\infty(\Omega)$ can be identified with the boson algebra $\Gamma_{1}(H)$. Our main theorems on Hodge-Dirac operators admit extensions in these cases, see Theorem \ref{Th-functional-calculus-bisector-Fourier} and Theorem \ref{prop-Hodge-Fourier-HI-on-Omega-psi}. We expect some differences of behaviour when $q$ varies. 

In Section \ref{new-quantum-metric-spaces}, we give a first elementary study of the triples associated to Markov semigroups of multipliers. In particular, we give sufficient conditions for the verification of axioms of noncommutative geometry and more generally the axioms of our Banach spectral triples. One of our main results in this part reads as follows and can be used with arbitrary amenable groups or free groups. Here $\mathrm{C}^*_r(G)$ denotes the reduced $\mathrm{C}^*$-algebra of the group $G$.

\begin{thm}[see Theorem \ref{First-spectral-triple-2}]
\label{First-spectral-triple-2-intro}
Suppose $1<p<\infty$ and $-1 \leq q < 1$. Let $G$ be a discrete group. Let $(T_t)_{t \geq 0}$ be a markovian semigroup of Fourier multipliers on the group von Neumann algebra $\VN(G)$. Consider the associated function $b_\psi \co G \to H$ from \eqref{liens-psi-bpsi}. Assume that the Hilbert space $H$ is finite-dimensional, that $b_\psi$ is injective and that
\[ 
\G_\psi 
\ov{\mathrm{def}}{=} \inf_{b_\psi(s) \neq b_\psi(t)} \norm{b_\psi(s) - b_\psi(t)}_H^2 > 0.
\]
Finally assume that the von Neumann crossed product $\Gamma_q(H) \rtimes_\alpha G$ has $\QWEP$. Let 
\[
\pi \co \mathrm{C}^*_r(G) \to \B(\L^p(\VN(G)) \oplus_p \L^p(\Gamma_q(H) \rtimes_\alpha G))
\]
be the Banach algebra homomorphism such that $\pi(a)(x,y) \ov{\mathrm{def}}{=} (ax,(1 \rtimes a)y)$ where $x \in \L^p(\VN(G))$ and $y \in \L^p(\Gamma_q(H) \rtimes_\alpha G)$. Let $\D_{\psi,q,p}$ be the Hodge-Dirac operator on the Banach space $\L^p(\VN(G)) \oplus_p \L^p(\Gamma_q(H) \rtimes_\alpha G)$ defined in \eqref{Def-Dirac-operator-Fourier-2}. Then $(\mathrm{C}^*_r(G),\L^p(\VN(G)) \oplus_p \L^p(\Gamma_q(H) \rtimes_\alpha G), \D_{\psi,q,p})$ is a Banach spectral triple in the sense of Definition \ref{Def-Banach-spectral-triple}.
\end{thm}

We are equally interested by the metric aspect \cite{Lat2} of noncommutative geometry, see Subsection \ref{subsec-prelims-quantum} for background on quantum (locally) compact metric spaces. We introduce new quantum (locally) compact metric spaces in the sense of \cite{Lat1}, \cite{Lat2}, \cite{Rie3} associated with these spectral triples. It relies on $\L^p$-variants of the seminorms of \cite[Section 1.2]{JM1}. Here, we check carefully the axioms taking into account all problems of domains required by this theory. Note that it is not clear how to do the same analysis at the level $p=\infty$ considered in \cite[Section 1.2]{JM1} since we cannot hope the boundedness of Riesz transform on $\L^\infty$ which is an important tool in this part. So, on this point, $\L^p$-seminorms seem more natural.   

We observe significant differences between the case of Fourier multipliers (Theorem \ref{thm-Fourier-quantum-compact-metric}) and the one of Schur multipliers (Theorem \ref{th-quantum-metric-space}). For example, we need to use \textit{$\mathrm{C}^*$-algebras} for semigroups of Fourier multipliers and which produces quantum \textit{compact} metric spaces contrarily to the case of semigroups of Schur multipliers which requires \textit{order-unit spaces} and that leads to quantum \textit{locally} compact metric spaces if $I$ is infinite. Furthermore, our analysis with quantum compact metric spaces relying on noncommutative $\L^p$-spaces makes appear a new phenomenon when the value of the parameter $p$ changes, see Remark \ref{rem-Fourier-quantum-compact-metric-exponent-restriction}. 

Finally, note that the combination of our spectral triples and our quantum (locally) compact metric spaces is in the spirit of the papers \cite{Lat8} and \cite{BMR1} (see also \cite{Con4}) but it is more subtle here since the link between the norms of the commutators and the seminorms of our quantum metric spaces is not as direct as the ones of \cite{Lat8} and \cite{BMR1}.

Finally, we are interested in other Dirac operators related to Riesz transforms. We study the associated noncommutative geometries in Subsection \ref{Sec-spectral-triples-Fourier-II} and Subsection \ref{subsec-new-spectral-triples-Schur-II} and their properties of functional calculus in Subsection \ref{Sec-functional-calculus}.

The paper is organized as follows. 

Section \ref{Preliminaries} gives background and preliminary results. In Subsection \ref{various}, we collect some elementary information on operator theory, functional calculus and semigroup theory. In Subsection \ref{Sec-q-gaussians}, we recall some important information on isonormal Gaussian processes and more generally $q$-Gaussian functors. These notions are fundamental for the construction of the noncommutative gradients. In Subsection \ref{Sec-sesquilinear-forms}, we develop a short theory of vector-valued unbounded bilinear forms on Banach spaces that we need for our carr\'e du champ with values in a noncommutative $\L^p$-space. In Subsection \ref{subsubsec-Fourier-mult-crossed-product}, we introduce a new transference result (Proposition \ref{prop-Fourier-mult-crossed-product}) between Fourier multipliers on crossed products and classical Fourier multipliers. In Subsection \ref{Sec-Hilbertian-valued}, we discuss Hilbertian valued noncommutative $\L^p$-spaces which are fundamental for us. We complete and clarify some technical points of the literature. In Subsection \ref{Sec-infos-on-Gamma-Fourier}, we introduce the carr\'e du champ $\Gamma$ and the gradients with value in a noncommutative $\L^p$-space for semigroups of Fourier multipliers.
In Subsection \ref{Sec-infos-on-Gamma}, we discuss the carr\'e du champ $\Gamma$ and the gradients associated to Markov semigroups of Schur multipliers.

In Section \ref{sec:Semigroups-selfadjoint-contractive-Schur-multipliers}, we investigate dimension free Riesz estimates for semigroups of Schur multipliers and we complement the results of \cite{JMP2} for Riesz transforms associated to semigroups of Markov Fourier multipliers. In Subsection \ref{twisted-Khintchine}, we obtain Khintchine inequalities for $q$-Gaussians in crossed products generalising some result of $\cite{JMP2}$. In Subsection \ref{sec-Fourier-Kato}, we extend the equivalences \eqref{Riesz-Parcet-commutatif} of \cite{JMP2} to $q$-deformed Gaussians variables on a larger (maximal) domain, see \eqref{Equivalence-square-root-domaine-Schur-sgrp}. We also examine the constants. In Subsection \ref{Sec-infos-on-Gamma-Fourier2}, we are interested by Meyer's equivalences \eqref{Meyer-estimates} in the context of semigroups of Fourier multipliers in the case $p \geq 2$, see Theorem \ref{cor-Riesz-equivalence-group-2}.
In Subsection \ref{Kato-Schur}, we show dimension free estimates for Riesz transforms associated to semigroups of Markov Schur multipliers, see Theorem \ref{Thm-Riesz-equivalence-Schur} and \eqref{Equivalence-square-root-domaine-Schur}. Moreover, we examine carefully the obtained constants in this equivalence. In Subsection \ref{Riesz-transforms-cocycle-Schur}, we turn again to the formulation of these equivalences in the spirit of Meyer's equivalence \eqref{Meyer-estimates} in connection with the carr\'e du champ. We also obtain concrete equivalences similar to the ones of Lust-Piquard \cite{Lus2}, \cite{Lus1} and related to \eqref{ine-intro-6}. Note that for the case $1 < p < 2$, the statement becomes more involved, see Theorem \ref{cor-Riesz-equivalence-Schur-MIfin}.

In the next Section \ref{Sec-Hodge-Dirac}, we show our main results on functional calculus of Hodge-Dirac operators. In Subsection \ref{subsec-HI-calculus-Hodge-Dirac-Fourier}, we start with the case of the Hodge-Dirac operator for semigroups of Fourier multipliers. We show the boundedness of the functional calculus on $\L^p(\VN(G)) \oplus_p \ovl{\Ran \partial_{\psi,q,p}}$, see Theorems \ref{Th-D-R-bisectorial-Fourier} and \ref{Th-functional-calculus-bisector-Fourier}. In Section \ref{Extension-full-Hodge-Fourier}, we extend this boundedness to the larger space $\L^p(\VN(G)) \oplus \L^p(\Gamma_q(H) \rtimes_\alpha G)$. In Subsection \ref{Sec-Omega-psi}, we obtain the boundedness on  $\L^p(\VN(G)) \oplus \Omega_{\psi,q,p}$. In Subsection \ref{Module-Omega}, we examine the bimodule properties of $\Omega_{\psi,q,p}$.
In Subsection \ref{Hodge-Schur}, we change the setting and we consider the Hodge-Dirac operator associated to semigroups of Schur multipliers and we show the boundedness of the functional calculus on $S^p_I \oplus_p \ovl{\Ran \partial_{\alpha,q,p}}$, see Theorem \ref{Th-functional-calculus-bisector-Schur}. In a second step, we extend in Subsection \ref{Extension-full-Hodge-Schur} this boundedness on $S^p_I \oplus \L^p(\Gamma_q(H) \otvn \B(\ell^2_I))$ in Corollary \ref{cor-key-Hodge-decomposition} and Theorem \ref{Thm-full-operator-bisectorial}.
In Subsection \ref{Sec-Indep-alpha}, we show that the constants are independent from $H$ and $\alpha$.

In Section \ref{new-quantum-metric-spaces}, we examine the noncommutative geometries induced by the Hodge-Dirac operators in both settings and another two Hodge-Dirac operators. In Subsection \ref{subsec-prelims-quantum}, we give background on quantum (locally) compact metric spaces. In Subsection \ref{subsec-new-compact-quantum-metric-spaces-Fourier}, we introduce $\L^p$-versions of the quantum compact metric spaces of \cite{JM1}. In Subsection \ref{Sec-some-estimates}, we establish fundamental estimates for the sequel and we introduce the constant $\G_\alpha$ which plays an important role. In Subsection \ref{Section-Leibniz-Schur}, we introduce the seminorm which allows us to define the new quantum (locally) compact metric spaces in the next Subsection \ref{Sec-I-finite}. In Subsection \ref{subsec-calculs-gap}, we compute the gaps of some explicit semigroups and we compare the notions in both of our settings (Fourier multipliers and Schur multipliers), see Proposition \ref{prop-lien-Gap-Fourier-Herz-Schur}. In Subsection \ref{Sec-Banach-spectral-triples-first}, we define a Banach space generalization of the notion of spectral triple. In Subsection \ref{Sec-spectral-triples-Fourier-I}, we give sufficient conditions in order that the triple induced by the Hodge-Dirac operator for a semigroup of Fourier multipliers give rise to Banach spectral triple satisfying axioms of noncommutative geometry and more generally the axioms of Subsection \ref{Sec-Banach-spectral-triples-first}. In Subsection \ref{Sec-spectral-triples-Fourier-II}, we introduce a second Hodge-Dirac operator for Fourier multipliers and we study the noncommutative geometry induced by this operator. In Subsection \ref{Sec-new-spectral-triples-first}, we give sufficient conditions in order that the triple induced by the Hodge-Dirac operator for a semigroup of Schur multipliers give rise to a spectral triple. In Subsection \ref{subsec-new-spectral-triples-Schur-II}, we introduce a second related Hodge-Dirac operator and we study the noncommutative geometry induced by this operator. In Subsection \ref{Sec-functional-calculus}, in the bosonic case $q=1$ we establish the bisectoriality and the boundedness of the functional calculus of Dirac operators introduced in Subsection \ref{Sec-spectral-triples-Fourier-II} and Subsection \ref{subsec-new-spectral-triples-Schur-II}.

Finally, in the short Section \ref{Sec-Levy-and-cocycles}, we discuss how the L\'evy measure of a continuous conditionally of negative type function $\psi$ on a locally compact abelian group $G$ induces a 1-cocycle.



\section{Preliminaries}
\label{Preliminaries}
\subsection{Operators, functional calculus and semigroups}
\label{various}

In this short subsection, we collect various objects that play a role in this article at different places.
\paragraph{Closed and weak* closed operators}
An operator $T \co \dom T \subseteq X \to Y$ is closed if 
\begin{align}
\label{Def-operateur-ferme}
\MoveEqLeft
\text{for any sequence $(x_n)$ of $\dom T$ with $x_n \to x$ and $T(x_n) \to y$ with $x \in X$ and $y \in Y$ }\\
&\text{implies that $x \in \dom T$ and $T(x)=y$.}   \nonumber         
\end{align}  
A linear subspace $C$ of $\dom T$ is a core of $T$ if $C$ is dense in $\dom T$ for the graph norm, that is 
\begin{equation}
\label{Def-core}
\text{for any $x \in \dom T$ there is a sequence $(x_n)$ of $C$ s. t. $x_n \to x$ in $X$ and $T(x_n) \to T(x)$ in $Y$.}
\end{equation}
If $T$ is closed, a subspace $C$ of $\dom T$ is a core for $T$ if and only if $T$ is the closure of its restriction $T|C$. Recall that if $T \co \dom T \subseteq X \to Y$ is a densely defined unbounded operator then $\dom T^*$ is equal to 
\begin{equation}
\label{Def-domaine-adjoint}
\big\{ y^* \in Y^* : \text{there exists } x^* \in X^* \text{ such that } \langle T(x),y^* \rangle_{Y,Y^*}=\langle x,x^* \rangle_{X,X^*} \text{ for all } x \in \dom T\big\}.
\end{equation}
If $y^* \in \dom T^*$, the above $x^* \in X^*$ is determined uniquely by $y^*$ and we let $T^*(y^*)=x^*$. For any $x \in \dom T$ and any $y^* \in \dom T^*$, we have
\begin{equation}
\label{crochet-duality}
\langle T(x),y^* \rangle_{Y,Y^*}
=\langle x,T^*(y^*) \rangle_{X,X^*}.
\end{equation}
Recall that the unbounded operator $T \co \dom T \subseteq X \to Y$ is closable \cite[p.~165]{Kat1} if and only if 
\begin{equation}
\label{Def-closable}
x_n \in \dom T,\, x_n \to 0 \text{ and } T(x_n) \to y \text{ imply } y=0. 
\end{equation}
If the unbounded operator $T \co \dom T \subseteq X \to Y$ is closable then 
\begin{equation}
\label{Domain-closure}
x \in \dom \ovl{T} \text{ iff there exists } (x_n) \subseteq \dom T \text{ such that } x_n \to x
\text{ and } T(x_n) \to y \text{ for some } y.  
\end{equation}
In this case $\ovl{T}(x)=y$. Finally, if $T$ is a densely defined closable operator then $T^*=\ovl{T}^*$ by \cite[Problem 5.24]{Kat1}. If $T$ is densely defined, by \cite[Problem 5.27 p. 168]{Kat1}, we have
\begin{equation}
\label{lien-ker-image}
\ker T
=(\Ran T)^\perp.
\end{equation}
If $T$ is a densely defined unbounded operator and if $T \subseteq R$, by \cite[Problem 5.25 p. 168]{Kat1} we have
\begin{equation}
\label{Inclusion-adjoint-unbounded}
R^* \subseteq T^*.
\end{equation}
If $TR$ and $T$ are densely defined by \cite[Problem 5.26 p. 168]{Kat1} we have
\begin{equation}
\label{Adjoint-product-unbounded}
R^*T^*\subseteq (TR)^*.
\end{equation}
In Section \ref{new-quantum-metric-spaces}, we need the following lemma. To this end, an operator $T \co \dom T \subseteq X \to Y$ between dual spaces $X$ and $Y$ is called weak* closed if \eqref{Def-operateur-ferme} holds, when the convergences are to understand in the weak* sense.

\begin{lemma}
\label{lem-Krein-Smulian-2}
Let $X,Y$ be Banach spaces and $T \co \dom T \subseteq X^* \to Y^*$ be a weak* closed opeator with a weak*-core\footnote{\thefootnote. For any $a \in \dom T$ there exists a net $(a_j)$ in $D$ such that $a_j \to a$ and $T(a_j) \to T(a)$ both for the weak* topology.} $D \subseteq \dom T$. Then any $a \in \dom T$ admits a bounded net $(a_j)$ of elements of $D$ such that the net $(T(a_j))$ is also bounded, $a_j \to a$ and $T(a_j) \to T(a)$ both for the weak* topology.
\end{lemma}

\begin{proof}
We consider the Banach space $X^* \oplus Y^*$ which is canonically the dual space of $X \oplus Y$ and thus carries a weak* topology. Let $N = \{ (a,T(a)) : \: a \in D \} \subseteq X^* \oplus Y^*$ and $N_1  = N \cap B_1$ where $B_1$ is the closed unit ball of $X^* \oplus Y^*$. 

We claim that $\C \ovl{N_1}^{\mathrm{w}^*} = \ovl{N}^{\mathrm{w}^*}$. The inclusion ``$\subset$'' is clear.  For the inclusion ``$\supset$'' it suffices to show that $\C \ovl{N_1}^{\mathrm{w}^*}$ contains $N$ and is weak* closed. That it contains $N$ is easy to see. 

For the weak*-closedness, we will employ again the Krein-Smulian Theorem \cite[Theorem 2.7.11]{Meg1}. So it suffices to show that $\C \ovl{N_1}^{\mathrm{w}^*} \cap B_r$ is weak* closed for any $r > 0$. Let $(x_j)$ be a net in $\C \ovl{N_1}^{\mathrm{w}^*} \cap B_r$ converging to some $x \in X^* \oplus Y^*$. Since $B_r$ is weak* closed (by Alaoglu Theorem), we have $x \in B_r$. Moreover, for any $j$ we can write $x_j = \lambda_j y_j$ with $y_j \in \ovl{N_1}^{\mathrm{w}^*}$ and $\lambda_j \in \C$. We have a freedom in choice of both factors and can take $\norm{y_j}_{X^* \oplus Y^*} = 1$. Thus, $|\lambda_j| \leq r$. By Alaoglu Theorem, $(y_j)$ admits a weak* convergent subnet $(y_{j_k})$ with a limit thus in $\ovl{N_1}^{\mathrm{w}^*}$. Moreover, we can suppose that the net $(\lambda_{j_k})$ is convergent. Thus, $x_{j_k} = \lambda_{j_k} y_{j_k}$ converges for the weak* topology with limit a fortiori equal to $x$. Thus $x \in \C \overline{N_1}^{\mathrm{w}^*}$. We infer that $x \in \C \ovl{N_1}^{\mathrm{w}^*} \cap B_r$, so $\C \ovl{N_1}^{\mathrm{w}^*} \cap B_r$ is weak*-closed and hence $\C \overline{N_1}^{\mathrm{w}^*}$ is weak*-closed.  So we have shown $\C \ovl{N_1}^{\mathrm{w}^*} = \ovl{N}^{\mathrm{w}^*}$. 

Since $D$ is a weak*-core of $T$, the closure $\ovl{N}^{\mathrm{w}^*}$ equals the graph of $T$.  Now if $a \in \dom T$, the elements $(a,T(a))$ belongs to $\ovl{N}^{\mathrm{w}^*} = \C \ovl{N_1}^{\mathrm{w}^*}$, so that we can choose a net $(a_j)$ of graph norm less than $1$ and $\lambda \in \C$ such that $a_j \in D$, $\lambda a_j \to a$ and $T(\lambda a_j) \to T(a)$ both for the weak* topology.
\end{proof}


\paragraph{Sectorial operators, semigroups and functional calculus}
If $(T_t)_{t \geq 0}$ is a strongly continuous semigroup on a Banach space $X$ with (negative) infinitesimal generator $A$, we have
\begin{equation}
\label{semigroupe-vers-generateur}
\frac{1}{t}	\big(T_t(x)-x)
\xra[t \to 0^+]{} -A(x), \qquad x \in \dom A.
\end{equation}
by \cite[Definition 1.2 p.~49]{EnN1}. Moreover, by \cite[Corollary 5.5 p.~223]{EnN1}, for any $t \geq 0$, we have
\begin{equation}
\label{Widder}
\bigg(\Id+\frac{t}{n}A\bigg)^{-n}x
\xra[n \to +\infty]{} T_t(x), \quad x \in X.
\end{equation}
Moreover, by \cite[(1.5) p.~50]{EnN1} if $x \in \dom A$ and $t \geq 0$ then $T_t(x)$ belongs to $\dom A$ and
\begin{equation}
\label{A-et-Tt-commute}
T_tA(x)
=AT_t(x).
\end{equation}
We refer to \cite{HvNVW2} for the background on (bi-)sectorial operators and their $\HI$ functional calculus, as well as $R$-boundedness.
In particular, for $\sigma \in (0,\pi)$, we let $\Sigma_\sigma^+ = \Sigma_{\sigma} \ov{\mathrm{def}}{=} \{ z \in \C \backslash \{ 0 \} : \: | \arg z | < \sigma\}$.
A sectorial operator $A \co \dom A \subseteq X \to X$ on a Banach space $X$ is a closed densely defined operator such that $\C \backslash \ovl{\Sigma_{\sigma}}$ belongs to its resolvent set for some $\sigma \in (0,\pi)$ and such that $\sup\{\norm{z R(z,A)} :\: z \not\in \ovl{\Sigma_{\sigma'}} \} \leq C$ for some $\sigma' \geq \sigma$.
One writes 
\begin{equation}
\label{equ-sectorial-operator}
\omega(A) = \inf \sigma',
\end{equation}
where $\sigma'$ is as above.
A sectorial operator $A$ is called $R$-sectorial if for some $\omega(A)< \sigma< \pi$ the family
\begin{equation}
\label{Def-R-sectoriel}
\big\{zR(z,A) : z \not\in \ovl{\Sigma_\sigma}\big\}
\end{equation}
is $R$-bounded. For any $x \in X$, we have by e.g. \cite[(3.2)]{JMX}
\begin{equation}
\label{Resolvent-Laplace}
(z-A)^{-1}x
=-\int_{0}^{\infty} \e^{zt}T_t(x) \d t.
\end{equation}
In a similar manner, for $\sigma \in (0,\frac{\pi}{2})$, we let $\Sigma_\sigma^\pm \ov{\mathrm{def}}{=} \Sigma_\sigma \cup (-\Sigma_\sigma)$ and call an operator $A$ ($R$-)bisectorial, if $\big\{z R(z,A) : z \not\in \ovl{\Sigma_\sigma^\pm} \big\}$ is ($R$-)bounded for some $\sigma \in (0,\frac{\pi}{2})$. By \cite[p.~447]{HvNVW2}, a linear operator is bisectorial if and only if 
\begin{equation}
\label{Def-R-bisectoriel}
\i \R-\{0\} \subseteq \rho(A)
\quad \text{and} \quad
\sup_{t \in \R^+-\{0\}} \norm{tR(\i t,A)} 
< \infty.
\end{equation}
The following is a particular case of \cite[Proposition 2.3]{NeV1}, see also \cite[Theorem 10.6.7]{HvNVW2}.

\begin{prop}
\label{Prop-liens-differents-calculs-fonctionnels}
Suppose that $A$ is an R-bisectorial operator on a Banach space $X$ of finite cotype. Then $A^2$ is R-sectorial and for each $\omega \in (0,\frac{\pi}{2})$ the following assertions are equivalent.
\begin{enumerate}
	\item The operator $A$ admits a bounded $\HI(\Sigma^\pm_\omega)$ functional calculus.
\item The operator $A^2$ admits a bounded $\HI(\Sigma_{2\omega})$ functional calculus.
\end{enumerate}
\end{prop}

If $A$ is bisectorial, by \cite[Proposition 10.6.2 (2)]{HvNVW2} we have
\begin{equation}
\label{Bisec-Ran-Ker}
\ovl{\Ran A^2}
=\ovl{\Ran A}
\quad \text{and} \quad
\ker A^2
=\ker A.
\end{equation}

The following is \cite[Proposition G.2.4]{HvNVW2}.

\begin{lemma}
\label{Lemma-core-semigroup}
Let $(T_t)_{t \geq 0}$ be a strongly continuous semigroup of bounded operators on a Banach space $X$ with (negative) generator $A$. If $Y$ is subspace of $\dom A$ which is dense in $X$ and invariant under each operator $T_t$, then $Y$ is a core of $A$.
\end{lemma}

\paragraph{Compact operators and complex interpolation} Recall the following result \cite[Theorem 9]{CwK1} (see also \cite[Theorem 5.5]{KaMS}) which allows to obtain compactness via complex interpolation.

\begin{thm}
\label{Th-interpolation-Kalton}
Suppose that $(X_0,X_1)$ and $(Y_0,Y_1)$ are Banach couples and that $X_0$ is a $\UMD$-space. Let $T \co X_0+X_1 \to Y_0+Y_1$ such that is $T_0$ compact and $T_1$ is bounded. Then for any $0 < \theta < 1$ the map $T \co (X_0,X_1)_\theta \to (Y_0,Y_1)_\theta$ is compact.
\end{thm}

\paragraph{Markovian semigroups of Fourier multipliers}

In this paragraph, we recall the basic theory of markovian semigroups of Fourier multipliers. The following definition and properties of a markovian semigroup are fundamental for us. Thus the assumptions and notations which follow these lines are standing for all the paper. A companion definition of equal importance is Definition \ref{def-Fourier-markovian} together with Proposition \ref{prop-Schoenberg} which follow. Thn otion of selfadjoint operator is defined in \cite[p.~43]{JMX}.

\begin{defi}
\label{def-markovian}
Let $M$ be a von Neumann algebra equipped with a normal semifinite faithful trace. We say that a weak* continuous semigroup $(T_t)_{t \geq 0}$ of operators on $M$ is a markovian semigroup if each $T_t$ is a weak* continuous selfadjoint completely positive unital contraction.
\end{defi}

For any $1 \leq p < \infty$, such semigroup induces a strongly continuous semigroup of operators on each $\L^p(M)$ satisfying
\begin{enumerate}
\item each $T_t$ is a contraction on $\L^p(M)$,
\item each $T_t$ is selfadjoint on $\L^2(M)$,
\item each $T_t$ is completely positive on $\L^p(M)$.
\item $T_t(1)=1$
\end{enumerate}

Let $G$ be a discrete group and consider its regular representation $\lambda$ over the Hilbert space $\ell^2_G$ given by left translations $\lambda_s \co \delta_r \mapsto \delta_{sr}$ where $r,s \in G$. Then the group von Neumann algebra $\VN(G)$ is the von Neumann subalgebra of $\B(\ell^2_G)$ generated by these left translations. We also recall that $\mathrm{C}^*_r(G)$ stands for the reduced group $\mathrm{C}^*$-algebra sitting inside $\VN(G)$ and generated by these $\lambda_s$. Let us write 
\begin{equation}
\label{equ-trigonometric-polynomials}
\P_G 
\ov{\mathrm{def}}{=} \VNGfin
\end{equation}
for the space of ``trigonometric polynomials''. The von Neumann algebra $\VN(G)$ is equipped with the tracial faithful normal state $\tau(\lambda_s) \ov{\mathrm{def}}{=} \delta_{s = e} = \langle \lambda_s \delta_e, \delta_e \rangle$. Now, we introduce the main class of multipliers which interest us.  


\begin{defi}
\label{def-Fourier-markovian}
Let $G$ be a discrete group. A Fourier multiplier on $\VN(G)$ is a weak* continuous linear map $T \co \VN(G) \to \VN(G)$ such that there exists a (unique) complex function $\phi \co G \to \C$ such that for any $s \in G$ we have $T(\lambda_s) = \phi_s \lambda_s$. In this case, we let $M_\phi=T$ and we say that $\phi$ is the symbol of $T$.
\end{defi}

We have the following folklore characterization of markovian semigroups of Fourier multipliers. Proposition \ref{prop-Schoenberg} is central for the remainder of the article and the mappings $\pi$ and $b_\psi = b$ will be used throughout tacitly.

Recall that a function $\psi \co G \to \C$ is conditionally negative definite if $\psi(s)=\ovl{\psi(s^{-1})}$ for any $s \in G$, $\psi(e) \geq 0$ and if the condition $\sum_{i=1}^{n} c_i = 0$ implies $\sum_{i,j=1}^n \ovl{c_i} c_{j} \psi(s_j^{-1}s_i) \leq 0$. If $H$ is a real Hilbert space, $\mathrm{O}(H)$ stands for the orthogonal group.

\begin{prop}
\label{prop-Schoenberg}
Let $G$ be a discrete group and $(T_t)_{t \geq 0}$ be a family of operators on $\VN(G)$. Then the following are equivalent.
\begin{enumerate}
\item $(T_t)_{t \geq 0}$ is a markovian semigroup of Fourier multipliers.

\item There exists a (unique) real-valued conditionally negative definite function $\psi \co G \to \R$ satisfying $\psi(e) = 0$ such that $T_t(\lambda_s) = \exp(-t \psi(s)) \lambda_s$ for any $t \geq 0$ and any $s \in G$.

\item There exists a real Hilbert space $H$ together with a mapping $b \co G \to \C$ and a homomorphism $\pi \co G \to \mathrm{O}(H)$ such that the $1$-cocycle law holds 
\begin{equation}
\label{Cocycle-law}
\pi_s(b_\psi(t))
=b_\psi(st)-b_\psi(s),
\quad \text{i.e.} \quad 
b_\psi(st)
=b_\psi(s)+\pi_s(b_\psi(t))
\end{equation}
for any $s,t \in G$ and such that 
\begin{equation}
\label{liens-psi-bpsi}
\psi(s)
=\norm{b_\psi(s)}_H^2, \quad s \in G.
\end{equation}
\end{enumerate}
\end{prop}
Under these conditions, we say that $\psi$ is a conditionally negative length. Note also the equality \cite[Exercise 2.14.1]{BHV}
\begin{equation}
\label{cocycle-b--1}
b(s^{-1})
=-\pi_{s^{-1}}b(s).
\end{equation}
We refer to Section \ref{Sec-Levy-and-cocycles} for a link between L\'evy measure and 1-cohomology strongly related to this result.

\paragraph{Truncations of matrices}

\begin{defi}
\label{def-Troncature}
We denote by $\M_{I,\fin}$ the space of matrices indexed by $I \times I$ with a finite number of non null entries. Given a set $I$, the set $\mathcal{P}_f(I)$ of all finite subsets of $I$ is directed with respect to set inclusion. For $J\in \mathcal{P}_f(I)$ and $A\in \Mat_{I}$ a matrix indexed by $I \times I$, we write $\mathcal{T}_{J}(A)$ for the matrix obtained from $A$ by setting each entry to zero if its row and column index are not both in $J$. We call $\big(\mathcal{T}_{J}(A)\big)_{J \in \mathcal{P}_f(I)}$ the net of finite submatrices of $A$. See \cite[p.~488]{Arh2}.
\end{defi}

If $J$ is a finite subset of $I$, recall the truncation $\Tron_J \co S^p_I \to S^p_I$ is a  complete contraction. It is well-known that $\Tron_J(x) \to x$ in $S^p_I$ as $J \to I$ where $1 \leq p \leq \infty$.

\paragraph{Markovian semigroups of Schur multipliers}

Let $I$ be any non-empty index set. Let $A=[a_{ij}]_{i,j\in I}$ be a matrix of $\mathbb{M}_I$. By definition, the Schur multiplier on $\B(\ell^2_I)$ associated with this matrix is the unbounded linear operator $M_A$ whose domain $\dom M_A$ is the space of all $B=[b_{ij}]_{i,j\in I}$ of $\B(\ell^2_I)$ such that $[a_{ij}b_{ij}]_{i,j\in I}$ belongs to $\B(\ell^2_I)$, and whose action on $B=[b_{ij}]_{i,j \in I}$ is given by $M_A(B) \ov{\mathrm{def}}{=} [a_{ij}b_{ij}]_{i,j \in I}$. Often, we write $A$ for $M_A$ and $A*x$ for the Schur product $[a_{ij}b_{ij}]_{i,j \in I}$. Note the following property.


The following description \cite[Proposition 5.4]{Arh1} of a markovian semigroups consisting of Schur multipliers is central for our paper. See also \cite{Arh6} for a generalization.

\begin{prop}
\label{def-Schur-markovian}
Let $I$ be some non-empty index set and $(T_t)_{t \geq 0}$ be a family of operators on $\B(\ell^2_I)$. Then the following are equivalent.
\begin{enumerate}
\item $(T_t)_{t \geq 0}$ is a markovian semigroup of Schur multipliers.

\item There exists a real Hilbert space $H$ and a family 
\begin{equation}
\label{equ-Schur-markovian-alpha}
\alpha = (\alpha_i)_{i \in I}
\end{equation}
of elements of $H$ such that each $T_t \co \B(\ell^2_I) \to \B(\ell^2_I)$ is the Schur multiplier associated with the matrix
\begin{equation}
\label{Semigroup-Schur}
\Big[\e^{-t\norm{\alpha_i-\alpha_j}_{H}^2}\Big]_{i,j \in I}.
\end{equation}
\end{enumerate}
In this case, the weak* (negative) infinitesimal generator $A$ acts by $A(e_{ij}) = \norm{\alpha_i-\alpha_j}_{H}^2 e_{ij}$.
\end{prop}

\subsection{$q$-Gaussian functors, isonormal processes and probability}
\label{Sec-q-gaussians}

The definitions \eqref{def-partial-psi} and \eqref{def-delta-alpha} of the noncommutative gradients associated with our markovian semigroups of Schur and Fourier multipliers  need $q$-deformed Gaussian variables from \cite{BKS}. We recall here several facts about the associated von Neumann algebras. We denote by $S_n$ the symmetric group. If $\sigma$ is a permutation of $S_n$ we denote by $|\sigma|$ the number $\card\big\{ (i,j)\ |\ 1\leq i,j \leq n, \sigma(i)>\sigma(j)\big\}$ of inversions of $\sigma$. Let $H$ be a real Hilbert space with complexification $H_{\C}$. If $-1 \leq q < 1 $ the $q$-Fock space over $H$ is
$$
\mathcal{F}_{q}(H)
\ov{\mathrm{def}}{=} \C \Omega \oplus \bigoplus_{n \geq 1}
H_{\C}^{\ot_{n}}
$$
where $\Omega$ is a unit vector, called the vacuum and where the scalar product on $H_{\C}^{\ot_{n}}$ is given by
$$
\langle h_1 \ot \dots \ot h_n ,k_1 \ot \cdots \ot k_n \rangle_{q}
=\sum_{\sigma \in S_n} q^{|\sigma|}\langle h_1,k_{\sigma(1)}\rangle_{H_{\C}}\cdots\langle h_n,k_{\sigma(n)}\rangle_{H_{\C}}.
$$
If $q=-1$, we must first divide out by the null space, and we obtain the usual antisymmetric Fock space. The creation operator $\ell(e)$ for $e \in H$ is given by
$$
\begin{array}{cccc}
   \ell(e) \co &   \mathcal{F}_{q}(H)  &  \longrightarrow   &  \mathcal{F}_{q}(H)  \\
         &  h_1 \ot \dots \ot h_n   &  \longmapsto       & e \ot h_1 \ot  \dots \ot h_n.   \\
\end{array}
$$
They satisfy the $q$-relation
$$
\ell(f)^*\ell(e)-q\ell(e)\ell(f)^* 
= \langle f,e\rangle_{H} \Id_{\mathcal{F}_{q}(H)}.
$$
We denote by $s_q(e) \co \mathcal{F}_{q}(H) \to \mathcal{F}_{q}(H)$ the selfadjoint operator $\ell(e)+\ell(e)^*$. The $q$-von Neumann algebra $\Gamma_q(H)$ is the von Neumann algebra over $\mathcal{F}_q(H)$ generated by the operators $s_q(e)$ where $e \in H$. It is a finite von Neumann algebra with the trace $\tau$ defined by $\tau(x)=\langle\Omega,x.\Omega\rangle_{\mathcal{F}_{q}(H)}$ where $x \in \Gamma_q(H)$.

Let $H$ and $K$ be real Hilbert spaces and $T \co H \to K$ be a contraction with complexification $T_{\C} \co H_{\C} \to K_{\C}$. We define the following linear map
$$
\begin{array}{cccc}
 \mathcal{F}_q(T) \co   &  \mathcal{F}_{q}(H)   &  \longrightarrow   & \mathcal{F}_{q}(K)  \\
    &  h_1 \ot \dots \ot h_n   &  \longmapsto       &    T_{\C}(h_1) \ot \dots \ot T_{\C}(h_n). \\
\end{array}
$$

Then there exists a unique map $\Gamma_q(T) \co \Gamma_q(H) \to \Gamma_q(K)$ such that for every $x \in \Gamma_q(H)$ we have
\begin{equation}
\label{SQq}
\big(\Gamma_q(T)(x)\big)\Omega
=\mathcal{F}_q(T)(x\Omega).
\end{equation}
This map is normal, unital, completely positive and trace preserving. If $T \co H \to K$ is an isometry, $\Gamma_q(T)$ is an injective $*$-homomorphism. If $1 \leq p < \infty$, it extends to a contraction $\Gamma_q^p(T) \co \L^p(\Gamma_q(H)) \to \L^p(\Gamma_q(K))$.

Moreover, we need the following Wick formula, (see \cite[p.~2]{Boz} and \cite[Corollary 2.1]{EfP}). In order to state this, we denote, if $k \geq 1$ is an integer, by $\mathcal{P}_2(2k)$ the set
of 2-partitions of the set $\{1,2,\ldots,2k\}$. If $\mathcal{V}\in \mathcal{P}_2(2k)$ we let $c(\mathcal{V})$ the number of crossings of $\mathcal{V}$, which is given, by the number of pairs of blocks
of $\mathcal{V}$ which cross (see \cite[p.~8630]{EfP} for a precise definition). Then, if $f_1,\ldots,f_{2k}\in H$ we have
\begin{equation}
\label{formule de Wick}
 \tau\big(s_q(f_1)s_q(f_2)\cdots s_q(f_{2k})\big)  
=\sum_{\mathcal{V}\in \mathcal{P}_2(2k)} q^{c(\mathcal{V})}
\prod_{(i,j)\in \mathcal{V}}\langle f_i,f_j\rangle_{H}
\end{equation}
and for an odd number of factors of $q$-Gaussians, the trace vanishes,
\begin{equation}
\label{formule de Wick odd}
\tau \big(s_q(f_1)s_q(f_2)\cdots  s_q(f_{2k-1}) \big) = 0.
\end{equation}
In particular, for any $e,f \in H$, we have
\begin{equation}
\label{petit-Wick}
\tau\big(s_q(e)s_q(f)\big) 
=\langle e,f \rangle_{H}.
\end{equation}

Recall that if $e \in H$ has norm 1, then the operator $s_{-1}(e)$ satisfies 
\begin{equation}
\label{fermion-carre}
s_{-1}(e)^2
=\Id_{\mathcal{F}_{-1}(H)}.
\end{equation}


If $q=1$, the $q$-Gaussian functor identifies to an $H$-isonormal process on a probability space $(\Omega,\mu)$ \cite[Definition 1.1.1]{Nua1}, \cite[Definition 6.5]{Neer1}, that is a linear mapping $\W \co H \to \L^0(\Omega)$ with the following properties:
\begin{flalign}
& \label{isonormal-gaussian} \text{for any $h \in H$ the random variable $\W(h)$ is a centred real Gaussian,} \\
&\label{esperance-isonormal} \text{for any } h_1, h_2 \in H \text{ we have } \E\big(\W(h_1) \W(h_2)\big)= \langle h_1, h_2\rangle_H. \\
&\label{density-isonormal} \text{The linear span of the products } \W(h_1)\W(h_2)\cdots \W(h_m), 
\text{ with } m \geq 0 \text{ and } h_1,\ldots, h_m \\
&\nonumber\text{in }H, \text{ is dense in the real Hilbert space $\L^2_\R(\Omega)$.}\end{flalign}  
Here $\L^0(\Omega)$ denote the space of measurable functions on $\Omega$ and we make the convention that the empty product, corresponding to $m=0$ in \eqref{density-isonormal}, is the constant function $1$. Recall that the span of elements $\e^{\i\W(h)}$ is dense in $\L^p(\Omega)$ by \cite[Theorem 2.12]{Jan1} if $1 \leq p<\infty$ and weak* dense if $p=\infty$.  

If $1 \leq p \leq \infty$ and if $u \co H \to H$ is a contraction, we denote by $\Gamma_p(u) \co \L^p(\Omega) \to \L^p(\Omega)$ the (symmetric) second quantization of $u$ acting on the \textit{complex} Banach space $\L^p(\Omega)$. Recall that the map $\Gamma_\infty(u) \co \L^\infty(\Omega) \to \L^\infty(\Omega)$ preserves the integral\footnote{\thefootnote. That means that for any $f \in \L^\infty(\Omega)$ we have $\int_{\Omega} \Gamma_\infty(u)f \d\mu=\int_{\Omega} f \d\mu$.}. If $u \co H_1 \to H_2$ is an isometry between Hilbert spaces then $\Gamma_\infty(u) \co \L^\infty(\Omega_{H_1}) \to \L^\infty(\Omega_{H_2})$ is a trace preserving injective normal unital $*$-homomorphism which is surjective if $u$ is surjective. Moreover, we have if $1 \leq p< \infty$
\begin{equation}
\label{SQ2}
\Gamma_p(u)\W(h)
=\W(u(h)), \quad h \in H,
\end{equation}
and
\begin{equation}
\label{SQ1}
\Gamma_\infty(u) \e^{\i\W(h)}
=\e^{\i\W(u(h))}, \quad h \in H.
\end{equation}
Finally, note that if $u \co H \to H$ is a surjective isometry then $\Gamma_\infty(u) \co \L^\infty(\Omega) \to \L^\infty(\Omega)$ is a $*$-automorphism of the von Neumann algebra $\L^\infty(\Omega)$.



\paragraph{Rosenthal's inequalities} Recall the following notion of conditional independence, which was introduced in \cite[Definition 6.1]{DPPS1} and which is similar to the one of \cite[(I) p.~233]{JuX3}. Let $M$ be a semifinite von Neumann algebra equipped with a normal semifinite faithful trace $\tau$. Suppose that $(N_k)$ is a family of von Neumann subalgebras of $M$ and that $N$ is a common von Neumann subalgebra of the $N_k$ such that $\tau|N$ is semifinite. Let $\E_N \co M \to N$ be the canonical trace preserving faithful normal conditional expectation with respect to $N$. We say that the family $(N_k)$ is independent\footnote{\thefootnote. The authors of \cite{JuZ2} say fully independent.} with respect to $N$ if for every $k$ we have
$$
\E_N(xy) 
=\E_N(x)\E_N(y),\quad x \in N_k, y \in \W^*((N_j)_{j \not=k})
$$
where $\W^*((N_j)_{j \not=k})$ denotes the von Neumann subalgebra generated by the $N_j$ with $j \not=k$. A sequence $(x_k)$ of elements of $\L^p(M)$ is said to be faithfully independent with respect to a von Neumann subalgebra $N$ if there exists a family $(N_k)$ with $x_k \in \L^p(N_k)$ for any $k$ such that $(N_k)$ is independent with respect to $N$.

Let $(N_k)_{k \geq 1}$ be an independent sequence of von Neumann subalgebras of $M$ with respect to $N$. In \cite[Theorem 0.4 and (3.1)]{JuZ2}, it is shown in the case where $M$ is finite that if $2 \leq p <\infty$ and if $(x_k)$ is a sequence such that $x_k \in \L^p(N_k)$ and $\E_N(x_k) = 0$ for all $k$, then
\begin{align}
\label{Inequality-Rosenthal-1-prime}
\Bgnorm{\sum_{k=1}^{n} x_k}_{\L^p(M)}
&\lesssim \max\Bigg\{ \sqrt{p}\Bgnorm{\bigg(\sum_{k=1}^n \E_N(|x_k|^2)\bigg)^{\frac{1}{2}}}_{\L^p(M)},\sqrt{p}\Bgnorm{\bigg(\sum_{k=1}^n \E_N(|x_k^*|^2)\bigg)^{\frac{1}{2}}}_{\L^p(M)},\Bigg.\\
&\Bigg. p\bigg(\sum_{k=1}^n \norm{x_k}_{\L^p(M)}^p\bigg)^{\frac{1}{p}}\Bigg\}. \nonumber
\end{align}
We refer to \cite[Theorem 6.3]{DPPS1} for a variant of this result. Using reduction \cite{HJX}, it is easy to extend  this inequality to the non-finite case.

\subsection{Vector-valued unbounded bilinear forms on Banach spaces}
\label{Sec-sesquilinear-forms}

We need some notions on vector-valued unbounded bilinear forms on Banach spaces.
The reason is that these notions give an appropriate framework for the carr\'e du champ operators defined in \eqref{Def-Gamma-sgrp} and \eqref{Def-Gamma} and their domains associated with the markovian semigroups of Schur multipliers or Fourier multipliers. In the case of scalar valued unbounded forms on Hilbert spaces, we refer to \cite{Kat1}, \cite{MaR1} and \cite{Ouh1}.   



\begin{defi}
Let $X$ be a Banach space and $Z$ be an ordered-$*$-quasi-Banach space so carrying a positive cone and an involution. Assume that $D \subseteq X$ is a subspace and $a \co D \times D \to Z$ is an $\R$-bilinear map. We say that
\begin{enumerate}
  \item $a$ is densely defined if $D$ is dense in $X$.
	\item $a$ is symmetric if for any $x,y \in D$ we have $a(x,y)^* = a(y,x)$.
	\item $a$ is positive if for any $x \in D$ we have $a(x,x) \geq 0$. 
	\item $a$ satisfies the Cauchy-Schwarz inequality if for any $x,y \in D$
	\begin{equation}
\label{Sector-condition}
\norm{a(x,y)}_{Z}
\leq \norm{a(x,x)}_{Z}^{\frac{1}{2}}\norm{a(y,y)}_{Z}^{\frac{1}{2}}.
\end{equation}
\item A net $(x_i)$ of elements of $D$ is called $a$-convergent to $x \in X$ if $x_i \to x$ in $X$ and if $a(x_i-x_j,x_i-x_j) \to 0$ for $i,j \to \infty$. We denote this by $x_i \xra[]{a} x$. 
\end{enumerate}
\end{defi}

%
%

%

\begin{remark}
\label{Remark-xra}
\normalfont
If $a$ satisfies the Cauchy-Schwarz inequality and if $Z$ is $p$-normed\footnote{\thefootnote. That means that $\norm{y+z}_{Z}^p \leq \norm{y}_{Z}^p+\norm{z}_{Z}^p$ for any $y,z \in Z$.}, then for any $x,y \in D$ we have
\begin{align}
\MoveEqLeft
\label{Inegalite-triangulaire-1}
\norm{a(x+y,x+y)}_{Z}^{\frac{p}{2}} 
=\big(\norm{a(x+y,x+y)}_{Z}^p\big)^{\frac{1}{2}}
=\big(\norm{a(x,x)+a(x,y)+a(y,x)+a(y,y)}_{Z}^p\big)^{\frac{1}{2}} \\         
&\leq \big(\norm{a(x,x)}_{Z}^p+\norm{a(x,y)}_{Z}^p+\norm{a(y,x)}_{Z}^p + \norm{a(y,y)}_{Z}^p\big)^{\frac{1}{2}}\nonumber\\
&\ov{\eqref{Sector-condition}}{\leq} \big(\norm{a(x,x)}_{Z}^p+2\norm{a(y,y)}_{Z}^{\frac{p}{2}}\norm{a(y,y)}_{Z}^{\frac{p}{2}}+\norm{a(y,y)}_{Z}^p\big)^{\frac{1}{2}}\nonumber
=\norm{a(x,x)}_{Z}^{\frac{p}{2}}+\norm{a(y,y)}_{Z}^{\frac{p}{2}}. 
\end{align}
Hence, if $x_i \xra[]{a} x$ and $y_j \xra[]{a} y$ and if $\alpha,\beta \in \R$ then $\alpha x_i+\beta y_j \xra[]{a} \alpha x+\beta y$.
\end{remark}

In the sequel, we suppose that $Z$ is an ordered $*$-Banach space.
\begin{defi}
Let $a \co D \times D \to Z$ be a densely defined $\R$-bilinear map satisfying the Cauchy-Schwarz inequality. We say that $a$ is closed if $D$ is complete for the norm
\begin{equation}
\label{Def-norm-a}
\norm{x}_{a}
\ov{\mathrm{def}}{=}\sqrt{\norm{a(x,x)}_Z+\norm{x}_X^2}.
\end{equation} 
\end{defi}


\begin{prop}
\label{Prop-carac-closed}
Let $a \co D \times D \to Z$ be a densely defined $\R$-bilinear map satisfying the Cauchy-Schwarz inequality. The following are equivalent.
\begin{enumerate}
	\item $a$ is closed 
	\item For any net $(x_i)$, $x_i \xra[]{a} x$ implies $x \in D$ and $a(x_i-x,x_i-x) \to 0$
	\item For any sequence $(x_n)$, $x_n \xra[]{a} x$ implies $x \in D$ and $a(x_n-x,x_n-x) \to 0$.
\end{enumerate}

\end{prop}

\begin{proof}
$1 \Rightarrow 2$: Suppose that $x_i \xra[]{a} x$. This implies that $a(x_i-x_j,x_i-x_j) \to 0$ and $\norm{x_i- x_j}_X \to 0$ for $i,j \to \infty$ so by \eqref{Def-norm-a} that $\norm{x_i-x_j}_a \to 0$ as $i,j \to \infty$. By the completeness of $D$ there is a $x_0 \in D$ such that $\norm{x_i-x}_a \to 0$. Hence by \eqref{Def-norm-a} we have $a(x_i-x_0,x_j-x_0) \to 0$ and $\norm{x_i-x_0}_X \to 0$. Hence we must have $x=x_0 \in D$ and $a(x_i-x,x_i-x) \to 0$.

$2\Rightarrow 3$: It is obvious.

$3 \Rightarrow 1$: Let $(x_n)$ be a Cauchy sequence in $D$, that is $\norm{x_n-x_m}_a \to 0$ for $n,m \to \infty$. By \eqref{Def-norm-a}, the sequence $(x_n)$ is also Cauchy in $X$, so that there is a $x \in X$ such that $x_n \to x$. Since $a(x_n -x_m,x_n-x_m) \to 0$  by \eqref{Def-norm-a}, the sequence $(x_n)$ is $a$-convergent and the assumption implies that $x \in D$ and $a(x_n-x,x_n-x) \to 0$. By \eqref{Def-norm-a}, we have $\norm{x_n-x}_a \to 0$, which shows that $D$ is complete.
\end{proof}

\begin{prop}
\label{Prop-existence-limit}
Let $a \co D \times D \to Z$ be a densely defined $\R$-bilinear form satisfying the Cauchy-Schwarz inequality. 
\begin{enumerate}
	\item Let $(x_i)$ and $(y_i)$ be nets in $D$. If $x_i \xra[]{a} x$ and $y_i \xra[]{a} y$ such that the nets $(a(x_i,x_i))$ and $(a(y_i,y_i))$ are bounded, then $\lim_{i} a(x_i,y_i)$ exists. If $a$ is closed, then $\lim_{i} a(x_i,y_i)=a(x,y)$.
	\item Let $(x_n)$ be a sequence in $D$. If $x_n \xra[]{a} x$ then the sequence $(a(x_n,x_n))$ is bounded. 
\end{enumerate}
\end{prop}

\begin{proof}
1. In order to see that $\lim_{i} a(x_i,y_i)$ exists, we write: 
\begin{align*}
\MoveEqLeft
\norm{a(x_i,y_i)-a(x_j,y_j)}_Z
=\norm{a(x_i-x_j,y_i) + a(x_j,y_i-y_j)}_Z \\
&\leq \norm{a(x_i-x_j,y_i)}_Z+\norm{a(x_j,y_i-y_j)}_Z \\
&\ov{\eqref{Sector-condition}}{\leq} \norm{a(x_i-x_j,x_i-x_j)}_Z^{\frac{1}{2}}\norm{a(y_i,y_i)}_Z^{\frac{1}{2}} 
+\norm{a(x_j,x_j)}_Z^{\frac{1}{2}}\norm{a(y_i-y_j,y_i-y_j)}_Z^{\frac{1}{2}}.            
\end{align*}
It follows from the previous inequality that $(a(x_i,y_i))$ is a Cauchy net, hence a convergent net.

If $a$ is in addition closed note by Proposition \ref{Prop-carac-closed} that $x,y \in D$, $a(x_i-x,x_i-x) \to 0$ and $a(y_i-y,y_i-y) \to 0$. We conclude with
\begin{align*}
\MoveEqLeft
\norm{a(x_i,y_i)-a(x,y)}_Z
=\norm{a(x_i-x,y_i) + a(x,y_i-y)}_Z \\
&\leq \norm{a(x_i-x,y_i)}_Z+\norm{a(x,y_i-y)}_Z \\
&\ov{\eqref{Sector-condition}}{\leq} \norm{a(x_i-x,x_i-x)}_Z^{\frac{1}{2}}\norm{a(y_i,y_i)}_Z^{\frac{1}{2}} 
+\norm{a(x,x)}_Z^{\frac{1}{2}}\norm{a(y_i-y,y_i-y)}_Z^{\frac{1}{2}}.            
\end{align*}
2. We have $a(x_n-x_m,x_n-x_m) \to 0$ and $\norm{x_n- x_m}_X \to 0$ for $n,m \to \infty$ so by \eqref{Def-norm-a} that $\norm{x_n-x_m}_a \to 0$ as $n,m \to \infty$. So $(x_n)$ is a Cauchy sequence for the norm \eqref{Def-norm-a}, hence a bounded sequence. Hence the sequence $(a(x_n,x_n))$ is bounded by \eqref{Def-norm-a}. 
\end{proof}

\begin{defi}
A densely defined $\R$-bilinear map $a$ satisfying the Cauchy-Schwarz inequality is called closable if there exists a closed map $b \co D' \times D'\to Z$ such that $D \subseteq D' \subseteq X$ with $a(x,y) = b(x,y)$ for any $x,y \in D$.
\end{defi}


If $a$ is closable, we let
\begin{equation}
\label{Domaine-closure}
\dom \ovl{a}
\ov{\mathrm{def}}{=} \big\{x \in X \text{ such that there exists a sequence } (x_n) \text{ of } D \text{ satisfying } x_n \xra[]{a} x\big\}.
\end{equation}

\begin{prop}
\label{Prop-closable}
Let $a \co D \times D \to Z$ be a densely defined (symmetric positive) $\R$-bilinear form satisfying the Cauchy-Schwarz inequality. Then $a$ is closable if and only if $x_n \xra[]{a} 0$ implies $a(x_n,x_n) \to 0$. In this case, if $x_n \xra[]{a} x$ and if $y_n \xra[]{a} y$ where $x,y \in \dom \ovl{a}$ then the limit
\begin{equation}
\label{Def-closure}
\ovl{a}(x,y) 
\ov{\mathrm{def}}{=} \lim_{n \to \infty} a(x_n,y_n)
\end{equation}
exists and $\ovl{a}$ is a well-defined, densely defined (symmetric positive) closed form and satisfies the Cauchy-Schwarz inequality. In addition, every closed extension of $a$ is also an extension of $\ovl{a}$.
\end{prop}

\begin{proof}
$\Rightarrow$: Assume that $a$ is closable and let $b \co D' \times D' \to Z$ be a closed extension. Suppose $x_n \xra[]{a} 0$. Then $x_n \xra[]{b} 0$. Hence $a(x_n,x_n)=b(x_n,x_n)=b(x_n-0,x_n-0) \to 0$ since $b$ is closed.

$\Leftarrow$: By Proposition \ref{Prop-existence-limit}, the limit \eqref{Def-closure} exists. We will prove that $\lim_{n \to \infty} a(x_n,y_n)$ is independent of the chosen sequences $(x_n)$ and $(y_n)$. Indeed, if in addition $x_n' \xra[]{a} x$ and $y_n' \xra[]{a} y$, then
\begin{align*}
\MoveEqLeft
\norm{a(x_n,y_n)- a(x_n',y_n')}_Z
=\norm{a(x_n-x_n',y_n) + a(x_n',y_n-y_n')}_Z\\
&\leq \norm{a(x_n-x_n',y_n)}_Z+\norm{a(x_n',y_n-y_n')}_Z \\
&\ov{\eqref{Sector-condition}}{\leq} \norm{a(x_n-x_n',x_n-x_n')}_Z^{\frac{1}{2}}\norm{a(y_n,y_n)}_Z^{\frac{1}{2}}  +\norm{a(x_n',x_n')}_Z^{\frac{1}{2}}\norm{a(y_n-y_n',y_n-y_n')}_Z^{\frac{1}{2}}.            
\end{align*}
Note that $x_n-x_n' \xra[]{a} 0$ and $y_n-y_n' \xra[a]{} 0$. So by the assumption $a(x_n-x_n',x_n-x_n') \to 0$ and $a(y_n-y_n',y_n-y_n') \to 0$. With the point 2 of Proposition \ref{Prop-existence-limit}, we obtain $\norm{a(x_n,y_n)- a(x_n',y_n')}_Z \to 0$ as $n \to \infty$. 

For the closedeness of $\ovl{a}$, consider a Cauchy sequence $(x_n)$ of elements of $D$ for the norm $\norm{x}_{\ovl{a}}$. It converges to some $x \in X$. From the definition \eqref{Domaine-closure} of $\dom \ovl{a}$, we infer that $x$ belongs to $\dom \ovl{a}$. We infer that
$$
\ovl{a}(x_n-x,x_n-x)
\ov{\eqref{Def-closure}}{=}\lim_{n \to \infty} \lim_{m \to \infty} a(x_n-x_m,x_n-x_m)
=0
$$
which means that the sequence $(x_n)$ is convergent for the norm $\norm{x}_{\ovl{a}}$. Finally, note that $D$ is dense in $\dom \ovl{a}$. So by a classical argument $\dom \ovl{a}$ is also complete. 

Now we will show that $\ovl{a}$ is the smallest closed extension of $a$. Let $b$ be a closed extension of $a$. Let $x \in \dom \ovl{a}$. By \eqref{Domaine-closure}, there exists a sequence $(x_n)$ of $D$ such that $x_n \xra[]{a} x$. Then $x_n \xra[]{b} x$. So $x$ is an element of $\dom b$. We conclude that $\dom \ovl{a} \subseteq \dom b$ and $ a(x,y) =b(x,y)$ if $x,y \in D$ by \eqref{Def-closure} and Proposition \ref{Prop-existence-limit}. 
\end{proof}

\begin{defi}
\label{defi-closure}
If the form a is closable, then $\ovl{a}$ defined by \eqref{Def-closure} with domain $D(\ovl{a})$ is called the closure of the  bilinear map $a$.
\end{defi}
\subsection{Transference of Fourier multipliers on crossed product von Neumann algebras}
\label{subsubsec-Fourier-mult-crossed-product}

In the statement of our results on Riesz transforms associated with markovian semigroups of Fourier multipliers of Section \ref{sec:Semigroups-selfadjoint-contractive-Schur-multipliers} and later, we will need to use crossed products which is a basic method of constructing a new von Neumann algebra $M \rtimes_\alpha G$ from a von Neumann algebra $M$ acted on by a discrete group $G$. In our work, $M$ will be a $q$-Gaussian von Neumann algebra $\Gamma_q(H)$ of Subsection \ref{Sec-q-gaussians} and the action will be constructed in \eqref{equ-markovian-automorphism-group} by the second quantization from the orthogonal representation of Proposition \ref{prop-Schoenberg}. 
We refer to \cite{Haa1}, \cite{Haa3}, \cite{Str}, \cite{Sun} and \cite{Tak2} for more information on crossed products. 

Let $M$ be a von Neumann algebra acting on a Hilbert space $H$. Let $G$ be a locally compact group equipped with some left Haar measure $\mu_G$. Let $\alpha \co G \to \Aut(M)$ be a trace preserving representation of $G$ on $M$ which is weak* continuous, i.e. for any $x \in M$ and any $y \in M_*$, the map $G \to M$, $s \mapsto \langle \alpha_{s}(x), y \rangle_{M,M_*}$ is continuous. For any $x \in M$, we define the operators $\pi(x)\co \L^2(G,H) \to \L^2(G,H)$ \cite[(2) p.~263]{Str} by
\begin{equation}
\label{}
\big(\pi(x) \xi \big)(s)
\ov{\mathrm{def}}{=}\alpha^{-1}_s(x) \xi(s),\quad
\quad  \xi \in \L^2(G, H), s \in G.
\end{equation}
These operators satisfy the following commutation relation \cite[(2) p.~292]{Str}:
\begin{equation}
\label{commutation-rules}
(\lambda_s \ot \Id_H) \pi(x) (\lambda_s \ot \Id_H)^*
= \pi(\alpha_{s}(x)),
\quad x \in M, s \in G.
\end{equation}
Recall that the crossed product of $M$ and $G$ with respect to $\alpha$ is the von Neumann algebra $M \rtimes_\alpha G=(\pi(M) \cup \{\lambda_s \ot \Id_{H}: s \in G\})''$ on the Hilbert space $\L^2(G,H)$ generated by the operators $\pi(x)$ and $\lambda_s \ot \Id_{H}$ where $x \in M$ and $s \in G$. By \cite[p.~263]{Str} or \cite[Proposition 2.5]{Dae1}, $\pi$ is a normal injective $*$-homomorphism from $M$ into $M \rtimes_\alpha G$ (hence $\sigma$-strong* continuous). 

In the sequel, we suppose that $G$ is discrete. For any $s \in G$ and any $x \in M$, we let $x \rtimes \lambda_s\ov{\mathrm{def}}{=}\pi(x)(\lambda_s \ot \Id_H)$. We recall the rules of product and adjoint:
\begin{equation}
\label{inverse-crossed-product}
(x \rtimes \lambda_s)^*
=\alpha_{s^{-1}}(x^*) \rtimes \lambda_{s^{-1}}	
\end{equation}
and
\begin{equation}
\label{Product-crossed-product}
(x \rtimes \lambda_s)(y \rtimes \lambda_t)
=x\alpha_s(y) \rtimes \lambda_{st}, \quad s,t \in G, \xi \in H.
\end{equation}
These relations will be used frequently, as well as the following definition.
Namely, for $1 \leq p \leq \infty$ and $M$ a semifinite von Neumann algebra, we let
\begin{equation}
\label{equ-P-rtimes-G}
\P_{p,\rtimes,G} 
=\vect \{ x \rtimes \lambda_s : \: x \in \L^p(M) , \: s \in G \} \subseteq \L^p(M \rtimes_\alpha G),
\end{equation}
and write in short $\P_{\rtimes,G} = \P_{\infty,\rtimes,G}$, which will be the most used of these spaces.
Finally, if $M$ is equipped with a normal finite faithful trace $\tau_M$, then $M \rtimes_\alpha G$ is equipped with the normal finite faithful trace defined by $\tau_{M \rtimes_\alpha G}(x \rtimes \lambda_s) = \tau_M(x) \delta_{s = e}$ if $\tau_M$ is $G$-invariant. 

\begin{lemma}
\label{lem-PG-dense}
Assume $M$ to be a finite von Neumann algebra.
Let $1 < p < \infty$.
Then $\P_{\rtimes,G}$ and $\P_{p,\rtimes,G}$ are dense subspaces of $\L^p(M \rtimes_\alpha G)$.
In particular, $\P_G$ is dense in $\L^p(\VN(G))$.
\end{lemma}

\begin{proof}
Since $M$ is finite, we have $M \subseteq \L^p(M)$, and thus, $\P_{\rtimes,G} \subseteq \P_{p,\rtimes,G}$.
It thus suffices to prove that $\P_{\rtimes,G}$ is dense.
Let $y \in \L^{p^*}(M \rtimes_\alpha G)$ such that $\langle x , y \rangle = 0$ for all $x \in \P_{\rtimes,G}$.
By the Hahn-Banach theorem, it suffices to show that $y = 0$.
But $y \in \L^1(M \rtimes_\alpha G) = (M \rtimes_\alpha G)_*$.
By the rule \eqref{Product-crossed-product}, $\P_{\rtimes,G}$ contains $\vect(\pi(M) \cup \{\lambda_s \ot \Id_{H}: s \in G\})$, and therefore by construction, $\P_{\rtimes,G}$ is weak* dense in $M \rtimes_\alpha G$.
So $\langle x , y \rangle = 0$ for all $x \in M \rtimes_\alpha G$.
Thus, $y = 0$.
\end{proof}

Given a discrete group $G$, we can consider the ``fundamental unitary'' $W_H \co \ell^2_G \ot_2 H \ot_2 \ell^2_G \to \ell^2_G \ot_2 H \ot_2 \ell^2_G$ defined by
\begin{equation}
\label{eq:reg-omega-rep-unitary}
W_H(\epsi_t \ot \xi \ot \epsi_r )
=\epsi_t \ot \xi \ot \epsi_{t r}
\end{equation}
(see also \cite[p.~8]{BoBr}).

\begin{lemma}
\label{Lemma-Twisted-coproduct}
Then, for any $s \in G$ and any $x \in M$, we have
$$
W_H \big(  (x \rtimes \lambda_{s})\ot \Id_{\ell^2_G}\big) W_H^* 
=(x \rtimes \lambda_{s}) \ot \lambda_{s}.
$$
\end{lemma}

\begin{proof}
On the one hand, for any $s,t,r \in G$, we have
\begin{align*}
\MoveEqLeft
  W_H\big(  (x \rtimes \lambda_{s})\ot \Id_{\ell^2_G} \big)(\epsi_t \ot \xi \ot \epsi_r) 
		=W_H\big( x(\lambda_s \ot \Id_{H})\ot \Id_{\ell^2_G}\big)(\epsi_t \ot \xi\ot \epsi_r ) \\
		&=W_H\big(x \ot \Id_{\ell^2_G}\big)(\epsi_{st} \ot \xi \ot \epsi_{r})
		=W_H(\epsi_{st} \ot \alpha_{st}^{-1}(x)\xi \ot \epsi_{r}) \\
	&=\epsi_{st} \ot \alpha_{st}^{-1}(x)\xi \ot \epsi_{str}.
\end{align*}
On the other hand, we have
\begin{align*}
\MoveEqLeft
\big((x \rtimes \lambda_{s}) \ot \lambda_{s}\big)  W_H(\epsi_t \ot \xi \ot  \epsi_r)
		=\big( (x(\lambda_s \ot \Id_{H})) \ot \lambda_{s} \big) (\epsi_t\ot \xi \ot \epsi_{tr})\\
		&=\big( x \ot \Id_{\ell^2_G} \big) (\epsi_{st} \ot \xi \ot \epsi_{str})
		=\epsi_{st}  \ot \alpha_{st}^{-1}(x)\xi \ot \epsi_{str}.
\end{align*}
We conclude by linearity and density.
\end{proof}

As a composition of an amplification and a spatial isomorphism, the map $M \rtimes_\alpha \VN(G) \to \VN(G) \otvn (M \rtimes_{\alpha} G)$, $z \mapsto W_H \big(z \ot \Id_{\ell^2_G}\big) W_H^*$ is a unital normal injective $*$-homomorphism. This map is a well-defined kind of ``crossed coproduct'' defined by 
\begin{equation}
\label{Twisted-coproduct}
\begin{array}{cccc}
  \Delta \co &  M \rtimes_\alpha \VN(G)   &  \longrightarrow   &  \VN(G) \otvn (M \rtimes_{\alpha} G)  \\
           &  x \rtimes \lambda_{s}    &  \longmapsto & \lambda_s \ot (x \rtimes \lambda_s)   \\
\end{array}
\end{equation}
This homormorphism is trace preserving since 
\begin{align*}
\MoveEqLeft
\tau_{\VN(G) \otvn (M \rtimes_\alpha G)}\big(\Delta\big(x \rtimes \lambda_s)\big) 
\ov{\eqref{Twisted-coproduct}}{=}\tau_{\VN(G) \otvn (M \rtimes_\alpha G)} \left(\lambda_s \ot (x \rtimes \lambda_s)\right) \\
&=\tau_M(x) \delta_{s = e} \tau_{\VN(G)} (\lambda_s) 
=\tau_M(x) \delta_{s = e} 
=\tau_{M \rtimes_\alpha \VN(G)}(x \rtimes \lambda_s).
\end{align*}
By \cite[Lemma 4.1]{Arh2}, it admits a completely isometric $\L^p$-extension $\Delta_p \co \L^p(M \rtimes_\alpha \VN(G)) \to \L^p(\VN(G) \otvn (M \rtimes_{\alpha} G))$ for any $1 \leq p <\infty$.

If $x \in M$ and $s \in G$ and if $\phi \co G \to \C$, we let
\begin{equation}
\label{def-crossed-multiplier}
\big(\Id_{M} \rtimes M_\phi\big)(x \rtimes \lambda_s)
\ov{\mathrm{def}}{=} \phi(s) x \rtimes \lambda_s.
\end{equation}
The main result of this subsection is the following transference result. The assumptions are satisfied in the case where $M$ has $\QWEP$ and where the action $\alpha \co G \to \Aut(M)$ is amenable, see \cite[Proposition 4.1 (vi)]{Oza}. See also \cite[Proposition 4.8]{Arh1}.

\begin{prop}
\label{prop-Fourier-mult-crossed-product}
Suppose $1 \leq p \leq \infty$. Let $\phi \co G \to \C$ a function which induces a completely bounded Fourier multiplier $M_\phi \co \L^p(\VN(G)) \to \L^p(\VN(G))$. If $1 \leq p<\infty$, assume in addition that $M \rtimes_\alpha G$ has $\QWEP$. Then \eqref{def-crossed-multiplier} induces a completely bounded map $\Id_{\L^p(M)} \rtimes M_\phi \co \L^p(M \rtimes_\alpha G) \to \L^p(M \rtimes_\alpha G)$ and
\begin{equation}
\label{twisted-multiplier-major}
\norm{\Id_{\L^p(M)} \rtimes M_\phi}_{\cb,\L^p(M \rtimes_\alpha G) \to \L^p(M \rtimes_\alpha G)} 
\leq \norm{M_\phi}_{\cb,\L^p(\VN(G)) \to \L^p(\VN(G))}.
\end{equation}
Moreover, we have
\begin{equation}
\label{entelace-crossed}
\Delta_p\big(\Id_{M} \rtimes M_\phi\big)
=\big(M_\phi \ot \Id_{\L^p(M \rtimes_\alpha G)}\big)\Delta_p.
\end{equation}

\end{prop}

\begin{proof}
On the one hand, for any $x \in \L^p(M)$ and any $s \in G$, we have
\begin{align*}
\MoveEqLeft
\Delta_p\big(\Id_{\L^p(M)} \rtimes M_\phi\big)(x \rtimes \lambda_s)            
\ov{\eqref{def-crossed-multiplier}}{=}\phi(s)\Delta_p(x \rtimes \lambda_s)
\ov{\eqref{Twisted-coproduct}}{=} \phi(s)\big(\lambda_s \ot (x \rtimes \lambda_s)\big).
\end{align*} 
On the other hand, we have
\begin{align*}
\MoveEqLeft
\big(M_\phi \ot \Id_{\L^p(M \rtimes_\alpha G)}\big)\Delta_p(x \rtimes \lambda_s)             
\ov{\eqref{Twisted-coproduct}}{=} \big(M_\phi \ot \Id_{\L^p(M \rtimes_\alpha G)}\big) \big(\lambda_s \ot ( x \rtimes \lambda_s )\big) 
=\phi(s)\big(\lambda_s \ot (x \rtimes \lambda_s)\big).
\end{align*}
We conclude that \eqref{entelace-crossed} is true on $\P_{p,\rtimes,G}$. Since $M \rtimes_\alpha G$ has $\QWEP$, the von Neumann algebra $\VN(G)$ is also $\QWEP$. So, by \cite[(iii) p.~984]{Jun}, the map $M_\phi \ot \Id_{\L^p(M \rtimes_\alpha G)} \co \L^p(\VN(G) \otvn (M \rtimes_\alpha G)) \to \L^p( \VN(G) \otvn(M \rtimes_\alpha G))$ is completely bounded (in the case $p=\infty$, note that this assumption is useless) with 
\begin{equation}
\label{Equa-inter-QWEP}
\norm{M_\phi \ot \Id_{\L^p(M \rtimes_\alpha G)}}_{\cb,\L^p(\VN(G) \otvn (M \rtimes_\alpha G)) \to \L^p(\VN(G) \otvn (M \rtimes_\alpha G))} 
\leq \norm{M_\phi}_{\cb,\L^p(\VN(G)) \to \L^p(\VN(G))}.
\end{equation}
Note that $\Id_{S^p} \ot \Delta_p \co \L^p(\B(\ell^2) \otvn (M \rtimes_\alpha G)) \to \L^p(\B(\ell^2) \otvn \VN(G) \otvn (M \rtimes_\alpha G))$ is a complete isometry. If $z \in S^p \ot \P_{p,\rtimes,G}$, we have

\begin{align*}
\MoveEqLeft
\norm{\big(\Id_{S^p} \ot (\Id_{\L^p(M)} \rtimes M_\phi)\big)(z)}_{}
=\norm{\big(\Id_{S^p} \ot \big(\Delta_p\big(\Id_{\L^p(M)} \rtimes M_\phi\big)\big)\big)(z)}_{} \\
&\ov{\eqref{entelace-crossed}}{=} \norm{\big(\Id_{S^p} \ot\big(M_\phi \ot \Id_{\L^p(M \rtimes_\alpha G)}\big)\Delta_p\big)(z)}_{} 
\leq \norm{M_\phi \ot \Id_{\L^p(M \rtimes_\alpha G)}}_{\cb} \norm{\Delta_p}_{\cb} \norm{z}_{} \\
&\ov{\eqref{Equa-inter-QWEP}}{\leq} \norm{M_\phi}_{\cb,\L^p(\VN(G)) \to \L^p(\VN(G))}\norm{z}_{}.
\end{align*}
So \eqref{entelace-crossed} is proved.
\end{proof}


Recall that a discrete group has AP, if there exists a net $(\varphi_j)$ of finitely supported functions on $G$ such that $M_{\varphi_j} \ot \Id_{\B(H)} \to \Id_{\mathrm{C}^*_\lambda(G) \ot_{} \B(H)}$ in the point-norm topology for any Hilbert space $H$ by \cite[Theorem 1.9]{HK}. A discrete group $G$ is weakly amenable if there exists a net $(\varphi_j)$ of finitely supported functions on $G$ such that $\varphi_j \to 1$ pointwise and $\sup_j \bnorm{M_{\varphi_j}}_{\cb,\VN(G) \to \VN(G)} < \infty$. Recall that a weakly amenable discrete group has AP by \cite[p.~677]{HK}. Finally note that amenable groups, free groups, $\mathrm{SL}(2,\Z)$ and hyperbolic groups are weakly amenable. The first three parts of the following are \cite[Lemma 2.1]{CaSm1}. The part 4 is a straightforward consequence of Proposition \ref{prop-Fourier-mult-crossed-product} and the proof of \cite[Theorem 1.12]{HK}.

\begin{prop}
\label{Prop-weak-amenable-crossed}
Let $G$ be a discrete group with $\mathrm{AP}$ and $\alpha \co G \to \Aut(M)$ be a weak* continuous representation of $G$ on a von Neumann algebra $M$. Then there exist a net $(\varphi_j)_{}$ of finitely supported functions on $G$ with the following properties:
\begin{enumerate}
\item For each $s \in G$, we have $\lim_j \varphi_j(s) =1$.

\item Each $M_{\varphi_j}$ induces a normal completely bounded map $\Id \rtimes M_{\varphi_j} \co M \rtimes_\alpha G \to M \rtimes_\alpha G$.

\item The net $(\Id \rtimes M_{\varphi_j})$ converges to $\Id_{M \rtimes_\alpha G}$ for the point weak* topology.

\item If $G$ is in addition weakly amenable then $\sup_j \norm{\Id \rtimes M_{\varphi_j}}_{\cb, M \rtimes_\alpha G \to M \rtimes_\alpha G} <  \infty$.
\end{enumerate}
\end{prop}

An operator space $E$ has CBAP \cite[p.~205]{ER} (see also \cite[p.~365]{BrO} for the particular case of $\mathrm{C}^*$-algebras) when there exists a net $(T_j)$ of finite-rank linear maps $T_j \co E \to E$ satisfying the properties:
\begin{enumerate}
	\item for any $x \in E$, we have $\lim_j \norm{T_j(x)-x}_{E}=0$,
	\item $\sup_j \norm{T_j}_{\cb,E \to E} <\infty$.
\end{enumerate}
 
Then $E$ has moreover CCAP if the supremum in 2. above can be chosen to be $1$. By \cite[Proposition 4.3]{JR1}, if $M$ is a finite von Neumann algebra then $\L^p(M)$ has $\CCAP$ if and only if $\L^p(M)$ has $\CBAP$. In the case where $E$ is a noncommutative $\L^p$-space associated to a group von Neumann algebra of a discrete group $G$, we can use an average construction of Haagerup \cite[proof of Lemma 2.5]{Haa3} to replace the net $(T_j)$ by a net of Fourier multipliers. Indeed, we have the following observation which is contained in the proof of \cite[Theorem 4.4]{JR1}.

\begin{lemma}
\label{Lemma-carac-CBAP-Lp-spaces}
Let $G$ be a discrete group such that $\VN(G)$ is $\QWEP$. Suppose $1 \leq p<\infty$. Then the operator space $\L^p(\VN(G))$ has $\CCAP$ if and only if there exists a net $(\varphi_j)$ of functions $\varphi_j \co G \to \C$ with finite support which converge pointwise to 1 such that the net $(M_{\varphi_j})$ converges to $\Id_{\L^p(\VN(G))}$ in the point-norm topology with
$$
\sup_j \norm{M_{\varphi_j}}_{\cb, \L^p(\VN(G)) \to \L^p(\VN(G))} 
\leq 1.
$$
\end{lemma}

Recall that by \cite[Proposition 3.5]{JR1} if $G$ is a weakly amenable discrete group then $\L^p(\VN(G))$ has $\CBAP$ for any $1<p<\infty$. A careful reading of the proofs of \cite[Proposition 3.2]{JR1} and \cite[Proposition 3.5]{JR1} shows that the existence of an approximating net $(M_{\varphi_j})$ as in Lemma \ref{Lemma-carac-CBAP-Lp-spaces} which converges to $\Id_{\L^p(\VN(G))}$ for \textit{all} $1<p<\infty$ \textit{and} in $\mathrm{C}^*_r(G)$ in the point-norm topologies with in addition $\sup \norm{M_{\varphi_j}}_{\cb,\mathrm{C}^*_r(G) \to \mathrm{C}^*_r(G)} <\infty$. 
If $M \rtimes_\alpha G$ has QWEP, then using \cite[Proposition C.18]{EFHN}, Proposition \ref{prop-Fourier-mult-crossed-product} and the density of $\mathcal{P}_{\rtimes,G}$ in $\L^p(M \rtimes_\alpha G)$, it is easy to check that the net $(\Id \rtimes M_{\varphi_j})$ of completely bounded maps converges to $\Id_{\L^p(M \rtimes_\alpha G)}$ for the point norm topology of $\L^p(M \rtimes_\alpha G)$.

Consider finally a discrete group $G$ with AP such that $\VN(G)$ has $\QWEP$. By \cite[Theorem 1.2]{JR1}, if $1<p<\infty$ then $\L^p(\VN(G))$ has CCAP. But in this case, it is unclear if the previous paragraph is true for this class of groups and it is an interesting open question. 

We will use the following result which is stated in \cite[Proposition 3.13 and its proof]{Dae1}, see also \cite[Theorem 4.1]{HJX} for complements. Then the rather easy proof of Lemma \ref{lem-van-Daele-Arhancet} follows the same lines.

\begin{lemma}
\label{lem-van-Daele-Arhancet}
Let $T \co M \to M$ be a completely bounded normal map. Let $G$ be a discrete group and $\alpha \co G \to \Aut(M)$ be a weak* continuous action. If $T \alpha(s) = \alpha(s) T$ for any $s \in G$, then there is a linear normal mapping $T \rtimes \Id_{\VN(G)} \co M \rtimes_\alpha G \to M \rtimes_\alpha G$ of completely bounded norm that of $T$ such that
$$
\big(T \rtimes \Id_{\VN(G)}\big) (x \rtimes \lambda_s) 
= T(x) \rtimes \lambda_s.
$$
If $T$ is in addition a (faithful) normal conditional expectation then $T \rtimes \Id_{\VN(G)}$ is a (faithful) normal conditional expectation.
\end{lemma}

\subsection{Hilbertian valued noncommutative $\L^p$-spaces}
\label{Sec-Hilbertian-valued}

We shall use various Hilbert-valued noncommutative $\L^p$-spaces. These are some kind of noncommutative Bochner spaces which will be invaluable when we will use noncommutative Khintchine inequalities and some variants. We refer to \cite[Chapter 2]{JMX} for more information on these spaces. However, note that our notations are slightly different since we want that our notations remain compatible with vector-valued noncommutative $\L^p$-spaces in the hyperfinite case in the spirit of the notations of \cite{JRZ1}. 

Let $M$ be a semifinite von Neumann algebra. If a complex vector space $X$ which is a right $M$-module is equipped with an $\L^{\frac{p}{2}}(M)$-valued inner product $\langle \cdot, \cdot \rangle_X$\footnote{\thefootnote. By an $\L^{\frac{p}{2}}$-valued inner product on $X$ we mean a sesquilinear mapping $\langle \cdot, \cdot \rangle \co X \times X \to \L^{\frac{p}{2}}(M)$, conjugate linear in the first variable, satisfying for any $x,y \in X$ and any $z \in M$
\begin{enumerate}
\item $\langle x, y z \rangle = \langle x, y \rangle z$
\item $\langle x, y \rangle = \langle y, x \rangle^*$
\item $\langle x, x \rangle \geq 0$
\item $\langle x, x \rangle = 0 \iff x = 0$.
\end{enumerate}} where $0 < p \leq \infty$, recall that for any $x,y \in X$, we have by \cite[(2.2)]{JuS1}
\begin{equation} 
\label{Cauchy-Schwarz-Lp-module}
\norm{\langle x, y \rangle_X}_{\L^{\frac{p}{2}}(M)} 
\leq \norm{x}_X \norm{y}_X
\end{equation}
where
\begin{equation}
\label{Norm-Lp-module}
\norm{x}_X 
\ov{\mathrm{def}}{=} \bnorm{\langle x, x \rangle_{X}}_{\L^{\frac{p}{2}}(M)}^{\frac{1}{2}}.
\end{equation}

Note that the inner product is antilinear in the first and linear in the second variable.
We keep the same convention for the duality bracket between (subspaces of) our noncommutative $\L^p$ spaces.
The following is \cite[Proposition 2.2]{JuS1}. We include a short proof for the sake of completeness.

\begin{lemma}
Suppose $0< p \leq \infty$. Let $X$ be a complex vector space which is a right $M$-module equipped with an $\L^{\frac{p}{2}}(M)$-valued inner product $\langle \cdot, \cdot \rangle_{X}$ on $X$. Then \eqref{Norm-Lp-module} defines a norm on $X$ if $2 \leq p \leq \infty$ and a $\frac{p}{2}$-norm if $0<p \leq 2$.
\end{lemma}

\begin{proof}
For any $x,y \in X$, we have
$$
\norm{\langle x, y \rangle_X}_{\L^{\frac{p}{2}}(M)}
\ov{\eqref{Cauchy-Schwarz-Lp-module}}{\leq} \norm{x}_X \norm{y}_X
\ov{\eqref{Norm-Lp-module}}{=}\bnorm{\langle x, x \rangle_X}_{\L^{\frac{p}{2}}(M)}^{\frac{1}{2}}\bnorm{\langle y, y \rangle_X}_{\L^{\frac{p}{2}}(M)}^{\frac{1}{2}}.
$$
Hence the inequality \eqref{Sector-condition} is valid with $Z=\L^{\frac{p}{2}}(M)$ and with $a(x,y)=\langle x, y \rangle_{X}$. By Remark \ref{Remark-xra}, the proof is finished.
\end{proof}

Suppose $1 \leq p <\infty$. Let $H$ be a Hilbert space. For any elements $\sum_{k=1}^{n} x_k \ot a_k, \sum_{j=1}^{m} y_j \ot b_j$ of $\L^p(M) \ot H$, we define the $\L^{\frac{p}{2}}(M)$-valued inner product
\begin{equation}
\label{def-scal-product-Lp-H}
\bigg\langle \sum_{k=1}^{n} x_k \ot a_k,\sum_{j=1}^{m} y_j \ot b_j \bigg\rangle_{\L^p(M,H_{c,p})}
\ov{\mathrm{def}}{=}\sum_{k,j=1}^{n,m} \langle a_k,b_j\rangle_H\, x_k^*y_j.
\end{equation}
For any element $\sum_{k=1}^{n} x_k \ot a_k$ of $\L^p(M) \ot H$, we set \cite[(2.9)]{JMX}
\begin{align}
	\label{JMX-equa-2.9}
	\norm{\sum_{k=1}^{n} x_k \ot a_k}_{\L^p(M,H_{c,p})}
	&\ov{\mathrm{def}}{=} \norm{\left(\bigg\langle \sum_{k=1}^{n} x_k \ot a_k, \sum_{j=1}^{n} x_j \ot a_j\bigg\rangle_{\L^p(M,H_{c,p})} \right)^{\frac{1}{2}}}_{\L^p(M)} \\
	&=\norm{\bigg(\sum_{k,j=1}^{n} \langle a_k,a_j\rangle_H\, x_k^*x_j\bigg)^{\frac{1}{2}}}_{\L^p(M)}. \nonumber
\end{align}
The space $\L^p(M,H_{c,p})$ is the completion of $\L^p(M) \ot H$ for this norm \cite[p.~10]{JMX}. It is folklore that $\L^p(M,H_{c,p})$ equipped with \eqref{def-scal-product-Lp-H} is a right $\L^p$-$M$-module \cite[Definition 3.3]{JuS1}, see also \cite[Section 1.3]{JuPe1}.


If $(e_1,\ldots,e_n)$ is an orthonormal family of $H$ and if $x_1,\ldots,x_n$ belong to $\L^p(M)$, it follows from \eqref{JMX-equa-2.9} (see \cite[(2.10)]{JMX}) that 
\begin{equation*}
\Bgnorm{\sum_{k=1}^n x_k\ot
e_k}_{\L^p(M,H_{c,p})}=\Bgnorm{\bigg(\sum_{k=1}^{n} |x_k|^2\bigg)^{\frac{1}{2}}}_{\L^p(M)}
=\norm{\sum_{k=1}^{n} e_{k1} \ot x_k}_{S^p(\L^p(M))}.
\end{equation*}
We define similarly $\L^p(M,H_{r,p})$ and we let $\L^p(M,H_{\rad,p}) \ov{\mathrm{def}}{=}\L^p(M,H_{r,p}) \cap \L^p(M,H_{c,p})$ if $p \geq 2$ and $\L^p(M,H_{\rad,p}) \ov{\mathrm{def}}{=} \L^p(M,H_{r,p}) + \L^p(M,H_{c,p})$ if $1 \leq p \leq 2$. Recall that \cite[(2.25)]{JMX}
\begin{equation}
\label{Duality-conditional}
\L^p(M,H_{\rad,p})^*
=\L^{p^*}(M,\ovl{H}_{\rad,p^*}).
\end{equation}

If $T \co \L^p(M) \to \L^p(M)$ is completely bounded then by \cite[p.~22]{JMX} the map $T \ot \Id_H$ induces a bounded operator on $\L^p(M,H_{c,p})$, on $\L^p(M,H_{r,p})$ and on $\L^p(M,H_{\rad,p})$ and we have
\begin{equation}
\label{extens-col-1}
\norm{T \ot \Id_H}_{\L^p(M,H_{c,p}) \to \L^p(M,H_{c,p})} 
\leq \norm{T}_{\cb,\L^p(M) \to \L^p(M)}.
\end{equation}
(and similarly for the others, and also similarly for $T \co \L^p(M) \to \L^p(N)$ completely bounded between distinct spaces). If $T \co H \to K$ between Hilbert spaces is bounded then the map $\Id_{\L^p(M)} \ot T$ induces a bounded operator from $\L^p(M,H_{\rad,p})$ into $\L^p(M,K_{\rad,p})$ (and on $\L^p(M,H_{c,p})$ and $\L^p(M,H_{r,p})$) \cite[p.~14]{JMX} and we have
\begin{equation}
\label{extens-rad-1}
\norm{\Id_{\L^p(M)} \ot T}_{\L^p(M,H_{\rad,p}) \to \L^p(M,K_{\rad,p})} 
= \norm{T}_{H \to K}
\end{equation}
(and similarly for the others). 

We need some semifinite variant of spaces introduced in \cite{Jun3} in the $\sigma$-finite case. Let $\E \co \mathcal{N} \to M$ be a trace preserving normal faithful conditional expectation between semifinite von Neumann algebras equipped with normal semifinite faithful traces. If $2 \leq p \leq \infty$, for any $f,g \in \L^p(\mathcal{N})$, using the boundedness\footnote{\thefootnote. In the case $0<p<2$, note that the conditional expectation $\E \co \L^{\frac{p}{2}}(\mathcal{N}) \to \L^{\frac{p}{2}}(M)$ is not bounded in general.} of the conditional expectation $\E \co \L^{\frac{p}{2}}(\mathcal{N}) \to \L^{\frac{p}{2}}(M)$ we let $\langle f,g \rangle_{\L^p_{c}(\E)} \ov{\mathrm{def}}{=}\E(f^*g)$. It is clear that this defines an $\L^{\frac{p}{2}}(M)$-valued inner product and the associated right $\L^p$-$M$-module\footnote{\thefootnote. In the case $p=\infty$, we obtain a right $\W^*$-module.} is denoted by $\L^p_{c}(\E)$. This means that $\L^p_{c}(\E)$ is the completion of $\L^p(\mathcal{N})$ with respect to the norm
\begin{equation}
\label{Norm-conditional-general}
\norm{f}_{\L^p_{c}(\E)}
\ov{\mathrm{def}}{=} \bnorm{\E(f^*f)}_{\L^{\frac{p}{2}}(M)}^\frac{1}{2}.
\end{equation}
We denote always by $\langle f,g \rangle_{\L^p_{c}(\E)}$ the extension of the bracket on this space. Similarly, we define $\norm{f}_{\L^p_{r}(\E)} \ov{\mathrm{def}}{=} \bnorm{\E(ff^*)}_{\L^{\frac{p}{2}}(M)}^\frac{1}{2}$. The proof of \cite[Remark 2.7]{Jun3}, gives the following result. For the sake of completeness, we give the short computation.

\begin{lemma}
\label{lem-contractive-crossed-product-Hilbert-valued}
Suppose $2 \leq p \leq \infty$. We have a contractive injective inclusion of $\L^p(\mathcal{N})$ into $\L^p_{c}(\E)$ and into $\L^p_{r}(\E)$.
\end{lemma}

\begin{proof}
For any $f \in \L^p(\mathcal{N})$, we have
\begin{align*}
\MoveEqLeft
\norm{f}_{\L^p_{c}(\E)}
\ov{\eqref{Norm-conditional-general}}{=} \norm{\big(\E(f^* f)\big) }_{\L^{\frac{p}{2}}(M)}^{\frac12}            
\leq \norm{f^* f }_{\L^{\frac{p}{2}}(\mathcal{N})}^{\frac12}=\norm{f}_{\L^{p}(\mathcal{N})}.
\end{align*} 
We conclude by density. The row result is similar.
\end{proof}

In the case $0<p<2$, the definition of these spaces is unclear in the \textit{semifinite} case\footnote{\thefootnote. In the finite case, we can use the conditional expectation at the level $p=\infty$ and define $\norm{\cdot}_{\L^p_{c}(\E)}$ on the vector space $\mathcal{N}$. See below.}. However, if $\mathcal{N}=M \otvn N$ is a product with $N$ finite, we can use the conditional expectation $\E=\Id_M \ot \tau \co M \otvn N \to M$ at the level $p=\infty$. Indeed, for any $f,g \in (M \cap \L^p(M)) \ot N$, we can consider the element $\langle f,g \rangle_{\L^p_{c}(\E)} \ov{\mathrm{def}}{=}\E(f^*g)$ of $M \cap \L^{\frac{p}{2}}(M)$. It is clear that this defines an $\L^{\frac{p}{2}}(M)$-valued inner product and we can consider the associated right $\L^p$-$M$-module $\L^p_{c}(\E)$.

\begin{lemma}
\label{Lema-egalite-fantas} 
Suppose $0 < p< \infty$. Assume that the von Neumann algebra $N$ is finite and that $M$ is semifinite. Consider the canonical trace preserving normal faithful conditional expectation $\E \co M \otvn N \to M$. If $f$ and $g$ are elements of $(M \cap \L^p(M)) \ot N$ then we have
\begin{equation}
\label{Egalite-fantastique-2}
\E(f^* g)
=\langle f , g \rangle_{\L^p(M,\L^2(N)_{c,p})}.
\end{equation}
and
\begin{equation}
\label{Egalite-fantastique}
\norm{f}_{\L^p_{c}(\E)}
=\norm{f}_{\L^p(M,\L^2(N)_{c,p})}.
\end{equation}
So, we have a canonical isometric isomorphism $\L^p_{c}(\E)=\L^p(M,\L^2(N)_{c,p})$. 
\end{lemma}

\begin{proof}
We can suppose that $f=\sum_{k=1}^{n} x_k \ot a_k$ and $g=\sum_{j=1}^{m} y_j \ot b_j$. Note that $N \subseteq \L^p(N)$ for any $0 < p \leq \infty$. If $\tau$ is the normalized trace of $N$, we have
\begin{align*}
\MoveEqLeft
\E(f^* g)
=(\Id_M \ot \tau) \bigg(\sum_{k,j=1}^{n,m} (x_k^* \ot a_k^*) (y_j \ot b_j)\bigg)  
=\sum_{k,j=1}^{n,m} \tau(a_k^*b_j) x_k^*y_j  \\
&=\sum_{k,j=1}^{n,m} \langle a_k,b_j\rangle_{\L^2(N)}x_k^* y_j  
\ov{\eqref{def-scal-product-Lp-H}}{=} \langle f , g \rangle_{\L^p(M,\L^2(N)_{c,p})}.
\end{align*}
Moreover, we have
\begin{align*}
\MoveEqLeft
\norm{f}_{\L^p_{c}(\E)}
\ov{\eqref{Norm-conditional-general}}{=} \norm{\big(\E(f^* f)\big) }_{\L^{\frac{p}{2}}(M)}^{\frac12}
&\ov{\eqref{Egalite-fantastique-2}}{=} \norm{ \left( \langle f , f \rangle_{\L^p(M,\L^2(N)_{c,p})} \right) }_{\L^{\frac{p}{2}}(M)}^{\frac12} 
\ov{\eqref{JMX-equa-2.9}}{=} \norm{f}_{\L^p(M,\L^2(N)_{c,p})}.
\end{align*}
Finally, if $0<p <2$ note that $(M \cap \L^p(M)) \ot N$ is dense in $\L^p_{c}(\E)$ by definition. If $2 \leq p \leq \infty$, note that $(M \cap \L^p(M)) \ot N$ is dense in $\L^p(M \otvn N)$, hence in $\L^p_{c}(\E)$ by Lemma \ref{lem-contractive-crossed-product-Hilbert-valued}. Moreover, $(M \cap \L^p(M)) \ot N$ is dense in $\L^p(M,\L^2(N)_{c,p})$ since $\norm{\cdot}_{\L^p(M,\L^2(N)_{c,p})}$ is a cross norm by \eqref{JMX-equa-2.9}.
\end{proof}

Let $\E \co \mathcal{N} \to M$ be a trace preserving normal faithful conditional expectation between finite von Neumann algebras equipped with normal finite faithful traces. For any $f,g \in \mathcal{N}$, we can consider the element $\langle f,g \rangle_{\L^p_{c}(\E)} \ov{\mathrm{def}}{=}\E(f^*g)$ of $M \subseteq \L^{\frac{p}{2}}(M)$. It is clear that this defines an $\L^{\frac{p}{2}}(M)$-valued inner product and we can consider the associated right $\L^p$-$M$-module $\L^p_{c}(\E)$. If $\mathcal{N}$ is a crossed product, we have the following result. See Subsection \ref{subsubsec-Fourier-mult-crossed-product} for background on crossed products.
We denote here in accordance with \eqref{equ-P-rtimes-G} below $\P_{\rtimes,G} = \vect\{x \rtimes \lambda_s : \: x \in N, \: s \in G \}$.
Note that in \eqref{Egalite-fantastique-group}, we interpret $f$ to be equal to $\sum_{k = 1}^n x_k \rtimes \lambda_{s_k}$ on the left hand side, and reinterpret this as $f = \sum_{k  = 1}^n x_k \ot \lambda_{s_k}$ on the right hand side.

\begin{lemma}
\label{Lema-egalite-fantas-group}
Suppose $0 < p <\infty$. Let $G$ be a discrete group. Assume that the von Neumann algebra $N$ is finite. Let $\alpha \co G \to N$ be an action of $G$ on $N$ by trace preserving $*$-automorphisms. Consider the canonical trace preserving normal faithful conditional expectation $\E \co N \rtimes_\alpha G \to \VN(G)$. If $f=\sum_{k=1}^{n} x_k \rtimes \lambda_{s_k}$ and $g=\sum_{j=1}^{m} y_j \rtimes \lambda_{s_j}$ are elements of $\P_{\rtimes,G}$, then we have
\begin{equation}
\label{Egalite-fantastique-2-group}
\E(f^* g)
=\bigg\langle \sum_{k=1}^{n} \lambda_{s_k} \ot x_k, \sum_{j=1}^{m} \lambda_{s_j} \ot y_j \bigg\rangle_{\L^p(\VN(G),\L^2(N)_{c,p})},
\end{equation}
and
\begin{equation}
\label{Egalite-fantastique-group}
\norm{f}_{\L^p_{c}(\E)}
=\norm{f}_{\L^p(\VN(G),\L^2(N)_{c,p})}.
\end{equation}
So, we have a canonical isometric isomorphism $\L^p_{c}(\E)=\L^p(\VN(G),\L^2(N)_{c,p})$. 
\end{lemma}

\begin{proof}
Note that $N \subseteq \L^p(N)$ for any $0 < p \leq \infty$. If $\tau$ is the normalized trace of $N$, we have
\begin{align*}
\MoveEqLeft
\E(f^* g)
=\sum_{k,j=1}^{n,m} \E\big( (x_k \rtimes \lambda_{s_k})^* (y_j \rtimes \lambda_{s_j})\big)  
\ov{\eqref{inverse-crossed-product}}{=} \sum_{k,j=1}^{n,m} \E\big( (\alpha_{s_k^{-1}}(x_k^*) \rtimes \lambda_{s_k}^{-1}) (y_j \rtimes \lambda_{s_j})\big) \\
&\ov{\eqref{Product-crossed-product}}{=} \sum_{k,j=1}^{n,m} \E\big( (\alpha_{s_k^{-1}}(x_k^*)\alpha_{s_k^{-1}}(y_j) \rtimes \lambda_{s_k}^{-1} \lambda_{s_j})\big)
=\sum_{k,j=1}^{n,m}  \tau\big(\alpha_{s_k^{-1}}(x_k^*y_j)\big) \rtimes \lambda_{s_k}^{-1} \lambda_{s_j} \\
&=\sum_{k,j=1}^{n,m} \tau(x_k^*y_j) \lambda_{s_k}^*\lambda_{s_j}  
\ov{\eqref{def-scal-product-Lp-H}}{=} \bigg\langle \sum_{k=1}^{n} \lambda_{s_k} \ot x_k, \sum_{j=1}^{m} \lambda_{s_j} \ot y_j \bigg\rangle_{\L^p(\VN(G),\L^2(N)_{c,p})}.
\end{align*}
Moreover, we have
\begin{align*}
\MoveEqLeft
\norm{f}_{\L^p_{c}(\E)}
\ov{\eqref{Norm-conditional-general}}{=} \norm{\big(\E(f^* f)\big) }_{\L^{\frac{p}{2}}(\VN(G))}^{\frac12} \\
&\ov{\eqref{Egalite-fantastique-2-group}}{=} \norm{ \left( \langle f , f \rangle_{\L^p(\VN(G),\L^2(N)_{c,p})} \right) }_{\L^{\frac{p}{2}}(\VN(G))}^{\frac12} 
\ov{\eqref{JMX-equa-2.9}}{=} \norm{f}_{\L^p(\VN(G),\L^2(N)_{c,p})}.
\end{align*}
Finally, if $0<p <2$ note that $\P_{\rtimes,G}$ is dense in $\L^p_{c}(\E)$ by definition. If $2 \leq p \leq \infty$, note that $\P_{\rtimes,G}$ is dense in $\L^p(N \rtimes_\alpha G)$, hence in $\L^p_{c}(\E)$. Moreover, $\VN(G) \ot N$ is dense in $\L^p(\VN(G),\L^2(N)_{c,p})$ since $\norm{\cdot}_{\L^p(\VN(G),\L^2(N)_{c,p})}$ is a cross norm by \eqref{JMX-equa-2.9}.
\end{proof}

\begin{remark} \normalfont
\label{rem-problem-row}
It is easy to see that \eqref{Egalite-fantastique} holds in the row case $\norm{f}_{\L^p_r(\E)} = \norm{f}_{\L^p(M,\L^2(N)_{r,p})}$, too.
Note however that it in general the row analogue of \eqref{Egalite-fantastique-group} does not hold true.
Indeed, an elementary calculation reveals that
\[ 
\norm{ \sum_k \lambda_{s_k} \rtimes x_k }_{\L^{p}_r(\E)} 
= \norm{ \left( \sum_{k,l} \tau(x_k \alpha_{s_k s_l^{-1}}(x_l^*)) \lambda_{s_ks_l^{-1}} \right)^{\frac12} }_{\L^{p}_r(\E)} , 
\]
whereas according to the row version of \eqref{JMX-equa-2.9}
\[ 
\norm{ \sum_k \lambda_{s_k} \ot x_k}_{\L^{p}(\VN(G),\L^2(N)_{r,p})} 
= \norm{ \left( \sum_{k,l} \tau(x_k x_l^*) \lambda_{s_k s_l^{-1}} \right)^{\frac12} }_{\L^{p}(\VN(G))} . 
\]
\end{remark}

For the rest of the subsection, we shall show that in the situations of Lemma \ref{Lema-egalite-fantas} and Lemma \ref{Lema-egalite-fantas-group}, the Banach spaces $\L^p_c(\E)$ and $\L^p_r(\E)$ define an interpolation couple and that $\L^p_c(\E)$ and $\L^p_r(\E)$ admit $\L^{p^*}_c(\E)$ and $\L^{p^*}_r(\E)$ as dual spaces. 

First we consider the easy situation of Lemma \ref{Lema-egalite-fantas}. According to Lemma \ref{Lema-egalite-fantas} (resp. Remark \ref{rem-problem-row}), the spaces $\L^p_c(\E)$ and $\L^p(M,\L^2(N)_{c,p})$ (resp. $\L^p_r(\E)$ and $\L^p(M,\L^2(N)_{r,p})$) are isometric. According to \cite[p.~13]{JMX}, the spaces $\L^p_c(\E)$ and $\L^p_r(\E)$ embedd into a common Banach space $Z$ (one can take the injective tensor product of $\L^p(M)$ and $\L^2(N)$). In the sequel, we will always consider these embeddings.


\begin{lemma}
\label{lem-compatibility-Schur}
Suppose $1 < p < \infty$. Let $N$ be a finite von Neumann algebra and $M$ be a semifinite von Neumann algebra. Consider the conditional expectation $\E \co M \otvn N \to M$ from Lemma \ref{Lema-egalite-fantas}. 
\begin{enumerate}
	\item The embeddings respect $(M \cap \L^p(M)) \ot \L^2(N)$, which is dense in $\L^p_c(\E)$, in $\L^p_r(\E)$ and in $\L^p_c(\E) \cap \L^p_r(\E)$.
	\item The dual space of $\L^p_{c}(\E)$ (resp. $\L^p_{r}(\E)$) is $\L^{p^*}_c(\E)$ (resp. $\L^{p^*}_r(\E)$) with duality bracket $\langle f, g \rangle = \tau_{M \otvn N}(f^*g)$ for $f \in (M \cap \L^p(M)) \ot \L^2(N)$ and $g \in (M \cap \L^{p^*}(M)) \ot \L^2(N)$.
\end{enumerate}
\end{lemma}

\begin{proof}
1. The above mentioned embedding into the injective tensor product respects the subspace $(M \cap \L^p(M)) \ot \L^2(N)$ which is dense in both $\L^p_c(\E)$ and $\L^p_r(\E)$. We prove the third density. Let $(p_i)$ be an increasing net of orthogonal projections over $\L^2(N)$ with finite dimensional range, converging strongly to the identity. Then by \eqref{extens-rad-1}, \eqref{Egalite-fantastique} and Remark \ref{rem-problem-row}, the net $(\Id_{\L^p(M)} \ot p_i)$ is uniformly bounded as operators from $\L^p_c(\E)$ into $\L^p_c(\E)$ and from $\L^p_r(\E)$ into $\L^p_r(\E)$. Moreover, by density of $(M \cap \L^p(M)) \ot \L^2(N)$ separately in both $\L^p_c(\E)$ and $\L^p_r(\E)$, the net $(\Id_M \ot p_i)$ converges strongly to the identity of both $\L^p_c(\E)$ and $\L^p_r(\E)$. Moreover, both mappings are compatible as restrictions of a common mapping $\Id_{\L^p(M)} \ot p_i \co Z \to Z$. But if $f \in \L^p_c(\E) \cap \L^p_r(\E)$, then $(\Id_M \ot p_i)(f) = \sum_k f_k \ot e_k$, where the sum is finite, $f_k \in \L^p(M)$, and $e_k \in \L^2(N)$ span the range of $p_i$. Taking $g_k \in M \cap \L^p(M)$ with $\norm{g_k - f_k}_{\L^p(M)} < \epsi$, we see that $(M \cap \L^p(M)) \ot \L^2(N)$ is dense in the intersection $\L^p_c(\E) \cap \L^p_r(\E)$.

2. For the duality result, see \cite[(2.14)]{JMX} (note that the duality bracket is bilinear there, so we have to reinterpret it).
\end{proof}

Now we consider the situation of Lemma \ref{Lema-egalite-fantas-group}. Concerning the compatibility, we have the following result.

\begin{prop}
\label{prop-compatibility-crossed}
Suppose $1 < p < \infty$. Let $N$ be a finite von Neumann algebra. Consider the conditional expectation $\E \co N \rtimes_\alpha G \to \VN(G)$ from Lemma \ref{Lema-egalite-fantas-group}. There exist continuous injective embeddings $w_c \co \L^p_c(\E) \hookrightarrow \ell^\infty_G(\L^2(N))$ and $w_r \co \L^p_r(\E) \hookrightarrow \ell^\infty_G(\L^2(N))$ such that $w_c(x) = w_r(x)$ for all $x \in \P_{\rtimes,G}$.
\end{prop}

\begin{proof}
For $x = \sum_s x_s \rtimes \lambda_s$ belonging to $\P_{\rtimes,G}$, we consider the map $w(x) \ov{\mathrm{def}}{=} (x_s)_{s \in G}$. Let us show that $w$ extends to continuous mappings $w_c \co\L^p_c(\E) \to \ell^\infty_G(\L^2(N))$ and $w_r \co \L^p_r(\E) \to \ell^\infty_G(\L^2(N))$. 

Fix some $t \in G$. Since the trace $\tau$ is unital, we have a completely\footnote{\thefootnote. Recall that a contractive linear form is completely contractive by \cite[Corollary 2.2.3]{ER} and that $\L^p(M) \subseteq \L^1(M)$.} contractive mapping $\phi_t \co \L^p(\VN(G)) \to \C, \: x \mapsto x_t \ov{\mathrm{def}}{=} \tau(x \lambda_{t^{-1}})$. Thus according to \eqref{extens-col-1}, $\phi_t \ot \Id_{\L^2(N)}$ extends to a contraction $\L^p(\VN(G),\L^2(N)_{c,p}) \to \L^p(\C,\L^2(N)_{c,p}) = \L^2(N)$ where the last equality is \cite[Remark 2.13]{JMX}. We deduce the inequality $\norm{x_t}_{\L^2(N)} \leq \norm{x}_{\L^p(\VN(G),\L^2(N)_{c,p})} \ov{\eqref{Egalite-fantastique-group}}{=} \norm{x}_{\L^p_c(\E)}$. Taking the supremum over all $t \in G$ on the left hand side, since $\P_{\rtimes,G}$ is dense in $\L^p_c(\E)$, this shows that $w$ extends to a contraction $w_c \co \L^p_c(\E) \to \ell^\infty_G(\L^2(N))$.

Next we show that $w$ extends to a continuous mapping $w_r \co \L^p_r(\E) \to \ell^\infty_G(\L^2(N))$. According to Remark \ref{rem-problem-row}, we have $\norm{x}_{\L^p_r(\E)} = \norm{y}_{\L^p(\VN(G),\L^2(N)_{r,p})}$ with $y = \sum_s  \lambda_s \ot \alpha_{s^{-1}}(x_s)$. For any $t \in G$, note that since $\phi_t \ot \Id_{\L^2(N)}$ equally extends to a contraction $\L^p(\VN(G),\L^2(N)_{r,p}) \to \L^p(\C,\L^2(N)_{r,p}) = \L^2(N)$ and $\alpha_{t^{-1}} \co \L^2(N) \to \L^2(N)$ is an isometry, we also have $\norm{x_t}_{\L^2(N)} = \norm{\alpha_{t^{-1}}(x_t)}_{\L^2(N)} \leq \norm{y}_{\L^p(\VN(G),\L^2(N)_{r,p})} = \norm{x}_{\L^p_r(\E)}$. Taking the supremum over all $t \in G$ on the left hand side, since $\P_{\rtimes,G}$ is dense in $\L^p_r(\E)$, this shows that $w$ extends to a contraction $w_r \co \L^p_r(\E) \to \ell^\infty_G(\L^2(N))$.

Next we show that $w_c$ is injective. Assume that $w_c(x) = 0$ for some $x \in \L^p_c(\E)$. If $\Psi \co \L^p_c(\E) \to \L^p(\VN(G),\L^2(N)_{c,p})$ is the isometry from \eqref{Egalite-fantastique-group}, the above proof of continuity shows that $w_c(x) = \big((\phi_s \ot \Id_{\L^2(N)}) \circ \Psi(x)\big)_{s \in G}$. We infer that $(\phi_s \ot \Id_{\L^2(N)}) \circ \Psi(x) = 0$ for all $s \in G$. Thus\footnote{\thefootnote. If $(\phi_s \ot \Id_{\L^2(N)})(y) = 0$ for some $y \in \L^p(\VN(G),\L^2(N)_{c,p})$ and all $s \in G$, then to conclude that $y = 0$, by duality, it suffices to check that $\langle y , y^* \rangle = 0$ for $y^*$ belonging to some total subset of $\L^{p^*}(\VN(G),\L^2(N)_{c,p})$. We can choose the total subset $\P_G \ot \L^2(N)$ and conclude by $\langle y, \lambda_s \ot h \rangle = \langle (\phi_s \ot \Id_{\L^2(N)})(y), h \rangle = 0$ (the first equality is true on dense subset, hence correct).}
, $\Psi(x) = 0$, hence $x = 0$.

Finally we show that $w_r$ is injective. Assume that $w_r(x) = 0$ for some $x \in \L^p_r(\E)$. If $\Phi \co \L^p_r(\E) \to \L^p(\VN(G),\L^2(N)_{r,p})$ is the isometry from Remark \ref{rem-problem-row}, we have $w_r(x) = \big(\alpha_{s} \circ ( \phi_s \ot \Id_{\L^2(N)}) \circ \Phi(x)\big)_{s \in G}$. We infer that $\alpha_s \circ( \phi_s \ot \Id_{\L^2(N)}) \circ \Phi(x) = 0$ for all $s \in G$. Consequently, we have $(\phi_s \ot \Id_{\L^2(N)}) \circ \Phi(x) = 0$ for all $s \in G$. Hence $\Phi(x) = 0$, so $x = 0$.
\end{proof}

We keep the notations from Proposition \ref{prop-compatibility-crossed} above.
As we have observed in the proof of Lemma \ref{Lema-egalite-fantas-group}, $\P_{\rtimes,G}$ is dense in $\L^p_c(\E)$, and thus also in $\L^p_r(\E)$, since $\P_{\rtimes,G}$ is stable by taking adjoints. We shall also need that $\P_{\rtimes,G}$ is dense in the intersection $\L^p_c(\E) \cap \L^p_r(\E)$. This is implicitly used in a similar situation in \cite[p.~541]{JMP2}\footnote{\thefootnote.\label{footnote-lem-density-intersection-CBAP} In line -3 of \cite[p.~541]{JMP2}, it is stated that $f$ is a finite sum $\sum_{g,h} (B(h) \ot f_{g,h}) \rtimes \lambda(g)$.} and \cite[p.~543]{JMP2}\footnote{\thefootnote. In the calculations line 3 and 4 of \cite[p.~543]{JMP2}, it is used that the duality brackets between $\L_{rc}^p(E_{\mathcal{M} \rtimes G})$ and $\L_{rc}^q(E_{\mathcal{M} \rtimes G})$, and between $\L^p(\mathcal{A})$ and $\L^q(\mathcal{A})$ are compatible, which a priori only makes sense on common subspaces. It is also used that the Gaussian projection is bounded $\widehat{Q} \co \L^q_{rc}(E_{\mathcal{M} \rtimes G}) \to G_q(\mathcal{M}) \rtimes G$ which is a priori only clear on the subspace of $f$ as in footnote \ref{footnote-lem-density-intersection-CBAP}.} in the case of a crossed product $\mathcal{M} \otvn \L^\infty(\Omega,\Sigma,\mu) \rtimes G$ where $\mathcal{M}$ is a finite von Neumann algebra. We provide a proof with some assumptions in Lemma \ref{lem-density-intersection-CBAP} which is needed in the proof of Theorem \ref{Khintchine-group-twisted} below. Note that the approach to define compatibility via \textit{non-constructive} Hilbert module theory seems \textit{ineffective }to obtain this kind of results.
We first record, also for later use, the following lemma.

\begin{lemma}
\label{lem-Fourier-mult-on-row-space}
Let $1 < p < \infty$ and $G$ be a discrete group acting on $N$ by normal trace preserving $*$-automorphisms, where $N$ is some finite von Neumann algebra. Let $M_\varphi$ be a completely bounded Fourier multiplier on $\L^p(\VN(G))$. Then
\begin{equation}
\label{equ-Fourier-mult-on-column-space} 
\bnorm{(\Id_{N} \rtimes M_{\varphi})(f)}_{\L^p_{c}(\E)} 
\lesssim_p \norm{f}_{\L^p_{c}(\E)}
\quad \text{and} \quad
\bnorm{(\Id_{N} \rtimes M_{\varphi})(f)}_{\L^p_{r}(\E)} 
\lesssim_p \norm{f}_{\L^p_{r}(\E)} .
\end{equation}
\end{lemma}

\begin{proof}
By the isometry in Lemma \ref{Lema-egalite-fantas-group} and the fact that $M_{\varphi}$ is completely bounded, the column part of \eqref{equ-Fourier-mult-on-column-space} follows immediately from \eqref{extens-col-1}. For the row part of \eqref{equ-Fourier-mult-on-column-space}, we need a little bit of work, due to Remark \ref{rem-problem-row}. Let $f = \sum_k x_k \rtimes \lambda_{s_k}$ be an element of $\P_{\rtimes,G}$ which is a dense subspace of $\L^p_r(\E)$.
Using the notation $\check{\varphi}(s) \ov{\mathrm{def}}{=} \ovl{\varphi(s^{-1})}$ and \cite[Lemma 6.4]{ArK1} in the last equality, we have
\begin{align*}
\MoveEqLeft
\bnorm{(\Id_{N} \rtimes M_{\varphi})(f) }_{\L^p_r(\E)} 
= \bnorm{\big((\Id_{N} \rtimes M_{\varphi})(f)\big)^* }_{\L^p_c(\E)} \nonumber 
= \norm{ \bigg( \sum_k \varphi(s_k) x_k \rtimes \lambda_{s_k} \bigg)^* }_{\L^p_c(\E)} \nonumber \\
& \overset{\eqref{inverse-crossed-product}}{=} \norm{ \sum_k \ovl{\varphi(s_k)} \alpha_{s_k^{-1}}(x_k^*) \rtimes \lambda_{s_k^{-1}} }_{\L^p_c(\E)} \nonumber 
= \norm{(\Id_{N} \rtimes M_{\check{\varphi}}) \bigg( \sum_k \alpha_{s_k^{-1}}(x_k^*) \rtimes \lambda_{s_k^{-1}} \bigg) }_{\L^p_c(\E)} \nonumber \\
&=\norm{(\Id_{N} \rtimes M_{\check{\varphi}})(f^*) }_{\L^p_c(\E)} \nonumber 
\overset{\eqref{equ-Fourier-mult-on-column-space}}{\leq} \norm{M_{\check{\varphi}}}_{\cb, \L^p(\VN(G)) \to \L^p(\VN(G))} \norm{f^*}_{\L^p_c(\E)} \nonumber \\
&=\norm{M_{\check{\varphi}}}_{\cb, \L^p(\VN(G)) \to \L^p(\VN(G))} \norm{f}_{\L^p_r(\E)} 
=\norm{M_\varphi}_{\cb, \L^p(\VN(G)) \to \L^p(\VN(G))} \norm{f}_{\L^p_r(\E)} \label{equ-proof-lem-Fourier-mult-on-row-space}.
\end{align*} 
\end{proof}

With this lemma at hand, we get the following.

\begin{lemma}
\label{lem-density-intersection-CBAP}
Let $1 < p < \infty$ and keep the assumptions of Proposition \ref{prop-compatibility-crossed}. Assume in addition that $\L^p(\VN(G))$ has $\CCAP$ and that $\VN(G)$ has $\QWEP$. Consider $\L^p_c(\E)$ and $\L^p_r(\E)$ injected into $\ell^\infty_G(\L^2(N))$ as in Proposition \ref{prop-compatibility-crossed}. Then $\P_{\rtimes,G}$ is dense in $\L^p_c(\E) \cap \L^p_r(\E)$.
\end{lemma}

\begin{proof}
Consider the approximating net $(M_{\varphi_j})$ of completely bounded Fourier multipliers provided by the $\CCAP$ assumption by Lemma \ref{Lemma-carac-CBAP-Lp-spaces} with $\norm{M_{\varphi_j}}_{\cb, \L^p(\VN(G)) \to \L^p(\VN(G))} \leq 1$. By \cite[Lemma 6.5, Lemma 6.6]{ArK1}, note that $(\varphi_j)$ is a bounded net of $\ell^\infty_G$.

A direct estimate shows that the $M_{\varphi_j} \ot \Id_{\L^2(N)} \co \ell^\infty_G(\L^2(N)) \to \ell^\infty_G(\L^2(N))$, $(x_s)_{s \in G} \mapsto (\varphi_j(s)x_s)_{s \in G}$'s define a bounded net of operators. Note that according to Lemma \ref{lem-Fourier-mult-on-row-space}, the extensions $\Id_{N} \rtimes M_{\varphi_j}$ are uniformly bounded as operators $\L^p_c(\E) \to \L^p_c(\E)$ and $\L^p_r(\E) \to \L^p_r(\E)$. 

Thus, by density of $\P_{\rtimes,G}$ in $\L^p_c(\E)$ and in $\L^p_r(\E)$, the net $(\Id_N \rtimes M_{\varphi_j})$ converges strongly to $\Id_{\L^p_c(\E)}$ resp. $\Id_{\L^p_r(\E)}$. By compatibility of this mapping on both spaces as a restriction of $M_{\varphi_j} \ot \Id_{\L^2(N)} \co \ell^\infty_G(\L^2(N)) \to \ell^\infty_G(\L^2(N))$, the net $(\Id_N \rtimes M_{\varphi_j})$ converges also to the identity on $\L^p_c(\E) \cap \L^p_r(\E)$. It thus suffices to show density of $\P_{\rtimes,G}$ in $\Ran(\Id_N \rtimes M_{\varphi_j} \co \L^p_c(\E) \cap \L^p_r(\E) \to \L^p_c(\E) \cap \L^p_r(\E))$ for some fixed $j$.

For $f \in \L^p_c(\E) \cap \L^p_r(\E)$, we have $(\Id_N \rtimes M_{\varphi_j})(f) = \sum_{s \in \supp \varphi_j} \varphi_j(s) f_s \rtimes \lambda_s$ for some $f_s \in \L^2(N)$.
Since the above sum is finite, we are left to show that $f \rtimes \lambda_s$ can be approximated by elements of $\P_{\rtimes,G}$ in $\L^p_c(\E) \cap \L^p_r(\E)$.
But this is immediate since both column and row norms are cross norms and $N$ is dense in $\L^2(N)$.
\end{proof}

By Lemma \ref{lem-compatibility-Schur} (resp. Proposition \ref{prop-compatibility-crossed}), $(\L^p_c(\E), \L^p_r(\E))$ is an interpolation couple. Therefore, it is meaningful to let $\L^p_{cr}(\E) \ov{\mathrm{def}}{=} \L^p_c(\E) + \L^p_r(\E)$ if $1 < p < 2$ and $\L^p_{cr}(\E) \ov{\mathrm{def}}{=} \L^p_c(\E) \cap \L^p_r(\E)$ if $2 \leq p < \infty$. In the situation of Lemma \ref{lem-compatibility-Schur}, we get immediately with \cite[2.7.1 Theorem]{BeL} that for $1 < p < \infty$, $\L^p_{cr}(\E)$ and $\L^{p^*}_{cr}(\E)$ are dual spaces with duality bracket $\langle f, g \rangle = \tau_{M \otvn N}(f^*g)$.

In the situation of Proposition \ref{prop-compatibility-crossed}, we have the following duality result, see \cite[Corollary 2.12]{Jun3} and \cite[2.7.1 Theorem]{BeL}. Note that there is a separability assumption in \cite[Corollary 2.12]{Jun3}. So our discrete groups considered in this paper are implicitly countable when we use this result. Nevertheless, this assumption seems removable with the transparent argument of \cite[Remark 6.11]{JuZ3} (see also the discussion before \cite[Lemma 3.4]{JuR}).

\begin{lemma}
\label{lem-expectation-duality}
For $1 < p < \infty$, $\L^p_c(\E)$ and $\L^{p^*}_c(\E)$ (resp. $\L^p_r(\E)$ and $\L^{p^*}_r(\E)$, resp. $\L^p_{cr}(\E)$ and $\L^{p^*}_{cr}(\E)$) are dual spaces with respect to each other. Here, the duality pairing is $\langle f , g \rangle = \tau_{\rtimes}(f^*g)$ for all $f,g \in N \rtimes_\alpha G$.
\end{lemma}

\subsection{Carr\'e du champ $\Gamma$ and first order differential calculus for Fourier multipliers}
\label{Sec-infos-on-Gamma-Fourier}

Let $G$ be a discrete group. Given a markovian semigroup $(T_t)_{t \geq 0}$ of Fourier multipliers on the von Neumann algebra $\VN(G)$ as given in Proposition \ref{prop-Schoenberg}, we shall introduce the noncommutative gradient $\partial_{\psi,q}$ and the carr\'e du champ $\Gamma$, and interrelate them together with the square root $A^{\frac12}$ of the negative generator $A$. For any $1 \leq p \leq \infty$, note  that the space $\P_G$ from \eqref{equ-trigonometric-polynomials} is a subspace of $\L^p(\VN(G))$.

For $x,y \in \P_G$, we define the element 
\begin{equation}
\label{Def-Gamma-sgrp}
\Gamma(x,y) 
\ov{\textrm{def}}{=}\frac{1}{2}\big[ A(x^*) y + x^* A(y) - A(x^* y)\big] 
\end{equation}
of $\VNGfin$. For any $s,t \in G$, note that
\begin{align}
\MoveEqLeft
\label{Def-Gamma-lambda-s}
\Gamma(\lambda_s,\lambda_t)
=\frac{1}{2}\big[ \psi(s^{-1}) + \psi(t) - \psi(s^{-1}t) \big]\lambda_{s^{-1}t}.
\end{align}
The form $\Gamma$ is usually called carr\'e du champ, associated with the markovian semigroup $(T_t)_{t \geq 0}$.
The first part of the following result shows how to construct the carr\'e du champ directly from the semigroup $(T_t)_{t \geq 0}$.

\begin{lemma}
\label{Lemma-gamma-infos-sgrp}
Let $G$ be a discrete group.
\begin{enumerate}
	\item Suppose $1 \leq p \leq \infty$. For any $x, y \in \P_G$, we have
\begin{equation}
	\label{Equa-Gamma-xy-sgrp}
\Gamma(x,y) 
=\lim_{t \to 0^+} \frac1{2t} \big( T_t(x^*y) - (T_t(x))^*T_t(y) \big)
\end{equation}
in $\L^p(\VN(G))$.
	\item For any $x \in \P_G$, we have $\Gamma(x,x) \geq 0$.
	\item For any $x, y \in \P_G$, the matrix 
$\begin{bmatrix} 
\Gamma(x,x) & \Gamma(x,y) \\ 
\Gamma(y,x) & \Gamma(y,y) 
\end{bmatrix}$ is positive.
\item For any $x,y \in \P_G$, we have
\begin{equation}
\label{CS-group}
\norm{\Gamma(x,y)}_{\L^p(\VN(G))} \leq \norm{\Gamma(x,x)}_{\L^p(\VN(G))}^{\frac12} \norm{\Gamma(y,y)}_{\L^p(\VN(G))}^{\frac12}.
\end{equation}
This inequality even holds for $0 < p \leq \infty$.
\end{enumerate}
\end{lemma}

\begin{proof}
1. For $x, y \in \VNGfin$, using the joint continuity of multiplication on bounded sets for the strong operator topology on $\L^p(\VN(G))$ \cite[Proposition C.19]{EFHN} (if $p = \infty$ we use that $x,y \in \VNGfin$, so all operators below belong to some finite dimensional space $\subseteq \L^p(\VN(G)) \cap \VN(G)$ on which convergence in $\VN(G)$ norm follows from convergence in $\L^p(\VN(G))$ norm), we have
\begin{align*}
\MoveEqLeft
\frac{1}{2t}\big( T_t(x^*y) - (T_t(x))^*T_t(y) \big) 
=\frac{1}{2t} \left(T_t(x^*y) - x^*y + x^*y - (T_t(x))^*T_t(y) \right)  \\
&=\frac{1}{2t} \big(T_t(x^*y) - x^*y + x^*y - (T_t(x))^*y + (T_t(x))^*y - (T_t(x))^*T_t(y) \big)  \\
&=\frac{1}{2t} \big(T_t(x^*y) - x^*y + (x^* - T_t(x^*))y + (T_t(x))^*(y-T_t(y)) \big)\\
&\xra[t \to 0^+]{\eqref{semigroupe-vers-generateur}} \frac12 \big(-A(x^*y) + A(x^*)y + x^*A(y) \big)
\ov{\eqref{Def-Gamma-sgrp}}{=}\Gamma(x,y).
\end{align*}

2. Each operator $T_t \co \VN(G) \to \VN(G)$ is completely positive and unital. Hence by the Schwarz inequality \cite[Proposition 3.3]{Pau} (or \cite[Proposition 9.9.4]{Pal2}), we have $T_t(x)^*T_t(x) \leq T_t(x^*x)$. Recall that the positive cone of $\VN(G)$ is weak* closed (see e.g. \cite[(2.3)]{ArK1}). Using the point 1, we conclude that $\Gamma(x,x) \geq 0$. 

3. For $x, y \in \VNGfin$, using the point 1 in the first equality and the point 2 in the last inequality, we obtain 
\begin{align*}
\MoveEqLeft
\begin{bmatrix} 
\Gamma(x,x) & \Gamma(x,y) \\ 
\Gamma(y,x) & \Gamma(y,y) 
\end{bmatrix} 
\ov{\eqref{Equa-Gamma-xy-sgrp}}{=}\lim_{t \to 0^+} \frac{1}{2t} 
\begin{bmatrix} 
 T_t(x^*x) - (T_t(x))^*T_t(x)  &  T_t(x^*y) - (T_t(x))^*T_t(y)  \\ 
 T_t(y^*x) - (T_t(y))^*T_t(x)  &   T_t(y^*y) - (T_t(y))^*T_t(y) 
\end{bmatrix}\\
&=\lim_{t \to 0^+} \frac{1}{2t} \left(\begin{bmatrix} 
T_t(x^*x) & T_t(x^*y) \\ 
T_t(y^*x) & T_t(y^*y) 
\end{bmatrix} - \begin{bmatrix} 
(T_t(x))^*T_t(x) & (T_t(x))^* T_t(y) \\ 
(T_t(y))^*T_t(x) & (T_t(y))^*T_t(y) 
\end{bmatrix} \right) \\
&=\lim_{t \to 0^+} \frac{1}{2t}\left( (\Id_{\M_2} \ot T_t) 
\left(\begin{bmatrix} 
x^*x & x^* y \\ 
y^* x & y^* y 
\end{bmatrix}\right)- \begin{bmatrix} 
(T_t(x))^* & 0\\ 
(T_t(y))^* & 0
\end{bmatrix} 
\begin{bmatrix} 
T_t(x) & T_t(y) \\ 
0 & 0 
\end{bmatrix} 
\right) \\
&=\lim_{t \to 0^+} \frac{1}{2t}\left( (\Id_{\M_2} \ot T_t) \left( 
\begin{bmatrix} 
x & y \\ 
0 & 0 
\end{bmatrix}^* 
\begin{bmatrix} 
x & y \\ 
0 & 0 
\end{bmatrix} \right) \right.\\
& \left.
-\left((\Id_{\M_2} \ot T_t)
\left(\begin{bmatrix} 
x & y \\ 
0 & 0 
\end{bmatrix}\right)\right)^* \left((\Id_{\M_2} \ot T_t)
\begin{bmatrix} 
x & y \\ 
0 & 0 \end{bmatrix} \right)\right)\\ 
&\geq 0.
\end{align*}

4. By \cite[Lemma 2.11]{ArK1}, 
$$ 
\norm{\Gamma(x,y)}_{\L^p(\VN(G))} 
\leq \frac{1}{2^{\frac1p}} \left( \norm{\Gamma(x,x)}_{\L^p(\VN(G))}^p + \norm{\Gamma(y,y)}_{\L^p(\VN(G))}^p \right)^{\frac1p}. 
$$
Since we have $\Gamma(\lambda x , \frac{1}{\lambda} y) = \Gamma(x,y)$ for $\lambda > 0$, we deduce

$$ 
\norm{\Gamma(x,y)}_{\L^p(\VN(G))} 
\leq \frac{1}{2^{\frac1p}} \inf_{\lambda > 0} \left(\lambda \norm{\Gamma(x,x)}^p_{\L^p(\VN(G))} + \frac{1}{\lambda}\norm{\Gamma(y,y)}^p_{\L^p(\VN(G))} \right)^{\frac1p}. 
$$
Ruling out beforehand the easy cases $\Gamma(x,x) = 0$ or $\Gamma(y,y) = 0$, we choose
$$
\lambda 
= \left(\norm{\Gamma(y,y)}^p_{\L^p(\VN(G))} / \norm{\Gamma(x,x)}^p_{\L^p(\VN(G))} \right)^{\frac12}
$$
and obtain
\begin{align}
\MoveEqLeft
\label{Inegalite-CS-Gamma-Lp-sgrp}
\norm{\Gamma(x,y)}_{\L^p(\VN(G))} 
\leq \frac{1}{2^{\frac1p}} \left( 2 \norm{\Gamma(x,x)}_{\L^p(\VN(G))}^{\frac{p}{2}} \norm{\Gamma(y,y)}_{\L^p(\VN(G))}^{\frac{p}{2}} \right)^{\frac1p} \\
&= \norm{\Gamma(x,x)}_{\L^p(\VN(G))}^{\frac12} \norm{\Gamma(y,y)}_{\L^p(\VN(G))}^{\frac12}. \nonumber  \end{align}
\end{proof}

Suppose that $(T_t)_{t \geq 0}$ is associated to the length $\psi \co G \to \R_+$ with associated cocycle $b = b_\psi \co G \to H$ and the orthogonal representation  $\pi \co G \to \B(H)$, $s \mapsto \pi_s$ of $G$ on $H$ from Proposition \ref{prop-Schoenberg}.  Suppose $-1 \leq q \leq 1$. For any $s \in G$, we will use the second quantization from \eqref{SQq} by letting
\begin{equation}
\label{equ-markovian-automorphism-group}
\alpha_s \ov{\mathrm{def}}{=} \Gamma^\infty_q(\pi_s) \co \Gamma_{q}(H) \to \Gamma_{q}(H)
\end{equation}
which is trace preserving. We obtain an action $\alpha \co G \to \Aut(\Gamma_{q}(H))$. So we can consider the crossed product $\Gamma_{q}(H) \rtimes_{\alpha} G$ as studied in Subsection \ref{subsubsec-Fourier-mult-crossed-product}, which comes equipped with its canonical normal finite faithful trace $\tau_\rtimes$.

Suppose $-1 \leq q \leq 1$ and $1 \leq p <\infty$. We introduce the map $\partial_{\psi,q} \co \P_G \to \L^p(\Gamma_{q}(H) \rtimes_\alpha G)$ defined by
\begin{equation}
\label{def-partial-psi}
\partial_{\psi,q}(\lambda_s)
=s_q(b_\psi(s)) \rtimes \lambda_s.
\end{equation}
which is a slight generalization of the map of \cite[p.~535]{JMP2}.

%

Note that $\L^p(\Gamma_{q}(H) \rtimes_\alpha G)$ is a $\VN(G)$-bimodule with left and right actions induced by
\begin{equation}
\label{def-bimodule-actions-VNG}
\lambda_s(z \rtimes \lambda_t)
\ov{\mathrm{def}}{=} \alpha_s(z) \rtimes \lambda_{st}
\quad \text{and} \quad
(z \rtimes \lambda_t)\lambda_s
\ov{\mathrm{def}}{=} z \rtimes \lambda_{ts}, \qquad z \in \Gamma_{q}(H),s,t \in G.	
\end{equation}
The following is stated in the particular case $q=1$ without proof in \cite[p.~544]{JMP2}. 
For the sake of completness, we give a short proof.

\begin{lemma}
\label{Lemma-Schur-Leibniz-nabla-sgrp}
Suppose $-1 \leq q \leq 1$. Let $G$ be a discrete group. For any $x,y \in \P_G$, we have
\begin{equation}
\label{Leibniz-Schur-gradient-mieux-sgrp}
\partial_{\psi,q}(xy)
=x\partial_{\psi,q}(y)+\partial_{\psi,q}(x)y.
\end{equation}
\end{lemma}

\begin{proof}
For any $s,t \in G$, we have
\begin{align*}
\MoveEqLeft
\lambda_s\partial_{\psi,q}(\lambda_t)+\partial_{\psi,q}(\lambda_s)\lambda_t     
\ov{\eqref{def-partial-psi}}{=} \lambda_s\big(s_q(b_\psi(t)) \rtimes \lambda_t\big)+\big(s_q(b_{\psi}(s)) \rtimes \lambda_s\big) \lambda_t\\
&\ov{\eqref{def-bimodule-actions-VNG}}{=} \alpha_s\big(s_q(b_\psi(t))\big) \rtimes \lambda_{st}+s_q(b_\psi(s)) \rtimes \lambda_{st}
=s_q\big(\pi_s(b_\psi(t))\big) \rtimes \lambda_{st}+s_q(b_\psi(s)) \rtimes \lambda_{st}\\
&\ov{\eqref{def-partial-psi}}{=} \big(s_q(\pi_s(b_{\psi}(t))+s_q(b_{\psi}(s))\big) \rtimes \lambda_{st}
\ov{\eqref{Cocycle-law}}{=} \big(s_q(b_\psi(st))\big) \rtimes \lambda_{st}=\partial_{\psi,q}(\lambda_{st})
=\partial_{\psi,q}(\lambda_s\lambda_t).
\end{align*}  
%
\end{proof}

The following is a slight generalization of \cite[Remark 1.3]{JMP2}. The proof is not difficult.\footnote{\thefootnote. On the one hand, for any $s,t \in G$, we have 
\begin{align*}
\MoveEqLeft
\Gamma(\lambda_s,\lambda_t)             
\ov{\eqref{Def-Gamma-lambda-s}}{=} \frac{1}{2}\big[ \psi(s^{-1}) + \psi(t) - \psi(s^{-1}t) \big]\lambda_{s^{-1}t}
\ov{\eqref{liens-psi-bpsi}}{=}\frac{1}{2}\Big[ \bnorm{b_\psi(s^{-1})}^2+\norm{b_\psi(t)}^2-\bnorm{b_\psi(s^{-1}t)}^2 \Big] \lambda_{s^{-1}t}\\
&\ov{\eqref{Cocycle-law}}{=}\frac{1}{2}\Big[ \bnorm{b_\psi(s^{-1})}^2+\norm{b_\psi(t)}^2-\bnorm{b_\psi(s^{-1})}^2- 2 \big\langle b_\psi(s^{-1}), \pi_{s^{-1}}(b_\psi(t)) \big\rangle_H-\norm{b_\psi(t)}^2 \Big] \lambda_{s^{-1}t}\\
&=-\big\langle b_\psi(s^{-1}), \pi_{s^{-1}}(b_\psi(t)) \big\rangle_H \lambda_{s^{-1}t}.
\end{align*}   
On the other hand, we have
\begin{align*}
\MoveEqLeft
\E\big[ \big(\partial_{\psi,q}(\lambda_s)\big)^* \partial_{\psi,q}(\lambda_t) \big]   
\ov{\eqref{def-partial-psi}}{=} \E \Big[ \big(s_q(b_\psi(s)) \rtimes \lambda_s\big)^*\big( s_q(b_\psi(t)) \rtimes \lambda_t\big)\Big] \\
&\ov{\eqref{inverse-crossed-product}}{=} \E \Big[ \big(\alpha_{s^{-1}}\big(s_q(b_\psi(s))\big) \rtimes \lambda_{s^{-1}}\big)\big(s_q(b_\psi(t)) \rtimes \lambda_t \big)\Big] 
=\E \Big[ \big(s_q(\pi_{s^{-1}}(b_\psi(s))) \rtimes \lambda_{s^{-1}}\big)\big(s_q(b_\psi(t)) \rtimes \lambda_t \big)\Big]\\
&\ov{\eqref{Cocycle-law}}{=} \E \Big[ \big(s_q(b_\psi(e)-b_\psi(s^{-1})) \rtimes \lambda_{s^{-1}}\big)\big(s_q(b_\psi(t)) \rtimes \lambda_t\big)\Big]
=-\E \Big[ \big(s_q(b_\psi(s^{-1})) \rtimes \lambda_{s^{-1}}\big)\big(s_q(b_\psi(t)) \rtimes \lambda_t\big)\Big]\\
&\ov{\eqref{Product-crossed-product}}{=}-\E \Big[ \Big(s_q(b_\psi(s^{-1}))\alpha_{s^{-1}}(s_q(b_\psi(t))) \Big)\rtimes \lambda_{s^{-1}t} \Big]
\ov{\eqref{Cocycle-law}}{=} -\E \Big[ \Big(s_q(b_\psi(s^{-1}))s_q\big( \pi_{s^{-1}}\big(b_\psi(t)\big)\big) \Big) \rtimes \lambda_{s^{-1}t} \Big] \\
&\ov{\eqref{petit-Wick}}{=}-\big\langle b_\psi(s^{-1}), \pi_{s^{-1}}\big(b_\psi(t)\big)\big\rangle_H \lambda_{s^{-1}t}.
\end{align*} } Here $\E \co \L^p(\Gamma_{q}(H) \rtimes_\alpha G) \to \L^p(\VN(G))$ denotes the canonical conditional expectation. 

\begin{prop}
\label{Prop-Fourier-grad-Gamma-Lp}
Suppose $-1 \leq q \leq 1$. For any $x,y \in \P_G$, we have 
\begin{equation}
\label{Equa-Fourier-grad-Gamma-Lp}
\Gamma(x,y)
=\E\big[ \big(\partial_{\psi,q}(x)\big)^* \partial_{\psi,q}(y) \big]. 	
\end{equation}
\end{prop}

Suppose $1 \leq p<\infty$. The following equalities are in \cite[pp.~930-931]{JRZ1} for $q=1$. For any $x,y \in \P_G$, we have
\begin{equation}
\label{lien-Gamma-partial-psi-junge}
\Gamma(x,y)
=\big\langle \partial_{\psi,q}(x),\partial_{\psi,q}(y) \big\rangle_{\L^p(\VN(G),\L^2(\Gamma_{q}(H))_{c,p})},
\end{equation}
and
\begin{equation}
\label{nabla-norm-Lp-gradient-psi-Junge}
\bnorm{\Gamma(x,x)^{\frac{1}{2}}}_{\L^p(\VN(G))}
=\bnorm{\partial_{\psi,q}(x)}_{\L^p(\VN(G),\L^2(\Gamma_{q}(H))_{c,p})}.
\end{equation}

We shall also need the following Riesz transform norm equivalence for markovian semigroup of Fourier multipliers from \cite[Theorem A2, Remark 1.3]{JMP2}.

\begin{prop}
\label{prop-Riesz-Fourier}
Let $G$ be a discrete group and $(T_t)_{t \geq 0}$ a markovian semigroup of Fourier multipliers with symbol $\psi$ of the negative generator $A$. Then for $2 \leq p < \infty$, we have the norm equivalence

\begin{equation}
\label{Equivalence-JMP2}
\bnorm{A^{\frac12}(x)}_{\L^p(\VN(G))} 
\cong_p \max \left\{ \bnorm{\Gamma(x,x)^{\frac12}}_{\L^p(\VN(G))}, \bnorm{\Gamma(x^*,x^*)^{\frac12} }_{\L^p(\VN(G))} \right\}, \quad x \in \P_G.
\end{equation}
\end{prop}

\subsection{Carr\'e du champ $\Gamma$ and first order differential calculus for Schur multipliers}
\label{Sec-infos-on-Gamma}

This subsection on markovian semigroups of Schur multipliers is the analogue of Subsection \ref{Sec-infos-on-Gamma-Fourier} on Fourier multipliers. We suppose here that we are given a markovian semigroup of Schur multipliers as in Proposition \ref{def-Schur-markovian}, and fix the associated family $(\alpha_i)_{i \in I}$ from \eqref{equ-Schur-markovian-alpha}.

For any $x,y \in \M_{I,\fin}$, we define the element 
\begin{equation}
\label{Def-Gamma}
\Gamma(x,y) 
\ov{\textrm{def}}{=} \frac{1}{2}\big[ A(x^*) y + x^* A(y) - A(x^* y)\big].	
\end{equation}
of $\M_{I,\fin}$. 
As in \eqref{Def-Gamma-sgrp} in the case of a markovian semigroup of Fourier multipliers, $\Gamma$ is called the carr\'e du champ, associated with the semigroup $(T_t)_{t \geq 0}$.
For any $i,j,k,l \in I$, note that
\begin{align}
\MoveEqLeft
\label{Def-Gamma-eij}
\Gamma(e_{ij},e_{kl})
=\frac{1}{2}\big[ A(e_{ji}) e_{kl} +  e_{ji}A(e_{kl}) - A(e_{ji} e_{kl})\big]             
=\frac{1}{2}\delta_{i=k}\big[ a_{ji} +  a_{kl} - a_{jl}\big]e_{jl}.
\end{align}
We refer to the books \cite{BGL1} and \cite{BoH1} for the general case of semigroups acting on classical $\L^p$-spaces and to \cite{Mey2}. For the semigroups acting on noncommutative $\L^p$-spaces, the carr\'e du champ and its applications were studied in the papers \cite[Section 9]{CiS1}, \cite{Cip1}, \cite{JM1}, \cite{Sau1} and \cite{JuZ3} mainly in the $\sigma$-finite case and for $\L^2$-spaces (see \cite{DaL1} for related things). Unfortunately, these papers do not suffice for our setting. So, in the following, we describe an elementary, independent and more concrete approach. However, some results of this section could be known, at least in the case $p=2$ and $I$ finite. 

The first part of the following result shows how to construct the carr\'e du champ directly from the semigroup $(T_t)_{t \geq 0}$.
The proof is similar to the one of Lemma \ref{Lemma-gamma-infos-sgrp} to which we refer.

\begin{lemma}
\label{Lemma-gamma-infos}
\begin{enumerate}
	\item Suppose $1 \leq p \leq \infty$. For any $x, y \in \M_{I,\fin}$, we have
\begin{equation}
	\label{Equa-Gamma-xy}
\Gamma(x,y) 
=\lim_{t \to 0^+} \frac1{2t} \big( T_t(x^*y) - (T_t(x))^*T_t(y) \big)
\end{equation}
in $S^p_I$.
	\item For any $x \in \M_{I,\fin}$, we have $\Gamma(x,x) \geq 0$.
	\item For any $x, y \in \M_{I,\fin}$, the matrix 
$\begin{bmatrix} 
\Gamma(x,x) & \Gamma(x,y) \\ 
\Gamma(y,x) & \Gamma(y,y) 
\end{bmatrix}$ is positive.
\item For any $x,y \in \M_{I,\fin}$, we have
\begin{equation}
\label{Inegalite-carre-Schur}
\norm{\Gamma(x,y)}_{S^p_I} 
\leq \norm{\Gamma(x,x)}_{S^p_I}^{\frac12} \norm{\Gamma(y,y)}_{S^p_I}^{\frac12}.
\end{equation}
This inequality even holds for $0 < p \leq \infty$.
\end{enumerate}
\end{lemma}

We shall link, more profoundly in Section \ref{sec:Semigroups-selfadjoint-contractive-Schur-multipliers}, the (square root of the) negative generator $A_p$, the carr\'e du champ $\Gamma$ and some noncommutative gradient from \eqref{def-delta-alpha}.
Suppose $1 \leq p< \infty$ and $-1 \leq q \leq 1$. Note that the Bochner space $\L^p(\Gamma_q(H) \otvn \B(\ell^2_I))$ is equipped with a canonical structure of $S^\infty_I$-bimodule whose operations are defined by $(f \ot x)y \ov{\textrm{def}}{=} f \ot xy$ and $y(f \ot x) \ov{\textrm{def}}{=} f \ot yx$ where $x \in S^p_I$, $y \in S^\infty_I$ and $f \in \L^p(\Gamma_q(H))$.
If $\Gamma_q$ is the $q$-Gaussian functor of Subsection \ref{Sec-q-gaussians} associated with the real Hilbert space $H$ stemming from Definition \ref{def-Schur-markovian}, we can consider the linear map $\partial_{\alpha,q} \co \M_{I,\fin} \to \Gamma_q(H) \ot \M_{I,\fin}$ (resp. $\partial_{\alpha,1} \co \M_{I,\fin} \to \L^0(\Omega) \ot \M_{I,\fin}$ if $q=1$) defined by 
\begin{equation}
\label{def-delta-alpha}
\partial_{\alpha,q}(e_{ij}) 
\ov{\mathrm{def}} 
=s_q(\alpha_i-\alpha_j) \ot e_{ij}, \qquad i,j \in I.
\end{equation} 

We have the following Leibniz rule.

\begin{lemma}
\label{Lemma-Schur-Leibniz-nabla}
Suppose $-1 \leq q\leq 1$. For any $x,y \in \M_{I,\fin}$, we have
\begin{equation}
\label{Leibniz-Schur-gradient-mieux}
\partial_{\alpha,q}(xy)
=x\partial_{\alpha,q}(y)+\partial_{\alpha,q}(x)y.
\end{equation}
\end{lemma}

\begin{proof}
On the one hand, for any $i,j,k,l \in I$, we have
$$
\partial_{\alpha,q}(e_{ij}e_{kl})
=\delta_{j=k}\partial_{\alpha,q}(e_{il})
\ov{\eqref{def-delta-alpha}}{=}\delta_{j=k}s_q(\alpha_i-\alpha_l) \ot e_{il}.
$$ 
On the other hand, for any $i,j,k,l \in I$, we have
\begin{align*}
\MoveEqLeft
 e_{ij}\partial_{\alpha,q}(e_{kl})+\partial_{\alpha,q}(e_{ij})e_{kl}      
		=e_{ij}\big(s_q(\alpha_k-\alpha_l) \ot e_{kl}\big)+\big(s_q(\alpha_i-\alpha_j)\ot e_{ij}\big) e_{kl}\\
		&=s_q(\alpha_k-\alpha_l) \ot e_{ij}e_{kl}+s_q(\alpha_i-\alpha_j) \ot e_{ij}e_{kl} \\
		&=s_q(\alpha_k-\alpha_l) \ot \delta_{j=k}e_{il}+s_q(\alpha_i-\alpha_j) \ot \delta_{j=k}e_{il}\\
		&=\delta_{j=k}\big(s_q(\alpha_k-\alpha_l)+s_q(\alpha_i-\alpha_j)\big) \ot e_{il}
		=\delta_{j=k}s_q(\alpha_i-\alpha_l) \ot e_{il}.
\end{align*}  
The result follows by linearity.
\end{proof}

Now, we describe a connection between the carr\'e du champ and the map $\partial_{\alpha}$ which is analogous to the equality of \cite[Section 1.4]{Sau1} (see also \cite{Sau2}). For that, we introduce the canonical trace preserving normal faithful conditional expectation $\E \co \Gamma_q(H) \otvn \B(\ell^2_I) \to \B(\ell^2_I)$.

\begin{prop}
\label{Prop-Schur-grad-Gamma}\label{Lemma-Hilbert-module}
Suppose $-1 \leq q\leq 1$ and $1 \leq p<\infty$. 
\begin{enumerate}
\item For any $x,y \in \M_{I,\fin}$, we have 
\begin{equation}
\label{Equa-Schur-grad-Gamma} 
\Gamma(x,y)
=\E\big(\big(\partial_{\alpha,q}(x)\big)^* \partial_{\alpha,q}(y)\big)
=\big\langle \partial_{\alpha,q}(x),\partial_{\alpha,q}(y) \big\rangle_{S^{p}_I(\L^2(\Gamma_q(H))_{c,p})}.
\end{equation}
\item For any $x \in \M_{I,\fin}$, we have
\begin{equation}
\label{nabla-norm-Lp-gradient}
\bnorm{\Gamma(x,x)^{\frac{1}{2}}}_{S^{p}_I}
=\bnorm{\partial_{\alpha,q}(x)}_{S^{p}_I(\L^2(\Gamma_q(H))_{c,p})},
\end{equation}
\end{enumerate}
\end{prop}

\begin{proof}
1. On the one hand, for any $i,j,k,l \in I$, we have 
\begin{align*}
\MoveEqLeft
\Gamma(e_{ij},e_{kl})             
\ov{\eqref{Def-Gamma-eij}}{=} \frac{1}{2}\delta_{i=k}\big[ a_{ji} +  a_{kl} - a_{jl}\big]e_{jl}
=\frac12 \delta_{i=k}\left( \norm{ \alpha_j - \alpha_i}_H^2 + \norm{\alpha_k - \alpha_l}_H^2 - \norm{\alpha_j - \alpha_l}_H^2 \right)e_{jl} \\
&=\frac12\delta_{i=k}\left( 2 \norm{\alpha_k}^2 - 2 \langle \alpha_i , \alpha_k \rangle - 2 \langle \alpha_j , \alpha_k \rangle + 2 \langle \alpha_i, \alpha_j \rangle \right)e_{jl} 
=\langle \alpha_k - \alpha_i, \alpha_k - \alpha_j \rangle_H e_{jl}.
\end{align*}   
On the other hand, we have
\begin{align*}
\MoveEqLeft
\E\big[ \big(\partial_{\alpha,q}(e_{ij})\big)^* \partial_{\alpha,q}(e_{kl}) \big]   
\ov{\eqref{def-delta-alpha}}{=} \E \Big[ \big(s_q(\alpha_i-\alpha_j) \ot e_{ij}\big)^*\big( s_q(\alpha_k-\alpha_l) \ot e_{kl}\big)\Big] \\
&=\E \Big[ s_q(\alpha_i-\alpha_j)s_q(\alpha_k-\alpha_l) \ot e_{ji}e_{kl}\Big] 
=\delta_{i=k}\tau( s_q(\alpha_i-\alpha_j)s_q(\alpha_k-\alpha_l))  e_{jl} \\
&=\delta_{i=k} \langle \alpha_i-\alpha_j, \alpha_k-\alpha_l\rangle_{H}  e_{jl}
\end{align*}
The second equality is a consequence of \eqref{Egalite-fantastique-2}.

2. If $x \in \M_{I,\fin}$, we have
\begin{align*}
\MoveEqLeft
 \bnorm{\partial_{\alpha,q}(x)}_{S^{p}_I(\L^2(\Gamma_q(H))_{c,p})}           
		\ov{\eqref{JMX-equa-2.9}}{=}\Bnorm{\big\langle\partial_{\alpha,q}(x),\partial_{\alpha,q}(x)\big\rangle_{S^{p}_I(\L^2(\Gamma_q(H))_{c,p})}^{\frac{1}{2}}}_{S^{p}_I}
		\ov{\eqref{Equa-Schur-grad-Gamma}}{=} \bnorm{\Gamma(x,x)^{\frac{1}{2}}}_{S^{p}_I}.
\end{align*}
  \end{proof}
\section{Riesz transforms associated to semigroups of Markov multipliers}
\label{sec:Semigroups-selfadjoint-contractive-Schur-multipliers}

\subsection{Khintchine inequalities for $q$-Gaussians in crossed products}
\label{twisted-Khintchine}

In this subsection, we consider a markovian semigroup of Fourier multipliers on $\VN(G)$, where $G$ is a discrete group satisfying Proposition \ref{prop-Schoenberg} with cocycle $(b_\psi,\pi,H)$. Moreover, we have seen in \eqref{equ-markovian-automorphism-group} that by second quantization we have an action $\alpha \co G \to \Aut(\Gamma_q(H))$. The aim of the subsection is to prove Theorem \ref{Khintchine-group-twisted} below which generalizes \cite[Theorem 1.1]{JMP2}.
In the following, the conditional expectation is $\E \co \Gamma_q(H) \rtimes_\alpha G \to \Gamma_q(H) \rtimes_\alpha G,\: x \rtimes \lambda_s \mapsto \tau_{\Gamma_q(H)}(x) \lambda_s$.
We let 
\begin{equation}
\label{Def-Gaussqp-group}
\Gauss_{q,p,\rtimes}(\L^p(\VN(G)))
\ov{\mathrm{def}}{=} \ovl{\Span \big\{ s_q(h) \rtimes x : h \in H,x \in \L^p(\VN(G)) \big\}}
\end{equation}
where the closure is taken in $\L^p(\Gamma_q(H) \rtimes_\alpha G)$ (for the weak* topology if $p=\infty$ and $-1 \leq q<1$).

\begin{lemma}
\label{lem-weakly-amenable-Gauss-density}
Let $1 < p < \infty$.
Let $G$ be a discrete group.
If $(e_k)_{k \in K}$ is an orthonormal basis of the Hilbert space $H$, then $\Gauss_{q,p,\rtimes}(\L^p(\VN(G)))$ equals the closure of the span of the $s_q(e_k) \rtimes \lambda_s$'s with $k \in K$ and $s \in G$.
\end{lemma}

\begin{proof}
We write temporarily $G' \subseteq \Gauss_{q,p}(\L^p(\VN(G)))$ for the closure of the span of the $s_q(e_k) \rtimes \lambda_s$. First note that it is trivial that $G' \subseteq \Gauss_{q,p}(\L^p(\VN(G)))$. For the reverse inclusion, by linearity and density and \eqref{Def-Gaussqp-group}, it suffices to show that $s_q(h) \rtimes \lambda_s$ belongs to $G'$ for any $h \in H$ and any $s \in G$. Let $\epsi > 0$ be fixed. Then there exists some finite subset $F$ of $K$ and $\alpha_k \in \C$ for $k \in F$ such that $\norm{h-\sum_{k \in F} \alpha_k e_k}_H < \epsi$. We have $\sum_{k \in F} \alpha_k s_q(e_k) \rtimes \lambda_s \in G'$ and\footnote{\thefootnote. Recall that $\alpha$ preserves the trace. If $s \in G$ and $k \in H$, we have 
\begin{align*}
\MoveEqLeft
\norm{s_q(h) \rtimes \lambda_s}_{\L^p(\Gamma_q(H) \rtimes_\alpha G)}^p            
=\tau_\rtimes\big((s_q(h) \rtimes \lambda_s)^*(s_q(h) \rtimes \lambda_s)\big)^\frac{p}{2})\ov{\eqref{inverse-crossed-product}}{=} \tau_\rtimes\big((\alpha_{s^{-1}}(s_q(h)) \rtimes \lambda_s)(s_q(h) \rtimes \lambda_s))^\frac{p}{2}\big)\\
&\ov{\eqref{Product-crossed-product}}{=} \tau_\rtimes\big((\alpha_{s^{-1}}(s_q(h)^2) \rtimes \lambda_e)^\frac{p}{2}\big) 
=\tau_\rtimes\big(\pi\big(\alpha_{s^{-1}}(s_q(h)^2)\big)^\frac{p}{2}\big)
=\tau\big(\alpha_{s^{-1}}\big((s_q(h)^{p})\big)\big) \\
&=\tau\big(s_q(h)^p\big)
=\norm{s_q(h)}_{\L^p(\Gamma_q(H))}^p.
\end{align*}.}
\begin{align*}
\MoveEqLeft
\norm{s_q(h) \rtimes \lambda_s - \sum_{k \in F} \alpha_k s_q(e_k) \rtimes \lambda_s}_{\L^p(\Gamma_q(H) \rtimes_\alpha G)} 
= \norm{s_q\bigg(h - \sum_{k \in F}\alpha_k e_k\bigg) \rtimes \lambda_s}_{\L^p(\Gamma_q(H) \rtimes_\alpha G)} \\
&=\norm{s_q\bigg(h - \sum_{k \in F}\alpha_k e_k\bigg)}_{\L^p(\Gamma_q(H))} 
\cong \norm{s_q\bigg(h - \sum_{k \in F}\alpha_k e_k\bigg)}_{\L^2(\Gamma_q(H))}
< \epsi.
\end{align*}
We conclude $\Gauss_{q,p,\rtimes}(\L^p(\VN(G))) \subseteq G'$ by closedness of $G'$.
\end{proof}

If $2 \leq p<\infty$, recall that we have a contractive injective inclusion 
\begin{equation}
\label{Inclusion-sympa-group}
\Gauss_{q,p,\rtimes}(\L^p(\VN(G))) \subseteq \L^p(\Gamma_{q}(H) \rtimes_\alpha G) \overset{\text{Lemma }\ref{lem-contractive-crossed-product-Hilbert-valued}}{\subset} \L^p_{cr}(\E)
\end{equation}
and that if $1<p<2$ we have an inclusion 
\begin{equation}
\label{inclusion-mechante-group}
\Span\{ s_q(h) \rtimes x : \: h \in H,\: x \in \L^p(\VN(G)) \}
\subseteq \L^p(\VN(G),\L^2(\Gamma_q(H)_{c,p})
\ov{\eqref{Egalite-fantastique-group}}{=} \L^p_{c}(\E) \subseteq \L^p_{cr}(\E).
\end{equation}
Here, the last but one equality follows from \eqref{Egalite-fantastique-group}. In Theorem \ref{Khintchine-group-twisted} below, we shall need several operations on $\L^p_c(\E),\:\L^p_r(\E)$ and $\L^p(\Gamma_q(H) \rtimes_\alpha G)$.
To this end, we have the following result.

\begin{lemma}
\label{lem-repair-row}
Let $-1 \leq q \leq 1$ and $1 < p < \infty$. Consider $P \co \L^2(\Gamma_q(H)) \to \L^2(\Gamma_q(H))$ the orthogonal projection onto $\Gauss_{q,2}(\C) = \Span \{ s_q(e): \: e \in H \}$. Let $G$ be a discrete group. Consider $Q_p =   P \rtimes \Id_{\L^p(\VN(G))}$ initially defined on $\P_{\rtimes,G}$. Then $Q_p$ induces well-defined contractions
\begin{equation}
\label{equ-2-repair-row}
Q_p  \co \L^p_c(\E) \to \L^p_c(\E)
\quad \text{and} \quad
Q_p  \co \L^p_r(\E) \to \L^p_r(\E).
\end{equation}
Consequently, according to Proposition \ref{prop-compatibility-crossed}, $Q_p$ equally extends to a contraction on $\L^p_{cr}(\E)$.
\end{lemma}

\begin{proof}
By \eqref{extens-col-1} and \eqref{Egalite-fantastique-group}, the column part of \eqref{equ-2-repair-row} follows immediately. We turn to row part of \eqref{equ-2-repair-row}. Note that the immediate analogue of \eqref{extens-col-1} to the row case holds true, but not for \eqref{Egalite-fantastique-group}, according to Remark \ref{rem-problem-row}. So we need to argue differently.

We start by showing that for any $h \in \L^2(\Gamma_q(H))$ and any $s \in G$ we have
\begin{equation}
\label{equ-repair-row}
\alpha_s\big( (P(h))^* \big) 
= P\big(\alpha_s (h^*)\big).
\end{equation}
Indeed, by anti-linearity and $\L^2(\Gamma_q(H))$ continuity of both sides, it suffices to show \eqref{equ-repair-row} for $h=w(e_{n_1} \ot \cdots \ot e_{n_N})$ the Wick word of $e_{n_1} \ot \cdots \ot e_{n_N}$, where $(e_k)$ is an orthonormal basis of $H$. But then we have\footnote{\thefootnote. In the last equality, we use that $w(e_{n_1} \ot \cdots \ot e_{n_N})=w(e_{n_N} \ot \cdots \ot e_{n_1})^*$ \cite[p.~21]{Nou2}.}
\begin{align*}
\MoveEqLeft
\alpha_s (P(h)^*) 
= \alpha_s (\delta_{N = 1} s_q(e_{n_1})^*) 
= \delta_{N = 1} \alpha_s(s_q(e_{n_1})) 
\ov{\eqref{equ-markovian-automorphism-group}}{=} \delta_{N = 1} s_q(\pi_s(e_{n_1})) 
= P(\alpha_s(h^*)).
\end{align*}
We deduce that
\begin{align}
\MoveEqLeft
\norm{ \big(P \rtimes \Id_{\L^{p}(\VN(G))}\big) \bigg(\sum_k  h_k \rtimes \lambda_{s_k}  \bigg) }_{\L^{p}_r(\E)}
=\norm{\bigg( \big(P \rtimes \Id_{\L^{p}(\VN(G))}\big) \bigg(\sum_k h_k \rtimes  \lambda_{s_k}\bigg) \bigg)^*}_{\L^{p}_c(\E)} \nonumber \\
&=\norm{ \bigg( \sum_k P(h_k) \rtimes \lambda_{s_k}\bigg)^* }_{\L^{p}_c(\E)} \nonumber 
\overset{\eqref{inverse-crossed-product}}{=} \norm{ \sum_k \alpha_{s_k^{-1}}(P(h_k)^*)  \rtimes \lambda_{s_k^{-1}}}_{\L^{p}_c(\E)} \nonumber \\
&\overset{\eqref{equ-repair-row}}{=} \norm{ \sum_k P\big(\alpha_{s_k^{-1}}(h_k^*)\big) \rtimes  \lambda_{s_k^{-1}} }_{\L^{p}_c(\E)} \nonumber 
=\norm{\big(P \rtimes \Id_{\L^p(\VN(G))}\big) \bigg( \sum_k \alpha_{s_k^{-1}}(h_k^*) \rtimes \lambda_{s_k^{-1}}\bigg)}_{\L^{p}_c(\E)} \nonumber \\
&\overset{\eqref{inverse-crossed-product}}{=} \norm{\big(P \rtimes \Id_{\L^{p}(\VN(G))}\big) \left( \bigg( \sum_k h_k \rtimes \lambda_{s_k}\bigg)^* \right)}_{\L^{p}_c(\E)} \nonumber 
\overset{\eqref{equ-2-repair-row}}{\leq} \norm{\bigg( \sum_k h_k \rtimes \lambda_{s_k}   \bigg)^*}_{\L^{p}_c(\E)} \nonumber \\
&=\norm{ \sum_k  h_k \rtimes \lambda_{s_k}}_{\L^{p}_r(\E)} \nonumber .
\end{align}
We have proved the row part of \eqref{equ-2-repair-row}.
\end{proof}

Now the noncommutative Khintchine inequalities can be rewritten under the following form.

\begin{thm} 
\label{Khintchine-group-twisted}
Consider $-1 \leq q \leq 1$ and $1 < p< \infty$.
Let $G$ be a discrete group.
	\begin{enumerate}
\item Suppose $1<p <2$. For any element $f=\sum_{s,h} f_{s,h} s_q(h) \rtimes \lambda_s$ of $\Span \big\{ s_q(h) : h \in H \} \rtimes \P_G$, we have
\begin{align}
\MoveEqLeft
\label{Khintchine-group-p<2}
\norm{f}_{\Gauss_{q,p,\rtimes}(\L^p(\VN(G)))}  
\approx_p \inf_{f=g+h} \bigg\{\norm{\big(\E(g^*g)\big)^{\frac12}}_{\L^p(\VN(G))},\norm{\big(\E(h h^*) \big)^{\frac12}}_{\L^p(\VN(G))}\bigg\}.
\end{align}
Here the infimum runs over all $g \in \L^p_c(\E)$ and $h \in \L^p_r(\E)$ such that $f = g + h$.
Further one can restrict the infimum to all $g \in \Ran Q_p|_{\L^p_c(\E)},h \in \Ran Q_p|_{\L^p_r(\E)}$, where $Q_p$ is the mapping from Lemma \ref{lem-repair-row}.

Finally, in case that $\L^p(\VN(G))$ has $\CCAP$ and $\VN(G)$ has $\QWEP$, the infimum can be taken over all $g,h \in \Span \big\{ s_q(e) : e \in H \} \rtimes \P_G$.

\item Suppose $2 \leq p <\infty$. For any element $f=\sum_{s,h} f_{s,h} s_q(h) \rtimes \lambda_s$ of $\Gauss_{q,p,\rtimes}(\L^p(\VN(G)))$ with $f_{s,h} \in \C$, we have
	\begin{align}
\MoveEqLeft
 	\label{Khintchine-group-p>2-Thm332}
\max\bigg\{\norm{\big(\E(f^* f)\big)^{\frac12}}_{\L^p(\VN(G))},\norm{\big(\E(f f^*) \big)^{\frac12}}_{\L^p(\VN(G))} \bigg\}
\leq \norm{f}_{\Gauss_{q,p,\rtimes}(\L^p(\VN(G)))}    \\
		&\lesssim \sqrt{p} \underbrace{\max\bigg\{\norm{\big(\E(f^* f)\big)^{\frac12}}_{\L^p(\VN(G))},\norm{\big( \E(f f^*) \big)^{\frac12}}_{\L^p(\VN(G))} \bigg\}}_{\norm{f}_{\L^p_{cr}(\E)}}. \nonumber
\end{align}	
\end{enumerate}

\end{thm}

\begin{proof}
2. Suppose $2 <p <\infty$. We begin by proving the upper estimate of \eqref{Khintchine-group-p>2-Thm332}. Fix $m \geq 1$. We have\footnote{\thefootnote. If $h \in H$, we have
$$
\norm{\frac{1}{\sqrt{m}} \sum_{l=1}^{m}  e_l \ot h}_{\ell^2_m(H)}
=\frac{1}{\sqrt{m}}\norm{\sum_{l=1}^{m}  e_l}_{\ell^2_m} \norm{h}_{H}
=\norm{h}_{H}.
$$} an isometric embedding $J_m  \co H \to \ell^2_m(H)$ defined by
\begin{equation}
\label{Def-Jm-group}
J_m(h)
\ov{\mathrm{def}}{=} \frac{1}{\sqrt{m}} \sum_{l=1}^{m}  e_l \ot h.
\end{equation}
Hence we can consider the associated operator $\Gamma_q^p(J_m) \co \L^p(\Gamma_q(H)) \to \L^p(\Gamma_q(\ell^2_m(H)))$ of second quantization \eqref{SQq} which is an isometric completely positive map. We define the maps $\pi^m_s \co G \to \B(\ell^2_m(H))$, $s \mapsto (e_l \ot h \mapsto e_l \ot \pi_s(h))$ and $\alpha^m_s\ov{\mathrm{def}}{=} \Gamma_q(\pi_s^m) \co \Gamma_q(\ell^2_m(H)) \to \Gamma_q(\ell^2_m(H))$. For any $h \in H$, any $s \in G$, note that
\begin{align}
\label{action-commute}
\MoveEqLeft
\pi_s^m \circ J_m(h)            
\ov{\eqref{Def-Jm-group}}{=} \pi_s^m\bigg(\frac{1}{\sqrt{m}} \sum_{l=1}^{m}  e_l \ot h\bigg)
=\frac{1}{\sqrt{m}} \sum_{l=1}^{m} \pi_s^m(e_l \ot h) 
=\frac{1}{\sqrt{m}} \sum_{l=1}^{m} e_l \ot \pi_s(h) \\
&\ov{\eqref{Def-Jm-group}}{=} J_m \circ \pi_s(h). 
\end{align} 
We deduce that
\begin{align*}
\MoveEqLeft
\alpha^m_s \circ \Gamma_q(J_m)          
=\Gamma_q(\pi_s^m) \circ \Gamma_q(J_m)
=\Gamma_q(\pi_s^m \circ J_m)
\ov{\eqref{action-commute}}{=} \Gamma_q(J_m \circ \pi_s)
=\Gamma_q(J_m) \circ \alpha_s.
\end{align*} 
By Lemma \ref{lem-van-Daele-Arhancet}, we obtain a trace preserving unital normal injective $*$-homomorphism $\Gamma_q(J_m) \rtimes \Id_{\VN(G)} \co \Gamma_q(H) \rtimes_{\alpha} G \to \Gamma_q(\ell^2_m(H)) \rtimes_{\alpha^m} G$. This map induces an isometric map $\rho \ov{\mathrm{def}}{=}\Gamma_q^p(J_m) \rtimes_{\alpha} \Id_{\L^p(\VN(G))} \co \L^p(\Gamma_q(H) \rtimes_{\alpha} G) \to \L^p(\Gamma_q(\ell^2_m(H)) \rtimes_{\alpha^m} G)$. For any finite sum $f=\sum_{s,h} f_{s,h} s_q(h) \rtimes \lambda_s$ of $\P_{\Gamma_q(H),G}$ with $f_{s,h} \in \C$, we obtain
\begin{align}
\MoveEqLeft
 \norm{f}_{\L^p(\Gamma_q(H) \rtimes_{\alpha} G)} 
=\norm{\sum_{s,h} f_{s,h} s_q(h) \rtimes \lambda_s}_{\L^p(\Gamma_q(H) \rtimes_{\alpha} G)} \\
&=\norm{\rho\bigg(\sum_{s,h} f_{s,h} s_q(h) \rtimes \lambda_s\bigg)}_{\L^p(\Gamma_q(\ell^2_m(H)) \rtimes_{\alpha^m} G)} \nonumber\\
&=\norm{\sum_{s,h} f_{s,h} \Gamma_q(J_m)(s_q(h)) \rtimes_{\alpha^m} \lambda_s}_{\L^p(\Gamma_q(\ell^2_m(H)) \rtimes_{\alpha^m} G)} \nonumber \\
&\ov{\eqref{SQq}}{=}\norm{\sum_{s,h} f_{s,h} s_{q,m}(J_m(h)) \rtimes \lambda_s}_{\L^p(\Gamma_q(\ell^2_m(H)) \rtimes_{\alpha^m} G)}\nonumber\\
&\ov{\eqref{Def-Jm-group}}{=}\norm{\sum_{s,h} f_{s,h} s_{q,m}\bigg(\frac{1}{\sqrt{m}} \sum_{l=1}^{m} e_l \ot h\bigg) \rtimes \lambda_s}_{\L^p(\Gamma_q(\ell^2_m(H)) \rtimes_{\alpha^m} G)} \nonumber \\
&=\frac{1}{\sqrt{m}} \Bgnorm{\sum_{s,h,l} f_{s,h} s_{q,m}(e_l \ot h) \rtimes \lambda_s}_{\L^p(\Gamma_q(\ell^2_m(H)) \rtimes_{\alpha^m} G)}.  \label{Equa-inter-4532-group}
\end{align}
For any $1 \leq l \leq m$, consider the element $f_l \ov{\mathrm{def}}{=} \sum_{s,h}  f_{s,h} s_{q,m}(e_l \ot h) \rtimes \lambda_s$ and the conditional expectation $\E \co \L^p(\Gamma_q(\ell^2_m(H)) \rtimes_{\alpha^m} G) \to \L^p(\VN(G))$. For any $1 \leq l \leq m$, note that
\begin{align*}
\MoveEqLeft
\E(f_l)   
=\E\bigg(\sum_{s,h} f_{s,h} s_{q,m}(e_l \ot h) \rtimes \lambda_s\bigg)
=\sum_{s,h} f_{s,h} \E\big(s_{q,m}(e_l \ot h)\rtimes \lambda_s\big)\\
&=\sum_{s,h} f_{s,h} \tau(s_{q,m}(e_l \ot h)) \lambda_s
\ov{\eqref{formule de Wick odd}}{=} 0.
\end{align*} 
We deduce that the random variables $f_l$ are mean-zero. Now, we prove in addition that the $f_l$ are independent over $\VN(G)$, so that we will be able to apply the noncommutative Rosenthal inequality \eqref{Inequality-Rosenthal-1-prime} to them.

\begin{lemma}
\label{lem-fl-faithfully-independent-Thm332}
The $f_l$ defined here above are independent over $\VN(G)$. More precisely, the von Neumann algebras generated by $f_l$ are independent with respect to $\VN(G)$.
\end{lemma}

\begin{proof}
Due to normality of $\E$, it suffices to prove 
\begin{equation}
\label{equ-1-lem-fl-faithfully-independent-Thm332}
\E(xy) 
=\E(x) \E(y)
\end{equation}
for $x$ (resp. $y$) belonging to a weak$^*$ dense subset of $\W^*(f_l)$ (resp. of $\W^*((f_k)_{k \neq l})$). Moreover, by bilinearity in $x,y$ of both sides of \eqref{equ-1-lem-fl-faithfully-independent-Thm332} and selfadjointness of $s_{q,m}(e_l \ot h)$, it suffices to take $x = s_{q,m}(e_l \ot h)^n \rtimes \lambda_s$ for some $n \in \N$ and $y = \prod_{t = 1}^T s_{q,m}(e_{k_t} \ot h)^{n_t} \rtimes \lambda_u$ for some $T \in \N$, $k_t \neq l$ and $n_t \in \N$. Since $\alpha_s^m=\Gamma_q(\pi_s^m)$ is trace preserving, we have 
\begin{align*}
\MoveEqLeft
\E(x) \E(y)
=\tau(s_{q,m}(e_l \ot h)^n)\lambda_s\, \tau\bigg(\prod_{t=1}^T s_{q,m}(e_{k_t} \ot h)^{n_t}\bigg) \lambda_u\\
&=\tau(s_{q,m}(e_l \ot h)^n)\tau\bigg(\prod_{t=1}^T s_{q,m}(e_{k_t} \ot h)^{n_t}\bigg) \lambda_{s u} \\
&=\tau(s_{q,m}(e_l \ot h)^n)\tau\bigg(\alpha_s^m\bigg(\prod_{t=1}^T s_{q,m}(e_{k_t} \ot h)^{n_t}\bigg)\bigg) \lambda_{s u}\\
&=\tau(s_{q,m}(e_l \ot h)^n)\tau\bigg(\prod_{t=1}^T s_{q,m}(e_{k_t} \ot \pi_s(h))^{n_t}\bigg) \lambda_{s u}
\end{align*} 
and 
\begin{align*}
\MoveEqLeft
\E(xy) 
=\E\bigg(\big(s_{q,m}(e_l \ot h)^n \rtimes \lambda_s\big)\bigg(\prod_{t = 1}^T s_{q,m}(e_{k_t} \ot h)^{n_t} \rtimes \lambda_u\bigg)\bigg) \\
&\ov{\eqref{Product-crossed-product}}{=} \E\bigg(\big(s_{q,m}(e_l \ot h)^n \alpha_s^m\bigg(\prod_{t=1}^T s_{q,m}(e_{k_t} \ot h)^{n_t}\bigg) \rtimes \lambda_{su}\big)\bigg)\\
&=\tau \left(s_{q,m}(e_l \ot h)^n \alpha_s^m\bigg(\prod_{t=1}^T s_{q,m}(e_{k_t} \ot h)^{n_t}\bigg) \right) \lambda_{su} \\
&=\tau\bigg(s_{q,m}(e_l \ot h)^n \prod_{t=1}^T s_{q,m}(e_{k_t} \ot \pi_s(h))^{n_t}\bigg) \lambda_{su}.
\end{align*} 
We shall now apply the Wick formulae \eqref{formule de Wick} and \eqref{formule de Wick odd} to the trace term above. 

Note that if $n + \sum_{t = 1}^T n_t$ is odd, then according to the Wick formula \eqref{formule de Wick odd} we have $\E(xy)=0$. On the other hand, then either $n$ or $\sum_{t = 1}^T n_t$ is odd, so according to the Wick formula \eqref{formule de Wick odd}, either $\E(x) = 0$ or $\E(y)=0$. Thus, \eqref{equ-1-lem-fl-faithfully-independent-Thm332} follows in this case. 

Now suppose that $n + \sum_{t = 1}^T n_t\ov{\mathrm{def}}{=} 2k$ is even. Consider a $2$-partition $\mathcal{V} \in \mathcal{P}_2(2k)$. If both $n$ and $\sum_{t = 1}^T n_t$ are odd, then we must have some $(i,j) \in \mathcal{V}$ such that one term $\langle f_i, f_j \rangle_{\ell^2_m(H)}$ in the Wick formula \eqref{formule de Wick} equals $\big\langle e_l \ot h, e_{k_t} \ot \pi_s(h) \big\rangle = 0$, since $k_t \neq l$. Thus, $\E(xy) = 0$, and since $n$ is odd, also $\E(x) = 0$, and therefore, \eqref{equ-1-lem-fl-faithfully-independent-Thm332} follows. If both $n$ and $\sum_{t = 1}^T n_t$ are even, then in the Wick formula \eqref{formule de Wick}, we only need to consider those $2$-partitions $\mathcal{V}$ without a mixed term as the $\langle f_i, f_j\rangle = 0$ above. Such a $\mathcal{V}$ is clearly the disjoint union $\mathcal{V} = \mathcal{V}_1 \cup \mathcal{V}_2$ of $2$-partitions corresponding to $n$ and to $\sum_{t = 1}^T n_t$. Moreover, we have for the number of crossings, $c(\mathcal{V}) = c(\mathcal{V}_1) + c(\mathcal{V}_2)$. With $(f_1,f_2,\ldots,f_{n+T}) = (\underbrace{e_{l} \ot h ,e_l \ot h,\ldots,e_l \ot h}_{n},e_{k_1} \ot \pi_s(h),\ldots, e_{k_T} \ot \pi_s(h))$, we obtain 
\begin{align*}
\MoveEqLeft
\tau \left(s_{q,m}(e_l \ot h)^n \prod_{t = 1}^T s_{q,m}(e_{k_t} \ot \pi_s(h))^{n_t} \right) 
\ov{\eqref{formule de Wick}}{=} \sum_{\mathcal{V} \in \mathcal{P}_2(2k)} q^{c(\mathcal{V})} \prod_{(i,j) \in \mathcal{V}} \langle f_i, f_j \rangle_{\ell^2_m(H)} \\
&=\sum_{\mathcal{V}_1,\mathcal{V}_2} q^{c(\mathcal{V}_1) + c(\mathcal{V}_2)} \prod_{(i_1,j_1) \in \mathcal{V}_1} \langle f_{i_1}, f_{j_1} \rangle_{\ell^2_m(H)} \prod_{(i_2,j_2) \in \mathcal{V}_2} \langle f_{i_2}, f_{j_2} \rangle_{\ell^2_m(H)} \\
&=\bigg(\sum_{\mathcal{V}_1} q^{c(\mathcal{V}_1)} \prod_{(i_1,j_1) \in \mathcal{V}_1} \langle f_{i_1}, f_{j_1} \rangle_{\ell^2_m(H)}\bigg)\bigg(\sum_{\mathcal{V}_2} q^{c(\mathcal{V}_2)} \prod_{(i_2,j_2) \in \mathcal{V}_2} \langle f_{i_2}, f_{j_2} \rangle_{\ell^2_m(H)}\bigg) \\
&\ov{\eqref{formule de Wick}}{=} \tau \left(s_{q,m}(e_l \ot h)^n \right) \tau \left( \prod_{t = 1}^T s_{q,m}(e_{k_t} \ot \pi_s(h))^{n_t} \right).
\end{align*}
Thus, also in this case \eqref{equ-1-lem-fl-faithfully-independent-Thm332} follows.
\end{proof}

Now we are able to apply the noncommutative Rosenthal inequality \eqref{Inequality-Rosenthal-1-prime} which yields
\begin{align}
\label{Eq-divers-3212-group}
\MoveEqLeft
\norm{f}_{\L^p(\Gamma_q(H) \rtimes_{\alpha} G)} 
\ov{\eqref{Equa-inter-4532-group}}{=}\frac{1}{\sqrt{m}} \norm{\sum_{l=1}^{m} f_l}_{\L^p(\Gamma_q(\ell^2_m(H)) \rtimes_{\alpha^m} G)} \\
& \ov{\eqref{Inequality-Rosenthal-1-prime}}{\lesssim} \frac{1}{\sqrt{m}} \Bigg[ p\bigg( \sum_{l=1}^m \norm{f_l}_{\L^p(\Gamma_q(\ell^2_m(H)) \rtimes_{\alpha^m} G)}^p \bigg)^{\frac{1}{p}} \\
&+ \sqrt{p}\norm{ \bigg(\sum_{l=1}^m \E(f_l^* f_l) \bigg)^{\frac12}}_{\L^p(\VN(G))} +\sqrt{p}\norm{ \bigg(\sum_{l=1}^m \E(f_l f_l^*) \bigg)^{\frac12}}_{\L^p(\VN(G))} \Bigg]. \nonumber
\end{align}
For any integer $1 \leq l \leq m$, note that
\begin{align*}
\MoveEqLeft
\E(f_l^* f_l)  
=\E\Bigg(\bigg(\sum_{s,h} f_{s,h} s_{q,m}(e_l \ot h) \rtimes \lambda_{s}\bigg)^* \bigg(\sum_{t,k} f_{t,k} s_{q,m}(e_l \ot k) \rtimes \lambda_{t}\bigg)\Bigg)\\
&\ov{\eqref{inverse-crossed-product}}{=} \sum_{s,h,t,k} \ovl{f_{s,h}} f_{t,k} \E\Big(\big(\alpha_{s^{-1}}^m(s_{q,m}(e_l \ot h))\big) \rtimes \lambda_{s^{-1}}\big) \big(s_{q,m}(e_l \ot k) \rtimes \lambda_{t}\big)\Big)\\
&\ov{\eqref{Product-crossed-product}}{=} \sum_{s,h,t,k} \ovl{f_{s,h}} f_{t,k} \E\Big(\alpha_{s^{-1}}^m(s_{q,m}(e_l \ot h))\alpha_{s^{-1}}^m(s_{q,m}(e_l \ot k)) \rtimes \lambda_{s^{-1}t} \Big)\\\
&=\sum_{s,h,t,k} \ovl{f_{s,h}} f_{t,k} \tau(\alpha_{s^{-1}}^m(s_{q,m}(e_l \ot h)(s_{q,m}(e_l \ot k))) \lambda_{s^{-1}t} \\
&=\sum_{s,h,t,k} \ovl{f_{s,h}}f_{t,k} \tau\big(s_{q,m}(e_l \ot h)s_{q,m}(e_l \ot k)\big) \lambda_{s^{-1}t} 
\ov{\eqref{petit-Wick}}{=} \sum_{i,j,h,s} \ovl{f_{s,h}}f_{t,k} \langle h,k \rangle_H \lambda_{s^{-1}t}.
\end{align*} 
and
\begin{align*}
\MoveEqLeft
\E(f^* f)   
=\E\Bigg(\bigg(\sum_{s,h} f_{s,h} s_q(h)  \rtimes \lambda_s\bigg)^* \bigg(\sum_{t,k} f_{t,k} s_q(k) \rtimes \lambda_{t}\bigg)\Bigg) \\
&\ov{\eqref{inverse-crossed-product}}{=} \sum_{s,h,t,k}\ovl{f_{s,h}}f_{t,k}\E\Big(\big( \alpha_{s^{-1}}( s_q(h)) \rtimes \lambda_{s^{-1}}\big) \big(  s_q(k) \rtimes \lambda_{t}\big)\Big) \\ 
&\ov{\eqref{Product-crossed-product}}{=} \sum_{s,h,t,k} \ovl{f_{s,h}}f_{t,k}\E\Big(\big(\alpha_{s^{-1}}( s_q(h))\alpha_{s^{-1}}(s_q(k)) \rtimes \lambda_{s^{-1}t} \Big) \\
&=\sum_{s,h,t,k} \ovl{f_{s,h}}f_{t,k}\tau\big(\alpha_{s^{-1}}(s_q(h)s_q(k))\big) \lambda_{s^{-1}t}\\
&=\sum_{s,h,t,k} \ovl{f_{s,h}} f_{t,k}  \tau\big( s_q(h)s_q(k) \big) \lambda_{s^{-1}t}
\ov{\eqref{petit-Wick}}{=} \sum_{i,j,h,s} \ovl{f_{s,h}}f_{t,k} \langle h,k \rangle_H \lambda_{s^{-1}t}
\end{align*}
and similarly for the row terms. We conclude that 
\begin{equation}
\label{Eq-divers-98089-group}
\E(f_l^* f_l)
=\E(f^* f)
\quad \text{and} \quad 
\E(f_l f_l^*) 
= \E(f f^*).
\end{equation}
Moreover, for $1 \leq l \leq m$, using the isometric map $\psi_l \co H \to \ell^2_m(H)$, $h \to e_l \ot h$, we can introduce the second quantization operator $\Gamma_q(\psi_l) \co \Gamma_q(H) \to \Gamma_q(\ell^2_m(H))$. 
For any $h \in H$, any $s \in G$, note that
\begin{align}
\label{action-commute-prime}
\MoveEqLeft
\pi_s^m \circ \psi_l(h)            
=\pi_s^m(e_l \ot h)
=e_l \ot \pi_s(h)
=\psi_l \circ \pi_s(h). 
\end{align} 
We deduce that
\begin{align*}
\MoveEqLeft
\alpha^m_s \circ \Gamma_q(\psi_l)          
=\Gamma_q(\pi_s^m) \circ \Gamma_q(\psi_l)
=\Gamma_q(\pi_s^m \circ \psi_l)
\ov{\eqref{action-commute-prime}}{=} \Gamma_q(\psi_l \circ \pi_s)
=\Gamma_q(\psi_l) \circ \alpha_s.
\end{align*} 
By Lemma \ref{lem-van-Daele-Arhancet}, we obtain a trace preserving unital normal injective $*$-homomorphism $\Gamma_q(\psi_l) \rtimes \Id_{\VN(G)} \co \Gamma_q(H) \rtimes_{\alpha} G \to \Gamma_q(\ell^2_m(H)) \rtimes_{\alpha^m} G$. This map induces an isometric map $\Gamma_q^p(\psi_l) \rtimes_{\alpha} \Id_{\L^p(\VN(G))} \co \L^p(\Gamma_q(H) \rtimes_{\alpha} G) \to \L^p(\Gamma_q(\ell^2_m(H)) \rtimes_{\alpha^m} G)$. We have
\begin{align}
\label{fl-f-group}
\MoveEqLeft
  \norm{f_l}_{\L^p(\Gamma_q(\ell^2_m(H)) \rtimes_{\alpha^m} G)}  
		=\norm{\sum_{s,h} f_{s,h} s_{q,m}(e_l \ot h) \rtimes \lambda_{s}}_{\L^p(\Gamma_q(\ell^2_m(H)) \rtimes_{\alpha^m} G)} \\
		&\ov{\eqref{SQq}}{=} \norm{\sum_{s,h} f_{s,h} \Gamma_q(\psi_l)(s_{q,m}(h)) \rtimes \lambda_{s}}_{\L^p(\Gamma_q(\ell^2_m(H)) \rtimes_{\alpha^m} G)} 
	=\norm{\sum_{s,h} f_{s,h} s_q(h) \rtimes \lambda_{s}}_{\L^p(\Gamma_q(H) \rtimes G)} \nonumber \\
	&=\norm{f}_{\L^p(\Gamma_q(H) \rtimes_{\alpha} G)}. \nonumber
\end{align} 
 We infer that
\begin{align*}
\MoveEqLeft
\norm{f}_{\L^p(\Gamma_q(H) \rtimes_{\alpha} G)} 
\ov{\eqref{Eq-divers-3212-group}\eqref{fl-f-group}\eqref{Eq-divers-98089-group}}{\lesssim} \frac{1}{\sqrt{m}} \Bigg[ p\bigg(\sum_{l=1}^m \norm{f}_{\L^p(\Gamma_q(H) \rtimes_{\alpha} G)}^p \bigg)^{\frac{1}{p}}  \Bigg.\\
&\qquad \qquad \qquad \Bigg.+ \sqrt{p}\norm{ \bigg(\sum_{l=1}^m \E(f^* f) \bigg)^{\frac12}}_{\L^p(\VN(G))}+\sqrt{p}\norm{ \bigg(\sum_{l=1}^m \E(f f^*) \bigg)^{\frac12}}_{\L^p(\VN(G))} \Bigg]\\
&=\frac{1}{\sqrt{m}} \Bigg[p m^{\frac{1}{p}} \norm{f}_{\L^p(\Gamma_q(H) \rtimes_{\alpha} G)} \\
& \qquad \qquad + \sqrt{pm}\norm{\big(\E(f^* f)\big)^{\frac12}}_{\L^p(\VN(G))} +\sqrt{pm}\norm{\big( \E(f f^*) \big)^{\frac12}}_{\L^p(\VN(G))} \Bigg]\\
&=p m^{\frac{1}{p}-\frac12} \norm{f}_{\L^p(\Gamma_q(H) \rtimes_{\alpha} G)}+\sqrt{p}\norm{\big(\E(f^* f)\big)^{\frac12}}_{\L^p(\VN(G))}+\sqrt{p}\norm{\big( \E(f f^*) \big)^{\frac12}}_{\L^p(\VN(G))}.
\end{align*}
Since $p>2$, passing to the limit when $m \to \infty$ and noting that Rosenthal's inequality comes with an absolute constant not depending on the von Neumann algebras under consideration \cite[p.~4303]{JuZ2}, we finally obtain 
$$
\norm{f}_{\L^p(\Gamma_q(H) \rtimes_{\alpha} G)}  
\lesssim \sqrt{p}\bigg[\norm{\big(\E(f^* f)\big)^{\frac12}}_{\L^p(\VN(G))}+\norm{\big(\E(f f^*) \big)^{\frac12}}_{\L^p(\VN(G))} \bigg].
$$
Using the equivalence $\ell^1_2 \approx \ell^\infty_2$, we obtain the upper estimate of \eqref{Khintchine-group-p>2-Thm332}.

The lower estimate of \eqref{Khintchine-group-p>2-Thm332} holds with constant $1$ from the contractivity of the conditional expectation $\E$ on $\L^{\frac{p}{2}}(\Gamma_q(H) \rtimes_{\alpha} G)$: 
\begin{align*}
\MoveEqLeft
\max\bigg\{\norm{\big(\E(f^* f)\big)^{\frac12}}_{\L^p(\VN(G))},\norm{\big(\E(f f^*) \big)^{\frac12}}_{\L^p(\VN(G))} \bigg\}  \\  
		&= \max\bigg\{\bnorm{\E(f^* f)}_{\L^{\frac{p}{2}}(\VN(G))}^{\frac12},\bnorm{\E(f f^*)}_{\L^{\frac{p}{2}}(\VN(G))}^{\frac12} \bigg\} \\
		&\leq \max\Bigg\{\norm{f^* f}_{\L^{\frac{p}{2}}(\Gamma_q(H) \rtimes_{\alpha} G)}^{\frac12},\norm{f f^*}_{\L^{\frac{p}{2}}(\Gamma_q(H) \rtimes_{\alpha} G)}^{\frac12} \Bigg\}
		=\norm{f}_{\L^p(\Gamma_q(H) \rtimes_{\alpha} G)}.
\end{align*}
1. Now, let us consider the case $1 < p < 2$. We will proceed by duality as follows. Consider the Gaussian projection $Q_p$ from Lemma \ref{lem-repair-row} which is a contraction on $\L^p_{cr}(\E)$. Note that $\tau_{\rtimes}(f^*g) = \tau_{\rtimes}(Q_p(f)^*g) = \tau_{\rtimes}(f^*Q_{p^*}(g))$ for any $g \in \L^{p^*}(\Gamma_q(H) \rtimes_\alpha G)$ and the fixed $f$ from the statement of the theorem.
Indeed, this can be seen from selfadjointness of $P$ and \eqref{equ-repair-row}.
Note also that $\L^p_{cr}(\E)$ and $\L^{p^*}_{cr}(\E)$ are dual spaces with respect to each other according to Lemma \ref{lem-expectation-duality},
with duality bracket $\langle f, g \rangle = \tau_\rtimes(f^*g)$, i.e. the same as the duality bracket between $\L^p(\Gamma_q(H) \rtimes_\alpha G)$ and $\L^{p^*}(\Gamma_q(H) \rtimes_\alpha G)$.

Using this in the first two equalities, the contractivity $Q_{p^*} \co \L^p_{cr}(\E) \to \L^p_{cr}(\E)$ from Lemma \ref{lem-repair-row} in the third equality, and the upper estimate of \eqref{Khintchine-group-p>2-Thm332} together with the density result of Lemma \ref{lem-density-intersection-CBAP} in the last inequality, we obtain for any $f \in \Gauss_{q,2}(\C) \rtimes \P_{G}$
\begin{align*}
\MoveEqLeft
\norm{f}_{\L^p_{cr}(\E)} = \sup_{\norm{g}_{\L^{p^*}_{cr}(\E)} \leq 1} \big| \tau_{\rtimes}( f^* g ) \big| = \sup_{\norm{g}_{\L^{p^*}_{cr}(\E)} \leq 1} \big| \tau_{\rtimes}( f^* Q_{p^*}(g) ) \big| \\
&=\sup_{\stackrel{\norm{g}_{\L^{p^*}_{cr}(\E)} \leq 1 }{g \in \Ran Q_{p^*}}} \big| \tau_{\rtimes}(f^*g) \big|\\ 
&\leq \norm{f}_{\L^p(\Gamma_q(H) \rtimes_{\alpha} G)} \sup_{\stackrel{\norm{g}_{\L^{p^*}_{cr}(\E)} \leq 1}{g \in \Ran Q_{p^*}}} \norm{g}_{\L^{p^*}(\Gamma_q(H) \rtimes_{\alpha} G)} 
\ov{\eqref{Khintchine-group-p>2-Thm332}}{\lesssim_p} \norm{f}_{\L^p(\Gamma_q(H) \rtimes_{\alpha} G)}.
\end{align*}
Note that in the definition of the $\L^p_{cr}(\E)$-norm, the infimum runs over all $g \in \L^p_c(\E)$ and $h \in \L^p_r(\E)$ such that $f = g + h$.
However, we can write $f = Q_p(f) = Q_p(g) + Q_p(h)$ and $\norm{Q_p(g)}_{\L^p_c(\E)} \overset{\eqref{equ-2-repair-row}}{\leq} \norm{g}_{\L^p_c(\E)}$ and $\norm{Q_p(h)}_{\L^p_r(\E)} \overset{\eqref{equ-2-repair-row}}{\leq} \norm{h}_{\L^p_r(\E)}$. Thus, in the definition of the $\L^p_{cr}(\E)$-norm of $f$, we can restrict to $g, h \in \Ran Q_p$.

Our next goal is to restrict in the lower estimate of \eqref{Khintchine-group-p<2} to $g,h \in \Span\left\{ s_q(e):\: e \in H \right\} \rtimes \P_G$ with $f = g+h$, in the case where $f \in \Span\left\{s_q(e):\: e \in H \right\} \rtimes \P_G$, under the supplementary assumptions that $\L^p(\VN(G))$ has $\CCAP$ and that $\VN(G)$ has $\QWEP$. We consider the approximating net $(M_{\varphi_j})$ of Fourier multipliers. Moreover, let $f = g + h$ be a decomposition with $g \in \L^p_c(\E)$ and $h \in \L^p_r(\E)$ such that $\norm{g}_{\L^p_c(\E)} + \norm{h}_{\L^p_r(\E)} \leq 2 \norm{f}_{\L^p_{cr}(\E)}$. According to the above, we can already assume that $g,h$ belong to $\Ran Q_p$. Then we newly decompose for each $j$
\begin{align*}
\MoveEqLeft
f = (\Id_{\Gamma_q(H)} \rtimes M_{\varphi_j})(f) + \left( f - \Id_{\Gamma_q(H)} \rtimes M_{\varphi_j}(f) \right) \\
& = \underset{\in \L^p_c(\E)}{\underbrace{(\Id_{\Gamma_q(H)} \rtimes M_{\varphi_j})(g) + \left(f - (\Id_{\Gamma_q(H)} \rtimes M_{\varphi_j})(f) \right)}} + \underset{\in \L^p_r(\E)}{\underbrace{(\Id_{\Gamma_q(H)} \rtimes M_{\varphi_j})(h)}} .
\end{align*}
Note that $(\Id_{\Gamma_q(H)} \rtimes M_{\varphi_j})(g)$ and $(\Id_{\Gamma_q(H)} \rtimes M_{\varphi_j})(h)$ belong to $\Span\left\{s_q(e):\: e \in H \right\} \rtimes \P_G$ since $\varphi_j$ is of finite support. Moreover, also $f - \Id_{\Gamma_q(H)} \rtimes M_{\varphi_j}(f)$ belongs to that space. Now, we will control the norms in this new decomposition.
Recall that $f =  \sum_{s,e} f_{s,e} s_q(e) \rtimes \lambda_s$ with finite sums.
Since $\varphi_j(s)$ converges pointwise to $1$, there is $j$ such that 
\begin{align*}
\MoveEqLeft
\norm{f - \Id_{\Gamma_q(H)} \rtimes M_{\varphi_j}(f)}_{\L^p_c(\E)} \leq \sum_{s,e} |f_{s,e}| \norm{s_q(e) \rtimes \lambda_s - \varphi_j(s) s_q(e) \rtimes \lambda_s}_{\L^p_c(\E)} \\
& \leq \sum_{s,e} |f_{s,e}| \, |1 - \varphi_j(s)| \norm{s_q(e) \rtimes \lambda_s}_{\L^p_c(\E)} 
\leq \epsi \leq \norm{f}_{\L^p_{cr}(\E)} \lesssim \norm{f}_{\L^p(\Gamma_q(H) \rtimes_\alpha G)}.
\end{align*}
We turn to the other two parts. With \eqref{equ-Fourier-mult-on-column-space}, we obtain
\begin{align*}
\MoveEqLeft
\norm{(\Id_{\Gamma_q(H)} \ot M_{\varphi_j})(g)}_{\L^p_c(\E)} + \norm{(\Id_{\Gamma_q(H)} \ot M_{\varphi_j})(h)}_{\L^p_r(\E)} \lesssim_p \norm{g}_{\L^p_c(\E)} + \norm{h}_{\L^p_r(\E)} \\
& \lesssim \norm{f}_{\L^p_{cr}(\E)} \lesssim_p \norm{f}_{\L^p(\Gamma_q(H)\rtimes_\alpha G)} .
\end{align*}
Together we have shown that also the infimum in \eqref{Khintchine-group-p<2} restricted to elements of the vector space $\Span \left\{ s_q(e) : \: e \in H \right\} \rtimes \P_G$ is also controlled by $\norm{f}_{\Gauss_{q,p,\rtimes}(\L^p(\VN(G)))}$.

Now, we will prove the remaining estimate, that is, $\norm{f}_{\Gauss_{q,p,\rtimes}(\L^p(\VN(G)))}$ is controlled by the second expression in \eqref{Khintchine-group-p<2}. Since $1<p<2$, the function $\R^+ \to \R^+$, $t \mapsto t^{\frac{p}{2}}$ is operator concave by \cite[p.~112]{Bha1}. Using \cite[Corollary 2.2]{HaP1} applied with the trace preserving positive map $\E$, we can write
\begin{align*}
\MoveEqLeft
\norm{f}_{\L^p(\Gamma_q(H) \rtimes_{\alpha} G)}
=\norm{|f|^2}_{\L^{\frac{p}{2}}(\Gamma_q(H) \rtimes_{\alpha} G)}^{\frac{1}{2}}
\leq \norm{\E(|f|^2)}_{\L^{\frac{p}{2}}(\Gamma_q(H) \rtimes_{\alpha} G)}^{\frac{1}{2}}
= \norm{f}_{\L^p_c(\E)}
\end{align*}
and similarly for the row term.
Thus, passing to the infimum over all decompositions $f = g+h$, we obtain $\norm{f}_{\L^p(\Gamma_q(H) \rtimes_\alpha G)} \leq \norm{f}_{\L^p_{cr}(\E)}$, which can be majorised in turn by the infimum of $\norm{g}_{\L^p_c(\E)} + \norm{h}_{\L^p_r(\E)}$, where $f = g + h$ and $g,h \in \Span\{s_q(h) : \: h \in H \} \rtimes \P_G$.
Hence, we have the last equivalence in the part 1 of the theorem. 


The case $p=2$ is obvious since we have isometrically $\L^2_{cr}(\E)=\L^2(\VN(G),\L^2(\Gamma_q(H))_{\rad,2})=\L^2(\Gamma_q(H)) \ot_2 \L^2(\VN(G))$ by Lemma \ref{Lema-egalite-fantas-group} and \cite[Remark 2.3 (1)]{JMX}.
\end{proof}

The remainder of the subsection is devoted to extend Theorem \ref{Khintchine-group-twisted} to the case of $f$ being a generic element of $\Gauss_{q,p,\rtimes}(\L^p(\VN(G)))$.
First we have the following lemma.

\begin{lemma}
\label{lem-Qp-bounded}
Let $-1 \leq q \leq 1$ and $1 < p < \infty$.
Consider again $Q_p = P \rtimes \Id_{\L^p(\VN(G))}$ the extension of the Gaussian projection from Lemma \ref{lem-repair-row}.
Then $Q_p$ extends to a bounded operator
\begin{equation}
\label{equ-4-repair-row}
Q_p \co \L^p(\Gamma_q(H) \rtimes_\alpha G) \to \L^p(\Gamma_q(H) \rtimes_\alpha G).
\end{equation}
\end{lemma}

\begin{proof}
First note that the case $G = \{e\}$ is contained in \cite[Theorem 3.5]{JuL1}, putting there $d = 1$.
Note that the closed space spanned by $\{ s_q(h) : h \in H \}$ coincides in this case with $\mathcal{G}^1_{p,q}$ there.
To see that the projection considered in this source is $Q_p$, we refer to \cite[Proof of Theorem 3.1]{JuL1}.

We turn to the case of general discrete group $G$.
Since $P$ is selfadjoint, $(Q_p)^* = Q_{p^*}$.
We obtain for $f = \sum_{s \in F} f_s \rtimes \lambda_s$ with $F \subseteq G$ finite and $f_s \in \L^p(\Gamma_q(H))$
\begin{align*}
\MoveEqLeft
\norm{Q_p(f) }_{\Gauss_{q,p,\rtimes}(\L^p(\VN(G)))} \ov{\text{Theorem }\ref{Khintchine-group-twisted}}{\lesssim} \norm{ Q_p(f) }_{\L^p_{cr}(\E)} 
= \sup \left\{ \left|\langle Q_p(f), g \rangle \right| : \: \norm{g}_{\L^{p^*}_{cr}(\E)} \leq 1 \right\} \\
& = \sup \left\{ \left| \langle f, Q_{p^*}(g) \rangle_{\L^p_{cr}(\E),\L^{p^*}_{cr}(\E)} \right| : \: g \right\} 
 = \sup \left\{ \left| \langle f, Q_{p^*}(g) \rangle_{\L^p(\Gamma_q(H) \rtimes_\alpha G),\L^{p^*}} \right| : \: g \right\}   \\
& \leq \norm{f}_{\L^p(\Gamma_q(H) \rtimes_\alpha G)} \sup \left\{ \norm{ Q_{p^*}(g) }_{\L^{p^*}(\Gamma_q(H) \rtimes_\alpha G)} : \: g \right\} \\
& \ov{\text{Theorem }\ref{Khintchine-group-twisted}}{\lesssim} \norm{f}_{\L^p(\Gamma_q(H) \rtimes_\alpha G)} \sup \left\{ \norm{ Q_{p^*}(g) }_{\L^{p^*}_{cr}(\E)} : \: \norm{g}_{\L^{p^*}_{cr}(\E)} \leq 1 \right\} \\
& \lesssim \norm{f}_{\L^p(\Gamma_q(H) \rtimes_\alpha G)} ,
\end{align*}
where in the last step we used that $Q_{p^*}$ is bounded on $\L^{p^*}_{cr}(\E)$ according to Lemma \ref{lem-repair-row}.
\end{proof}

Now we obtain the following extension of Theorem \ref{Khintchine-group-twisted}.

\begin{prop}
\label{Khintchine-group-twisted-full-space}
Let $-1 \leq q \leq 1$ and $1 < p < \infty$. Assume that $\L^p(\VN(G))$ has $\CCAP$ and that $\Gamma_q(H) \rtimes_\alpha G$ has $\QWEP$. Let $f$ be an element of $\Ran Q_p = \Gauss_{q,p,\rtimes}(\L^p(\VN(G)))$\footnote{\thefootnote. This equality follows from the fact that $\P_{\rtimes,G}$ is dense in $\L^p(\Gamma_q(H) \rtimes_\alpha G)$ according to Lemma \ref{lem-PG-dense}, so $Q_p(\P_{\rtimes,G}) \subseteq \Gauss_{q,p,\rtimes}(\L^p(\VN(G)))$ is dense in $\Ran Q_p$. The other inclusion $\Gauss_{q,p,\rtimes}(\L^p(\VN(G))) \subseteq \Ran Q_p$ follows from the fact that the span of the $s_q(h) \rtimes x$ with $x \in \L^p(\VN(G))$ is dense in $\Gauss_{q,p,\rtimes}(\L^p(\VN(G)))$ and obviously lies in $\Ran Q_p$ which in turn is closed.}.
\begin{enumerate}
\item
Suppose $1 < p < 2$.
Then we have
\begin{align}
\MoveEqLeft
\label{Khintchine-group-p<2-full-space}
\norm{f}_{\Gauss_{q,p,\rtimes}(\L^p(\VN(G)))}  
\approx_p \inf_{f=g+h} \bigg\{\norm{\big(\E(g^*g)\big)^{\frac12}}_{\L^p(\VN(G))},\norm{\big(\E(h h^*) \big)^{\frac12}}_{\L^p(\VN(G))}\bigg\}.
\end{align}
Here the infimum runs over all $g \in \L^p_c(\E)$ and $h \in \L^p_r(\E)$ such that $f = g + h$.
Further one can restrict the infimum to all $g,h \in \Ran Q_p$.
\item Suppose $2 \leq p <\infty$. Then we have
	\begin{align}
\MoveEqLeft
 	\label{Khintchine-group-p>2-Thm332-full-space}
\max\bigg\{\bnorm{\big(\E(f^* f)\big)^{\frac12}}_{\L^p(\VN(G))},\norm{\big(\E(f f^*) \big)^{\frac12}}_{\L^p(\VN(G))} \bigg\}
\leq \norm{f}_{\Gauss_{q,p,\rtimes}(\L^p(\VN(G)))}    \\
		&\lesssim \sqrt{p} \max\bigg\{\norm{\big(\E(f^* f)\big)^{\frac12}}_{\L^p(\VN(G))},\norm{\big( \E(f f^*) \big)^{\frac12}}_{\L^p(\VN(G))} \bigg\}. \nonumber
\end{align}
\end{enumerate}
\end{prop}

\begin{proof}
We let $(M_{\varphi_j})_j$ be the approximating net of finitely supported Fourier multipliers guaranteed by $\CCAP$ assumption. Observe that $\Id_{\Gamma_q(H)} \rtimes M_{\varphi_j}$ approximates the identity on $\P_{\rtimes,G}$. Moreover, according to Proposition \ref{prop-Fourier-mult-crossed-product}, $(\Id_{\Gamma_q(H)} \rtimes M_{\varphi_j})_j$ is a bounded net in $\B(\L^p(\Gamma_q(H)\rtimes_\alpha G))$. We conclude by density of $\P_{\rtimes,G}$ in $\L^p(\Gamma_q(H) \rtimes_\alpha G)$ that $(\Id_{\Gamma_q(H)} \rtimes M_{\varphi_j})_j$ converges in the point norm topology of $\L^p(\Gamma_q(H) \rtimes_\alpha G)$ to the identity. Moreover, replacing in this argument Proposition \ref{prop-Fourier-mult-crossed-product} by Lemma \ref{lem-Fourier-mult-on-row-space}, the same argument yields that $(\Id_{\Gamma_q(H)} \rtimes M_{\varphi_j})_j$ converges to the identity in the point norm topology of $\L^p_c(\E)$ and of $\L^p_r(\E)$. Then again by the same argument, we infer that $\norm{f}_{\L^p_{cr}(\E)} = \lim_j \norm{(\Id_{\Gamma_q(H)} \rtimes M_{\varphi_j})(f)}_{\L^p_{cr}(\E)}$.
Note that for fixed $j$, $(\Id_{\Gamma_q(H)} \rtimes M_{\varphi_j})(f) = (\Id_{\Gamma_q(H)} \rtimes M_{\varphi_j})(Q_pf) = Q_p(\Id_{\Gamma_q(H)} \rtimes M_{\varphi_j})(f)$ belongs to $\Span \{ s_q(e):\:e \in H \} \rtimes \P_G$. Thus, Theorem \ref{Khintchine-group-twisted} applies to $f$ replaced by $\Id_{\Gamma_q(H)} \rtimes M_{\varphi_j}(f)$, and
therefore,
\begin{align*}
\MoveEqLeft
\norm{f}_{\Gauss_{q,p,\rtimes}(\L^p(\VN(G)))} = \lim_j \norm{(\Id_{\Gamma_q(H)} \rtimes M_{\varphi_j})(f)}_{\Gauss_{q,p,\rtimes}(\L^p(\VN(G)))} \\
& \ov{\text{Theorem }\ref{Khintchine-group-twisted}}{\cong} \lim_j \norm{(\Id_{\Gamma_q(H)} \rtimes M_{\varphi_J})(f)}_{\L^p_{cr}(\E)} 
= \norm{f}_{\L^p_{cr}(\E)}.
\end{align*}
Finally, the fact that one can restrict the infimum to all $g,h \in \Ran Q_p$ can be proved in the same way as that in Theorem \ref{Khintchine-group-twisted}.
\end{proof}

\subsection{Kato's square root problem for semigroups of Fourier multipliers}
\label{sec-Fourier-Kato}

Throughout this subsection, we consider a discrete group $G$, and fix a markovian semigroup of Fourier multipliers $(T_t)_{t \geq 0}$ acting on $\VN(G)$ with negative generator $A$ and representing objects $b_\psi \co G \to H,\: \pi \co G \to \mathrm{O}(H),\: \alpha \co G \to \Aut(\Gamma_q(H))$, see Proposition \ref{prop-Schoenberg} and \eqref{equ-markovian-automorphism-group}. We also have the noncommutative gradient $\partial_{\psi,q} \co \P_G \subseteq \L^p(\VN(G)) \to \L^p(\Gamma_q(H) \rtimes_\alpha G)$ from \eqref{def-partial-psi}. The aim of this subsection is to compare $A^{\frac12}(x)$ and $\partial_{\psi,q}(x)$ in the $\L^p$-norm. We shall also extend $\partial_{\psi,q}$ to a closed operator and identify the domain of its closure. To this end, we need some facts from the general theory of $C_0$-semigroups.
The following is a straightforward extension of \cite[Lemma 4.2]{Bak1}.

\begin{prop}
\label{prop-Bakry-estimate-Markov}
Let $(T_t)_{t \geq 0}$ be a strongly continuous bounded semigroup acting on a Banach space $X$ with (negative) generator $A$. We have
\begin{equation}
\label{Bakry-estimate-Markov}
\bnorm{\big(\Id_{X}+A\big)^{\frac{1}{2}}(x)}_{X}	
\approx \norm{x}_{X}+\bnorm{A^{\frac{1}{2}}(x)}_{X}, \quad x \in X.
\end{equation}
\end{prop}

\begin{proof}
It is well-known \cite[Lemma 2.3]{Bak1} that the function $f_1 \co t \mapsto (1+t)^{\frac{1}{2}}(1+t^{\frac{1}{2}})^{-1}$ is the Laplace transform $\mathcal{L}(\mu_1)$ of some bounded measure $\mu_1$. By \cite[Proposition 3.3.2]{Haa1} (note that $f_1 \in \HI(\Sigma_{\pi - \epsi})$ for any $\epsi \in (0,\frac{\pi}{2})$ and $f_1$ has finite limits $\lim_{z \in \Sigma_{\pi - \epsi},\: z \to 0} f_1(z)$ and $\lim_{z \in \Sigma_{\pi - \epsi}, \: |z| \to \infty} f_1(z)$) (see also \cite[Lemma 2.12]{LM1}), we have
\begin{align*}
\MoveEqLeft
 \int_{0}^{\infty} T_t \d \mu_1(t) 
=\mathcal{L}(\mu_1)(A)
=f_1(A)
=(\Id_X+A)^{\frac{1}{2}}\Big(\Id_X+A^{\frac{1}{2}}\Big)^{-1}.        
\end{align*}
Hence for any $x \in \dom(\Id_X + A^{\frac12}) = \dom A^{\frac12}$, we have
\begin{equation}
\label{formula-Id+Ap-1-2}
(\Id_X+A)^{\frac{1}{2}}x
=\int_{0}^{\infty} T_t\big(x+A^{\frac{1}{2}}x\big) \d \mu_1(t).	
\end{equation}
By \cite[Remark 3.3.3]{Haa1}, we know that $\int_{0}^{\infty} T_t \d\mu_1(t)$ is a bounded operator on $X$ of norm $\leq \norm{\mu_1}$. Using the triangular inequality in the last inequality, we conclude that
\begin{align*}
\MoveEqLeft
\bnorm{(\Id_X+A)^{\frac{1}{2}}x}_{X}
\ov{\eqref{formula-Id+Ap-1-2}}{=}\norm{\int_{0}^{\infty} T_t\big(x+A^{\frac{1}{2}}x\big) \d \mu_1(t)}_{X}
\leq \norm{\mu_1} \bnorm{x + A^{\frac{1}{2}}x}_{X} \\
&\leq \norm{\mu_1}\Big(\norm{x}_{X} + \bnorm{A^{\frac{1}{2}}x}_{X}\Big).
\end{align*}
It is easy to see \cite[Lemma 2.3]{Bak1} that the functions $f_2 \co t \mapsto (1+t^{\frac{1}{2}})(1+t)^{-\frac{1}{2}}$ and $f_3 \co t \mapsto (1+t)^{-\frac{1}{2}}$ are the Laplace transforms $\mathcal{L}(\mu_2)$ and $\mathcal{L}(\mu_3)$ of some bounded measures $\mu_2$ and $\mu_3$. So we have
\begin{equation*}
\label{}
\int_{0}^{\infty} T_t \d \mu_2(t)
=\Big(\Id_X+A^{\frac{1}{2}}\Big)(\Id+A)^{-\frac{1}{2}}
\quad \text{and} \quad
\int_{0}^{\infty} T_t \d \mu_3(t)
=(\Id_X+A)^{-\frac{1}{2}}.
\end{equation*}
Following the same argument as above, we obtain
\begin{equation}
\label{estimate-mu-2}
\bnorm{x + A^{\frac{1}{2}}x}_{X}
\leq \norm{\mu_2} \bnorm{(\Id_X+A)^{\frac{1}{2}}x}_{X}	
\text{ and }
\norm{x}_{X}
\leq \norm{\mu_3} \bnorm{(\Id_X+A)^{\frac{1}{2}}x}_{X}.
\end{equation}
Note that
\begin{equation}
\label{Estimate-divers-Ap-3}
\bnorm{A^{\frac{1}{2}}x}_{X} 
= \bnorm{-x+x+A^{\frac{1}{2}}x}_{X}
\leq \bnorm{x}_{X}+\bnorm{x+A^{\frac{1}{2}}x}_{X}.	
\end{equation}
We conclude that
$$
\norm{x}_{X} + \bnorm{A^{\frac{1}{2}}x}_{X} 
\ov{\eqref{Estimate-divers-Ap-3}}{\leq} 2\norm{x}_{X} + \bnorm{x +A^{\frac{1}{2}}x}_{X} 
\ov{\eqref{estimate-mu-2}}{\lesssim} \bnorm{(\Id_X+A)^{\frac{1}{2}}x}_{X}.
$$
\end{proof}


We start to observe that the estimates in \eqref{Estimations-Riesz-group-q=1} below come with a constant independent of the group $G$ and the cocycle $(\alpha,H)$. This is essentially the second proof of the appendix \cite[pp.~574-575]{JMP2}. Note that we are unfortunately unable\footnote{\thefootnote. More precisely, with the notations of \cite{JMP2} we are unable to check that ``$H \rtimes \Id_G$ extends to a bounded operator on $L_p$''. A part of the very concise explanation given in \cite[p.~544]{JMP2} is ``$H$ is $\mathrm{G}$-equivariant''. But it seems to be strange. Indeed we have an action $\alpha \co G \to \Aut(L_\infty(\R^n_{\mathrm{bohr}}))$, $f \mapsto \big[x \mapsto \alpha_g(f)(x) = f(\pi_g(x))\big]$ for some map $\pi_g \co \R^n_{\mathrm{bohr}} \to \R^n_{\mathrm{bohr}}$ where $g \in G$ and an induced action $\alpha$ from $G$ on $L_\infty(\R^n_{\mathrm{bohr}} \times \R^n,\nu \times \gamma)$. Now, note that
$$
(H(\alpha_gf))(x,y)
=\bigg(\pv \int_\R \beta_t \alpha_g f \frac{\d t}{t}\bigg)(x,y) 
=\bigg(\pv \int_\R \beta_t (f \circ \pi_g) \frac{\d t}{t} \bigg)(x,y)
=\pv \int_\R f(\pi_g(x+ty)) \frac{\d t}{t}
$$
and
$$
(\alpha_g(H(f))(x,y)
=\alpha_g\bigg(\pv \int_\R \beta_t f \frac{\d t}{t}\bigg)(x,y)
=\pv \int_\R f(\pi_g(x) + ty) \frac{\d t}{t}
$$
which could be different if $\pi$ is not trivial.} to check the original proof given in \cite[p.~ 544]{JMP2}. Here, $\L^\infty(\Omega) = \Gamma_1(H)$ is the Gaussian space from Subsection \ref{Sec-q-gaussians}.

\begin{lemma}
\label{lem-Riesz-Fourier-dimension-free-q=1}
Suppose $1 < p < \infty$. For any $x \in \dom \P_G$, we have
\begin{equation}
\label{Estimations-Riesz-group-q=1}
\frac{1}{K  \max(p,p^*)} \bnorm{A_p^{\frac12}(x)}_{\L^p(\VN(G))} 
\leq \bnorm{\partial_{\psi,1,p}(x)}_{\L^p(\L^\infty(\Omega) \rtimes_\alpha G)}
\leq K  \max(p,p^*)^{\frac{3}{2}} \bnorm{A_p^{\frac12}(x)}_{\L^p(\VN(G))}
\end{equation}
with an absolute constant $K$ not depending on $G$ nor the cocycle $(\alpha,H)$.
\end{lemma}

We define the densely defined unbounded operator $\partial_{\psi,q}^* \co \P_{\rtimes,G} \subseteq \L^p(\Gamma_q(H) \rtimes_\alpha G) \to \L^p(\Gamma_q(H))$ by
\begin{equation}
\label{Adjoint-partial-Fourier}
\partial_{\psi,q}^*(f \rtimes \lambda_s)
=\big\langle  s_q(b_\psi(s)),f\big\rangle_{\L^{p^*}(\Gamma_{q}(H)),\L^{p}(\Gamma_{q}(H))} \lambda_s, \quad s \in G, f \in \L^p(\Gamma_q(H)).
\end{equation}

The following lemma is left to the reader.

\begin{lemma}
\label{Lemma-adjoint-prtial-Fourier}
The operators $\partial_{\psi,q}$ and $\partial_{\psi,q}^*$ are formal adjoints.
\end{lemma}

The next proposition extends \eqref{Estimations-Riesz-group-q=1} to the case of $q$-Gaussians.

\begin{prop}
\label{prop-Riesz-Fourier-dimension-free}
Suppose $1 < p < \infty$ and $-1 \leq q \leq 1$. For any $x \in \dom P_G$, we have
\begin{equation}
\label{Estimations-Riesz-group-q}
\frac{1}{K  \max(p,p^*)^{\frac{3}{2}}} \bnorm{A_p^{\frac12}(x)}_{\L^p(\VN(G))} 
\leq \bnorm{\partial_{\psi,q,p}(x)}_{\L^p(\Gamma_q(H) \rtimes_\alpha G)} 
\leq K  \max(p,p^*)^{2} \bnorm{A_p^{\frac12}(x)}_{\L^p(\VN(G))} 
\end{equation}
with an absolute constant $K$ not depending on $G$ nor the cocycle $(\alpha,b_\psi,H)$.
\end{prop}

\begin{proof}
Start with the case $2 \leq p < \infty$. The case $q = 1$ is covered by Lemma \ref{lem-Riesz-Fourier-dimension-free-q=1}. Consider now the case $-1 \leq q < 1$. Pick some element $x = \sum_s x_s \lambda_s$ of $\P_G$. We recall from Subsection \ref{twisted-Khintchine} that we have a conditional expectation $\E  \co \Gamma_q(H) \rtimes_\alpha G \to \Gamma_q(H) \rtimes_\alpha G,\: x \rtimes \lambda_s \mapsto \tau_{\Gamma_q(H)}(x) \lambda_s$. In the following calculations, we consider be an orthonormal basis $(e_k)$ of $H$. Therefore, using the orthonormal systems $(\W(e_k))$ and $(s_q(e_k))$ in $\L^2(\Omega)$ and $\L^2(\Gamma_q(H))$ in the third equality, we obtain
\begin{align*}
\MoveEqLeft
\norm{ \sum_{s} x_s \W(b_\psi(s)) \rtimes \lambda_s }_{\L^p_c(\E)} 
= \norm{  \sum_{s,k} x_s \langle e_k, b_\psi(s) \rangle \W(e_k) \rtimes \lambda_s}_{\L^p_c(\E)} \\
&\ov{\eqref{Egalite-fantastique-group}}{=} \norm{ \sum_{s,k} x_s \langle e_k, b_\psi(s) \rangle \lambda_s \ot \W(e_k) }_{\L^p(\VN(G),\L^2(\Omega)_{c,p})} \\ 
&= \norm{ \sum_{s,k} x_s \langle e_k, b_\psi(s) \rangle \lambda_s \ot s_q(e_k) }_{\L^p(\VN(G),\L^2(\Gamma_q(H))_{c,p})} 
\ov{\eqref{Egalite-fantastique-group}}{=} \norm{ \sum_s x_s s_q(b_\psi(s)) \rtimes \lambda_s }_{\L^p_c(\E)}.
\end{align*}
We claim that this equality holds also for the row space in place of column space. Note that according to Remark \ref{rem-problem-row}, we have to argue differently. We have with $x_{s,k} \ov{\mathrm{def}}{=} x_s \langle e_k , b_\psi(s) \rangle_H$,
\begin{align*}
\MoveEqLeft
\norm{ \sum_s x_s s_q(b_\psi(s)) \rtimes \lambda_s }_{\L^p_r(\E)}^2 
= \norm{ \sum_{s,k,t,l} \sum_{t,l}x_{s,k} \ovl{x_{t,l}} \E\big((s_q(e_k) \rtimes \lambda_s) (s_q(e_l) \rtimes \lambda_t)^* \big) }_{\frac{p}{2}} \\
& \ov{\eqref{inverse-crossed-product}}{=} \norm{ \sum_{s,k,t,l} x_{s,k} \ovl{x_{t,l}} \E \big( (s_q(e_k) \rtimes \lambda_s)( \alpha_{t^{-1}}(s_q(e_l)) \rtimes \lambda_{t^{-1}} )\big) }_{\frac{p}{2}} \\
& \ov{\eqref{Product-crossed-product}}{=} \norm{ \sum_{s,k,t,l} x_{s,k} \ovl{x_{t,l}} \E \big(   s_q(e_k) \alpha_{st^{-1}}(s_q(e_l)) \rtimes \lambda_{st^{-1}} \big) }_{\frac{p}{2}} \\
& = \norm{ \sum_{s,k,t,l} x_{s,k} \ovl{x_{t,l}} \tau_q\big(s_q(e_k) s_q(\pi_{st^{-1}}(e_l))\big) \lambda_{st^{-1}} }_{\frac{p}{2}} 
\ov{\eqref{petit-Wick}}{=} \norm{ \sum_{s,k,t,l} x_{s,k} \ovl{x_{t,l}} \big\langle e_k, \pi_{st^{-1}}(e_l) \big\rangle_H \lambda_{st^{-1}} }_{\frac{p}{2}}.
\end{align*}
The point is that this last quantity does not depend on $q$. We infer that
\begin{equation}
\label{equ-1-Estimation-Riesz-group-q}
\norm{\sum_s x_s \W(b_\psi(s)) \rtimes \lambda_s }_{\L^p_{cr}(\E)} = \norm{ \sum_s x_s s_q(b_\psi(s)) \rtimes \lambda_s}_{\L^p_{cr}(\E)} .
\end{equation}
Then we have
\begin{align*}
\MoveEqLeft
\bnorm{A_p^{\frac12}(x)}_{\L^p(\VN(G))} \ov{\eqref{Estimations-Riesz-group-q=1}}{\leq} K p \bnorm{\partial_{\psi,1,p}(x)}_{\L^p(\L^\infty(\Omega) \rtimes_\alpha G)} \\
&= Kp \norm{\sum_s x_s \W(b_\psi(s)) \rtimes \lambda_s}_p = Kp \norm{\sum_{s,k}  x_s \langle e_k,b_\psi(s) \rangle \W(e_k) \rtimes \lambda_s}_p \\
& \ov{\eqref{Khintchine-group-p>2-Thm332}}{\leq} K' p \cdot p^{\frac12} \norm{\sum_{s,k} x_s \langle e_k, b_\psi(s) \rangle \W(e_k) \rtimes \lambda_s}_{\L^p_{cr}(\E)} \\
&\ov{\eqref{equ-1-Estimation-Riesz-group-q}}{=} K' p^{\frac32} \norm{\sum_{s,k} x_s \big\langle e_k, b_\psi(s) \big\rangle s_q(e_k) \rtimes \lambda_s}_{\L^p_{cr}(\E)} \\
& \ov{ \eqref{Khintchine-group-p>2-Thm332}}{\leq} K' p^{\frac32} \norm{\sum_{s,k} x_s \big\langle e_k,b_\psi(s)\big\rangle s_q(e_k) \rtimes \lambda_s}_{\L^p(\Gamma_q(H) \rtimes_\alpha G)} \\
&= K' p^{\frac32} \bnorm{\partial_{\psi,q,p}(x)}_{\L^p(\Gamma_q(H) \rtimes_\alpha G)}.
\end{align*}

We pass to the converse inequality. We have
\begin{align*}
\MoveEqLeft
\bnorm{A_p^{\frac12}(x)}_{\L^p(\VN(G))} 
\overset{\eqref{Estimations-Riesz-group-q}}{\geq} \frac{1}{K p^{\frac32}} \bnorm{\partial_{\psi,1,p}(x)}_{\L^p(\L^\infty(\Omega) \rtimes_\alpha G)} \\
& = \frac{1}{Kp^{\frac32}} \norm{\sum_s x_s \W(b_\psi(s)) \rtimes \lambda_s}_p = \frac{1}{Kp^{\frac32}} \norm{\sum_{s,k} x_s \langle e_k, b_\psi(s) \rangle \W(e_k) \rtimes \lambda_s}_p \\
& \overset{\eqref{Khintchine-group-p>2-Thm332}}{\geq} \frac{1}{K' p^{\frac32}} \norm{\sum_{s,k} x_s \langle e_k, b_\psi(s) \rangle \W(e_k) \rtimes \lambda_s}_{\L^p_{cr}(E)} \\
&\ov{\eqref{equ-1-Estimation-Riesz-group-q}}{ = } \frac{1}{K' p^{\frac32}} \norm{\sum_{s,k} x_s \langle e_k, b_\psi(s) \rangle s_q(e_k) \rtimes \lambda_s}_{\L^p_{cr}(E)} \\
& \overset{\eqref{Khintchine-group-p>2-Thm332}}{\geq} \frac{1}{K''p^{\frac32} \cdot p^{\frac12}} \norm{\sum_{s,k} x_s \langle e_k,b_\psi(s)\rangle s_q(e_k) \rtimes \lambda_s}_{\L^p(\Gamma_q(H) \rtimes_\alpha G)} \\
&= \frac{1}{K'' p^2} \norm{\partial_{\psi,q,p}(x)}_{\L^p(\Gamma_q(H) \rtimes_\alpha G)}.
\end{align*}
Altogether we have shown \eqref{Estimations-Riesz-group-q} in the case $2 \leq p < \infty$.

We turn to the case $1 < p < 2$.
Note that \eqref{equ-1-Estimation-Riesz-group-q} still holds in this case.
Indeed, for elements $f \in \Gauss_{q,p,\rtimes}(\L^p(\VN(G)))$, the norm $\norm{f}_{\L^p_{cr}(\E)} = \inf\{\norm{g}_{\L^p_c(\E)} + \norm{h}_{\L^p_r(\E)} :\: f = g + h\}$ remains unchanged if $g,h$ are choosen in $\Ran Q_p|_{\L^p_c(\E)}$ and $\Ran Q_p|_{\L^p_r(\E)}$, see Theorem \ref{Khintchine-group-twisted}.
But for those elements $g,h$, the $\L^p_c(\E)$ and $\L^p_r(\E)$ norms remain unchanged upon replacing classical Gaussian variables $\W(e_k)$ by $q$-Gaussian variables $s_q(e_k)$, see the beginning of the proof.
Then the proof in the case $1 < p < 2$ can be executed as in the case $2 \leq p < \infty$, noting that the additional factor ${p^{*}}^{\frac12}$ will appear at another step of the estimate, in accordance with Theorem \ref{Khintchine-group-twisted} parts 1. and 2.
\end{proof}


\begin{prop}
\label{Prop-derivation-closable-sgrp}
\label{prop-fermable-sgrp}
Let $G$ be a discrete group. Suppose $1 < p<\infty$ and $-1 \leq q \leq 1$.
\begin{enumerate}
	\item The operator $\partial_{\psi,q} \co \P_G \subseteq \L^p(\VN(G)) \to \L^p(\Gamma_q(H) \rtimes_{\alpha} G)$ is closable as a densely defined operator on $\L^p(\VN(G))$ into $\L^p(\Gamma_q(H) \rtimes_{\alpha} G)$. We denote by $\partial_{\psi,q,p}$ its closure. So $\P_G$ is a core of $\partial_{\psi,q,p}$.
	
	\item $(\Id_{\L^p(\VN(G))} +  A_p)^{\frac{1}{2}}(\P_G)$ is a dense subspace of $\L^p(\VN(G))$.
	
	\item $\P_G$ is a core of $A_p^{\frac{1}{2}}$.
	
	\item We have $\dom \partial_{\psi,q,p}=\dom A_p^{\frac{1}{2}}$. Moreover, for any $x \in \dom A_p^{\frac{1}{2}}$, we have
\begin{equation} 
\label{Equivalence-square-root-domaine-Schur-sgrp}
\bnorm{A_p^{\frac12}(x)}_{\L^p(\VN(G))} 
\approx_{p} \bnorm{\partial_{\psi,q,p}(x)}_{\L^p(\Gamma_q(H) \rtimes_{\alpha} G)}. 
\end{equation}	
Finally, for any $x \in \dom A_p^{\frac{1}{2}}$, there exists a sequence $(x_n)$ of elements of $\P_G$ such that $x_n \to x$, $A_p^{\frac{1}{2}}(x_n) \to A_p^{\frac{1}{2}}(x)$ and $\partial_{\psi,q,p}(x_n) \to \partial_{\psi,q,p}(x)$.

\item If $x \in \dom \partial_{\psi,q,p}$, we have $x^* \in \dom \partial_{\psi,q,p}$ and 
\begin{equation}
\label{relation-partial-star-group}
(\partial_{\psi,q,p}(x))^* 
=-\partial_{\psi,q,p}(x^*).
\end{equation}
\item Let $M_\psi \co \L^p(\VN(G)) \to \L^p(\VN(G))$ be a finitely supported bounded Fourier multiplier such that the map $\Id \rtimes M_\psi \co \L^p(\Gamma_q(H) \rtimes_{\alpha} G) \to \L^p(\Gamma_q(H) \rtimes_{\alpha} G)$ is a well-defined bounded operator. For any $x \in \dom \partial_{\psi,q,p}$, the element $M_{\psi}(x)$ belongs to $\dom \partial_{\psi,q,p}$ and we have 
\begin{equation}
\label{commute-troncature}
\partial_{\psi,q,p} M_{\psi}(x)
=(\Id \rtimes M_{\psi}) \partial_{\psi,q,p}(x) .
\end{equation}
\end{enumerate} 
\end{prop}

\begin{proof}
1. This is a consequence of \cite[Theorem 5.28 p.~168]{Kat1} together with Lemma \ref{Lemma-adjoint-prtial-Fourier} above.

2. We consider the semigroup $(R_t)_{t \geq 0}$ associated to the negative generator $(\Id_{\L^p(\VN(G))} +  A_p)^{\frac{1}{2}}$. Let $x \in \L^{p^*}(\VN(G))$ orthogonal to $(\Id_{\L^p(\VN(G))} +  A_p)^{\frac{1}{2}}(\P_G)$. For an element $y \in \P_G$, using \cite[Lemma 1.3 (iv) p.~50]{EnN1} in the second equality, we have
\begin{align*}
\MoveEqLeft
\langle R_t(x), y \rangle_{}
=\langle x, R_t(y) \rangle_{}            
=\bigg\langle x,(\Id_{\L^p(\VN(G))} +  A_p)^{\frac{1}{2}}\int_{0}^{t} R_s(y) \d s+y \bigg\rangle_{}\\
&=\bigg\langle x,(\Id_{\L^p(\VN(G))} +  A_p)^{\frac{1}{2}}\int_{0}^{t} R_s(y) \d s\bigg\rangle +\langle x,y \rangle_{}
=\langle x,y \rangle_{}.
\end{align*}  
By density, we deduce that $x=R_t(x)$. We know that $(R_t(x))$ converges to 0 when $t$ goes to $+\infty$. Indeed, using \cite[2.2 p.~60]{EnN1} in the second equality, a standard computation gives
\begin{align*}
\MoveEqLeft
\bnorm{R_t(x)}_{\L^p(\VN(G))}
\ov{\eqref{subordinated-formula}}{=}\norm{\int_{0}^{+\infty} f_t(s) \e^{-s(\Id+A_p)}(x) \d s}_{\L^p(\VN(G))}
=\norm{\int_{0}^{\infty} f_t(s)\e^{-s} T_s(x) \d s}_{\L^p(\VN(G))} \\          
&\leq \int_{0}^{\infty} f_t(s)\e^{-s} \bnorm{T_s(x)}_{\L^p(\VN(G))} \d s
\leq \bigg(\int_{0}^{\infty} \e^{-s}f_t(s) \d s\bigg)\norm{x}_{\L^p(\VN(G))} \\
&\ov{\eqref{integrale-utile-1}}{=} \e^{-t}\norm{x}_{\L^p(\VN(G))}
\xra[t \to +\infty]{} 0.
\end{align*}  
We infer that $x=0$. We conclude that 0 is the only element orthogonal to $(\Id_{\L^p(\VN(G))} +  A_p)^{\frac{1}{2}}(\P_G)$ as a subset of $\L^p(\VN(G))$. It follows that the subspace $(\Id_{\L^p(\VN(G))} +  A_p)^{\frac{1}{2}}(\P_G)$ is dense in $\L^p(\VN(G))$.

3. Let $x \in \dom A_p^{\frac{1}{2}}$. By \cite[Proposition 3.8.2]{ABHN1}, we have $\dom A_p^{\frac{1}{2}}=\dom (\Id_{\L^p(\VN(G))} + A_p )^{\frac{1}{2}}$. Hence $x \in \dom (\Id_{\L^p(\VN(G))} + A_p )^{\frac{1}{2}}$. From the point 2, there exists a sequence $(x_n)$ of $\P_G$ such that $(\Id_{\L^p(\VN(G))} +  A_p)^{\frac{1}{2}} x_n \to (\Id_{\L^p(\VN(G))} + A_p)^{\frac{1}{2}}x$ in $\L^p(\VN(G))$. We obtain that
\begin{align*}
\MoveEqLeft
\norm{x_n - x}_{\dom A_p^{\frac{1}{2}}}
\ov{\eqref{Def-graph-norm}}{=}\norm{x_n - x}_{\L^p(\VN(G))} + \bnorm{ A_p^{\frac{1}{2}} (x_n - x)}_{\L^p(\VN(G))} \\
&\ov{\eqref{Bakry-estimate-Markov}}{\lesssim_p} \bnorm{(\Id_{\L^p(\VN(G))} + A_p )^{\frac{1}{2}}(x_n - x)}_{\L^p(\VN(G))}
\xra[n \to +\infty]{} 0.            
\end{align*}
Consequently the sequence $(x_n)$ converges to $x$ in $\dom A_p^{\frac{1}{2}}$.

4. Let $x \in \dom A_p^{\frac{1}{2}}$. By the point 3, $\P_G$ is dense in $\dom A_p^{\frac{1}{2}}$ equipped with the graph norm. Hence we can find a sequence $(x_n)$ of $\P_G$ such that $x_n \to x$ and $A_p^{\frac{1}{2}}(x_n) \to A_p^{\frac{1}{2}}(x)$. For any integers $n,m$, we obtain
\begin{align*}
\MoveEqLeft
\norm{x_n-x_m}_{\L^p(\VN(G))} + \norm{\partial_{\psi,q,p}(x_n) - \partial_{\psi,q,p}(x_m)}_{\L^p(\Gamma_q(H) \rtimes_\alpha G)} \\
&\ov{\eqref{Estimations-Riesz-group-q}}{\lesssim_p} \norm{x_n-x_m}_{\L^p(\VN(G))} + \bnorm{A_p^{\frac{1}{2}}(x_n)-A_p^{\frac{1}{2}}(x_m)}_{\L^p(\VN(G))}.
\end{align*} 
which shows that $(x_n)$ is a Cauchy sequence in $\dom \partial_{\psi,q,p}$. By the closedness of $\partial_{\psi,q,p}$, we infer that this sequence converges to some $x' \in \dom \partial_{\psi,q,p}$ equipped with the graph norm. Since $\dom \partial_{\psi,q,p}$ is continuously embedded into $\L^p(\VN(G))$, we have $x_n \to x'$ in $\L^p(\VN(G))$, and therefore $x=x'$ since $x_n \to x$. It follows that $x \in \dom \partial_{\psi,q,p}$. This proves the inclusion $\dom A_p^{\frac{1}{2}} \subseteq \dom \partial_{\psi,q,p}$. Moreover, for any integer $n$, we have 
$$
\norm{\partial_{\psi,q,p}(x_n)}_{\L^p(\Gamma_q(H) \rtimes_\alpha  G)} 
\ov{\eqref{Estimations-Riesz-group-q}}{\lesssim_p}\bnorm{A_p^{\frac{1}{2}}(x_n)}_{\L^p(\VN(G))}.
$$
Since $x_n \to x$ in $\dom\partial_{\psi,q,p}$ and in $\dom A_p^{\frac{1}{2}}$ both equipped with the graph norm, we conclude that
$$
\bnorm{\partial_{\psi,q,p}(x)}_{\L^p(\Gamma_q(H) \rtimes_\alpha  G)} 
\lesssim_p\bnorm{A_p^{\frac{1}{2}}(x)}_{\L^p(\VN(G))}.
$$
The proof of the reverse inclusion and of the reverse estimate are similar. Indeed, by part 1, $\P_G$ is a dense subspace of $\dom\partial_{\psi,q,p}$ equipped with the graph norm.

5. Recall that $b_\psi(e)=0$. For any $s \in G$, we have
\begin{align*}
\MoveEqLeft
\big(\partial_{\psi,q,p}(\lambda_s)\big)^*
\ov{\eqref{def-partial-psi}}{=} \big(s_q(b_\psi(s)) \rtimes \lambda_s\big)^*
\ov{\eqref{inverse-crossed-product}}{=} \alpha_{s^{-1}}(s_q(b_\psi(s))) \rtimes \lambda_{s^{-1}} 
\ov{\eqref{Cocycle-law}}{=} s_q(b_\psi(e)-b_\psi(s^{-1})) \rtimes \lambda_{s^{-1}} \\
&=-s_q(b_\psi(s^{-1})) \rtimes \lambda_{s^{-1}} 
\ov{\eqref{def-partial-psi}}{=} -\partial_{\psi,q,p}(\lambda_{s^{-1}}) 
=-\partial_{\psi,q,p}(\lambda_s^*).            
\end{align*}
Let $x \in \dom \partial_{\psi,q,p}$. By the point 1, $\P_G$ is core of $\partial_{\psi,q,p}$. Hence there exists a sequence $(x_n)$ of $\P_G$ such that $x_n \to x$ and $\partial_{\psi,q,p}(x_n) \to \partial_{\psi,q,p}(x)$. We have $x_n^* \to x^*$ and by the first part of the proof $\partial_{\psi,q,p}(x_n^*)=-(\partial_{\psi,q,p}(x_n))^* \to  -(\partial_{\psi,q,p}(x))^*$. By \eqref{Domain-closure}, we conclude that $x^* \in \dom \partial_{\psi,q,p}$ and that $\partial_{\psi,q,p}(x^*)=-(\partial_{\psi,q,p}(x))^*$.

6. If $s \in G$, we have
\begin{align*}
\MoveEqLeft
(\Id \rtimes M_{\psi}) \partial_{\psi,q,p}(\lambda_s) 
\ov{\eqref{def-partial-psi}}{=} (\Id \rtimes M_{\psi})(s_q(b_\psi(s)) \rtimes \lambda_s)     
=\psi_s s_q(b_\psi(s)) \rtimes \lambda_s \\
&\ov{\eqref{def-partial-psi}}{=} \psi_s\partial_{\psi,q,p}(\lambda_s)
=\partial_{\psi,q,p} M_{\psi}(\lambda_s).
\end{align*} 
By linearity, \eqref{commute-troncature} is true for elements of $\P_G$. Let $x \in \dom \partial_{\psi,q,p}$. There exists a sequence $(x_n)$ of $\P_G$ such $x_n \to x$ and $\partial_{\psi,q,p}(x_n) \to \partial_{\psi,q,p}(x)$. We have $M_{\psi}(x_n) \to M_{\psi}(x)$ and
$$
\partial_{\psi,q,p} M_{\psi}(x_n)
=(\Id \rtimes M_{\psi}) \partial_{\psi,q,p}(x_n)
\xra[n \to +\infty]{}  (\Id \rtimes M_{\psi}) \partial_{\psi,q,p}(x).
$$ 
By \eqref{Def-operateur-ferme}, we deduce \eqref{commute-troncature}.
\end{proof}

For the property $\mathrm{AP}$ in the next proposition, we refer to the preliminary Subsection \ref{subsubsec-Fourier-mult-crossed-product}.

\begin{prop}
\label{Prop-closable-deriv-p-infty}
Assume $-1 \leq q < 1$.
Let $G$ be a discrete group with $\mathrm{AP}$. The operator $\partial_{\psi,q} \co \P_G \subseteq \VN(G) \to \Gamma_q(H) \rtimes_{\alpha} G$ is weak* closable\footnote{\thefootnote. That is, if $(x_n)$ is a sequence in $\P_G$ such that $x_n \to 0$ and $\partial_{\psi,q}(x_n) \to y$ for some $y \in \Gamma_q(H) \rtimes_\alpha G$ both for the weak* topology, then $y = 0$.}. We denote by $\partial_{\psi,q,\infty}$ its weak* closure. 
\end{prop}

\begin{proof}
Suppose that $(x_i)$ is a net of $\P_G$ which converges to 0 for the weak* topology with $x_i=\sum_{s \in G} x_{i,s} \lambda_s$ such that the net $(\partial_{\psi,q}(x_i))$ converges  for the weak* topology to some $y$ belonging to $\Gamma_q(H) \rtimes_{\alpha} G$. Let $(M_{\varphi_j})$ be the net of Fourier multipliers approximating the identity from Proposition \ref{Prop-weak-amenable-crossed}. For any $j$, we have 
\begin{align*}
\MoveEqLeft
\partial_{\psi,q}(M_{\varphi_j}x_i)
=\partial_{\psi,q}\bigg(M_{\varphi_j}\bigg(\sum_{s \in G} x_{i,s} \lambda_s\bigg)\bigg)
=\partial_{\psi,q}\bigg(\sum_{s \in \supp \varphi_j} x_{i,s}\lambda_s\bigg)
=\sum_{s \in \supp \varphi_j} x_{i,s}\partial_{\psi,q}(\lambda_s)\\  
&=(\Id_{} \rtimes M_{\varphi_j}) \partial_{\psi,q}\bigg(\sum_{s \in G} x_{i,s} \lambda_s\bigg)
=(\Id_{} \rtimes M_{\varphi_j}) \partial_{\psi,q} (x_i) 
\xra[i]{} (\Id_{} \rtimes M_{\varphi_j})(y).            
\end{align*}
On the other hand, for all $s \in \supp \varphi_j$ we have $x_{i,s} \to 0$ as $i \to \infty$. Hence
$$
\partial_{\psi,q}(M_{\varphi_j} x_i) 
= \partial_{\psi,q} \left( \sum_{s \in \supp(\varphi_j)} \varphi_j(s) x_{i,s} \lambda_s \right) 
\ov{\eqref{def-partial-psi}}{=} \sum_{s \in \supp(\varphi_j)} \varphi_j(s) x_{i,s} s_q(b_\psi(s)) \rtimes \lambda_s 
\xra[i]{} 0.
$$ 
This implies by uniqueness of the limit that $(\Id_{} \rtimes M_{\varphi_j})(y) = 0$ for any $j$. By the point 3 of Proposition \ref{Prop-weak-amenable-crossed}, we deduce that $y = \weakstar\lim_j (\Id_{\Gamma_q(H)} \rtimes M_{\varphi_j}) (y) = 0$.
\end{proof}

\begin{remark} \normalfont
\label{Remark-AP-necessary}
We do not know if the assumption ``AP'' in Proposition \ref{Prop-closable-deriv-p-infty} is really necessary. 
\end{remark}

The following generalizes an observation of \cite{JMP2} (in the case $q=1$).

\begin{lemma}
\label{lem-formule-trace-trace-Fourier}
Let $-1 \leq q \leq 1$. Suppose $1 < p < \infty$. For any $x \in \dom A_p^{\frac12}$ and any $y \in \dom A_{p^*}^{\frac12}$, we have
\begin{equation}
\label{Formula-trace-trace-group}
\tau_{\Gamma_q(H) \rtimes_\alpha G}\big(\partial_{\psi,q,p}(x) (\partial_{\psi,q,p^*}(y))^* \big) 
=\tau_G \Big(A_p^{\frac12}(x) \big(A_{p^*}^{\frac12}(y)\big)^* \Big).
\end{equation}
\end{lemma}

\begin{proof}
By the fact that $\P_G$ is a core for $\partial_{\psi,q,p},\partial_{\psi,q,p^*},A_p^{\frac12}$ and $A_{p^*}$ and the norm equivalence \eqref{Equivalence-square-root-domaine-Schur-sgrp}, we can assume $x,y \in \P_G$. Consider some elements $x = \sum_{s \in G} x_s \lambda_s$ and $y = \sum_{r \in G} y_r \lambda_r$ of $\P_G$ where both sums are finite. On the one hand, 
we have 
\begin{align*}
\MoveEqLeft
\tau_{\Gamma_q(H) \rtimes_\alpha G}\big(\partial_{\psi,q}(x) (\partial_{\psi,q}(y))^* \big) 
=\sum_{s,r \in G} x_s \ovl{y_r}\tau_{\Gamma_q(H) \rtimes_\alpha G}\big(\partial_{\psi,q}(\lambda_s) (\partial_{\psi,q}(\lambda_r))^* \big) \\
&\ov{\eqref{def-partial-psi}}{=} \sum_{s,r \in G} x_s \ovl{y_r}\tau_{\Gamma_q(H) \rtimes_\alpha G} \big((s_q(b_\psi(s)) \rtimes \lambda_s) (s_q(b_\psi(r)) \rtimes \lambda_r)^* \big) \\
&\ov{\eqref{inverse-crossed-product}}{=} \sum_{s,r \in G} x_s \ovl{y_r}\tau_{\Gamma_q(H) \rtimes_\alpha G} \big((s_q(b_\psi(s)) \rtimes \lambda_s) (\alpha_{r^{-1}}(s_q(b_\psi(r))) \rtimes \lambda_{r^{-1}} \big)\\
&\ov{\eqref{Product-crossed-product}}{=} \sum_{s,r \in G} x_s \ovl{y_r}\tau_{\Gamma_q(H) \rtimes_\alpha G} \big((s_q(b_\psi(s))\alpha_{sr^{-1}}(s_q(b_\psi(r))) \rtimes \lambda_{sr^{-1}} \big)\\
&=\sum_{s \in G} x_s \ovl{y_s} \tau_{\Gamma_q(H)} \big(s_q(b_\psi(s)) s_q(b_\psi(s))\big) 
\ov{\eqref{petit-Wick}}{=} \sum_{s \in G} x_s \ovl{y_s} \norm{b_\psi(s)}_H^2.
\end{align*}
On the other hand, we have
\begin{align*}
\MoveEqLeft
\tau_G\left(A^{\frac12}(x) \big(A^{\frac12}(y)\big)^* \right) 
=\sum_{s,r \in G} x_s \ovl{y_r}\tau_G\left(A^{\frac12}(\lambda_s) \big(A^{\frac12}(\lambda_r)\big)^* \right)  \\
&=\sum_{s,r \in G} x_s \ovl{y_r} \norm{b_\psi(s)}_H\norm{b_\psi(r)}_H\tau_G\big(  \lambda_s  \lambda_r^* \big) 
=\sum_{s \in G} x_s \ovl{y_s} \norm{b_\psi(s)}_H^2.
\end{align*}
\end{proof}

\subsection{Extension of the carr\'e du champ $\Gamma$ for Fourier multipliers}
\label{Sec-infos-on-Gamma-Fourier2}

In this subsection, we consider again a markovian semigroup $(T_t)_{t \geq 0}$ of Fourier multipliers acting on the von Neumann algebra $\VN(G)$ where $G$ is a discrete group. We shall extend the carr\'e du champ $\Gamma$ associated with $(T_t)_{t \geq 0}$ to a closed form and identify its domain. It will be easier to consider simultaneously $x \mapsto \Gamma(x,x)$ and $x \mapsto \Gamma(x^*,x^*)$ and the case $2 \leq p < \infty$ throughout the subsection. We will also link the carr\'e du champ with $A_p^{\frac12}$ and the gradient $\partial_{\psi,q,p}$. In most of the results of this subsection, we need approximation properties of $\L^p(\VN(G))$ and $\L^p(\Gamma_q(H) \rtimes_\alpha G)$. Note that by \cite[Theorem 1.2]{JR1}, if $G$ is a discrete group with AP such that $\VN(G)$ has $\QWEP$, then $\L^p(\VN(G))$ has the completely contractive approximation property $\mathrm{CCAP}$ for any $1<p<\infty$. With the following lemma, one can extend the definition of $\Gamma(x,y)$ to a larger domain. 

\begin{lemma}
\label{lem-Gamma-closure-group}
Suppose $2 \leq p <\infty$. Let $G$ be a discrete group. The forms $a \co \P_G \times  \P_G \to \L^{\frac{p}{2}}(\VN(G)) \oplus_\infty \L^{\frac{p}{2}}(\VN(G))$, $(x,y) \mapsto \Gamma(x,y) \oplus \Gamma(x^*,y^*)$ and $\Gamma \co \P_G \times  \P_G \to \L^{\frac{p}{2}}(\VN(G))$, $(x,y) \mapsto \Gamma(x,y)$ are symmetric, positive and satisfy the Cauchy-Schwarz inequality. The first is in addition closable and the domain of the closure $\ovl{a}$ is $\dom A^{\frac12}_p$.
\end{lemma}

\begin{proof}
According to the point 3 of Lemma \ref{Lemma-gamma-infos-sgrp}, we have $\Gamma(x,y)^* = \Gamma(y,x)$, so $a$ is symmetric. Moreover, again according to Lemma \ref{Lemma-gamma-infos-sgrp}, $a$ is positive. For any $x,y \in \P_G$, we have
\begin{align*}
\MoveEqLeft
\norm{a(x,y)}_{\L^{\frac{p}{2}}(\VN(G)) \oplus_\infty \L^{\frac{p}{2}}(\VN(G))}            
=\norm{\Gamma(x,y) \oplus \Gamma(x^*,y^*)}_{\L^{\frac{p}{2}}(\VN(G)) \oplus_\infty \L^{\frac{p}{2}}(\VN(G))} \\
&=\max\big\{\norm{\Gamma(x,y)}_{\frac{p}{2}},\norm{\Gamma(x^*,y^*)}_{\frac{p}{2}}\big\}\\
&\ov{\eqref{CS-group}}{\leq} \max\Big\{\norm{\Gamma(x,x)}_{\frac{p}{2}}^{\frac{1}{2}} \norm{\Gamma(y,y)}_{\frac{p}{2}}^{\frac{1}{2}},\norm{\Gamma(x^*,x^*)}_{\frac{p}{2}}^{\frac{1}{2}}\norm{\Gamma(y^*,y^*)}_{\frac{p}{2}}^{\frac{1}{2}}\Big\}\\
&\leq \max\big\{\norm{\Gamma(x,x)}_{\frac{p}{2}}^{\frac{1}{2}},\norm{\Gamma(x^*,x^*)}_{\frac{p}{2}}^{\frac{1}{2}}\big\}\max\big\{\norm{\Gamma(y,y)}_{\frac{p}{2}}^{\frac{1}{2}},\norm{\Gamma(y^*,y^*)}_{\frac{p}{2}}^{\frac{1}{2}}\Big\}\\
&=\norm{a(x,x)}_{\L^{\frac{p}{2}}(\VN(G)) \oplus_\infty \L^{\frac{p}{2}}(\VN(G))}^{\frac{1}{2}}\norm{a(y,y)}_{\L^{\frac{p}{2}}(\VN(G)) \oplus_\infty \L^{\frac{p}{2}}(\VN(G))}^{\frac{1}{2}}.
\end{align*}
So $a$ satisfies the Cauchy-Schwarz inequality. The assertions concerning $\Gamma$ are similar. Suppose $x_n \xra[]{a} 0$ that is $x_n \in \P_G$, $x_n \to 0$ and $a(x_n-x_m,x_n-x_m) \to 0$. For any integer $n,m$, we have
\begin{align*}
\MoveEqLeft
\bnorm{A_p^{\frac12}(x_n-x_m)}_{\L^p(\VN(G))} \\
&\ov{\eqref{Equivalence-JMP2}}{\lesssim_p} 
\max \left\{ \bnorm{\Gamma(x_n-x_m,x-n-x_m)^{\frac12}}_{\L^p}, \bnorm{\Gamma((x_n-x_m)^*,(x_n-x_m)^*)^{\frac12}}_{\L^p} \right\} \\
&=\norm{a(x_n-x_m,x_n-x_m)}_{\L^{\frac{p}{2}}(\VN(G)) \oplus_\infty \L^{\frac{p}{2}}(\VN(G))}  
\to 0.
\end{align*}  
We infer that $\big(A_p^{\frac12}(x_n)\big)$ is a Cauchy sequence, hence a convergent sequence. Since $x_n \to 0$, by the closedness of $A_p^{\frac12}$, we deduce that $A_p^{\frac12}(x_n) \to 0$. Now, we have
\begin{align*}
\MoveEqLeft
\norm{a(x_n,x_n)}_{\L^{\frac{p}{2}}(\VN(G)) \oplus_\infty \L^{\frac{p}{2}}(\VN(G))}            
=\max \left\{ \bnorm{\Gamma(x_n,x_n)^{\frac12}}_{\L^p(\VN(G))}, \bnorm{\Gamma(x_n^*,x_n^*))^{\frac12}}_{\L^p(\VN(G))} \right\} \\
&\ov{\eqref{Equivalence-JMP2}}{\lesssim_p} \bnorm{A_p^{\frac12}(x_n)}_{\L^p(\VN(G))}
\to 0.
\end{align*}  
By Proposition \ref{Prop-closable}, we obtain the closability of the form $a$.

Let $x \in \L^p(\VN(G))$. By \eqref{Domaine-closure}, we have $x \in \dom \ovl{a}$ if and only if there exists a sequence $(x_n)$ of $\P_G$ satisfying $x_n \xra[]{a} x$, that is satisfying $x_n \to x$ and $\Gamma(x_n-x_m,x_n-x_m) \to 0$, $\Gamma((x_n-x_m)^*,(x_n-x_m)^*) \to 0 \text{ as } n,m \to \infty$. By \eqref{Equivalence-JMP2}, this is equivalent to the existence of a sequence $x_n \in \P_G$, such that $x_n \to x$, $A_p^{\frac12}(x_n-x_m) \to 0$ as $n,m \to \infty$. Now recalling that $\P_G$ is a core of the operator $A_p^{\frac12}$, we conclude that this is equivalent to $x \in \dom A_p^{\frac12}$.
\end{proof}

\begin{remark} \normalfont
\label{rem-after-lem-Gamma-closure-group}
The closability and biggest reasonable domain of the form $\Gamma$ are unclear.
\end{remark}

If $2 \leq p < \infty$, and $x,y \in \dom A_p^{\frac12}$, then we let $\Gamma(x,y)$ be the first component of $\ovl{a}(x,y)$, where $a \co \P_G \times \P_G \to \L^{\frac{p}{2}}(\VN(G)) \oplus_\infty \L^{\frac{p}{2}}(\VN(G)),\: (u,v) \mapsto \Gamma(u,v) \oplus \Gamma(u^*,v^*)$ is the form in Lemma \ref{lem-Gamma-closure-group}.

\begin{lemma}
\label{lem-Gamma-cut-off-contractive-sgrp}
Suppose $2 \leq p < \infty$. Let $G$ be a discrete group and assume that $\L^p(\VN(G))$ has $\CCAP$ and that $\VN(G)$ has $\QWEP$. Let $(\varphi_j)$ be a net of functions $\varphi_j \co G \to \C$ with finite support such that the net $(M_{\varphi_j})$ converges to $\Id_{\L^p(\VN(G))}$ in the point-norm topology with $\sup_j \norm{M_{\varphi_j}}_{\cb,\L^p(\VN(G) \to \L^p(\VN(G)} \leq 1$. If $x \in \dom A^{\frac12}_{p}$ then for any $j$
$$
\norm{\Gamma\big(M_{\varphi_j}(x),M_{\varphi_j}(x)\big)}_{\L^{\frac{p}{2}}(\VN(G))} 
\leq \bnorm{\partial_{\psi,q,p}(x)}_{\L^{p}(\VN(G),\L^2(\Gamma_q(H))_{c,p})}.
$$
\end{lemma}

\begin{proof}
Let $x \in \dom A^{\frac12}_{p}=\dom \partial_{\psi,q,p}$. Since $\P_G$ is a core of $\partial_{\psi,q,p}$ by Proposition \ref{Prop-derivation-closable-sgrp}, there exists a sequence $(x_n)$ of elements of $\P_G$ such that $x_n \to x$ and $\partial_{\psi,q,p}(x_n) \to \partial_{\psi,q,p}(x)$. Note that $M_{\varphi_j} \co \L^p(\VN(G)) \to \L^p(\VN(G))$ is a complete contraction. The linear map $M_{\varphi_j} \ot \Id_{\L^2(\Gamma_q(H))} \co \L^{p}(\VN(G),\L^2(\Gamma_q(H))_{c,p}) \to \L^{p}(\VN(G),\L^2(\Gamma_q(H))_{c,p})$ is also a contraction according to \eqref{extens-col-1}. We deduce that
\begin{align*}
\MoveEqLeft
\label{Equa-inter-1000}
\partial_{\psi,q,p} M_{\varphi_j}(x_n) 
=\sum_{s \in \supp \varphi_j} x_{ns} \partial_{\psi,q,p}(\lambda_s)
=\big(M_{\varphi_j} \ot \Id_{\L^2(\Gamma_q(H))}\big)\bigg(\sum_{s \in G} x_{ns} \partial_{\psi,q,p}(\lambda_s)\bigg) \\
&=\big(M_{\varphi_j} \ot \Id_{\L^2(\Gamma_q(H))} \big) (\partial_{\psi,q,p} x_n)
\xra[n \to +\infty]{} \big(M_{\varphi_j} \ot \Id_{\L^2(\Gamma_q(H))}\big) (\partial_{\psi,q,p} x).            
\end{align*}  
Since $\partial_{\psi,q,p} M_{\varphi_j}$ is bounded by \cite[Problem 5.22]{Kat1}, we deduce that $\partial_{\psi,q,p} M_{\varphi_j}(x)=\big(M_{\varphi_j} \ot \Id_{\L^2(\Omega)}\big) (\partial_{\psi,q,p} x)$. Now, we have
\begin{align*}
\MoveEqLeft
\bnorm{\Gamma(M_{\varphi_j}x,M_{\varphi_j}x)}_{\L^{\frac{p}{2}}}^{\frac{1}{2}} 
=\bnorm{\Gamma(M_{\varphi_j}x,M_{\varphi_j}x)^\frac{1}{2}}_{\L^p}
\ov{\eqref{nabla-norm-Lp-gradient-psi-Junge}}{=} \bnorm{\partial_{\psi,q,p}M_{\varphi_j}(x)}_{\L^{p}(\VN(G),\L^2(\Gamma_q(H))_{c,p})} \\
&=\bnorm{(M_{\varphi_j} \ot \Id_{\L^2(\Gamma_q(H))}) (\partial_{\psi,q,p} x)}_{\L^{p}(\VN(G),\L^2(\Gamma_q(H))_{c,p})}
\leq \bnorm{\partial_{\psi,q,p}(x)}_{\L^{p}(\VN(G),\L^2(\Gamma_q(H))_{c,p})}. 
\end{align*} 
\end{proof}

Now, we give a very concrete way to approximate the carr\'e du champ for a large class of groups.

\begin{lemma}
\label{lemma-agrandir-domaine-1-sgrp}
Let $2 \leq p < \infty$. Assume that $\L^p(\VN(G))$ has $\CCAP$ and that $\VN(G)$ has $\QWEP$. Let $(\varphi_j)$ be a net of functions $\varphi_j \co G \to \C$ with finite support such that the net $(M_{\varphi_j})$ converges to $\Id_{\L^p(\VN(G))}$ in the point-norm topology with $\sup_j \norm{M_{\varphi_j}}_{\cb,\L^p(\VN(G) \to \L^p(\VN(G)}\leq 1$. For any $x,y \in \dom A^{\frac12}_{p}$, we have in $\L^{\frac{p}{2}}(\VN(G))$
\begin{equation}
\label{Def-Gamma-group}
\Gamma(x,y)
=\lim_{j} \Gamma\big(M_{\varphi_j}(x),M_{\varphi_j}(y)\big).
\end{equation}
\end{lemma}

\begin{proof}
If $x \in \dom A^{\frac12}_{p}$, we have for any $j,k$
\begin{align*}
\MoveEqLeft
\norm{\Gamma(M_{\varphi_j}x-M_{\varphi_k}x,M_{\varphi_j}x-M_{\varphi_k}x)}_{\L^{\frac{p}{2}}(\VN(G)}^{\frac{1}{2}}  
=\bnorm{\Gamma(M_{\varphi_j}x-M_{\varphi_k}x,M_{\varphi_j}x-M_{\varphi_k}x)^{\frac{1}{2}}}_{\L^p(\VN(G))}   \\         
&\ov{\eqref{Equivalence-JMP2}}{\lesssim_p} \bnorm{A^{\frac12}(M_{\varphi_j}x-M_{\varphi_k}x)}_{\L^p(\VN(G))} 
= \bnorm{M_{\varphi_j} A_p^{\frac12} x-M_{\varphi_k} A_p^{\frac12} x}_{\L^p(\VN(G))}.
\end{align*}
The same inequality holds with $\Gamma((M_{\varphi_j}x-M_{\varphi_k}x)^*,(M_{\varphi_j}x-M_{\varphi_k}x)^*)$ on the left hand side.
Note that since $A^{\frac12}_p(x)$ belongs to $\L^{p}(\VN(G))$, $(M_{\varphi_j}A_p^{\frac12}(x))_p$ is a Cauchy net of $\L^p(\VN(G))$. Since $M_{\varphi_j}(x) \to x$, we infer that $M_{\varphi_j}(x) \xra{a} x$, where $a$ is the first form from Lemma \ref{lem-Gamma-closure-group}. Then \eqref{Def-Gamma-group} is a consequence of Lemma \ref{lem-Gamma-cut-off-contractive-sgrp} and Proposition \ref{Prop-existence-limit}. 
\end{proof}

\begin{lemma}
\label{Lemma-Hilbert-module-groups} \label{Prop-Schur-grad-Gamma62}
Suppose $2 \leq p<\infty$ and $-1 \leq q \leq 1$. 

\begin{enumerate}
	\item For any $x,y \in \dom A^{\frac12}_p=\dom \partial_{\psi,q,p}$, we have
\begin{equation}
\label{lien-Gamma-partial-psi}
\Gamma(x,y)
=\big\langle \partial_{\psi,q,p}(x),\partial_{\psi,q,p}(y) \big\rangle_{\L^p(\VN(G),\L^2(\Gamma_q(H))_{c,p})}
=\E\big(\big(\partial_{\psi,q,p}(x)\big)^* \partial_{\psi,q,p}(y)\big).
\end{equation}

\item For any $x \in \dom A^{\frac12}_p$, we have
\begin{equation}
\label{nabla-norm-Lp-gradient-psi}
\bnorm{\Gamma(x,x)^{\frac{1}{2}}}_{\L^p(\VN(G))}
=\bnorm{\partial_{\psi,q,p}(x)}_{\L^p(\VN(G),\L^2(\Gamma_q(H))_{c,p})}.
\end{equation}
\end{enumerate}
\end{lemma}

\begin{proof}
1. Consider some elements $x,y \in \dom A^{\frac12}_p=\dom \partial_{\psi,q,p}$. According to the point 3. of Proposition \ref{Prop-derivation-closable-sgrp}, $\P_G$ is a core of $\dom \partial_{\psi,q,p}$. So by \eqref{Def-core} there exist sequences $(x_n)$ and $(y_n)$ of $\P_G$ such that $x_n \to x$, $y_n \to y$, $\partial_{\psi,q,p}(x_n) \to \partial_{\psi,q,p}(x)$ and $\partial_{\psi,q,p}(y_n) \to \partial_{\psi,q,p}(y)$. By Lemma \ref{lem-contractive-crossed-product-Hilbert-valued}, we have, since $p \geq 2$, a contractive inclusion 
$$
\L^p(\Gamma_q(H) \rtimes_{\alpha} G) 
\ov{}{\hookrightarrow}  \L^p_{c}(\E)
\ov{\eqref{Egalite-fantastique-group}}{=} \L^p(\VN(G),\L^2(\Gamma_q(H))_{c,p}).
$$ 
We deduce that $\partial_{\psi,q,p}(x_n) \to \partial_{\psi,q,p}(x)$ and $\partial_{\psi,q,p}(y_n) \to \partial_{\psi,q,p}(y)$ in $\L^p(\VN(G),\L^2(\Gamma_q(H))_{c,p})$. We obtain
\begin{align*}
\MoveEqLeft
\norm{\Gamma(x_n-x_m,x_n-x_m)}_{\L^{\frac{p}{2}}(\VN(G))} \\
&\ov{\eqref{nabla-norm-Lp-gradient-psi-Junge}}{=} \norm{\big\langle \partial_{\psi,q,p}(x_n-x_m), \partial_{\psi,q,p}(x_n-x_m) \big\rangle_{S^p_I(\L^2(\Gamma_q(H))_{c,p})}}_{\L^{\frac{p}{2}}(\VN(G))}\\
&\leq \norm{\partial_{\psi,q,p}(x_n-x_m)}_{\L^p(\L^2(\Gamma_q(H))_{c,p})} \norm{\partial_{\psi,q,p}(x_n-x_m)}_{\L^p(\L^2(\Gamma_q(H))_{c,p})} 
\xra[n,m \to +\infty]{} 0
\end{align*}
and similarly for $x_n^*,y_n$ and $y_n^*$. We deduce that
\begin{align*}
\Gamma(x,y)
&\ov{\eqref{Def-closure}}{=}\lim_{n \to +\infty} \Gamma(x_n,y_n)
\ov{\eqref{lien-Gamma-partial-psi-junge}}{=} \lim_{n \to +\infty}  \big\langle \partial_{\psi,q,p}(x_n), \partial_{\psi,q,p}(y_n) \big\rangle_{\L^p(\VN(G),\L^2(\Gamma_q(H))_{c,p})}\\
&=\big\langle \partial_{\psi,q,p}(x), \partial_{\psi,q,p}(y) \big\rangle_{\L^p(\VN(G),\L^2(\Gamma_q(H))_{c,p})}.            
\end{align*}


\noindent

2. If $x \in \dom A^{\frac12}_p$, we have
\begin{align*}
\MoveEqLeft
\bnorm{\partial_{\psi,q,p}(x)}_{\L^p(\VN(G),\L^2(\Gamma_q(H))_{c,p})}          
\ov{\eqref{Norm-Lp-module}}{=}\Bnorm{\big\langle\partial_{\psi,q,p}(x),\partial_{\psi,q,p}(x)\big\rangle_{\L^p(\VN(G),\L^2(\Gamma_q(H))_{c,p})}^{\frac{1}{2}}}_{\L^p(\VN(G))} \\
&\ov{\eqref{lien-Gamma-partial-psi}}{=} \bnorm{\Gamma(x,x)^{\frac{1}{2}}}_{\L^p(\VN(G))}.
\end{align*}
\end{proof}

%

Now, we can extend Proposition \ref{prop-Riesz-Fourier}.

\begin{thm}
\label{cor-Riesz-equivalence-group-2}
Suppose that $\L^p(\VN(G))$ has $\CCAP$ and that $\VN(G)$ has $\QWEP$. Let $2 \leq p < \infty$. For any $x \in \dom A_p^{\frac12}$, we have
\begin{equation}
\label{max-Gamma-2-group}
\bnorm{A_p^{\frac12}(x)}_{\L^p(\VN(G))} 
\approx_p 
\max \left\{ \bnorm{\Gamma(x,x)^{\frac12}}_{\L^p(\VN(G))}, \bnorm{\Gamma(x^*,x^*)^{\frac12}}_{\L^p(\VN(G))} \right\}.
\end{equation}		
\end{thm}

\begin{proof}
Pick any $-1 \leq q \leq 1$. For any $x \in \dom A_p^{\frac12}$, first note that
\begin{equation}
\label{equa-inter-1000}
\E\big(\partial_{\psi,q,p}(x) (\partial_{\psi,q,p}(x)^*)\big)\big)
\ov{\eqref{relation-partial-star-group}}{=} \E\big(\partial_{\psi,q,p}(x^*)^* (\partial_{\psi,q,p}(x^*))\big)\big)
\ov{\eqref{lien-Gamma-partial-psi}}{=} \Gamma(x^*,x^*).
\end{equation}
We conclude that
\begin{align*}
\MoveEqLeft         
\bnorm{A_p^{\frac12}(x)}_{\L^p(\VN(G))} 
\ov{\eqref{Equivalence-square-root-domaine-Schur-sgrp}}{\approx_{p}} \bnorm{\partial_{\psi,q,p}(x)}_{\L^p(\Gamma_q(H) \rtimes_{\alpha} G)} \\
&\ov{\text{Proposition }\ref{Khintchine-group-twisted-full-space}}{\approx_p} \max \bigg\{ \norm{\big(\E((\partial_{\psi,q,p}(x))^* \partial_{\psi,q,p}(x))\big)^{\frac12}}_{p},\norm{\big(\E(\partial_{\psi,q,p}(x) (\partial_{\psi,q,p}(x)^*))\big)^{\frac12}}_{p} \bigg\}\\
&\ov{\eqref{lien-Gamma-partial-psi}\eqref{equa-inter-1000}}{=} \max \left\{ \bnorm{\Gamma(x,x)^{\frac12}}_{\L^p(\VN(G))}, \bnorm{\Gamma(x^*,x^*)^{\frac12}}_{\L^p(\VN(G))} \right\}.
\end{align*} 
\end{proof}

\subsection{Kato's square root problem for semigroups of Schur multipliers}
\label{Kato-Schur}


In this subsection, we shall consider the Kato square root problem for markovian semigroups of Schur multipliers. Thus, we fix for the whole subsection such a markovian semigroup from Proposition \ref{def-Schur-markovian}, with its generator $A_p$ on $S^p_I$ and also the gradient type operator $\partial_{\alpha,q}$ from \eqref{def-delta-alpha}. Kato's square root problem is then the question whether $A^{\frac12}_p(x)$ and $\partial_{\alpha,q}(x)$ are comparable in $\L^p$ norm. The main results in this subsection, answering affirmatively to this question, are then Theorem \ref{Thm-Riesz-equivalence-Schur} (case of classical Gaussians) and Proposition \ref{Prop-equiv-q-gaussians} (case of $q$-deformed Gaussians), together with Lemma \ref{Khintchine-Schur} on a Khintchine type equivalence in some Gaussian subspace of $\L^p(\Gamma_q(H) \otvn \B(\ell^2_I))$. Finally, the problem of exact description of the domain of the closure of the gradients is investigated in Proposition \ref{Prop-derivation-closable}. Throughout this subsection, if $H$ is a Hilbert space, we denote by $H_\disc$ the abelian group $(H,+)$ equipped with the \textit{discrete} topology. We will use the trace preserving normal unital injective $*$-homomorphism map  
\begin{equation}
\label{equ-trivial-map-J}
J \co \VN(H_\disc) \to \L^\infty(\Omega) \otvn \VN(H_\disc),\: \lambda_h  \mapsto  1 \ot \lambda_h 
\end{equation}
and the unbounded Fourier multiplier $(-\Delta)^{-\frac12} \co \P_{H_\disc} \subseteq \L^p(\VN(H_{\disc})) \to \L^p(\VN(H_{\disc}))$ defined by
\begin{equation}
\label{Map-moins-Delta-1}
(-\Delta)^{-\frac12}(\lambda_h)
\ov{\mathrm{def}}{=} \frac{1}{2\pi\norm{h}_H}\lambda_h.
\end{equation}
The following is inspired by \cite{Arh1} \cite{Arh4}.

\begin{prop}
\label{Prop-3-1}
Let $H$ be a Hilbert space. There exists a unique weak* continuous group $(U_t)_{t \in \R}$ of $*$-automorphisms of $\L^\infty(\Omega) \otvn \VN(H_\disc)$ such that
\begin{equation}
\label{Def-U_t}
U_t(f \ot \lambda_h)
=\e^{\sqrt{2}\i t \W(h)}f \ot \lambda_h, \qquad t \in \R,f \in \L^\infty(\Omega), h \in H.
\end{equation}
Moreover, each $U_t$ is trace preserving.
\end{prop}

\begin{proof}
For any $t \in \R$, we consider the (continuous) function $u_t \co H_\disc \to \mathrm{U}(\L^\infty(\Omega))$, $h \mapsto \e^{-\sqrt{2}\i t \W(h)}$. For any $t \in \R$ and any $h_1,h_2 \in H_\disc$, note that  
\begin{align*}
\MoveEqLeft
u_t(h_1+h_2)
=\e^{-\sqrt{2}\i t \W(h_1+h_2)}
=\e^{-\sqrt{2}\i t \W(h_1)}\e^{-\sqrt{2}\i t \W(h_2)}
=u_t(h_1)u_t(h_2).            
\end{align*}
By \cite[Proposition 2.4]{Arh7} applied with $M=\L^\infty(\Omega)$, $G=H_\disc$ and by considering the trivial action $\alpha$, for any $t \in \R$, we have a unitary $V_t \co \L^2(H_\disc,\L^2(\Omega)) \to \L^2(H_\disc,\L^2(\Omega))$, $\xi \mapsto (h \mapsto u_t(-h)(\xi(h)))$ and a $*$-isomorphism $U_t \co \L^\infty(\Omega) \otvn \VN(H_\disc)  \to \L^\infty(\Omega) \otvn \VN(H_\disc)$, $x \mapsto V_txV_t^*$ satisfying \eqref{Def-U_t}. The uniqueness is clear by density.

For any $t,t' \in \R$ any $\xi \in \L^2(H_\disc,\L^2(\Omega))$, note that almost everywhere
\begin{align*}
\MoveEqLeft
\big(V_{t}V_{t'}(\xi)\big)(h)            
		=u_t(-h)\big((u_{t'}(-h)(\xi(h))\big)
		=\e^{-\sqrt{2}\i t \W(-h)}\e^{-\sqrt{2}\i t' \W(-h)}\xi(h)\\
		&=\e^{-\sqrt{2}\i (t+t') \W(-h)}\xi(h)
		=u_{t+t'}(-h)(\xi(h))
		=\big(V_{t+t'}(\xi)\big)(h).
\end{align*} 
We conclude that $V_{t}V_{t'}=V_{t+t'}$. Moreover, for any $\xi,\eta \in \L^2(H_\disc,\L^2(\Omega))=\L^2(H_\disc \times \Omega)$, using dominated convergence theorem, we obtain
\begin{align*}
\MoveEqLeft
\big\langle V_t(\xi), \eta \big\rangle_{\L^2(H_\disc,\L^2(\Omega))}           
		=\int_{H_\disc \times \Omega} \ovl{V_t(\xi)(h)(\omega) \eta(h,\omega)} \d\mu_{H_\disc}(h)\d \mu(\omega) \\
		&=\int_{H_\disc \times \Omega} \ovl{u_t(-h)(\xi(h))(\omega)} \eta(h,\omega) \d\mu_{H_\disc}(h)\d \mu(\omega) \\
		&=\int_{H_\disc \times \Omega} \e^{\sqrt{2}\i t \W(-h)(\omega)}\ovl{\xi(h,\omega)} \eta(h,\omega)\d\mu_{H_\disc}(h)\d \mu(\omega).
		\\
		&\xra[t \to 0]{} \int_{H_\disc \times \Omega} \ovl{\xi(h,\omega)} \eta(h,\omega)\d\mu_{H_\disc}(s)\d \mu(\omega)
		=\langle \xi, \eta \rangle_{\L^2(H_\disc,\L^2(\Omega))}.	
\end{align*} 
So $(V_t)_{t \in \R}$ is a weakly continuous group of unitaries hence a strongly continuous group by \cite[Lemma 13.4]{Str} or \cite[p.~239]{Tak2}. By \cite[p.~238]{Tak2}, we conclude that $(U_t)_{t \in \R}$ is a weak* continuous group of $*$-automorphisms. 

Finally, for any $t \geq 0$, any $f \in \L^\infty(\Omega)$ and any $h \in H$, we have
\begin{align*}
\MoveEqLeft
\bigg(\int_\Omega \cdot \ot \tau_{\VN(H_{\disc})}\bigg)\big(U_t(f \ot \lambda_h)\big)            
=\bigg(\int_\Omega \cdot \ot \tau_{\VN(H_{\disc})}\bigg)\Big(\e^{\sqrt{2}\i t \W(h)}f \ot \lambda_h\Big)\\
&=\bigg(\int_\Omega \e^{\sqrt{2}\i t \W(h)}f \d \mu \bigg) \tau_{\VN(H_{\disc})}(\lambda_h)
=\bigg(\int_\Omega \e^{\sqrt{2}\i t \W(h)}f \d \mu\bigg) \delta_{0,h}\\
&=\bigg(\int_\Omega f \d \mu \bigg) \tau_{\VN(H_{\disc})}(\lambda_h)
=\bigg(\int_\Omega \cdot \ot \tau_{\VN(H_{\disc})}\bigg)(f \ot \lambda_h).
\end{align*} 
\end{proof}

We consider the unbounded operator 
$$
\delta \co \Span\left\{ \lambda_s : \: s \in H_{\disc} \right\} \subseteq \L^p(\VN(H_{\disc})) \to \L^p(\Omega,\L^p(\VN(H_{\disc})))
$$ 
defined by
\begin{equation}
\label{Map-delta-1}
\delta(\lambda_h)
\ov{\mathrm{def}}{=} 2\pi\i  \W(h) \ot \lambda_h.
\end{equation}

Recall the classical transference principle \cite[Theorem 2.8]{BGM2}. Let $G$ be a locally compact abelian group and $G \to \B(X)$, $t \to \pi_t$ be a strongly continuous representation of $G$ on a Banach space $X$
such that $c =\sup \{\norm{\pi_t} : t \in G \} < \infty$. Let $k \in \L^1(G)$ and let $T_k \co X \to X$ be the operator defined by $T_k(x)=\int_G k(t)\pi_{-t}(x) \d \mu_G(t)$. Then 
\begin{equation}
\label{Transference}
\norm{T_k}_{X \to X} 
\leq c^2\norm{k* \cdot}_{\L^p(G,X)\to \L^p(G,X)}.
\end{equation}
We will use the function $k_{\epsi,R}(t)=\frac{1}{\pi t} 1_{\epsi < |t|<R}$ \cite[p.~388]{HvNVW1}.

Recall that we say that a function $f \in \L^1_{\mathrm{loc}}(\R^*,X)$ admits a Cauchy principal value if the limit $\lim_{\epsi \to 0^+}\big(\int_{-\frac{1}{\epsi}}^{-\epsi} f(t)\d t+\int_{\epsi}^{\frac{1}{\epsi}} f(t)\d t\big)$ exists and we let 
$$
\pv \int_\R f(t) \d t 
\ov{\mathrm{def}}{=} \lim_{\epsi \to 0^+} \bigg(\int_{-\frac{1}{\epsi}}^{-\epsi} f(t) \d t+\int_{\epsi}^{\frac{1}{\epsi}}f(t) \d t\bigg).
$$
 
\begin{prop}
\label{prop-H-boundedness}

Suppose $1< p <\infty$. For any Hilbert space $H$, the map
\begin{equation}
\label{Def-calH}
\begin{array}{cccc}
\mathcal{H}_{\disc} \co & \L^p(\Omega,\L^p(\VN(H_\disc))) &  \longrightarrow & \L^p(\Omega,\L^p(\VN(H_\disc))) \\
       & z  &  \longmapsto & \dsp  \pv \frac{1}{\pi} \int_\R U_t(z) \frac{\d t}{t}   \\
\end{array}.
\end{equation}
is well-defined and completely bounded.
\end{prop}

\begin{proof}
For any $\epsi>0$ large enough, using the fact that the space $\L^p(\Omega,\L^p(\VN(H_\disc)),S^p_I)$ is UMD and \cite[p.~485]{AKM} \cite[Theorem 13.5]{Haa2} in the last inequality, we have by transference 
\begin{align*}
\MoveEqLeft
\norm{\frac{1}{\pi} \int_{\epsi < |t|<\frac{1}{\epsi}} (U_t \ot \Id_{S^p_I}) \frac{\d t}{t}}_{\L^p(\Omega,\L^p(\VN(H_\disc)),S^p_I) \to \L^p(\Omega,\L^p(\VN(H_\disc)),S^p_I)} \\
&=\norm{\int_{\R} k_{\epsi,\frac{1}{\epsi}}(t) (U_t \ot \Id_{S^p_I}) \d t}_{\L^p(\Omega,\L^p(\VN(H_\disc)),S^p_I) \to \L^p(\Omega,\L^p(\VN(H_\disc)),S^p_I)}\\           
&\ov{\eqref{Transference}}{\leq} \norm{\big(k_{\epsi,\frac{1}{\epsi}} *\cdot \big) \ot \Id_{\L^p(\Omega,\L^p(\VN(H_\disc)),S^p_I)}}_{\L^p(\R \times \Omega,\L^p(\VN(H_\disc)),S^p_I)) \to \L^p(\R \times \Omega,\L^p(\VN(H_\disc)),S^p_I))} \\
&\lesssim_p 1.
\end{align*}  
If $z=f \ot \lambda_h$, we have
\begin{align*}
\MoveEqLeft
\int_{-\frac{1}{\epsi}}^{-\epsi} U_t(z) \frac{\d t}{t} +\int_{\epsi}^{\frac{1}{\epsi}} U_t(z) \frac{\d t}{t}          
=\int_{\epsi}^{\frac{1}{\epsi}} (U_t(z)-U_{-t}(z)) \frac{\d t}{t}\\
&\ov{\eqref{Def-U_t}}{=} 2\i \int_{\epsi}^{\frac{1}{\epsi}} \sin (\sqrt{2}\W(h)t)(f \ot \lambda_h) \frac{\d t}{t}
=2\i\bigg(\int_{\epsi}^{\frac{1}{\epsi}} \sin (\sqrt{2}\W(h)t) \frac{\d t}{t}\bigg)(f \ot \lambda_h)
\end{align*}  
which admits a limit when $\epsi \to 0$. We have the existence of the principal value by linearity and density. 
\end{proof}

The following is a variation of Pisier formula. Here $Qf = \sum_k \left(\int_\Omega f \gamma_k\right) \cdot \gamma_k$ is the Gaussian projection where $(\gamma_k)$ is a family of independent standard Gaussian variables given by $\gamma_k = \W(e_k)$, here $e_k$ running through an orthonormal basis of $H$.
(Note that if $H$ were to be non-separable, then for any fixed $f \in \L^2(\Omega)$, only countably many $\int_\Omega f \gamma_k$ would be non-zero.)
Then $Q$ is independent of the particular choice of the $e_k$\footnote{\thefootnote. If $(e_k)$ and $(f_j)$ are two such orthonormal bases of $H$, then $f_j = \sum_k u_{jk} e_k$, where the sum is over at most countably many $k$ and the coefficients $u_{jk}$ form an orthogonal matrix, $\sum_k u_{jk} u_{ik} = \delta_{ji}$.
We have 
\begin{align*}
\MoveEqLeft
Q_{(f_j)}(f) 
= \sum_j \langle \W(f_j) , f\rangle \W(f_j) 
= \sum_j \bigg\langle  \sum_k u_{jk} \W(e_k), f \bigg\rangle \sum_l u_{jl} \W(e_l) 
= \sum_{k,l} \sum_j u_{jk} u_{jl} \langle \W(e_k), f \rangle \W(e_l) \\
&= \sum_{k,l} \delta_{kl} \langle \W(e_k), f \rangle \W(e_l) 
= \sum_l \langle \W(e_k), f \rangle \W(e_k) 
= Q_{(e_k)}(f).
\end{align*}}.
Moreover, $J$ is given in \eqref{equ-trivial-map-J} and $(-\Delta)^{-\frac12}$ is given in \eqref{Map-moins-Delta-1} at the beginning of this subsection.

\begin{prop}
\label{Prop-variation-pisier}
Let $1 < p < \infty$.
For any $h \in H$, we have 
\begin{equation}
\label{Pisier-formula-Schur-2}
\sqrt{\frac{2}{\pi}} \delta (-\Delta)^{-\frac{1}{2}}(\lambda_h)
=\big(Q \ot \Id_{\L^p(\VN(H_{\disc}))}\big)\mathcal{H_{\disc}}J(\lambda_h).
\end{equation}
\end{prop}

\begin{proof}
For any $h \in H$, we have
\begin{align*}
\MoveEqLeft
 \sqrt{\frac{2}{\pi}}  \delta (-\Delta)^{- \frac12} (\lambda_h)
		\ov{\eqref{Map-moins-Delta-1}}{=}\frac{1}{\sqrt{2\pi^3}\norm{h}_H} \delta(\lambda_h)
		\ov{\eqref{Map-delta-1}}{=}\frac{\sqrt{2} \i}{\sqrt{\pi}\norm{h}_H}  \W(h) \ot \lambda_h.
\end{align*}   
Recall that for any $\alpha \in \R$ we have $\int_{0}^{+\infty} \frac{\sin(\alpha t)}{t} \d t=\sign(\alpha)\frac{\pi}{2}$. Consequently, we have
\begin{align*}
\MoveEqLeft
  \big(Q \ot \Id_{\L^p(\VN(H_{\disc}))}\big)\mathcal{H}_{\disc}J(\lambda_h)  
		\ov{\eqref{Def-calH}}{=}\big(Q \ot \Id_{\L^p(\VN(H_{\disc}))}\big)\bigg(\pv \frac{1}{\pi} \int_\R \e^{\i \sqrt{2}t \W(h)}\frac{\d t}{t} \ot \lambda_h\bigg)\\
		&=Q\bigg(\pv \frac{1}{\pi} \int_\R \e^{\i \sqrt{2}t \W(h)}\frac{\d t}{t}\bigg) \ot \lambda_h
		=\frac{1}{\pi}Q\bigg(\pv\int_\R \e^{\i \sqrt{2} t \W(h)}\frac{\d t}{t}\bigg) \ot \lambda_h\\
		&=\frac{2\i}{\pi}Q\bigg(\int_{0}^{+\infty} \sin\big(\sqrt{2}t\W(h)\big) \frac{\d t}{t}\bigg) \ot \lambda_h   
		=\i Q\big(\mathrm{sign} \circ \W(h)\big) \ot \lambda_h.
\end{align*}
The Gaussian projection $Q$ is independent of the choice of the family $(\gamma_k)$ of $\L^2(\Omega)$ consisting of independent standard Gaussian variables stemming from the process $\W(e),\:e \in H$. Choosing $\gamma_1 \ov{\mathrm{def}}{=} \frac{\W(h)}{\norm{\W(h)}_{\L^2(\Omega)}} \ov{\eqref{esperance-isonormal}}{=} \frac{\W(h)}{\norm{h}_H}$, for any $k \geq 2$ the random variables $\gamma_k$ and $\W(h)$ are independent. Thus $\gamma_k$ and $\mathrm{sign} \circ \W(h)$ are also independent. So $\int_\Omega \mathrm{sign} \circ \W(h) \cdot \gamma_k = \int_\Omega \mathrm{sign} \circ \W(h) \int_\Omega \gamma_k = 0$. Using \cite[p.~100]{DeF1} or \cite[(E.2)]{HvNVW2} in the last equality, we infer that
\begin{align*}
\MoveEqLeft
Q (\mathrm{sign} \circ \W(h)) 
= \left( \int_\Omega \mathrm{sign} \circ \W(h) \cdot \gamma_1  \right) \gamma_1 = \left( \int_\Omega \mathrm{sign} \circ \gamma_1 \cdot \gamma_1 \right) \gamma_1 
= \left( \int_\Omega |\gamma_1| \right) \gamma_1 
=\sqrt{\frac{2}{\pi}} \frac{\W(h)}{\norm{h}_H}.
\end{align*}
We conclude that
$$
\big(Q \ot \Id_{\L^p(\VN(H_{\disc}))}\big)\mathcal{H}_{\disc}J(\lambda_h)  
=\frac{\i\sqrt{2}}{\sqrt{\pi}\norm{h}_H} \W(h) \ot \lambda_h.
$$
\end{proof}

\begin{lemma} 
\label{boundedness-Schur-multipliers}
If $1 < p< \infty$, the operator
\begin{equation}
\label{bounded-Riesz}
\delta (-\Delta)^{- \frac12} \ot \Id_{S^p_I} \co \L^p\big(\VN(H_{\disc}), S^p_I \big) \to \L^p\big(\Omega,\L^p(\VN(H_{\disc}),S^p_I)\big)
\end{equation}
is well-defined and bounded.
\end{lemma}

\begin{proof}
Note that the map $J \co \L^p(\VN(H_\disc)) \to \L^p(\Omega,\L^p(\VN(H_\disc))$ is completely bounded. By Proposition \ref{prop-H-boundedness}, the operator $\mathcal{H}_{\disc} \ot \Id_{S^p_I} \co \L^p(\Omega,\L^p(\VN(H),S^p_I)) \to \L^p(\Omega,\L^p(\VN(H),S^p_I))$ is well-defined and bounded. Moreover, since the Banach space $\L^p(\VN(H_{\mathrm{disc}}),S^p_I)$ is K-convex, the map $Q \ot \Id_{\L^p(\VN(H_{\mathrm{disc}}),S^p_I)}$ induces a bounded operator on $\L^p(\Omega,\L^p(\VN(H_{\mathrm{disc}}),S^p_I))$. Using the equality \eqref{Pisier-formula-Schur-2}, we obtain the result by composition.
\end{proof}

For the remainder of this subsecion, we consider and fix a markovian semigroup $(T_t)_{t \geq 0}$ of Schur multipliers on $\B(\ell^2_I)$ from Definition \ref{def-Schur-markovian}.
Recall the real Hilbert space $H$ and the family $(\alpha_i)_{i \in I}$ in $H$ from \eqref{equ-Schur-markovian-alpha}.


Note that $u=\mathrm{diag}(\lambda_{\alpha_i}: i \in I)$ is a unitary of the von Neumann algebra $\M_I(\VN(H_{\disc}))=\VN(H_{\disc}) \otvn \B(\ell^2_I)$. Hence we have a trace preserving normal unital injective homomorphism $\pi \co \B(\ell^2_I) \to \VN(H_{\disc}) \otvn \B(\ell^2_I)$, $x \mapsto u(1 \ot x)u^*$. For any $i,j \in I$, note that
\begin{equation}
\label{Map-pi-1}
\pi(e_{ij})
=\lambda_{\alpha_i -\alpha_j} \ot e_{ij}.
\end{equation}
Moreover, if $1 \leq p \leq \infty$ and if $y \in \L^p(\Omega,S^p_I)$, we have
\begin{equation}
\label{pi-isometry}
\norm{(\Id_{\L^p(\Omega)} \ot \pi_p)(y)}_{\L^p(\Omega,\L^p(\VN(H_{\disc}),S^p_I))}
=\norm{y}_{\L^p(\Omega,S^p_I)}.
\end{equation}

%
%

The next two lemmas prepare our first main result in this subsection, which is Theorem \ref{Thm-Riesz-equivalence-Schur}.
In the following lemma, we consider the operator $A_p^{-\frac12}$ on $\M_{I,\fin} \subseteq S^p_I$, given as the Schur multiplier $A_p^{-\frac12}e_{ij} = a_{ij}^{-\frac12} e_{ij}$.
Note that this makes sense when $a_{ij} \neq 0$, so for $e_{ij} \in \Ran(A_p)$, and we interpret $A_p^{-\frac12}x$ to be $0$ if $x \in \ker(A_p)$.
Note that by theory of  markovian semigroups, in general we have a direct sum decomposition $S^p_I = \ovl{\Ran{A_p}} \oplus \ker(A_p)$ \cite[p.~361]{HvNVW2}.
Then the statement of \eqref{intertwining-formula} below then includes that the right hand side vanishes for $x \in \ker(A_p)$, and both sides are bounded operators on $\Ran(A_p)$, and thus on $\ovl{\Ran(A_p)}$.

\begin{lemma}[intertwining formula]
\label{Lemma-magic-formula-1}
Suppose $1 \leq p <\infty$. For any $x \in \M_{I,\fin}$, we have
\begin{equation}
\label{intertwining-formula}
\i \big(\Id_{\L^p(\Omega)} \ot \pi_p \big) \partial_\alpha A_p^{-\frac12}
=\big(\delta (-\Delta)^{-\frac12} \ot \Id_{S^p_I} \big) \pi_p.
\end{equation} 
\end {lemma}

\begin{proof}
For any $i,j \in I$, we have
\begin{align*}
\MoveEqLeft
\i \big(\Id_{\L^p(\Omega)} \ot \pi_p \big) \partial_{\alpha,1} A_p^{-\frac12} (e_{ij})    
=\frac{\i}{\norm{\alpha_i-\alpha_j}_{H}}\big(\Id_{\L^p(\Omega)} \ot \pi_p \big) \partial_\alpha (e_{ij}) \\
&\ov{\eqref{def-delta-alpha}}{=}\frac{\i}{\norm{\alpha_i-\alpha_j}_{H}}\big(\Id_{\L^p(\Omega)} \ot \pi_p \big) \big(\W(\alpha_i-\alpha_j) \ot e_{ij}\big)  
=\frac{\i}{\norm{\alpha_i-\alpha_j}_{H}}\W(\alpha_i-\alpha_j) \ot \pi_p(e_{ij}) \\
&\ov{\eqref{Map-pi-1}}{=} \frac{\i}{\norm{\alpha_i-\alpha_j}_{H}}\W(\alpha_i-\alpha_j) \ot \lambda_{\alpha_i -\alpha_j} \ot e_{ij}
\end{align*}
and
\begin{align*}
\MoveEqLeft
\big( \delta (-\Delta)^{-\frac12} \ot \Id_{S^p_I} \big) \pi_p(e_{ij}) 
\ov{\eqref{Map-pi-1}}{=} \big( \delta (-\Delta)^{-\frac12} \ot \Id_{S^p_I} \big) (\lambda_{\alpha_i -\alpha_j} \ot e_{ij})
= \delta (-\Delta)^{-\frac12} (\lambda_{\alpha_i -\alpha_j}) \ot e_{ij}) \\ 
&\ov{\eqref{Map-moins-Delta-1}}{=}\frac{1}{2 \pi \norm{\alpha_i-\alpha_j}_{H}} \delta (\lambda_{\alpha_i -\alpha_j}) \ot e_{ij}
\ov{\eqref{Map-delta-1}}{=}\frac{\i}{\norm{\alpha_i-\alpha_j}_{H}}\big(\W(\alpha_i-\alpha_j) \ot \lambda_{\alpha_i -\alpha_j} \big) \ot e_{ij}.
\end{align*}
We conclude by linearity.
\end{proof}

%

For future use, we state the next lemma equally for the case of $q$-deformed algebras.

\begin{lemma}
\label{lem-formule-trace-trace}
For $-1 \leq q \leq 1$ and any $x,y \in \M_{I,\fin}$, we have
\begin{equation}
\label{formule-trace-trace}
(\tau_{\Gamma_q(H)} \ot \tr_{\B(\ell^2_I)})\big((\partial_{\alpha,q}(x))^* \partial_{\alpha,q}(y) \big)
=\tr_{\B(\ell^2_I)} \big((A^\frac12(x))^* A^\frac12(y)\big).
\end{equation}
\end{lemma}

\begin{proof}
For any $i,j,k,l \in I$, we have
\begin{align*}
\MoveEqLeft
\tr\big((A^\frac12(e_{ij}))^* A^\frac12(e_{kl})\big)   
=\norm{\alpha_i-\alpha_j}_H\norm{\alpha_k-\alpha_l}_H\tr(e_{ij}^* e_{kl}) 
=\delta_{i=k} \delta_{j=l}\norm{\alpha_i-\alpha_j}_H^2.
\end{align*}
Moreover, we have
\begin{align*}
\MoveEqLeft
(\tau \ot \tr)\big((\partial_{\alpha,q}(e_{ij}))^* \partial_{\alpha,q}(e_{kl}) \big) 
\ov{\eqref{def-delta-alpha}}{=} (\tau \ot \tr)\big(( s_q(\alpha_i - \alpha_j) \ot e_{ij})^*  (s_q(\alpha_k - \alpha_l) \ot e_{kl}) \big) \\
&= \sum_{i,j,k,l} \tau(s_q(\alpha_i - \alpha_j)s_q(\alpha_k - \alpha_l))\tr(e_{ij}^*e_{kl}) 
\ov{\eqref{petit-Wick}, \eqref{esperance-isonormal}}{=} \delta_{i=k} \delta_{j=l} \norm{\alpha_i - \alpha_j}^2_H.
\end{align*}
We conclude by linearity.
\end{proof}

Suppose $1 \leq p< \infty$. We denote by $A_p^{\frac12} \co \dom A^{\frac12}_{p} \subseteq S^{p}_I\to S^{p}_I$ the square root of the sectorial operator $A_p \co \dom A_{p} \subseteq S^{p}_I\to S^{p}_I$.
It is again a Schur multiplier, associated with symbol $\norm{\alpha_i-\alpha_j}$.

\begin{thm}
\label{Thm-Riesz-equivalence-Schur}
Suppose $1<p<\infty$. For any $x \in \M_{I,\fin}$, we have
\begin{equation} 
\label{OneLineDimFree}
\bnorm{A_p^{\frac12}(x)}_{S^p_I} 
\approx_{p} \norm{\partial_{\alpha,1}(x)}_{\L^p(\Omega,S^p_I)}. 
\end{equation}
\end{thm}

\begin{proof} 
We write in short $\partial_{\alpha,1} = \partial_\alpha$ in this proof.
Using Lemma \ref{Lemma-magic-formula-1}, for any $x \in \M_{I,\fin}$, we obtain the inequality 
\begin{align*}
\MoveEqLeft
\norm{\partial_\alpha(x)}_{\L^p(\Omega,S^p_I)} 
=\bnorm{\partial_\alpha A_p^{-\frac12} A_p^{\frac12}(x)}_{\L^p(\Omega,S^p_I)} \\
&\ov{\eqref{pi-isometry}}{=} \bnorm{\big(\Id_{\L^p(\Omega)} \ot \pi_p\big) \partial_\alpha A_p^{-\frac12} A_p^{\frac12} (x)}_{\L^p(\Omega,\L^p(\VN(H_{\disc}),S^p_I))} \\
&\ov{\eqref{intertwining-formula}}{=}\bnorm{ \big(\delta (-\Delta)^{-\frac12} \ot \Id_{S^p_I} \big) \pi_p A_p^{\frac12}(x)}_{\L^p(\Omega,\L^p(\VN(H_{\disc}),S^p_I))}\\
&\leq \bnorm{ \big(\delta (-\Delta)^{-\frac12} \ot \Id_{S^p_I} \big)}_{\L^p(\VN(H_{\disc}),S^p_I) \to \L^p(\Omega,\L^p(\VN(H_{\disc}),S^p_I))} \norm{\pi_p} \bnorm{A_p^{\frac12}(x)}_{S^p_I}\\
&\ov{\eqref{bounded-Riesz}}{\lesssim_p} \bnorm{A_p^{\frac12}(x)}_{S^p_I}. 
\end{align*}  
The reverse estimate follows with essentially the same constant from a duality argument. Indeed, if we fix $x$ to be an element of $S^p_I$, by Hahn-Banach theorem, there exists an element $y$ of $S^{p^*}_I$ with $\norm{y}_{S^{p^*}_I} = 1$ satisfying 
\begin{equation}
\label{Riesz-dual-y}
\bnorm{A_p^{\frac{1}{2}}(x)}_{S^p_I}
=\tr\big(y^* A_p^\frac12(x) \big).
\end{equation}
We can suppose that $y \in \M_{I,\fin}$. Note that $A^{-\frac{1}{2}}$ is well-defined on $e_{ij}$ only if $\norm{\alpha_i-\alpha_j} \not=0$. We consider on the set $I$ the equivalence relation $i \cong j \ov{\mathrm{def}}{\Longleftrightarrow} \alpha_i = \alpha_j$. Then we can write $I = \bigcup_{k \in K}J_k$ partitioned into pairwise disjoint subsets $J_k$ according to this equivalence relation. Consider the subalgebra $M = \ovl{\bigoplus_{k \in K} \B(\ell^2_{J_k})}^{\mathrm{w}^*}$ of $\B(\ell^2_I)$ and the conditional expectation $\E_I \co S^p_I \to \L^p(M,\tr|M)$. We obtain an element $y_0 \ov{\mathrm{def}}{=} y - \E_I(y)$ with 
$$
\norm{y_0}_{S^{p^*}_I}
=\norm{y- \E_I(y)}_{S^{p^*}_I} 
\leq \norm{y}_{S^{p^*}_I}+\norm{\E_I(y)}_{S^{p^*}_I} 
\leq 2.
$$
Note that $y_{0ij} = 0$ if $\alpha_i = \alpha_j$, hence $A_p^{-\frac12}(y_0)$ is well-defined. For any $i \in I$, we have 
\begin{align*}
\MoveEqLeft
\Big[(\E_I(y))^* A_p^\frac12(x)\Big]_{ii}  
=\sum_{k \in I} (\E_I(y)^*)_{ik} \norm{\alpha_k - \alpha_i} x_{ki} 
=\sum_{j=1}^{|K|} \sum_{l \in J_j} (\E_I(y)^*)_{il} \norm{\alpha_l - \alpha_i} x_{li} \\
&=\sum_{j=1}^{|K|} \sum_{l \in J_j} \delta_{i \in J_j} y_{li} \norm{\alpha_l - \alpha_i} x_{li}
=\sum_{j=1}^{|K|} \sum_{l \in J_j} \delta_{i \in J_j} y_{li} \cdot 0 \cdot x_{li}
=0.
\end{align*}
Hence, $\tr\big((\E_I(y))^* A_p^\frac12(x)\big) = 0$. Using the first part of the proof in the last estimate, we conclude that
\begin{align*}
\MoveEqLeft
\bnorm{A_p^\frac{1}{2}(x)}_{S^p_I}
\ov{\eqref{Riesz-dual-y}}{=}\tr \big(y^* A_p^\frac12(x) \big) 
=\tr \big(y^* A_p^\frac12(x) \big) - \tr \big( (\E_I(y))^* A_p^\frac12 (x)\big) 
=\tr \big(y_0^* A_p^\frac12(x) \big) \\
&=\tr \big((A_{p^*}^\frac12 A_{p^*}^{-\frac12}(y_0))^* A^\frac12(x) \big) 
\ov{\eqref{formule-trace-trace}}{=}\tr_{\L^\infty(\Omega) \ot \B(\ell^2_I)} \big((\partial_\alpha A_{p^*}^{-\frac12}(y_0))^* \partial_\alpha(x) \big) \\
&\leq \bnorm{\partial_\alpha A_{p^*}^{-\frac12}(y_0)}_{\L^{p^*}(\Omega,S^{p^*}_I)} \norm{\partial_\alpha(x)}_{\L^p(\Omega,S^p_I)} 
\lesssim_p \norm{\partial_\alpha(x)}_{\L^p(\Omega,S^p_I)}.
\end{align*}
\end{proof}

\begin{remark} \normalfont
\label{rem-constants-Riesz-equivalence}
Keeping track of the constants in the two-sided estimate of Theorem \ref{Thm-Riesz-equivalence-Schur}, we obtain
\begin{equation}
\label{equ-constants-Riesz-equivalence}
\frac{1}{K \max(p,p^*)^{\frac32}} \bnorm{\partial_{\alpha,1}(x)}_{\L^p(\Omega,S^p_I)} 
\leq \bnorm{A^{\frac12}_p(x)}_{S^p_I} \leq K \max(p,p^*)^{\frac32} \bnorm{\partial_{\alpha,1}(x)}_{\L^p(\Omega,S^p_I)},
\end{equation}
where $K$ is an absolute constant.
Indeed, for the upper estimate, looking into the proof of Theorem \ref{Thm-Riesz-equivalence-Schur}, we have a control by
\begin{align*}
\MoveEqLeft
\norm{\delta(-\Delta)^{-\frac12} \ot \Id_{S^p_I} }_{\L^p(\VN(H_\disc),S^p_I) \to \L^p(\Omega,\L^p(\VN(H_\disc),S^p))} \norm{\pi_p}_{p \to p} \\
& \overset{\text{Proposition }\ref{Prop-variation-pisier}}{=} \sqrt{\frac{\pi}{2}} \norm{Q \ot \Id_{\L^p(\VN(H_\disc))} \mathcal{H}_\disc J}_{\cb} \cdot 1 \\
& \lesssim \norm{Q \ot \Id_{\L^p(\VN(H_\disc),S^p)}} \norm{\mathcal{H}_\disc}_\cb \norm{J}_\cb.
\end{align*}
Let us estimate the three factors. For the first, we recall that $Q$ is the Gaussian projection $Q \co \L^p(\Omega) \to \L^p(\Omega),\:f \mapsto \sum_k \langle f,\W(e_k)\rangle \W(e_k)$. Its norm is controlled by $C\sqrt{\max(p,p^*)}$ according to Lemma \ref{lem-K-convexity-constant} below. Furthermore, going into the proof of Proposition \ref{prop-H-boundedness}, we see that the norm $\norm{\mathcal{H}_\disc}_\cb$ is no more than the truncated Hilbert transform norm 
$$
\sup_{\epsi \in (0,1)} \bnorm{k_{\epsi,\frac{1}{\epsi}} \ast \cdot}_{\B(\L^p(\R,\L^p(\Omega,\L^p(\VN(H_\disc)),S^p)))}.
$$
The last expression in turn is controlled by an absolute constant times the norm of the Hilbert transform $\norm{H}_{\B(\L^p(\R,\L^p(\Omega,\L^p(\VN(H_\disc)),S^p)))}$ according to \cite[p.~222-223]{Haa3}. This in turn is controlled by $\frac{8}{\pi} \max(p,p^*)$ \cite[Proposition 5.4.2]{HvNVW1} and \cite[Theorem 4.3]{OsY}. Finally, $J$ is a complete contraction.
In all, we obtain
\[ 
\norm{Q \ot \Id_{\L^p(\VN(H_\disc),S^p)}} \norm{\mathcal{H}_\disc}_\cb \norm{J}_\cb 
\leq K \sqrt{\max(p,p^*)} \cdot \max(p,p^*). 
\]
The lower estimate in \eqref{equ-constants-Riesz-equivalence} follows in the same way, since in the proof of Theorem \ref{Thm-Riesz-equivalence-Schur}, it was obtained by duality.
\end{remark}

\begin{lemma}
\label{lem-K-convexity-constant}
Let $1 < p < \infty$ and $X$ be a $\sigma$-finite measure space.
We have a control of the Gaussian projection
\begin{equation}
\label{equ-K-convexity-constant}
\norm{Q \ot \Id_{\L^p(X,S^p)}}_{\L^p(\Omega,\L^p(X,S^p)) \to \L^p(\Omega,\L^p(X,S^p))} \leq C \sqrt{\max(p,p^*)} ,
\end{equation}
where $C$ is some absolute constant.
\end{lemma}

\begin{proof}
First some remarks are in order. Note that the Gaussian projection was defined as $Q(f) \ov{\mathrm{def}}{=} \sum_{k \in K} \langle f, \gamma_k \rangle \gamma_k$, where $(\gamma_k)_{k \in K}$ is an orthonormal family of standard Gaussians. If $K$ is at most countable, then the main part of the proof below applies directly. If $K$ happened to be larger, then consider for a subset $K_0$ of $K$ the modified Gaussian projection $Q_{K_0}f = \sum_{k \in K_0} \langle f, \gamma_k \rangle \gamma_k$. Note that for fixed $f \in D \ov{\mathrm{def}}{=} \L^p(\Omega,\L^p(X,S^p)) \cap (\L^2(\Omega) \ot \L^p(X,S^p))$, the sum over $K_0$ has at most countably many non-zero entries. This also implies by density of $D$ that
\begin{align*}
\MoveEqLeft
\norm{Q}_{\B(\L^p(\Omega,\L^p(X,S^p)))} 
= \sup\left\{ \frac{\norm{Q(f)}_p}{\norm{f}_p} :\: f \in D \right\} 
= \sup\left\{ \frac{\norm{Q_{K_0}(f)}_p}{\norm{f}_p} : \: f \in D,\: K_0 \subseteq K \text{ countable} \right\} \\
& = \sup\left\{ \norm{Q_{K_0}}_{\L^p(\Omega,\L^p(X,S^p)) \to \L^p(\Omega,\L^p(X,S^p))} :\: K_0 \subseteq K \text{ countable} \right\} .
\end{align*}
Then the proof below first gives an estimate for $\norm{Q_{K_0}}$, which is uniform in $K_0$, and consequently, by the above, gives an estimate for $Q_K = Q$.

Now we proceed to the proof. We recall constants of noncommutative Khintchine inequalities.
Namely, for $2 \leq p < \infty$ and a sequence of independent standard Gaussians $(\gamma_k)$ over $\Omega$, we have\footnote{\thefootnote.Indeed, according to 
\cite[p. 193]{Pis7} \cite[p. 271]{Pis15}, \eqref{equ-1-proof-K-convexity-constant} holds for $S^p$ instead of $\L^p(X,S^p)$. Then \eqref{equ-1-proof-K-convexity-constant} follows with $\L^p(X,S^p)$ by the following Fubini argument. For any finite family $(x_k)$ of $\L^p(X,S^p)$, we have
\begin{align*}
\MoveEqLeft
\norm{\sum_k \gamma_k \ot x_k}_{\L^p(\Omega,\L^p(X,S^p))}^p 
= \int_X \norm{\sum_k \gamma_k \ot x_k(t)}_{\L^p(\Omega,S^p)}^p \d t 
\leq C \sqrt{p}^p \int_X \norm{\sum_k x_k(t) \ot e_k}_{S^p(H_{\rad,p})}^p \d t \\
& = C \sqrt{p}^p \int_X \max \left\{ \norm{\left(\sum_k |x_k(t)|^2\right)^{\frac12}}_{S^p}, \norm{\left( \sum_k |x_k(t)^*|^2 \right)^{\frac12}}_{S^p} \right\}^p \d t \\
& \leq C \sqrt{p}^p \int_X \norm{\left(\sum_k |x_k(t)|^2\right)^{\frac12}}_{S^p}^p + \norm{\left( \sum_k |x_k(t)^*|^2 \right)^{\frac12}}_{S^p}^p \d t \\
& = C\sqrt{p}^p \left(\norm{\left(\sum_k |x_k|^2\right)^{\frac12}}_{\L^p(X,S^p)}^p + \norm{\left( \sum_k |x_k^*|^2 \right)^{\frac12}}_{\L^p(X,S^p)}^p \right)\\
& \leq 2C \sqrt{p}^p \max \left\{ \norm{\left(\sum_k |x_k|^2\right)^{\frac12}}_{\L^p(X,S^p)},\norm{\left( \sum_k |x_k^*|^2 \right)^{\frac12}}_{\L^p(X,S^p)} \right\}^p \\
&=2C \sqrt{p}^p \norm{\sum_k x_k \ot e_k}_{\L^p(X,S^p(H_{\rad,p}))}^p.
\end{align*} .} 
\begin{equation}
\label{equ-1-proof-K-convexity-constant}
\norm{\sum_k \gamma_k \ot x_k}_{\L^p(\Omega,\L^p(X,S^p))} \leq C \sqrt{p} \norm{\sum_k x_k \ot e_k}_{\L^p(X,S^p,H_{\rad,p})} ,
\end{equation}
with an absolute constant $C$, where the latter space was defined in Subsection \ref{Sec-Hilbertian-valued} and $(e_k)$ is any orthonormal sequence in $H$. Moreover, we claim that for $1 \leq p \leq 2$ we have
\begin{equation}
\label{equ-2-proof-K-convexity-constant}
\norm{\sum_k \gamma_k \ot x_k}_{\L^p(\Omega,\L^p(X,S^p))} \leq C \norm{\sum_k x_k \ot e_k}_{\L^p(X,S^p(H_{\rad,p}))}
\end{equation}
with some universal constant $C$. Indeed, this follows since \eqref{equ-2-proof-K-convexity-constant} holds with constant $C_1$ (resp. $C_2$) for $p = 1$ (resp. $p = 2$) according to noncommutative Khintchine inequalities. Moreover, both $\L^p(X,S^p,H_{\rad,p})$ and $\L^p(\Omega,\L^p(X,S^p))$ form a complex interpolation scale, so that we have $C_p \leq C_1^\theta C_2^{1-\theta} \leq C$ with the correct interpolation parameter $\theta \in [0,1]$.

Now the end of the argument for \eqref{equ-K-convexity-constant} is as follows.
Write $(Q \ot \Id) (f) = \sum_k \gamma_k \ot x_k$ with $x_k = \int_\Omega f \gamma_k$.
Then for $1 < p \leq 2$, by duality $(\L^p(X,S^p,H_{\rad,p}))^* \ov{\eqref{Duality-conditional}}{=} \L^{p^*}(X,S^{p^*},\ovl{H}_{\rad,p^*})$,
\begin{align*}
\MoveEqLeft
\bnorm{(Q \ot \Id) (f)}_p 
= \norm{\sum_k \gamma_k \ot x_k}_{\L^p(\Omega,\L^p(X,S^p))} 
\overset{\eqref{equ-2-proof-K-convexity-constant}}{\leq} C \norm{ \sum_k x_k \ot e_k}_{\L^p(X,S^p(H_{\rad,p}))} \\
& \overset{\eqref{Duality-conditional}}{=} C \sup\left\{ \left|\sum_k \langle x_k,y_k \rangle \right| : \: {\norm{\sum_k y_k \ot e_k}_{\L^{p^*}(X,S^{p^*}(H_{\rad,p^*}))}\leq 1} \right\}\\
& = C \sup \left\{  \left| \int_\Omega \sum_k \langle x_k,y_k \rangle \gamma_k^2 \right| : \: {\norm{\sum_k y_k \ot e_k}_{\L^{p^*}(X,S^{p^*}(H_{\rad,p^*}))}\leq 1} \right\}\\
& = C \sup \left\{   \left| \int_\Omega \big\langle(Q \ot \Id) (f), (Q \ot \Id)(y) \big\rangle \right| : \: {y \in \L^{p^*}, \norm{\sum_k \bigg(\int_\Omega y \gamma_k \bigg) \ot e_k}\leq 1} \right\} \\
& = C \sup \left\{ \left| \int_\Omega f (Q \ot \Id)(y) \right| : \: {y \in \L^{p^*}, \norm{\sum_k \bigg(\int_\Omega y \gamma_k \bigg) \ot e_k}\leq 1} \right\} \\
& \leq C \norm{f}_{p} \sup \left\{  \norm{(Q \ot \Id)(y)}_{\L^{p^*}(\Omega,\L^{p^*}(S^{p^*}))} : \: {y \in \L^{p^*}, \norm{\sum_k \bigg(\int_\Omega y \gamma_k \bigg) \ot e_k}\leq 1} \right\} \\
& \overset{\eqref{equ-1-proof-K-convexity-constant}}{\leq} C' \sqrt{p^*} \norm{f}_{\L^p(\Omega,\L^p(X,S^p))}.
\end{align*}
The case $2 \leq p < \infty$ follows by duality, since $Q$ is selfadjoint on $\L^2(\Omega)$.
\end{proof}

The $q$-Gaussian gradients equally satisfy the equivalence with $A_p^{\frac12}$.
To this end, the following Lemma \ref{Khintchine-Schur} (which is probably folklore) will be useful. We denote here the canonical conditional expectation $\E \co \L^p(\Gamma_q(H) \otvn \B(\ell^2_I)) \to S^p_I$.
We let 
\begin{equation}
\label{Def-Gaussqp}
\Gauss_{q,p}(S^p_I)
\ov{\mathrm{def}}{=} \ovl{\Span \big\{ s_q(h) \ot x : h \in H,x \in S^p_I \big\}}
\end{equation}
where the closure is taken in $\L^p(\Gamma_q(H) \otvn \B(\ell^2_I))$ (for the weak* topology if $p=\infty$ and $-1 \leq q<1$). If $(e_k)_{k \in K}$ is an orthonormal basis of the Hilbert space $H$, note that $\Gauss_{q,p}(S^p_I)$ is also the closure of the span of the $s_q(e_k) \ot x$'s with $x \in S^p_I$, closure which we denote temporarily by $G'$. Indeed, it is trivial that $G' \subseteq \Gauss_{q,p}(S^p_I)$. For the reverse inclusion, it suffices by linearity and density to show that $s_q(h) \ot x$ belongs to $G'$ for any $h \in H$ and any $x \in S^p_I$. Let $\epsi > 0$ be fixed. Then there exists some finite subset $F$ of $K$ and $\alpha_k \in \C$ for $k \in F$ such that $\norm{h-\sum_{k \in F} \alpha_k e_k}_H < \epsi$. We have $\sum_{k \in F} \alpha_k s_q(e_k) \ot x \in G'$ and
\begin{align*}
\MoveEqLeft
\norm{s_q(h) \ot x - \sum_{k \in F} \alpha_k s_q(e_k) \ot x}_{\L^p(\Gamma_q(H) \otvn \B(\ell^2_I))} 
= \norm{s_q\bigg(h - \sum_{k \in F}\alpha_k e_k\bigg) \ot x}_{\L^p(\Gamma_q(H) \otvn \B(\ell^2_I))} \\
&=\norm{s_q\bigg(h - \sum_{k \in F}\alpha_k e_k\bigg)}_{\L^p(\Gamma_q(H))} \norm{x}_{S^p_I}
\cong \norm{s_q\bigg(h - \sum_{k \in F}\alpha_k e_k\bigg)}_{\L^2(\Gamma_q(H))} \norm{x}_{S^p_I} 
< \epsi \norm{x}_{S^p_I}.
\end{align*}
We conclude $\Gauss_{q,p}(S^p_I) \subseteq G'$ by closedness of $G'$.


If $2 \leq p<\infty$, recall that we have a contractive injective inclusion 
\begin{equation}
\label{Inclusion-sympa}
\Gauss_{q,p}(S^p_I) \subseteq \L^p(\Gamma_{q}(H) \otvn \B(\ell^2_I)) \subseteq \L^p_{cr}(\E)
\end{equation}
and that if $1<p<2$ we have an inclusion 
\begin{equation}
\label{inclusion-mechante}
\Gauss_{q,2}(\C) \ot S^p_I
\subseteq \L^2(\Gamma_q(H)) \ot S^p_I
\subseteq S^p_{I}(\L^2(\Gamma_q(H))_{\rad,p})
\ov{\eqref{Egalite-fantastique}}{=}\L^p_{cr}(\E).
\end{equation}

We introduce the orthogonal projection $P \co \L^{2}(\Gamma_q(H)) \to \L^{2}(\Gamma_q(H))$ onto the closed span $\Gauss_{q,2}(\C)$ of the $s_q(e_k)$.
We let $Q_{p} \ov{\mathrm{def}}{=} \Id_{S^{p}_I} \ot P$ initially defined on $\M_{I,\fin} \ot \L^2(\Gamma_q(H))$.
We recall some properties of $Q_p$.

\begin{lemma}
\label{lem-gaussian-projection-Schur}
Let $1 < p < \infty$ and $-1 \leq q \leq 1$. Then $Q_p$ induces well-defined contractions
\begin{align}
\label{equ-1-gaussian-projection-Schur}
Q_p  \co \L^p_c(\E) \to \L^p_c(\E) 
\quad \text{and} \quad
Q_p \co \L^p_r(\E) \to \L^p_r(\E) 
\end{align}
Thus, according to Lemma \ref{lem-compatibility-Schur}, $Q_p$ also extends to a contraction $Q_p \co \L^p_{cr}(\E) \to \L^p_{cr}(\E)$.
\end{lemma}

\begin{proof}
Since $P$ is a contraction, by \eqref{extens-rad-1} and \eqref{Egalite-fantastique}, the column part of \eqref{equ-1-gaussian-projection-Schur} follows.
Then by Remark \ref{rem-problem-row}, also the row part of \eqref{equ-1-gaussian-projection-Schur} follows.
\end{proof}

Note that the noncommutative Khintchine inequalities can be rewritten under the following form.

\begin{lemma} 
\label{Khintchine-Schur}
Consider $-1 \leq q \leq 1$.
	\begin{enumerate}
\item Suppose $1<p <2$. For any element $f=\sum_{i,j,h} f_{i,j,h} s_q(h) \ot e_{ij}$ of $ \Span \big\{ s_q(h) : h \in H \} \ot \M_{I,\fin}$, we have
\begin{equation}
\label{Khintchine-Schur-p<2}
\norm{f}_{\Gauss_{q,p}(S^p_I)}  
\approx_p \norm{f}_{S^p_I(\L^2(\Gamma_q(H))_{\rad,p})}
\approx_p \inf_{f=g+h} \bigg\{\norm{\big(\E(g^*g)\big)^{\frac12}}_{S^p_I},\norm{\big(\E(h h^*) \big)^{\frac12}}_{S^p_I}\bigg\}	
\end{equation}
where the infimum can equally be taken over all $g,h \in \Span \big\{ s_q(e) : e \in H \} \ot \M_{I,\fin}$. 
		\item Suppose $2 \leq p <\infty$. For any element $f=\sum_{i,j,h} f_{i,j,h} s_q(h) \ot e_{ij}$ of $\Gauss_{q,p}(S^p_I)$ with $f_{i,j,h} \in \C$
	\begin{align}
\MoveEqLeft
 	\label{Khintchine-Schur-p>2}
\max\bigg\{\norm{\big(\E(f^* f)\big)^{\frac12}}_{S^p_I},\norm{\big(\E(f f^*) \big)^{\frac12}}_{S^p_I} \bigg\}
\leq \norm{f}_{\Gauss_{q,p}(S^p_I)}    \\
		&\lesssim \sqrt{p} \underbrace{\max\bigg\{\norm{\big(\E(f^* f)\big)^{\frac12}}_{S^p_I},\norm{\big( \E(f f^*) \big)^{\frac12}}_{S^p_I} \bigg\}}_{\norm{f}_{\L^p_{cr}(\E)}}. \nonumber
\end{align}	
\end{enumerate}

\end{lemma}

\begin{proof}
2. Suppose $2 <p <\infty$. We begin by proving the upper estimate of \eqref{Khintchine-Schur-p>2}. Fix $m \geq 1$.
We have\footnote{\thefootnote. If $h \in H$, we have
$$
\norm{\frac{1}{\sqrt{m}} \bigg(\sum_{l=1}^{m}  e_l\bigg) \ot h}_{\ell^2_m(H)}
=\frac{1}{\sqrt{m}}\norm{\sum_{l=1}^{m}  e_l}_{\ell^2_m} \norm{h}_{H}
=\norm{h}_{H}.
$$} an isometric embedding $J_m  \co \H \to \ell^2_m(H)$ defined by
\begin{equation}
\label{Def-Jm-Lemma310}
J_m(h)
\ov{\mathrm{def}}{=} \frac{1}{\sqrt{m}} \bigg(\sum_{l=1}^{m}  e_l\bigg) \ot h.
\end{equation}
Hence we can consider the associated operator $\Gamma_p(J_m) \co \L^p(\Gamma_q(H)) \to \L^p(\Gamma_q(\ell^2_m(H)))$ of second quantization \eqref{SQq} which is an isometric completely positive map.
By tensorizing with the identity $\Id_{S^p_I}$, we obtain an isometric map $\pi_I\ov{\mathrm{def}}{=}\Gamma_p(J_m) \ot \Id_{S^p_I} \co \L^p(\Gamma_q(H) \otvn \B(\ell^2_I)) \to \L^p(\Gamma_q(\ell^2_m(H)) \otvn \B(\ell^2_I))$. For any finite sum $f=\sum_{i,j,h} f_{i,j,h} s_q(h) \ot e_{ij}$ of $\L^p(\Gamma_q(H) \otvn \B(\ell^2_I))$ with $f_{i,j,h} \in \C$, we obtain
\begin{align}
\label{Equa-inter-4532}
\MoveEqLeft
 \norm{f}_{\L^p(\Gamma_q(H) \otvn \B(\ell^2_I))} 
=\norm{\sum_{i,j,h} f_{i,j,h} s_q(h) \ot e_{ij}}_{\L^p(\Gamma_q(H) \otvn \B(\ell^2_I))} \\
&=\norm{\pi_I\bigg(\sum_{i,j,h} f_{i,j,h} s_q(h) \ot e_{ij}\bigg)}_{\L^p(\Gamma_q(\ell^2_m(H)) \otvn \B(\ell^2_I)))} \nonumber\\
&=\norm{\sum_{i,j,h} f_{i,j,h} \Gamma_p(J_m)(s_q(h))  \ot e_{ij}}_{\L^p(\Gamma_q(\ell^2_m(H)) \otvn \B(\ell^2_I)))} \nonumber\\
&\ov{\eqref{SQq}}{=}\norm{\sum_{i,j,h} f_{i,j,h} s_{q,m}(J_m(h)) \ot e_{ij}}_{\L^p(\Gamma_q(\ell^2_m(H)) \otvn \B(\ell^2_I)))}\nonumber\\
&\ov{\eqref{Def-Jm-Lemma310}}{=}\norm{\sum_{i,j,h} f_{i,j,h} s_{q,m}\bigg(\frac{1}{\sqrt{m}} \bigg(\sum_{l=1}^{m}  e_l\bigg) \ot h\bigg) \ot e_{ij}}_{\L^p(\Gamma_q(\ell^2_m(H)) \otvn \B(\ell^2_I))} \nonumber\\
&=\frac{1}{\sqrt{m}} \Bgnorm{\sum_{i,j,h,l} f_{i,j,h} s_{q,m}(e_l \ot h) \ot e_{ij}}_{\L^p(\Gamma_q(\ell^2_m(H)) \otvn \B(\ell^2_I))}.  \nonumber
\end{align}
For any $1 \leq l \leq m$, consider the element $f_l \ov{\mathrm{def}}{=} \sum_{i,j,h}  f_{i,j,h} s_{q,m}(e_l \ot h) \ot e_{ij}$ and the conditional expectation $\E \co \L^p(\Gamma_q(\ell^2_m(H)) \otvn \B(\ell^2_I)) \to S^p_I$. For any $1 \leq l \leq m$, note that
\begin{align*}
\MoveEqLeft
 \E(f_l)   
		=\E\bigg(\sum_{i,j,h} f_{i,j,h} s_{q,m}(e_l \ot h) \ot e_{ij}\bigg)
		=\sum_{i,j,h}f_{i,j,h} \E\big(s_{q,m}(e_l \ot h)\ot e_{ij}\big)\\
		&=\sum_{i,j,h} f_{i,j,h} \tau( s_{q,m}(e_l \ot h)) e_{ij}
		\ov{\eqref{formule de Wick odd}}{=} 0.
\end{align*} 
We deduce that the random variables $f_l$ are mean-zero. Now, we prove in addition that the $f_l$ are  independent over $\B(\ell^2_I)$, so that we will be able to apply the noncommutative Rosenthal inequality \eqref{Inequality-Rosenthal-1-prime} to them.

\begin{lemma}
\label{lem-fl-faithfully-independent}
The $f_l$ defined here above are independent over $\B(\ell^2_I)$.
More precisely, the von Neumann algebras generated by $f_l$ are independent with respect to $\B(\ell^2 _I)$.
\end{lemma}

\begin{proof}
Due to normality of $\E$, it suffices to prove 
\begin{equation}
\label{equ-1-lem-fl-faithfully-independent-Lemma310}
\E(xy) 
=\E(x) \E(y)
\end{equation}
for $x$ (resp. $y$) belonging to a weak$^*$ dense subset of $\W^*(f_l)$ (resp. of $\W^*((f_k)_{k \neq l})$). Moreover, by bilinearity in $x,y$ of both sides of \eqref{equ-1-lem-fl-faithfully-independent-Lemma310} and selfadjointness of $s_{q,m}(e_l \ot h)$, it suffices to take $x = s_{q,m}(e_l \ot h)^n \ot e_{ij}$ for some $n \in \N$ and $y = \prod_{t = 1}^T s_{q,m}(e_{k_t} \ot h)^{n_t} \ot e_{rs}$ for some $T \in \N$, $k_t \neq l$ and $n_t \in \N$. We have $\E(x) \E(y)=\tau(s_{q,m}(e_l \ot h)^n)\tau\Big(\prod_{t = 1}^T s_{q,m}(e_{k_t} \ot h)^{n_t}\Big) \ot e_{ij} e_{rs}$ and $\E(xy) = \tau \left(s_{q,m}(e_l \ot h)^n \prod_{t = 1}^T s_{q,m}(e_{k_t} \ot h)^{n_t} \right) \ot e_{ij} e_{rs}$. We shall now apply the Wick formulae \eqref{formule de Wick} and \eqref{formule de Wick odd} to the trace term above. 

Note that if $n + \sum_{t = 1}^T n_t$ is odd, then according to the Wick formula \eqref{formule de Wick odd} we have $\E(xy) = 0$. On the other hand, then either $n$ or $\sum_{t = 1}^T n_t$ is odd, so according to the Wick formula \eqref{formule de Wick odd}, either $\E(x) = 0$ or $\E(y) = 0$. Thus, \eqref{equ-1-lem-fl-faithfully-independent-Lemma310} follows in this case. 

Now suppose that $2k \ov{\mathrm{def}}{=} n + \sum_{t = 1}^T n_t$ is even. Consider a $2$-partition $\mathcal{V} \in \mathcal{P}_2(2k)$. If both $n$ and $\sum_{t = 1}^T n_t$ are odd, then we must have some $(i,j) \in \mathcal{V}$ such that one term $\langle f_i, f_j \rangle_{\ell^2_m(H)}$ in the Wick formula \eqref{formule de Wick} equals $\langle e_l \ot h, e_{k_t} \ot h \rangle = 0$, since $k_t \neq l$. Thus, $\E(xy) = 0$, and since $n$ is odd, also $\E(x) = 0$, and therefore, \eqref{equ-1-lem-fl-faithfully-independent-Lemma310} follows. If both $n$ and $\sum_{t = 1}^T n_t$ are even, then in the Wick formula \eqref{formule de Wick}, we only need to consider those $2$-partitions $\mathcal{V}$ without a mixed term as the $\langle f_i, f_j\rangle = 0$ above. Such a $\mathcal{V}$ is clearly the disjoint union $\mathcal{V} = \mathcal{V}_1 \cup \mathcal{V}_2$ of $2$-partitions corresponding to $n$ and to $\sum_{t = 1}^T n_t$.
Moreover, we have for the number of crossings, $c(\mathcal{V}) = c(\mathcal{V}_1) + c(\mathcal{V}_2)$. With $(f_1,f_2,\ldots,f_{n+T}) = (\underbrace{e_{l} \ot h ,e_l \ot h,\ldots,e_l \ot h}_{n},e_{k_1} \ot h,\ldots, e_{k_T} \ot h)$, we obtain 
\begin{align*}
\MoveEqLeft
\tau \left(s_{q,m}(e_l \ot h)^n \prod_{t = 1}^T s_{q,m}(e_{k_t} \ot h)^{n_t} \right) 
\ov{\eqref{formule de Wick}}{=} \sum_{\mathcal{V} \in \mathcal{P}_2(2k)} q^{c(\mathcal{V})} \prod_{(i,j) \in \mathcal{V}} \langle f_i, f_j \rangle_{\ell^2_m(H)} \\
&=\sum_{\mathcal{V}_1,\mathcal{V}_2} q^{c(\mathcal{V}_1) + c(\mathcal{V}_2)} \prod_{(i_1,j_1) \in \mathcal{V}_1} \langle f_{i_1}, f_{j_1} \rangle_{\ell^2_m(H)} \prod_{(i_2,j_2) \in \mathcal{V}_2} \langle f_{i_2}, f_{j_2} \rangle_{\ell^2_m(H)} \\
&=\bigg(\sum_{\mathcal{V}_1} q^{c(\mathcal{V}_1)} \prod_{(i_1,j_1) \in \mathcal{V}_1} \langle f_{i_1}, f_{j_1} \rangle_{\ell^2_m(H)}\bigg)\bigg(\sum_{\mathcal{V}_2} q^{c(\mathcal{V}_2)} \prod_{(i_2,j_2) \in \mathcal{V}_2} \langle f_{i_2}, f_{j_2} \rangle_{\ell^2_m(H)}\bigg) \\
&\ov{\eqref{formule de Wick}}{=} \tau \left(s_{q,m}(e_l \ot h)^n \right) \tau \left( \prod_{t = 1}^T s_{q,m}(e_{k_t} \ot h)^{n_t} \right).
\end{align*}
Thus, also in this case \eqref{equ-1-lem-fl-faithfully-independent-Lemma310} follows.
\end{proof}

Now we are able to apply the noncommutative Rosenthal inequality \eqref{Inequality-Rosenthal-1-prime} which yields
\begin{align}
\label{Eq-divers-3212}
\MoveEqLeft
\norm{f}_{\L^p(\Gamma_q(H) \otvn \B(\ell^2_I))} 
\ov{\eqref{Equa-inter-4532}}{=}\frac{1}{\sqrt{m}} \norm{\sum_{l=1}^{m} f_l}_{\L^p(\Gamma_q(\ell^2_m(H)) \otvn \B(\ell^2_I))} \\
& \ov{\eqref{Inequality-Rosenthal-1-prime}}{\lesssim} \frac{1}{\sqrt{m}} \Bigg[ p\bigg( \sum_{l=1}^m \norm{f_l}_{\L^p(\Gamma_q(\ell^2_m(H)) \otvn \B(\ell^2_I))}^p \bigg)^{\frac{1}{p}} \\
&+ \sqrt{p}\norm{ \bigg(\sum_{l=1}^m \E(f_l^* f_l) \bigg)^{\frac12}}_{S^p_I} +\sqrt{p}\norm{ \bigg(\sum_{l=1}^m \E(f_l f_l^*) \bigg)^{\frac12}}_{S^p_I} \Bigg]. \nonumber
\end{align}
For any integer $1 \leq l \leq m$, note that
\begin{align*}
\MoveEqLeft
  \E(f_l^* f_l)  
		=\E\Bigg(\bigg(\sum_{i,j,h} f_{i,j,h} s_{q,m}(e_l \ot h) \ot e_{ij}\bigg)^* \bigg(\sum_{r,s,k} f_{r,s,k} s_{q,m}(e_l \ot k) \ot e_{rs}\bigg)\Bigg)\\
		&=\sum_{i,j,h,r,s,k} \ovl{f_{i,j,h}}f_{r,s,k} \bigg( \tau (s_{q,m}(e_l \ot h)s_{q,m}(e_l \ot k))\bigg) e_{ji} e_{rs} \\
		&\ov{\eqref{petit-Wick}}{=} 
		\sum_{i,j,h,s,k} \ovl{f_{i,j,h}}f_{i,s,k} \langle h,k\rangle_He_{js}.
\end{align*} 
and
\begin{align*}
\MoveEqLeft
\E(f^* f)   
=\E\Bigg(\bigg(\sum_{i,j,h} f_{i,j,h} s_q(h)  \ot e_{ij}\bigg)^* \bigg(\sum_{r,s,k} f_{r,s,k} s_q(k)  \ot e_{rs}\bigg)\Bigg)   \\
&=\sum_{i,j,h,r,s,k}\bigg( \ovl{f_{i,j,h}} f_{r,s,k}  \tau\big( s_q(h)s_q(k) \big)\bigg) e_{ji}e_{rs} 
\ov{\eqref{petit-Wick}}{=} \sum_{i,j,h,s,k} \ovl{f_{i,j,h}} f_{i,s,k}\langle h,k\rangle_H e_{js}.
\end{align*}
and similarly for the row terms. We conclude that 
\begin{equation}
\label{Eq-divers-98089}
\E(f_l^* f_l)
=\E(f^* f)
\quad \text{and} \quad 
\E(f_l f_l^*) 
= \E(f f^*).
\end{equation}
Moreover, for $1 \leq l \leq m$, using the isometric map $\psi_l \co H \to \ell^2_m(H)$, $h \to e_l \ot h$, we can introduce the second quantization operator $\Gamma_p(\psi_l) \co \L^p(\Gamma_q(\ell^2_m(H))) \to \L^p(\Gamma_q(\ell^2_m(H)))$. We have
\begin{align}
\label{fl-f}
\MoveEqLeft
  \norm{f_l}_{\L^p(\Gamma_q(\ell^2_m(H)) \otvn \B(\ell^2_I))}  
		=\norm{\sum_{i,j,h} f_{i,j,h} s_{q,m}(e_l \ot h) \ot e_{ij}}_{\L^p(\Gamma_q(\ell^2_m(H)) \otvn \B(\ell^2_I))} \\
		&\ov{\eqref{SQq}}{=} \norm{\sum_{i,j,h} f_{i,j,h} \Gamma_p(\psi_l)(s_{q,m}(h)) \ot e_{ij}}_{\L^p(\Gamma_q(\ell^2_m(H)) \otvn \B(\ell^2_I))} \nonumber\\
	&=\norm{\sum_{i,j,h} f_{i,j,h} s_q(h) \ot e_{ij}}_{\L^p(\Gamma_q(H) \otvn \B(\ell^2_I))}  \nonumber\\
	&=\norm{f}_{\L^p(\Gamma_q(H) \otvn \B(\ell^2_I))}. \nonumber
\end{align} 
 We infer that
\begin{align*}
\MoveEqLeft
\norm{f}_{\L^p(\Gamma_q(H) \otvn \B(\ell^2_I))} 
\ov{\eqref{Eq-divers-3212}\eqref{fl-f}\eqref{Eq-divers-98089}}{\lesssim} \frac{1}{\sqrt{m}} \Bigg[ p\bigg(\sum_{l=1}^m \norm{f}_{\L^p(\Gamma_q(H) \otvn \B(\ell^2_I))}^p \bigg)^{\frac{1}{p}} \Bigg.\\
&\qquad \qquad \Bigg. + \sqrt{p}\norm{ \bigg(\sum_{l=1}^m \E(f^* f) \bigg)^{\frac12}}_{S^p_I}+\sqrt{p}\norm{ \bigg(\sum_{l=1}^m \E(f f^*) \bigg)^{\frac12}}_{S^p_I} \Bigg]\\
&=\frac{1}{\sqrt{m}} \Bigg[p m^{\frac{1}{p}} \norm{f}_{\L^p(\Gamma_q(H) \otvn \B(\ell^2_I))}+ \sqrt{pm}\norm{\big(\E(f^* f)\big)^{\frac12}}_{S^p_I} +\sqrt{pm}\norm{\big( \E(f f^*) \big)^{\frac12}}_{S^p_I} \Bigg]\\
&=p m^{\frac{1}{p}-\frac12} \norm{f}_{\L^p(\Gamma_q(H) \otvn \B(\ell^2_I))}+\sqrt{p}\norm{\big(\E(f^* f)\big)^{\frac12}}_{S^p_I}+\sqrt{p}\norm{\big( \E(f f^*) \big)^{\frac12}}_{S^p_I}.
\end{align*}
Since $p>2$, passing to the limit when $m \to \infty$, we finally obtain 
$$
\norm{f}_{\L^p(\Gamma_q(H) \otvn \B(\ell^2_I))}  
\lesssim \sqrt{p}\bigg[\norm{\big(\E(f^* f)\big)^{\frac12}}_{S^p_I}+\norm{\big(\E(f f^*) \big)^{\frac12}}_{S^p_I} \bigg].
$$
Using the equivalence $\ell^1_2 \approx \ell^\infty_2$, we obtain the upper estimate of \eqref{Khintchine-Schur-p>2}.

The lower estimate of \eqref{Khintchine-Schur-p>2} holds with constant $1$ from the contractivity of the conditional expectation $\E$ on $\L^{\frac{p}{2}}(\Gamma_q(H) \otvn \B(\ell^2_I))$: 
\begin{align*}
\MoveEqLeft
\max\bigg\{\norm{\big(\E(f^* f)\big)^{\frac12}}_{S^p_I},\norm{\big(\E(f f^*) \big)^{\frac12}}_{S^p_I} \bigg\}    
		= \max\bigg\{\bnorm{\E(f^* f)}_{S^{\frac{p}{2}}_I}^{\frac12},\bnorm{\E(f f^*)}_{S^{\frac{p}{2}}_I}^{\frac12} \bigg\} \\
		&\leq \max\Bigg\{\norm{f^* f}_{\L^{\frac{p}{2}}(\Gamma_q(H) \otvn \B(\ell^2_I))}^{\frac12},\norm{f f^*}_{\L^{\frac{p}{2}}(\Gamma_q(H) \otvn \B(\ell^2_I))}^{\frac12} \Bigg\}
		=\norm{f}_{\L^p(\Gamma_q(H) \otvn \B(\ell^2_I))}.
\end{align*}
1. Let us now consider the case $1 < p < 2$. We will proceed by duality as follows. 
Recall the Gaussian projection $Q_p$ from Lemma \ref{lem-gaussian-projection-Schur}.

Using in the second equality that $Q_p^* = Q_{p^*}$ and that $Q_{p^*}$ extends to a contraction on $S^{p^*}_I(\L^2(\Gamma_q(H))_{\rad,p^*})$ according to Lemma \ref{lem-gaussian-projection-Schur}, and using the upper estimate of \eqref{Khintchine-Schur-p>2} and the density of $\Span \{ s_q(h):\: h \in H \} \ot \M_{I,\fin}$ in $S^{p^*}_I(\L^2(\Gamma_q(H))_{\rad,p^*})$ in the last inequality, we obtain for any $f =Q_p(f)\in \Gauss_{q,2}(\C) \ot S^p_I$
\begin{align*}
\MoveEqLeft
\norm{f}_{S^{p}_I(\L^2(\Gamma_q(H))_{\rad,p})} 
\ov{\eqref{Duality-conditional}}{=} \sup_{\norm{g}_{S^{p^*}_I(\L^2(\Gamma_q(H))_{\rad,p^*})} \leq 1} \big| \langle f,g \rangle_{} \big|
=\sup_{\stackrel{\norm{g}_{S^{p^*}_I(\L^2(\Gamma_q(H))_{\rad,p^*})} \leq 1 }{g \in \Ran Q_{p^*}}} \big| \langle f,g \rangle_{} \big|\\ 
&\leq \norm{f}_{\L^p(\Gamma_q(H) \otvn \B(\ell^2_I))} \sup_{\stackrel{\norm{g}_{S^{p^*}_I(\L^2(\Gamma_q(H)_{\rad,p^*})} \leq 1}{g \in \Ran Q_{p^*}}} \norm{g}_{\L^{p^*}(\Gamma_q(H) \otvn \B(\ell^2_I))} 
\ov{\eqref{Khintchine-Schur-p>2}}{\lesssim_p} \norm{f}_{\L^p(\Gamma_q(H) \otvn \B(\ell^2_I))}.
\end{align*} 

Note that in the definition of the norm of $S^p_I(\L^2(\Gamma_q(H))_{\rad,p})$ above, the infimum runs over $g,h$ belonging to $\L^p_c(\E)$ and $\L^p_r(\E)$.
Our next goal is to restrict to those $g,h$ belonging to $\Span\{ s_q(e) : \: e \in H\} \ot \M_{I,\fin}$.
To this end, consider a decomposition $f=g+h$ with $g \in \L^p_{c}(\E)$ and $h \in \L^p_{r}(\E)$ such that
$$
\norm{g}_{\L^p_{c}(\E)} + \norm{h}_{\L^p_{r}(\E)}
\leq 2\norm{f}_{S^{p}_I(\L^2(\Gamma_q(H))_{\rad,p})}.
$$ 
Then for some large enough $J$ and with as before $Q_p$ the Gaussian projection, we have 
\[ 
f 
= Q_p(\Id_{\Gamma_q(H)} \ot \Tron_J)(f) 
= Q_p(\Id_{\Gamma_q(H)} \ot \Tron_J)(g) + Q_p(\Id_{\Gamma_q(H)} \ot \Tron_J)(h) 
\]
and $Q_p(\Id_{\Gamma_q(H)} \ot \Tron_J)(g)$, $Q_p(\Id_{\Gamma_q(H)} \ot \Tron_J)(h)$ belong to $\Span \big\{ s_q(e) : e \in H \} \ot \M_{I,\fin}$. We claim that we have
$$
\bnorm{Q_p(\Id_{\Gamma_q(H)} \ot \Tron_J)(g)}_{\L^p_{c}(\E)}
\leq \norm{g}_{\L^p_{c}(\E)}
\quad \text{and} \quad
\bnorm{Q_p(\Id_{\Gamma_q(H)} \ot \Tron_J)(h)}_{\L^p_{r}(\E)}
\leq \norm{h}_{\L^p_{r}(\E)}.
$$
Indeed, first note that $Q_p \co \L^p_c(\E) \to \L^p_c(\E)$ and $Q_p \co \L^p_r(\E) \to \L^p_r(\E)$) are contractive according to \eqref{equ-1-gaussian-projection-Schur}. Moreover, since $\Tron_J \co S^p_I \to S^p_I$ is a complete contraction, according to \eqref{extens-col-1}, the linear map $\Tron_J \ot \Id_{\L^2(\Gamma_q(H))} \co S^{p}_I(\L^2(\Gamma_q(H))_{c,p}) \to S^{p}_I(\L^2(\Gamma_q(H))_{c,p})$ is also a contraction. So we have
\begin{align*}
\MoveEqLeft
\bnorm{Q_p(\Id_{\Gamma_q(H)} \ot \Tron_J)(g)}_{\L^p_{c}(\E)}          
\ov{\eqref{Egalite-fantastique}}{=} \norm{Q_p(\Tron_J \ot \Id_{\L^2(\Gamma_q(H))})(g)}_{S^p_I(\L^2(\Gamma_q(H))_{c,p})} \\
&\ov{\eqref{extens-col-1}}{\leq} \norm{g}_{S^p_I(\L^2(\Gamma_q(H))_{c,p})}
\ov{\eqref{Egalite-fantastique}}{=}\norm{g}_{\L^p_{c}(\E)}.
\end{align*} 
The row estimate is similar.
Together we have shown that the third expression in \eqref{Khintchine-Schur-p<2} is controlled by $\norm{f}_{\Gauss_{q,p}(S^p_I)}$.

Now, we will prove the remaining estimate, that is, $\norm{f}_{\Gauss_{q,p}(S^p_I)}$ is controlled by the second expression in \eqref{Khintchine-Schur-p<2}. Since $1<p<2$, the function $\R^+ \to \R^+$, $t \mapsto t^{\frac{p}{2}}$ is operator concave by \cite[p.~112]{Bha1}. Using \cite[Corollary 2.2]{HaP1} applied with the trace preserving positive map $\E$, we can write
\begin{align*}
\MoveEqLeft
\norm{f}_{\L^p(\Gamma_q(H) \otvn \B(\ell^2_I))}
=\norm{|f|^2}_{\L^{\frac{p}{2}}(\Gamma_q(H) \otvn \B(\ell^2_I))}^{\frac{1}{2}}
\leq \norm{\E(|f|^2)}_{\L^{\frac{p}{2}}(\Gamma_q(H) \otvn \B(\ell^2_I))}^{\frac{1}{2}}
\ov{}{=} \norm{f}_{\L^p_c(\E)}
\end{align*}
and similarly for the row term.
Thus, passing to the infimum over all decompositions $f = g+h$, we obtain $\norm{f}_{\L^p(\Gamma_q(H)\otvn \B(\ell^2_I))} \leq \norm{f}_{S^p_I(\L^2(\Gamma_q(H))_{\rad,p})}$, which can be majorised in turn by the infimum of $\norm{g}_{\L^p_c(\E)} + \norm{h}_{\L^p_r(\E)}$, where $f = g + h$ and $g,h \in \Span\{ s_q(e) :\: e \in H\} \ot \M_{I,\fin}$.
Hence, we have the last equivalence in the part 1 of the theorem. 

The case $p=2$ is obvious since we have isometrically $\L^2_{cr}(\E)=S_I^2(\L^2(\Gamma_q(H))_{\rad,2})=\L^2(\Gamma_q(H)) \ot_2S^2_I$ by Lemma \ref{Lema-egalite-fantas} and \cite[Remark 2.3 (1)]{JMX}.  
\end{proof}

We can equally extend Lemma \ref{Khintchine-Schur} to the case that $f \in \L^p(\Gamma_q(H) \otvn \B(\ell^2_I))$ belongs to $\Ran Q_p$, where we recall that $Q_p = \Id_{S^p_I} \ot P \co \L^p(\Gamma_q(H) \otvn \B(\ell^2_I)) \to \L^p(\Gamma_q(H) \otvn \B(\ell^2_I))$ is the Gaussian projection from Lemma \ref{lem-gaussian-projection-Schur}.

\begin{lemma}
\label{lem-Qp-bounded-Schur}
Let $-1 \leq q \leq 1$. Suppose $1 < p < \infty$. Then the map $Q_p \co \L^p(\Gamma_q(H) \otvn \B(\ell^2_I)) \to \L^p(\Gamma_q(H) \otvn \B(\ell^2_I))$ is completely bounded.
\end{lemma}

\begin{proof}
Since $P$ is selfadjoint, we have $Q_p^* = Q_{p^*}$.
We obtain for $f \in \Gamma_q(H) \otvn \M_{I,\fin}$,
\begin{align*}
\MoveEqLeft
\norm{Q_p(f) }_{\Gauss_{q,p}(S^p_I)} \ov{\text{Lemma }\ref{Khintchine-Schur}}{\lesssim} \norm{ Q_p(f) }_{\L^p_{cr}(\E)} 
= \sup \left\{ \left|\langle Q_p(f), g \rangle \right| : \: \norm{g}_{\L^{p^*}_{cr}(\E)} \leq 1 \right\} \\
& = \sup \left\{ \left| \langle f, Q_{p^*}(g) \rangle_{\L^p_{cr}(\E),\L^{p^*}_{cr}(\E)} \right| : \: g \right\} \\
& = \sup \left\{ \left| \langle f, Q_{p^*}(g) \rangle_{\L^p(\Gamma_q(H) \otvn \B(\ell^2_I)),\L^{p^*}} \right| : \: g \right\}   \\
& \leq \norm{f}_{\L^p(\Gamma_q(H) \otvn \B(\ell^2_I))} \sup \left\{ \norm{ Q_{p^*}(g) }_{\L^{p^*}(\Gamma_q(H) \otvn \B(\ell^2_I))} : \: g \right\} \\
& \ov{\text{Lemma }\ref{Khintchine-Schur}}{\lesssim} \norm{f}_{\L^p(\Gamma_q(H) \otvn \B(\ell^2_I))} \sup \left\{ \norm{ Q_{p^*}(g) }_{\L^{p^*}_{cr}(\E)} : \: \norm{g}_{\L^{p^*}_{cr}(\E)} \leq 1 \right\} \\
& \leq \norm{f}_{\L^p(\Gamma_q(H) \otvn \B(\ell^2_I))} ,
\end{align*}
where in the last step we used that $Q_{p^*}$ is contractive on $\L^{p^*}_{cr}(\E)$ according to Lemma \ref{lem-gaussian-projection-Schur}.
This shows that $Q_p$ is bounded.
Then the fact that $Q_p$ is completely bounded follows from a standard matrix amplification argument since we can replace $I$ by $I \times \{1, \ldots, N\}$.
\end{proof}

\begin{lemma} 
\label{Khintchine-Schur-full-space}
Consider $-1 \leq q \leq 1$ and let $f \in \Ran Q_p = \Gauss_{q,p}(S^p_I)$.
	\begin{enumerate}
\item Suppose $1<p <2$. We have
\begin{equation}
\label{Khintchine-Schur-p<2-full-space}
\norm{f}_{\Gauss_{q,p}(S^p_I)}  
\approx_p \norm{f}_{S^p_I(\L^2(\Gamma_q(H))_{\rad,p})}
\approx_p \inf_{f=g+h} \bigg\{\norm{\big(\E(g^*g)\big)^{\frac12}}_{S^p_I},\norm{\big(\E(h h^*) \big)^{\frac12}}_{S^p_I}\bigg\}	
\end{equation}
where the infimum can be taken over all $g,h \in \Ran Q_p$.
		\item Suppose $2 \leq p <\infty$. We have
	\begin{align}
\MoveEqLeft
 	\label{Khintchine-Schur-p>2-full-space}
\max\bigg\{\norm{\big(\E(f^* f)\big)^{\frac12}}_{S^p_I},\norm{\big(\E(f f^*) \big)^{\frac12}}_{S^p_I} \bigg\}
\leq \norm{f}_{\Gauss_{q,p}(S^p_I)}    \\
		&\lesssim \sqrt{p} \max\bigg\{\norm{\big(\E(f^* f)\big)^{\frac12}}_{S^p_I},\norm{\big( \E(f f^*) \big)^{\frac12}}_{S^p_I} \bigg\}. \nonumber
\end{align}	
\end{enumerate}
\end{lemma}

\begin{proof}
Observe that $\Id_{\Gamma_q(H)} \ot \Tron_J$ approximates the identity on $\Gamma_q(H) \otvn \M_{\I,\fin}$. Moreover, according to Lemma \ref{lem-block-Schur}, the net $(\Id_{\Gamma_q(H)} \ot \Tron_J)$ is bounded in $\B(\L^p(\Gamma_q(H) \otvn \B(\ell^2_I)))$. We conclude by density of $\Gamma_q(H) \otvn \M_{I,\fin}$ in $\L^p(\Gamma_q(H) \otvn \B(\ell^2_I))$ that the net $(\Id_{\Gamma_q(H)} \ot \Tron_J)$ converges in the point norm topology of $\L^p(\Gamma_q(H) \otvn \B(\ell^2_I))$ to the identity. Moreover, replacing in this argument boundedness in $\B(\L^p(\Gamma_q(H) \otvn \B(\ell^2_I)))$ by boundedness in $\B(\L^p_c(\E))$ (resp. $\B(\L^p_r(\E)),\: \B(\L^p_{cr}(\E))$) according to \eqref{extens-col-1} and \eqref{Egalite-fantastique} (resp. Remark \ref{rem-problem-row}), we obtain that the net $(\Id_{\Gamma_q(H)} \ot \Tron_J)$ converges in the point norm topology of $\L^p_c(\E)$ (resp. $\L^p_r(\E),\: \L^p_{cr}(\E)$) to the identity.
Note that for fixed $J$, $(\Id_{\Gamma_q(H)} \ot \Tron_J)(f) = (\Id_{\Gamma_q(H)} \ot \Tron_J)(Q_pf) = Q_p(\Id_{\Gamma_q(H)} \ot \Tron_J)(f)$ belongs to $\Span \{ s_q(e) : \: e \in H \} \ot \M_{\I,\fin}$. Thus, Lemma \ref{Khintchine-Schur} applies to $f$ replaced by $(\Id_{\Gamma_q(H)} \ot \Tron_J)(f)$ and therefore,
\begin{align*}
\MoveEqLeft
\norm{f}_{\Gauss_{q,p}(S^p_I)} 
= \lim_J \norm{(\Id_{\Gamma_q(H)} \ot \Tron_J)(f)}_{\Gauss_{q,p}(S^p_I)} \\
& \ov{\text{Lemma }\ref{Khintchine-Schur}}{\cong} \lim_J \norm{(\Id_{\Gamma_q(H)} \ot \Tron_J)(f)}_{\L^p_{cr}(\E)} 
= \norm{f}_{\L^p_{cr}(\E)}.
\end{align*}
Finally, the fact that one can restrict the infimum to all $g,h \in \Ran Q_p$ can be proved in the same way as that in Lemma \ref{Khintchine-Schur}.
\end{proof}

Now we can state the Kato square root problem for the case of the gradient taking values in a $q$-deformed algebra.

\begin{prop}
\label{Prop-equiv-q-gaussians}
Suppose $-1 \leq q \leq 1$ and $1 < p < \infty$.
For any $x \in \M_{I,\fin}$, we have
\begin{equation} 
\label{equiv-free}
\bnorm{A_p^{\frac12}(x)}_{S^p_I} 
\approx_{p} \bnorm{\partial_{\alpha,q}(x)}_{\L^p(\Gamma_q(H) \otvn \B(\ell^2_I))}. 
\end{equation} 
\end{prop}

\begin{proof}
Note that if $q = 1$, then this is a consequence of Theorem \ref{Thm-Riesz-equivalence-Schur}.
We show  the remaining case $-1 \leq q < 1$ by reducing it to the case $q = 1$. Let $x = \sum_{i,j} x_{ij} e_{ij} \in \M_{I,\fin}$. Choose an orthonormal basis $(e_k)$ of $\Span \{ \alpha_i \}$, where we consider only those indices $i$ appearing in the above double sum describing $x$. Thus, $\alpha_i-\alpha_j = \sum_k  \langle\alpha_i-\alpha_j, e_k \rangle_H e_k$, where the sum is finite. Then
\begin{align*}
\MoveEqLeft
\bnorm{A_p^{\frac12}(x)}_{S^p_I}
\ov{\eqref{OneLineDimFree}}{\approx_{p}} \norm{\partial_{\alpha,1}(x)}_{\L^p(\Omega,S^p_I)} 
=\norm{\partial_{\alpha,1}(x)}_{\Gauss_{1,p}(S^p_I)} \\
& \ov{\eqref{def-delta-alpha}}{=}\norm{\sum_{i,j} x_{ij} \W(\alpha_i-\alpha_j) \ot e_{ij}}_{\Gauss_{1,p}(S^p_I)}    \\
&=\norm{\sum_{i,j} x_{ij} \W\bigg(\sum_k \langle e_k,\alpha_i-\alpha_j \rangle_H e_k\bigg) \ot e_{ij}}_{\Gauss_{1,p}(S^p_I)} \\
&=\norm{\sum_{i,j} x_{ij}\sum_k  \langle e_k, \alpha_i-\alpha_j \rangle_H \W(e_k) \ot e_{ij}}_{\Gauss_{1,p}(S^p_I)}\\
&=\norm{\sum_k \W(e_k) \ot \bigg(\sum_{i,j} x_{ij}  \langle e_k,\alpha_i-\alpha_j\rangle_H e_{ij}\bigg)}_{\Gauss_{1,p}(S^p_I)}\\
& \overset{\cite[(2.21), (2.22) \mathrm{p}.~12]{JMX}}{\approx}\norm{\sum_k \bigg(\sum_{i,j} x_{ij}  \langle e_k, \alpha_i-\alpha_j \rangle_H e_{ij}\bigg) \ot e_k}_{S^p_I(H_{rc})}\\
& \overset{\text{Lemma} \, \ref{Khintchine-Schur}} \approx\norm{\sum_k s_q(e_k) \ot \bigg(\sum_{i,j} x_{ij}  \langle e_k,\alpha_i-\alpha_j\rangle_H e_{ij}\bigg)}_{\Gauss_{q,p}(S^p_I)}\\
&=\norm{\sum_{i,j} x_{ij} s_q(\alpha_i-\alpha_j) \ot e_{ij}}_{\Gauss_{q,p}(S^p_I)}
\ov{\eqref{def-delta-alpha}}{=} \norm{\partial_{\alpha,q}(x)}_{\Gauss_{q,p}(S^p_I)}.
\end{align*}
\end{proof}

\begin{remark} \normalfont
\label{rem-constants-Riesz-equivalence-q}
Assume $-1 \leq q \leq 1$.
In Proposition \ref{Prop-equiv-q-gaussians} above, we obtain again constants depending on $p$ as in \eqref{equ-constants-Riesz-equivalence} of a slightly different form, that is, for some absolute constant $K > 0$, we have for all $x \in \M_{I,\fin}$,
\begin{align}
\frac{1}{K p^2} \bnorm{\partial_{\alpha,q}(x)}_{\L^p(\Gamma_q(H) \otvn \B(\ell^2_I))} 
\leq \bnorm{A^{\frac12}_p(x)}_{S^p_I} \leq K p^{2} \bnorm{\partial_{\alpha,q}(x)}_{\L^p(\Gamma_q(H) \otvn \B(\ell^2_I))} & \quad (2 \leq p < \infty) \label{equ-1-constants-Riesz-equivalence-q}\\
\intertext{and}
\frac{1}{K (p^*)^2} \bnorm{\partial_{\alpha,q}(x)}_{\L^p(\Gamma_q(H) \otvn \B(\ell^2_I))} 
\leq \bnorm{A^{\frac12}_p(x)}_{S^p_I} \leq K (p^*)^2 \bnorm{\partial_{\alpha,q}(x)}_{\L^p(\Gamma_q(H) \otvn \B(\ell^2_I))} & \quad (1 < p \leq 2). \label{equ-2-constants-Riesz-equivalence-q}
\end{align}

Let us prove this and start with the case $2 \leq p < \infty$.
Note that $\W(e_k)$ and $s_q(e_k)$ are orthonormal systems in their spaces $\L^2(\Omega)$ and $\L^2(\Gamma_q(H))$ respectively.
Thus, we have isometrically, for $x_k \in \M_{I,\fin}$,
\begin{equation}
\label{equ-2-rem-constants-Riesz-equivalence-q}
\norm{ \sum_k x_k \ot \W(e_k)}_{S^p_I(\L^2(\Omega)_{\rad,p})} = \norm{\sum_k x_k \ot s_q(e_k)}_{S^p_I(\L^2(\Gamma_q(H))_{\rad,p})} .
\end{equation}
Thus, for the upper estimate in \eqref{equ-1-constants-Riesz-equivalence-q}, we have
\begin{align*}
\MoveEqLeft
\bnorm{A^{\frac12}_p(x)}_{S^p_I} 
\ov{\eqref{equ-constants-Riesz-equivalence}}{\leq}  K p^{\frac32} \bnorm{\partial_{\alpha,1}(x)}_{\Gauss_{1,p}(S^p_I)} \\
& = K p^{\frac32} \Bnorm{\sum_k \W(e_k) \ot \bigg( \sum_{i,j} x_{ij} \langle e_k, \alpha_i - \alpha_j \rangle_H e_{ij} \bigg) }_{\Gauss_{1,p}(S^p_I)} \\
& \overset{\text{Lemma }\ref{Khintchine-Schur}}{\leq} K' p^{\frac32} \cdot p^{\frac12} \Bnorm{ \sum_k \bigg( \sum_{i,j} x_{ij} \langle e_k, \alpha_i - \alpha_j \rangle_H e_{ij} \bigg) \ot \W(e_k) }_{S^p_I(\L^2(\Omega)_{\rad,p})} \\
& \overset{\eqref{equ-2-rem-constants-Riesz-equivalence-q}}{=}
K' p^{\frac32} \cdot p^{\frac12} \Bnorm{ \sum_k \bigg( \sum_{i,j} x_{ij} \langle e_k, \alpha_i - \alpha_j \rangle_H e_{ij} \bigg) \ot s_q(e_k) }_{S^p_I(\L^2(\Gamma_q(H))_{\rad,p})} \\
& \overset{\text{Lemma }\ref{Khintchine-Schur}}{\leq} K' p^2 \Bnorm{ \sum_k s_q(e_k) \ot \bigg( \sum_{i,j} x_{ij} \langle e_k, \alpha_i - \alpha_j \rangle_H e_{ij} \bigg) }_{\Gauss_{q,p}(S^p_I)} \\
& =K' p^2 \bnorm{\partial_{\alpha,q}(x)}_{\Gauss_{q,p}(S^p_I)}.
\end{align*}
In the other direction, we have, again for $2 \leq p < \infty$,
\begin{align*}
\MoveEqLeft
\bnorm{A^{\frac12}_p(x)}_{S^p_I} \ov{\eqref{equ-constants-Riesz-equivalence}}{\geq} \frac{1}{K p^{\frac32}} \bnorm{\partial_{\alpha,1}(x)}_{\Gauss_{1,p}(S^p_I)} \\
& = \frac{1}{K p^{\frac32}} \Bnorm{\sum_k \W(e_k) \ot \bigg( \sum_{i,j} x_{ij} \langle e_k, \alpha_i - \alpha_j \rangle_H e_{ij} \bigg) }_{\Gauss_{1,p}(S^p_I)} \\
& \overset{\text{Lemma }\ref{Khintchine-Schur}}{\geq} \frac{1}{K p^{\frac32}} \Bnorm{ \sum_k \bigg( \sum_{i,j} x_{ij} \langle e_k, \alpha_i - \alpha_j \rangle_H e_{ij} \bigg) \ot \W(e_k) }_{S^p_I(\L^2(\Omega)_{\rad,p})} \\
& \overset{\eqref{equ-2-rem-constants-Riesz-equivalence-q}}{=}
\frac{1}{K p^{\frac32}} \Bnorm{ \sum_k \bigg( \sum_{i,j} x_{ij} \langle e_k, \alpha_i - \alpha_j \rangle_H e_{ij} \bigg) \ot s_q(e_k) }_{S^p_I(\L^2(\Gamma_q(H))_{\rad,p})} \\
& \overset{\text{Lemma }\ref{Khintchine-Schur}}{\geq} \frac{1}{K' p^{\frac32} \cdot p^{\frac12}} \Bnorm{ \sum_k s_q(e_k) \ot \bigg( \sum_{i,j} x_{ij} \langle e_k, \alpha_i - \alpha_j \rangle_H e_{ij} \bigg) }_{\Gauss_{q,p}(S^p_I)} \\
& = \frac{1}{K' p^2} \bnorm{\partial_{\alpha,q}(x)}_{\Gauss_{q,p}(S^p_I)}.
\end{align*}
We turn to the case $1 < p \leq 2$. Note that \eqref{equ-2-rem-constants-Riesz-equivalence-q} equally holds in this case.
Then we argue in the same way as in the case $p \geq 2$, noting that the additional factor ${p^*}^{\frac12}$ appears at another step of the estimate, in accordance with Lemma \ref{Khintchine-Schur} parts 1. and 2.
\end{remark}

\begin{remark} \normalfont
\label{rem-forthcoming-Riesz-norm-Schur}
In a forthcoming paper, we shall improve the dependence on $p$ of the constants in the above Remark \ref{rem-constants-Riesz-equivalence-q} with a different method.
\end{remark}

In the remainder of this subsection, we shall extend the gradient $\partial_{\alpha,q}$ to a closed (weak* closed if $p = \infty$) operator $S^p_I \to \L^p(\Gamma_q(H) \otvn \B(\ell^2_I))$, and in particular identify its domain in terms of the generator of the markovian semigroup. All this will be achieved in Proposition \ref{Prop-derivation-closable}. 

We refer to \cite{Kat1} for information on closable and closed operators acting on Banach spaces. By \cite[p.~165]{Kat1}, an operator $T \co \dom T \subseteq X \to Y$ is closed if and only if its domain $\dom T$ is a complete space with respect to the graph norm 
\begin{equation}
\label{Def-graph-norm}
\norm{x}_{\dom T}
\ov{\mathrm{def}}{=} \norm{x}_X+\norm{T(x)}_{Y}.	
\end{equation}
For $t,s>0$, we let $f_{t}(s)\ov{\mathrm{def}}{=}\frac{t}{\sqrt{4\pi s^3}} \e^{-\frac{t^2}{4s}}$. It is well-known \cite[Lemma 1.6.7]{ABHN1} that
\begin{equation}
\label{integrale-utile-1}
\int_{0}^{+\infty} \e^{-s \lambda} f_{t}(s) \d s
=\e^{-t \sqrt{\lambda}}.
\end{equation}
Let $B$ be the negative generator of a semigroup $(\e^{-tB})_{t \geq 0}$ acting on a Banach space $X$. By \cite[Example 3.4.6]{Haa1}, for any $x \in X$ we have
\begin{equation}
	\label{subordinated-formula}
	\e^{-t B^\frac{1}{2}}x
	=\int_{0}^{+\infty} f_t(s) \e^{-sB}(x) \d s.
\end{equation}

We define the densely defined unbounded operator 
$$
\partial_{\alpha,q}^* \co \L^p(\Gamma_q(H)) \ot \M_{I,\fin} \subseteq \L^p(\Gamma_q(H)) \otvn \B(\ell^2_I)) \to S^p_I
$$ 
by
\begin{equation}
\label{Adjoint-partial-Schur}
\partial_{\alpha,q}^*(f \ot e_{ij})
=\big\langle s_q(\alpha_i-\alpha_j),f\big\rangle_{\L^{p^*},\L^p} e_{ij} =\tau(s_q(\alpha_i-\alpha_j)f)e_{ij}, \quad i,j \in I, f \in \L^p(\Gamma_q(H)).
\end{equation}

\begin{lemma}
\label{Lemma-adjoint-prtial-Schur}
The operators $\partial_{\alpha,q}$ and $\partial_{\alpha,q}^*$ are formal adjoints (with respect to duality brackets $\langle x ,y \rangle = \tau(x^*y)$).
\end{lemma}

\begin{proof}
For any $i,j,k,l \in I$ and any $f \in \L^p(\Gamma_q(H))$, we have 
\begin{align}
\MoveEqLeft
\label{Calcul-100000}
\big\langle \partial_{\alpha,q}(e_{kl}), f \ot e_{ij} \big\rangle_{\L^{p^*}(\Gamma_q(H) \otvn \B(\ell^2_I)),\L^p(\Gamma_q(H) \otvn \B(\ell^2_I))} \\
&\ov{\eqref{def-delta-alpha}}{=}\big\langle s_q(\alpha_k-\alpha_l) \ot e_{kl}, f \ot e_{ij}\big\rangle_{\L^{p^*}(\Gamma_q(H) \otvn \B(\ell^2_I)),\L^p(\Gamma_q(H) \otvn \B(\ell^2_I))}  \nonumber \\
&=\delta_{ki} \delta_{lj} \big\langle  s_q(\alpha_i-\alpha_j),f\big\rangle_{\L^{p^*}(\Gamma_q(H)),\L^{p}(\Gamma_q(H))} \nonumber \\
&=\Big\langle e_{kl}, \big\langle  s_q(\alpha_i-\alpha_j),f\big\rangle_{\L^{p^*}(\Gamma_q(H)),\L^{p}(\Gamma_q(H))} e_{ij}\Big\rangle_{S^{p^*}_I,S^{p}_I} \nonumber=\Big\langle e_{kl}, \partial_{\alpha,q}^*(f \ot e_{ij})\Big\rangle_{S^{p^*}_I,S^{p}_I}. \nonumber
\end{align}
\end{proof}

\begin{prop}
\label{Prop-derivation-closable} 
Suppose $1 < p<\infty$ and $-1 \leq q \leq 1$.
\label{prop-fermable}
\begin{enumerate}
	\item The operator $\partial_{\alpha,q} \co \M_{I,\fin} \subseteq S^p_I \to \L^p(\Gamma_q(H) \otvn \B(\ell^2_I))$ is closable as a densely defined operator on $S^p_I$ into $\L^p(\Gamma_q(H) \otvn \B(\ell^2_I))$. We denote by $\partial_{\alpha,q,p}$ its closure. So $\M_{I,\fin}$ is a core of $\partial_{\alpha,q,p}$.
	
	\item $(\Id_{S^p_I} +  A_p)^{\frac{1}{2}}(\M_{I,\fin})$ is a dense subspace of $S^p_I$.
	
	\item $\M_{I,\fin}$ is a core of $\dom A_p^{\frac{1}{2}}$.
	
	\item We have $\dom \partial_{\alpha,q,p}=\dom A_p^{\frac{1}{2}}$. Moreover, for any $x \in \dom A_p^{\frac{1}{2}}$, we have
\begin{equation} 
\label{Equivalence-square-root-domaine-Schur}
\bnorm{A_p^{\frac12}(x)}_{S^p_I} 
\approx_{p} \bnorm{\partial_{\alpha,q,p}(x)}_{\Gauss_{q,p}(S^p_I)}. 
\end{equation}	
Finally, for any $x \in \dom A_p^{\frac{1}{2}}$ there exists a sequence $(x_n)$ of elements of $\M_{I,\fin}$ such that $x_n \to x$, $A_p^{\frac{1}{2}}(x_n) \to A_p^{\frac{1}{2}}(x)$ and $\partial_{\alpha,q,p}(x_n) \to \partial_{\alpha,q,p}(x)$.

\item If $x \in \dom \partial_{\alpha,q,p}$, we have $x^* \in \dom \partial_{\alpha,q,p}$ and 
\begin{equation}
\label{relation-partial-star}
(\partial_{\alpha,q,p}(x))^* 
=-\partial_{\alpha,q,p}(x^*).
\end{equation}
\item Suppose that $-1 \leq q < 1$.
Then the operator $\partial_{\alpha,q} \co \M_{I,\fin} \subseteq \B(\ell^2_I) \to \Gamma_q(H) \otvn \B(\ell^2_I)$ is weak* closable. We denote by $\partial_{\alpha,q,\infty}$ its weak* closure.

\item Let $M_B \co S^p_I \to S^p_I$ be a finitely supported bounded Schur multiplier such that the map $\Id \ot M_B \co \L^p(\Gamma_q(H) \otvn \B(\ell^2_I)) \to \L^p(\Gamma_q(H) \otvn \B(\ell^2_I))$ is a well-defined bounded operator. For any $x \in \dom \partial_{\alpha,q,p}$, the element $M_B(x)$ belongs to $\dom \partial_{\alpha,q,p}$ and we have 
\begin{equation}
\label{commute-troncature-Schur}
\partial_{\alpha,q,p} M_{B}(x)
=(\Id \ot M_{B}) \partial_{\alpha,q,p}(x).
\end{equation}
\end{enumerate}
\end{prop}

\begin{proof}
The proofs of 1., 2., 3., 4. and 7. are identical to the proof of Proposition \ref{prop-fermable-sgrp} replacing $\P_G$ by $\M_{I,\fin}$, $\L^p(\VN(G))$ by $S^p_I$ and $\L^p(\Gamma_q(H) \rtimes_\alpha  G)$ by $\L^p(\Gamma_q(H) \ot \B(\ell^2_I))$.

5. For any $i,j \in I$, we have
\begin{align*}
\MoveEqLeft
\big(\partial_{\alpha,q,p}(e_{ij})\big)^*
\ov{\eqref{def-delta-alpha}}{=} \big(s_q(\alpha_i-\alpha_j) \ot e_{ij}\big)^*
=s_q(\alpha_i-\alpha_j) \ot e_{ji}
=-s_q(\alpha_j-\alpha_i) \ot e_{ji} \\
&\ov{\eqref{def-delta-alpha}}{=} -\partial_{\alpha,q,p}(e_{ji}) 
=-\partial_{\alpha,q,p}(e_{ij}^*).            
\end{align*}
Let $x \in \dom \partial_{\alpha,q,p}$. By the point 1, $\M_{I,\fin}$ is core of $\partial_{\alpha,q,p}$. Hence there exists a sequence $(x_n)$ of $\M_{I,\fin}$ such that $x_n \to x$ and $\partial_{\alpha,q,p}(x_n) \to \partial_{\alpha,q,p}(x)$. We have $x_n^* \to x^*$ and by the first part of the proof $\partial_{\alpha,q,p}(x_n^*)=-(\partial_{\alpha,q,p}(x_n))^* \to  -(\partial_{\alpha,q,p}(x))^*$. By \eqref{Domain-closure}, we conclude that $x^* \in \dom \partial_{\alpha,q,p}$ and that $\partial_{\alpha,q,p}(x^*)=-(\partial_{\alpha,q,p}(x))^*$.

6. 
Note that $\big(\Id_{\Gamma_q(H) \otvn \B(\ell^2_I)} \ot \Tron_J\big)$ converges for the point weak* topology to $\Id_{\Gamma_q(H) \otvn \B(\ell^2_I)}$. Suppose that $(x_k)$ is a net of $\M_{I,\fin}$ with $x_k=\sum_{i,j \in I} x_{k,i,j}e_{ij}$ which converges to 0 for the weak* topology such that the net $(\partial_{\alpha,q}(x_k))$ converges for the weak* topology to some $y$ belonging to $\Gamma_q(H) \otvn \B(\ell^2_I)$. For any finite subset $J$ of $I$, we have
\begin{align*}
\MoveEqLeft
\partial_{\alpha,q}\Tron_J(x_k)
=\partial_{\alpha,q}\bigg(\sum_{i,j \in J} x_{k,i,j}e_{ij}\bigg) 
=\sum_{i,j \in J} x_{k,i,j} \partial_{\alpha,q}(e_{ij}) \\
&=\big(\Id_{\Gamma_q(H) \otvn \B(\ell^2_I)} \ot \Tron_J\big) \partial_{\alpha,q} \bigg(\sum_{i,j \in I} x_{k,i,j}e_{ij}\bigg) 
=\big(\Id_{\Gamma_q(H) \otvn \B(\ell^2_I)} \ot \Tron_J\big) \partial_{\alpha,q} (x_k) \\
&\xra[k]{}
\big(\Id_{\Gamma_q(H) \otvn \B(\ell^2_I)} \ot \Tron_J\big)(y)            
\end{align*}  
for the weak* topology. On the other hand, for any $i,j \in I$ we have $x_{k,ij} \to 0$ as $k \to \infty$. Hence for any finite subset $J$ of $I$
$$
\partial_{\alpha,q}\Tron_J (x_k) 
= \partial_{\alpha,q} \left( \sum_{i,j \in J} x_{k,i,j} e_{ij} \right) 
\ov{\eqref{def-delta-alpha}}{=}  \sum_{i,j \in J} x_{k,i,j} s_q(\alpha_i - \alpha_j) \ot e_{ij} 
\xra[k]{} 0.
$$ 
This implies by uniqueness of the $\Gamma_q(H) \otvn \B(\ell^2_I)$-limit that $\big(\Id_{\Gamma_q(H) \otvn \B(\ell^2_I)} \ot \Tron_J\big)(y) = 0$. We deduce that $y = \weakstar\lim_J \big(\Id_{\Gamma_q(H) \otvn \B(\ell^2_I)} \ot \Tron_J\big)(y) = 0$.
\end{proof}

\begin{remark} \normalfont
\label{Natural-remark}
It would be more natural to replace the definition \eqref{def-delta-alpha} by the formula 
\[\partial_{\alpha,q}(e_{ij}) 
\ov{\mathrm{def}} 
=2\pi \i s_q(\alpha_i-\alpha_j) \ot e_{ij}.\]
With this new definition, the derivation is symmetric, i.e. $\partial_{\alpha,q}(x^*)=\partial_{\alpha,q}(x)^*$. We refer to \cite{Wea1} for more information on weak* closed derivations on von Neumann algebras.
\end{remark}

Finally, note the extension of Lemma \ref{lem-formule-trace-trace} which can be proved with the point 4 of Proposition \ref{Prop-derivation-closable}. 

\begin{lemma}
\label{lem-formule-trace-trace-2}
Suppose $1 < p < \infty$ $-1 \leq q \leq 1$. For any $x \in \dom A_p^{\frac12}$ and any $y \in \dom A_{p^*}^{\frac12}$, we have
\begin{equation}
\label{formule-trace-trace-2}
(\tau_{\Gamma_q(H)} \ot \tr_{\B(\ell^2_I)})\big((\partial_{\alpha,q,p}(x))^* \partial_{\alpha,q,p^*}(y) \big)
=\tr_{\B(\ell^2_I)} \big((A_p^\frac12(x))^* A_{p^*}^\frac12(y)\big).
\end{equation}
\end{lemma}

\subsection{Meyer's problem for semigroups of Schur multipliers}
\label{Riesz-transforms-cocycle-Schur}

In this subsection, we again fix a markovian semigroup of Schur multipliers, and thus we have some Hilbert space $H$ and a family $(\alpha_i)_{i \in I}$ in $H$ representing the semigroup in the sense of Proposition \ref{def-Schur-markovian}. We investigate the so-called Meyer's problem which consists in expressing the $\L^p$ norm of the gradient form \eqref{Equivalence-square-root-domaine-Schur} in terms of the carr\'e du champ $\Gamma$. Here our results split into the cases $1 < p < 2$ and $2 \leq p < \infty$, and are more satisfactory in the second one. The first main result will be Theorem \ref{cor-Riesz-equivalence-Schur}. In the second part of this subsection, we shall have a look at Riesz transforms \eqref{Riesz-directional-def} associated with the semigroup, carrying a direction vector in the Hilbert space $H$. Then these directional Riesz transforms have some sort of square functions expressing the $S^p_I$ norm, see Theorem \ref{ThA-Schur}.

In the case of finite sized matrices, we are already in the position to state the essence of Meyer's problem for semigroups of Schur multipliers. Namely, the following is an immediate consequence of Lemma \ref{Khintchine-Schur}, Proposition \ref{Prop-equiv-q-gaussians} and \eqref{Equa-Schur-grad-Gamma}. We recall that we have a canonical conditional expectation $\E \co \L^p(\Gamma_q(H) \otvn \B(\ell^2_I)) \to S^p_I$.

\begin{cor}
\label{cor-Riesz-equivalence-Schur-MIfin}
\begin{enumerate}
\item Suppose $1 < p <2$. Let $-1 \leq q \leq 1$. 
For any $x \in \M_{I,\fin}$ we have
$$
\bnorm{A_p^{\frac12}(x)}_{S^p_I} 
\approx_p 
\inf_{\partial_{\alpha,q,p}(x) = y + z}  \norm{y}_{\L_{r}^p(\E)}+\norm{z}_{\L_{c}^p(\E)}
$$
where the infimum is taken over all $y \in \Gauss_{q,2}(\C) \ot \M_{I,\fin}$ and all $z \in \Gauss_{q,2}(\C) \ot \M_{I,\fin}$.

\item Suppose $2 \leq p < \infty$. For any $x \in \M_{I,\fin}$, we have
\begin{equation}
\label{max-Gamma-2}
\bnorm{A_p^{\frac12}(x)}_{S^p_I} 
\approx_p 
\max \left\{ \bnorm{\Gamma(x,x)^{\frac12}}_{S^p_I}, \bnorm{\Gamma(x^*,x^*)^{\frac12}}_{S^p_I} \right\}.
\end{equation}		
\end{enumerate}
\end{cor}

Now, we will extend the definition of $\Gamma(x,y)$ to a larger domain. 

\begin{lemma}
\label{lem-Gamma-closure}
Suppose $2 \leq p <\infty$. The forms $a \co \M_{I,\fin} \times  \M_{I,\fin} \to S^{\frac{p}{2}}_I \oplus_\infty S^{\frac{p}{2}}_I$, $(x,y) \mapsto \Gamma(x,y) \oplus \Gamma(x^*,y^*)$ and $\Gamma \co \M_{I,\fin} \times \M_{I,\fin} \to S^{\frac{p}{2}}_I$, $(x,y) \mapsto \Gamma(x,y)$ are symmetric, positive and satisfy the Cauchy-Schwarz inequality. The first is in addition closable and the domain of the closure $\ovl{a}$ is $\dom A^{\frac12}_p$.
\end{lemma}

\begin{proof}
According to the point 3 of Lemma \ref{Lemma-gamma-infos}, we have $\Gamma(x,y)^* = \Gamma(y,x)$, so $a$ is symmetric. Moreover, again according to Lemma \ref{Lemma-gamma-infos}, $a$ is positive. For any $x,y \in \M_{I,\fin}$, we have
\begin{align*}
\MoveEqLeft
\norm{a(x,y)}_{S^{\frac{p}{2}}_I \oplus_\infty S^{\frac{p}{2}}_I}            
=\norm{\Gamma(x,y) \oplus \Gamma(x^*,y^*)}_{S^{\frac{p}{2}}_I \oplus_\infty S^{\frac{p}{2}}_I} 
=\max\big\{\norm{\Gamma(x,y)}_{\frac{p}{2}},\norm{\Gamma(x^*,y^*)}_{\frac{p}{2}}\big\}\\
&\ov{\eqref{Inegalite-carre-Schur}}{\leq} \max\Big\{\norm{\Gamma(x,x)}_{\frac{p}{2}}^{\frac{1}{2}} \norm{\Gamma(y,y)}_{\frac{p}{2}}^{\frac{1}{2}},\norm{\Gamma(x^*,x^*)}_{\frac{p}{2}}^{\frac{1}{2}}\norm{\Gamma(y^*,y^*)}_{\frac{p}{2}}^{\frac{1}{2}}\Big\}\\
&\leq \max\big\{\norm{\Gamma(x,x)}_{\frac{p}{2}}^{\frac{1}{2}},\norm{\Gamma(x^*,x^*)}_{\frac{p}{2}}^{\frac{1}{2}}\big\}\max\big\{\norm{\Gamma(y,y)}_{\frac{p}{2}}^{\frac{1}{2}},\norm{\Gamma(y^*,y^*)}_{\frac{p}{2}}^{\frac{1}{2}}\Big\}\\
&=\norm{a(x,x)}_{S^{\frac{p}{2}}_I \oplus_\infty S^{\frac{p}{2}}_I}^{\frac{1}{2}}\norm{a(y,y)}_{S^{\frac{p}{2}}_I \oplus_\infty S^{\frac{p}{2}}_I}^{\frac{1}{2}}.
\end{align*}
So $a$ satisfies the Cauchy-Schwarz inequality. The assertions concerning $\Gamma$ are similar. Suppose $x_n \xra[]{a} 0$ that is $x_n \in \M_{I,\fin}$, $x_n \to 0$ and $a(x_n-x_m,x_n-x_m) \to 0$. For any integer $n,m$, we have
\begin{align*}
\MoveEqLeft
\bnorm{A_p^{\frac12}(x_n-x_m)}_{S^p_I} 
\ov{\eqref{max-Gamma-2}}{\lesssim_p} 
\max \left\{ \bnorm{\Gamma(x_n-x_m,x_n-x_m)^{\frac12}}_{S^p_I}, \bnorm{\Gamma((x_n-x_m)^*,(x_n-x_m)^*)^{\frac12}}_{S^p_I} \right\} \\
&=\norm{a(x_n-x_m,x_n-x_m)}_{S^{\frac{p}{2}}_I \oplus_\infty S^{\frac{p}{2}}_I}  
\to 0.
\end{align*}  
We infer that $\big(A_p^{\frac12}(x_n)\big)$ is a Cauchy sequence, hence a convergent sequence. Since $x_n \to 0$, by the closedness of $A_p^{\frac12}$, we deduce that $A_p^{\frac12}(x_n) \to 0$. Now, we have
\begin{align*}
\MoveEqLeft
\norm{a(x_n,x_n)}_{S^{\frac{p}{2}}_I \oplus_\infty S^{\frac{p}{2}}_I}            
=\max \left\{ \bnorm{\Gamma(x_n,x_n)^{\frac12}}_{S^p_I}, \bnorm{\Gamma(x_n^*,x_n^*))^{\frac12}}_{S^p_I} \right\}
\ov{\eqref{max-Gamma-2}}{\lesssim_p} \bnorm{A_p^{\frac12}(x_n)}_{S^p_I}
\to 0.
\end{align*}  
By Proposition \ref{Prop-closable}, we obtain the closability of the form $a$.

Let $x \in S^p_I$. By \eqref{Domaine-closure}, we have $x \in \dom \ovl{a}$ if and only if there exists a sequence $(x_n)$ of $\M_{I,\fin}$ satisfying $x_n \xra[]{a} x$, that is satisfying $x_n \to x$ and $\Gamma(x_n-x_m,x_n-x_m) \to 0$, $\Gamma((x_n-x_m)^*,(x_n-x_m)^*) \to 0 \text{ as } n,m \to \infty$. By the equivalence \eqref{max-Gamma-2}, this is equivalent to the existence of a sequence $x_n \in \M_{I,\fin}$, such that $x_n \to x$, $A_p^{\frac12}(x_n-x_m) \to 0$ as $n,m \to \infty$. Now recalling that $\M_{I,\fin}$ is a core of the operator $A_p^{\frac12}$, we conclude that this is equivalent to $x \in \dom A_p^{\frac12}$.
\end{proof}

\begin{remark} \normalfont
\label{rem-after-lem-Gamma-closure}
The closability of $\Gamma$ and the biggest reasonable domain of the form $\Gamma$ are unclear.
\end{remark}

If $2 \leq p < \infty$, and $x,y \in \dom A_p^{\frac12}$, then we let $\Gamma(x,y)$ be the first component of $\ovl{a}(x,y)$, where $a \co \M_{I,\fin} \times \M_{I,\fin} \to S^{\frac{p}{2}}_I \oplus_\infty S^{\frac{p}{2}}_I,\: (u,v) \mapsto \Gamma(u,v) \oplus \Gamma(u^*,v^*)$ is the form in Lemma \ref{lem-Gamma-closure}.

%
	%

\begin{lemma}
\label{Lemma-Hilbert-module2}
Suppose $-1 \leq q \leq 1$ and $2 \leq p<\infty$\footnote{\thefootnote. In the proof, we recall that $\L^p(\Gamma_q(H) \otvn \B(\ell^2_I)) \subseteq S^p_I(\L^2(\Gamma_q(H))_{c,p})$.}.
\begin{enumerate}
\item For any $x,y \in \dom A^{\frac12}_p = \dom \partial_{\alpha,q,p}$, we have 
\begin{equation}
\label{Equa-Schur-grad-Gamma-2}
\Gamma(x,y) 
=\big\langle \partial_{\alpha,q,p} (x) , \partial_{\alpha,q,p} (y) \big\rangle_{S^p_I(\L^2(\Omega)_{c,p})}
=\E\big((\partial_{\alpha,q,p}(x))^* \partial_{\alpha,q,p}(y)\big).
\end{equation}
\item For any $x \in \dom A^{\frac12}_p = \dom \partial_{\alpha,q,p}$, we have
\begin{equation}
\label{Last-equality}
\bnorm{\Gamma(x,x)^{\frac12}}_{S^p_I} 
=\bnorm{\partial_{\alpha,q,p}(x)}_{S^p_I(\L^2(\Gamma_q(H))_{c,p})}.
\end{equation}
\end{enumerate}
\end{lemma}


\begin{proof}
Consider some elements $x,y \in \dom A^{\frac12}_p=\dom \partial_{\alpha,q,p}$. According to Proposition \ref{Prop-derivation-closable}, $\M_{I,\fin}$ is a core of $\dom \partial_{\alpha,q,p}$. So by \eqref{Def-core}, there exist sequences $(x_n)$ and $(y_n)$ of $\M_{I,\fin}$ such that $x_n \to x$, $y_n \to y$, $\partial_{\alpha,q,p}(x_n) \to \partial_{\alpha,q,p}(x)$ and $\partial_{\alpha,q,p}(y_n) \to \partial_{\alpha,q,p}(y)$. Since $p \geq 2$, we have a contractive inclusion $\L^p(\Gamma_q(H) \otvn \B(\ell^2_I)) \hookrightarrow \L^p_c(\E) \ov{\eqref{Egalite-fantastique}}{=} S^p_I(\L^2(\Gamma_q(H))_{c,p})$.  We deduce that $\partial_{\alpha,q,p}(x_n) \to \partial_{\alpha,q,p}(x)$ and $\partial_{\alpha,q,p}(y_n) \to \partial_{\alpha,q,p}(y)$ in the space $S^p_I(\L^2(\Gamma_q(H))_{c,p})$. Note that
\begin{align*}
\MoveEqLeft
\norm{\Gamma(x_n-x_m,x_n-x_m)}_{S^{\frac{p}{2}}_I}
\ov{\eqref{Equa-Schur-grad-Gamma}}{=} \norm{\big\langle \partial_{\alpha,q,p}(x_n-x_m), \partial_{\alpha,q,p}(x_n-x_m) \big\rangle_{S^p_I(\L^2(\Gamma_q(H))_{c,p})}}_{S^{\frac{p}{2}}_I}\\
&\ov{\eqref{Cauchy-Schwarz-Lp-module}}{\leq} \norm{\partial_{\alpha,q,p}(x_n-x_m)}_{S^p_I(\L^2(\Gamma_q(H))_{c,p})} \norm{\partial_{\alpha,q,p}(x_n-x_m)}_{S^p_I(\L^2(\Gamma_q(H))_{c,p})} 
\xra[n,m \to +\infty]{} 0
\end{align*}
and similarly for $y_n$. Hence $x_n \xra[]{\Gamma} x$ and $y_n \xra[]{\Gamma} y$. We deduce that
\begin{align*}
\MoveEqLeft
\Gamma(x,y)
\ov{\eqref{Def-closure}}{=}\lim_{n \to +\infty} \Gamma(x_n,y_n)
\ov{\eqref{Equa-Schur-grad-Gamma}}{=}\lim_{n \to +\infty}  \big\langle \partial_{\alpha,q,p}(x_n), \partial_{\alpha,q,p}(y_n) \big\rangle_{S^p_I(\L^2(\Gamma_q(H))_{c,p})} \\
&=\big\langle \partial_{\alpha,q,p}(x), \partial_{\alpha,q,p}(y) \big\rangle_{S^p_I(\L^2(\Gamma_q(H))_{c,p})}.            
\end{align*}  

2. If $x \in \dom A^{\frac12}_p$, we have
\begin{align*}
\MoveEqLeft
\bnorm{\partial_{\alpha,q,p}(x)}_{S^{p}_I(\L^2(\Gamma_q(H))_{c,p})}           
\ov{\eqref{Norm-Lp-module}}{=}\Bnorm{\big\langle\partial_{\alpha,q,p}(x),\partial_{\alpha,q,p}(x)\big\rangle_{S^{p}_I(\L^2(\Gamma_q(H))_{c,p})}^{\frac{1}{2}}}_{S^{p}_I}
\ov{\eqref{Equa-Schur-grad-Gamma-2}}{=} \bnorm{\Gamma(x,x)^{\frac{1}{2}}}_{S^{p}_I}.
\end{align*}
\end{proof}

The following is our first main theorem of this subsection.

\begin{thm}
\label{cor-Riesz-equivalence-Schur}
Suppose $2 \leq p < \infty$. For any $x \in \dom A_p^{\frac12}$, we have
\begin{equation}
\label{equivalence-Ap-Schur-dom}
\bnorm{A_p^{\frac12}(x)}_{S^p_I} 
\approx_p 
\max \left\{ \bnorm{\Gamma(x,x)^{\frac12}}_{S^p_I}, \bnorm{\Gamma(x^*,x^*)^{\frac12}}_{S^p_I} \right\}.
\end{equation}		
\end{thm}

\begin{proof}
Pick any $-1 \leq q \leq 1$.
For any $x \in \dom A_p^{\frac12}$, first note that
\begin{equation}
\label{equa-inter-10001}
\E\big(\partial_{\alpha,q,p}(x) (\partial_{\alpha,q,p}(x))^*\big)
\ov{\eqref{relation-partial-star}}{=} \E\big((\partial_{\alpha,q,p}(x^*))^* \partial_{\alpha,q,p}(x^*)\big)
\ov{\eqref{Equa-Schur-grad-Gamma-2}}{=} \Gamma(x^*,x^*).
\end{equation}
We conclude that
\begin{align*}
\MoveEqLeft         
\bnorm{A_p^{\frac12}(x)}_{S^p_I} 
\ov{\eqref{Equivalence-square-root-domaine-Schur}}{\approx_{p}} \bnorm{\partial_{\alpha,q,p}(x)}_{\Gauss_{q,p}(S^p_I)} \\
&\ov{\eqref{Khintchine-Schur-p>2-full-space} }{\approx_p}\max \bigg\{ \norm{\big(\E((\partial_{\alpha,q,p}(x))^* \partial_{\alpha,q,p}(x))\big)^{\frac12}}_{S^p_I},\norm{\big(\E(\partial_{\alpha,q,p}(x) (\partial_{\alpha,q,p}(x)^*))\big)^{\frac12}}_{S^p_I} \bigg\}\\
&\ov{\eqref{Equa-Schur-grad-Gamma-2}\eqref{equa-inter-10001}}{=} \max \left\{ \bnorm{\Gamma(x,x)^{\frac12}}_{S^p_I}, \bnorm{\Gamma(x^*,x^*)^{\frac12}}_{S^p_I} \right\}.
\end{align*} 
\end{proof}

In Lemma \ref{lem-Gamma-closure}, we have shown that the carr\'e du champ $\Gamma$ is a closable form. In fact, one can give a concrete way how to reach all of its domain by approximation with $\M_{I,\fin}$ matrices.
This is the content of the next two lemmas.
Recall the truncations $\Tron_J$ from Definition \ref{def-Troncature}.

\begin{lemma}
\label{lem-Gamma-cut-off-contractive}
Suppose $-1 \leq q \leq 1$ and $2 \leq p < \infty$. Let $J$ be a finite subset of $I$. If $x \in \dom A^{\frac12}_{p}$, then
$$
\norm{\Gamma\big(\Tron_J(x),\Tron_J(x)\big)}_{S^{\frac{p}{2}}_I} 
\leq \bnorm{\partial_{\alpha,q,p}(x)}_{S^{p}_I(\L^2(\Gamma_q(H))_{c,p})}.
$$
\end{lemma}

\begin{proof}
Let $x \in \dom A^{\frac12}_{p}=\dom \partial_{\alpha,q,p}$. Since $\M_{I,\fin}$ is a core of $\partial_{\alpha,q,p}$ by Proposition \ref{Prop-derivation-closable}, there exists a sequence $(x_n)$ of $\M_{I,\fin}$ such that $x_n \to x$ and $\partial_{\alpha,q,p}(x_n) \to \partial_{\alpha,q,p}(x)$. Note that since $\Tron_J \co S^{p}_I \to S^{p}_I$ is a complete contraction, the linear map $\Tron_J \ot \Id_{\L^2(\Gamma_q(H))} \co S^{p}_I(\L^2(\Gamma_q(H))_{c,p}) \to S^{p}_I(\L^2(\Gamma_q(H))_{c,p})$ is also a contraction according to \eqref{extens-col-1}. We deduce that
\begin{align*}
\MoveEqLeft
\label{Equa-inter-1000}
\partial_{\alpha,q,p} \Tron_J(x_n) 
=\sum_{i,j \in J} x_{nij} \partial_{\alpha,q,p}(e_{ij})
=\big(\Id_{\L^2(\Gamma_q(H))} \ot \Tron_J\big)\bigg(\sum_{i,j \in I} x_{nij} \partial_{\alpha,q,p}(e_{ij})\bigg) \\
&=\big(\Id_{\L^2(\Gamma_q(H))}  \ot \Tron_J\big) (\partial_{\alpha,q,p} x_n)
\xra[n \to +\infty]{} \big(\Id_{\L^2(\Gamma_q(H))}  \ot \Tron_J\big) (\partial_{\alpha,q,p} x).            
\end{align*}  
Since $\partial_{\alpha,q,p} \Tron_J$ is bounded by \cite[Problem 5.22]{Kat1}, we deduce that $\partial_{\alpha,q,p} \Tron_J(x)=\big(\Id_{\L^2(\Gamma_q(H))}  \ot \Tron_J\big) (\partial_{\alpha,q,p} x)$. Now, we have
\begin{align*}
\MoveEqLeft
\bnorm{\Gamma(\Tron_Jx,\Tron_Jx)}_{S^{\frac{p}{2}}_I}^{\frac{1}{2}} 
=\bnorm{\Gamma(\Tron_Jx,\Tron_Jx)^\frac{1}{2}}_{S^{p}_I}
\ov{\eqref{nabla-norm-Lp-gradient}}{=} \bnorm{\partial_{\alpha,q,p}\Tron_J(x)}_{{S^{p}_I(\L^2(\Gamma_q(H))_{c,p})}} \\
&=\bnorm{(\Id_{\L^2(\Gamma_q(H))} \ot \Tron_J) (\partial_{\alpha,q,p} x)}_{{S^{p}_I(\L^2(\Gamma_q(H))_{c,p})}}
\ov{\eqref{extens-col-1}}{\leq} \bnorm{\partial_{\alpha,q,p}(x)}_{{S^{p}_I(\L^2(\Gamma_q(H))_{c,p})}}. 
\end{align*} 
\end{proof}

Now, we give a very concrete way to approximate the carr\'e du champ.

\begin{lemma}
\label{lemma-agrandir-domaine-1}
Let $2 \leq p < \infty$. For any $x,y \in \dom A^{\frac12}_{p}$, we have in $S^{\frac{p}{2}}_I$ 
\begin{equation}
\label{Def-Gamma-schur}
\Gamma(x,y)
=\lim_{J} \Gamma\big(\Tron_J(x),\Tron_J(y)\big).
\end{equation} 
\end{lemma}

\begin{proof}
If $x \in \dom A^{\frac12}_{p}$, we have for any finite subsets $J,K$ of $I$
\begin{align*}
\MoveEqLeft
\norm{\Gamma(\Tron_Jx-\Tron_Kx,\Tron_Jx-\Tron_Kx)}_{S^{\frac{p}{2}}_I}^{\frac{1}{2}}  
=\bnorm{\Gamma(\Tron_Jx-\Tron_Kx,\Tron_Jx-\Tron_Kx)^{\frac{1}{2}}}_{S^p_I}   \\         
&\ov{\eqref{equivalence-Ap-Schur-dom}}{\lesssim_p} \bnorm{A_p^{\frac12}(\Tron_Jx-\Tron_Kx)}_{S^p_I} 
= \bnorm{\Tron_J A_p^{\frac12} x-\Tron_KA_p^{\frac12} x}_{S^p_I}.
\end{align*}
Note that since $A^{\frac12}_p(x)$ belongs to $S^p_I$, $(\Tron_JA_p^{\frac12}(x))_p$ is a Cauchy net of $S^p_I$. Since $\Tron_J(x) \to x$ in $S^p_I$, we infer that $\Tron_J(x) \xra{\Gamma} x$. Then \eqref{Def-Gamma-schur} is a consequence of Lemma \ref{lem-Gamma-cut-off-contractive} and Proposition \ref{Prop-existence-limit}. 
\end{proof}



In the second half of this subsection, we shall show that directional Riesz transforms associated with markovian semigroups of Schur multipliers decompose the $S^p_I$ norm.
We recall that since $A$ generates a markovian semigroup, we have 
\begin{equation}
\label{equa-2-Sp0}
\SpI
= \Big\{x \in S^p_I  :  \lim_{t \to +\infty} T_{t}(x)= 0 \Big\}.
\end{equation}
Since $T_{t}(e_{ij}) = \e^{-t\norm{\alpha_i - \alpha_j}_H^2} e_{ij}$ for any $i,j \in I$ and any $t \geq 0$, we have
\begin{equation}
\label{equa-Sp0}
\SpI
=\big\{x \in S^p_I  :  x_{ij} = 0  \text{ for all }  i,j \text{ with } \alpha_i=\alpha_j \big\}.	
\end{equation}
Note that $(T_t)_t$ is strongly continuous on $S^\infty_I$ and that it is not difficult to show that \eqref{equa-2-Sp0} and \eqref{equa-Sp0} also hold for $\SinftyI$, where $A_\infty$ denotes the generator of $(T_t)_t$ on $S^\infty_I$.
If $h \in H$, we define the $h$-directional Riesz transform $R_{\alpha,h}$ defined on $\M_{I,\fin}$ by
\begin{equation}
\label{Riesz-directional-def}
R_{\alpha,h}(e_{ij})
\ov{\mathrm{def}}{=} \frac{\langle \alpha_i-\alpha_j, h \rangle_{H}}{\norm{\alpha_i-\alpha_j}_H} e_{ij}
\text{ if } i,j \text{ satisfy }  \alpha_i \not=\alpha_j	
\end{equation}
and $R_{\alpha,h}(e_{ij})=0$ if it is not the case. If $(e_k)_{k \in K}$ is an orthonormal basis of the Hilbert space $H$, we let
\begin{equation}
\label{Riesz-def-k-Schur}
R_{\alpha,k}
\ov{\mathrm{def}}{=} R_{\alpha,e_k}.
\end{equation} 

Suppose $1 < p<\infty$. Consider the contractive linear map 
$$
U \co \L^2(\Gamma_q(H)) \to \ell^2_K, g \mapsto \sum_{k} \big\langle s_q(e_k),g \big\rangle_{\L^2(\Gamma_q(H))} e_k.
$$ 
By \eqref{extens-rad-1}, we have a bounded map $u \ov{\mathrm{def}}{=}\Id_{S^p_I} \ot U \co S^p_I(\L^2(\Gamma_q(H))_{c,p}) \to S^p_I(\ell^2_{K,c,p}) \subseteq S^p_I(S^p_K)$. For any $g \in \L^2(\Gamma_q(H))$ and any $i,j \in I$ we have
\begin{equation}
\label{Def-de-u}
u(e_{ij} \ot g)
=e_{ij} \ot \sum_k \big\langle s_q(e_k),g \big\rangle_{\L^2(\Gamma_q(H))}e_{k1}
=\sum_k \big\langle s_q(e_k),g \big\rangle_{\L^2(\Gamma_q(H))} e_{ij} \ot e_{k1}.
\end{equation}
With the linear form $\langle s_q(e_k), \cdot\rangle \co \L^2(\Gamma_q(H)) \to \C$ and \eqref{extens-rad-1}, we can introduce the map 
\[u_k \ov{\mathrm{def}}{=}\Id_{S^p_I} \ot \langle s_q(e_k), \cdot\rangle \co S^p_I(\L^2(\Gamma_q(H))_{c,p}) \to S^p_I.\]
 For any $g \in \L^2(\Gamma_q(H))$, any $i,j \in I$ and any $k \in K$, we have
\begin{equation}
\label{Def-de-upk}
u_{k}(e_{ij} \ot g)
=\big\langle s_q(e_k),g \big\rangle_{\L^2(\Gamma_q(H))} e_{ij}.
\end{equation}
For any $k \in I$ and any $f \in \L^2(\Gamma_q(H)) \ot S^p_I$, we have
\begin{equation}
\label{Def-de-up}
u(f)
=\sum_{k} u_{k}(f) \ot e_{k1}.
\end{equation}


In the following proposition, recall again that we have a canonical conditional expectation $\E \co \L^p(\Gamma_q(H) \otvn \B(\ell^2_I)) \to S^p_I$.

\begin{prop}
\label{Prop-2-pi}
Let $-1 \leq q \leq 1$.
\begin{enumerate}
	\item Suppose $2 \leq p <\infty$. For any $f,h \in \Gauss_{q,p}(S^p_I)$ we have 
\begin{equation}
\label{Facto-Esperance-Schur}
\E(f^*h) 
=(u(f))^* u(h)
\end{equation}
where we identify $(u(f))^* u(h)$ as an element of $S^{\frac{p}{2}}_I$ (with its unique non-zero entry).

\item Suppose $1<p<2$. For any $f,h \in \Gauss_{q,2}(\C) \ot S^p_I$ we have 
\begin{equation}
\label{Facto-Esperance-Schur-2}
\langle f,h \rangle_{\L^p_c(\E)}
=(u(f))^* u(h)
\end{equation}
where we identify $(u(f))^* u(h)$ as an element of $S^{\frac{p}{2}}_I$ (with its unique non-zero entry). 

\item Suppose $1<p <\infty$. If $x \in \MIfinzero$, we have\footnote{\thefootnote. Recall that if $2 \leq p <\infty$, we have $\L^p(\Gamma_q(H) \otvn \B(\ell^2_I)) \subseteq S^p_I(\L^2(\Gamma_q(H))_{c,p})$.}
\begin{equation}
\label{u-de-part}
u_k\partial_{\alpha,q,p} A_p^{-\frac12}(x) 
=R_{\alpha,k}(x).
\end{equation}

\item Suppose $1<p <\infty$. If $h \in S^p_I(\L^2(\Gamma_q(H))_{c,p})$, for any $k \in K$ we have 
\begin{equation}
\label{star-passe-a-travers}
 \big(u_k(h^*))^*    
=u_k(h). 
\end{equation} 
\end{enumerate}
\end{prop}

\begin{proof}
1 and 2. For any $i,j,r,s \in I$ and any $g,w \in \Gauss_{q,2}(\C)$ (this assumption is used in a crucial way in the last equality), we have
\begin{align*}
\MoveEqLeft
(u(e_{ij} \ot g))^* u(e_{rs} \ot w)   
\ov{\eqref{Def-de-u}}{=}\bigg(\sum_k \big\langle s_q(e_k),g \big\rangle_{\L^2} e_{ij} \ot e_{k1}\bigg)^*\bigg(\sum_l \big\langle s_q(e_l),w \big\rangle_{\L^2} e_{rs} \ot e_{l1}\bigg)
\\
&=\bigg(\sum_k \ovl{\big\langle s_q(e_k),g \big\rangle_{\L^2(\Gamma_q(H))}} e_{ij}^* \ot e_{1k}\bigg)\bigg(\sum_l \big\langle s_q(e_l),w \big\rangle_{\L^2(\Gamma_q(H))} e_{rs} \ot e_{l1}\bigg)\\
&=\sum_{k,l}\bigg( \ovl{\big\langle s_q(e_k),g \big\rangle_{\L^2(\Gamma_q(H)}} \big\langle s_q(e_l),w \big\rangle_{\L^2(\Gamma_q(H))} e_{ij}^*e_{rs} \ot e_{1k}e_{l1} \bigg)\\
&=\sum_{k} \ovl{\big\langle s_q(e_k),g \big\rangle_{\L^2(\Gamma_q(H))}} \big\langle s_q(e_k),w \big\rangle_{\L^2(\Gamma_q(H))} e_{ij}^*e_{rs} \ot e_{11} 
=\sum_{k}\langle g, w\rangle_{\L^2(\Gamma_q(H))} e_{ij}^*e_{rs}\ot e_{11}.
\end{align*}    
On the other hand, we have
\begin{align*}
\MoveEqLeft
\langle e_{ij} \ot g, e_{rs} \ot w \rangle_{S^p_I(\L^2(\Gamma_q(H))_{c,p})} 
\ov{\eqref{def-scal-product-Lp-H}}{=}
\langle g, w \rangle_{\L^2(\Gamma_q(H))} e_{ij}^*e_{rs}.
\end{align*}
We conclude by bilinearity and density.  

3. For any $i,j \in I$ such that $\alpha_i \not= \alpha_j$, we have 
\begin{align*}
\label{Equa-u-delta-alpha}
\MoveEqLeft
u_k\partial_{\alpha,q,p} A_p^{-\frac12}(e_{ij})
=u_k\partial_{\alpha,q,p} \bigg( \frac{1}{\norm{\alpha_i-\alpha_j}_{H}} e_{ij}\bigg) 
=\frac{1}{\norm{\alpha_i-\alpha_j}_{H}} u_k\partial_{\alpha,q,p}(e_{ij})\\ 
&\ov{\eqref{def-delta-alpha}}{=} \frac{1}{\norm{\alpha_i-\alpha_j}_{H}} u_k\big(e_{ij} \ot s_q(\alpha_i-\alpha_j)\big) 
\ov{\eqref{Def-de-upk}}{=} \frac{1}{\norm{\alpha_i-\alpha_j}_{H}} \big\langle s_q(e_k),s_q(\alpha_i-\alpha_j) \big\rangle_{\L^2(\Gamma_q(H))} e_{ij}\\ 
&\ov{\eqref{petit-Wick}}{=} \frac{\langle e_k,\alpha_i-\alpha_j \rangle_{H}}{\norm{\alpha_i-\alpha_j}_{H}} e_{ij}
\ov{\eqref{Riesz-directional-def}\eqref{Riesz-def-k-Schur}}{=} R_{\alpha,k}(e_{ij}). 
\end{align*}

4. Since $s_q(e_k)$ is selfadjoint, we have for any $i,j \in I$ any $g \in \L^2(\Gamma_q(H))$
\begin{align*}
\MoveEqLeft
\big[u_{k}((e_{ij} \ot g)^*)\big]^* 
=\big[u_{k}(e_{ji} \ot g^*)\big]^* 
\ov{\eqref{Def-de-upk}}{=} \Big(\big\langle s_q(e_k),g^* \big\rangle_{\L^2(\Gamma_q(H))} e_{ji}\Big)^* \\
&=\ovl{\big\langle s_q(e_k), g^* \big\rangle_{\L^2(\Gamma_q(H))}} e_{ij} 
=\big\langle g^*,s_q(e_k) \big\rangle_{\L^2(\Gamma_q(H))} e_{ij}
=\tau(g s_q(e_k)) e_{ij}\\
&=\tau(s_q(e_k)g) e_{ij}
=\big\langle s_q(e_k),g \big\rangle_{\L^2(\Gamma_q(H))} e_{ij}
\ov{\eqref{Def-de-upk}}{=} u_{k}(e_{ij} \ot g). 
\end{align*} 
We conclude by linearity and density.
\end{proof}

\begin{lemma}
\label{lem-proof-ThA-Schur}
Suppose $1<p \leq 2$. The restriction of the map $u \co S^p_I(\L^2(\Gamma_q(H))_{c,p}) \to S^p_I(S^p_K)$ on $S^p_I(\Gauss_{q,2}(\C)_{c,p})$ induces an isometric map and the range of this restriction is $S^p_I(\ell^2_{K,c,p})$.
\end{lemma}

\begin{proof}
For any $f \in S^p_I \ot \Gauss_{q,2}(\C)$, we have
\begin{align*}
\MoveEqLeft
\norm{f}_{\L^p_{c}(\E)}
\ov{\eqref{Norm-Lp-module}}{=} \bnorm{\langle f,f\rangle_{\L^p_{c}(\E)}}_{S^{\frac{p}{2}}_I}^\frac{1}{2} 
\ov{\eqref{Facto-Esperance-Schur-2}}{=} \bnorm{(u(f))^* u(f)}_{S^{\frac{p}{2}}_I}^\frac{1}{2}
=\bnorm{u(f)}_{S^p_I(S^p_K)}.
\end{align*} 
and $u(e_{ij} \ot s_q(e_l)) \ov{\eqref{Def-de-u}}{=} e_{ij} \ot e_{l}$ for any $i,j \in I$ and any $k \in K$. Alternatively, note that the linear map $\Gauss_{q,2}(\C) \to \ell^2_K$, $g \mapsto \sum_{k \in K} \big\langle s_q(e_k), g \big\rangle_{\L^2(\Gamma_q(H))} e_k$ is a surjective isometry. So the associated map $S^p_I(\Gauss_{q,2}(\C)_{c,p}) \to S^p_I(\ell^2_{K,c,p})$ is also a surjective isometry.
\end{proof}


\begin{thm}
\label{ThA-Schur}
Suppose $1<p < \infty$. 
\begin{enumerate}
\item If $1 < p \leq 2$ and if $x \in \MIfinzero$ we have
\begin{equation}
\label{inf-de-dingue}
\norm{x}_{S^p_I} 
\approx_{p} \inf_{R_{\alpha,k}(x)=a_k+b_k} \bgnorm{\bigg(\sum_{k \in K} |a_k|^2\bigg)^{\frac12}}_{S^p_I} + \bgnorm{\bigg(\sum_{k \in K} |b_k^*|^2\bigg)^{\frac12}}_{S^p_I}
\end{equation}
where the infimum is taken over all $(a_k), (b_k) \in S^p_I(\ell^2_{K,c})$.
\item If $2 \leq p < \infty$ and if $x \in \MIfinzero$, we have 
$$
\norm{x}_{S^p_I} 
\approx_{p} \max \left\{\bgnorm{\bigg(\sum_{k \in K} |R_{\alpha,k}(x)|^2 \bigg)^{\frac12} }_{S^p_I}, \bgnorm{\bigg(\sum_{k \in K} |\big(R_{\alpha,k}(x)\big)^*|^2\bigg)^{\frac12}}_{S^p_I} \right\}.
$$
\end{enumerate}
\end{thm}


\begin{proof}
Here we use some fixed $-1 \leq q \leq 1$ and we drop the index in the notation $\partial_{\alpha,(q,)p}$.

1. Suppose $1<p<2$. For any $x \in \MIfinzero$, we have
\begin{align*}
\MoveEqLeft
\norm{x}_{S^p_I}
=\bnorm{A_p^\frac12 A_p^{-\frac12}(x)}_{S^p_I}
\ov{\eqref{Equivalence-square-root-domaine-Schur}}{\approx_{p}} \bnorm{\partial_{\alpha,p} A_p^{-\frac12}(x)}_{\L^p(\Gamma_q(H) \otvn \B(\ell^2_I))} \\
&\ov{\eqref{Khintchine-Schur-p<2}}{\approx_p} \inf_{\partial_{\alpha,p}A_p^{-\frac12}(x)=g+h} \norm{g}_{\L_{c}^p(\E)}+\norm{h}_{\L_{r}^p(\E)}
\end{align*}   
where the infimum is taken over all $g,h \in \Gauss_{2}(\C) \ot \M_{I,\fin}$. By Proposition \ref{Prop-2-pi}, we infer that 
\begin{align*}
\MoveEqLeft
\langle g,g \rangle_{\L^p_{c}(\E)}^{\frac{1}{2}}
\ov{\eqref{Facto-Esperance-Schur-2}}{=}\big(u(g)^* u(g)\big)^{\frac{1}{2}}
\ov{\eqref{Def-de-up}}{=}\Bigg(\bigg(\sum_{k} u_{k}(g) \ot e_{k1} \bigg)^*\bigg(\sum_{l} u_{l}(g) \ot e_{l1} \bigg)\Bigg)^{\frac12}\\
&=\bigg(\sum_{k,l} \big(u_{k}(g))^* u_{l}(g) \ot e_{1k}e_{l1}\bigg)^{\frac12}
=\bigg(\sum_{k} \big(u_{k}(g))^* u_{k}(g)\bigg)^{\frac{1}{2}}.
\end{align*}  
Similarly, we have
\begin{align*}
\MoveEqLeft
\langle h^*,h^* \rangle_{\L^p_{c}(\E)}^{\frac{1}{2}}  
\ov{\eqref{Facto-Esperance-Schur-2}}{=}\big((u(h^*))^* u(h^*)\big)^{\frac{1}{2}}
\ov{\eqref{Def-de-up}}{=}\Bigg(\bigg(\sum_{k} u_{k}(h^*) \ot e_{k1}\bigg)^*\bigg(\sum_{l} u_{l}(h^*) \ot e_{l1}\bigg)\Bigg)^{\frac12}\\
&\ov{\eqref{star-passe-a-travers}}{=}\Bigg(\bigg(\sum_{k} u_{k}(h) \ot e_{1k}\bigg)\bigg(\sum_{l} u_{l}(h^*) \ot e_{l1} \bigg)\Bigg)^{\frac12}
=\bigg(\sum_{k} u_{k}(h) u_{k}(h^*) \bigg)^{\frac{1}{2}}.
\end{align*} 
By Lemma \ref{lem-proof-ThA-Schur}, the restriction of the map $u \co S^p_I(\L^2(\Gamma_q(H))_{c,p}) \to S^p_I(S^p_I)$ on $S^p_I(\Gauss_{q,2}(\C)_{c,p})$ is injective. So we have $\partial_{\alpha,p}A_p^{-\frac12}(x)=g+h$ if and only if $u\partial_{\alpha,p}A_p^{-\frac12}(x)=u(g)+u(h)$. By \eqref{u-de-part} and \eqref{Def-de-up}, this is equivalent to $R_{\alpha,k}(x)=u_{k}(g) + u_{k}(h)$ for any $k$. Hence this computation gives that $\norm{x}_{S^p_{I}}$ is comparable to the norm
\begin{equation}
\label{Premier-infimum-R-A-k}
\inf_{R_{\alpha,k}(x)=u_{k}(g) + u_{k}(h)} \bgnorm{\bigg( \sum_{k \in K} \big(u_{k}(g)\big)^* u_{k}(g) \bigg)^\frac12}_{S^p_I} + \bgnorm{\bigg(\sum_{k \in K} \big(u_{k}(h) u_{k}(h^*) \bigg)^\frac12}_{S^p_I}.	
\end{equation}
Therefore, using again Lemma \ref{lem-proof-ThA-Schur}, this infimum is finally equal to the infimum \eqref{inf-de-dingue}.



2. Now, suppose $2 \leq p < \infty$. For any $x \in \MIfinzero$, we have 
\begin{align*}
\MoveEqLeft
\norm{x}_{S^p_{I}}
=\bnorm{A_p^\frac12A_p^{-\frac12}(x)}_{S^p_{I}}
\ov{\eqref{Equivalence-square-root-domaine-Schur}}{\approx_{p}} \bnorm{\partial_{\alpha,p} A_p^{-\frac12}(x)}_{\Gauss_{q,p}(S^p_I)}
\ov{\eqref{Khintchine-Schur-p>2}}{\approx_p} \bnorm{\partial_{\alpha,p} A_p^{-\frac12}(x)}_{\L_{rc}^p(\E)}\\
&=\max\bigg\{\norm{\Big(\E\Big(\big(\partial_{\alpha,p} A_p^{-\frac12}(x)\big)^*\partial_{\alpha,p} A_p^{-\frac12}(x)\Big)\Big)^{\frac12}}_{p},\norm{\Big( \E\Big(\partial_\alpha A_p^{-\frac12}(x) \big(\partial_{\alpha,p} A_p^{-\frac12}(x)\big)^*\Big) \Big)^{\frac12}}_{p} \bigg\}.
\end{align*}  
We obtain
\begin{align*}
\MoveEqLeft
\E\Big(\big(\partial_{\alpha,p} A_p^{-\frac12} x\big)^* \partial_{\alpha,p} A_p^{-\frac12} x \Big) 
\ov{\eqref{Facto-Esperance-Schur}}{=}\big[u\big(\partial_{\alpha,p} A_p^{-\frac12}(x)\big) \big]^* u\big(\partial_{\alpha,p} A_p^{-\frac12}(x)\big) \\
&\ov{\eqref{u-de-part}}{=}\bigg(\sum_{k \in K} R_{\alpha,k}(x) \ot e_{k1}\bigg)^*\bigg(\sum_{l \in K} R_{\alpha,l}(x) \ot e_{l1}\bigg)\\
&=\bigg(\sum_{k \in K} \big(R_{\alpha,k}(x)\big)^* \ot e_{1k}\bigg)\bigg(\sum_{l \in K} R_{\alpha,l}(x) \ot e_{l1}\bigg) \\
&=\sum_{k \in K} \big(R_{\alpha,k}(x)\big)^* R_{\alpha,k}(x)
=\sum_{k \in K} |R_{\alpha,k}(x)|^2.   
\end{align*} 
Similarly, using \eqref{star-passe-a-travers} in addition, we obtain
\begin{align*}
\MoveEqLeft
\E\Big(\partial_{\alpha,p} A_p^{-\frac12}(x) \big(\partial_{\alpha,p} A_p^{-\frac12} x\big)^* \Big) 
\ov{\eqref{Facto-Esperance-Schur}}{=}\big(u\big((\partial_{\alpha,p} A_p^{-\frac12} x)^*\big) \big)^* u\big((\partial_{\alpha,p} A_p^{-\frac12}x)^*\big)\\
&\ov{\eqref{Def-de-up}}{=} \bigg(\sum_{k} u_{k}\big((\partial_{\alpha,p} A_p^{-\frac12} x)^*\big) \ot e_{k1}\bigg)^*\bigg(\sum_{l} u_{l}\big((\partial_{\alpha,p} A_p^{-\frac12} x)^*\big) \ot e_{l1}\bigg)\\
&\ov{\eqref{star-passe-a-travers}}{=} \bigg(\sum_{k} \big(u_{k}\big(\partial_{\alpha,p} A_p^{-\frac12} x\big)\big)^* \ot e_{k1}\bigg)^*\bigg(\sum_{l} \big(u_{l}(\partial_{\alpha,p} A_p^{-\frac12} x)\big)^* \ot e_{l1}\bigg)\\
&=\bigg(\sum_{k} u_{k}\big(\partial_{\alpha,p} A_p^{-\frac12} x\big) \ot e_{1k}\bigg)\bigg(\sum_{l} \big(u_{l}(\partial_{\alpha,p} A_p^{-\frac12} x)\big)^* \ot e_{l1}\bigg)\\
&\ov{\eqref{u-de-part}}{=}\bigg(\sum_{k} R_{\alpha,k}(x) \ot e_{1k}\bigg)\bigg(\sum_{l}(R_{\alpha,l}(x))^*  \ot e_{l1}\bigg)
=\sum_{k \in K} \big|\big(R_{\alpha,k}(x)\big)^*\big|^2.
\end{align*} 
The proof is complete. 
\end{proof}

\section{Boundedness of $\HI$ functional calculus of Hodge-Dirac operators}
\label{Sec-Hodge-Dirac}

\subsection{Boundedness of functional calculus of Hodge-Dirac operators for Fourier multipliers}
\label{subsec-HI-calculus-Hodge-Dirac-Fourier}

In this subsection, we consider a discrete group $G$ and a semigroup $(T_t)_{t \geq 0}$ of Markov Fourier multipliers satisfying Proposition \ref{prop-Schoenberg}. If $1 \leq p <\infty$, we denote by $A_p$ the (negative) infinitesimal generator on $\L^p(\VN(G))$.
By \cite[(5.2)]{JMX}, we have $(A_{p})^*=A_{p^*}$ if $1<p<\infty$. If $-1 \leq q \leq 1$, recall that by Proposition \ref{Prop-derivation-closable-sgrp}, we have a closed operator 
$$
\partial_{\psi,q,p} \co \dom \partial_{\psi,q,p} \subseteq \L^p(\VN(G)) \to \L^p(\Gamma_{q}(H) \rtimes_\alpha G),\quad 
\lambda_s \mapsto s_q(b_\psi(s)) \rtimes \lambda_s	
$$ 
and a closed operator $(\partial_{\psi,q,p^*})^* \co \dom (\partial_{\psi,q,p^*})^* \subseteq \L^p(\Gamma_{q}(H) \rtimes G) \to \L^p(\VN(G))$.  
Now, we will define a Hodge-Dirac operator $D_{\psi,q,p}$ in \eqref{Def-Dirac-operator-Fourier}, from $\partial_{\psi,q,p}$ and its adjoint. 
Then the main topic of this subsection will be to show that $D_{\psi,q,p}$ is $R$-bisectorial (Theorem \ref{Th-D-R-bisectorial-Fourier}) and has a bounded $\HI$ functional calculus on $\L^p(\VN(G)) \oplus \ovl{\Ran \partial_{\psi,q,p}}$ (Theorem \ref{Th-functional-calculus-bisector-Fourier}). By Remark \ref{remark-sgn-Fourier}, this extends the Kato square root equivalence from \eqref{Equivalence-square-root-domaine-Schur-sgrp}.

Note that we use Proposition \ref{prop-Fourier-mult-crossed-product}, so we need some assumption on the group $G$. First, with exception of Proposition \ref{Prop-core-1-group}, we assume that $\Gamma_q(H) \rtimes_\alpha G$ has $\QWEP$ (e.g. $G$ is amenable or is a free group and $q = \pm 1$ by an adaptation of \cite[Proposition 4.8]{Arh1}). Second, since we need approximating compactly supported Fourier multipliers, from Lemma \ref{Petit-lemma} on, we also assume that $G$ is weakly amenable. 

We start with the following intuitive formula which show that $\partial_{\psi,q,p}$ can be seen as a ``gradient'' for $A_p$ in the spirit of the link between the classical laplacian and the classical gradient.

\begin{prop}
\label{Prop-lien-generateur-partial-Fourier}
Suppose $1<p<\infty$ and $-1 \leq q \leq 1$. As unbounded operators, we have
\begin{equation}
\label{eq-lien-generateur-partial-Fourier}
A_p
=(\partial_{\psi,q,p^*})^*\partial_{\psi,q,p}. 
\end{equation}
\end{prop}

\begin{proof}
By Lemma \ref{Lemma-adjoint-prtial-Fourier} and \cite[p.~167]{Kat1}, $\partial_{\psi,q,p}(\P_G)$ is a subspace of $\dom (\partial_{\psi,q,p^*})^*$. For any $s \in G$, we have
\begin{align}
\MoveEqLeft
\label{Equa-10001}
(\partial_{\psi,q,p^*})^*\partial_{\psi,q,p}(\lambda_s)           
\ov{\eqref{def-partial-psi}}{=}(\partial_{\psi,q,p^*})^*\big(s_q(b_\psi(s)) \rtimes \lambda_s\big) 
\ov{\eqref{Adjoint-partial-Fourier}}{=} \tau\big( s_q(b_\psi(s))s_q(b_\psi(s))\big) \lambda_s \\
&\ov{\eqref{petit-Wick}}{=}\norm{b_\psi(s)}_H^2\, \lambda_s
=A_p(\lambda_s). \nonumber
\end{align}
Hence for any $x,y \in \P_G$, by linearity we have 
\begin{align*}
\MoveEqLeft
\big\langle A_p^{\frac{1}{2}}(x), A_{p^*}^{\frac{1}{2}}(y) \big\rangle_{\L^p(\VN(G)),\L^{p^*}(\VN(G))}
=\big\langle A_p(x), y \big\rangle_{\L^p(\VN(G)),\L^{p^*}(\VN(G))} \\
&\ov{\eqref{Equa-10001}}{=} \big\langle (\partial_{\psi,q,p^*})^*\partial_{\psi,q,p}(x), y \big\rangle_{\L^p(\VN(G)),\L^{p^*}(\VN(G))}
\ov{\eqref{crochet-duality}}{=}\big\langle \partial_{\psi,q,p}(x),\partial_{\psi,q,p^*}(y) \big\rangle_{\L^p,\L^{p^*}}.            
\end{align*} 
Using the part 4 of Proposition \ref{Prop-derivation-closable-sgrp}, it is not difficult to see that this identity extends to elements $x \in \dom A_p$. For any $x \in \dom A_p$ and any $y \in \P_G$, we obtain
$$
\big\langle A_p(x), y \big\rangle_{\L^p(\VN(G)),\L^{p^*}(\VN(G))}
=\big\langle \partial_{\psi,q,p}(x),\partial_{\psi,q,p^*}(y) \big\rangle_{\L^p,\L^{p^*}}.
$$
Recall that $\P_G$ is a core of $\partial_{\psi,q,p^*}$ by the part 1 of Proposition \ref{Prop-derivation-closable-sgrp}. So using \eqref{Def-core}, it is easy to check that this identity remains true for elements $y$ of $\dom \partial_{\psi,q,p^*}$. By \eqref{Def-domaine-adjoint}, this implies that $\partial_{\psi,q,p}(x) \in \dom (\partial_{\psi,q,p^*})^*$ and that $(\partial_{\psi,q,p^*})^*\partial_{\psi,q,p}(x) = A_p(x)$. We conclude that $A_p \subseteq (\partial_{\psi,q,p^*})^*\partial_{\psi,q,p}$.

To prove the other inclusion we consider some $x \in \dom \partial_{\psi,q,p}$ such that $\partial_{\psi,q,p}(x)$ belongs to $\dom (\partial_{\psi,q,p^*})^*$. By \cite[Theorem 5.29 p.~168]{Kat1}, we have $(\partial_{\psi,q,p^*})^{**}=\partial_{\psi,q,p^*}$. We infer that $(\partial_{\psi,q,p})^*\partial_{\psi,q,p^*} \ov{\eqref{Adjoint-product-unbounded}}{\subseteq} \big((\partial_{\psi,q,p^*})^*\partial_{\psi,q,p}\big)^* \ov{\eqref{Inclusion-adjoint-unbounded}}{\subseteq} A_p^*$. For any $y \in \P_G$, using $\partial_{\psi,q,p}(x) \in \dom (\partial_{\psi,q,p^*})^*$ in the last equality, we deduce that
\begin{align*}
\MoveEqLeft
\big\langle A_p^*(y),x \big\rangle_{\L^{p^*}(\VN(G)),\L^{p}(\VN(G))}
=\big\langle (\partial_{\psi,q,p})^*\partial_{\psi,q,p^*}(y),x \big\rangle_{\L^{p^*}(\VN(G)),\L^{p}(\VN(G))} \\
&\ov{\eqref{crochet-duality}}{=} \big\langle \partial_{\psi,q,p^*}(y),\partial_{\psi,q,p}(x) \big\rangle_{\L^{p^*},\L^{p}}
\ov{\eqref{crochet-duality}}{=}\big\langle y,(\partial_{\psi,q,p^*})^*\partial_{\psi,q,p}(x) \big\rangle_{\L^{p^*}(\VN(G)),\L^{p}(\VN(G))}.
\end{align*}
Since $\P_G$ is a core for $A_p^*=A_{p^*}$ by definition, this implies \cite[Problem 5.24 p. 168]{Kat1} that $x \in \dom A_p^{**}=A_p$ and that $A_p(x) =  (\partial_{\psi,q,p^*})^*\partial_{\psi,q,p}(x)$.
\end{proof}

Now, we show how the noncommutative gradient $\partial_{\psi,q,p}$ commutes with the semigroup and the resolvents of its generator.

\begin{lemma}
\label{lem-T_t-et-derivations-Fourier}
Let $G$ be a discrete group such that $\Gamma_q(H) \rtimes_\alpha G$ has $\QWEP$. Suppose $1<p<\infty$ and $-1 \leq q \leq 1$. If $x \in \dom \partial_{\psi,q,p}$ and $t \geq 0$, then $T_{t,p}(x)$ belongs to $\dom \partial_{\psi,q,p}$ and we have 
\begin{equation}
\label{T_t-et-derivations-Fourier}
\big(\Id_{\L^p(\Gamma_{q}(H))} \rtimes T_{t,p}\big) \partial_{\psi,q,p}(x) 
=\partial_{\psi,q,p} T_{t,p}(x).
\end{equation}
\end{lemma}

\begin{proof}
For any $s \in G$, we have
\begin{align*}
\MoveEqLeft
\big(\Id_{\L^p(\Gamma_{q}(H))} \rtimes T_{t,p}\big) \partial_{\psi,q,p}(\lambda_s)           
\ov{\eqref{def-partial-psi}}{=} \big(\Id_{\L^p(\Gamma_{q}(H))} \rtimes T_{t,p}\big)\big(s_q(b_\psi(s)) \rtimes \lambda_s\big) \\
&=\e^{-t\norm{b_\psi(s)}^2} s_q(b_\psi(s)) \rtimes \lambda_s 
\ov{\eqref{def-partial-psi}}{=} \e^{-t\norm{b_\psi(s)}^2}\partial_{\psi,q,p}(\lambda_s)
=\partial_{\psi,q,p}\big(\e^{-t\norm{b_\psi(s)}^2}\lambda_s\big) \\
&=\partial_{\psi,q,p} T_{t,p}(\lambda_s).
\end{align*}
So by linearity the equality \eqref{T_t-et-derivations-Fourier} is true for elements of $\P_G$. Now consider some $x \in \dom \partial_{\psi,q,p}$. By \cite[p.~166]{Kat1}, since $\partial_{\psi,q,p}$ is the closure of $\partial_{\psi,q} \co \P_G \subseteq \L^p(\VN(G)) \to \L^p(\Gamma_q(H) \rtimes_\alpha G)$, there exists a sequence $(x_n)$ of elements of $\P_G$ converging to $x$ in $\L^p(\VN(G))$ such that the sequence $(\partial_{\psi,q,p}(x_n))$ converges to $\partial_{\psi,q,p}(x)$. The complete boundedness of $T_{t,p} \co \L^p(\VN(G)) \to \L^p(\VN(G))$ implies by Proposition \ref{prop-Fourier-mult-crossed-product}  that we have a (completely) bounded linear operator $\Id_{\L^p(\Gamma_q(H))} \rtimes T_{t,p} \co \L^p(\Gamma_q(H) \rtimes_\alpha G) \to \L^p(\Gamma_q(H) \rtimes_\alpha G)$. We infer that in $\L^p(\VN(G))$ and $\L^p(\Gamma_q(H) \rtimes_\alpha G)$ we have 
$$
T_{t,p}(x_n) \xra[n \to +\infty]{} T_{t,p}(x)
\ \ \text{and} \ \ 
\big(\Id_{\L^p(\Gamma_q(H))} \rtimes T_{t,p}\big) \partial_{\psi,q,p}(x_n)
\xra[n \to +\infty]{} \big(\Id_{\L^p(\Gamma_q(H))} \rtimes T_{t,p}\big)\partial_{\psi,q,p}(x).
$$ 
For any integer $n$, we have $\big(\Id_{\L^p(\Gamma_q(H))} \rtimes T_{t,p}\big) \partial_{\psi,q,p}(x_n) =\partial_{\psi,q,p} T_{t,p}(x_n)$ by the first part of the proof. Since the left-hand side converges, we obtain that the sequence $(\partial_{\psi,q,p} T_{t,p}(x_n))$ converges to $(\Id_{\L^p(\Gamma_q(H))} \rtimes T_{t,p})\partial_{\psi,q,p}(x)$ in $\L^p(\Gamma_q(H) \rtimes_\alpha G)$. Since each $T_{t,p}(x_n)$ belongs to $\dom \partial_{\psi,q,p}$, the closedness of $\partial_{\psi,q,p}$ and \eqref{Def-operateur-ferme} shows that $T_{t,p}(x)$ belongs to $\dom \partial_{\psi,q,p}$ and that $\partial_{\psi,q,p} T_{t,p}(x)=\big(\Id_{\L^p(\Gamma_q(H))} \rtimes T_{t,p}\big) \partial_{\psi,q,p}(x)$.
\end{proof}


\begin{prop}
\label{Prop-commuting-Fourier}
Let $G$ be a discrete group such that $\Gamma_q(H) \rtimes_\alpha G$ has $\QWEP$. Suppose $1 < p < \infty$ and $-1 \leq q \leq 1$. For any $s \geq 0$ and any $x \in \dom \partial_{\psi,q,p}$, we have $\big(\Id_{}+sA_p\big)^{-1}(x) \in \dom \partial_{\psi,q,p}$ and 
\begin{equation}
\label{commuting-deriv-Fourier}
\big(\Id_{} \rtimes (\Id_{}+sA_p)^{-1}\big)\partial_{\psi,q,p}(x)
=\partial_{\psi,q,p}\big(\Id_{}+sA_p\big)^{-1}(x).
\end{equation}
\end{prop}

\begin{proof}
Note that for any $s>0$ and any $x \in \L^p(\VN(G))$ (resp. $x \in \dom \partial_{\psi,q,p}$) the continuous functions $\R^+ \to \L^p(\VN(G))$, $t \mapsto \e^{-s^{-1}t}T_{t,p}(x)$ and $\R^+ \to \L^p(\Gamma_q(H) \rtimes_\alpha G)$, $t \mapsto \e^{-s^{-1}t}(\Id_{} \rtimes T_{t,p})\partial_{\psi,q,p}(x)$ are Bochner integrable. If $t > 0$ and if $x \in \dom \partial_{\psi,q,p}$, taking Laplace transforms on both sides of \eqref{T_t-et-derivations-Fourier} and using \cite[Theorem 1.2.4]{HvNVW2} and the closedness of $\partial_{\psi,q,p}$ in the penultimate equality, we obtain that $\int_{0}^{\infty} \e^{-ts^{-1}}T_{t,p}(x) \d t$ belongs to $\dom \partial_{\psi,q,p}$ and that
\begin{align*}
\MoveEqLeft
\big(\Id_{} \rtimes (s^{-1}\Id_{}+A_p)^{-1}\big)\partial_{\psi,q,p}(x) 
=-\big(-s^{-1}\Id_{}-(\Id_{} \rtimes A_p)\big)^{-1}\partial_{\psi,q,p}(x)\\
&\ov{\eqref{Resolvent-Laplace}}{=} \int_{0}^{\infty} \e^{-s^{-1}t}(\Id_{} \rtimes T_{s,p})\partial_{\psi,q,p}(x) \d t 
\ov{\eqref{T_t-et-derivations-Fourier}}{=} \int_{0}^{\infty} \e^{-s^{-1}t}\partial_{\psi,q,p}T_{t,p}(x) \d t\\
&=\partial_{\psi,q,p}\bigg(\int_{0}^{\infty} \e^{-s^{-1}t}T_{t,p}(x) \d t\bigg)
\ov{\eqref{Resolvent-Laplace}}{=} 
\partial_{\psi,q,p}\big(s^{-1}\Id_{}+A_p\big)^{-1}(x).    
\end{align*}
We deduce the desired identity by multiplying by $s^{-1}$.
\end{proof}

Now, we prove a result which gives some $R$-boundedness.

\begin{prop}
\label{Prop-R-gradient-bounds-Fourier} 
Suppose $1<p<\infty$ and $-1 \leq q \leq 1$. Let $G$ be a discrete group. The family 
\begin{equation}
\label{R-gradient-bounds-groups}
\Big\{t\partial_{\psi,q,p}(\Id+t^2A_p)^{-1}: t>0 \Big\}
\end{equation} 
of operators of $\B(\L^p(\VN(G)),\L^p(\Gamma_{q}(H) \rtimes_\alpha G))$ is $R$-bounded.
\end{prop}

\begin{proof}
Note that the operator $\partial_{\psi,q,p}A_p^{-\frac{1}{2}} \co \ovl{\Ran A_p} \to \L^p(\Gamma_{q}(H) \rtimes_\alpha G)$ is bounded by \eqref{Equivalence-square-root-domaine-Schur-sgrp}. Suppose $t >0$. A standard functional calculus argument gives
\begin{align}
\label{Divers-987}
\MoveEqLeft
t\partial_{\psi,q,p}(\Id+t^2A_p)^{-1}            
=\partial_{\psi,q,p}A_p^{-\frac{1}{2}}\Big((t^2A_p)^{\frac{1}{2}}(\Id+t^2A_p)^{-1} \Big).
\end{align} 
By \cite{Arh7}, note that $A_p$ has a bounded $\HI(\Sigma_\theta)$ functional calculus for some $0< \theta < \frac{\pi}{2}$. Moreover, the Banach space $\L^p(\VN(G))$ is $\UMD$ by \cite[Corollary 7.7]{PiX}, hence has the triangular contraction property $(\Delta)$ by \cite[Theorem 7.5.9]{HvNVW2}. We deduce by \cite[Theorem 10.3.4 (2)]{HvNVW2} that the operator $A_p$ is $R$-sectorial. By \cite[Example 10.3.5]{HvNVW2} applied with $\alpha=\frac{1}{2}$ and $\beta=1$, we infer that the set
$$
\big\{(t^2A_p)^{\frac{1}{2}}(\Id+t^2A_p)^{-1}: t>0\big\}
$$
of operators of $\B(\L^p(\VN(G)))$ is $R$-bounded. Recalling that a singleton is $R$-bounded by \cite[Example 8.1.7]{HvNVW2}, we obtain by composition \cite[Proposition 8.1.19 (3)]{HvNVW2} that the set
$$
\Big\{ \partial_{\psi,q,p}A_p^{-\frac{1}{2}}\Big((t^2A_p)^{\frac{1}{2}}(\Id+t^2A_p)^{-1} \Big) : t>0 \Big\}
$$
of operators of $\B(\L^p(\VN(G)),\L^p(\Gamma_{q}(H) \rtimes_\alpha G))$ is $R$-bounded. Hence with \eqref{Divers-987} we conclude that the subset \eqref{R-gradient-bounds-groups} is $R$-bounded.
\end{proof}

Our Hodge-Dirac operator in \eqref{Def-Dirac-operator-Fourier} below will be constructed out of $\partial_{\psi,q,p}$ and the unbounded operator $(\partial_{\psi,q,p^*})^*|\ovl{\Ran \partial_{\psi,q,p}}$. Note that the latter is by definition an unbounded operator on the Banach space $\ovl{\Ran \partial_{\psi,q,p}}$ with values in $\L^p(\VN(G))$ having domain $\dom (\partial_{\psi,q,p^*})^* \cap \ovl{\Ran \partial_{\psi,q,p}}$.

\begin{lemma}
\label{Lemma-rest-closed-dense-Fourier}
Let $G$ be a discrete group. Suppose $1 < p  < \infty$ and $-1 \leq q \leq 1$.
The operator $(\partial_{\psi,q,p^*})^*|\ovl{\Ran \partial_{\psi,q,p}}$ is densely defined and is closed. More precisely, the subspace $\partial_{\psi,q,p}(\P_G)$ of $\dom (\partial_{\psi,q,p^*})^*$ is dense in $\ovl{\Ran \partial_{\psi,q,p}}$.
\end{lemma}

\begin{proof}
Let $y \in \ovl{\Ran \partial_{\psi,q,p}}$. Let $\epsi >0$. There exists $x \in \dom \partial_{\psi,q,p}$ such that $\norm{y-\partial_{\psi,q,p}(x)} < \epsi$. By Proposition \ref{prop-fermable-sgrp}, there exist $x_\fin \in \P_G$ such that $\norm{x-x_\fin}_{\L^p(\VN(G))} < \epsi$ and 
$$
\bnorm{\partial_{\psi,q,p}(x)-\partial_{\psi,q,p}(x_\fin)}_{\L^p(\Gamma_q(H) \rtimes_\alpha G)} 
< \epsi.
$$ 
We deduce that $\norm{y-\partial_{\psi,q,p}(x_\fin)}_{\L^p(\Gamma_q(H) \rtimes_\alpha G} < 2 \epsi$. By Proposition \ref{Prop-lien-generateur-partial-Fourier}, $\partial_{\psi,q,p}(\P_G)$ is a subspace of $\dom (\partial_{\psi,q,p^*})^*$. So $\partial_{\psi,q,p}(x_\fin)$ belongs to $\dom (\partial_{\psi,q,p^*})^*$. 

Since $(\partial_{\psi,q,p^*})^*$ is closed, the assertion on the closedness is (really) obvious.
\end{proof}


According to Lemma \ref{lem-T_t-et-derivations-Fourier}, $\Id \rtimes T_t$ leaves $\Ran \partial_{\psi,q,p}$ invariant for any $t \geq 0$, so by continuity of $\Id \rtimes T_t$ also leaves $\ovl{\Ran \partial_{\psi,q,p}}$ invariant. By \cite[pp.~60-61]{EnN1}, we can consider the generator
\begin{equation}
\label{Def-Bp-Fourier}
B_p 
\ov{\mathrm{def}}{=} (\Id_{\L^p(\Gamma_q(H))} \rtimes A_p)|\ovl{\Ran \partial_{\psi,q,p}}.
\end{equation}
of the restriction of $(\Id \rtimes T_t)_{t \geq 0}$ on $\ovl{\Ran \partial_{\psi,q,p}}$.

\begin{lemma}
\label{lem-Bp-is-sectorial-Fourier}
Let $G$ be a discrete group such that $\Gamma_q(H) \rtimes_\alpha G$ has $\QWEP$. Suppose $1 < p < \infty$ and $-1 \leq q \leq 1$. The operator $B_p$ is injective and sectorial on $\ovl{\Ran \partial_{\psi,q,p}}$.
\end{lemma}

\begin{proof}
The operator $B_p$ is sectorial of type $\frac{\pi}{2}$ by e.g. \cite[p.~25]{JMX}.

For the injectivity, we note that once $B_p$ is known to be sectorial on a reflexive space, we have the projection \cite[(10.5) and p.~367]{HvNVW2} onto the kernel of $B_p$ essentially given by the strong limit 
$$
P
= \lim_{\lambda \to 0^+} \lambda  (\lambda + \Id \rtimes A_p)^{-1}|_{\ovl{\Ran \partial_{\psi,q,p}}}.
$$ 
It is easy to check that $P(y) = 0$ for any $y \in \partial_{\psi,q,p}(\P_G)$. Indeed, by linearity, we can assume that $y = \partial_{\psi,q}(\lambda_s) \ov{\eqref{def-partial-psi}}{=} s_q(b_\psi(s)) \rtimes \lambda_s$ for some $s \in G$ such that in addition $b_\psi(s) \not=0$. We claim that
\begin{equation}
\label{equ-claim-projection-null-space-Fourier}
 (\lambda + \Id \rtimes A_p)^{-1}|_{\ovl{\Ran \partial_{\psi,q,p}}} \big(s_q(b_\psi(s)) \rtimes \lambda_s \big)
= \frac{1}{\lambda + \norm{b_\psi(s)}^2} s_q(b_\psi(s)) \rtimes \lambda_s.
\end{equation}
To see \eqref{equ-claim-projection-null-space-Fourier}, it suffices to calculate
\begin{align*}
\MoveEqLeft
(\lambda + \Id \rtimes A_p)\bigg(\frac{1}{\lambda + \norm{b_\psi(s)}^2} s_q(b_\psi(s)) \rtimes \lambda_s\bigg) \\
&= \frac{1}{\lambda + \norm{b_\psi(s)}^2}\big(\lambda  s_q(b_\psi(s)) \rtimes \lambda_s + \norm{b_\psi(s)}^2 s_q(b_\psi(s)) \rtimes \lambda_s\big) 
= s_q(b_\psi(s)) \rtimes \lambda_s.            
\end{align*} 
Now $P(y)= 0$ follows from $\frac{\lambda}{\lambda + \norm{b_\psi(s)}^2} \to \frac{0}{\norm{b_\psi(s)}^2} = 0$ when $\lambda \to 0^+$ since $\norm{b_\psi(s)}^2 \neq 0$. Since $\P_G$ is a core for $\partial_{\psi,q,p}$, it is not difficult to see that $P(y) = 0$ also for $y \in \Ran \partial_{\psi,q,p}$. Finally by continuity of $P$, we deduce $P= 0$. Thus, $B_p$ is injective.
\end{proof}

Recall that if $G$ is a weakly amenable discrete group, then by the considerations at the end of Subsection \ref{subsubsec-Fourier-mult-crossed-product}, there exists a net $(\Id_{\L^p(\Gamma_q(H))} \rtimes M_{\varphi_j})$ of finitely supported crossed product Fourier multipliers converging strongly to the identity. We also recall from \eqref{equ-P-rtimes-G} that $\P_{p,\rtimes,G}$ is the span of the $x \rtimes \lambda_s$'s, where $x \in \L^p(\VN(G))$ and $s \in G$.

We will use the following result.
 
\begin{lemma}
\label{Petit-lemma}
Let $G$ be a weakly amenable discrete group such that $\Gamma_q(H) \rtimes_\alpha G$ has $\QWEP$. If $z \in \P_{p,\rtimes,G} \cap \ovl{\Ran \partial_{\psi,q,p}}$ then $z \in \partial_{\psi,q,p}(\P_G)$.
\end{lemma}

\begin{proof}
Let $(z_n)$ be a sequence in $\dom \partial_{\psi,q,p}$ such that $\partial_{\psi,q,p}(z_n) \to z$. Let $M_{\varphi}$ be a compactly supported completely bounded Fourier multiplier such that $(\Id \rtimes M_{\varphi})(z) = z$. Then 
\begin{align*}
\MoveEqLeft
\norm{z - \partial_{\psi,q,p}M_{\varphi}(z_n)}_{\L^p} 
\ov{\eqref{commute-troncature}}{=}  \norm{z - (\Id \rtimes M_{\varphi}) \partial_{\psi,q,p}(z_n)}_{\L^p} \\
& \leq \norm{z - (\Id \rtimes M_{\varphi})(z)}_{\L^p} + \norm{(\Id \rtimes M_{\varphi})(z - \partial_{\psi,q,p}(z_n))}_{\L^p}  \\
& \leq \norm{z - (\Id \rtimes M_{\varphi})(z)}_{\L^p}+\norm{z - \partial_{\psi,q,p}(z_n)}_{\L^p}  
\xra[n \to +\infty]{} 0.
\end{align*}
Therefore, the sequence $(\partial_{\psi,q,p} M_{\varphi}(z_n))$ is convergent to $z$ for $n \to \infty$. Write $M_{\varphi}(z_n) = \sum_{s \in F} z_{n,s} \lambda_s$ with some $z_{n,s} \in \C$ and a finite subset $F$ of $G$. Since for any $s \in G$, the Fourier multiplier $M_s$ associated with the symbol $t \mapsto \delta_s(t)$ is completely bounded, $(\Id \rtimes M_s) \partial_{\psi,q,p} M_{\varphi}(z_n) = z_{n,s} s_q(b_\psi(s)) \rtimes \lambda_s$ is convergent for $n \to \infty$. Thus, either $b_\psi(s) = 0$ or $(z_{n,s})$ is a convergent scalar sequence with limit, say, $z^{s} \in \C$. We infer that 
$$
\partial_{\psi,q,p}M_{\varphi}(z_n) 
=\sum_{s \in F} z_{n,s} s_q(b_\psi(s)) \rtimes \lambda_s
\xra[n \to +\infty]{} \sum_{s \in F} z^{s} s_q(b_\psi(s)) \rtimes \lambda_s 
\ov{\eqref{def-partial-psi}}{=} \partial_{\psi,q,p}\bigg(\sum_{s \in F} z^{s} \lambda_s\bigg).
$$ 
We conclude that $z = \partial_{\psi,q,p}(\sum_{s \in F} z^{s} \lambda_s)$ belongs to $\partial_{\psi,q,p}(\P_G)$.
\end{proof}

\begin{prop}
\label{Prop-core-1-group}
Let $G$ be a discrete group such that $\Gamma_q(H) \rtimes_\alpha G$ has $\QWEP$. In the first three points below, assume in addition that $G$ is weakly amenable, so that the approximating Fourier multipliers $M_{\varphi_j}$ exist. Suppose $1 < p < \infty$ and $-1 \leq q \leq 1$.
\begin{enumerate}
\item $\P_{p,\rtimes,G}$ is a core of $(\partial_{\psi,q,p^*})^*$. Furthermore, if $y \in \dom(\partial_{\psi,q,p^*})^*$, then $(\Id \rtimes M_{\varphi_j})(y)$ belongs to $\dom (\partial_{\psi,q,p^*})^*$ for any $j$ and
\begin{equation}
\label{Commute-troncatures-duals}
(\partial_{\psi,q,p^*})^*(\Id \rtimes M_{\varphi_j})(y) 
= M_{\varphi_j} (\partial_{\psi,q,p^*})^*(y).
\end{equation}

\item $\partial_{\psi,q,p}(\P_G)$ is a core of $(\partial_{\psi,q,p^*})^*|\ovl{\Ran \partial_{\psi,q,p}}$.

\item $\partial_{\psi,q,p}(\P_G)$ is equally a core of $\partial_{\psi,q,p}(\partial_{\psi,q,p^*})^*|\ovl{\Ran \partial_{\psi,q,p}}$.

\item $\partial_{\psi,q,p}(\P_G)$ is equally a core of $B_p$.
\end{enumerate}
\end{prop}

\begin{proof}
1. Let $y \in \dom(\partial_{\psi,q,p^*})^*$. Then $(\Id_{\L^p(\Gamma_q(H))} \rtimes M_{\varphi_j})(y)$ belongs to $\P_{p,\rtimes,G}$. 
It remains to show that $(\Id_{\L^p(\Gamma_q(H))} \rtimes M_{\varphi_j})(y)$ converges to $y$ in the graph norm. Recall that $(\Id_{\L^p(\Gamma_q(H))} \rtimes M_{\varphi_j})(y)$ converges to $y$ in $\L^p(\Gamma_q(H) \rtimes_\alpha G)$ according to the assumptions. For any $s \in G$, we have
\begin{align*}
\MoveEqLeft
\big\langle (\partial_{\psi,q,p^*})^* (\Id \rtimes M_{\varphi_j})(y), \lambda_s \big\rangle 
\ov{\eqref{crochet-duality}}{=} \big\langle (\Id \rtimes M_{\varphi_j})(y), \partial_{\psi,q,p^*}(\lambda_s) \big\rangle 
=\big\langle y , (\Id \rtimes M_{\ovl{\varphi_j}})(\partial_{\psi,q,p^*}(\lambda_s)) \big\rangle \\
&\ov{\eqref{commute-troncature}}{=} \big\langle y ,\partial_{\psi,q,p^*} M_{\ovl{\varphi_j}}(\lambda_s) \big\rangle 
=\big\langle (\partial_{\psi,q,p^*})^*(y) , M_{\ovl{\varphi_j}}(\lambda_s) \big\rangle 
=\big\langle M_{\varphi_j} (\partial_{\psi,q,p^*})^*(y) , \lambda_s \big\rangle.
\end{align*}
By linearity and density, we deduce the relation \eqref{Commute-troncatures-duals}, i.e.
$$
(\partial_{\psi,q,p^*})^*(\Id \rtimes M_{\varphi_j})(y) 
= M_{\varphi_j} (\partial_{\psi,q,p^*})^*(y)
$$ 
which converges to $(\partial_{\psi,q,p^*})^*(y)$ in $\L^p(\VN(G))$.

2. Let $y \in \dom (\partial_{\psi,q,p^*})^*|\ovl{\Ran \partial_{\psi,q,p}}$ and $\epsi > 0$. By the proof of the first point, we already know that for some $j$ large enough, there is $z \ov{\mathrm{def}}{=} (\Id \rtimes M_{\varphi_j})(y)$ belonging to $\P_{p,\rtimes,G}$ such that 
\[ 
\norm{y-z}_{\L^p} < \epsi
\quad \text{and} \quad
\norm{(\partial_{\psi,q,p^*})^*(y) - (\partial_{\psi,q,p^*})^*(z)}_{\L^p} 
< \epsi. 
\]
We claim that $z$ belongs to $\ovl{\Ran \partial_{\psi,q,p}}$. Indeed, for $\delta > 0$, since $y \in \ovl{\Ran \partial_{\psi,q,p}}$, there is some $a \in \dom \partial_{\psi,q,p}$ such that $\norm{y - \partial_{\psi,q,p}(a)}_{\L^p} < \delta$. But then $\norm{z - (\Id \rtimes M_{\varphi_j}) \partial_{\psi,q,p}(a)}_{\L^p} = \norm{(\Id \rtimes M_{\varphi_j})(y- \partial_{\psi,q,p}(a))}_{\L^p} \lesssim \norm{y - \partial_{\psi,q,p}(a)}_{\L^p} < \delta$ and $(\Id \rtimes M_{\varphi_j}) \partial_{\psi,q,p}(a) \ov{\eqref{commute-troncature}}{=} \partial_{\psi,q,p} M_{\varphi_j}(a)$ belongs again to $\Ran \partial_{\psi,q,p}$. Letting $\delta \to 0$, it follows that $z$ belongs to $\ovl{\Ran \partial_{\psi,q,p}}$. 

By Lemma \ref{Petit-lemma}, we conclude that $z$ belongs to $\partial_{\psi,q,p}(\P_G)$. 

3. Take $y \in \dom \partial_{\psi,q,p} (\partial_{\psi,q,p^*})^*|\ovl{\Ran \partial_{\psi,q,p}}$, that is, $y \in \dom(\partial_{\psi,q,p^*})^* \cap \ovl{\Ran \partial_{\psi,q,p}}$ such that $(\partial_{\psi,q,p^*})^*(y) \in \dom \partial_{\psi,q,p}$. We have
$(\Id \rtimes M_{\varphi_j})(y) \to y$ and $(\Id \rtimes M_{\varphi_j})(y) \in \partial_{\psi,q,p}(\P_G)$ by the proof of part 2. Moreover, we have
\begin{align*}
\MoveEqLeft
\partial_{\psi,q,p} (\partial_{\psi,q,p^*})^*( \Id \rtimes M_{\varphi_j})(y)           
\ov{\eqref{Commute-troncatures-duals}}{=} \partial_{\psi,q,p} M_{\varphi_j}(\partial_{\psi,q,p^*})^*(y) 
\ov{\eqref{commute-troncature}}{=} (\Id \rtimes M_{\varphi_j})\partial_{\psi,q,p} (\partial_{\psi,q,p^*})^*(y) \\
& \xra[\ j \ ]{} \partial_{\psi,q,p} (\partial_{\psi,q,p^*})^*(y).
\end{align*}   

4. Note that $\partial_{\psi,q,p}(\P_G)$ is a dense subspace of $\ovl{\Ran \partial_{\psi,q,p}}$ which is clearly a subspace of $\dom B_p$ and invariant under each operator $(\Id \rtimes T_t)|\ovl{\Ran \partial_{\psi,q,p}}$ by Lemma \ref{lem-T_t-et-derivations-Fourier}. By Lemma \ref{Lemma-core-semigroup}, we deduce that $\partial_{\psi,q,p}(\P_G)$ is a core of $(\Id \rtimes A_p)|\ovl{\Ran \partial_{\psi,q,p}} \ov{\eqref{Def-Bp-Fourier}}{=} B_p$.
\end{proof}

\begin{prop}
\label{prop-commuting-deriv-Fourier}
Let $G$ be a weakly amenable discrete group such that $\Gamma_q(H) \rtimes_\alpha G$ has $\QWEP$. Suppose $1 < p< \infty$ and $-1 \leq q \leq 1$.
\begin{enumerate}
\item For any $s> 0$, the operator $\big(\Id_{}+sA_p\big)^{-1}(\partial_{\psi,q,p^*})^*$
induces a bounded operator on $\ovl{\Ran \partial_{\psi,q,p}}$.

\item For any $s \geq 0$ and any $y \in \ovl{\Ran \partial_{\psi,q,p}} \cap \dom(\partial_{\psi,q,p^*})^*$, the element $\big(\Id \rtimes (\Id+sA_p)^{-1}\big)(y)$ belongs to $\dom (\partial_{\psi,q,p^*})^*$ and
\begin{equation}
\label{commuting-deriv2-Fourier}
\big(\Id+sA_p\big)^{-1}(\partial_{\psi,q,p^*})^*(y) 
=(\partial_{\psi,q,p^*})^*\big(\Id \rtimes (\Id+sA_p)^{-1}\big)(y).
\end{equation}

\item For any $t \geq 0$ and any $y \in \ovl{\Ran \partial_{\psi,q,p}} \cap \dom(\partial_{\psi,q,p^*})^*$, the element $(\Id \rtimes T_t)(y)$ belongs to $\dom ( \partial_{\psi,q,p^*})^*$ and
$$
T_t(\partial_{\psi,q,p^*})^*(y) 
=(\partial_{\psi,q,p^*})^*\big(\Id \rtimes T_t\big)(y).
$$
\end{enumerate}
\end{prop}

\begin{proof}
1. Note that $(\Id + s A_p)^{-1}(\partial_{\psi,q,p^*})^* \ov{\eqref{Adjoint-product-unbounded}}{\subseteq} \big(\partial_{\psi,q,p^*}(\Id+sA_{p^*})^{-1}\big)^*$. By Proposition \ref{Prop-R-gradient-bounds-Fourier}, the operator $\big(\partial_{\psi,q,p^*}(\Id+sA_{p^*})^{-1}\big)^*$ is bounded. By Lemma \ref{Lemma-rest-closed-dense-Fourier}, the subspace $\partial_{\psi,q,p}(\P_G)$ of $\dom (\partial_{\psi,q,p^*})^*$ is dense in $\ovl{\Ran \partial_{\psi,q,p}}$. Now, the conclusion is immediate.

2. By Proposition \ref{Prop-lien-generateur-partial-Fourier}, for any $x \in \dom A_p$ we have $x \in \dom \partial_{\psi,q,p}$ and $\partial_{\psi,q,p}(x) \in \dom (\partial_{\psi,q,p^*})^*$. Moreover, for all $t > 0$ we have 
\begin{align*}
\MoveEqLeft
T_t (\partial_{\psi,q,p^*})^* \partial_{\psi,q,p}(x) 
\ov{\eqref{eq-lien-generateur-partial-Fourier}}{=} T_tA_p(x)
\ov{\eqref{A-et-Tt-commute}}{=} A_p T_t(x)
\ov{\eqref{eq-lien-generateur-partial-Fourier}}{=} (\partial_{\psi,q,p^*})^* \partial_{\psi,q,p} T_t(x) \\
&\ov{\eqref{T_t-et-derivations-Fourier}}{=} (\partial_{\psi,q,p^*})^* (\Id \rtimes T_t) \partial_{\psi,q,p}(x).
\end{align*}   
By taking Laplace transforms with \eqref{Resolvent-Laplace} and using the closedness of $(\partial_{\psi,q,p^*})^*$, we deduce that the element $(\Id + s \Id \rtimes A_p)^{-1} \partial_{\psi,q,p}(x)$ belongs to $\dom (\partial_{\psi,q,p^*})^*$ for any $s \geq 0$ and that
\begin{equation}
\label{eq:sec-1-Fourier}
(\Id + s A_p)^{-1} (\partial_{\psi,q,p^*})^* \partial_{\psi,q,p}(x)
=(\partial_{\psi,q,p^*})^* (\Id \rtimes (\Id +s A_p))^{-1} \partial_{\psi,q,p}(x).
\end{equation}

Let $y \in \ovl{\Ran \partial_{\psi,q,p}} \cap \dom(\partial_{\psi,q,p^*})^*$. Then according to Lemma \ref{Lemma-rest-closed-dense-Fourier}, there exists a sequence $(x_n)$ of $\P_G$ such that $\partial_{\psi,q,p}(x_n) \to y$. We have $(\Id \rtimes (\Id +sA_p))^{-1}\partial_{\psi,q,p}(x_n) \to (\Id \rtimes (\Id + s A_p))^{-1}(y)$. Since each $x_n$ belongs to $\dom A_p$, using the first point in the passage to the limit we deduce that
\begin{align*}
\MoveEqLeft
(\partial_{\psi,q,p^*})^* (\Id \rtimes (\Id +s A_p))^{-1} \partial_{\psi,q,p}(x_n)   
\ov{\eqref{eq:sec-1-Fourier}}{=} (\Id + s A_p)^{-1} (\partial_{\psi,q,p^*})^* \partial_{\psi,q,p}(x_n) \\
&\xra[n \to +\infty]{} (\Id + s A_p)^{-1} (\partial_{\psi,q,p^*})^*(y).
\end{align*}
Since $(\partial_{\psi,q,p^*})^*$ is closed, we infer by \eqref{Def-operateur-ferme} that $(\Id \rtimes (\Id +sA_p))^{-1}(y)$ belongs to $\dom (\partial_{\psi,q,p^*})^*$ and that
\begin{equation}
\label{Commutation-Fourier}
(\partial_{\psi,q,p^*})^*(\Id \rtimes (\Id +s A_p))^{-1}(y)
=(\Id + s A_p)^{-1} (\partial_{\psi,q,p^*})^*(y).
\end{equation}
Thus \eqref{commuting-deriv2-Fourier} follows.

3. Let $y \in \ovl{\Ran \partial_{\psi,q,p}} \cap \dom(\partial_{\psi,q,p^*})^*$. If $t \geq 0$, note that 
$$
\bigg(\Id \rtimes \bigg(\Id + \frac{t}{n} A_p\bigg)\bigg)^{-n}(y)
\xra[n \to +\infty]{\eqref{Widder}} (\Id \rtimes T_t)(y).
$$
Repeating the above commutation relation \eqref{Commutation-Fourier} together with the observation that $(\Id \rtimes (\Id + s A_p))^{-1}$ maps $\partial_{\psi,q,p}(\P_G)$ into itself hence by continuity $\ovl{\Ran \partial_{\psi,q,p}}$ into itself yields\footnote{\thefootnote. Note that we replace $s$ by $\frac{t}{n}$.} for any integer $n \geq 1$ and any $t \geq 0$
$$
(\partial_{\psi,q,p^*})^*\bigg(\Id \rtimes \bigg(\Id + \frac{t}{n} A_p\bigg)\bigg)^{-n}(y)
=\bigg(\Id + \frac{t}{n} A_p\bigg)^{-n} (\partial_{\psi,q,p^*})^*(y)
\xra[n \to +\infty]{\eqref{Widder}} T_t (\partial_{\psi,q,p^*})^*(y).
$$
Then by the closedness of $(\partial_{\psi,q,p^*})^*$, we deduce that $(\Id \rtimes T_t)(y)$ belongs to $\dom (\partial_{\psi,q,p^*})^*$ and that
$$
(\partial_{\psi,q,p^*})^*(\Id \rtimes T_t)(y) 
= T_t (\partial_{\psi,q,p^*})^*(y).
$$ 
Thus the third point follows. 
\end{proof}

Proposition \ref{Prop-core-1-group} enables us to identify $B_p$ in terms of $\partial_{\psi,q,p}$ and its adjoint.

\begin{prop}
\label{Prop-fundamental-Fourier}
Let $G$ be a weakly amenable discrete group such that $\Gamma_q(H) \rtimes_\alpha G$ has $\QWEP$. Let $1 < p < \infty$ and $-1 \leq q \leq 1$. As unbounded operators, we have
\begin{equation}
\label{relation-dur-fund-Fourier}
B_p
=\partial_{\psi,q,p}(\partial_{\psi,q,p^*})^*|\ovl{\Ran \partial_{\psi,q,p}}.
\end{equation}
\end{prop}

\begin{proof}
For any $s \in G$, we have
\begin{align*}
\MoveEqLeft
\partial_{\psi,q,p}(\partial_{\psi,q,p^*})^*  \partial_{\psi,q,p}(\lambda_s)
\ov{\eqref{eq-lien-generateur-partial-Fourier}}{=}\partial_{\psi,q,p}A_p(\lambda_s)          
=\norm{b_\psi(s)}_H^2 \partial_{\psi,q,p}(\lambda_s) \\
&\ov{\eqref{def-partial-psi}}{=}\norm{b_\psi(s)}_H^2 s_q(b_\psi(s)) \rtimes \lambda_s
=(\Id_{} \rtimes A_p)\big(s_q(b_\psi(s)) \rtimes \lambda_s\big) 
\ov{\eqref{def-partial-psi}}{=} (\Id_{} \rtimes A_p)\big(\partial_{\psi,q,p}(\lambda_s)\big).
\end{align*}  
We deduce that the operators $\partial_{\psi,q,p}(\partial_{\psi,q,p^*})^*|\ovl{\Ran \partial_{\psi,q,p}}$ and $(\Id_{\L^p(\Gamma_q(H))} \rtimes A_p)$ coincide on $\partial_{\psi,q,p}(\P_G)$. By Proposition \ref{Prop-core-1-group}, $\partial_{\psi,q,p}(\P_G)$ is a core for each operator. We conclude that they are equal.
\end{proof}

In the proof of Theorem \ref{crossed-Ap-functional-calculus}, we shall use the following folklore lemma, see e.g. \cite[Proposition 2.13]{LM1}.

\begin{lemma}
\label{lem-technical-HI-Laplace-transform}
Let $(T_t)_{t \geq 0}$ be a bounded strongly continuous semigroup on a Banach space $X$. Let $-A$ be its infinitesimal generator. For any $\frac{\pi}{2}<\theta<\pi$ the following are equivalent
\begin{enumerate}
	\item $A$ admits a bounded $\H^\infty(\Sigma_\theta)$ functional calculus. 
	\item There exists a constant $C>0$ such that for any $b \in \L^1(\R^+)$ whose Laplace transform $\scr{L}(b)$ belongs to $\H_0^\infty(\Sigma_\theta)$, we have
	\begin{equation}
\label{}
\norm{\int_{0}^{+\infty} b(t) T_t \d t}_{X \to X}
\leq C\norm{\scr{L}(b)}_{\H^\infty(\Sigma_\theta)}.
\end{equation}
\end{enumerate}
\end{lemma}
Recall that the above integral is defined in the strong operator topology sense.

\begin{thm}
\label{crossed-Ap-functional-calculus}
Let $G$ be a discrete group. Suppose $1<p<\infty$ and $-1 \leq q \leq 1$. Let $G$ be a weakly amenable discrete group such that $\Gamma_q(H) \rtimes_\alpha G$ has $\QWEP$. The operators $A_p$ and $\Id \rtimes A_p$ have a bounded $\H^\infty(\Sigma_\theta)$ functional calculus of angle $\theta$ for any $\theta>\pi|\frac{1}{p}-\frac{1}{2}|$. 
\end{thm}

\begin{proof}
For the $\QWEP$ property in case of $G$ being a free group, we refer to \cite[Proposition 4.8]{Arh1}.
According to Proposition \ref{prop-Fourier-mult-crossed-product}, $\Id \rtimes T_t$ extends to a (completely) contractive operator on $\L^p(\Gamma_q(H) \rtimes_\alpha G)$.
Moreover, since $\P_{\rtimes,G}$ is dense in $\L^p(\Gamma_q(H) \rtimes_\alpha G)$, the property of $(T_t)_{t \geq 0}$ being a strongly continuous semigroup carries over to $(\Id \rtimes T_t)_{t \geq 0}$. According to \cite[Theorem 4.1]{Arh7} (see also \cite{Arh8}), $A_p$ has a (completely) bounded $\H^\infty(\Sigma_\theta)$ functional calculus on $\L^p(\VN(G))$ for any $\frac{\pi}{2} <\theta < \pi$. By Lemma \ref{lem-technical-HI-Laplace-transform}, for any $b \in \L^1(\R_+)$ such that $\scr{L}(b)$ belongs to $\H_0^\infty(\Sigma_\theta)$ we have
\begin{equation}
\label{Equa-5001}
\norm{ \int_0^\infty b(t)T_t \d t}_{\cb,\L^p(\VN(G)) \to \L^p(\VN(G))} 
\lesssim \norm{\scr{L}(b)}_{\H^\infty(\Sigma_\theta)}.
\end{equation}
Note that for any $\lambda_s \in \P_{G}$, we have 
\begin{align*}
\MoveEqLeft
\bigg(\int_0^\infty b(t) T_t \d t\bigg)(\lambda_s) 
=\int_0^\infty b(t) T_t(\lambda_s) \d t 
=\int_0^\infty b(t) \e^{-t \norm{b_\psi(s)}^2}\lambda_s \d t \\
&=\bigg(\int_0^\infty b(t) \e^{-t \norm{b_\psi(s)}^2} \d t\bigg)(\lambda_s).
\end{align*}
So the map $\int_0^\infty b(t) T_t \d t$ is also a Fourier multiplier. Thus, by Proposition \ref{prop-Fourier-mult-crossed-product}
\begin{align*}
\MoveEqLeft
\norm{\int_0^\infty b(t) \Id \rtimes T_t(x) \d t}_{\L^p(\Gamma_q(H) \rtimes_\alpha G) \to \L^p(\Gamma_q(H) \rtimes_\alpha G)} 
=\norm{\Id \rtimes \bigg(\int_0^\infty b(t)  T_t(x) \d t\bigg)}_{\L^p \to \L^p} \\
&\ov{\eqref{twisted-multiplier-major}}{\leq} \norm{ \int_0^\infty b(t)T_t(x) \d t}_{\cb,\L^p(\VN(G)) \to \L^p(\VN(G))} 
\ov{\eqref{Equa-5001}}{\lesssim} \norm{\scr{L}(b)}_{\H^\infty(\Sigma_\theta)}.           
\end{align*}
By Lemma \ref{lem-technical-HI-Laplace-transform}, $\Id \rtimes A_p$ admits a bounded $\H^\infty(\Sigma_\theta)$ functional calculus for any $\theta >\frac{\pi}{2}$. Now, we reduce the angle and conclude with \cite[ Proposition 5.8]{JMX}.
\end{proof}

Suppose $1 < p <\infty$. 
We introduce the unbounded operator
\begin{equation}
\label{Def-Dirac-operator-Fourier}
D_{\psi,q,p} 
\ov{\mathrm{def}}{=}\begin{bmatrix} 
0 & (\partial_{\psi,q,p^*})^* \\ 
\partial_{\psi,q,p}& 0 
\end{bmatrix}
\end{equation}
on the Banach space $\L^p(\VN(G)) \oplus_p \ovl{\Ran\, \partial_{\psi,q,p}}$ defined by
\begin{equation}
\label{Def-D-psi}
D_{\psi,q,p}(x,y)
\ov{\mathrm{def}}{=}
\big((\partial_{\psi,q,p^*})^*(y),
\partial_{\psi,q,p}(x)\big), \quad x \in \dom \partial_{\psi,q,p},\  y \in \dom (\partial_{\psi,q,p^*})^* \cap \ovl{\Ran \partial_{\psi,q,p}}.
\end{equation}
We call it the Hodge-Dirac operator of the semigroup. This operator is a closed operator and can be seen as a differential square root of the generator of the semigroup $(T_{t,p})_{t \geq 0}$ since we have Proposition \ref{Prop-carre-Dirac-Fourier}.

\begin{thm}
\label{Th-D-R-bisectorial-Fourier}
Let $G$ be a weakly amenable discrete group such that $\Gamma_q(H) \rtimes_\alpha G$ has $\QWEP$. Suppose $1<p<\infty$ and $-1 \leq q \leq 1$. The Hodge-Dirac operator $D_{\psi,q,p}$ is R-bisectorial on $\L^p(\VN(G)) \oplus_p \ovl{\Ran \partial_{\psi,q,p}}$. 
\end{thm}

\begin{proof}
We will start by showing that the set $\{\i t : t \in \R, t \not=0\}$ is contained in the resolvent set of $D_{\psi,q,p}$. We will do this by showing that $\Id-\i tD_{\psi,q,p}$ has a two-sided bounded inverse $(\Id-\i tD_{\psi,q,p})^{-1}$ given by
\begin{equation}
\label{Resolvent-Fourier}
\begin{bmatrix} 
(\Id_{}+t^2A_p)^{-1} & \i t(\Id_{}+t^2A_p)^{-1}(\partial_{\psi,q,p^*})^* \\ 
\i t\partial_{\psi,q,p}(\Id_{}+t^2A_p)^{-1} & \Id_{} \rtimes (\Id+t^2A_p)^{-1} 
\end{bmatrix} 
\co \L^p \oplus_p \ovl{\Ran \partial_{\psi,q,p}} \to \L^p \oplus_p \ovl{\Ran \partial_{\psi,q,p}}.
\end{equation}
By Proposition \ref{Prop-R-gradient-bounds-Fourier} and since the operators $A_p$ and $\Id_{\L^p} \rtimes A_p$ satisfy the property \eqref{Def-R-sectoriel} of $R$-sectoriality (see Theorem \ref{crossed-Ap-functional-calculus}), the four entries are bounded. It only remains to check that this matrix defines a two-sided inverse of $\Id-\i t D_{\psi,q,p}$. We have the following equalities of operators acting on $\dom D_{\psi,q,p}$.
\begin{align*}
\MoveEqLeft
\begin{bmatrix} 
(\Id_{}+t^2A_p)^{-1} & \i t(\Id_{\L^p}+t^2A_p)^{-1}(\partial_{\psi,q,p^*})^* \\ 
\i t\partial_{\psi,q,p}(\Id_{}+t^2A_p)^{-1} & \Id_{} \rtimes (\Id+t^2A_p)^{-1} 
\end{bmatrix}(\Id-\i tD_{\psi,q,p})   \\       
&\ov{\eqref{Def-Dirac-operator-Fourier}}{=}
\begin{bmatrix} 
(\Id+t^2A_p)^{-1} & \i t(\Id_{}+t^2A_p)^{-1}(\partial_{\psi,q,p^*})^* \\ 
\i t\partial_{\psi,q,p}(\Id_{}+t^2A_p)^{-1} & \Id_{} \rtimes (\Id+t^2A_p)^{-1} 
\end{bmatrix}
\begin{bmatrix} 
\Id_{\L^p} & -\i t(\partial_{\psi,q,p^*})^* \\ 
-\i t\partial_{\psi,q,p} &\Id_{\ovl{\Ran \partial_{\psi,q,p}}}
\end{bmatrix}\\
&=\left[\begin{matrix} 
(\Id+t^2A_p)^{-1}+t^2(\Id+t^2A_p)^{-1}(\partial_{\psi,q,p^*})^*\partial_{\psi,q,p}\\ 
\i t\partial_{\psi,q,p}(\Id+t^2A_p)^{-1}-\i t(\Id_{} \rtimes (\Id+t^2A_p)^{-1})\partial_{\psi,q,p}
\end{matrix}\right.\\
&\qquad \qquad \qquad \qquad \qquad \qquad \left.\begin{matrix} 
-\i t(\Id+t^2A_p)^{-1}(\partial_{\psi,q,p^*})^*+\i t(\Id+t^2A_p)^{-1}(\partial_{\psi,q,p^*})^* \\ 
t^2\partial_{\psi,q,p}(\Id+t^2A_p)^{-1}(\partial_{\psi,q,p^*})^*+\Id_{} \rtimes(\Id+t^2A_p)^{-1}
\end{matrix}\right]\\
&\ov{\eqref{eq-lien-generateur-partial-Fourier} \eqref{commuting-deriv-Fourier} \eqref{commuting-deriv2-Fourier}}{=}
\left[\begin{matrix} 
(\Id+t^2A_p)^{-1}+t^2(\Id+t^2A_p)^{-1}A_p  \\ 
\i t\partial_{\psi,q,p}(\Id+t^2A_p)^{-1}-\i t\partial_{\psi,q,p}\big(\Id+t^2A_p\big)^{-1} 
\end{matrix}\right. \\
&\qquad \qquad \qquad \qquad \qquad \qquad \qquad \qquad\left.\begin{matrix} 0 \\
(t^2\partial_{\psi,q,p}(\partial_{\psi,q,p^*})^*+\Id_{})(\Id_{} \rtimes (\Id+t^2A_p)^{-1})
\end{matrix}\right] \\
&=\ov{\eqref{relation-dur-fund-Fourier}\eqref{Def-Bp-Fourier}}{=}  
\begin{bmatrix} 
\Id_{\L^p(\VN(G))} & 0 \\ 
0 & \Id_{\ovl{\Ran \partial_{\psi,q,p}}}
\end{bmatrix}
\end{align*} 
and similarly
\begin{align*}
\MoveEqLeft
(\Id-\i tD_{\psi,p}) 
\begin{bmatrix} 
(\Id_{}+t^2A_p)^{-1} & \i t(\Id_{}+t^2A_p)^{-1}(\partial_{\psi,q,p^*})^* \\ 
\i t\partial_{\psi,q,p}(\Id_{}+t^2A_p)^{-1} & \Id_{} \rtimes (\Id_{}+t^2A_p)^{-1} 
\end{bmatrix}  \\          
&=\begin{bmatrix} 
\Id_{} & -\i t(\partial_{\psi,q,p^*})^* \\ 
-\i t\partial_{\psi,q,p} &\Id_{\ovl{\Ran \partial_{\psi,q,p}}}
\end{bmatrix}
\begin{bmatrix} 
(\Id_{}+t^2A_p)^{-1} & \i t(\Id_{}+t^2A_p)^{-1}(\partial_{\psi,q,p^*})^* \\ 
\i t\partial_{\psi,q,p}(\Id_{}+t^2A_p)^{-1} & \Id_{} \rtimes (\Id_{}+t^2A_p)^{-1} 
\end{bmatrix} 
\\
&=\left[\begin{matrix} 
(\Id_{}+t^2A_p)^{-1}+ t^2(\partial_{\psi,q,p^*})^*\partial_{\psi,q,p}(\Id_{}+t^2A_p)^{-1} &  \\ 
-\i t\partial_{\psi,q,p}(\Id_{}+t^2A_p)^{-1}+\i t\partial_{\psi,q,p}(\Id_{}+t^2A_p)^{-1} &
\end{matrix}\right.\\
&\qquad \qquad \qquad \qquad \qquad\left.
\begin{matrix}
\i t(\Id_{}+t^2A_p)^{-1}(\partial_{\psi,q,p^*})^*-\i t(\partial_{\psi,q,p^*})^*\big(\Id_{} \rtimes (\Id_{}+t^2A_p)^{-1}\big)\\
t^2\partial_{\psi,q,p}(\Id_{}+t^2A_p)^{-1}(\partial_{\psi,q,p^*})^*+\Id_{} \rtimes (\Id_{}+t^2A_p)^{-1}
\end{matrix}\right]\\
&= \begin{bmatrix} 
\Id_{\L^p} & 0 \\ 
0 & \Id_{\ovl{\Ran \partial_{\psi,q,p}}}
\end{bmatrix}.
\end{align*} 
It remains to show that the set $\{\i t(\i t-D_{\psi,q,p})^{-1} : t \not=0\}=\{(\Id -\i t D_{\psi,q,p})^{-1} : t \not=0\}$ is $R$-bounded. For this, observe that the diagonal entries of \eqref{Resolvent-Fourier} are $R$-bounded by the $R$-sectoriality of $A_p$ and $\Id_{\L^p} \rtimes A_p$. The $R$-boundedness of the other entries follows from the $R$-gradient bounds of Proposition \ref{Prop-R-gradient-bounds-Fourier}. Since a set of operator matrices is $R$-bounded precisely when each entry is $R$-bounded, we conclude that \eqref{Def-R-bisectoriel} is satisfied, i.e. that $D_{\psi,q,p}$ is $R$-bisectorial.
\end{proof}

\begin{prop}
\label{Prop-carre-Dirac-Fourier}
Let $G$ be a weakly amenable discrete group such that $\Gamma_q(H) \rtimes_\alpha G$ has $\QWEP$. Suppose $1<p<\infty$ and $-1 \leq q \leq 1$. As densely defined closed operators on $\L^p(\VN(G)) \oplus_p \ovl{\Ran\, \partial_{\psi,q,p}}$, we have 
\begin{equation}
\label{D-alpha-carre-egal-Fourier}
D_{\psi,q,p}^2 
=\begin{bmatrix} 
A_p & 0 \\ 
0 & (\Id_{\L^p(\Gamma_{q}(H))} \rtimes A_p)|\ovl{\Ran \partial_{\psi,q,p}}
\end{bmatrix}.	
\end{equation}
\end{prop}

\begin{proof}
By Proposition \ref{Prop-fundamental-Fourier}, we have
\begin{align*}
\MoveEqLeft
C_p
\ov{\mathrm{def}}{=}\begin{bmatrix} 
A_p & 0 \\ 
0 & (\Id_{\L^p(\Gamma_q(H))} \rtimes A_p)|\ovl{\Ran \partial_{\psi,q,p}}
\end{bmatrix} \\
&\ov{\eqref{eq-lien-generateur-partial-Fourier}\eqref{relation-dur-fund-Fourier}}{=} 
\begin{bmatrix} 
(\partial_{\psi,q,p^*})^*\partial_{\psi,q,p} & 0 \\ 
0 & \partial_{\psi,q,p}(\partial_{\psi,q,p^*})^*|\ovl{\Ran \partial_{\psi,q,p}}\,
\end{bmatrix} 
=\begin{bmatrix} 
0 & (\partial_{\psi,q,p^*})^*|\ovl{\Ran \partial_{\psi,q,p}}\, \\ 
\partial_{\psi,q,p}& 0 
\end{bmatrix}^2 \\
&\ov{\eqref{Def-Dirac-operator-Fourier}}{=} D_{\psi,q,p}^2.            
\end{align*}
\end{proof}

Now, we can state the following main result of this subsection.

\begin{thm}
\label{Th-functional-calculus-bisector-Fourier}
Suppose $1<p<\infty$ and $-1 \leq q \leq 1$. Let $G$ be a weakly amenable discrete group such that $\Gamma_q(H) \rtimes_\alpha G$ has $\QWEP$. The Hodge-Dirac operator $D_{\psi,q,p}$ is R-bisectorial on $\L^p(\VN(G)) \oplus_p \ovl{\Ran \partial_{\psi,q,p}}$ and admits a bounded $\HI(\Sigma^\pm_\omega)$ functional calculus on a bisector. 
\end{thm}

\begin{proof}
By Theorem \ref{crossed-Ap-functional-calculus}, the operator $
D_{\psi,q,p}^2 
\ov{\eqref{D-alpha-carre-egal-Fourier}}{=}\begin{bmatrix} 
A_p & 0 \\ 
0 & B_p
\end{bmatrix}$ has a bounded $\HI$ functional calculus of angle $2 \omega < \frac{\pi}{2}$. Since $D_{\psi,q,p}$ is $R$-bisectorial by Theorem \ref{Th-D-R-bisectorial-Fourier}, we deduce by Proposition \ref{Prop-liens-differents-calculs-fonctionnels} that the operator $D_{\psi,q,p}$ has a bounded $\HI(\Sigma^\pm_\omega)$ functional calculus on a bisector.
\end{proof}

\begin{remark} \normalfont
\label{remark-sgn-Fourier}
The boundedness of the $\HI$ functional calculus of the operator $D_{\psi,q,p}$ implies the boundedness of the Riesz transforms and this result may be thought of as a strengthening of the equivalence \eqref{Equivalence-square-root-domaine-Schur-sgrp}. Indeed, consider the function $\sgn \in \H^\infty(\Sigma_\omega)$ defined by $\sgn(z) \ov{\mathrm{def}}{=} 1_{\Sigma_\omega^+}(z)-1_{\Sigma_\omega^-}(z)$. By Theorem \ref{Th-functional-calculus-bisector-Fourier}, the operator $D_{\psi,q,p}$ has a bounded $\H^\infty(\Sigma^\pm_\omega)$ functional calculus on $\L^p(\VN(G)) \oplus_p \ovl{\Ran \partial_{\psi,q,p}}$. Hence the operator $\sgn(D_{\psi,q,p})$ is bounded. This implies that
\begin{equation}
\label{lien-sign-abs-Fourier}
\big(D_{\psi,q,p}^2\big)^{\frac{1}{2}}=\sgn(D_{\psi,q,p}) D_{\psi,q,p}
\quad \text{and} \quad
D_{\psi,q,p}
=\sgn(D_{\psi,q,p})\big(D_{\psi,q,p}^2\big)^{\frac{1}{2}}.
\end{equation}
For any $x \in \dom D_{\psi,q,p} = \dom \big(D_{\psi,q,p}^2\big)^{\frac{1}{2}}$, we deduce that
\begin{align*}
\MoveEqLeft
\bnorm{D_{\psi,q,p}(x)}_{\L^p(\VN(G)) \oplus_p \L^p(\Gamma_q(H) \rtimes_\alpha G)} 
\ov{\eqref{lien-sign-abs-Fourier}}{=} \bnorm{\sgn(D_{\psi,q,p}) (D_{\psi,q,p}^2)^{\frac{1}{2}}(x)}_{\L^p(\VN(G)) \oplus_p \L^p(\Gamma_q(H) \rtimes_\alpha G)} \\
&\lesssim_p \bnorm{(D_{\psi,q,p}^2)^{\frac{1}{2}}(x)}_{\L^p(\VN(G)) \oplus_p \L^p(\Gamma_q(H) \rtimes_\alpha G)}
\end{align*}
and
\begin{align*}
\MoveEqLeft
\bnorm{(D_{\psi,q,p}^2)^{\frac{1}{2}}(x)}_{\L^p(\VN(G)) \oplus_p \L^p(\Gamma_q(H) \rtimes_\alpha G)} 
\ov{\eqref{lien-sign-abs-Fourier}}{=} \bnorm{\sgn(D_{\psi,q,p}) D_{\psi,q,p}(x)}_{\L^p(\VN(G)) \oplus_p \L^p(\Gamma_q(H) \rtimes_\alpha G)} \\
&\lesssim_p \bnorm{D_{\psi,q,p}(x)}_{\L^p(\VN(G)) \oplus_p  \L^p(\Gamma_q(H) \rtimes_{\alpha} G)}.            
\end{align*}
Recall that on $\L^p(\VN(G)) \oplus_p \ovl{\Ran \partial_{\psi,q,p}}$, we have
\begin{equation}
\label{Eq-inter-sans-fin1}
(D_{\psi,q,p}^2)^{\frac{1}{2}}
\ov{\eqref{D-alpha-carre-egal-Fourier}}{=}\begin{bmatrix}
A_p^{\frac{1}{2}} & 0 \\ 
0 & \Id_{\L^p(\Gamma_q(H))} \rtimes A_p^{\frac{1}{2}}|\ovl{\Ran \partial_{\psi,q,p}}
\end{bmatrix}.
\end{equation}
By restricting to elements of the form $(y,0)$ with $y \in \dom A_p^{\frac{1}{2}}$, we obtain the desired result.
\end{remark}

\begin{remark} \normalfont
\label{remark-amplified-Riesz-Fourier}
In a similar way to Remark \ref{remark-sgn-Fourier}, we also obtain for any element $(x,y)$ of $\dom \partial_{\psi,q,p} \oplus \dom (\partial_{\psi,q,p^*})^*|\ovl{\Ran \partial_{\psi,q,p}}$ that
\begin{align*}
\MoveEqLeft
\bnorm{(D^2_{\psi,q,p})^{\frac12}(x,y)}_p 
\ov{\eqref{Eq-inter-sans-fin1}}{\cong} \bnorm{A_p^{\frac12}(x)}_p + \bnorm{\big(\Id_{\L^p(\Gamma_q(H))} \rtimes A_p^{\frac12}\big)(y)}_p 
\cong \norm{\partial_{\psi,q,p}(x)}_p + \norm{(\partial_{\psi,q,p^*})^*(y)}_p \\
&\ov{\eqref{Def-Dirac-operator-Fourier}}{\cong} \norm{D_{\psi,q,p}(x,y)}_p.
\end{align*}
Similarly, we have
\begin{equation}
\label{equ-amplified-Riesz-Fourier}
\bnorm{\big(\Id_{\L^p(\Gamma_q(H))} \rtimes A_p^{\frac12}\big)(y)}_p 
\cong  \norm{(\partial_{\psi,q,p^*})^*(y)}_p, \quad y \in \dom(\partial_{\psi,q,p^*})^* \cap \ovl{\Ran \partial_{\psi,q,p}}.
\end{equation}
\end{remark}

\begin{prop}
\label{Prop-liens-ranges-group}
Let $G$ be a weakly amenable discrete group such that $\Gamma_q(H) \rtimes_\alpha G$ has $\QWEP$. Suppose $1 < p < \infty$ and $-1 \leq q \leq 1$. We have $\ovl{\Ran  A_p}=\ovl{\Ran (\partial_{\psi,q,p^*})^*}$, $\ovl{\Ran B_p}=\ovl{\Ran \partial_{\psi,q,p}}$, $\ker A_p=\ker \partial_{\psi,q,p}$, $\ker B_p=\ker (\partial_{\psi,q,p^*})^*=\{0\}$ and
\begin{equation}
\label{Decompo-sympa-group}
\L^p(\VN(G))
=\ovl{\Ran (\partial_{\psi,q,p^*})^*} \oplus \ker \partial_{\psi,q,p}.
\end{equation}
Here, by $(\partial_{\psi,q,p^*})^*$ we understand its restriction to $\ovl{\Ran \partial_{\psi,q,p}}$.
However, we have $\Ran (\partial_{\psi,q,p^*})^* = \Ran(\partial_{\psi,q,p^*})^*|{\ovl{\Ran \partial_{\psi,q,p}}}$.
\end{prop}

\begin{proof}
By \eqref{Bisec-Ran-Ker}, we have $\ovl{\Ran D_{\psi,q,p}^2}=\ovl{\Ran D_{\psi,q,p}}$ and $\ker D_{\psi,q,p}^2=\ker D_{\psi,q,p}$. It is not difficult to prove the first four equalities using \eqref{D-alpha-carre-egal-Fourier} and \eqref{Def-D-psi}. The last one is a consequence of the definition of $A_p$ and of \cite[p.~361]{HvNVW2}.
\end{proof}

\subsection{Extension to full Hodge-Dirac operator and Hodge decomposition}
\label{Extension-full-Hodge-Fourier}

We keep the standing assumptions of the preceding subsection and thus have a markovian semigroup $(T_t)_{t \geq 0}$ of Fourier multipliers with generator $A_p$, the noncommutative gradient $\partial_{\psi,q,p}$ and its adjoint $(\partial_{\psi,q,p^*})^*$, together with the Hodge-Dirac operator $D_{\psi,q,p}$. We shall now extend the operator $D_{\psi,q,p}$ to a densely defined bisectorial operator $\D_{\psi,q,p}$ on $\L^p(\VN(G)) \oplus \L^p(\Gamma_q(H) \rtimes_\alpha G)$ which will also be bisectorial and will have a $\HI(\Sigma_\omega^\pm)$ functional calculus on a bisector. The key will be Corollary \ref{cor-key-Hodge-decomposition-Fourier} below. We let
\begin{equation}
\label{equ-full-Hodge-Dirac-operator-Fourier}
\D_{\psi,q,p} 
\ov{\mathrm{def}}{=} \begin{bmatrix} 
0 & (\partial_{\psi,q,p^*})^* \\ 
\partial_{\psi,q,p} & 0 
\end{bmatrix}
\end{equation}
along the decomposition $\L^p(\VN(G)) \oplus \L^p(\Gamma_q(H) \rtimes_\alpha G)$, with natural domains for $\partial_{\psi,q,p}$ and $(\partial_{\psi,q,p^*})^*$. Again, except in Lemma \ref{lem-Hodge-decomposition-sum-Fourier}, we need in this subsection approximation properties of $G$, see before Lemma \ref{lem-eigenspace-q-OU-Fourier}.

Consider the sectorial operator $A_p^{\frac12}$ on $\L^p(\VN(G))$. According to \cite[(10.1) p.~361]{HvNVW2}, we have the \textit{topological} direct sum decomposition $\L^p(\VN(G))= \ovl{\Ran A_p^{\frac12}} \oplus \ker A_p^{\frac12}$. We define the operator $R_p \ov{\mathrm{def}}{=} \partial_{\psi,q,p} A_p^{-\frac12} \co \Ran A_p^{\frac12} \to \L^p(\Gamma_q(H) \rtimes_\alpha G)$. According to  the point 4 of Proposition \ref{Prop-derivation-closable-sgrp}, $R_p$ is bounded on $\Ran A_p^{\frac12}$, so extends to a bounded operator on $\ovl{\Ran A_p^{\frac12}}$. We extend it to a bounded operator $R_p \co \L^p(\VN(G)) \to \L^p(\Gamma_q(H) \rtimes_{\alpha} G)$, called Riesz transform,  by putting $R_p|\ker A_p^{\frac12} = 0$ along the above decomposition of $\L^p(\VN(G))$. We equally let $R_{p^*}^* \ov{\mathrm{def}}{=} (R_{p^*})^*$.

\begin{lemma}
\label{lem-Hodge-decomposition-sum-Fourier}
Let $-1 \leq q \leq 1$ and $1 < p < \infty$. Then we have the decomposition
\begin{equation}
\label{Sum-plus-Fourier}
\L^p(\Gamma_q(H) \rtimes_{\alpha} G) 
=\ovl{\Ran \partial_{\psi,q,p}} + \ker (\partial_{\psi,q,p^*})^*.
\end{equation}
\end{lemma}

\begin{proof}
Let $y \in \L^p(\Gamma_q(H) \rtimes_\alpha G)$ be arbitrary. We claim that $y = R_pR_{p^*}^*(y) + (\Id - R_pR_{p^*}^*)(y)$ is the needed decomposition for \eqref{Sum-plus-Fourier}. Note that $R_p$ maps $\Ran A_p^{\frac12}$ into $\Ran \partial_{\psi,q,p}$, so by boundedness, $R_p$ maps $\ovl{\Ran A_p^{\frac12}}$ to $\ovl{\Ran \partial_{\psi,q,p}}$. Thus, we indeed have $R_pR_{p^*}^*(y) \in \ovl{\Ran \partial_{\psi,q,p}}$. Next we claim that for any $z \in \L^p(\VN(G))$ and any $x \in \dom \partial_{\psi,q,p^*}$, we have
\begin{equation}
\label{equ-1-proof-lem-Hodge-decomposition-sum-Fourier}
\big\langle R_p(z), \partial_{\psi,q,p^*}(x)\big\rangle_{\L^p(\Gamma_q(H) \rtimes_\alpha G),\L^{p^*}(\Gamma_q(H) \rtimes_\alpha G)} 
=\big\langle z, A_{p^*}^{\frac12}(x) \big\rangle_{\L^p(\VN(G)),\L^{p^*}(\VN(G))}.
\end{equation}
According to the decomposition $\L^p(\VN(G)) = \ovl{\Ran A_p^{\frac12}} \oplus \ker A_p^{\frac12}$ above, we can write $z = \lim_{n \to +\infty} A_p^{\frac12}(z_n) + z_0$ with $z_n \in \dom A_p^{\frac12}$ and $z_0 \in \ker A_p^{\frac12}$. Then using Lemma \ref{lem-formule-trace-trace-Fourier} in the third equality, we have
\begin{align*}
\MoveEqLeft
\big\langle R_p(z), \partial_{\psi,q,p^*}(x) \big\rangle
=\lim_{n \to +\infty}\big\langle R_p\big(A_p^{\frac12}(z_n) + z_0\big), \partial_{\psi,q,p^*}(x) \big\rangle 
=\lim_{n \to +\infty} \big\langle \partial_{\psi,q,p}(z_n) , \partial_{\psi,q,p^*}(x) \big\rangle \\
&=\lim_{n \to +\infty} \big\langle A_p^{\frac12}(z_n), A_{p^*}^{\frac12} (x) \big\rangle 
=\big\langle z - z_0 , A_{p^*}^{\frac12}(x) \big\rangle 
=\big\langle z, A_{p^*}^{\frac12}(x) \big\rangle-\big\langle z_0 , A_{p^*}^{\frac12}(x) \big\rangle 
=\big\langle z, A_{p^*}^{\frac12}(x) \big\rangle.
\end{align*}
Thus, \eqref{equ-1-proof-lem-Hodge-decomposition-sum-Fourier} is proved. Now, for any $x \in \dom \partial_{\psi,q,p^*}$, we have
\begin{align*}
\MoveEqLeft
\big\langle (\Id - R_p R_{p^*}^*)(y) , \partial_{\psi,q,p^*}(x) \big\rangle 
=\big\langle y , \partial_{\psi,q,p^*}(x) \big\rangle  - \big\langle R_pR_{p^*}^*(y), \partial_{\psi,q,p^*}(x) \big\rangle \\
&\ov{\eqref{equ-1-proof-lem-Hodge-decomposition-sum-Fourier}}{=} \big\langle y , \partial_{\psi,q,p^*}(x) \big\rangle  - \big\langle R_{p^*}^*(y), A_{p^*}^{\frac12}(x) \big\rangle 
=\big\langle y, \partial_{\psi,q,p^*}(x) \big\rangle - \big\langle y, R_{p^*} A_{p^*}^{\frac12}(x) \big\rangle \\
&=\big\langle y, \partial_{\psi,q,p^*}(x) \big\rangle - \big\langle y, \partial_{\psi,q,p^*} A_{p^*}^{-\frac12} A_{p^*}^\frac12(x) \big\rangle 
= 0.
\end{align*}
By \eqref{lien-ker-image}, we conclude that $\big(\Id - R_pR_{p^*}^*\big)(y)$ belongs to $\ker (\partial_{\psi,q,p^*})^*$.
\end{proof}


From now on, we suppose the discrete group $G$ to be weakly amenable such that $\Gamma_q(H) \rtimes_\alpha G$ has $\QWEP$ (e.g. if $G$ is amenable). In the proof of Proposition \ref{prop-Hodge-decomposition-intersection-Fourier} below, we shall need some information on the Wiener-Ito chaos decomposition for $q$-Gaussians. This is collected in the following lemma.

\begin{lemma}
\label{lem-eigenspace-q-OU-Fourier}
Let $-1 \leq q \leq 1$ and $1 < p < \infty$. Let $G$ be a weakly amenable discrete group such that $\Gamma_q(H) \rtimes_\alpha G$ has $\QWEP$. Consider an approximating net $(M_{\varphi_j})$ of finitely supported Fourier multipliers.
\begin{enumerate}
\item There exists a completely bounded projection $\mathcal{P} \co \L^p(\Gamma_q(H)) \to \L^p(\Gamma_q(H))$ onto the closed space spanned by $\{s_q(h) : h \in H\}$. Moreover, the projections are compatible for different values of $p$. The mapping $\mathcal{P} \rtimes \Id_{\L^p(\VN(G))}$ extends to a bounded operator on $\L^p(\Gamma_q(H) \rtimes_\alpha G)$.

\item For any $j$ and any $y \in \L^p(\Gamma_q(H) \rtimes_\alpha G)$, the element $(\mathcal{P} \rtimes M_{\varphi_j})(y)$ can be written as $\sum_{s \in \supp \varphi_j} s_q(h_{s}) \rtimes \lambda_s$ for some $h_{s} \in H$.

\item Denoting temporarily by $\mathcal{P}_p$ and $M_{\varphi_j,p}$ the operator $\mathcal{P}$ and $M_{\varphi_j}$ on the $p$-level, the identity mapping on $\Gauss_{q,p}(\C) \rtimes \mathrm{span} \{ \lambda_s : s \in \supp \varphi_j\}$ extends to an isomorphism 
\begin{equation}
\label{equ-1-proof-prop-Hodge-decomposition-intersection-Fourier}
J_{p,2,j} \co \Ran(\mathcal{P}_p \rtimes M_{\varphi_j,p}) 
\subseteq \L^p(\Gamma_q(H) \rtimes_\alpha G) \to \Ran(\mathcal{P}_2 \rtimes M_{\varphi_j,2}) \subseteq \L^2(\Gamma_q(H) \rtimes_\alpha G).
\end{equation}
\end{enumerate}
\end{lemma}

\begin{proof}
1. We simply take for Wick words $w(h)$ with $h \in \mathcal{F}_q(H)$, $\mathcal{P}(w(h)) = \delta_{h \in H \subseteq \mathcal{F}_q(H)} w(h)$. Then 1. follows from Lemma \ref{lem-Qp-bounded}.

2. Once we know that $\mathcal{P} \rtimes M_{\varphi_j} = (\mathcal{P} \rtimes \Id) \circ (\Id \rtimes M_{\varphi_j})$ is bounded according to point 1. and Proposition \ref{prop-Fourier-mult-crossed-product}, this is easy and left to the reader.

3. We have 
$$
\norm{\sum_{s \in \supp \varphi_j} s_q(h_{s}) \rtimes \lambda_s}_{\L^p(\Gamma_q(H) \rtimes G)} 
\leq \sum_{s \in \supp \varphi_j} \norm{s_q(h_{s})}_{\L^p(\Gamma_q(H))} 
\lesssim_{q,p} \sum_{s \in \supp \varphi_j} \norm{h_{s}}_H
$$ 
(the same estimate on the $\L^2$-level). Note that for any $s_0 \in \supp \varphi_j$ fixed, we have a completely contractive projection of $\L^p(\Gamma_q(H) \rtimes_\alpha G)$ onto $\Span \{ x \rtimes \lambda_{s_0} : \:  x \in \L^p(\Gamma_q(H)) \}$.
So we have 
$$
\norm{\sum_{s \in \supp \varphi_j} s_q(h_{s}) \rtimes \lambda_s}_{\L^p(\Gamma_q(H) \rtimes_\alpha G)} \gtrsim \norm{s_q(h_{s_0})}_{\L^p(\Gamma_q(H))} \gtrsim_{q,p} \norm{h_{s_0}}_H. 
$$
It follows that 
$$
\norm{\sum_{s \in \supp \varphi_j} s_q(h_{s}) \rtimes \lambda_s}_{\L^p(\Gamma_q(H) \rtimes_\alpha G)} \cong \sum_{s \in \supp \varphi_j} \norm{h_{s}}_H \cong \norm{\sum_{s \in \supp \varphi_j} s_q(h_s) \rtimes \lambda_s}_{\L^2(\Gamma_q(H) \rtimes_\alpha G)}.
$$
Thus, \eqref{equ-1-proof-prop-Hodge-decomposition-intersection-Fourier} follows.
\end{proof}

\begin{prop}
\label{prop-Hodge-decomposition-intersection-Fourier}
Let $-1 \leq q \leq 1$ and $1 < p < \infty$. Let $G$ be a weakly amenable discrete group such that $\Gamma_q(H) \rtimes_\alpha G$ is $\QWEP$. Then the subspaces from Lemma \ref{lem-Hodge-decomposition-sum-Fourier} have trivial intersection, i.e. $\ovl{\Ran \partial_{\psi,q,p}} \cap \ker (\partial_{\psi,q,p^*})^* = \{0\}$.
\end{prop}

\begin{proof}
We begin with the case $p = 2$. According to Theorem \ref{First-spectral-triple-2}, the unbounded operator $\D_{\psi,q,2}$ is selfadjoint on $\L^2(\VN(G)) \oplus \L^2(\Gamma_q(H) \rtimes_\alpha G)$. We thus have the orthogonal sum $\ovl{\Ran \D_{\psi,q,2}} \oplus \ker \D_{\psi,q,2} =\L^2(\VN(G)) \oplus \L^2(\Gamma_q(H) \rtimes_\alpha G)$. Considering vectors in the second component, that is, in $\L^2(\Gamma_q(H) \rtimes_\alpha G)$, we deduce that $\ovl{\Ran \partial_{\psi,q,2}}$ and $\ker (\partial_{\psi,q,2})^*$ are orthogonal, hence have trivial intersection.

We turn to the case $1 < p < \infty$. Consider an approximating net $(M_{\varphi_j})$. According to Lemma \ref{lem-eigenspace-q-OU-Fourier} point 1. and Proposition \ref{prop-Fourier-mult-crossed-product}, we have for any $j$ a completely bounded mapping $\mathcal{P} \rtimes M_{\varphi_j} = (\mathcal{P} \rtimes \Id_{}) \circ (\Id_{} \rtimes M_{\varphi_j}) \co \L^p(\Gamma_q(H) \rtimes_\alpha G) \to \L^p(\Gamma_q(H) \rtimes_\alpha G)$. We claim that
\begin{align}
 \text{the subspace } \ovl{\Ran \partial_{\psi,q,p}}\text{ is invariant under  } \Id_{} \rtimes M_{\varphi_j} & \label{equ-2-proof-propHodge-decomposition-intersection-Fourier} \\
 \text{the subspace } \ker (\partial_{\psi,q,p^*})^* \text{ is invariant under  } \Id_{} \rtimes M_{\varphi_j} & \label{equ-3-proof-propHodge-decomposition-intersection-Fourier} \\
\text{the restriction of } \mathcal{P} \rtimes \Id_{} \text{ on } \ovl{\Ran \partial_{\psi,q,p}}\text{ is the identity mapping}. & \label{equ-5-proof-propHodge-decomposition-intersection-Fourier}
\end{align}
For \eqref{equ-2-proof-propHodge-decomposition-intersection-Fourier}, for any $s \in G$, note that 
$$
(\Id \rtimes M_{\varphi_j}) \partial_{\psi,q,p}(\lambda_s) 
\ov{\eqref{def-partial-psi}}{=} (\Id \rtimes M_{\varphi_j})\big(s_q(b_\psi(s)) \rtimes \lambda_s\big) 
= \varphi_j(s) s_q(b_\psi(s)) \rtimes \lambda_s.
$$
This element belongs to $\ovl{\Ran \partial_{\psi,q,p}}$. By linearity and since $\P_G$ is a core for $\partial_{\psi,q,p}$ according to Proposition \ref{Prop-derivation-closable-sgrp}, we deduce that $\Id \rtimes M_{\varphi_j}$ maps $\Ran \partial_{\psi,q,p}$ into $\ovl{\Ran \partial_{\psi,q,p}}$. Now \eqref{equ-2-proof-propHodge-decomposition-intersection-Fourier} follows from the continuity of $\Id \rtimes M_{\varphi_j}$. 

For \eqref{equ-3-proof-propHodge-decomposition-intersection-Fourier}, note that if $x \in \dom \partial_{\psi,q,p^*}$ and $f \in \ker (\partial_{\psi,q,p^*})^*$, then
\begin{align*}
\MoveEqLeft
\big\langle (\Id \rtimes M_{\varphi_j})(f), \partial_{\psi,q,p^*}(x) \big\rangle 
=\big\langle f, (\Id \rtimes M_{\ovl{\varphi_j}})\partial_{\psi,q,p^*}(x) \big\rangle \\
&=\big\langle f, \partial_{\psi,q,p^*} M_{\ovl{\varphi_j}}(x) \big\rangle 
=\big\langle (\partial_{\psi,q,p^*})^*(f),M_{\ovl{\varphi_j}}(x) \big\rangle 
= 0.
\end{align*}
By \cite[Problem 5.27, p.~168]{Kat1}, we conclude that $(\Id \rtimes M_{\varphi_j})(f)$ belongs to $\ker (\partial_{\psi,q,p^*})^*$ and \eqref{equ-3-proof-propHodge-decomposition-intersection-Fourier} follows. For \eqref{equ-5-proof-propHodge-decomposition-intersection-Fourier}, for any $s \in G$ we have 
$$
(\mathcal{P} \rtimes \Id) (\partial_{\psi,q,p}(\lambda_s)) 
\ov{\eqref{def-partial-psi}}{=} (\mathcal{P} \rtimes \Id)\big(s_q(b_\psi(s)) \rtimes \lambda_s\big) 
=s_q(b_\psi(s)) \rtimes \lambda_s 
\in \Ran \partial_{\psi,q,p}.
$$
Now use in a similar manner as before linearity, the fact that $\P_G$ is a core of $\partial_{\psi,q,p}$ and the continuity of $\mathcal{P} \rtimes \Id$. 

Now, let $z \in \ovl{\Ran \partial_{\psi,q,p}} \cap \ker (\partial_{\psi,q,p^*})^*$. Then according to \eqref{equ-2-proof-propHodge-decomposition-intersection-Fourier} -- \eqref{equ-5-proof-propHodge-decomposition-intersection-Fourier}, we infer that $(\mathcal{P} \rtimes M_{\varphi_j})(z)$ belongs\footnote{\thefootnote. Note that $(\mathcal{P} \rtimes \Id_{})(z)=z$.} again to $\ovl{\Ran \partial_{\psi,q,p}} \cap \ker (\partial_{\psi,q,p^*})^*$. We claim that $J_{p,2,j} (\mathcal{P} \rtimes M_{\varphi_j})(z)$ belongs to $\ovl{\Ran \partial_{\psi,q,2}} \cap \ker (\partial_{\psi,q,2})^*$, where the mapping $J_{p,2,j}$ was defined in Lemma \ref{lem-eigenspace-q-OU-Fourier}. First $(\mathcal{P} \rtimes M_{\varphi_j})(z)$ belongs to $\ovl{\Ran \partial_{\psi,q,p}}$, so that there exists a sequence $(x_n)$ in $\dom \partial_{\psi,q,p}$ such that 
\begin{equation}
\label{Divers-1253}
(\mathcal{P} \rtimes M_{\varphi_j})(z) 
= \lim_{n \to \infty} \partial_{\psi,q,p}(x_n).
\end{equation}
Then we have in $\L^p(\Gamma_q(H) \rtimes_\alpha G)$, with $\psi_j = 1_{\supp \varphi_j}$, which is of finite support and thus induces a completely bounded multiplier $M_{\psi_j}$\footnote{\thefootnote. If $x \in \dom \partial_{\alpha,q,p}$, it is really easy to check that $(\Id \rtimes M_{\psi_j}) \partial_{\alpha,q,p}(x) =\partial_{\alpha,q,p} M_{\psi_j}(x)$.} ,
\begin{align*}
\MoveEqLeft
(\mathcal{P} \rtimes M_{\varphi_j})(z) 
=(\Id \rtimes M_{\psi_j}) \cdot (\mathcal{P} \rtimes M_{\varphi_j})(z) 
\ov{\eqref{Divers-1253}}{=} \lim_{n \to \infty} (\Id \rtimes M_{\psi_j}) \partial_{\psi,q,p}(x_n) \\
&=\lim_{n \to \infty} \partial_{\psi,q,p} M_{\psi_j}(x_n) 
=\lim_{n \to \infty} \sum_{s \in \supp \varphi_j} x_{n,s} \partial_{\psi,q,p}(\lambda_s).
\end{align*}
Since $J_{p,2,j}$ is an isomorphism, this limit also holds in $\L^2(\Gamma_q(H) \rtimes_\alpha G)$ and the element $J_{p,2,j} (\mathcal{P} \rtimes M_{\varphi_j})(z) = \lim_{n \to \infty} \sum_{s \in \supp \varphi_j} x_{n,s} \partial_{\psi,q,2}(\lambda_s)$ belongs to $\ovl{\Ran \partial_{\psi,q,2}}$. Furthermore, for some family $(h_{s})$ of elements of $H$, we have
\begin{equation}
\label{Equa-divers-1001}
(\mathcal{P} \rtimes M_{\varphi_j})(z) 
=\sum_{s \in \supp \varphi_j} s_q(h_{s}) \rtimes \lambda_s.
\end{equation} 
Then using that $(\mathcal{P} \rtimes M_{\varphi_j})(z)$ belongs again to $\ker (\partial_{\psi,q,p^*})^*$ in the last equality, we obtain
\begin{align*}
\MoveEqLeft
(\partial_{\psi,q,2})^* J_{p,2,j} (\mathcal{P} \rtimes M_{\varphi_j})(z) 
\ov{\eqref{Equa-divers-1001}}{=} (\partial_{\psi,q,2})^* \bigg(\sum_{s \in \supp \varphi_j} s_q(h_s) \rtimes \lambda_s\bigg) \\
&=(\partial_{\psi,q,p^*})^*\bigg(\sum_{s \in \supp \varphi_j} s_q(h_s) \rtimes \lambda_s\bigg) 
=(\partial_{\psi,q,p^*})^*(\mathcal{P} \rtimes M_{\varphi_j})(z)
=0.             
\end{align*}
We have shown that $J_{p,2,j} (\mathcal{P} \rtimes M_{\varphi_j})(z)$ belongs to $\ovl{\Ran \partial_{\psi,q,2}} \cap \ker (\partial_{\psi,q,2})^*$. 

According to the beginning of the proof, the last intersection is trivial. It follows that $J_{p,2,j} (\mathcal{P} \rtimes M_{\varphi_j}) (z) = 0$. Since $J_{p,2,j}$ is an isomorphism, we infer that $(\mathcal{P} \rtimes M_{\varphi_j})(z) = 0$ for any $j$. Since $G$ is weakly amenable and $\Gamma_q(H) \rtimes_\alpha G$ has $\QWEP$, the net $(\Id \rtimes M_{\varphi_j})$ converges to $\Id_{\L^p(\Gamma_q(H) \rtimes_\alpha G)}$ for the point norm topology of $\L^p(\Gamma_q(H) \rtimes_\alpha G)$. We deduce that $(\mathcal{P} \rtimes \Id)(z) = 0$. But we had seen in \eqref{equ-5-proof-propHodge-decomposition-intersection-Fourier} that $(\mathcal{P} \rtimes \Id)(z)=z$, so that $z=0$ and we are done.
\end{proof}

Combining Lemma \ref{lem-Hodge-decomposition-sum-Fourier} and Proposition \ref{prop-Hodge-decomposition-intersection-Fourier}, we can now deduce the following corollary.

\begin{cor}
\label{cor-key-Hodge-decomposition-Fourier}
Let $-1 \leq q \leq 1$ and $1 < p < \infty$. Let $G$ be a weakly amenable discrete group such that $\Gamma_q(H) \rtimes_\alpha G$ is $\QWEP$ (e.g. $G$ amenable, or $G$ is a free group and $q = \pm 1$). Then we have a topological direct sum decomposition 
\begin{equation}
\label{Cor-decompo1-Fourier}
\L^p(\Gamma_q(H) \rtimes_\alpha G) 
= \ovl{\Ran \partial_{\psi,q,p}} \oplus \ker (\partial_{\psi,q,p^*})^*.
\end{equation}
where the associated first bounded projection is $R_pR_{p^*}^*$. In particular, we have
$$
\Ran (\partial_{\psi,q,p^*})^* 
=\Ran (\partial_{\psi,q,p^*})^*|{\ovl{\Ran \partial_{\psi,q,p}}}.
$$
\end{cor}

\begin{proof}
According to Lemma \ref{lem-Hodge-decomposition-sum-Fourier}, the above subspaces add up to $\L^p(\Gamma_q(H) \rtimes_\alpha G)$, and according to Proposition \ref{prop-Hodge-decomposition-intersection-Fourier}, the sum is direct. By \cite[Theorem 1.8.7]{KaR1}, we conclude that the decomposition is topological. In the course of the proof of Lemma \ref{lem-Hodge-decomposition-sum-Fourier}, we have seen that for any $y \in \L^p(\Gamma_q(H) \rtimes_\alpha G)$, we have the suitable decomposition $y = R_pR_{p^*}^*(y) + (\Id - R_pR_{p^*}^*)(y)$. So the associated first bounded projection is $R_pR_{p^*}^*$.
\end{proof}

\begin{thm}
\label{Thm-full-operator-bisectorial-Fourier}
Let $-1 \leq q \leq 1$ and $1 < p < \infty$. Let $G$ be a weakly amenable discrete group such that $\Gamma_q(H) \rtimes_\alpha G$ is $\QWEP$. Consider the operator $\D_{\psi,q,p}$ from \eqref{equ-full-Hodge-Dirac-operator-Fourier}. Then $\D_{\psi,q,p}$ is bisectorial and has a bounded $\HI(\Sigma_\omega^\pm)$ functional calculus.
\end{thm}

\begin{proof}
According to Corollary \ref{cor-key-Hodge-decomposition-Fourier}, the space $\L^p(\VN(G)) \oplus \L^p(\Gamma_q(H) \rtimes_\alpha G)$ admits the topological direct sum decomposition $\L^p(\VN(G)) \oplus \L^p(\Gamma_q(H) \rtimes_\alpha G) \ov{\eqref{Cor-decompo1-Fourier}}{=} \L^p(\VN(G)) \oplus \ovl{\Ran \partial_{\psi,q,p}} \oplus \ker (\partial_{\psi,q,p^*})^*$ into a sum of three subspaces. Along this decomposition, we can write\footnote{\thefootnote. Here the notation $(\partial_{\psi,q,p^*})^*$ is used for the restriction of $(\partial_{\psi,q,p^*})^*$ on the subspace $\ovl{\Ran \partial_{\psi,q,p}}$.}
$$
\D_{\psi,q,p}
\ov{\eqref{equ-full-Hodge-Dirac-operator-Fourier}}{=}\begin{bmatrix} 
0 & (\partial_{\psi,q,p^*})^* & 0 \\ 
\partial_{\psi,q,p} & 0 & 0 \\ 
0 & 0 & 0 
\end{bmatrix}
\ov{\eqref{Def-Dirac-operator-Fourier}}{=} \begin{bmatrix} 
D_{\psi,q,p} & 0 \\ 
0 & 0 
\end{bmatrix}.
$$
According to Theorem \ref{Th-functional-calculus-bisector-Fourier}, the operator $D_{\psi,q,p}$ is bisectorial and does have a bounded $\HI(\Sigma_\omega^\pm)$ functional calculus. So we conclude tthe same thing for $\D_{\psi,q,p}$. See also Proposition \ref{prop-resolvent-Hodge-Fourier} below.
\end{proof}

\begin{thm}[Hodge decomposition]
\label{Th-Hodge-decomposition-Fourier}
Suppose $1<p<\infty$ and $-1 \leq q \leq 1$. Let $G$ be a weakly amenable discrete group such that $\Gamma_q(H) \rtimes_\alpha G$ is $\QWEP$. If we identify $\ovl{\Ran \partial_{\psi,q,p}}$ and $\ovl{\Ran (\partial_{\psi,q,p^*})^*}$ as the closed  subspaces $\{0\} \oplus \ovl{\Ran \partial_{\psi,q,p}}$ and $\ovl{\Ran (\partial_{\psi,q,p^*})^*} \oplus \{0\}$ of $\L^p(\VN(G)) \oplus \L^p(\Gamma_q(H) \rtimes_\alpha G)$, we have 
\begin{equation}
\label{Hodge-decomposition}
\L^p(\VN(G)) \oplus \L^p(\Gamma_q(H) \rtimes_\alpha G) 
=\ovl{\Ran \partial_{\psi,q,p}} \oplus \ovl{\Ran (\partial_{\psi,q,p^*})^*} \oplus \ker \D_{\psi,q,p}. 
\end{equation}  
\end{thm}

\begin{proof}
From the definition \eqref{equ-full-Hodge-Dirac-operator-Fourier}, it is obvious that $\ker \D_{\psi,q,p}= (\ker \partial_{\psi,q,p} \oplus \{0\}) \oplus (\{0\} \oplus \ker (\partial_{\psi,q,p^*})^* )$. We deduce that
\begin{align*}
\MoveEqLeft
\L^p(\VN(G)) \oplus \L^p(\Gamma_q(H) \rtimes_\alpha G) \\            
&\ov{\eqref{Decompo-sympa-group} \eqref{Cor-decompo1-Fourier}}{=} \bigl(\ovl{\Ran (\partial_{\psi,q,p^*})^*} \oplus \ker \partial_{\psi,q,p}\bigr)
\oplus \bigl(\ker (\partial_{\psi,q,p^*})^* \oplus \ovl{\Ran \partial_{\psi,q,p}}\bigr) \\
&=(\{0\} \oplus \ovl{\Ran \partial_{\psi,q,p}}) \oplus (\ovl{\Ran (\partial_{\psi,q,p^*})^*} \oplus \{0\}) \oplus (\ker \partial_{\psi,q,p} \oplus \{  0 \}) \oplus (\{ 0 \} \oplus \ker (\partial_{\psi,q,p^*})^*) \\
&=(\{0\} \oplus \ovl{\Ran \partial_{\psi,q,p}}) \oplus (\ovl{\Ran (\partial_{\psi,q,p^*})^*} \oplus \{0\}) \oplus\ker \D_{\psi,q,p}.
\end{align*}  
\end{proof}

\begin{remark} \normalfont
\label{rem-dimension-free-Hodge-Fourier}
An inspection in all the steps of the proof of Theorem \ref{Thm-full-operator-bisectorial-Fourier} shows that the angle of the $\HI(\Sigma_\omega^\pm)$ calculus can be chosen $\omega > \frac{\pi}{2} | \frac1p - \frac12 |$ and that the norm of the calculus is bounded by a constant $K_\omega$ not depending on $G$ nor the cocycle $(b_\psi,H)$, in particular it is independent of the dimension of $H$.
\end{remark}

\begin{proof}
First note that since $A_p$ has a (sectorial) completely bounded $\HI(\Sigma_{2\omega})$ calculus with angle $2\omega > \pi | \frac1p - \frac12 |$ \cite[Theorem 4.1]{Arh7}, by the representation \eqref{D-alpha-carre-egal-Fourier} together with the fact that the spectral multipliers of $A_p$ are Fourier multipliers and with Proposition \ref{prop-Fourier-mult-crossed-product}, $\D_{\psi,q,p}^2$ also has a sectorial $\HI(\Sigma_{2\omega})$ calculus.
According to \cite[Proof of Theorem 10.6.7, Theorem 10.4.4 (1) and (3), Proof of Theorem 10.4.9]{HvNVW2}, the operator $\D_{\psi,q,p}$ has then an $\HI(\Sigma_\omega^\pm)$ bisectorial calculus to the angle $\omega > \frac{\pi}{2} | \frac1p - \frac12 |$ with a norm control
\begin{equation}
\label{equ-1-rem-dimension-free-Hodge-Fourier}
\norm{f(\D_{\psi,q,p})} \leq K_\omega \left( M_{2\omega,\D_{\psi,q,p}^2}^\infty \right)^2 \left( M_{\omega,\D_{\psi,q,p}}^R \right)^2 \norm{f}_{\infty,\omega},
\end{equation}
where $K_\omega$ is a constant only depending on $\omega$ (and not on $G$ nor the cocycle $(b_\psi,H)$).
Here $M_{2\omega,\D_{\psi,q,p}^2}^\infty$ is the $\HI(\Sigma_{2 \omega})$ calculus norm of $\D_{\psi,q,p}^2$ and
\begin{equation}
\label{equ-2-rem-dimension-free-Hodge-Fourier}
M_{\omega,\D_{\psi,q,p}}^R  = R\left(\left\{ \lambda (\lambda - \D_{\psi,q,p})^{-1} :\: \lambda \in \C \backslash \{ 0 \},\: \left|\,|\arg(\lambda)| - \frac{\pi}{2} \right| < \frac{\pi}{2} - \omega \right\}\right).
\end{equation}
Thus it remains to show that both $M_{2\omega,\D_{\psi,q,p}^2}^\infty$ and $M_{\omega,\D_{\psi,q,p}}^R$ can be chosen independently of the Hilbert space $H$ and the cocycle $b_\psi$.
Let us start with $M_{2\omega,\D_{\psi,q,p}^2}^\infty$.
It is controlled according to the above reasoning and \eqref{D-alpha-carre-egal-Fourier}, by $M_{2\omega,\Id_{S^p} \ot A_p}^\infty$, that is, the completely bounded $\HI$ calculus norm of $A_p$. Moreover, an application of \cite[Proposition 5.8]{JMX} shows that it suffices to consider only $2\omega > \frac{\pi}{2}$. According to \cite[Theorem 4.1]{Arh7}, we have a certain decomposition of the semigroup $(T_{t,p})_{t \geq 0}$ generated by $A_p$, given by
\[ 
\Id_{S^p} \ot T_{t,p} 
= (\Id_{S^p} \ot \E_p)(\Id_{S^p} \ot U_{t,p})(\Id_{S^p} \ot J_p) . 
\]
Here, $\Id_{S^p} \ot \E_p$ and $\Id_{S^p} \ot J_p$ are contractions and $\Id_{S^p} \ot U_{t,p}$ is a group of isometries. An inspection of the proof of Lemma \ref{lem-technical-HI-Laplace-transform} shows that the constant in the second condition there is a bound of the $\HI$ calculus in the first condition there, so that it suffices to show that the generator of $\Id_{S^p} \ot U_{t,p}$ has a bounded $\HI(\Sigma_{2\omega})$ calculus with a norm controlled by a constant depending only on $2 \omega > \frac{\pi}{2}$. By \cite[Proof of Theorem 10.7.10]{HvNVW2}, the norm of the calculus of $U_{t,p}$ is controlled by $c_{2\omega} \beta^2_{p,X} h_{p,X}$, $c_{2\omega}$ denoting a constant depending only on $2\omega$, $\beta_{p,X}$ denoting the UMD constant of $X = S^p(\L^p(M))$ and $h_{p,X}$ denoting the Hilbert transform norm, i.e. on the space $\L^p(\R,S^p(\L^p(M)))$.

These are controlled by a constant depending only on $p$ but not on $M$. Indeed, using \cite[Proposition 4.2.15]{HvNVW2}, it suffices to control the UMD constant of $S^p(\L^p(M))$. It is known that the constant is finite. By the following direct sum argument, one can see that it does not depend on $M$, needed for our dimension freeness. Suppose that there are noncommutative $\L^p$ spaces $\L^p(M_i)$ such that the Hilbert transform norm 
$$
\norm{H}_{\L^p(\R,S^p(\L^p(M_i))) \to \L^p(\R,S^p(\L^p(M_i)))} 
\xra[i \to +\infty]{} +\infty.
$$
We have a direct sum $\oplus_i^p S^p(\L^p(M_i)) = S^p(\L^p(\oplus_i M_i))$ isometrically.
But considering a sequence $(x_i)$ of elements of $\oplus_i^p S^p(\L^p(M_i))$ such that $x_i = 0$ for $i \neq i_0$ one sees that the Hilbert transform norm
$$ 
\norm{H}_{\L^p(\R,S^p(\L^p(\oplus_i M_i))) \to \L^p(\R,S^p(\L^p(\oplus_i M_i)))} 
\geq \norm{H}_{\L^p(\R,S^p(\L^p(M_{i_0}))) \to \L^p(\R,S^p(\L^p(M_{i_0})))} . 
$$
On the one hand, the left hand side is finite, on the other hand, the right hand side tends to infinity. By contradiction we conclude the claimed uniform control on the UMD constant of $S^p(\L^p(M))$. So we have the desired control of $M_{2\omega,\Id_{S^p} \ot A_p}^\infty$, and thus of $M_{2\omega,\D_{\psi,q,p}^2}^\infty$.

We turn to the control of $M_{\omega,\D_{\psi,q,p}}^R$ from \eqref{equ-2-rem-dimension-free-Hodge-Fourier}. According to \eqref{Resolvent-Fourier} extended to complex times $z$ belonging to some bisector $\Sigma_{\sigma}^\pm$ with $\sigma = \frac{\pi}{2} - \omega$,
it suffices to control the following $R$-bounds
\begin{align}
& R \left( \left\{ (\Id + z^2 A_p)^{-1} :\: z \in \Sigma_\sigma^\pm \right\} \right) \label{equ-3-rem-dimension-free-Hodge-Fourier}  \\
& R \left( \left\{ z (\Id + z^2 A_p)^{-1} (\partial_{\psi,q,p^*})^* : \: z \in \Sigma_\sigma^\pm \right\} \right) \label{equ-4-rem-dimension-free-Hodge-Fourier} \\
& R \left( \left\{ z \partial_{\psi,q,p} (\Id + z^2 A_p)^{-1} : \: z \in \Sigma_\sigma^\pm \right\} \right) \label{equ-5-rem-dimension-free-Hodge-Fourier} \\
& R \left( \left\{ \Id \rtimes (\Id + z^2 A_p)^{-1} : \: z \in \Sigma_\sigma^\pm \right\} \right) \label{equ-6-rem-dimension-free-Hodge-Fourier}.
\end{align}
Indeed, an operator matrix family is $R$-bounded if and only if all operator entries in the matrix are $R$-bounded.
According to \cite[Theorem 10.3.4 (1)]{HvNVW2}, \eqref{equ-3-rem-dimension-free-Hodge-Fourier} is $R$-bounded since we can write $(\Id + z^2 A_p)^{-1} = \Id - f(z^2A_p)$ with $f(\lambda) = \lambda(1+\lambda)^{-1}$.
Moreover, by the same reference, its $R$-bound is controlled by $M_{2\omega - \epsi,A_p}^\infty$, which in turn, by the above argument of dilation can be controlled independently of the Hilbert space $H$ and the cocycle $b_\psi$.
The same argument shows that also \eqref{equ-6-rem-dimension-free-Hodge-Fourier} is $R$-bounded.
Since the operator family in \eqref{equ-4-rem-dimension-free-Hodge-Fourier} consists of the family of the adjoints in \eqref{equ-5-rem-dimension-free-Hodge-Fourier} and $R$-boundedness is preserved under adjoints, it suffices to prove that \eqref{equ-5-rem-dimension-free-Hodge-Fourier} is $R$-bounded.
To this end, we decompose for $z \in \Sigma_\omega$ the positive part of the bisector (similarly if $z$ belongs to the negative part of the bisector)
\[ z \partial_{\psi,q,p} (\Id + z^2 A_p)^{-1} = \left[\partial_{\psi,q,p} A_p^{-\frac12}\right] \left(z^2 A_p \right)^{\frac12}(\Id + z^2 A_p)^{-1} = \left[\partial_{\psi,q,p} A_p^{-\frac12}\right] f(z^2A_p) \]
with $f(\lambda) = \sqrt{\lambda}(1+\lambda)^{-1}$.
Again \cite[Theorem 10.3.4 (1)]{HvNVW2} shows that the term $f(z^2A_p)$ is $R$-bounded with $R$-bound controlled by some constant independent of the cocycle.
Finally, we are left to show that the Riesz transform is bounded by a constant independent of the cocycle, that is,
\begin{equation}
\label{equ-7-rem-dimension-free-Hodge-Fourier}
\norm{\partial_{\psi,q,p} (x) }_{\L^p(\Gamma_q(H) \rtimes_\alpha G)} 
\leq C \bnorm{A_p^{\frac12}(x) }_{\L^p(\VN(G))}.
\end{equation}
For this in turn we refer to Proposition \ref{prop-Riesz-Fourier-dimension-free}.
\end{proof}

\subsection{Hodge-Dirac operator on $\L^p(\VN(G)) \oplus \Omega_{\psi,q,p}$}
\label{Sec-Omega-psi}

We keep the standing assumptions of the two preceding subsections and thus have a markovian semigroup $(T_t)_{t \geq 0}$ of Fourier multipliers with generator $A_p$, the noncommutative gradient $\partial_{\psi,q,p}$ and its adjoint $(\partial_{\psi,q,p^*})^*$. We also fix the parameters $1 < p < \infty$ and $-1 \leq q \leq 1$, and assume that the discrete group $G$ is weakly amenable such that $\Gamma_q(H) \rtimes_\alpha G$ has $\QWEP$, so that the main results from the preceding subsection are valid. For the rest of this section we consider the Hodge-Dirac operator
$$
\D_{\psi,q,p} 
\ov{\mathrm{def}}{=} 
\begin{bmatrix} 0 & (\partial_{\psi,q,p^*})^* & 0 \\ 
\partial_{\psi,q,p} & 0 & 0 \\ 
0 & 0 & 0 
\end{bmatrix}
$$
on the bigger space $\L^p(\VN(G)) \oplus \ovl{\Ran \partial_{\psi,q,p}} \oplus \ker (\partial_{\psi,q,p^*})^*$, with domain 
\[
\dom  \D_{\psi,q,p} 
\ov{\mathrm{def}}{=}  \dom \partial_{\psi,q,p} \oplus \left( \dom (\partial_{\psi,q,p^*})^* \cap \ovl{\Ran \partial_{\psi,q,p}} \right) \oplus \ker (\partial_{\psi,q,p^*})^*.
\]
In the following, we consider the bounded operator

\begin{equation}
\label{equ-resolvent-Hodge-Fourier}
T 
\ov{\mathrm{def}}{=} \begin{bmatrix} (\Id + t^2 A_p)^{-1} & \i t (\Id + t^2 A_p)^{-1} (\partial_{\psi,q,p^*})^* & 0 \\
\i t \partial_{\psi,q,p} (\Id + t^2 A_p)^{-1} & \Id \rtimes (\Id + t^2 A_p)^{-1} & 0 \\
0 & 0 & \Id_{\ker (\partial_{\psi,q,p^*})^*}
\end{bmatrix}
\end{equation}
on the space $\L^p(\VN(G)) \oplus \ovl{\Ran \partial_{\psi,q,p}} \oplus \ker (\partial_{\psi,q,p^*})^*$. Here, we interpret $\i t \partial_{\psi,q,p} (\Id + t^2 A_p)^{-1}$ and $\i t (\Id + t^2 A_p)^{-1} (\partial_{\psi,q,p^*})^*$ as the bounded extensions $\L^p(\VN(G)) \to \ovl{\Ran \partial_{\psi,q,p}} \subseteq \L^p(\Gamma_q(H) \rtimes_\alpha G)$ resp. $\ovl{\Ran \partial_{\psi,q,p}} \subseteq \L^p(\Gamma_q(H) \rtimes_\alpha G) \to \L^p(\VN(G))$ guaranteed by Proposition \ref{Prop-R-gradient-bounds-Fourier}. Then this proposition and Theorem \ref{Th-Hodge-decomposition-Fourier} yield that $T$ is a bounded operator on $\L^p(\VN(G)) \oplus \L^p(\Gamma_q(H) \rtimes_\alpha G)$.

\begin{prop}
\label{prop-resolvent-Hodge-Fourier}
Let $1 < p< \infty$, $-1 \leq q \leq 1$ and $G$ be a weakly amenable discrete group such that $\Gamma_q(H) \rtimes_\alpha G$ has $\QWEP$.
We have
\begin{equation}
\label{equ-1-prop-resolvent-Hodge-Fourier}
T (\Id - \i t \D_{\psi,q,p}) 
= \Id_{\dom \D_{\psi,q,p}}
\end{equation}
and
\begin{equation}
\label{equ-2-prop-resolvent-Hodge-Fourier}
(\Id - \i t \D_{\psi,q,p}) T 
= \Id_{\L^p(\VN(G)) \oplus \L^p(\Gamma_q(H) \rtimes_\alpha G)} .
\end{equation}
\end{prop}

\begin{proof}
If $R \ov{\mathrm{def}}{=} (\Id + t^2 A_p)^{-1}$, a straightforward calculation shows that
\begin{align*}
T (\Id - \i t \D_{\psi,q,p}) & = \begin{bmatrix} R(\Id + t^2 (\partial_{\psi,q,p^*})^* \partial_{\psi,q,p}) & 0 & 0 \\
\i t \partial_{\psi,q,p} R - \i t \Id \rtimes R \partial_{\psi,q,p} & t^2 \partial_{\psi,q,p} R (\partial_{\psi,q,p^*})^* + \Id \rtimes R & 0 \\
0 & 0 & \Id_{\ker (\partial_{\psi,q,p^*})^*} \end{bmatrix} \\
& \ov{\mathrm{def}}{=}
 \begin{bmatrix} 
(I) & 0 & 0 \\ 
(II) & (III) & 0 \\ 
0 & 0 & \Id_{\ker (\partial_{\psi,q,p^*})^*} 
\end{bmatrix}.
\end{align*}
We check the expressions $(I),(II),(III)$. For $(I)$, note that according to Proposition \ref{Prop-lien-generateur-partial-Fourier}, we have $(I) = (\Id + t^2 A_p)^{-1} + t^2(\Id + t^2 A_p)^{-1} (\partial_{\psi,q,p^*})^* \partial_{\psi,q,p} = \Id_{\dom \partial_{\psi,q,p}}$, since we recall that $(\Id + t^2 A_p)^{-1} (\partial_{\psi,q,p^*})^*$ is interpreted as a bounded operator $\ovl{\Ran \partial_{\psi,q,p^*}} \to \L^p(\VN(G))$. Then $(II) = 0$ on $\dom \partial_{\psi,q,p}$ according to Proposition \ref{Prop-commuting-Fourier}. Finally, we note that on the one hand, it is easy to check with Proposition \ref{Prop-commuting-Fourier} that $(III) = \Id$ on $\P_{\rtimes,G}$.
On the other hand, $\P_{\rtimes,G}$ is a core of $(\partial_{\psi,q,p^*})^*$ according to Proposition \ref{Prop-core-1-group} and $\partial_{\psi,q,p}(\Id + t^2 A_p)^{-1}$ is bounded, the two of which imply easily that $(III) = \Id_{\ovl{\Ran \partial_{\psi,q,p}} \cap \dom (\partial_{\psi,q,p^*})^*}$. Altogether, we have shown \eqref{equ-1-prop-resolvent-Hodge-Fourier}.
We turn to \eqref{equ-2-prop-resolvent-Hodge-Fourier}.
Again a straightforward calculation shows that
\begin{align*}
\MoveEqLeft
(\Id - \i t \D_{\psi,q,p})T \\
 & = \begin{bmatrix} (\Id + t^2 (\partial_{\psi,q,p^*})^* \partial_{\psi,q,p}) R  & \i t R (\partial_{\psi,q,p^*})^* - \i t (\partial_{\psi,q,p^*})^* \Id \rtimes R & 0 \\
0 & t^2 \partial_{\psi,q,p} R (\partial_{\psi,q,p^*})^* + \Id \rtimes R & 0 \\
0 & 0 & \Id_{\ker (\partial_{\psi,q,p^*})^*} \end{bmatrix} \\
& \ov{\mathrm{def}}{=} \begin{bmatrix} (I) & (II) & 0 \\ 0 & (III) & 0 \\ 0 & 0 & \Id_{\ker (\partial_{\psi,q,p^*})^*} \end{bmatrix}.
\end{align*}
As for \eqref{equ-1-prop-resolvent-Hodge-Fourier}, one shows that $(I) = \Id_{\L^p(\VN(G))}$.
Again with Proposition \ref{Prop-commuting-Fourier}, one shows that $(II) = 0$ on $\dom (\partial_{\psi,q,p^*})^*$.
But $(II)$ is closed, since the product $AB$ of a closed operator $A$ and a bounded operator $B$ is closed, so since $\dom (\partial_{\psi,q,p^*})^* \cap \ovl{\Ran \partial_{\psi,q,p}}$ is dense in $\ovl{\Ran \partial_{\psi,q,p}}$, $(II) = 0$ on $\ovl{\Ran \partial_{\psi,q,p}}$.
Finally, we have already shown above that $(III) = \Id$ on $\dom (\partial_{\psi,q,p^*})^*$, and as before, $(III)$ is closed.
We infer that $(III) = \Id_{\ovl{\Ran \partial_{\psi,q,p}}}$.
\end{proof}


We recall from Theorem \ref{Thm-full-operator-bisectorial-Fourier} that $\D_{\psi,q,p}$ has a bounded $\HI(\Sigma_\omega^\pm)$ functional calculus on a bisector on the space $\L^p(\VN(G)) \oplus \L^p(\Gamma_q(H) \rtimes_\alpha G)$.
Next we show that an appropriate restriction of $\D_{\psi,q,p}$ to the space $\L^p(\VN(G)) \oplus \Omega_{\psi,q,p}$, where
$$
\Omega_{\psi,q,p} 
\ov{\mathrm{def}}{=}  \ovl{\vect}^{\L^p} \left\{ s_q(\xi) \rtimes \lambda_s : \: \xi \in H_\psi, \: s \in G \right\}
$$
and $H_\psi \ov{\mathrm{def}}{=} \ovl{\vect} \{ b_\psi(s) :\: s \in G \}$, is still bisectorial and admits an $\HI(\Sigma_\omega^\pm)$ functional calculus on a bisector.

\begin{lemma}
\label{lem-projection-on-Omega-psi}
Assume $-1 \leq q \leq 1$ and $1 < p < \infty$ and $G$ be a weakly amenable discrete group such that $\Gamma_q(H) \rtimes_\alpha G$ has $\QWEP$. There is a bounded projection $W \co \L^p(\VN(G)) \oplus \L^p(\Gamma_q(H) \rtimes_\alpha G) \to \L^p(\VN(G)) \oplus \L^p(\Gamma_q(H) \rtimes_\alpha G)$ on $\Omega_{\psi,q,p}$ such that $\L^p(\VN(G)) \subseteq \ker W$.
\end{lemma}

\begin{proof}
We remind the reader that the bounded projection $\mathcal{P} \co \L^p(\Gamma_q(H)) \to \L^p(\Gamma_q(H))$ from the part 1 of Lemma \ref{lem-eigenspace-q-OU-Fourier} is given in the following way. If $(e_k)_{k \geq 1}$ denotes an orthonormal basis of $H$ and for a multi-index $\underline{i} = (i_1,i_2,\ldots,i_n)$ of length $|\underline{i}|=n \in \N$ we let $e_{\underline{i}} \ov{\mathrm{def}}{=} e_{i_1} \ot e_{i_2} \ot \cdots \ot e_{i_n} \in \Fc_q(H)$, then the Wick word $w(e_{\underline{i}}) \in \Gamma_q(H)$ is determined by $w(e_{\underline{i}})\Omega = e_{\underline{i}}$, where $\Omega$ denotes as usual the vacuum vector. Then a careful inspection of \cite{JuL1}, in particular Theorem 3.5 there, shows that we have $\mathcal{P}(\sum_{\underline{i}} \alpha_{\underline{i}} e_{\underline{i}}) = \sum_{|\underline{i}| = 1} \alpha_{\underline{i}} w(e_{\underline{i}})$ for any finite sum and $\alpha_{\underline{i}} \in \C$. Here, $w(e_{\underline{i}}) = s_q(e_{\underline{i}})$ in case $|\underline{i}| = 1$. According to Lemma \ref{lem-Qp-bounded}, $\mathcal{P} \rtimes \Id_{\L^p(\VN(G))} \co \L^p(\Gamma_q(H) \rtimes_\alpha G) \to \L^p(\Gamma_q(H) \rtimes_\alpha G)$ is a bounded mapping.

Now note that if $Q \co H \to H$ denotes the orthogonal projection onto the closed subspace $H_\psi$ of $H$, then by \cite[Theorem 2.11]{BKS} there exists a trace preserving conditional expectation $\E \co \L^p(\Gamma_q(H)) \to \L^p(\Gamma_q(H))$ such that $\E(w(e_{\underline{i}})) = w(\Fc_q(Q)e_{\underline{i}})$. Note that by \eqref{Cocycle-law}, for any $s \in G$, each $\pi_s$ induces an operator $\pi_s \co H_\psi \to H_\psi$, and thus by orthogonality of $\pi_{s^{-1}}$ also $\pi_s \co H_\psi^\perp \to H_\psi^\perp$. Thus each $\pi_s$ and $Q$ commute, whence $\alpha_s \ov{\eqref{equ-markovian-automorphism-group}}{=} \Gamma_q(\pi_s)$ and $\E = \Gamma_q(Q)$ commute. We can use Lemma \ref{lem-van-Daele-Arhancet} and deduce that $\E$ extends to a normal complete contraction $\E \rtimes \Id_{\VN(G)} \co \Gamma_q(H) \rtimes_\alpha G \to \Gamma_q(H) \rtimes_\alpha G$. Again by Lemma \ref{lem-van-Daele-Arhancet}, we infer that $\E \rtimes \Id_{\VN(G)}$ is a conditional expectation which is trace preserving, so extends to a contraction on $\L^p(\Gamma_q(H) \rtimes_\alpha G)$ for $1 \leq p \leq \infty$.

We claim that $\mathcal{P}\E \rtimes \Id_{\L^p(\VN(G))}$ is a projection. To this end, it suffices\footnote{\thefootnote. Recall that the product of two commuting projections on a Banach space is a projection.} to check that $\mathcal{P}$ and $\E$ commute. We assume that the orthonormal basis $(e_k)_{k \in \N}$ is chosen in such a way that $(e_k)_{k \in \N_\psi}$ is an orthonormal basis of $H_\psi$ for some $\N_\psi \subseteq \N$. From the foregoing, for a Wick word $w(e_{\underline{i}})$ we have
\[
\mathcal{P}\E w(e_{\underline{i}}) 
= \mathcal{P} w(\Fc_q(Q) e_{\underline{i}}) 
= \delta_{|\underline{i}|=1} w(\Fc_q(Q) e_{\underline{i}}) 
= \delta_{|\underline{i}|=1} \delta_{\underline{i} \in \N_\psi} w(e_{\underline{i}})
\]
and
\[
\E \mathcal{P} w(e_{\underline{i}}) 
= \E \delta_{|\underline{i}| = 1} w(e_{\underline{i}}) 
= \delta_{|\underline{i}| = 1} w(\Fc_q(Q)e_{\underline{i}}) 
= \delta_{|\underline{i}|=1} \delta_{\underline{i} \in \N_\psi} w(e_{\underline{i}}).
\]
Thus $\mathcal{P}$ and $\E$ commute on a total set, so commute on all of $\L^p(\Gamma_q(H))$. It suffices to consider the projection $W' \ov{\mathrm{def}}{=} \mathcal{P}\E \rtimes \Id_{\L^p(\VN(G))}$ and to finally extend it to $W  \co \L^p(\VN(G)) \oplus \L^p(\Gamma_q(H) \rtimes_\alpha G) \to \Omega_\psi$ by setting $W(x,y) \ov{\mathrm{def}}{=} 0 \oplus W'(y)$ and observe by a standard density and continuity argument that $\Ran W = \Omega_{\psi,q,p}$.
\end{proof}

The proof of the following elementary lemma is left to the reader.

\begin{lemma}
\label{lem-elementary-projections}
Let $X$ be a Banach space and let $Q_1,Q_2,Q_3,Q_4 \co X \to X$ be bounded projections. Assume that $Q_1 + Q_2 + Q_3 + Q_4 = \Id_X$ and that $Q_i Q_j = 0$ for $i < j$. Then we have a direct sum decomposition
$$
X 
= \Ran(Q_1) \oplus \Ran(Q_2) \oplus \Ran(Q_3) \oplus \Ran(Q_4).
$$
\end{lemma}

\begin{lemma}
\label{lem-Hodge-Fourier-stabilises-Omega-psi}
Let $1 < p <\infty$ and $-1 \leq q \leq 1$. Assume that the discrete group $G$ is weakly amenable and that $\Gamma_q(H) \rtimes_\alpha G$ has $\QWEP$. The subspace $\L^p(\VN(G)) \oplus \Omega_{\psi,q,p}$ is invariant under the resolvents $T$ from \eqref{equ-resolvent-Hodge-Fourier} of the Hodge-Dirac operator $\D_{\psi,q,p}$.
\end{lemma}

\begin{proof}
With respect to the decomposition projections from \eqref{Cor-decompo1-Fourier}, we can decompose the identity $\Id_{\L^p(\VN(G)) \oplus \L^p(\Gamma_q(H) \rtimes_\alpha G)}$ as a sum
\begin{align*}
\Id_{\L^p(\VN(G)) \oplus \L^p(\Gamma_q(H) \rtimes_\alpha G)} & = P_1 & \oplus P_2 & \oplus P_3 \\
\L^p(\VN(G)) \oplus \L^p(\Gamma_q(H) \rtimes_\alpha G) & = \L^p(\VN(G)) & \oplus \ovl{\Ran \partial_{\psi,q,p}}& \oplus \ovl{\Ran (\partial_{\psi,q,p^*})^*}.
\end{align*}
Then we claim that
\begin{align*}
\Id_{\L^p(\VN(G)) \oplus \L^p(\Gamma_q(H) \rtimes_\alpha G)} & = P_1 & \oplus P_2W \oplus P_3W & \oplus V \\
\L^p(\VN(G)) \oplus \L^p(\Gamma_q(H) \rtimes_\alpha G) & = \L^p(\VN(G)) & \oplus \underset{= \Omega_{\psi,q,p}}{\underbrace{\ovl{\Ran \partial_{\psi,q,p}}\oplus X_\psi }} & \oplus Y_\psi
\end{align*} 
for some subspaces $X_\psi$ and $Y_\psi$, coming with projections $P_1, P_2W$, $P_3W$ and $V \ov{\mathrm{def}}{=} (P_2+P_3)(\Id - W)$. To this end, we apply the auxiliary Lemma \ref{lem-elementary-projections}. Note first that $P_1 + P_2W + P_3W + V = P_1 + (P_2 + P_3)W + (P_2+P_3)(\Id - W) = P_1 + (P_2 + P_3) = \Id$. Then $P_2W$, $P_3W$ and $V$ are projections. Indeed, $P_2WP_2W = P_2 P_2W$, since $\Ran P_2W \subseteq \ovl{\Ran\partial_{\psi,q,p}} \subseteq \Omega_{\psi,q,p} = \Ran W$, and $P_2P_2W = P_2W$, since $P_2$ is a projection. Moreover, $P_3WP_3W = P_3W(\Id - P_1 - P_2)W = P_3W - P_3WP_1W - P_3WP_2W = P_3W - P_3\cdot 0 \cdot W - P_3P_2 W = P_3W - 0 - 0\cdot W = P_3W$. Thus, $P_2W$ and $P_3W$ are projections. Moreover, $V = (P_2+P_3)(\Id - W) = (P_1+P_2+P_3)(\Id-W) - P_1 = \Id - P_1 - W$, and $V^2 = (\Id-P_1-W)^2 = \Id + P_1 + W - 2 P_1 - 2W +P_1W + WP_1 = \Id - P_1 - W + 0 + 0 = V$. Thus, also $V$ is a projection. Now we check that some products of the four projections vanish as needed to apply Lemma \ref{lem-elementary-projections}. We choose the order $(Q_1,Q_2,Q_3,Q_4) = (P_1,P_3W,P_2W,V)$. First note that this is clear if one of the factors is $P_1$. Then $P_3WP_2W = 0$ since $\Ran(WP_2W) = \Ran(P_2W) \subseteq \ovl{\Ran \partial_{\psi,q,p}} \subseteq \ker P_3$. Moreover, $P_3WV = P_3W(P_2+P_3)(\Id-W) = P_3W(\Id-W) = 0$ and also $P_2WV = P_2W(P_2+P_3)(\Id-W) = 0$. We have shown the claim and thus have a direct sum decomposition of the space into four closed subspaces.

Now write the resolvent
\begin{equation}
\label{equ-1-proof-stabilises-Omega-psi}
T = 
\begin{bmatrix} 
A & B & 0 \\ 
C & D & 0 \\ 
0 & 0 & \Id_{\ker (\partial_{\psi,q,p^*})^*} \end{bmatrix}
\end{equation}
along the Hodge decomposition $\Id = P_1 + P_2 + P_3$. If $x \in \L^p(\VN(G)) \oplus \L^p(\Gamma_q(H) \rtimes_\alpha G)$, then $x$ belongs to $\L^p(\VN(G)) \oplus \Omega_{\psi,q,p}$ if and only if $V(x) = 0$. For such an $x$, we have
\begin{align*}
T(x) 
& = T(P_1x+P_2Wx+P_3Wx) \\
& = P_1AP_1x + P_2WCP_1x + P_3WCP_1x + VCP_1x + P_1BP_2Wx \\
& + P_2WDP_2Wx + P_3WDP_2Wx + VDP_2Wx + P_3Wx.
\end{align*}
The summands starting with $P_1$ and $P_2$ lie in $\L^p(\VN(G))$ and $\ovl{\Ran \partial_{\psi,q,p}} \subseteq \Omega_{\psi,q,p}$. The remaining summands are $P_3WCP_1x = P_3CP_1x = 0$ since $CP_1x \in \ovl{\Ran \partial_{\psi,q,p}}$; $VCP_1x = 0$ since $CP_1x \in \ovl{\Ran \partial_{\psi,q,p}} \subseteq \Omega_{\psi,q,p}$; $P_3WDP_2Wx = P_3DP_2Wx = 0$; $VDP_2Wx = 0$ since $DP_2Wx \in \ovl{\Ran \partial_{\psi,q,p}} \subseteq \Omega_{\psi,q,p}$.
Finally, $P_3Wx \in X_\psi \subseteq \Omega_{\psi,q,p}$. We conclude that $T(x) \in \L^p(\VN(G)) \oplus \Omega_{\psi,q,p}$.
\end{proof}

\begin{thm}
\label{prop-Hodge-Fourier-HI-on-Omega-psi}
Let $1 < p< \infty$ and $-1 \leq q \leq 1$.
Assume that the discrete group $G$ is weakly amenable and that $\Gamma_q(H) \rtimes_\alpha G$ has $\QWEP$. Consider the part $\D_{\psi,q,p}'$ of the Hodge-Dirac operator $\D_{\psi,q,p} \co \dom \D_{\psi,q,p} \subseteq \L^p(\VN(G)) \oplus \L^p(\Gamma_q(H) \rtimes_\alpha G ) \to \L^p(\VN(G)) \oplus \L^p(\Gamma_q(H) \rtimes_\alpha G)$ on the closed subspace $\L^p(\VN(G)) \oplus \Omega_{\psi,q,p}$. Then $\D_{\psi,q,p}'$ is bisectorial and has a bounded $\HI(\Sigma_\omega^\pm)$ functional calculus on a bisector.
\end{thm}

\begin{proof}
We want to apply \cite[Proposition 3.2.15]{Ege2}. To this end, note that we have proved above that $\D_{\psi,q,p}$ is bisectorial on $X = \L^p(\VN(G)) \oplus \L^p(\Gamma_q(H) \rtimes_\alpha G)$. Moreover, for any $t \in \R$, $(\Id - \i t \D_{\psi,q,p})^{-1}$ leaves invariant $Y = \L^p(\VN(G)) \oplus \Omega_{\psi,q,p}$ according to the above Lemma \ref{lem-Hodge-Fourier-stabilises-Omega-psi}. Note that then the same holds for $t$ belonging to any bisector $\Sigma^\pm$ to which $\D_{\psi,q,p}$ is bisectorial. Indeed, $z \mapsto (\Id - \i z \D_{\psi,q,p})^{-1}$ is an analytic function on the subset of $\C$ where it is defined. Then also $P_i (\Id - \i z \D_{\psi,q,p})^{-1} P_j$ is analytic for $i = 1,2,3$. Therefore by the uniqueness theorem of analytic functions, $(\Id - \i z \D_{\psi,q,p})^{-1}$ has the same form as \eqref{equ-1-proof-stabilises-Omega-psi} at least for $z$ belonging to such a bisector $\Sigma^\pm$.
But then the proof of Lemma \ref{lem-Hodge-Fourier-stabilises-Omega-psi} goes through for such $z$ in place of $t$. Now the theorem follows from an application of \cite[Proposition 3.2.15]{Ege2} together with Proposition \ref{prop-resolvent-Hodge-Fourier}.
\end{proof}

\subsection{Bimodule $\Omega_{\psi,q,p,c}$}
\label{Module-Omega}

In this short subsection, we continue to consider a markovian semigroup $(T_t)_{t \geq 0}$ of Fourier multipliers as in  Proposition \ref{prop-Schoenberg}. This time, we do not need approximation properties on the discrete group $G$. We shall clarify and generalize some result of \cite[pp.~585-586]{JMP2}. We need the following notion of bimodule which is different from the notion of \cite[Definition 5.4]{JuS1} and is inspired by the well-known theory of Hilbert bimodules. 
Recall that the notion of right $\L^p$-$M$-module is defined in Subsection \ref{Sec-Hilbertian-valued}.

\begin{defi}
\label{Def-bimodule}
Let $M$ and $N$ be von Neumann algebras. Suppose $1 \leq p < \infty$. An $\L^p$-$N$-$M$-bimodule is a right $\L^p$-$M$-module $X$ equipped with a structure of left-$N$-module such that the associated $\L^{\frac{p}{2}}(M)$-valued inner product $\langle \cdot, \cdot \rangle_X$ satisfies 
\begin{equation}
\label{Compatibility-bimodule}
\langle a^* x,y\rangle_X
=\langle x, ay \rangle_X, \quad x,y \in X, a \in N. 
\end{equation}
An $\L^p$-$M$-bimodule is an $\L^p$-$M$-$M$-bimodule.
\end{defi}

Suppose $-1 \leq q \leq 1$ and $2 \leq p < \infty$ (the case $p<2$ is entirely left to the reader). If $\E \co \L^p(\Gamma_q(H) \rtimes_\alpha G) \to \L^p(\VN(G))$ is the canonical conditional expectation, it is obvious that the formula 
\begin{equation}
\label{Bracket-bimodule}
\langle \omega, \eta \rangle_{}
\ov{\mathrm{def}}{=} \E(\omega^* \eta)
\end{equation}
defines an $\L^{\frac{p}{2}}(\VN(G))$-valued inner product on $\L^p(\Gamma_q(H) \rtimes_\alpha G)$. We can consider the associated right $\L^p$-$\VN(G)$-module $\L^p_c(\E)$. It is easy to see that $\L^p_c(\E)$ is an $\L^p$-$\VN(G)$-bimodule. We consider the closed subspace
\begin{equation}
\label{Def-Z-q-p}
\Omega_{\psi,q,p,c}' 
\ov{\mathrm{def}}{=} \ovl{\mathrm{span}} \big\{ \partial_{\psi,q,p}(x)a:  x \in \dom \partial_{\psi,q,p}, a \in \VN(G) \big\}
\end{equation}
of $\L^p(\Gamma_q(H) \rtimes_\alpha G)$. For any $a,b \in \VN(G)$ and any $x \in \dom \partial_{\psi,q,p}$, note that 
$$
(\partial_{\psi,q,p}(x)a)b
=\partial_{\psi,q,p}(x)ab.
$$ 
Thus by linearity and density, $\Omega_{\psi,q,p,c}'$ is a right $\VN(G)$-module.
Moreover, for any $a \in \VN(G)$ and any $x,b \in \P_G$, we have
\begin{align*}
\MoveEqLeft
b\partial_{\psi,q,p}(x)a
\ov{\eqref{Leibniz-Schur-gradient-mieux-sgrp}}{=} \big[\partial_{\psi,q,p}(b x)-\partial_{\psi,q,p}(b)x\big]a 
=\partial_{\psi,q,p}(bx)a-\partial_{\psi,q,p}(b)xa.            
\end{align*}
Thus $b\partial_{\psi,q,p}(x)a$ belongs to $\Omega_{\psi,q,p,c}'$.
Since $\P_G$ is a core for $\partial_{\psi,q,p}$ according to Proposition \ref{Prop-derivation-closable-sgrp}, the same holds for $x \in \dom \partial_{\psi,q,p}$.
If $b \in \VN(G)$ is a general element, we approximate it in the strong operator topology by a bounded net in $\P_G$ and obtain again that the same holds for $b \in \VN(G)$.
By linearity and density, we deduce that $\Omega_{\psi,q,p,c}'$ is a left $\VN(G)$-module, so finally a $\VN(G)$-bimodule. It is obvious that the restriction of the bracket \eqref{Bracket-bimodule} defines an $\L^{\frac{p}{2}}(\VN(G))$-valued inner product on this subspace. We can consider the associated right $\L^p$-$\VN(G)$-module $\Omega_{\psi,q,p,c}$ which is also an $\L^p$-$\VN(G)$-bimodule and which identifies canonically to a closed subspace of $\L^p_{c}(\E)\ov{\eqref{Egalite-fantastique-group}}{=} \L^p(\VN(G),\L^2(\Gamma_q(H))_{c,p})$. Finally, we recall that $H_\psi$ is the real Hilbert space generated by the $b_\psi(s)$'s where $s \in G$.

\begin{lemma}
Let $G$ be a discrete group.
\begin{enumerate}
	\item If $\xi \in H_\psi$ and if $s \in G$ then $s_q(\xi) \rtimes \lambda_s$ belongs to $\Omega_{\psi,q,p,c}$. 
	\item Moreover, we have
\begin{equation}
\label{Description-de-Omega-psi}
\Omega_{\psi,q,p,c}
=\ovl{\mathrm{span}}^{\Omega_{\psi,q,p,c}} \big\{ s_q(\xi) \rtimes \lambda_s : \xi \in H_\psi, s \in G \big\}.
\end{equation}
\end{enumerate}
\end{lemma}

\begin{proof}
1. For any $s \in G$, we have
$$
s_q(b_\psi(s)) \rtimes 1
\ov{\eqref{Product-crossed-product}}{=} \big(s_q(b_\psi(s)) \rtimes \lambda_s\big)(1 \rtimes \lambda_{s^{-1}})
\ov{\eqref{def-partial-psi}}{=} \partial_{\psi,q,p}(\lambda_s)(1 \rtimes \lambda_{s^{-1}}).
$$
Hence $s_q(b_\psi(s)) \rtimes 1$ belongs to \eqref{Def-Z-q-p}. If $\xi$ belongs to the span of the $b_\psi(s)$'s where $s \in G$, we deduce by linearity that $s_q(\xi) \rtimes 1$ belongs to $\Omega_{\psi,q,p,c}$. Now, for $\xi \in H_\psi$, there exists a sequence $(\xi_n)$ of elements of the previous span such that $\xi_n \to \xi$ in $H$. By \eqref{petit-Wick}, we infer that $s_q(\xi_n) \to s_q(\xi)$ in $\L^2(\Gamma_q(H))$. Hence $s_q(\xi_n) \rtimes 1 \to s_q(\xi) \rtimes 1$ in $\L^p(\VN(G),\L^2(\Gamma_q(H))_{c,p})$. We conclude that $s_q(\xi) \rtimes 1$ belongs to $\Omega_{\psi,q,p,c}$. Since $\Omega_{\psi,q,p,c}$ is a right $\VN(G)$-module, we conclude that $s_q(\xi) \rtimes \lambda_s=(s_q(\xi) \rtimes 1)\lambda_s$ belongs to $\Omega_{\psi,q,p,c}$.

2. For any $s,t \in G$, we have 
$$
\partial_{\psi,q,p}(\lambda_s)(1 \rtimes \lambda_t)
\ov{\eqref{def-partial-psi}}{=} (s_q(b_\psi(s)) \rtimes \lambda_s)(1 \rtimes \lambda_t)
\ov{\eqref{Product-crossed-product}}{=}s_q(b_\psi(s)) \rtimes \lambda_{st}
$$ 
which belongs to $\ovl{\mathrm{span}}^{\Omega_{\psi,q,p,c}} \big\{ s_q(\xi_s) \rtimes \lambda_s : \xi_s \in H_\psi, s \in G \big\}$. Since the closed span of the $\partial_{\psi,q,p}(\lambda_s)(1 \rtimes \lambda_t)$'s is $\Omega_{\psi,q,p,c}$\footnote{\thefootnote. It is clear that elements of the form $\partial_{\psi,q,p}(\lambda_s)(1 \rtimes \lambda_t)$ belong to $\Omega_{\psi,q,p,c}' \subseteq \Omega_{\psi,q,p,c}$.
On the other hand, for the density, since $\Omega_{\psi,q,p,c}'$ is by definition dense in $\Omega_{\psi,q,p,c}$, it suffices to approximate $\partial_{\psi,q,p}(x)(1 \rtimes a)$ by elements of the span of the $\partial_{\psi,q,p}(\lambda_s) (1 \rtimes \lambda_t)$ for $x \in \dom \partial_{\psi,q,p}$ and $a \in \VN(G)$.
The approximation needs to be in $\L^p_c(\E)$ norm, but according to Lemma \ref{lem-contractive-crossed-product-Hilbert-valued}, it suffices to approximate in $\L^p(\Gamma_q(H) \rtimes_\alpha G)$ norm.
By definition of $\VN(G)$ and Kaplansky's density theorem, there exists a bounded net $(a_\alpha)$ in $\P_G$ converging in the strong operator topology to $a$.
Then it is not difficult to see that the net $(1 \rtimes a_\alpha)$ converges in the strong operator topology to $1 \rtimes a$.
Thus we can approximate $\partial_{\psi,q,p}(x)(1 \rtimes a)$ in $\L^p$ norm by elements in the span of the $\partial_{\psi,q,p}(x) (1 \rtimes \lambda_t)$.
Since $\P_G$ is a core of $\partial_{\psi,q,p}$, we can then approximate in turn by elements in the span of the $\partial_{\psi,q,p}(\lambda_s)(1 \rtimes\lambda_t)$.}
 , the proof is complete.
\end{proof}

\subsection{Hodge-Dirac operators associated with semigroups of Markov Schur multipliers}
\label{Hodge-Schur}

In this subsection, we consider some markovian semigroup $(T_t)_{t \geq 0}$ of Schur multipliers acting on $\B(\ell^2_I)$ that we defined in Proposition \ref{def-Schur-markovian}.  If $1 \leq p <\infty$, we denote by $A_p$ the (negative) infinitesimal generator on $S^p_I$ which is defined as the closure of the unbounded operator $A \co \M_{I,\fin} \to S^p_I$, $e_{ij} \mapsto \norm{\alpha_i-\alpha_j}_H^2e_{ij}$. So $\M_{I,\fin}$ is a core of $A_p$. By \cite[(5.2)]{JMX}, we have $(A_{p})^*=A_{p^*}$ if $1<p<\infty$. If $-1 \leq q \leq 1$, recall that by Proposition \ref{Prop-derivation-closable}, we have a closed operator 
$$
\partial_{\alpha,q,p} \co \dom \partial_{\alpha,q,p} \subseteq S^p_I \to \L^p(\Gamma_q(H) \otvn \B(\ell^2_I)),\quad 
e_{ij} \mapsto s_q(\alpha_i-\alpha_j) \ot e_{ij}.	
$$ 
Note that the adjoint operator $(\partial_{\alpha,q,p^*})^* \co \dom (\partial_{\alpha,q,p^*})^* \subseteq \L^p(\Gamma_q(H) \otvn \B(\ell^2_I)) \to S^p_I$ is closed by \cite[p.~168]{Kat1}.  
We will now define a Hodge-Dirac operator $D_{\alpha,q,p}$ in \eqref{Def-Dirac-operator-Schur}, relying on $\partial_{\alpha,q,p}$ and its adjoint. Then the main topic of this subsection will be to show that $D_{\alpha,q,p}$ is $R$-bisectorial (Theorem \ref{Th-D-R-bisectorial-Schur}) and has a bounded $\HI$ functional calculus on $S^p_I \oplus \ovl{\Ran \partial_{\alpha,q,p}}$ (Theorem \ref{Th-functional-calculus-bisector-Schur}). By Remark \ref{remark-sgn-Schur}, this extends the Kato square root equivalence from Proposition \ref{Prop-equiv-q-gaussians}.


\begin{prop}
\label{Prop-lien-generateur-partial}
Suppose $1<p<\infty$ and $-1 \leq q \leq 1$. As unbounded operators, we have
\begin{equation}
\label{eq-lien-generateur-partial}
A_p
=(\partial_{\alpha,q,p^*})^*\partial_{\alpha,q,p}. 
\end{equation}
\end{prop}

\begin{proof}
By Lemma \ref{Lemma-adjoint-prtial-Schur} and \cite[p.~167]{Kat1}, $\partial_{\alpha,q,p}(\M_{I,\fin})$ is a subspace of $\dom (\partial_{\alpha,q,p^*})^*$. For any $i,j \in I$ we have
\begin{align}
\MoveEqLeft
\label{Equa-10000}
(\partial_{\alpha,q,p^*})^*\partial_{\alpha,q,p}(e_{ij})           
\ov{\eqref{def-delta-alpha}}{=}(\partial_{\alpha,q,p^*})^*\big(s_q(\alpha_i-\alpha_j) \ot e_{ij}\big) 
\ov{\eqref{Adjoint-partial-Schur}}{=} \tau\big( s_q(\alpha_i-\alpha_j)s_q(\alpha_i-\alpha_j)\big) e_{ij} \\
&\ov{\eqref{petit-Wick}}{=}\norm{\alpha_i-\alpha_j}_H^2\, e_{ij}
=A_p(e_{ij}). \nonumber
\end{align}
Hence for any $x,y \in \M_{I,\fin}$, by linearity we have 
\begin{align*}
\MoveEqLeft
\big\langle A_p^{\frac{1}{2}}(x), A_{p^*}^{\frac{1}{2}}(y) \big\rangle_{S^p_I,S^{p^*}_I}
=\big\langle A_p(x), y \big\rangle_{S^p_I,S^{p^*}_I}
\ov{\eqref{Equa-10000}}{=} \big\langle (\partial_{\alpha,q,p^*})^*\partial_{\alpha,q,p}(x), y \big\rangle_{S^p_I,S^{p^*}_I} \\
&\ov{\eqref{crochet-duality}}{=}\big\langle \partial_{\alpha,q,p}(x),\partial_{\alpha,q,p^*}(y) \big\rangle_{\L^p,\L^{p^*}}.            
\end{align*} 
Using the part 4 of Proposition \ref{Prop-derivation-closable}, it is not difficult to see that this identity extends to elements $x \in \dom A_p$. For any $x \in \dom A_p$ and any $y \in \M_{I,\fin}$, we obtain
$$
\big\langle A_p(x), y \big\rangle_{S^p_I,S^{p^*}_I}
=\big\langle \partial_{\alpha,q,p}(x),\partial_{\alpha,q,p^*}(y) \big\rangle_{\L^p,\L^{p^*}}.
$$
Recall that $\M_{I,\fin}$ is a core of $\partial_{\alpha,q,p^*}$ by the part 1 of Proposition \ref{Prop-derivation-closable}. So using \eqref{Def-core}, it is easy to check that this identity remains true for elements $y$ of $\dom \partial_{\alpha,q,p^*}$. By \eqref{Def-domaine-adjoint}, this implies that $\partial_{\alpha,q,p}(x) \in \dom (\partial_{\alpha,q,p^*})^*$ and that $(\partial_{\alpha,q,p^*})^*\partial_{\alpha,q,p}(x) = A_p(x)$. We conclude that $A_p \subseteq (\partial_{\alpha,q,p^*})^*\partial_{\alpha,q,p}$.

To prove the other inclusion we consider some $x \in \dom \partial_{\alpha,q,p}$ such that $\partial_{\alpha,q,p}(x)$ belongs to $\dom (\partial_{\alpha,q,p^*})^*$. By \cite[Theorem 5.29]{Kat1}, we have $(\partial_{\alpha,q,p^*})^{**}=\partial_{\alpha,q,p^*}$. We infer that $(\partial_{\alpha,q,p})^*\partial_{\alpha,q,p^*} \ov{\eqref{Adjoint-product-unbounded}}{\subseteq} \big((\partial_{\alpha,q,p^*})^*\partial_{\alpha,q,p}\big)^* \ov{\eqref{Inclusion-adjoint-unbounded}}{\subseteq} A_p^*$. For any $y \in \M_{I,\fin}$, using $\partial_{\alpha,q,p}(x) \in \dom (\partial_{\alpha,q,p^*})^*$ in the last equality, we deduce that
\begin{align*}
\MoveEqLeft
\big\langle A_p^*(y),x \big\rangle_{S^{p^*}_I,S^p_I}
=\big\langle (\partial_{\alpha,q,p})^*\partial_{\alpha,q,p^*}(y),x \big\rangle_{S^{p^*}_I,S^p_I}\\
&\ov{\eqref{crochet-duality}}{=} \big\langle \partial_{\alpha,q,p^*}(y),\partial_{\alpha,q,p}(x) \big\rangle_{S^{p^*}_I,S^p_I}
\ov{\eqref{crochet-duality}}{=}\big\langle y,(\partial_{\alpha,q,p^*})^*\partial_{\alpha,q,p}(x) \big\rangle_{S^{p^*}_I,S^p_I}.
\end{align*}
Since $\M_{I,\fin}$ is a core for $A_p^*=A_{p^*}$ by definition, this implies \cite[Problem 5.24]{Kat1} that $x \in \dom A_p^{**}=A_p$ and that $A_p(x) =  (\partial_{\alpha,q,p^*})^*\partial_{\alpha,q,p}(x)$.  
\end{proof}

%
Now, we show that the noncommutative gradient $\partial_{\alpha,q,p}$ commutes with the semigroup and the resolvents of its generator.

\begin{lemma}
\label{Lemma-commuting-Schur-1}
Suppose $1<p<\infty$ and $-1 \leq q \leq 1$. If $x \in \dom \partial_{\alpha,q,p}$ and $t \geq 0$, then $T_{t,p}(x)$ belongs to $\dom \partial_{\alpha,q,p}$ and we have 
\begin{equation}
\label{T_t-et-derivations}
\big(\Id_{\L^p(\Gamma_q(H))} \ot T_{t,p}\big) \partial_{\alpha,q,p}(x) 
=\partial_{\alpha,q,p} T_{t,p}(x).
\end{equation}
\end{lemma}

\begin{proof}
For any $i,j \in I$, we have
\begin{align*}
\MoveEqLeft
\big(\Id_{\L^p(\Gamma_q(H))} \ot T_{t,p}\big) \partial_{\alpha,q,p}(e_{ij})           
\ov{\eqref{def-delta-alpha}}{=} \big(\Id_{\L^p(\Gamma_q(H))} \ot T_{t,p}\big)\big(s_q(\alpha_i-\alpha_j) \ot e_{ij}\big) \\
&=\e^{-t\norm{\alpha_i -\alpha_j}^2} s_q(\alpha_i-\alpha_j) \ot e_{ij}
\ov{\eqref{def-delta-alpha}}{=} \e^{-t\norm{\alpha_i -\alpha_j}^2}\partial_{\alpha,q,p}(e_{ij})
=\partial_{\alpha,q,p}\big(\e^{-t\norm{\alpha_i -\alpha_j}^2}e_{ij}\big) \\
&=\partial_{\alpha,q,p} T_{t,p}(e_{ij}).
\end{align*}
So by linearity, the equality \eqref{T_t-et-derivations} is true for elements of $\M_{I,\fin}$. Now consider some $x \in \dom \partial_{\alpha,q,p}$. By \cite[p.~166]{Kat1}, since $\partial_{\alpha,q,p}$ is the closure of $\partial_{\alpha,q} \co \M_{I,\fin} \subseteq S^p_I \to \L^p(\Gamma_q(H) \otvn \B(\ell^2_I))$, there exists a sequence $(x_n)$ of elements of $\M_{I,\fin}$ converging to $x$ in $S^p_I$ such that the sequence $(\partial_{\alpha,q,p}(x_n))$ converges to $\partial_{\alpha,q,p}(x)$. The complete boundedness of $T_{t,p} \co S^p_I\to S^p_I$ implies by  \cite[p.~984]{Jun} that we have a bounded operator $\Id_{\L^p(\Gamma_q(H))} \ot T_{t,p} \co \L^p(\Gamma_q(H) \otvn \B(\ell^2_I)) \to \L^p(\Gamma_q(H) \otvn \B(\ell^2_I))$. We infer that in $S^p_I$ and $\L^p(\Gamma_q(H) \otvn \B(\ell^2_I))$, we have 
$$
T_{t,p}(x_n) \xra[n \to +\infty]{} T_{t,p}(x)
\quad \text{and} \quad 
\big(\Id_{\L^p(\Gamma_q(H))} \ot T_{t,p}\big) \partial_{\alpha,q,p}(x_n)
\xra[n \to +\infty]{} \big(\Id_{\L^p(\Gamma_q(H))} \ot T_{t,p}\big)\partial_{\alpha,q,p}(x).
$$ 
For any integer $n$, we have $\big(\Id_{\L^p(\Gamma_q(H))} \ot T_{t,p}\big) \partial_{\alpha,q,p}(x_n) 
=\partial_{\alpha,q,p} T_{t,p}(x_n)$ by the first part of the proof. Since the left-hand side converges, we obtain that the sequence $(\partial_{\alpha,q,p} T_{t,p}(x_n))$ converges to $(\Id_{\L^p(\Gamma_q(H))} \ot T_{t,p})\partial_{\alpha,q,p}(x)$ in $\L^p(\Gamma_q(H) \otvn \B(\ell^2_I))$. Since each $T_{t,p}(x_n)$ belongs to $\dom \partial_{\alpha,q,p}$, the closedness of $\partial_{\alpha,q,p}$ and \eqref{Def-operateur-ferme} shows that $T_{t,p}(x)$ belongs to $\dom \partial_{\alpha,q,p}$ and that $\partial_{\alpha,q,p} T_{t,p}(x)=\big(\Id_{\L^p(\Gamma_q(H))} \ot T_{t,p}\big) \partial_{\alpha,q,p}(x)$. 
\end{proof}

\begin{prop}
\label{Prop-commuting}
Let $1 < p < \infty$ and $-1 \leq q \leq 1$.
For any $s \geq 0$ and any $x \in \dom \partial_{\alpha,q,p}$, we have $\big(\Id_{S^p_I}+sA_p\big)^{-1}x \in \dom \partial_{\alpha,q,p}$ and 
\begin{equation}
\label{commuting-deriv}
\big(\Id \ot (\Id+sA_p)^{-1}\big)\partial_{\alpha,q,p}(x)
=\partial_{\alpha,q,p}\big(\Id+sA_p\big)^{-1}(x).
\end{equation}
\end{prop}

\begin{proof}
Note that for any $s>0$ and any $x \in S^p_I$ (resp. $x \in \dom \partial_{\alpha,q,p}$) the continuous functions $\R^+ \to S^p_I$, $t \mapsto \e^{-ts^{-1}}T_{t,p}(x)$ and $\R^+ \to \L^p(\Gamma_q(H) \otvn \B(\ell^2_I))$, $t \mapsto \e^{-ts^{-1}}(\Id_{\L^p(\Gamma_q(H))} \ot T_{t,p})\partial_{\alpha,q,p}(x)$ where $x \in \dom \partial_{\alpha,q,p}$ are Bochner integrable. If $t > 0$ and if $x \in \dom \partial_{\alpha,q,p}$, taking Laplace transforms on both sides of \eqref{T_t-et-derivations} and using \cite[Theorem 1.2.4]{HvNVW2} and the closedness of $\partial_{\alpha,q,p}$ in the penultimate equality, we obtain that $\int_{0}^{\infty} \e^{-ts^{-1}}T_{t,p}(x) \d t$ belongs to $\dom \partial_{\alpha,q,p}$ and that
\begin{align*}
\MoveEqLeft
\big(\Id \ot (s^{-1}\Id+A_p)^{-1}\big)\partial_{\alpha,q,p}(x) 
=-\big(-s^{-1}\Id-(\Id \ot A_p)\big)^{-1}\partial_{\alpha,q,p}(x)\\
&\ov{\eqref{Resolvent-Laplace}}{=} \int_{0}^{\infty} \e^{-s^{-1}t}(\Id \ot T_{s,p})\partial_{\alpha,q,p}(x) \d t 
\ov{\eqref{T_t-et-derivations}}{=} \int_{0}^{\infty} \e^{-s^{-1}t}\partial_{\alpha,q,p}T_{t,p}(x) \d t\\
&=\partial_{\alpha,q,p}\bigg(\int_{0}^{\infty} \e^{-s^{-1}t}T_{t,p}(x) \d t\bigg)
\ov{\eqref{Resolvent-Laplace}}{=} 
=\partial_{\alpha,q,p}\big(s^{-1}\Id+A_p\big)^{-1}(x).    
\end{align*}
We deduce the desired identity by multiplying by $s^{-1}$.
\end{proof}

%

We have the following analogue of Proposition \ref{Prop-R-gradient-bounds-Fourier} which can be proved in a similar manner.

\begin{prop}
\label{Prop-R-gradient bounds} 
Suppose $1<p<\infty$ and $-1 \leq q \leq 1$. The family 
\begin{equation}
\label{R-gradient-bounds}
\Big\{t\partial_{\alpha,q,p}(\Id+t^2A_p)^{-1}: t>0 \Big\}
\end{equation} 
of operators of $\B(S^p_I,\L^p(\Gamma_q(H) \otvn \B(\ell^2_I)))$ is $R$-bounded.
\end{prop}

\begin{proof}
Note that the operator $\partial_{\alpha,q,p}A_p^{-\frac{1}{2}} \co \SpI \to \L^p(\Gamma_q(H) \otvn \B(\ell^2_I))$ is bounded by \eqref{Equivalence-square-root-domaine-Schur}. Suppose $t >0$. A standard functional calculus argument gives
\begin{align}
\label{Divers-987-Proposition44}
\MoveEqLeft
t\partial_{\alpha,q,p}(\Id+t^2A_p)^{-1}            
=\partial_{\alpha,q,p}A_p^{-\frac{1}{2}}\Big((t^2A_p)^{\frac{1}{2}}(\Id+t^2A_p)^{-1} \Big).
\end{align} 
By \cite{Arh1} \cite{Arh4}, note that $A_p$ has a bounded $\HI(\Sigma_\theta)$-functional calculus for some $0< \theta < \frac{\pi}{2}$. Moreover, by \cite[p.~136]{HvNVW2}, the Banach space $S^p_I$ has the triangular contraction property $(\Delta)$. We deduce by \cite[Theorem 10.3.4 (2)]{HvNVW2} that the operator $A_p$ is $R$-sectorial. By \cite[Example 10.3.5]{HvNVW2} applied with $\alpha=\frac{1}{2}$ and $\beta=1$, we infer that the set
$$
\big\{(t^2A_p)^{\frac{1}{2}}(\Id+t^2A_p)^{-1}: t>0\big\}
$$
of operators of $\B(S^p_I)$ is $R$-bounded. Recalling that a singleton is $R$-bounded by \cite[Example 8.1.7]{HvNVW2}, we obtain by composition \cite[Proposition 8.1.19 (3)]{HvNVW2} that the set
$$
\Big\{ \partial_{\alpha,q,p}A_p^{-\frac{1}{2}}\Big((t^2A_p)^{\frac{1}{2}}(\Id+t^2A_p)^{-1} \Big) : t>0 \Big\}
$$
of operators of $\B(\S^p_I,\L^p(\Gamma_q(H) \otvn \B(\ell^2_I)))$ is $R$-bounded. Hence with \eqref{Divers-987-Proposition44} we conclude that the subset \eqref{R-gradient-bounds} is $R$-bounded.
\end{proof}

Our Hodge-Dirac operator in \eqref{Def-Dirac-operator-Schur} below will be constructed with $\partial_{\alpha,q,p}$ and the unbounded operator $(\partial_{\alpha,q,p^*})^*|\ovl{\Ran \partial_{\alpha,q,p}}$. Note that the latter is by definition an unbounded operator on the Banach space $\ovl{\Ran \partial_{\alpha,q,p}}$ having domain $\dom (\partial_{\alpha,q,p^*})^* \cap \ovl{\Ran \partial_{\alpha,q,p}}$.

\begin{lemma}
\label{Lemma-rest-closed-dense}
Let $1 < p < \infty$ and $-1\leq q \leq 1$. The operator $(\partial_{\alpha,q,p^*})^*|\ovl{\Ran \partial_{\alpha,q,p}}$ is densely defined and is closed. More precisely, the subspace $\partial_{\alpha,q,p}(\M_{I,\fin})$ of $\dom (\partial_{\alpha,q,p^*})^*$ is dense in $\ovl{\Ran \partial_{\alpha,q,p}}$.
\end{lemma}

\begin{proof}
Let $y \in \ovl{\Ran \partial_{\alpha,q,p}}$. Let $\epsi >0$. There exists $x \in \dom \partial_{\alpha,q,p}$ such that $\norm{y-\partial_{\alpha,q,p}(x)} < \epsi$. By Proposition \ref{Prop-derivation-closable}, there exists $x_\fin \in \M_{I,\fin}$ such that $\norm{x-x_\fin}_{S^p_I} < \epsi$ and 
$$
\norm{\partial_{\alpha,q,p}(x)-\partial_{\alpha,q,p}(x_\fin)}_{\L^p(\Gamma_q(H) \otvn \B(\ell^2_I))} < \epsi.
$$ 
We deduce that $\norm{y-\partial_{\alpha,q,p}(x_\fin)}_{\L^p(\Gamma_q(H) \otvn \B(\ell^2_I))} < 2 \epsi$. By Proposition \ref{Prop-lien-generateur-partial}, $\partial_{\alpha,q,p}(\M_{I,\fin})$ is a subspace of $\dom (\partial_{\alpha,q,p^*})^*$. So $\partial_{\alpha,q,p}(x_\fin)$ belongs to $\dom (\partial_{\alpha,q,p^*})^*$. 

Since $(\partial_{\alpha,q,p^*})^*$ is closed, the assertion on closedness is (really) obvious.
\end{proof}

According to Lemma \ref{Lemma-commuting-Schur-1}, $\Id \ot T_t$ leaves $\Ran \partial_{\alpha,q,p}$ invariant for any $t \geq 0$, so by continuity of $\Id \ot T_t$ also leaves $\ovl{\Ran \partial_{\alpha,q,p}}$ invariant. By \cite[pp.~60-61]{EnN1}, we can consider the generator
\begin{equation}
\label{Def-Bp-Schur}
B_p 
\ov{\mathrm{def}}{=} (\Id_{\L^p(\Gamma_q(H))} \ot A_p)|\ovl{\Ran \partial_{\alpha,q,p}}.
\end{equation}
of the restriction of $(\Id \ot T_t)_{t \geq 0}$ on $\ovl{\Ran \partial_{\alpha,q,p}}$.

\begin{lemma}
\label{lem-Bp-is-sectorial}
Let $1 < p < \infty$ and $-1 \leq q \leq 1$. The operator $B_p$ is injective and sectorial on $\ovl{\Ran \partial_{\alpha,q,p}}$.
\end{lemma}

\begin{proof}
The operator $B_p$ is sectorial of type $\frac{\pi}{2}$ by e.g. \cite[p.~25]{JMX}. 

For the injectivity, we note that once $B_p$ is known to be sectorial on a reflexive space, we have the projection \cite[(10.5)]{HvNVW2} onto the kernel of $B_p$ given by the strong limit 
$$
P
= \lim_{\lambda \to 0^+} \lambda  (\lambda + \Id \ot A_p)^{-1}|_{\ovl{\Ran \partial_{\alpha,q,p}}}.
$$ 
It will show that $P(y) = 0$ for $y \in \partial_{\alpha,q,p}(\M_{I,\fin})$. Indeed, by linearity, we can assume that $y = \partial_{\alpha,q}(e_{ij}) = s_q(\alpha_i - \alpha_j) \ot e_{ij}$ for some $i,j \in I$ with in addition $\alpha_i \neq \alpha_j$. We claim that
\begin{equation}
\label{equ-claim-projection-null-space}
 (\lambda + \Id \ot A_p)^{-1}|_{\ovl{\Ran \partial_{\alpha,q,p}}} \big(s_q(\alpha_i - \alpha_j) \ot e_{ij} \big)
= \frac{1}{\lambda + \norm{\alpha_i - \alpha_j}^2} s_q(\alpha_i - \alpha_j) \ot e_{ij}.
\end{equation}
To see \eqref{equ-claim-projection-null-space}, it suffices to calculate 
\begin{align*}
\MoveEqLeft
(\lambda + \Id \ot A_p)\bigg(\frac{1}{\lambda + \norm{\alpha_i - \alpha_j}^2} s_q(\alpha_i - \alpha_j) \ot e_{ij} \bigg) \\
&= \frac{1}{\lambda + \norm{\alpha_i - \alpha_j}^2}\big( \lambda s_q(\alpha_i - \alpha_j) \ot e_{ij} + \norm{\alpha_i - \alpha_j}^2 s_q(\alpha_i - \alpha_j) \ot e_{ij} \big)
= s_q(\alpha_i - \alpha_j) \ot e_{ij}.            
\end{align*} 
Now $P(y) = 0$ follows from $\frac{\lambda}{\lambda + \norm{\alpha_i - \alpha_j}^2} \xra[\lambda \to 0^+]{} \frac{0}{\norm{\alpha_i - \alpha_j}^2} = 0$ since $\norm{\alpha_i - \alpha_j}^2 \neq 0$. Since $\M_{I,\fin}$ is a core for $\partial_{\alpha,q,p}$, it is (really) not difficult to prove that $P(y) = 0$ is also true for any $y \in \Ran \partial_{\alpha,q,p}$. Finally by continuity of $P$, we deduce $P=0$. Thus the operator $B_p$ is injective.
\end{proof}

\begin{lemma}
\label{Petit-lemma-Schur}
If $z$ belongs to $\L^p(\Gamma_q(H)) \ot \M_{I,\fin} \cap \ovl{\Ran \partial_{\alpha,q,p}}$ then $z \in \partial_{\alpha,q,p}(\M_{I,\fin})$.
\end{lemma}

\begin{proof}
Let $(z_n)$ be a sequence in $\dom \partial_{\alpha,q,p}$ such that $\partial_{\alpha,q,p}(z_n) \to z$. Then for $J$ sufficiently large such that $(\Id \ot \Tron_J)(z) = z$,
\begin{align*}
\MoveEqLeft
\norm{z - \partial_{\alpha,q,p}\Tron_J(z_n)}_{\L^p} 
\ov{\eqref{commute-troncature-Schur}}{=} \norm{z - (\Id \ot \Tron_J) \partial_{\alpha,q,p}(z_n)}_{\L^p}  \\
&= \norm{(\Id \ot \Tron_J)(z - \partial_{\alpha,q,p}(z_n))}_{\L^p}  
\leq \norm{z - \partial_{\alpha,q,p}(z_n)}_{\L^p}  
\xra[n \to +\infty]{} 0.
\end{align*}
Therefore, the sequence $(\partial_{\alpha,q,p} \Tron_J(z_n))$ is convergent to $z$ for $n \to \infty$. Write $\Tron_J(z_n) = \sum_{i,j \in J} z_{n,ij} e_{ij}$ for some $z_{n,ij} \in \C$. Since for any $i,j \in J$, the Schur multiplier $\Tron_{ij}$ associated with the matrix $e_{ij}$ is completely bounded, $(\Id \ot \Tron_{ij}) \partial_{\alpha,q,p} \Tron_J (z_n) = z_{n,ij} s_q(\alpha_i - \alpha_j) \ot e_{ij}$ is convergent for $n \to \infty$. Thus, either $\alpha_i - \alpha_j = 0$ or $(z_{n,ij})_n$ is a convergent scalar sequence with limit, say, $z_{ij} \in \C$. We infer that 
$$
\partial_{\alpha,q,p}\Tron_J(z_n) 
=\sum_{i,j \in J} z_{n,ij} s_q(\alpha_i - \alpha_j) \ot e_{ij}
\xra[n \to +\infty]{} \sum_{i,j \in J} z_{ij} s_q(\alpha_i-\alpha_j) \ot e_{ij} 
\ov{\eqref{def-delta-alpha}}{=} \partial_{\alpha,q,p}\bigg(\sum_{i,j \in J} z_{ij} e_{ij}\bigg).
$$ 
Thus $z=\partial_{\alpha,q,p}(\sum_{i,j \in J} z_{ij} e_{ij})$ belongs to $\partial_{\alpha,q,p}(\M_{I,\fin})$. 
\end{proof}

\begin{prop}
\label{Prop-core-1}
Let $1 < p < \infty$ and $-1 \leq q \leq 1$.
\begin{enumerate}
\item $\L^p(\Gamma_q(H)) \ot \M_{I,\fin}$ is a core of $(\partial_{\alpha,q,p^*})^*$. Furthermore, if $y \in \dom(\partial_{\alpha,q,p^*})^*$ and if $J$ is a subset of $I$, we have
\begin{equation}
\label{Interminable-1245}
(\partial_{\alpha,q,p^*})^*(\Id \ot \Tron_J)(y) 
= \Tron_J (\partial_{\alpha,q,p^*})^*(y).
\end{equation}

\item $\partial_{\alpha,q,p}(\M_{I,\fin})$ is a core of $(\partial_{\alpha,q,p^*})^*|\ovl{\Ran \partial_{\alpha,q,p}}$.
\item $\partial_{\alpha,q,p}(\M_{I,\fin})$ is equally a core of $\partial_{\alpha,q,p}(\partial_{\alpha,q,p^*})^*|\ovl{\Ran \partial_{\alpha,q,p}}$.
\item $\partial_{\alpha,q,p}(\M_{I,\fin})$ is equally a core of $B_p$.
\end{enumerate}
\end{prop}

\begin{proof}
1. Recall that $\L^p(\Gamma_q(H)) \ot \M_{I,\fin}$ is a subset of $\dom(\partial_{\alpha,q,p^*})^*$, see \eqref{Adjoint-partial-Schur}. Let $y \in \dom(\partial_{\alpha,q,p^*})^*$. Then $(\Id_{\L^p(\Gamma_q(H))} \ot \Tron_J)(y)$ belongs to $\L^p(\Gamma_q(H)) \ot \M_{I,\fin}$. It remains to show that $(\Id_{\L^p(\Gamma_q(H))} \ot \Tron_J)(y)$ converges to $y$ in the graph norm. Since $\Tron_J$ converges strongly to $\Id_{S^p_I}$ and $\Id_{\L^p(\Gamma_q(H))} \ot \Tron_J$ is contractive, we deduce that $(\Id_{\L^p(\Gamma_q(H))} \ot \Tron_J)(y)$ converges to $y$ in $\L^p(\Gamma_q(H) \otvn \B(\ell^2_I))$. For any $i,j \in I$, we have 
\begin{align*}
\MoveEqLeft
\big\langle (\partial_{\alpha,q,p^*})^* (\Id \ot \Tron_J)(y), e_{ij} \big\rangle 
=\big\langle (\Id \ot \Tron_J)(y), \partial_{\alpha,q,p^*}(e_{ij}) \big\rangle 
=\big\langle y , (\Id \ot \Tron_J)(\partial_{\alpha,q,p^*}(e_{ij})) \big\rangle \\
&\ov{\eqref{commute-troncature-Schur}}{=} \big\langle y , \partial_{\alpha,q,p^*}\Tron_J(e_{ij}) \big\rangle 
=\big\langle (\partial_{\alpha,q,p^*})^*(y) , \Tron_J(e_{ij}) \big\rangle 
=\big\langle \Tron_J (\partial_{\alpha,q,p^*})^*(y) , e_{ij} \big\rangle.
\end{align*}
By linearity and density, we infer that
$$
(\partial_{\alpha,q,p^*})^*(\Id_{\L^p(\Gamma_q(H))} \ot \Tron_J)(y) 
= \Tron_J (\partial_{\alpha,q,p^*})^*(y)
$$ 
which converges to $(\partial_{\alpha,q,p^*})^*(y)$ in $S^p_I$.

2. Let $y \in \dom (\partial_{\alpha,q,p^*})^*|\ovl{\Ran \partial_{\alpha,q,p}}$ and $\epsi > 0$. By the proof of the first point, we already know that for some large enough finite subset $J$ of $I$, the element $z \ov{\mathrm{def}}{=} (\Id \ot \Tron_J)(y)$ belonging to $\L^p(\Gamma_q(H)) \ot \M_{I,\fin}$ satisfies 
$$
\norm{y-z}_{\L^p} < \epsi 
\quad \text{and} \quad 
\norm{(\partial_{\alpha,q,p^*})^*(y) - (\partial_{\alpha,q,p^*})^*(z)}_{S^p_I} < \epsi.
$$ 
We claim that $z$ belongs to $\ovl{\Ran \partial_{\alpha,q,p}}$. Indeed, for $\delta > 0$, since $y \in \ovl{\Ran \partial_{\alpha,q,p}}$, there is some $a \in \dom \partial_{\alpha,q,p}$ such that $\norm{y - \partial_{\alpha,q,p}(a)}_{\L^p} < \delta$. But then $\norm{z - (\Id \ot \Tron_J)\partial_{\alpha,q,p}(a)}_{\L^p} = \norm{(\Id \ot \Tron_J)(y - \partial_{\alpha,q,p}(a))}_{\L^p} \leq \norm{y - \partial_{\alpha,q,p}(a)}_{\L^p} < \delta$ and $(\Id \ot \Tron_J) \partial_{\alpha,q,p}(a) = \partial_{\alpha,q,p} \Tron_J(a)$ belongs again to $\Ran \partial_{\alpha,q,p}$. Letting $\delta \to 0$, it follows that $z$ belongs to $\ovl{\Ran \partial_{\alpha,q,p}}$. 
By Lemma \ref{Petit-lemma-Schur}, we conclude that $z \in \partial_{\alpha,q,p}(\M_{I,\fin})$.

3. Take $y \in \dom \partial_{\alpha,q,p} (\partial_{\alpha,q,p^*})^*|\ovl{\Ran \partial_{\alpha,q,p}}$, that is, $y \in \dom(\partial_{\alpha,q,p^*})^* \cap \ovl{\Ran \partial_{\alpha,q,p}}$ such that $(\partial_{\alpha,q,p^*})^*(y) \in \dom \partial_{\alpha,q,p}$. We have $(\Id \ot \Tron_J)(y) \to y$ and $(\Id \ot \Tron_J)(y) \in \partial_{\alpha,q,p}(\M_{I,\fin})$ by the proof of the second part. Moreover, we have
\begin{align*}
\MoveEqLeft
\partial_{\alpha,q,p} (\partial_{\alpha,q,p^*})^*( \Id \ot \Tron_J)(y)           
\ov{\eqref{Interminable-1245}}{=} \partial_{\alpha,q,p} \Tron_J(\partial_{\alpha,q,p^*})^*(y) 
\ov{\eqref{commute-troncature-Schur}}{=} (\Id \ot \Tron_J)\partial_{\alpha,q,p} (\partial_{\alpha,q,p^*})^*(y) \\
& \xra[\ j \ ]{} \partial_{\alpha,q,p} (\partial_{\alpha,q,p^*})^*(y).
\end{align*}   

4. Note that $\partial_{\alpha,q,p}(\M_{I,\fin})$ is a dense subspace of $\ovl{\Ran \partial_{\alpha,q,p}}$ which is clearly a subspace of $\dom B_p$ and invariant under each operator $(\Id \ot T_t)|\ovl{\Ran \partial_{\alpha,q,p}}$ by Lemma \ref{Lemma-commuting-Schur-1}. By Lemma \ref{Lemma-core-semigroup}, we deduce that $\partial_{\alpha,q,p}(\M_{I,\fin})$ is a core of $(\Id \ot A_p)|\ovl{\Ran \partial_{\alpha,q,p}} \ov{\eqref{Def-Bp-Schur}}{=} B_p$.

\end{proof}

\begin{prop}
\label{Prop-commuting-2}
Let $1 < p < \infty$ and $-1 \leq q \leq 1$.
\begin{enumerate}
\item For any $s> 0$, the operator $\big(\Id_{}+sA_p\big)^{-1}(\partial_{\alpha,q,p^*})^*$
induces a bounded operator on $\ovl{\Ran \partial_{\alpha,q,p}}$.

\item For any $s \geq 0$ and any $y \in \ovl{\Ran \partial_{\alpha,q,p}} \cap  \dom(\partial_{\alpha,q,p^*})^*$, the element $\big(\Id \ot (\Id+sA_p)^{-1}\big)(y)$ belongs to $\dom (\partial_{\alpha,q,p^*})^*$ and
\begin{equation}
\label{commuting-deriv2}
\big(\Id+sA_p\big)^{-1}(\partial_{\alpha,q,p^*})^*(y) 
=(\partial_{\alpha,q,p^*})^*\big(\Id \ot (\Id+sA_p)^{-1}\big)(y).
\end{equation}

\item For any $t \geq 0$ and any $y \in \ovl{\Ran \partial_{\alpha,q,p}} \cap \dom(\partial_{\alpha,q,p^*})^*$, the element $(\Id \ot T_t)(y)$ belongs to $\dom ( \partial_{\alpha,q,p^*})^*$ and
$$
T_t(\partial_{\alpha,q,p^*})^*(y) 
=(\partial_{\alpha,q,p^*})^*(\Id \ot T_t)(y).
$$
\end{enumerate}
\end{prop}

\begin{proof}
1. Note that $(\Id + s A_p)^{-1}(\partial_{\alpha,q,p^*})^* \ov{\eqref{Adjoint-product-unbounded}}{\subseteq} \big(\partial_{\alpha,q,p^*}(\Id+sA_{p^*})^{-1}\big)^*$. By Proposition \ref{Prop-R-gradient bounds}, the operator $\big(\partial_{\alpha,q,p^*}(\Id+sA_{p^*})^{-1}\big)^*$ is bounded. By Lemma \ref{Lemma-rest-closed-dense}, the subspace $\partial_{\alpha,q,p}(\M_{I,\fin})$ of $\dom (\partial_{\alpha,q,p^*})^*$ is dense in $\ovl{\Ran \partial_{\alpha,q,p}}$. Now, the conclusion is immediate.

2. By Proposition \ref{Prop-lien-generateur-partial}, for any $x \in \dom A_p$ we have $x \in \dom \partial_{\alpha,q,p}$ and $\partial_{\alpha,q,p}(x) \in \dom (\partial_{\alpha,q,p^*})^*$. Moreover, for all $t > 0$ we have 
\begin{align*}
\MoveEqLeft
T_t (\partial_{\alpha,q,p^*})^* \partial_{\alpha,q,p}(x) 
\ov{\eqref{eq-lien-generateur-partial}}{=} T_tA_p(x)
\ov{\eqref{A-et-Tt-commute}}{=} A_p T_t(x) \\
&\ov{\eqref{eq-lien-generateur-partial}}{=} (\partial_{\alpha,q,p^*})^* \partial_{\alpha,q,p} T_t(x)
\ov{\eqref{T_t-et-derivations}}{=} (\partial_{\alpha,q,p^*})^* (\Id \ot T_t) \partial_{\alpha,q,p}(x).
\end{align*}
By taking Laplace transforms and using the closedness of $(\partial_{\alpha,q,p^*})^*$, we deduce that the element $(\Id \ot (\Id + s A_p))^{-1} \partial_{\alpha,q,p}(x)$ belongs to $\dom (\partial_{\alpha,q,p^*})^*$ for any $s \geq 0$ and that
\begin{equation}
\label{eq:sec-1-Proposition46}
(\Id + s A_p)^{-1} (\partial_{\alpha,q,p^*})^* \partial_{\alpha,q,p}(x)
=(\partial_{\alpha,q,p^*})^* (\Id \ot (\Id +s A_p))^{-1} \partial_{\alpha,q,p}(x).
\end{equation}

Let $y \in \ovl{\Ran \partial_{\alpha,q,p}} \cap \dom(\partial_{\alpha,q,p^*})^*$. By Lemma \ref{Lemma-rest-closed-dense}, there exists a sequence $(x_n)$ of $\M_{I,\fin}$ such that $\partial_{\alpha,q,p}(x_n) \to y$. 
We have $(\Id \ot (\Id +sA_p))^{-1}\partial_{\alpha,q,p}(x_n) \to (\Id \ot (\Id + s A_p))^{-1}(y)$. Since each $x_n$ belongs to $\dom A_p$, using the first point in the passage to the limit we deduce that
\begin{align*}
\MoveEqLeft
(\partial_{\alpha,q,p^*})^* (\Id \ot (\Id +s A_p))^{-1} \partial_{\alpha,q,p}(x_n)   
\ov{\eqref{eq:sec-1-Proposition46}}{=} (\Id + s A_p)^{-1} (\partial_{\alpha,q,p^*})^* \partial_{\alpha,q,p}(x_n) \\
&\xra[n \to +\infty]{} (\Id + s A_p)^{-1} (\partial_{\alpha,q,p^*})^*(y).
\end{align*}
Since $(\partial_{\alpha,q,p^*})^*$ is closed, we infer by \eqref{Def-operateur-ferme} that $(\Id \ot (\Id +sA_p))^{-1}(y)$ belongs to $\dom (\partial_{\alpha,q,p^*})^*$ and that
\begin{equation}
\label{Commute-Schur-435}
(\partial_{\alpha,q,p^*})^*(\Id \ot (\Id +s A_p))^{-1}(y)
=(\Id + s A_p)^{-1} (\partial_{\alpha,q,p^*})^*(y)
\end{equation}
Thus \eqref{commuting-deriv2} follows.

3. Let $y \in \ovl{\Ran \partial_{\alpha,q,p}} \cap \dom(\partial_{\alpha,q,p^*})^*$. If $t \geq 0$, note that 
$$
\bigg(\Id \ot \bigg(\Id + \frac{t}{n} A_p\bigg)\bigg)^{-n}(y)
\xra[n \to +\infty]{\eqref{Widder}} (\Id \ot T_t)(y).
$$
Repeating the above commutation relation \eqref{commuting-deriv2} together with the observation that $(\Id \ot (\Id + s A_p))^{-1}$ maps $\partial_{\alpha,q,p}(\M_{I,\fin})$ into itself hence by continuity $\ovl{\Ran \partial_{\alpha,q,p}}$ into itself yields\footnote{\thefootnote. Note that we replace $s$ by $\frac{t}{n}$.} for any integer $n \geq 1$ and any $t \geq 0$
$$
(\partial_{\alpha,q,p^*})^*\bigg(\Id \ot \bigg(\Id + \frac{t}{n} A_p\bigg)\bigg)^{-n}(y)
=\bigg(\Id + \frac{t}{n} A_p\bigg)^{-n} (\partial_{\alpha,q,p^*})^*(y)
\xra[n \to +\infty]{\eqref{Widder}} T_t (\partial_{\alpha,q,p^*})^*(y).
$$
Then by the closedness of $(\partial_{\alpha,q,p^*})^*$ and again by \eqref{Widder}, $(\Id \ot T_t)(y)$ belongs to $\dom ( \partial_{\alpha,q,p^*})^*$ and $(\partial_{\alpha,q,p^*})^*(\Id \ot T_t)(y) = T_t (\partial_{\alpha,q,p^*})^*(y)$. Thus the third point follows.
\end{proof}

Proposition \ref{Prop-core-1} enables us to identify $B_p$ in terms of $\partial_{\alpha,q,p}$ and its adjoint.

\begin{prop}
\label{Prop-fundamental}
Let $1 < p < \infty$ and $-1 \leq q \leq 1$. As unbounded operators, we have
\begin{equation}
\label{relation-dur-fund}
B_p
=\partial_{\alpha,q,p}(\partial_{\alpha,q,p^*})^*|\ovl{\Ran \partial_{\alpha,q,p}}.
\end{equation}
\end{prop}

\begin{proof}
For any $i,j \in I$, we have
\begin{align*}
\MoveEqLeft
\partial_{\alpha,q,p}(\partial_{\alpha,q,p^*})^*  \partial_{\alpha,q,p}(e_{ij})
\ov{\eqref{eq-lien-generateur-partial}}{=}\partial_{\alpha,q,p}A_p(e_{ij})    
=\norm{\alpha_i-\alpha_j}_H^2 \partial_{\alpha,q,p}(e_{ij}) \\
&=\norm{\alpha_i-\alpha_j}_H^2 (s_q(\alpha_i-\alpha_j) \ot e_{ij})
=(\Id_{\L^p(\Gamma_q(H))} \ot A_p)\big(s_q(\alpha_i-\alpha_j) \ot e_{ij}\big) \\
&=(\Id_{\L^p(\Gamma_q(H))} \ot A_p)\big(\partial_{\alpha,q,p}(e_{ij})\big).
\end{align*}  
We deduce that the operators $\partial_{\alpha,q,p}(\partial_{\alpha,q,p^*})^*|\ovl{\Ran \partial_{\alpha,q,p}}$ and $(\Id_{\L^p(\Gamma_q(H))} \ot A_p)$ coincide on $\partial_{\alpha,q,p}(\M_{I,\fin})$. By Proposition \ref{Prop-core-1}, $\partial_{\alpha,q,p}(\M_{I,\fin})$ is a core for each operator. We conclude that they are equal.
\end{proof}

Note that the results of \cite{Arh1} \cite{Arh4} gives the following result.

\begin{prop}
\label{Th-functional-calculus-Schur-Ap}
Suppose $1<p<\infty$. The operators $A_p$ and $B_p$ have a bounded $\HI(\Sigma_\omega)$ functional calculus of angle $\omega$ for some $\omega <\frac{\pi}{2}$. 
\end{prop}


Suppose $1 < p <\infty$ and $-1 \leq q \leq 1$. 
We introduce the unbounded operator
\begin{equation}
\label{Def-Dirac-operator-Schur}
D_{\alpha,q,p} 
\ov{\mathrm{def}}{=}
\begin{bmatrix} 
0 & (\partial_{\alpha,q,p^*})^*|\ovl{\Ran \partial_{\alpha,q,p}}\, \\ 
\partial_{\alpha,q,p}& 0 
\end{bmatrix}
\end{equation}
on the Banach space $S^p_I \oplus_p \ovl{\Ran \partial_{\alpha,q,p}}$ defined by
\begin{equation}
\label{Def-Dirac-Schur-2}
D_{\alpha,q,p}(x,y)
\ov{\mathrm{def}}{=}
\big((\partial_{\alpha,q,p^*})^*(y),\partial_{\alpha,q,p}(x)\big), \quad x \in \dom \partial_{\alpha,q,p},\  y \in \dom (\partial_{\alpha,q,p^*})^* \cap \ovl{\Ran \partial_{\alpha,q,p}}.
\end{equation}
By Lemma \ref{Lemma-rest-closed-dense} and Proposition \ref{Prop-derivation-closable}, this operator is a closed operator and can be seen as a differential square root of the generator of the semigroup $(T_{t,p})_{t \geq 0}$. We call it the Hodge-Dirac operator of the semigroup since we have Proposition \ref{Dirac-carre-Dirac} below.

\begin{thm}
\label{Th-D-R-bisectorial-Schur}
Suppose $1<p<\infty$ and $-1 \leq q \leq 1$. The Hodge-Dirac operator $D_{\alpha,q,p}$ is R-bisectorial on $S^p_I \oplus_p \ovl{\Ran \partial_{\alpha,q,p}}$. 
\end{thm}

\begin{proof}
We will start by showing that the set $\{\i t : t \in \R, t \not=0\}$ is contained in the resolvent set of $D_{\alpha,q,p}$. We will do this by showing that $\Id-\i tD_{\alpha,q,p}$ has a two-sided bounded inverse $(\Id-\i tD_{\alpha,q,p})^{-1}$ given by
\begin{equation}
\label{Resolvent}
\begin{bmatrix} 
(\Id_{S^p_I}+t^2A_p)^{-1} & \i t(\Id_{S^p_I}+t^2A_p)^{-1}(\partial_{\alpha,q,p^*})^* \\ 
\i t\partial_{\alpha,q,p}(\Id_{S^p_I}+t^2A_p)^{-1} & \Id_{} \ot (\Id_{S^p_I}+t^2A_p)^{-1} 
\end{bmatrix} \co S^p_I \oplus_p \ovl{\Ran \partial_{\alpha,q,p}} \to S^p_I \oplus_p \ovl{\Ran \partial_{\alpha,q,p}}.
\end{equation}
By Proposition \ref{Prop-R-gradient bounds} and since the operators $A_p$ and $\Id \ot A_p$ satisfy the property \eqref{Def-R-sectoriel} of $R$-sectoriality, the four  entries are bounded. It only remains to check that this matrix defines a two-sided inverse of $\Id-\i t D_{\alpha,q,p}$. We have the following equalities of operators on $\dom D_{\alpha,q,p}$
\begin{align*}
\MoveEqLeft
\begin{bmatrix} 
(\Id_{S^p_I}+t^2A_p)^{-1} & \i t(\Id_{S^p_I}+t^2A_p)^{-1}(\partial_{\alpha,q,p^*})^* \\ 
\i t\partial_{\alpha,q,p}(\Id_{S^p_I}+t^2A_p)^{-1} & \Id_{} \ot (\Id_{S^p_I}+t^2A_p)^{-1} 
\end{bmatrix}(\Id-\i tD_{\alpha,q,p})   \\       
&\ov{\eqref{Def-Dirac-operator-Schur}}{=}
\begin{bmatrix} 
(\Id_{S^p_I}+t^2A_p)^{-1} & \i t(\Id_{S^p_I}+t^2A_p)^{-1}(\partial_{\alpha,q,p^*})^* \\ 
\i t\partial_{\alpha,q,p}(\Id_{S^p_I}+t^2A_p)^{-1} & \Id_{} \ot (\Id_{S^p_I}+t^2A_p)^{-1} 
\end{bmatrix}\begin{bmatrix} 
\Id_{S^p_I} & -\i t(\partial_{\alpha,q,p^*})^* \\ 
-\i t\partial_{\alpha,q,p} &\Id_{}
\end{bmatrix}\\
&=\left[\begin{matrix} 
(\Id+t^2A_p)^{-1}+t^2(\Id+t^2A_p)^{-1}(\partial_{\alpha,q,p^*})^*\partial_{\alpha,q,p}\\ 
\i t\partial_{\alpha,q,p}(\Id+t^2A_p)^{-1}-\i t(\Id_{} \ot (\Id+t^2A_p)^{-1})\partial_{\alpha,q,p}
\end{matrix}\right.\\
&\qquad \qquad \qquad \qquad \qquad \qquad \left.\begin{matrix} 
-\i t(\Id+t^2A_p)^{-1}(\partial_{\alpha,q,p^*})^*+\i t(\Id+t^2A_p)^{-1}(\partial_{\alpha,q,p^*})^* \\ 
t^2\partial_{\alpha,q,p}(\Id+t^2A_p)^{-1}(\partial_{\alpha,q,p^*})^*+\Id_{} \ot(\Id+t^2A_p)^{-1}
\end{matrix}\right]\\
&\ov{\eqref{eq-lien-generateur-partial}\eqref{commuting-deriv}\eqref{commuting-deriv2}}{=}
\left[\begin{matrix} 
(\Id+t^2A_p)^{-1}+t^2(\Id+t^2A_p)^{-1}A_p  \\ 
\i t\partial_{\alpha,q,p}(\Id+t^2A_p)^{-1}-\i t\partial_{\alpha,q,p}\big(\Id+t^2A_p\big)^{-1} 
\end{matrix}\right. \\
&\qquad \qquad \qquad \qquad \qquad \qquad \qquad \qquad\left.\begin{matrix} 0 \\
(t^2\partial_{\alpha,q,p}(\partial_{\alpha,q,p^*})^*+\Id_{})(\Id_{} \ot(\Id+t^2A_p)^{-1})
\end{matrix}\right] \\
&= \begin{bmatrix} 
\Id_{S^p_I} & 0 \\ 
0 & \Id_{\ovl{\Ran \partial_{\alpha,q,p}}}
\end{bmatrix}
\end{align*} 
and similarly
\begin{align*}
\MoveEqLeft
(\Id-\i tD_{\alpha,p}) 
\begin{bmatrix} 
(\Id_{S^p_I}+t^2A_p)^{-1} & \i t(\Id_{S^p_I}+t^2A_p)^{-1}(\partial_{\alpha,p^*})^* \\ 
\i t\partial_{\alpha,q,p}(\Id_{S^p_I}+t^2A_p)^{-1} & \Id_{} \ot (\Id_{S^p_I}+t^2A_p)^{-1} 
\end{bmatrix}  \\          
&=\begin{bmatrix} 
\Id & -\i t(\partial_{\alpha,q,p^*})^* \\ 
-\i t\partial_{\alpha,q,p} &\Id
\end{bmatrix}
\begin{bmatrix} 
(\Id_{S^p_I}+t^2A_p)^{-1} & \i t(\Id_{S^p_I}+t^2A_p)^{-1}(\partial_{\alpha,q,p^*})^* \\ 
\i t\partial_{\alpha,q,p}(\Id_{S^p_I}+t^2A_p)^{-1} & \Id_{} \ot (\Id_{S^p_I}+t^2A_p)^{-1} 
\end{bmatrix} 
\\
&=\left[\begin{matrix} 
(\Id_{S^p_I}+t^2A_p)^{-1}+ t^2(\partial_{\alpha,q,p^*})^*\partial_{\alpha,q,p}(\Id_{S^p_I}+t^2A_p)^{-1} &  \\ 
-\i t\partial_{\alpha,q,p}(\Id_{S^p_I}+t^2A_p)^{-1}+\i t\partial_{\alpha,q,p}(\Id_{S^p_I}+t^2A_p)^{-1} &
\end{matrix}\right.\\
&\qquad \qquad \qquad \qquad \qquad\left.
\begin{matrix}
\i t(\Id_{S^p_I}+t^2A_p)^{-1}(\partial_{\alpha,q,p^*})^*-\i t(\partial_{\alpha,q,p^*})^*\big(\Id_{} \ot (\Id_{S^p_I}+t^2A_p)^{-1}\big)\\
t^2\partial_{\alpha,q,p}(\Id_{S^p_I}+t^2A_p)^{-1}(\partial_{\alpha,q,p^*})^*+\Id_{} \ot (\Id_{S^p_I}+t^2A_p)^{-1}
\end{matrix}\right]\\
&= \begin{bmatrix} 
\Id_{S^p_I} & 0 \\ 
0 & \Id_{\ovl{\Ran \partial_{\alpha,q,p}}}
\end{bmatrix}.
\end{align*} 
It remains to show that the set $\{\i t(\i t-D_{\alpha,q,p})^{-1} : t \not=0\}=\{(\Id - \i t D_{\alpha,q,p})^{-1} : t \not=0\}$ is $R$-bounded. For this, observe that the diagonal entries of \eqref{Resolvent} are $R$-bounded by the $R$-sectoriality of $A_p$ and that of $\Id \ot A_p$. The $R$-boundedness of the other entries follows from the $R$-gradient bounds of Proposition \ref{Prop-R-gradient bounds}. Since a set of operator matrices is $R$-bounded precisely when each entry is $R$-bounded, we conclude that \eqref{Def-R-bisectoriel} is satisfied, i.e. that $D_{\alpha,q,p}$ is $R$-bisectorial.
\end{proof}

\begin{prop}
\label{Dirac-carre-Dirac}
Suppose $1<p<\infty$ and $-1 \leq q \leq 1$. We have 
\begin{equation}
\label{D-alpha-carre-egal}
D_{\alpha,q,p}^2
=\begin{bmatrix} 
A_p & 0 \\ 
0 & (\Id_{\L^p(\Gamma_q(H))} \ot A_p)|\ovl{\Ran \partial_{\alpha,q,p}}
\end{bmatrix}.	
\end{equation}
\end{prop}

\begin{proof}
By Proposition \ref{Prop-fundamental}, we have
\begin{align*}
\MoveEqLeft
\begin{bmatrix} 
A_p & 0 \\ 
0 & (\Id_{\L^p(\Gamma_q(H))} \ot A_p)|\ovl{\Ran \partial_{\alpha,q,p}}
\end{bmatrix} \\
&\ov{\eqref{eq-lien-generateur-partial}\eqref{relation-dur-fund}}{=} 
\begin{bmatrix} 
(\partial_{\alpha,q,p^*})^*\partial_{\alpha,q,p} & 0 \\ 
0 & \partial_{\alpha,q,p}(\partial_{\alpha,q,p^*})^*|\ovl{\Ran \partial_{\alpha,q,p}}\,
\end{bmatrix} 
=\begin{bmatrix} 
0 & (\partial_{\alpha,q,p^*})^*|\ovl{\Ran \partial_{\alpha,q,p}}\, \\ 
\partial_{\alpha,q,p}& 0 
\end{bmatrix}^2 \\
&\ov{\eqref{Def-Dirac-operator-Schur}}{=} D_{\alpha,q,p}^2.            
\end{align*} 
\end{proof}

Now, we can state the following main result of this subsection.

\begin{thm}
\label{Th-functional-calculus-bisector-Schur}
Suppose $1<p<\infty$. The Hodge-Dirac operator $D_{\alpha,q,p}$ has a bounded $\HI(\Sigma^\pm_\omega)$ functional calculus on a bisector, on the Banach space $S^p_I \oplus_p \ovl{\Ran \partial_{\alpha,q,p}}$. 
\end{thm}

\begin{proof}
By Proposition \ref{Th-functional-calculus-Schur-Ap}, the operator $
D_{\alpha,q,p}^2 
\ov{\eqref{D-alpha-carre-egal}}{=}\begin{bmatrix} 
A_p & 0 \\ 
0 & B_p
\end{bmatrix}$ has a bounded $\HI$ functional calculus of angle $2 \omega < \frac{\pi}{2}$. Since $D_{\alpha,q,p}$ is $R$-bisectorial by Theorem \ref{Th-D-R-bisectorial-Schur}, we deduce by Proposition \ref{Prop-liens-differents-calculs-fonctionnels} that the operator $D_{\alpha,q,p}$ has a bounded $\HI(\Sigma^\pm_\omega)$ functional calculus on a bisector.
\end{proof}

\begin{remark} \normalfont
\label{remark-sgn-Schur}
Similarly to Remark \ref{remark-sgn-Fourier}, the boundedness of the $\HI$ functional calculus of the operator $D_{\alpha,q,p}$ implies the boundedness of the Riesz transforms and this result may be thought of as a strengthening of the equivalence \eqref{Equivalence-square-root-domaine-Schur}. 
\end{remark}

\begin{remark} \normalfont
\label{remark-amplified-Riesz}
In a similar way to Remark \ref{remark-sgn-Schur}, we also obtain for $(x,y) \in \dom\partial_{\alpha,q,p} \oplus \dom(\partial_{\alpha,q,p^*})^*|\ovl{\Ran \partial_{\alpha,q,p}}$ that
\begin{align*}
\MoveEqLeft
\bnorm{(D^2_{\alpha,q,p})^{\frac12}(x,y)}_p 
\cong \bnorm{A_p^{\frac12}(x)}_p + \bnorm{(\Id_{\L^p(\Gamma_q(H))} \ot A_p^{\frac12})(y)}_p \cong \bnorm{\partial_{\alpha,q,p}(x)}_p + \bnorm{(\partial_{\alpha,q,p^*})^*(y)}_p \\
&\cong \norm{D_{\alpha,q,p}(x,y)}_p.           
\end{align*}
Moreover, we have
\begin{equation}
\label{equ-amplified-Riesz}
\bnorm{(\Id_{\L^p(\Gamma_q(H))} \ot A_p^{\frac12})(y)}_p 
\cong  \bnorm{(\partial_{\alpha,q,p^*})^*(y)}_p, \quad y \in \dom(\partial_{\alpha,q,p^*})^* \cap \ovl{\Ran \partial_{\alpha,q,p}}.
\end{equation}
\end{remark}

\begin{prop}
\label{Prop--liens-ranges}
We have $\ovl{\Ran  A_p}=\ovl{\Ran (\partial_{\alpha,q,p^*})^*}$, $\ovl{\Ran B_p}=\ovl{\Ran \partial_{\alpha,q,p}}$, $\ker A_p=\ker \partial_{\alpha,q,p}$, $\ker B_p=\ker (\partial_{\alpha,q,p^*})^*=\{0\}$ and
\begin{equation}
\label{Decompo-sympa}
S^p_I
=\ovl{\Ran (\partial_{\alpha,q,p^*})^*} \oplus \ker \partial_{\alpha,q,p}.
\end{equation}
Here, by $(\partial_{\alpha,q,p^*})^*$ we understand its restriction to $\ovl{\Ran \partial_{\alpha,q,p}}$. However, we shall see in Corollary \ref{cor-key-Hodge-decomposition} that $\Ran (\partial_{\alpha,q,p^*})^* = \Ran(\partial_{\alpha,q,p^*})^*|{\ovl{\Ran \partial_{\alpha,q,p}}}$.
\end{prop}

\begin{proof}
By \eqref{Bisec-Ran-Ker}, we have $\ovl{\Ran D_{\alpha,q,p}^2}=\ovl{\Ran D_{\alpha,q,p}}$ and $\ker D_{\alpha,q,p}^2=\ker D_{\alpha,q,p}$. It is not difficult to prove the first four equalities using \eqref{Dirac-carre-Dirac} and \eqref{Def-Dirac-Schur-2}. The last one is a consequence of the definition of $A_p$ and of \cite[p.~361]{HvNVW2}.
\end{proof}

\subsection{Extension to full Hodge-Dirac operator and Hodge decomposition}
\label{Extension-full-Hodge-Schur}

We keep the standing assumptions of the preceding subsection and thus we have a markovian semigroup $(T_t)_{t \geq 0}$ of Schur multipliers with generator $A_p$, the noncommutative gradient $\partial_{\alpha,q,p}$ and its adjoint $(\partial_{\alpha,q,p^*})^*$, together with the Hodge-Dirac operator $D_{\alpha,q,p}$. We shall now extend the operator $D_{\alpha,q,p}$ to a densely defined bisectorial operator $\D_{\alpha,q,p}$ on $S^p_I \oplus \L^p(\Gamma_q(H) \otvn \B(\ell^2_I))$ which will also be bisectorial and will have an $\HI(\Sigma_\omega^\pm)$ functional calculus on a bisector. The key will be Corollary \ref{cor-key-Hodge-decomposition} below. We let
\begin{equation}
\label{equ-full-Hodge-Dirac-operator}
\D_{\alpha,q,p} 
\ov{\mathrm{def}}{=} \begin{bmatrix} 
0 & (\partial_{\alpha,q,p^*})^* \\ 
\partial_{\alpha,q,p} & 0 
\end{bmatrix}
\end{equation}
along the decomposition $S^p_I \oplus \L^p(\Gamma_q(H) \otvn \B(\ell^2_I))$, with natural domains for $\partial_{\alpha,q,p}$ and $(\partial_{\alpha,q,p^*})^*$.

Consider the sectorial operator $A_p^{\frac12}$ on $S^p_I$. According to \cite[(10.1) p.~361]{HvNVW2}, we have the \textit{topological} direct sum decomposition $S^p_I= \ovl{\Ran A_p^{\frac12}} \oplus \ker A_p^{\frac12}$. We define the operator $R_p \ov{\mathrm{def}}{=} \partial_{\alpha,q,p} A_p^{-\frac12} \co \Ran A_p^{\frac12} \to \L^p(\Gamma_q(H) \otvn \B(\ell^2_I))$. According to  the point 4 of Proposition \ref{Prop-derivation-closable}, $R_p$ is bounded on $\Ran A_p^{\frac12}$, so extends to a bounded operator on $\ovl{\Ran A_p^{\frac12}}$. We extend it to a bounded operator $R_p \co S^p_I \to \L^p(\Gamma_q(H) \otvn \B(\ell^2_I))$, called Riesz transform, by putting $R_p|\ker A_p^{\frac12} = 0$ along the above decomposition of $S^p_I$. We equally let $R_{p^*}^* \ov{\mathrm{def}}{=} (R_{p^*})^*$.

\begin{lemma}
\label{lem-Hodge-decomposition-sum}
Let $-1 \leq q \leq 1$ and $1 < p < \infty$. Then we have the subspace sum 
\begin{equation}
\label{Decompo-inter-100}
\L^p(\Gamma_q(H) \otvn \B(\ell^2_I)) 
= \ovl{\Ran \partial_{\alpha,q,p}} + \ker (\partial_{\alpha,q,p^*})^*.
\end{equation}
\end{lemma}

\begin{proof}
Let $y \in \L^p(\Gamma_q(H) \otvn \B(\ell^2_I))$ be arbitrary. We claim that $y = R_pR_{p^*}^*(y) + (\Id - R_pR_{p^*}^*)(y)$ is the needed decomposition for \eqref{Decompo-inter-100}. Note that $R_p$ maps $\Ran A_p^{\frac12}$ into $\Ran \partial_{\alpha,q,p}$, so by boundedness, $R_p$ maps $\ovl{\Ran A_p^{\frac12}}$ to $\ovl{\Ran \partial_{\alpha,q,p}}$. Thus, we indeed have $R_pR_{p^*}^*(y) \in \ovl{\Ran \partial_{\alpha,q,p}}$. Next we claim that for any $z \in S^p_I$ and $x \in \dom \partial_{\alpha,q,p^*}$, we have
\begin{equation}
\label{equ-1-proof-lem-Hodge-decomposition-sum}
\big\langle R_p(z), \partial_{\alpha,q,p^*}(x)\big\rangle_{\L^p(\Gamma_q(H) \otvn \B(\ell^2_I)),\L^{p^*}(\Gamma_q(H) \otvn \B(\ell^2_I))} 
=\big\langle z, A_{p^*}^{\frac12}(x) \big\rangle_{S^p_I,S^{p^*}_I}.
\end{equation}
According to the decomposition $S^p_I = \ovl{\Ran A_p^{\frac12}} \oplus \ker A_p^{\frac12}$, we can write $z = \lim_{n \to +\infty} A_p^{\frac12}(z_n) + z_0$ with $z_n \in \dom A_p^{\frac12}$ and $z_0 \in \ker A_p^{\frac12}$. Then we have
\begin{align*}
\MoveEqLeft
\big\langle R_p(z), \partial_{\alpha,q,p^*}(x) \big\rangle 
=\lim_{n \to +\infty} \big\langle R_p\big( A_p^{\frac12}(z_n) + z_0\big), \partial_{\alpha,q,p^*}(x) \big\rangle
=\lim_{n \to +\infty} \big\langle \partial_{\alpha,q,p}(z_n) , \partial_{\alpha,q,p^*}(x) \big\rangle \\
&\overset{\eqref{formule-trace-trace}}{=} \lim_{n \to +\infty} \big\langle A^{\frac12}_p(z_n), A^{\frac12}_{p^*} (x) \big\rangle 
=\big\langle z - z_0 , A_{p^*}^{\frac12}(x) \big\rangle 
=\big\langle z, A^{\frac12}_{p^*}(x) \big\rangle-\big\langle z_0 , A^{\frac12}_{p^*}(x) \big\rangle \\
&=\big\langle z, A_{p^*}^{\frac12}(x) \big\rangle.
\end{align*}
Thus, \eqref{equ-1-proof-lem-Hodge-decomposition-sum} is proved. We can now conclude the proof, since for $x \in \dom \partial_{\alpha,q,p^*}$, we have
\begin{align*}
\MoveEqLeft
\big\langle (\Id - R_p R_{p^*}^*)(y) , \partial_{\alpha,q,p^*}(x) \big\rangle 
=\big\langle y , \partial_{\alpha,q,p^*}(x) \big\rangle  - \big\langle R_pR_{p^*}^*(y), \partial_{\alpha,q,p^*}(x) \big\rangle \\
&\ov{\eqref{equ-1-proof-lem-Hodge-decomposition-sum}}{=} \big\langle y , \partial_{\alpha,q,p^*}(x) \big\rangle  - \big\langle R_{p^*}^*(y), A_{p^*}^{\frac12}(x) \big\rangle 
=\big\langle y, \partial_{\alpha,q,p^*}(x) \big\rangle - \big\langle y, R_{p^*} A_{p^*}^{\frac12}(x) \big\rangle \\
&=\big\langle y, \partial_{\alpha,q,p^*}(x) \big\rangle - \big\langle y, \partial_{\alpha,q,p^*} A_{p^*}^{-\frac12} A_{p^*}^\frac12(x) \big\rangle 
= 0.
\end{align*}
By \eqref{lien-ker-image}, we infer that $\big(\Id - R_pR_{p^*}^*\big)(y)$ belongs to $\ker (\partial_{\alpha,q,p^*})^*$.
\end{proof}

In the proof of Proposition \ref{prop-Hodge-decomposition-intersection} below, we shall need some information on the Wiener-Ito chaos decomposition for $q$-Gaussians.
This is collected in the following lemma.

\begin{lemma}
\label{lem-eigenspace-q-OU}
Let $-1 \leq q \leq 1$ and $1 < p < \infty$.
\begin{enumerate}
\item There exists a completely bounded projection $\mathcal{P} \co \L^p(\Gamma_q(H)) \to \L^p(\Gamma_q(H))$ onto the closed space spanned by $\{s_q(h) : h \in H\}$. Moreover, the projections are compatible for different values of $p$.

\item For any finite subset $J$ of $I$ and $y \in \L^p(\Gamma_q(H) \otvn \B(\ell^2_I))$, $(\mathcal{P} \ot \Tron_J)(y)$ can be written as $\sum_{i,j \in J} s_q(h_{ij}) \ot e_{ij}$ for some $h_{ij} \in H$.

\item Let $J$ be a finite subset of $I$. Denoting temporarily $\mathcal{P}_p$ and $\Tron_{J,p}$ the operators $\mathcal{P}$ and $\Tron_J$ on the $p$-level, the identity mapping on $\Gauss_{q,p}(\C) \ot S^p_J$ extends to an isomorphism 
\begin{equation}
\label{equ-1-proof-prop-Hodge-decomposition-intersection}
J_{p,2,J} \co \Ran(\mathcal{P}_p \ot \Tron_{J,p}) \subseteq \L^p(\Gamma_q(H) \otvn \B(\ell^2_I)) \to \Ran(\mathcal{P}_2 \ot \Tron_{J,2}) \subseteq \L^2(\Gamma_q(H) \otvn \B(\ell^2_I)).
\end{equation}
\end{enumerate}
\end{lemma}

\begin{proof}
1. This is contained in \cite[Theorem 3.5]{JuL1}, putting there $d = 1$.
Note that the closed space spanned by $\{s_q(h),\: h \in H\}$ coincides in this case with $\mathcal{G}^1_{p,q}$ there. For the fact that the projections are compatible for different values of $p$, we refer to \cite[Proof of Theorem 3.1]{JuL1}.
See also the mapping $Q_p$ from Lemma \ref{lem-Qp-bounded-Schur}.

2. This is easy and left to the reader.

3. We have 
$$
\norm{\sum_{i,j \in J} s_q(h_{ij}) \ot e_{ij}}_{\L^p(\Gamma_q(H) \otvn \B(\ell^2_I))} \leq \sum_{i,j \in J} \norm{s_q(h_{ij})}_{\L^p(\Gamma_q(H))} 
\lesssim_{q,p} \sum_{i,j \in J} \norm{h_{ij}}_H
$$ 
(the same estimate on the $\L^2$-level). Note that for any $i_0, \: j_0 \in J$, we have a completely contractive Schur multiplier projecting onto the span of $e_{i_0j_0}$. So we have 
$$
\norm{\sum_{i,j \in J} s_q(h_{ij}) \ot e_{ij}}_{\L^p(\Gamma_q(H) \ot \B(\ell^2_I))} \gtrsim \norm{s_q(h_{i_0j_0})}_{\L^p(\Gamma_q(H))} \gtrsim_{q,p} \norm{h_{i_0j_0}}_H. 
$$
It follows that 
$$
\norm{\sum_{i,j \in J} s_q(h_{ij}) \ot e_{ij}}_{\L^p(\Gamma_q(H) \otvn \B(\ell^2_I))} \cong \sum_{i,j \in J} \norm{h_{ij}}_H \cong \norm{\sum_{i,j \in J} s_q(h_{ij}) \ot e_{ij}}_{\L^2(\Gamma_q(H) \otvn \B(\ell^2_I))}.
$$
Thus, \eqref{equ-1-proof-prop-Hodge-decomposition-intersection} follows.
\end{proof}

\begin{prop}
\label{prop-Hodge-decomposition-intersection}
Let $-1 \leq q \leq 1$ and $1 < p < \infty$. Then the subspaces from Lemma \ref{lem-Hodge-decomposition-sum} have trivial intersection, $\ovl{\Ran \partial_{\alpha,q,p}} \cap \ker (\partial_{\alpha,q,p^*})^* = \{ 0 \}$.
\end{prop}

\begin{proof}
We begin with the case $p = 2$.
According to Proposition \ref{First-spectral-triple} below, the unbounded operator $\D_{\alpha,q,2}$ is selfadjoint on $S^2_I \oplus \L^2(\Gamma_q(H) \otvn \B(\ell^2_I))$.
We thus have the orthogonal sum $\ovl{\Ran \D_{\alpha,q,2}} \oplus \ker \D_{\alpha,q,2} = S^2_I \oplus \L^2(\Gamma_q(H) \otvn \B(\ell^2_I))$. Considering vectors in the second component, that is, in $\L^2(\Gamma_q(H) \otvn \B(\ell^2_I))$, we deduce that $\ovl{\Ran \partial_{\alpha,q,2}}$ and $\ker (\partial_{\alpha,q,2})^*$ are orthogonal, hence have trivial intersection.

We turn to the case $1 < p < \infty$. We recall that for a finite subset $J$ of $I$, we have a completely contractive Schur multiplier $\Tron_J \co S^p_I \to S^p_I$ projecting onto the finite-dimensional subspace $S^p_J$ of $S^p_I$. We thus have according to Lemma \ref{lem-eigenspace-q-OU} for any finite subset $J$ of $I$ a completely bounded mapping $\mathcal{P} \ot \Tron_J = (\mathcal{P} \ot \Id_{S^p_I}) \circ (\Id_{\L^p(\Gamma_q(H))} \ot \Tron_J) \co \L^p(\Gamma_q(H) \otvn \B(\ell^2_I)) \to \L^p(\Gamma_q(H) \otvn \B(\ell^2_I))$. We claim that
\begin{align}
 \text{the subspace } \ovl{\Ran \partial_{\alpha,q,p}}\text{ is invariant under  } \Id_{} \ot \Tron_J & \label{equ-2-proof-propHodge-decomposition-intersection} \\
 \text{the subspace } \ker (\partial_{\alpha,q,p^*})^* \text{ is invariant under  } \Id_{} \ot \Tron_J & \label{equ-3-proof-propHodge-decomposition-intersection} \\
\text{the restriction of } \mathcal{P} \ot \Id_{S^p_I} \text{ on } \ovl{\Ran \partial_{\alpha,q,p}}\text{ is the identity mapping}. & \label{equ-5-proof-propHodge-decomposition-intersection}
\end{align}
For \eqref{equ-2-proof-propHodge-decomposition-intersection}, for any $i,j \in I$, note that 
$$
(\Id \ot \Tron_J) \partial_{\alpha,q,p}(e_{ij}) 
\ov{\eqref{def-delta-alpha}}{=} (\Id \ot \Tron_J)\big(s_q(\alpha_i-\alpha_j)  \ot e_{ij}\big) 
= 1_J(i)1_J(j) s_q(\alpha_i-\alpha_j) \ot e_{ij}.
$$
This element belongs to $\ovl{\Ran \partial_{\alpha,q,p}}$. By linearity and since $\M_{I,\fin}$ is a core for $\partial_{\alpha,q,p}$ according to Proposition \ref{Prop-derivation-closable}, we deduce that $\Id \ot \Tron_J$ maps $\Ran \partial_{\alpha,q,p}$ into $\ovl{\Ran \partial_{\alpha,q,p}}$. Now \eqref{equ-2-proof-propHodge-decomposition-intersection} follows from boundedness of $\Id \ot \Tron_J$. 

For \eqref{equ-3-proof-propHodge-decomposition-intersection}, note that if $x \in \M_{I,\fin}$ and $f \in \ker (\partial_{\alpha,q,p^*})^*$, then
\begin{align*}
\MoveEqLeft
\big\langle (\Id \ot \Tron_J) (f), \partial_{\alpha,q,p^*}(x) \big\rangle 
=\big\langle f, (\Id \ot \Tron_J) \partial_{\alpha,q,p^*}(x) \big\rangle 
=\big\langle f, \partial_{\alpha,q,p^*} \Tron_J(x) \big\rangle \\
&=\big\langle (\partial_{\alpha,q,p^*})^* f, \Tron_J(x) \big\rangle 
= 0.
\end{align*}
Thus $(\Id \ot \Tron_J)(f)$ belongs to $\ker (\partial_{\alpha,q,p^*})^*$ and \eqref{equ-3-proof-propHodge-decomposition-intersection} follows. For \eqref{equ-5-proof-propHodge-decomposition-intersection}, for any $i,j \in I$ we have  
$$
(\mathcal{P} \ot \Id) (\partial_{\alpha,q,p}(e_{ij})) 
\ov{\eqref{def-delta-alpha}}{=} (\mathcal{P} \ot \Id) (s_q(\alpha_i - \alpha_j) \ot e_{ij}) 
=s_q(\alpha_i -\alpha_j) \ot e_{ij} \in \Ran \partial_{\alpha,q,p}.
$$
Now use in a similar manner as before linearity, the fact that $\M_{I,\fin}$ is a core of $\partial_{\alpha,q,p}$ and boundedness of $\mathcal{P} \ot \Id$.

Now, let $y \in \ovl{\Ran \partial_{\alpha,q,p}} \cap \ker (\partial_{\alpha,q,p^*})^*$. Let $J$ be a finite subset of $I$. Then according to \eqref{equ-2-proof-propHodge-decomposition-intersection} -- \eqref{equ-5-proof-propHodge-decomposition-intersection}, we infer that $(\mathcal{P} \ot \Tron_J)(y)$ belongs\footnote{\thefootnote. Note that $(\mathcal{P} \ot \Id_{S^p_I})(y)=y$.} again to $\ovl{\Ran \partial_{\alpha,q,p}} \cap \ker (\partial_{\alpha,q,p^*})^*$. We claim that $J_{p,2,J} (\mathcal{P} \ot \Tron_J) (y)$ belongs to $\ovl{\Ran \partial_{\alpha,q,2}} \cap \ker (\partial_{\alpha,q,2})^*$, where the mapping $J_{p,2,J}$ was defined in Lemma \ref{lem-eigenspace-q-OU}. First we have $(\mathcal{P} \ot \Tron_J)(y) \in \ovl{\Ran \partial_{\alpha,q,p}}$, so that there exists a sequence $(x_n)$ in $\dom \partial_{\alpha,q,p}$ such that $(\mathcal{P} \ot \Tron_J)(y) = \lim_n \partial_{\alpha,q,p}(x_n)$. Then
\begin{align*}
\MoveEqLeft
(\mathcal{P} \ot \Tron_J)(y) 
=(\Id \ot \Tron_J) (\mathcal{P} \ot \Tron_J) (y) 
=\lim_n (\Id \ot \Tron_J) \partial_{\alpha,q,p}(x_n) \\
&\overset{\eqref{commute-troncature-Schur}}{=}\lim_n \partial_{\alpha,q,p} \Tron_J(x_n) 
=\lim_n \sum_{i,j \in J} x_{n,ij} \partial_{\alpha,q,p}(e_{ij}).
\end{align*}
Since $J_{p,2,J}$ is an isomorphism, this limit also holds in $\L^2(\Gamma_q(H) \otvn \B(\ell^2_I))$ and the element $J_{p,2,J} (\mathcal{P} \ot \Tron_J) (y) = \lim_n \sum_{i,j \in J} x_{n,ij} \partial_{\alpha,q,2} (e_{ij})$ belongs to $\ovl{\Ran \partial_{\alpha,q,2}}$. Furthermore, we have
\begin{equation}
\label{Une-der-equa}
(\mathcal{P} \ot \Tron_J)(y) 
= \sum_{i,j \in J} s_q(h_{ij}) \ot e_{ij}
\end{equation}
for some $h_{ij} \in H$. Then
$$
(\partial_{\alpha,q,2})^* J_{p,2,J} (\mathcal{P}   \ot \Tron_J)(y) 
\ov{\eqref{Une-der-equa}}{=}(\partial_{\alpha,q,2})^* \bigg(\sum_{i,j \in J} s_q(h_{ij}) \ot e_{ij}\bigg) 
=(\partial_{\alpha,q,p^*})^*\bigg(\sum_{i,j \in J} s_q(h_{ij}) \ot e_{ij}\bigg) 
=0. 
$$
We have shown that $J_{p,2,J} (\mathcal{P} \ot \Tron_J) (y)$ belongs to $\ovl{\Ran \partial_{\alpha,q,2}} \cap \ker (\partial_{\alpha,q,2})^*$. According to the beginning of the proof, the last intersection is trivial. It follows that $J_{p,2,J} (\mathcal{P} \ot \Tron_J)(y) = 0$, and since $J_{p,2,J}$ is an isomorphism, that $(\mathcal{P} \ot \Tron_J)(y) = 0$. Thus, for any finite subset $J$ of $I$ we have $(\mathcal{P} \ot \Tron_J)(y) = 0$. As $\Tron_J$ converges strongly and completely boundedly on $S^p_I$ to the identity, we infer that $(\mathcal{P} \ot \Id)(y) = 0$. But we had seen in \eqref{equ-5-proof-propHodge-decomposition-intersection} that $(\mathcal{P} \ot \Id)(y) = y$, so that $y = 0$ and we are done.
\end{proof}

Combining Lemma \ref{lem-Hodge-decomposition-sum} and Proposition \ref{prop-Hodge-decomposition-intersection}, we can now deduce the following corollary.

\begin{cor}
\label{cor-key-Hodge-decomposition}
Let $-1 \leq q \leq 1$ and $1 < p < \infty$. Then we have a topological direct sum decomposition  
\begin{equation}
\label{Cor-decompo1}
\L^p(\Gamma_q(H) \otvn \B(\ell^2_I)) 
= \ovl{\Ran \partial_{\alpha,q,p}} \oplus \ker (\partial_{\alpha,q,p^*})^*
\end{equation}
where the associated first bounded projection is $R_pR_{p^*}^*$. In particular, we have
$$
\Ran (\partial_{\alpha,q,p^*})^* 
=\Ran (\partial_{\alpha,q,p^*})^*|{\ovl{\Ran \partial_{\alpha,q,p}}}.
$$
\end{cor}

\begin{proof}
According to Lemma \ref{lem-Hodge-decomposition-sum}, the above subspaces add up to $\L^p(\Gamma_q(H)  \otvn \B(\ell^2_I))$, and according to Proposition \ref{prop-Hodge-decomposition-intersection}, the sum is direct. By \cite[Theorem 1.8.7]{KaR1}, we conclude that the decomposition is topological. In the course of the proof of Lemma \ref{lem-Hodge-decomposition-sum}, we have seen that for any $y \in \L^p(\Gamma_q(H)  \otvn \B(\ell^2_I))$, we have the suitable decomposition $y = R_pR_{p^*}^*(y) + (\Id - R_pR_{p^*}^*)(y)$. So the associated first bounded projection is $R_pR_{p^*}^*$.
\end{proof}

\begin{thm}
\label{Thm-full-operator-bisectorial}
Let $-1 \leq q \leq 1$ and $1 < p < \infty$. Consider the operator $\D_{\alpha,q,p}$ from \eqref{equ-full-Hodge-Dirac-operator}. Then $\D_{\alpha,q,p}$ is bisectorial and has a bounded $\HI(\Sigma_\omega^\pm)$ calculus.
\end{thm}

\begin{proof}
According to Corollary \ref{cor-key-Hodge-decomposition}, the space $S^p_I \oplus \L^p(\Gamma_q(H) \otvn \B(\ell^2_I))$ admits the topological direct sum decomposition $S^p_I \oplus \L^p(\Gamma_q(H) \otvn \B(\ell^2_I)) \ov{\eqref{Cor-decompo1}}{=} S^p_I \oplus \ovl{\Ran \partial_{\alpha,q,p}} \oplus \ker (\partial_{\alpha,q,p^*})^*$ into a sum of three subspaces. Along this decomposition, we can write\footnote{\thefootnote. Here the notation $(\partial_{\alpha,q,p^*})^*$ is used for the restriction of $(\partial_{\alpha,q,p^*})^*$ on the subspace $\ovl{\Ran \partial_{\alpha,q,p}}$.}
$$
\D_{\alpha,q,p}
\ov{\eqref{equ-full-Hodge-Dirac-operator}}{=}\begin{bmatrix} 
0 & (\partial_{\alpha,q,p^*})^* & 0 \\ 
\partial_{\alpha,q,p} & 0 & 0 \\ 
0 & 0 & 0 
\end{bmatrix}
\ov{\eqref{Def-Dirac-operator-Schur}}{=} \begin{bmatrix} 
D_{\alpha,q,p} & 0 \\ 
0 & 0 
\end{bmatrix}.
$$
According to Theorem \ref{Th-functional-calculus-bisector-Schur}, the operator $D_{\alpha,q,p}$ is bisectorial and has a bounded $\HI(\Sigma_\omega^\pm)$ functional calculus. So we conclude that the same thing for $\D_{\alpha,q,p}$.
\end{proof}

\begin{thm}[Hodge decomposition]
\label{Th-Hodge-decomposition}
Suppose $1<p<\infty$ and $-1 \leq q \leq 1$. If we identify $\ovl{\Ran \partial_{\alpha,q,p}}$ and $\ovl{\Ran (\partial_{\alpha,q,p^*})^*}$ as the closed  subspaces $\{0\} \oplus \ovl{\Ran \partial_{\alpha,q,p}}$ and $\ovl{\Ran (\partial_{\alpha,q,p^*})^*} \oplus \{0\}$ of $S^p_I \oplus \L^p(\Gamma_q(H) \otvn \B(\ell^2_I))$, we have 
\begin{equation}
\label{Hodge-decomposition}
S^p_I \oplus \L^p(\Gamma_q(H) \otvn \B(\ell^2_I)) 
=\ovl{\Ran \partial_{\alpha,q,p}} \oplus \ovl{\Ran (\partial_{\alpha,q,p^*})^*} \oplus \ker \D_{\alpha,q,p}. 
\end{equation}  
\end{thm}

\begin{proof}
From the definition \eqref{equ-full-Hodge-Dirac-operator}, it is obvious that $\ker \D_{\alpha,q,p}= (\ker \partial_{\alpha,q,p} \oplus \{ 0 \}) \oplus (\{ 0 \} \oplus \ker (\partial_{\alpha,q,p^*})^* )$. We deduce that
\begin{align*}
\MoveEqLeft
S^p_I \oplus \L^p(\Gamma_q(H) \otvn \B(\ell^2_I))             
\ov{\eqref{Decompo-sympa} \eqref{Cor-decompo1}}{=} (\ovl{\Ran (\partial_{\alpha,q,p^*})^*} \oplus \ker \partial_{\alpha,q,p})
\oplus (\ker (\partial_{\alpha,q,p^*})^* \oplus \ovl{\Ran \partial_{\alpha,q,p}})\\
&=(\{0\} \oplus \ovl{\Ran \partial_{\alpha,q,p}}) \oplus (\ovl{\Ran (\partial_{\alpha,q,p^*})^*} \oplus \{0\}) \oplus (\ker \partial_{\alpha,q,p} \oplus \{  0 \}) \oplus (\{ 0 \} \oplus \ker (\partial_{\alpha,q,p^*})^*) \\
&=(\{0\} \oplus \ovl{\Ran \partial_{\alpha,q,p}}) \oplus (\ovl{\Ran (\partial_{\alpha,q,p^*})^*} \oplus \{0\}) \oplus\ker \D_{\alpha,q,p}.
\end{align*}  
\end{proof}

\subsection{Independence from $H$ and $\alpha$}
\label{Sec-Indep-alpha}

In this short subsection, we again keep the standing assumptions from the two preceding subsections and have a markovian semigroup of Schur multipliers. The main topic will be the proof of Theorem \ref{prop-explicitly-dimension-free} showing that the bound of the bisectorial $\HI(\Sigma_\omega^\pm)$ functional calculus in Theorem \ref{Th-functional-calculus-bisector-Schur} and Theorem \ref{Thm-full-operator-bisectorial} comes with a constant only depending on $\omega,q,p$, but not on the markovian semigroup nor the associated Hilbert space $H$ nor $\alpha \co I \to H$. We show several intermediate lemmas before.

\begin{lemma}
Let $-1 \leq q \leq 1$.
Let $H$ be a Hilbert space and $i_1 \co H_1 \subseteq H$ an embedding of a sub Hilbert space $H_1$. Moreover, let $I$ be an index set and $I_1 \subseteq I$ a subset. Consider the mapping
$$
J_1 \co 
\begin{cases} \Gamma_q(H_1) \otvn \B(\ell^2_{I_1}) & \to \Gamma_q(H) \otvn \B(\ell^2_I) \\
a \ot x & \mapsto \Gamma_q(i_1)(a) \ot j_1(x) 
\end{cases},
$$
where $j_1 \co \B(\ell^2_{I_1}) \to \B(\ell^2_I),\: x \mapsto P_1^* x P_1$ and $P_1 \co \ell^2_I \to \ell^2_{I_1}, (\xi_i)_{i \in I} \mapsto (\xi_i \delta_{i \in I_1})_{i \in I}$ is the canonical orthogonal projection. Then $J_1$ is a normal faithful trace preserving $*$-homomorphism and thus extends to a complete isometry on the $\L^p$ level, $1 \leq p \leq \infty$.
\end{lemma}

\begin{proof}
According to \cite[Theorem 2.11]{BKS}, since $i_1$ is an isometric embedding, $\Gamma_q(i_1)$ is a faithful $*$-homomorphism which preserves the traces. According to \cite[p.~97]{JMX}, $\Gamma_q(i_1)$ is normal. Moreover, it is easy to check that $j_1$ is also a normal faithful $^*$-homomorphism. Thus also $J_1$ is a normal faithful $^*$-homomorphism, see also \cite[p.~32]{Pis7}. Since $\Gamma_q(i_1)$ and $j_1$ preserve the traces, also $J_1$ preserves the trace. Now $J_1$ is a (complete) $\L^p$ isometry according to \cite[p.~92]{JMX}.
\end{proof}

In the following, we let $(H_n)_{n \in \N}$ be a sequence of mutually orthogonal sub Hilbert spaces of some big Hilbert space $H = \bigoplus_{n \in \N} H_n$ and let $I = \bigsqcup_{n \in \N} I_n$ be a partition of a big index set $I$ into smaller pieces $I_n$.

\begin{lemma}
\label{lem-ellp-sum-isometry}
Let $1 < p < \infty$ and $-1 \leq q \leq 1$.
Consider the mappings $\Psi_n \co S^p_{I_n} \to S^p_I, \: x \mapsto P_n^* x P_n$, where $P_n \co \ell^2_I \to \ell^2_{I_n}$ is the canonical orthogonal projection as above,
\[ 
\Psi \co \bigoplus_{n \in \N}^p S^p_{I_n} \to S^p_I,
\: (x_n) \mapsto \sum_n \Psi_n(x_n) 
\]
and
\[ 
J \co \bigoplus_{n \in \N}^p \L^p(\Gamma_q(H_n) \otvn \B(\ell^2_{I_n})) \to \L^p(\Gamma_q(H) \otvn \B(\ell^2_I)), 
\: (y_n) \mapsto \sum_n J_n(y_n).
\]
\begin{enumerate}
\item On the $\L^\infty$ level, $\Ran(J_n) \cdot \Ran(J_m) = \{ 0 \}$ for $n \neq m$.
\item If the above domains of $\Psi$ and $J$ that are exterior Banach space sums, are as indicated equipped with the $\ell^p$-norms, then $\Psi$ and $J$ are isometries.
\end{enumerate}
\end{lemma}

\begin{proof}
1. By a density and normality argument, it suffices to pick $a \ot x \in \Gamma_q(H_n) \ot \B(\ell^2_{I_n})$ and $b \ot y \in \Gamma_q(H_m) \ot \B(\ell^2_{I_m})$ and calculate the product $J_n(a \ot x) J_m(b \ot y)$.
We have
\[ 
J_n(a \ot x) J_m(b \ot y) = \Gamma_q(i_n)(a) \Gamma_q(i_m)(b) \ot P_n^* x P_n P_m^* y P_m = 0 , 
\] 
since $P_n P_m^* \co \ell^2_{I_m} \to \ell^2_{I_n}$ equals $0$.

2. According to the first point, we have for any $x_n \in \Gamma_q(H_n) \otvn \B(\ell^2_{I_n})$,
\begin{align*}
\left( \sum_{n = 1}^N J_n(x_n) \right)^* \left( \sum_{m = 1}^N J_m(x_m) \right) 
& = \left(\sum_{n = 1}^N J_n(x_n^*)\right) \left( \sum_{m = 1}^N J_m(x_m) \right) 
= \sum_{n = 1}^N J_n(x_n^*x_n) .
\end{align*}
In the same way, by functional calculus of $\left|\sum_{n =1}^N J_n(x_n) \right|^2$ and $J_n(|x_n|^2)$, we obtain
\[ 
\left| \sum_{n =1}^N J_n(x_n) \right|^p = \sum_{n = 1}^N |J_n(x_n)|^p, 
\]
and thus, taking traces, we obtain
\[ 
\norm{ \sum_{n = 1}^N J_n(x_n) }_p = \left( \sum_{n = 1}^N \norm{J_n(x_n)}_p^p \right)^{\frac1p} = \left( \sum_{n =1}^N \norm{x_n}_p^p \right)^{\frac1p} . 
\]
By a density argument, we infer that $J$ is an isometry.
The proof for for $\Psi$ is easier and left to the reader.
\end{proof}

In the following, we let moreover $\alpha_n \co I_n \to H_n$ be mappings and associate with it $\alpha \co I \to H$ given by $\alpha(i) = \alpha_n(i)$ if $i \in I_n \subseteq \bigsqcup_{k \in \N}I_k$.
We thus have noncommutative gradients $\partial_{\alpha_n,q,p}$ and $\partial_{\alpha,q,p}$.

\begin{lemma}
\label{lem-exchange-gradient}
Let $1 < p < \infty$ and $-1 \leq q \leq 1$.
Recall the mappings $J_n$ and $\Psi_n$ from Lemma \ref{lem-ellp-sum-isometry}.
\begin{enumerate}
\item For any $n \in \N$ and $x_n \in \M_{I_n,\fin}$, we have $\Psi_n(x_n) \in \M_{I,\fin}$ and
\[ 
\partial_{\alpha,q,p} \Psi_n(x_n) 
= J_n \partial_{\alpha_n,q,p}(x_n) . 
\]
\item For any $n \in \N$, $y_n \in \Gamma_q(H_n) \ot \M_{I_n,\fin}$
, we have
$J_n(y_n) \in \dom(\partial_{\alpha,q,p^*})^*$
and
\[ 
(\partial_{\alpha,q,p^*})^* J_n(y_n) 
= \Psi_n (\partial_{\alpha_n,q,p^*})^*(y_n) . 
\]
\end{enumerate}
\end{lemma}

\begin{proof}
1. We calculate with $x_n = \sum_{i,j \in I_n} x_{ij} e_{ij}$ and explicit embedding $k_n \co I_n \hookrightarrow I$,
\begin{align*}
\MoveEqLeft
\partial_{\alpha,q,p} \Psi_n(x_n) 
=\partial_{\alpha,q,p} (P_n^* x_n P_n) 
=\partial_{\alpha,q,p} \bigg(P_n^* \sum_{i,j} x_{ij} e_{ij} P_n \bigg) \\
&=\sum_{i,j} x_{ij} \partial_{\alpha,q,p}(e_{k_n(i)k_n(j)}) 
=\sum_{i,j} x_{ij} s_q(\alpha(k_n(i)) - \alpha( k_n(j))) \ot e_{k_n(i)k_n(j)}.
\end{align*}
On the other hand, we have
\begin{align*}
\MoveEqLeft
J_n \partial_{\alpha_n,q,p}(x_n) 
= \sum_{i,j} x_{ij} J_n(s_q(\alpha_n(i) - \alpha_n(j)) \ot e_{ij}) \\
&=\sum_{i,j} x_{ij}\Gamma_q(i_n)(s_q(\alpha_n(i) - \alpha_n(j))) \ot P_n^* e_{ij} P_n 
=\sum_{i,j} x_{ij} s_q(\alpha(k_n(i)) - \alpha(k_n(j))) \ot e_{k_n(i)k_n(j)}.
\end{align*}

2. First we note that $J_n(\Gamma_q(H_n) \ot \M_{I_n,\fin}) \subseteq \dom(\partial_{\alpha,q,p^*})^*$ and that $(\partial_{\alpha,q,p^*})^*$ maps $J_n(\Gamma_q(H_n) \ot \M_{I_n,\fin})$ into $\Psi_n(S^p_{I_n})$. Indeed for $s_q(h_n) \ot e_{ij} \in \Gamma_q(H_n) \ot \M_{I_n,\fin}$ and $x = \sum_{k,l} x_{kl} e_{kl} \in S^{p^*}_I$, we have
\begin{align*}
x &\mapsto  \tau_{\Gamma_q(H) \otvn \B(\ell^2_I)} \left( (s_q(i_n(h_n)) \ot e_{k_n(i)k_n(j)})^* \cdot \partial_{\alpha,q,p^*}(x) \right) \\
& = \sum_{k,l} x_{kl} \tau_{\Gamma_q(H) \otvn \B(\ell^2_I)} \big(s_q(i_n(h_n)) \ot e_{k_n(j)k_n(i)} \cdot s_q(\alpha(k) - \alpha(l)) \ot e_{kl}\big) \\
& = \tau_{\Gamma_q(H)}\big(s_q(i_n(h_n)) s_q(\alpha(k_n(j))-\alpha(k_n(i)))\big) x_{k_n(i)k_n(j)}.
\end{align*}
This defines clearly a linear form on $S^{p^*}_I$, so indeed $J_n(\Gamma_q(H_n) \ot \M_{I_n,\fin}) \subseteq \dom(\partial_{\alpha,q,p^*})^*$.
Moreover, we see above that if $x = e_{kl}$ and $k,l$ not both in $k_n(I_n)$, then $\langle (\partial_{\alpha,q,p^*})^*(J_n(y_n)), e_{kl} \rangle = 0$.
Therefore, $(\partial_{\alpha,q,p^*})^*$ maps $J_n(\Gamma_q(H_n) \ot \M_{I_n,\fin})$ into $\Psi_n(S^p_{I_n})$. Now we check the claimed equality in the statement of the lemma by applying a dual element $x_n$ to both sides. By density of $\M_{I,\fin}$ in $S^{p^*}_I$, we can assume $x_n \in \M_{I,\fin}$. Moreover, by the above, we can assume $x_n \in \Psi_n(\M_{I_n,\fin})$, so that $x_n = \Psi_n \Psi_n^*(x_n)$. Then we have according to the first point of the lemma
\begin{align*}
\MoveEqLeft
\tr_I(\Psi_n (\partial_{\alpha_n,q,p^*})^*( y_n) x_n^*) \\
&=\tr_{I_n}( (\partial_{\alpha_n,q,p^*})^*(y_n)(\Psi_n^*(x_n))^*) 
=\tau_{\Gamma_q(H_n) \otvn \B(\ell^2_{I_n})} (y_n (\partial_{\alpha_n,q,p^*} \Psi_n^*(x_n))^*) \\
&=\tau_{\Gamma_q(H) \otvn \B(\ell^2_I)}(J_n(y_n (\partial_{\alpha_n,q,p^*} \Psi_n^*(x_n))^*)) 
=\tau_{\Gamma_q(H) \otvn \B(\ell^2_I)}(J_n(y_n) (J_n \partial_{\alpha_n,q,p^*} \Psi_n^*(x_n))^*) \\
&=\tau_{\Gamma_q(H) \otvn \B(\ell^2_I)}(J_n(y_n) (\partial_{\alpha,q,p^*} \Psi_n \Psi_n^*(x_n))^*) \\
&=\tr_{I}((\partial_{\alpha,q,p^*})^*J_n(y_n) (\Psi_n \Psi_n^*(x_n))^*) 
=\tr_{I}((\partial_{\alpha,q,p^*})^*J_n(y_n) x_n^*).
\end{align*}
Here we have also used that $J_n$ is trace preserving and multiplicative. The lemma is proved.
\end{proof}

\begin{lemma}
\label{lem-Hodge-Dirac-transference}
Let $1 < p < \infty$ and $-1 \leq q \leq 1$.
Recall the mappings $\Psi$ and $J$ from Lemma \ref{lem-ellp-sum-isometry}.
\begin{enumerate}
\item
Then for any $x = (x_n)$ with $x_n \in \M_{I_n,\fin} \subseteq S^p_{I_n}$ and any $y = (y_n)$ with $y_n \in \Gamma_q(H_n) \ot \M_{I_n,\fin} \cap \ovl{\Ran(\partial_{\alpha_n,q,p})}$ such that only finitely many $x_n$ and $y_n$ are non-zero, we have
\begin{equation}
\label{equ-Hodge-Dirac-transference}
D_{\alpha,q,p} \circ \begin{bmatrix} \Psi & 0 \\ 0 & J \end{bmatrix} (x,y) = \begin{bmatrix} \Psi & 0 \\ 0 & J \end{bmatrix} \circ \begin{bmatrix} D_{\alpha_1,q,p} & 0 & \ldots & \\ 0 & D_{\alpha_2,q,p} & 0 & \ldots \\ \vdots & 0 & \ddots & \end{bmatrix} (x,y)
\end{equation}
\item Let $\omega \in (0,\frac{\pi}{2})$ such that $D_{\alpha,q,p}$ and all $D_{\alpha_n,q,p}$ have an $\HI(\Sigma_\omega^\pm)$ calculus, e.g. according to Theorem \ref{Th-functional-calculus-bisector-Schur}, $\omega = \frac{\pi}{4}$.
Then for any $m \in \HI(\Sigma_\omega^\pm)$, we have
\begin{equation}
\label{equ-2-Hodge-Dirac-transference}
m(D_{\alpha,q,p}) \circ \begin{bmatrix} \Psi & 0 \\ 0 & J \end{bmatrix} = \begin{bmatrix} \Psi & 0 \\ 0 & J \end{bmatrix} \circ \begin{bmatrix} m(D_{\alpha_1,q,p}) & 0 & \ldots \\ 0 & m(D_{\alpha_2,q,p}) & 0  \\ \vdots & 0 & \ddots \end{bmatrix}
\end{equation}
as bounded operators $\bigoplus_{n \in \N}^p \left( S^p_{I_n} \oplus_p \ovl{\Ran(\partial_{\alpha_n,q,p})} \right) \to S^p_I \oplus_p \ovl{\Ran(\partial_{\alpha,q,p})}$.
\end{enumerate}
\end{lemma}

\begin{proof}
1. We have according to Lemma \ref{lem-exchange-gradient}
\begin{align*}
\MoveEqLeft
D_{\alpha,q,p} (\Psi(x),J(y)) 
=D_{\alpha,q,p}\bigg(\sum_n \Psi_n(x_n), \sum_n J_n(y_n)\bigg)  \\ 
&=\bigg((\partial_{\alpha,q,p^*})^*\sum_n J_n(y_n), \partial_{\alpha,q,p} \sum_n \Psi_n(x_n)\bigg) 
=\bigg(\sum_n \Psi_n (\partial_{\alpha_n,q,p^*})^* y_n, \sum_n J_n \partial_{\alpha_n,q,p} x_n\bigg) \\
&=\big(\Psi \text{diag}((\partial_{\alpha_n,q,p^*})^*: \: n \in \N)y, J\text{diag}(\partial_{\alpha_n,q,p}: \: n \in \N)x \big).
\end{align*}
This shows \eqref{equ-Hodge-Dirac-transference}.

2. According to \eqref{equ-Hodge-Dirac-transference}, for any $\lambda \in \rho(D_{\alpha,q,p}) \cap \bigcap_{n \in \N} \rho(D_{\alpha_n,q,p})$, we have for $(x,y)$ as in the first part of the lemma,
\begin{align*}
\MoveEqLeft (\lambda - D_{\alpha,q,p})^{-1} 
\begin{bmatrix} 
\Psi & 0 \\ 0 & J 
\end{bmatrix} 
\begin{bmatrix} 
(\lambda - D_{\alpha_1,q,p}) & 0 \\ 
0 & \ddots 
\end{bmatrix} (x,y)
= (\lambda - D_{\alpha,q,p})^{-1}(\lambda - D_{\alpha,q,p}) 
\begin{bmatrix} 
\Psi & 0 \\ 
0 & J 
\end{bmatrix}(x,y) \\
& = \begin{bmatrix} 
\Psi & 0 \\ 
0 & J 
\end{bmatrix}(x,y) 
=\begin{bmatrix} 
\Psi & 0 \\ 
0 & J 
\end{bmatrix} 
\begin{bmatrix} (\lambda - D_{\alpha_1,q,p})^{-1} & 0 \\ 
0 & \ddots \end{bmatrix} 
\begin{bmatrix} (\lambda - D_{\alpha_1,q,p}) & 0 \\ 
0 & \ddots 
\end{bmatrix}(x,y).
\end{align*}
Thus, we obtain
\begin{equation}
\label{equ-1-proof-Hodge-Dirac-transference}
(\lambda-D_{\alpha,q,p})^{-1} 
\begin{bmatrix} 
\Psi & 0 \\ 
0 & J \end{bmatrix} 
= 
\begin{bmatrix} 
\Psi & 0 \\ 0 & J 
\end{bmatrix} 
\begin{bmatrix} (\lambda - D_{\alpha_1,q,p})^{-1} & 0 \\ 0 & \ddots 
\end{bmatrix}
\end{equation}
on a dense subspace of $\bigoplus_{n \in \N}^p \left( S^p_{I_n} \oplus_p \ovl{\Ran(\partial_{\alpha_n,q,p})} \right)$.
By boundedness of both sides of \eqref{equ-1-proof-Hodge-Dirac-transference}, we obtain equality in \eqref{equ-1-proof-Hodge-Dirac-transference} as bounded operators $\bigoplus_{n \in \N}^p \left( S^p_{I_n} \oplus_p \ovl{\Ran(\partial_{\alpha_n,q,p})} \right) \to S^p_I \oplus_p \ovl{\Ran(\partial_{\alpha,q,p})}$.
Now by the Cauchy integral formula, we obtain \eqref{equ-2-Hodge-Dirac-transference} for $m \in \HI_0(\Sigma_\omega^\pm)$, and by the $\HI$ convergence lemma (see \cite[Lemma 2.1]{CDMY} for the sectorial case) applied first for fixed $(x,y)$ as in the first part of the lemma, also for $m \in \HI(\Sigma_\omega^\pm)$.
\end{proof}

Now, we can give an answer to a variant of \cite[Problem C.5]{JMP2}.

\begin{thm}
\label{prop-explicitly-dimension-free}
Let $1 < p < \infty$, $-1 \leq q \leq 1$, $\omega = \frac{\pi}{4}$ 
 and $D_{\alpha,q,p}$ as in Theorem \ref{Th-functional-calculus-bisector-Schur}.
Then the $\HI(\Sigma_\omega^\pm)$ calculus norm of $D_{\alpha,q,p}$ is controlled independently of $H$ and $\alpha$, that is, $\norm{m(D_{\alpha,q,p})} \leq C_{q,p} \norm{m}_{\HI(\Sigma_\omega^\pm)}$.
In particular, $\sgn(D_{\alpha,q,p}) = D_{\alpha,q,p} (D_{\alpha,q,p}^2)^{-\frac12}$ is bounded by a constant not depending on $H$ or its dimension nor the mapping $\alpha \co I \to H$.
\end{thm}

\begin{proof}
Suppose that the statement of the theorem is false. For any $n \in \N$, then there exists a sequence $(D_{\alpha_n,q,p})$ of Hodge-Dirac operators such that the $\HI(\Sigma_{\omega}^\pm)$ functional calculus norm is bigger than $n$. We let $I \ov{\mathrm{def}}{=} \bigsqcup_{n \in \N} I_n$, $H \ov{\mathrm{def}}{=} \bigoplus_{n \in \N} H_n$ and consider $\alpha \co I \to H$, $i \in I_n \mapsto \alpha_n(i)$. Let $f_n$ be a function of $\HI(\Sigma_\omega^\pm)$ such that $\norm{f_n}_{\HI(\Sigma_{\omega}^\pm)} = 1$ and $\norm{f_n(D_{\alpha_n,q,p})} \geq n$. According to Theorem \ref{Th-functional-calculus-bisector-Schur}, the operator $D_{\alpha,q,p}$ has a bounded $\HI(\Sigma_\omega^\pm)$ functional calculus. Thus, for some constant $C < \infty$ and any $f_n$ as above, according to Lemma \ref{lem-Hodge-Dirac-transference},
\begin{align*}
\MoveEqLeft
C \geq \norm{f_n(D_{\alpha,q,p})} 
\geq \norm{f_n(D_{\alpha,q,p}) 
\circ \begin{bmatrix} 
\Psi & 0 \\ 
0 & J 
\end{bmatrix} } 
=\norm{ \begin{bmatrix} 
\Psi & 0 \\ 
0 & J \end{bmatrix} 
\circ \begin{bmatrix} 
f_n(D_{\alpha_1,q,p}) & 0 \\ 
0 & \ddots \end{bmatrix} } \\
&= \norm{ \begin{bmatrix} 
f_n(D_{\alpha_1,q,p}) & 0 \\ 
0 & \ddots \end{bmatrix} } 
\geq \norm{f_n(D_{\alpha_n,q,p})} 
\geq n.
\end{align*}
Taking $n \to \infty$ yields a contradiction. The last sentence of Theorem \ref{prop-explicitly-dimension-free} follows from taking $f(z) = 1_{\Sigma_\omega^+}(z) - 1_{\Sigma_\omega^-}(z)$.
\end{proof}

\begin{remark} \normalfont
\label{rem-dimension-free-Hodge-Schur}
We have the following alternative proof of Theorem \ref{prop-explicitly-dimension-free} using dimension-free estimates of the Riesz transform associated with $A_p$. It will show that the angle of the $\HI(\Sigma_\omega^\pm)$ calculus of $D_{\alpha,q,p}$ and $\D_{\alpha,q,p}$ can be chosen $\omega > \frac{\pi}{2} | \frac1p - \frac12 |$ and that the norm of the calculus is bounded by a constant $K_\omega$ not depending on $I$ nor the representation $(\alpha,H)$, in particular it is independent of the dimension of $H$.
\end{remark}

\begin{proof}
First note that since $A_p$ has a (sectorial) completely bounded $\HI(\Sigma_{2\omega})$ calculus with angle $2\omega > \pi | \frac1p - \frac12 |$ (see the dilation argument below), by the representation \eqref{Dirac-carre-Dirac}, $\D_{\alpha,q,p}^2$ also has a sectorial $\HI(\Sigma_{2\omega})$ calculus.
According to \cite[Proof of Theorem 10.6.7, Theorem 10.4.4 (1) and (3), Proof of Theorem 10.4.9]{HvNVW2}, the operator $\D_{\alpha,q,p}$ has then an $\HI(\Sigma_\omega^\pm)$ bisectorial calculus to the angle $\omega > \frac{\pi}{2} | \frac1p - \frac12 |$ with a norm control
\begin{equation}
\label{equ-1-rem-dimension-free-Hodge-Schur}
\norm{f(\D_{\alpha,q,p})} \leq K_\omega \left( M_{2\omega,\D_{\alpha,q,p}^2}^\infty \right)^2 \left( M_{\omega,\D_{\alpha,q,p}}^R \right)^2 \norm{f}_{\infty,\omega},
\end{equation}
where $K_\omega$ is a constant only depending on $\omega$ (and not on $I$ nor the representation $(\alpha,H)$).
Here $M_{2\omega,\D_{\alpha,q,p}^2}^\infty$ is the $\HI(\Sigma_{2 \omega})$ calculus norm of $\D_{\alpha,q,p}^2$ and
\begin{equation}
\label{equ-2-rem-dimension-free-Hodge-Schur}
M_{\omega,\D_{\alpha,q,p}}^R  = R\left(\left\{ \lambda (\lambda - \D_{\alpha,q,p})^{-1} :\: \lambda \in \C \backslash \{ 0 \},\: \left|\,|\arg(\lambda)| - \frac{\pi}{2} \right| < \frac{\pi}{2} - \omega \right\}\right).
\end{equation}
Thus it remains to show that both $M_{2\omega,\D_{\alpha,q,p}^2}^\infty$ and $M_{\omega,\D_{\alpha,q,p}}^R$ can be chosen independently of the Hilbert space $H$ and the representation $\alpha$.
Let us start with $M_{2\omega,\D_{\alpha,q,p}^2}^\infty$.
It is controlled according to the above reasoning and Proposition \ref{Dirac-carre-Dirac}, by $M_{2\omega,\Id_{S^p} \ot A_p}^\infty$, that is, the completely bounded $\HI$ calculus norm of $A_p$.
An application of \cite[Proposition 5.8]{JMX} shows that it suffices to consider only $2\omega > \frac{\pi}{2}$.
According to \cite{Arh1} \cite{Arh4}, we have a decomposition of the semigroup $(T_{t,p})_{t \geq 0}$ generated by $A_p$, given by
\[ 
\Id_{S^p} \ot T_{t,p} 
= (\Id_{S^p} \ot \E_p)(\Id_{S^p} \ot U_{t,p})(\Id_{S^p} \ot J_p) . 
\]
Here, $\E_p \co \L^p(\Omega,S^p_I) \to S^p_I$ and $J_p$ are complete contractions and $U_{t,p} \co \L^p(\Omega,S^p_I) \to \L^p(\Omega,S^p_I)$ is a group of isometries. An inspection of the proof of Lemma \ref{lem-technical-HI-Laplace-transform} shows that the constant in the second condition there is a bound of the $\HI$ calculus in the first condition there, so that it  suffices to show that the generator of $\Id_{S^p} \ot U_{t,p}$ has a bounded $\HI(\Sigma_{2\omega})$ calculus with a norm controlled by a constant depending only on $2 \omega > \frac{\pi}{2}$. By \cite[Proof of Theorem 10.7.10]{HvNVW2}, the norm of the calculus of $U_{t,p}$ is controlled by $c_{2\omega} \beta^2_{p,X} h_{p,X}$, $c_{2\omega}$ denoting a constant depending only on $2\omega$, $\beta_{p,X}$ denoting the UMD constant of $X = S^p(\L^p(M))$ and $h_{p,X}$ denoting the Hilbert transform norm, i.e. on the space $\L^p(\R,S^p(\L^p(M)))$. Here, $M = \L^\infty(\Omega) \otvn \B(\ell^2_I)$.

These are controlled by a constant depending only on $p$ but not on $M$. Indeed, using \cite[Proposition 4.2.15]{HvNVW2}, it suffices to control the UMD constant of $S^p(\L^p(M))$ and thus uniformly the one of $S^p_n(\L^p(M))$. Now $\M_n(M)$ is a $\QWEP$ finite von Neumann algebra, so $\M_n(M)$ admits a trace preserving embedding in an ultrapower $\B(\ell^2_J)^\ul$ of $\B(\ell^2_J)$ for some index set $J$ and some ultrafilter $\ul$. So $S^p_n(\L^p(M))$ embeds in $(S^p_J)^\ul$. By \cite[Theorem 8.13]{DJT}, $(S^p_J)^\ul$ is finitely representable in $S^p_J$. We conclude since UMD is a super-property by \cite[p.~308]{HvNVW2}.
So we have the desired control of $M_{2\omega,\Id_{S^p} \ot A_p}^\infty$, and thus of $M_{2\omega,\D_{\alpha,q,p}^2}^\infty$.
We turn to the control of $M_{\omega,\D_{\alpha,q,p}}^R$ from \eqref{equ-2-rem-dimension-free-Hodge-Schur}.
According to \eqref{Resolvent} extended to complex times $z$ belonging to some bisector $\Sigma_{\sigma}^\pm$ with $\sigma = \frac{\pi}{2} - \omega$,
it suffices to control the following $R$-bounds
\begin{align}
& R \left( \left\{ (\Id + z^2 A_p)^{-1} :\: z \in \Sigma_\sigma^\pm \right\} \right) \label{equ-3-rem-dimension-free-Hodge-Schur}  \\
& R \left( \left\{ z (\Id + z^2 A_p)^{-1} (\partial_{\alpha,q,p^*})^* : \: z \in \Sigma_\sigma^\pm \right\} \right) \label{equ-4-rem-dimension-free-Hodge-Schur} \\
& R \left( \left\{ z \partial_{\alpha,q,p} (\Id + z^2 A_p)^{-1} : \: z \in \Sigma_\sigma^\pm \right\} \right) \label{equ-5-rem-dimension-free-Hodge-Schur} \\
& R \left( \left\{ \Id \ot (\Id + z^2 A_p)^{-1} : \: z \in \Sigma_\sigma^\pm \right\} \right) \label{equ-6-rem-dimension-free-Hodge-Schur}.
\end{align}
Indeed, an operator matrix family is $R$-bounded if and only if all operator entries in the matrix are $R$-bounded.
According to \cite[Theorem 10.3.4 (1)]{HvNVW2}, \eqref{equ-3-rem-dimension-free-Hodge-Schur} is $R$-bounded since we can write $(\Id + z^2 A_p)^{-1} = \Id - f(z^2A_p)$ with $f(\lambda) = \lambda(1+\lambda)^{-1}$.
Moreover, by the same reference, its $R$-bound is controlled by $M_{2\omega - \epsi,A_p}^\infty$, which in turn, by the above argument of dilation can be controlled independently of the Hilbert space $H$ and the representation $\alpha$.
The same argument shows that also \eqref{equ-6-rem-dimension-free-Hodge-Schur} is $R$-bounded.
Since the operator family in \eqref{equ-4-rem-dimension-free-Hodge-Schur} consists of the family of the adjoints in \eqref{equ-5-rem-dimension-free-Hodge-Schur} and $R$-boundedness is preserved under adjoints, it suffices to prove that \eqref{equ-5-rem-dimension-free-Hodge-Schur} is $R$-bounded.
To this end, we decompose for $z \in \Sigma_\omega$ the positive part of the bisector (similarly if $z$ belongs to the negative part of the bisector)
\[ z \partial_{\alpha,q,p} (\Id + z^2 A_p)^{-1} = \left[\partial_{\alpha,q,p} A_p^{-\frac12}\right] \left(z^2 A_p \right)^{\frac12}(\Id + z^2 A_p)^{-1} = \left[\partial_{\alpha,q,p} A_p^{-\frac12}\right] f(z^2A_p) \]
with $f(\lambda) = \sqrt{\lambda}(1+\lambda)^{-1}$.
Again \cite[Theorem 10.3.4 (1)]{HvNVW2} shows that the term $f(z^2A_p)$ is $R$-bounded with $R$-bound controlled by some constant independent of the representation $(\alpha,H)$.
Finally, we are left to show that the Riesz transform is bounded by a constant independent of the cocycle, that is,
\begin{equation}
\label{equ-7-rem-dimension-free-Hodge-Schur}
\norm{\partial_{\alpha,q,p} (x) }_{\L^p(\Gamma_q(H) \otvn \B(\ell^2_I))} 
\leq C \bnorm{A_p^{\frac12}(x) }_{S^p_I} .
\end{equation}
For this in turn we refer to Remark \ref{rem-constants-Riesz-equivalence-q}.
\end{proof}

\section{Locally compact quantum metric spaces and spectral triples}
\label{new-quantum-metric-spaces}
\subsection{Background on quantum compact metric spaces}
\label{subsec-prelims-quantum}

We recall definitions and characterizations of the notions that we need in this section. The main notion is that of quantum compact metric space. This concept has its origins in Connes' paper \cite{Con4} of 1989 (see also \cite[Chapter 6]{Con3} and \cite[Chapter 3]{Var1}), in which he shows that we can recover the geodesic distance $d$ of a compact riemannian spin manifold $M$ using the Dirac operator $D$ by the formula
$$
d(p,q)
=\sup_{f \in \mathrm{C}(M), \norm{[D,f]} \leq 1} |f(p)-f(q)|, \quad p,q \in M
$$ 
where the commutator $[D,f]=Df-fD$ extends to a bounded operator\footnote{\thefootnote. Recall that $D$ is an unbounded operator acting on the Hilbert space of $\L^2$-spinors and that the functions of $\mathrm{C}(M)$ act on the same Hilbert space by multiplication operators.}. Indeed, it is known that the commutator $[D,f]$ induces a bounded operator if and only if $f$ is a Lipschitz function and in this case the Lipschitz norm of $f$ is equal to $\norm{[D,f]}$. Moreover, this space of functions is norm dense in $\mathrm{C}(M)$. If we identify the points $p,q$ as pure states $\omega_p$ and $\omega_q$ on the algebra $\mathrm{C}(M)$, this formula can be seen as 
$$
d(\omega_p,\omega_q)
=\sup_{f \in \mathrm{C}(M), \norm{[D,f]} \leq 1} |\omega_p(f)-\omega_q(f)|, \quad p,q \in M.
$$

Afterwards, Rieffel \cite{Rie3} and Latr\'emoli\`ere \cite{Lat1} axiomatized this formula replacing $\mathrm{C}(M)$ by a unital $\mathrm{C}^*$-algebra $A$ (or even by an order-unit space $\mathcal{A}$), $f \mapsto \norm{[D,f]}$ by a seminorm $\norm{\cdot}$ defined on a subspace of $A$ and $\omega_p,\omega_q$ by arbitrary states obtaining essentially the formula \eqref{Distance-Monge} below and giving rise to a theory of \textit{quantum} compact metric spaces. With this notion, Rieffel was able to define a quantum analogue of Gromov-Hausdorff distance and to give a meaning to many approximation found in the physics literature. We refer to the surveys \cite{Lat2} and \cite{Rie3} and references therein for more information.

Recall that an order-unit space \cite[p.~69]{Alf1}, \cite[Definition 1.8]{AlS1} is an ordered $\R$-vector space $\A$ with a closed positive cone and an element $1_\A$, satisfying $\norm{a}_\A=\inf\{\lambda>0 : -\lambda1_\A  \leq a \leq \lambda 1_\A\}$. The element $1_\A$ is called the distinguished order unit. The definition of an order-unit space is due to Kadison \cite{Kad}. Important examples of order-unit spaces are given by real linear subspaces of selfadjoint elements containing the unit element in a unital $\mathrm{C}^*$-algebra.

The following is a slight generalization of \cite[Definition 2.3]{Lat2} and \cite[Definition 2.2]{Lat6}. 
\begin{defi}
\label{Def-Lipschitz-pair}
A unital Lipschitz pair $(\A,\norm{\cdot})$ is a pair where $\A$ is an order-unit-space and where $\norm{\cdot}$ is a seminorm defined on a dense subspace $\dom \norm{\cdot}$ of $\A$ such that
\begin{equation}
\label{equa=1A}
\big\{a \in \dom \norm{\cdot} : \norm{a} = 0\big\} 
=\R1_{\A}.
\end{equation}
\end{defi}

\begin{remark} \normalfont
\label{Extension-quantum-metric-space}
Note that if a seminorm $\norm{\cdot}$ is defined on some subspace $\dom \norm{\cdot}$ of a unital $\mathrm{C}^*$-algebra $A$ such that $A_\sa \cap \dom \norm{\cdot}$ is dense in $A_\sa$ and such that $\big\{a \in \dom \norm{\cdot} : \norm{a} = 0\big\} =\C1_{A}$, then its restriction on $A_\sa \cap \dom \norm{\cdot}$ defines a unital Lipschitz pair. In this case, we also say that $(A,\norm{\cdot})$ is a unital Lipschitz pair (or a compact quantum metric space if Definition \ref{Def-quantum-compact-metric-space} is satisfied).
\end{remark}

If $(X,\dist)$ is a compact metric space, a fundamental example \cite[Example 2.6]{Lat2}, \cite[Example 2.9]{Lat7} is given by $(\mathrm{C}(X)_\sa,\Lip)$ where $\mathrm{C}(X)$ is the commutative $\mathrm{C}^*$-algebra of continuous functions and where $\Lip$ is the Lipschitz seminorm, defined for any Lipschitz function $f \co X \to \C$ by 
\begin{equation}
\label{Lip-eq}
\Lip(f) 
\ov{\mathrm{def}}{=} \sup\left\{\frac{|f(x)-f(y)|}{\dist(x,y)} : x,y\in X,x\not=y \right\}.
\end{equation}
It is immediate that a function $f$ has zero Lipschitz constant if and only if it is constant on $X$. Moreover, the set of real Lipschitz functions is norm-dense in $\mathrm{C}(X)_\sa$ by the Stone-Weierstrass theorem.

Now, following \cite[Definition 2.2]{Rie5} (see also \cite[Theorem 2.42]{Lat2}, \cite[Definition 1.2]{Lat3} or \cite[Definition 2.6]{Lat6}), we introduce a notion of quantum compact metric space. Recall that a linear functional $\varphi$ on an order-unit space $\A$ is a state \cite[p.~72]{Alf1} if $\norm{\varphi}=\varphi(1_\A)=1$.

\begin{defi}
\label{Def-quantum-compact-metric-space}
A quantum compact metric space $(\A,\norm{\cdot})$ is a unital Lipschitz pair whose associated Monge-Kantorovich metric
\begin{equation}
\label{Distance-Monge}
\dist_{\mk}(\varphi,\psi) 
\ov{\mathrm{def}}{=}\sup \big\{|\varphi(a)-\psi(a)| : a \in \A, \norm{a} \leq 1 \big\}, \quad \varphi, \psi \in \S(\A)
\end{equation}
metrizes the weak* topology restricted to the state space $\S(\A)$ of $\A$. When a Lipschitz pair $(\A,\norm{\cdot})$ is a quantum compact metric space, the seminorm $\norm{\cdot}$ is referred to as a Lip-norm.
\end{defi}

In the case of $(\mathrm{C}(X)_\sa,\Lip)$, we obtain the dual formulation of the classical Kantorovich-Rubinstein metric \cite[Remark 6.5]{Vil1} for Borel probability measures $\mu$ and $\nu$ on $X$
\begin{equation}
\label{Distance-Monge-commutatif}
\dist_{\mk}(\mu,\nu) 
\ov{\mathrm{def}}{=} \sup \left\{\left|\int_X f \d\mu- \int_X f\d \nu \right| : f \in \mathrm{C}(X)_\sa, \Lip(f) \leq 1 \right\}.
\end{equation}
which is is a basic concept in optimal transport theory \cite{Vil1}. Considering the Dirac measures $\delta_x$ and $\delta_y$ at points $x,y \in X$ instead of $\mu$ and $\nu$, we recover the distance $\dist(x,y)$ with the formula \eqref{Distance-Monge-commutatif}.


The compatibility of Monge-Kantorovich metric with the weak* topology is hard to check directly in general. Fortunately, there exists a condition \cite[Proposition 1.3]{OzR}, \cite[Theorem 2.5]{AgL}, \cite[Theorem 2.43]{Lat2} and \cite[Theorem 2.10]{Lat6} which is more practical. This condition is inspired by the fact that Arz\'ela-Ascoli's theorem shows that for any $x \in X$ the set
$$
\big\{f \in \mathrm{C}(X)_\sa :  \Lip(f) \leq 1, f(x)= 0 \big\}
$$
is norm relatively compact and it is known that this property implies that \eqref{Distance-Monge-commutatif} metrizes the weak* topology on the space of Borel probability measures on $X$. The following is a slight generalization of \cite[Proposition 1.3]{OzR} with a similar proof left to the reader.

\begin{prop}
\label{Prop-norm-precompact}
Let $(\A,\norm{\cdot})$ be a unital Lipschitz pair. Let $\mu \in \S(\A)$ be a state of $\A$. If the set $\{a \in \dom \norm{\cdot}: \norm{a} \leq 1, \mu(a) = 0\}$ is norm relatively compact in $\A$ then $(\A,\norm{\cdot})$ is a quantum compact metric space.
\end{prop}

The Lipschitz seminorm $\mathrm{Lip}$ associated to a compact metric space $(X,\dist)$ enjoys a natural property with respect to the multiplication of functions in $\mathrm{C}(X)$, called the Leibniz property for any Lipschitz functions $f,g \co X \to \C$:
\begin{equation}
\label{equa-Leibniz}
\Lip(fg)
\leq \norm{f}_{\mathrm{C}(X)}\Lip(g) + \Lip(f)\norm{g}_{\mathrm{C}(X)}.
\end{equation}

Moreover, the Lipschitz seminorm is lower-semicontinuous with respect to the $\mathrm{C}^*$-norm of $\mathrm{C}(X)$, i.e. the uniform convergence norm on $X$. These two additional properties were not assumed in the above Definition \ref{Def-quantum-compact-metric-space}, yet they are quite natural.  However, as research in noncommutative metric geometry progressed, the need for a noncommutative analogue of these properties for some developments became evident. So, sometimes, some additional conditions are often added to Definition \ref{Def-quantum-compact-metric-space} which brings us to the following definition.

The following is a slight generalization of \cite[Definition 2.21]{Lat2} for order-unit spaces embedding in unital $\mathrm{C}^*$-algebras (see also \cite[Definition 1.3]{Lat5}) and \cite[Definition 2.45]{Lat2}, \cite[Definition 2.19]{Lat7}, \cite[Definition 2.19]{Lat6} and \cite[Definition 2.2.2]{Lat4}. Here $a \circ b \ov{\mathrm{def}}{=} \frac{1}{2}(ab+ba)$ and $\{a,b\} \ov{\mathrm{def}}{=} \frac{1}{2 \i} (ab-ba)$.

\begin{defi}
\label{Def-Leibniz-pair} \label{Def-Leibniz-quantum-compact-metric-space}
\begin{enumerate}
	\item A unital Leibniz pair $(\A,\norm{\cdot})$ is a unital Lipschitz pair where $\A$ is a real linear subspace of selfadjoint elements containing the unit element in a unital $\mathrm{C}^*$-algebra $A$ such that:
\begin{enumerate}
	\item the domain of $\norm{\cdot}$ is a Jordan-Lie subalgebra of $A_\sa$,
	\item for any $a, b \in \dom \norm{\cdot}$, we have:
$$
\norm{a \circ b} 
\leq \norm{a}_A \norm{b} + \norm{a}\norm{b}_A
\quad \text{and} \quad
\norm{\{a,b\}}
\leq \norm{a}_A \norm{b} +\norm{a}\norm{b}_A.
$$
\end{enumerate}

\item A unital Leibniz pair $(\A,\norm{\cdot})$ is a Leibniz quantum compact
metric space when $\norm{\cdot}$ is a lower semicontinuous Lip-norm.
\end{enumerate}
\end{defi}

Note that neither $\A$ nor $\A \cap \dom \norm{\cdot}$ are in general a Jordan or Lie or a usual algebra. 
We continue with a useful observation \cite[Theorem 2.17]{Lat7}, \cite[p.~18]{Lat2} in the spirit of Remark \ref{Extension-quantum-metric-space}.

\begin{prop}
\label{Prop-passagea-se}
Let $A$ be a unital $\mathrm{C}^*$-algebra and $\norm{\cdot}$ be a seminorm defined on a dense $\C$-subspace $\dom \norm{\cdot}$ of $A$, such that $\dom \norm{\cdot}$ is closed under the adjoint operation, such that $\{a \in \dom \norm{\cdot}: \norm{a} = 0\} = \C 1_\A$ and, for all $a, b \in \dom \norm{\cdot}$, we have:
$$
\norm{ab} 
\leq  \norm{a}_A \norm{b} + \norm{a} \norm{b}_A.
$$
If $\norm{\cdot}_\sa$ is the restriction of $\norm{\cdot}$ to $A_\sa \cap \dom\norm{\cdot}$, then $(A_\sa,\norm{\cdot}_\sa)$ is a unital Leibniz pair.
\end{prop}


In the case of semigroups Schur multipliers with $I$ infinite, we also need a notion of quantum locally compact metric spaces. The paper \cite{Lat1} gives such a definition in the case of a Lipschitz pair $(\A,\norm{\cdot})$ where $\A$ is a $\mathrm{C}^*$-algebra. However, we need a version for order-unit-spaces not covered by \cite{Lat1}. 

The following is a variant of \cite[Definition 2.3]{Lat1} and \cite[Definition 2.27]{Lat1}. Here $uA$ is the unitization of $A$.

\begin{defi}
A Lipschitz pair $(\A,\norm{\cdot})$ is a closed subspace $\A$ of selfadjoint elements of a non-unital $\mathrm{C}^*$-algebra $A$ and a seminorm $\norm{\cdot}$ defined on a dense subspace of $\A \oplus \R 1_{u A}$ such that $\{ x \in \A \oplus \R 1_{u A} : \norm{x}=0 \}=\R 1_{u A}$.

A Lipschitz triple $(\A,\norm{\cdot},\mathfrak{M})$ is a Lipschitz pair $(\A,\norm{\cdot})$ and an abelian $\mathrm{C}^*$-algebra $\mathfrak{M}$ of $A$ such that $\mathfrak{M}$ contains an approximate unit of $A$.
\end{defi}

The following is \cite[Definition 2.23]{Lat1}.

\begin{defi}
\label{def-local-state}
Let $A$ be a non-unital $\mathrm{C}^*$-algebra and $\mathfrak{M}$ be an abelian $\mathrm{C}^*$-subalgebra of $A$ containing an approximate unit of $A$.
Let $\mu \co A \to \C$ be a state of $A$.
We call $\mu$ a local state (of $(A,\mathfrak{M})$) provided that there exists a projection $e$ in $\mathfrak{M}$ of compact support\footnote{\thefootnote. It is not clear if the support must be in addition open in \cite{Lat1} since the indicator of a subset $A$ is continuous if and only if $A$ is both open and closed.}  (in the Gelfand spectrum) such that $\mu(e) = 1$.
\end{defi}

Inspired by \cite[Theorem 3.10]{Lat1}, \cite[Theorem 2.73]{Lat2}, we introduce the following definition. Of course, we recognize that this definition is a bit artificial. A better choice could be to generalize the results and definitions of \cite{Lat1} to a larger context. Here $\mu$ is extended on the unitization as in \cite[Notation 2.2]{Lat1}.

\begin{defi}
\label{Defi-locally-order-unit}
We say that a Lipschitz triple $(\A,\norm{\cdot},\mathfrak{M})$ is a quantum locally compact metric space if for a local state $\mu$ of $(A,\mathfrak{M})$ and any compactly supported $a,b \in \mathfrak{M}$, the set 
$$
a \big\{x \in \A \oplus \R 1_{uA} : \norm{x} \leq 1, \mu(x)=0 \big\} b
$$
is totally bounded in the norm topology of $uA$.
\end{defi}

Now, we give a variant of Definition \ref{Def-Leibniz-pair}.

\begin{defi}
\label{def-Leibniz-local}
Let $(\A,\norm{\cdot},\mathfrak{M})$ be a quantum locally compact metric space.
We call it a Leibniz quantum locally compact metric space if $(\A \oplus \R 1_{uA},\norm{\cdot})$ is a unital Leibniz pair in the sense of Definition \ref{Def-Leibniz-pair} and if the seminorm $\norm{\cdot}$ is lower semicontinuous with respect to the $\mathrm{C}^*$-norm $\norm{\cdot}_{uA}$.
\end{defi}

In this definition, note in a similar way to Definition \ref{Def-Leibniz-pair} that neither $\A$ nor $\A \cap \dom \norm{\cdot}$ are supposed to be a Jordan or Lie or a usual algebra.
In a similar way to Proposition \ref{Prop-passagea-se}, we obtain that if $(\A,\norm{\cdot})$ is a Lipschitz pair such that $\norm{ab} \leq \norm{a}_{uA} \norm{b} + \norm{a} \norm{b}_{uA}$ for any $a,b \in \A \oplus \R 1_{uA}$, then $(\A \oplus \R 1_{uA},\norm{\cdot})$ is a unital Leibniz pair.

%
%

Finally, following \cite[Condition 4.3]{Lat1} (see also \cite{BMR1} for a related discussion), we introduce the following definition.

\begin{defi}
\label{Defi-bounded-locally}
We say that a locally compact metric space $(\A,\norm{\cdot},\mathfrak{M})$ in the sense of Definition \ref{Defi-locally-order-unit} is bounded if the Lipschitz ball $\big\{ a \in \A: \norm{a} \leq 1\big\}$ is norm bounded.
\end{defi}

\subsection{Quantum compact metric spaces associated to semigroups of Fourier multipliers}
\label{subsec-new-compact-quantum-metric-spaces-Fourier}

In this subsection, we consider a markovian semigroup $(T_t)_{t \geq 0}$ of Fourier multipliers on a von Neumann algebra $\VN(G)$ where $G$ is a discrete group as in Proposition \ref{prop-Schoenberg}. We introduce new compact quantum metric spaces in the spirit of the ones of \cite{JM1}. We also add to the picture the lower semicontinuity and a careful examination of the domains. By Lemma \ref{lem-Gamma-closure-group}, the following definition is correct. It is not clear how to do the same analysis at the level $p=\infty$ considered in \cite[Section 1.2]{JM1} since the domain is unclear in this case.

\begin{defi}
\label{defi-gamma-psi-p}
Suppose $2 \leq p < \infty$. Let $G$ be a discrete group. Let $A_p$ denote the $\L^p$ realization of the (negative) generator of $(T_t)_{t \geq 0}$. For any $x \in \dom A^{\frac12}_p$ we let
\begin{equation}
\label{def-gamma-psi-p}
\norm{x}_{\Gamma,p}
\ov{\mathrm{def}}{=}
\max \Big\{ \bnorm{\Gamma(x,x)^\frac12}_{\L^p(\VN(G))}, \bnorm{ \Gamma(x^*,x^*)^\frac12}_{\L^p(\VN(G))} \Big\}.
\end{equation}
\end{defi}

We start with an elementary fact.

\begin{prop}
\label{Prop-C1-Schur-infinite-sgrp}
Suppose $2 \leq p <\infty$. Let $G$ be a discrete group. Then $\norm{\cdot}_{\Gamma,p}$ is a seminorm on the subspace $\dom A^{\frac12}_p$ of $\L^p(\VN(G))$.
\end{prop}

\begin{proof}
1. If $x$ belongs to $\dom A^{\frac12}_p$ and $k \in \C$, we have
\begin{align*}
\MoveEqLeft
\norm{k x}_{\Gamma,p}
\ov{\eqref{def-gamma-psi-p}}{=} \max \Big\{ \bnorm{\Gamma(k x,k x)^\frac12}_{\L^p(\VN(G))}, \bnorm{ \Gamma((k x)^*,(k x)^*)^\frac12}_{\L^p(\VN(G))} \Big\} \\  
&=\max \Big\{ \bnorm{\Gamma(k x,k x)^\frac12}_{\L^p(\VN(G))}, \bnorm{ \Gamma(\ovl{k} x^*,\ovl{k} x^*)^\frac12}_{\L^p(\VN(G))}\Big\}\\         
&=\max \Big\{ |k| \bnorm{\Gamma(x,x)^\frac12}_{\L^p(\VN(G))}, |k|\bnorm{ \Gamma(x^*,x^*)^\frac12}_{\L^p(\VN(G))} \Big\}\\
&=|k|\max \Big\{ \bnorm{\Gamma(x,x)^\frac12}_{\L^p(\VN(G))}, \bnorm{ \Gamma(x^*,x^*)^\frac12}_{\L^p(\VN(G))} \Big\}
\ov{\eqref{def-gamma-psi-p}}{=}|k|\norm{x}_{\Gamma,p}.
\end{align*} 
Let us turn to the triangular inequality. Assume that $x,y \in \P_G$. Note that $\L^{\frac{p}{2}}(\VN(G))$ is a normed space since $p \geq 2$. According to the part 4 of Lemma \ref{Lemma-gamma-infos-sgrp} and \eqref{Inegalite-triangulaire-1} applied with $Z=\L^{\frac{p}{2}}(\VN(G))$, we have
$$
\norm{\Gamma(x+y,x+y)}_{\L^{\frac{p}{2}}(\VN(G))}^{\frac{1}{2}}           
\ov{\eqref{Inegalite-triangulaire-1}}{\leq} \norm{\Gamma(x,x)}_{\L^{\frac{p}{2}}(\VN(G))}^{\frac{1}{2}}+\norm{\Gamma(y,y)}_{\L^{\frac{p}{2}}(\VN(G))}^{\frac{1}{2}}.
$$
We can rewrite this inequality under the form
\begin{equation}
\label{ine-6543}
\bnorm{\Gamma(x+y,x+y)^{\frac12}}_{\L^p(\VN(G))} 
\leq \bnorm{\Gamma(x,x)^{\frac12}}_{\L^p(\VN(G))} + \bnorm{\Gamma(y,y)^{\frac12}}_{\L^p(\VN(G))}.
\end{equation}
Then, for the general case where $x,y \in \dom A^{\frac12}_p$, consider some sequences $(x_n)$ and $(y_n)$ of $\P_G$ such that $x_n \xra[]{a} x$ and $y_n \xra[]{a} y$ with Lemma \ref{lem-Gamma-closure-group} and \eqref{Domaine-closure}. By Remark \ref{Remark-xra}, we have $x_n+y_n \xra[]{a} x+y$ and thus
\begin{align*}
\MoveEqLeft
\bnorm{\Gamma(x+y,x+y)^{\frac12}} _{\L^p(\VN(G))}
=\bnorm{\Gamma(x+y,x+y)}_{\L^{\frac{p}{2}}(\VN(G))}^{\frac12} \\
&\ov{\eqref{Def-closure}}{=} \lim_n \bnorm{\Gamma(x_n+y_n,x_n+y_n)}_{\L^{\frac{p}{2}}(\VN(G))}^{\frac12} 
=\lim_n \bnorm{\Gamma(x_n+y_n,x_n+y_n)^{\frac12}} _{\L^p(\VN(G))} \\
&\ov{\eqref{ine-6543}}{\leq} \lim_n \bnorm{\Gamma(x_n,x_n)^{\frac12}}_{\L^p(\VN(G))} + \lim_n \bnorm{\Gamma(y_n,y_n)^{\frac12}}_{\L^p(\VN(G))} \\
&\ov{\eqref{Def-closure}}{=}\bnorm{\Gamma(x,x)^{\frac12}}_{\L^p(\VN(G))} + \bnorm{\Gamma(y,y)^{\frac12}}_{\L^p(\VN(G))}.
\end{align*}
A similar reasoning for $x,y$ replaced by $x^*,y^*$ gives 
\begin{equation}
\label{Equa-divers-100000}
\bnorm{\Gamma(x^*+y^*,x^* + y^*)^{\frac12}}_{\L^p(\VN(G))} \leq \bnorm{\Gamma(x^*,x^*)^{\frac12}}_{\L^p(\VN(G))} + \bnorm{\Gamma(y^*,y^*)^{\frac12}}_{\L^p(\VN(G))}.
\end{equation}
Now we have
\begin{align*}
\MoveEqLeft
\norm{x+y}_{\Gamma,p} 
\ov{\eqref{def-gamma-psi-p}}{=} \max\left\{ \bnorm{\Gamma(x+y,x+y)^{\frac12}}_{\L^p(\VN(G))}, \bnorm{\Gamma(x^*+y^*,x^* + y^*)^{\frac12}}_{\L^p(\VN(G))} \right\} \\
&\ov{\eqref{ine-6543}\eqref{Equa-divers-100000}}{\leq} \max \left\{ \bnorm{\Gamma(x,x)^{\frac12}}_{\L^p(\VN(G))} + \bnorm{\Gamma(y,y)^{\frac12}}_{\L^p(\VN(G))},\right. \nonumber\\
&\qquad \qquad \qquad \qquad \qquad \qquad \left. \bnorm{\Gamma(x^*,x^*)^{\frac12}}_{\L^p(\VN(G))} + \bnorm{\Gamma(y^*,y^*)^{\frac12}}_{\L^p(\VN(G))} \right\} \nonumber\\
&\leq \max \left\{ \bnorm{\Gamma(x,x)^{\frac12}}_{\L^p(\VN(G))}, \bnorm{\Gamma(x^*,x^*)^{\frac12}}_{\L^p(\VN(G))} \right\} \nonumber \\
&\qquad \quad+ \max \left\{ \bnorm{\Gamma(y,y)^{\frac12}}_{\L^p(\VN(G))}, \bnorm{\Gamma(y^*,y^*)^{\frac12}}_{\L^p(\VN(G))} \right\}
\ov{\eqref{def-gamma-psi-p}}{=} \norm{x}_{\Gamma,p} + \norm{y}_{\Gamma,p}. \nonumber
\end{align*}
\end{proof}

Now, we prove the Leibniz property of these seminorms.
We recall that $\mathrm{C}^*_r(G)$ is the reduced group $\mathrm{C}^*$ algebra containing $\P_G$ as a dense subspace.

\begin{prop}
\label{Prop-Leibniz-groupes-Lp}
Suppose $2 \leq p<\infty$.  Let $G$ be a discrete group. For any $x,y \in \P_G$, we have
\begin{equation}
\label{Leibniz-rule-groupes}
\norm{xy}_{\Gamma,p} 
\leq \norm{x}_{\mathrm{C}^*_r(G)} \norm{y}_{\Gamma,p} +  \norm{x}_{\Gamma,p}\norm{y}_{\mathrm{C}^*_r(G)}.
\end{equation}
\end{prop}

\begin{proof}
Let $x,y \in \P_G$. Using the structure of bimodule of $\L^p(\VN(G),\L^2(\Omega)_{c,p})$ in the second inequality, we see that
\begin{align*}
\MoveEqLeft
\bnorm{\Gamma(xy,xy)^{\frac{1}{2}}}_{\L^p(\VN(G))} \\
&\ov{\eqref{nabla-norm-Lp-gradient-psi}}{=}\bnorm{\partial_\psi(xy)}_{\L^p(\VN(G),\L^2(\Omega)_{c,p})}
\ov{\eqref{Leibniz-Schur-gradient-mieux-sgrp}}{=} \bnorm{x\partial_\psi(y)+\partial_\psi(x)y}_{\L^p(\VN(G),\L^2(\Omega)_{c,p})}\\
&\leq \bnorm{x\partial_\psi(y)}_{\L^p(\VN(G),\L^2(\Omega)_{c,p})}+\bnorm{\partial_\psi(x)y}_{\L^p(\VN(G),\L^2(\Omega)_{c,p})}\\
&\leq \norm{x}_{\mathrm{C}^*_r(G)} \bnorm{\partial_\psi(y)}_{\L^p(\VN(G),\L^2(\Omega)_{c,p})}+\bnorm{\partial_\psi(x)}_{\L^p(\VN(G),\L^2(\Omega)_{c,p})} \norm{y}_{\mathrm{C}^*_r(G)}\\
&\ov{\eqref{nabla-norm-Lp-gradient-psi}}{=} \norm{x}_{\mathrm{C}^*_r(G)} \bnorm{\Gamma(y,y)^{\frac{1}{2}}}_{\L^p(\VN(G))}+\bnorm{\Gamma(x,x)^{\frac{1}{2}}}_{\L^p(\VN(G))} \norm{y}_{\mathrm{C}^*_r(G)}.
\end{align*}
Similarly, we have
\begin{align*}
\MoveEqLeft
\bnorm{\Gamma(y^*x^*,y^*x^*)^{\frac{1}{2}}}_{\L^p(\VN(G))} \\
&\leq \norm{x^*}_{\mathrm{C}^*_r(G)} \bnorm{\Gamma(y^*,y^*)^{\frac{1}{2}}}_{\L^p(\VN(G))} + \bnorm{\Gamma(x^*,x^*)^{\frac{1}{2}}}_{\L^p(\VN(G))} \norm{y^*}_{\mathrm{C}^*_r(G)} \\
&=\norm{x}_{\mathrm{C}^*_r(G)} \bnorm{\Gamma(y^*,y^*)^{\frac{1}{2}}}_{\L^p(\VN(G))} + \bnorm{\Gamma(x^*,x^*)^{\frac{1}{2}}}_{\L^p(\VN(G))} \norm{y}_{\mathrm{C}^*_r(G)}.
\end{align*} 
Therefore, $\norm{\cdot}_{\Gamma,p}$ satisfies the Leibniz property.
\end{proof}

For the proof of Theorem \ref{thm-Fourier-quantum-compact-metric}, we shall need the following short lemma. Note that $\dom \partial_{\psi,q,p}$ was defined in Proposition \ref{Prop-derivation-closable-sgrp}.
\begin{lemma}
\label{lem-partial-closed-rad-Fourier}
Let $2 \leq p < \infty$ and $-1 \leq q \leq 1$. Let $G$ be a discrete group. Suppose that $\L^p(\VN(G))$ has $\CCAP$ and that $\VN(G)$ has $\QWEP$. The map $\partial_{\psi,q,p} \co \dom \partial_{\psi,q,p} \subseteq \L^p(\VN(G)) \to \L^p(\VN(G),\L^2(\Gamma_q(H))_{\rad,p})$ is closed.
\end{lemma}

\begin{proof}
Let $(x_n)$ be a sequence in $\dom \partial_{\psi,q,p}$ such that $x_n \to x$ where $x$ is some element of $\L^p(\VN(G))$ and $\partial_{\psi,q,p}(x_n) \to z$ where $z \in \L^p(\VN(G),\L^2(\Gamma_q(H))_{\rad,p})$. We have $\bnorm{\partial_{\psi,q,p}(y)}_{\L^p(\VN(G),\L^2(\Gamma_q(H))_{\rad,p})} \ov{\eqref{nabla-norm-Lp-gradient-psi}\eqref{relation-partial-star-group}}{\cong } \norm{y}_{\Gamma,p} \ov{\eqref{max-Gamma-2-group}}{\cong} \bnorm{A^{\frac12}_p(y)}_{\L^p(\VN(G))}$.  Thus the sequence $(A^{\frac12}_p(x_n))$ converges in $\L^p(\VN(G))$. Since $A^{\frac12}_p$ is closed, $x$ belongs to $\dom A^{\frac12}_p = \dom \partial_{\psi,q,p}$ and $A^{\frac12}_p(x) = \lim_n A^{\frac12}_p(x_n)$. According to \eqref{Equivalence-square-root-domaine-Schur-sgrp}, the sequence $(\partial_{\psi,q,p}(x_n))$ is then convergent in $\L^p(\Gamma_q(H) \rtimes_\alpha G)$. Then as $\partial_{\psi,q,p}$ is closed, seen as an unbounded operator $\dom \partial_{\psi,q,p} \subseteq \L^p(\VN(G)) \to \L^p(\Gamma_q(H) \rtimes_\alpha G)$ by Proposition \ref{Prop-derivation-closable-sgrp}, $x$ belongs to $\dom \partial_{\psi,q,p}$ and $\partial_{\psi,q,p}(x) = \lim_n \partial_{\psi,q,p}(x_n)$, limit in $\L^p(\Gamma_q(H) \rtimes_\alpha G)$. Since $p \geq 2$, by \cite[Proposition 2.5 (2)]{JMX} the latter space embeds into the space $\L^p(\VN(G),\L^2(\Gamma_q(H))_{\rad,p})$, and the above limit also holds in $\L^p(\VN(G),\L^2(\Gamma_q(H))_{\rad,p})$.
\end{proof}

In the next theorem, recall that if $A$ is the generator of a markovian semigroup of Fourier multipliers, then there exists a real Hilbert space $H$ together with a mapping $b_\psi \co G \to H$ such that the symbol $\psi \co G \to \C$ of $A$ satisfies $\psi(s) = \norm{b_\psi(s)}_H^2$. Following \cite[p.~1962]{JMP1}, we define
\begin{equation}
\label{equ-Delta-A-prime-sgrp}
\G_\psi 
\ov{\mathrm{def}}{=} \inf_{b_\psi(s) \neq b_\psi(t)} \norm{b_\psi(s) - b_\psi(t)}_H^2.
\end{equation}
By \cite[Proposition 2.10.2]{BHV}, note that $\G_\psi$ is independent of $b_\psi$, that is, if $b_\psi \co I \to H_b$ and $c_\psi \co I \to H_c$ are two cocycles such that $\psi(s) = \norm{b_\psi(s)}_{H_b}^2 = \norm{c_\psi(s)}_{H_c}^2$ for all $s \in G$, then $\G_\psi$ takes the same value in \eqref{equ-Delta-A-prime-sgrp}, whether it is defined via $b_\psi$ or $c_\psi$.

We denote by $\L^p_0(\VN(G))$ the range of the bounded projection $P \co \L^p(\VN(G)) \to \L^p(\VN(G))$, $\lambda_s \mapsto (1-\delta_{s = e})\lambda_s = (1-\tau(\lambda_s))\lambda_s$ (see the proof of Lemma \cite[Lemma 7.27]{ArK1} for a proof of the boundedness). Recall that $\mathrm{C}^*_r(G)$ is the reduced group $\mathrm{C}^*$-algebra associated with the discrete group $G$. We also write $\mathrm{C}^*_r(G)_{0}$ for the range of the bounded projection $\mathrm{C}^*_r(G) \to \mathrm{C}^*_r(G) \co \lambda_s \mapsto (1-\delta_{s = e}) \lambda_s$, on the space of elements of $\mathrm{C}^*_r(G)$ with vanishing trace (see again the proof of Lemma \cite[Lemma 7.27]{ArK1} for a proof of the boundedness). In the sequel, we denote by $\norm{\cdot}_{\Gamma,p}$ the restriction of \eqref{def-gamma-psi-p} on $\dom A^{\frac12}_p \cap \mathrm{C}^*_r(G)$. Note Remark \ref{Extension-quantum-metric-space}.

\begin{thm}
\label{thm-Fourier-quantum-compact-metric}
Let $2 \leq p < \infty$. Let $G$ be a discrete group. Let $b_\psi \co G \to H$ be an injective cocycle with values in a finite-dimensional Hilbert space of dimension $n<p$. Assume that $\G_\psi > 0$. 
\begin{enumerate}
	\item If $\L^p(\VN(G))$ has $\CCAP$ and $\VN(G)$ has $\QWEP$ then $(\mathrm{C}^*_r(G),\norm{\cdot}_{\Gamma,p})$ is a  quantum compact metric space.
	\item If $G$ is in addition weakly amenable then $(\mathrm{C}^*_r(G),\norm{\cdot}_{\Gamma,p})$ is a Leibniz quantum compact metric space.
\end{enumerate}
\end{thm}

\begin{proof}
We will use Proposition \ref{Prop-passagea-se}.

1. Since $\dom \norm{\cdot}_{\Gamma,p}=\dom A^{\frac12}_p \cap \mathrm{C}^*_r(G)$ contains the dense subspace $\P_G$ of $\mathrm{C}^*_r(G)$, the $\C$-subspace $\dom \norm{\cdot}_{\Gamma,p}$ is dense in $\mathrm{C}^*_r(G)$. By Lemma \ref{Lemma-Hilbert-module-groups} and Proposition \ref{Prop-derivation-closable-sgrp} $\dom \norm{\cdot}_{\Gamma,p}$ is closed under the adjoint operation. According to Proposition \ref{Prop-C1-Schur-infinite-sgrp}, $\norm{\cdot}_{\Gamma,p}$ is a seminorm. We will show that it is a Lip-norm, i.e. 
\begin{equation}
\label{equ-2-proof-prop-Fourier-quantum-compact-metric}
\left\{ x \in \dom \norm{\cdot}_{\Gamma,p} : \: \norm{x}_{\Gamma,p} = 0 \right\} 
= \C 1 .
\end{equation}
Indeed, we have $A(1) = \psi(e) 1 = 0$, so that $\Gamma(1,1)\ov{\eqref{Def-Gamma-sgrp}}{=}
\frac{1}{2}\big[ A(1^*) 1 + 1^* A(1) - A(1^* 1)\big] = 0$ and consequently $\norm{1}_{\Gamma,p} \ov{\eqref{def-gamma-psi-p}}{=} \max \big\{ \bnorm{\Gamma(1,1)^\frac12}_{\L^p(\VN(G))}, \bnorm{ \Gamma(1^*,1^*)^\frac12}_{\L^p(\VN(G))} \big\}= 0$. In the other direction, if $\norm{x}_{\Gamma,p} = 0$, then according to \eqref{max-Gamma-2-group}, we have $A_p^{\frac12}(x) = 0$. For any $s \in G$, we infer that 
$$
0
=\tau_G \big(A_p^{\frac12}(x) \lambda_s^*\big) 
=\tau_G \big(x A_p^{\frac12}((\lambda_s)^*)\big) 
=\psi(s^{-1})^{\frac12} \tau_G(x \lambda_s^*).
$$
Since $\psi(s^{-1}) \neq 0$ for $s \neq e$, we deduce that $\tau_G(x \lambda_s^*) = 0$ for these $s$, so that finally, $x \in \C 1$. Thus, $\norm{\cdot}_{\Gamma,p}$ is a Lip-norm. 

Since $\tau_G \co \mathrm{C}^*_r(G) \to \C$ is a state, with Proposition \ref{Prop-norm-precompact} it suffices to show that
\begin{equation}
\label{equ-proof-prop-Fourier-quantum-compact-metric}
\left\{ x  \in \dom \norm{\cdot}_{\Gamma,p} : \norm{x}_{\Gamma,p} \leq 1 ,\: \tau_G(x) = 0 \right\} \text{ is relatively compact in }\mathrm{C}^*_r(G).
\end{equation}
Note that $A_2^{-1} \co \L^2_0(\VN(G)) \to \L^2_0(\VN(G))$ is compact. Indeed, $(\lambda_s)_{s \in G \backslash \{ e \}}$ is an orthonormal basis of $\L^2_0(\VN(G))$ consisting of eigenvectors of $A_2^{-1}$. Moreover, the corresponding eigenvalues are $\norm{b_\psi(s)}_H^{-2}$ vanish at infinity\footnote{\thefootnote. Recall that a family $(x_s)_{s \in I}$ vanishes at infinity means that for any $\epsi >0$, there exists a finite subset $J$ of $I$ such that for any $s \in I-J$ we have $|x_s| \leq \epsi$.}. For the latter, note that the condition $\G_\psi = \inf_{b_\psi(s) \neq b_\psi(t)} \norm{b_\psi(s) - b_\psi(t)}_H^2 > 0$ together with the injectivity of $b_\psi$ imply that any compact subset of $H$ meets the $b_\psi(s)$ only for a finite number of $s \in G$. As $H$ is finite-dimensional by assumption, the closed balls of $H$ are compact, so that $\norm{b_\psi(s)}_H^{-2} \to 0$. We have shown that $A_2^{-1}$ is compact. Now according to \cite[Lemma 5.8]{JMP1}, \cite[Theorem 1.1.7]{JM1} applied with $z=\frac{1}{2}$ and $q=\infty$, the operator $A^{-\frac12} \co \L^p_0(\VN(G)) \to \VN(G)_0$ is compact. So the set
$$
\Big\{ x \in \mathrm{C}^*_r(G)_{0} \cap \dom A_p^{\frac12} : \bnorm{A_p^{\frac12}(x)}_{\L^p(\VN(G))} \leq 1 \Big\} 
$$ 
is relatively compact in $\mathrm{C}^*_r(G)$. Note that for any $x \in \dom A_p^{\frac12} $ we have
\begin{align*}
\MoveEqLeft
\bnorm{A_p^{\frac12}(x)}_{\L^p(\VN(G))} 
\ov{\eqref{max-Gamma-2-group}}{\lesssim_{p}} \max \Big\{ \bnorm{\Gamma(x,x)^{\frac12} }_{\L^p(\VN(G))}, \bnorm{\Gamma(x^*,x^*)^\frac12}_{\L^p(\VN(G))} \Big\} 
\ov{\eqref{def-gamma-psi-p}}{=} \norm{x}_{\Gamma,p}.           
\end{align*} 
We deduce from this inequality that 
$$
\big\{x \in \mathrm{C}^*_r(G)_{0} \cap \dom A_p^{\frac12} : \norm{x}_{\Gamma,p} \leq 1  \big\}
$$
is relatively compact in $\mathrm{C}^*_r(G)$. We have
$$
\big\{x \in \mathrm{C}^*_r(G)_{0} \cap \dom A_p^{\frac12}: \norm{x}_{\Gamma,p} \leq 1  \big\} 
= \big\{ x \in \mathrm{C}^*_r(G) \cap \dom A_p^{\frac12}: \norm{x}_{\Gamma,p} \leq 1, \: \tau_G(x) = 0 \big\}.
$$ 
We deduce \eqref{equ-proof-prop-Fourier-quantum-compact-metric}.

%

2. Let $x \in \mathrm{C}^*_r(G)$ and $(x_n)$ be a sequence of elements of $\dom A^{\frac12}_p \cap \mathrm{C}^*_r(G)$ such that $\norm{x_n-x}_{\mathrm{C}^*_r(G)} \to 0$ and $\norm{x_n}_{\Gamma,p} \ov{\eqref{nabla-norm-Lp-gradient-psi}}{=} \bnorm{\partial_{\psi,1,p}(x_{n})}_{\L^p(\VN(G),\L^2(\Omega)_{\rad,p})} \leq 1$. We have to show that $x$ belongs to $\dom A^{\frac12}_p$ and that $\bnorm{\partial_{\psi,1,p}(x)}_{\L^p(\VN(G),\L^2(\Omega)_{\rad,p})} \leq 1$. Note that 
\[\norm{x_n-x}_{\L^p(\VN(G))} \to 0.\] 
Then $\big((x_n,\partial_{\psi,1,p}(x_n)\big)$ is a sequence of elements of the graph of $\partial_{\psi,1,p}$ (the operator $\partial_{\psi,1,p}$ is seen as in Lemma \ref{lem-partial-closed-rad-Fourier}), which is bounded for the graph norm. Note that this closed graph is convex, hence weakly closed. Since bounded sets are weakly relatively compact, there exists a subnet $\big((x_{n_i},\partial_{\psi,1,p}(x_{n_i})\big)$ which converges weakly to an element $\big(z,\partial_{\psi,1,p}(z)\big)$ of the graph of $\partial_{\psi,1,p}$. In particular, $(x_{n_i})$ converges weakly to $z$ and $\partial_{\psi,1,p}(x_{n_i})$ converges weakly to $\partial_{\psi,1,p}(z)$. Then necessarily $x=z$. We conclude that $x$ belongs to $\dom \partial_{\psi,1,p} = \dom A_p^{\frac{1}{2}}$. Moreover, with \cite[Theorem 2.5.21]{Meg1}, we have by weak convergence
\begin{align*}
\MoveEqLeft
\bnorm{\partial_{\psi,1,p}(x)}_{\L^p(\VN(G),\L^2(\Omega)_{\rad,p})} 
\leq \liminf_{i} \bnorm{\partial_{\psi,1,p}(x_{n_i})}_{\L^p(\VN(G),\L^2(\Omega)_{\rad,p})}
\leq 1.
\end{align*}

Now, suppose that $G$ is weakly amenable. We will extend the Leibniz property \eqref{Leibniz-rule-groupes} to elements $x,y$ of $\dom A^{\frac12}_p \cap \mathrm{C}^*_r(G)$. Let $(\varphi_j)$ be a net of functions $\varphi_j \co G \to \C$ with finite support such that the net $(M_{\varphi_j})$ converges to $\Id_{\L^p(\VN(G))}$ and to $\Id_{\mathrm{C}^*_\lambda(G)}$ in the point-norm topologies with $\sup_j \norm{M_{\varphi_j}}_{\cb,\L^p(\VN(G)) \to \L^p(\VN(G))} < \infty$ (see the end of Subsection \ref{subsubsec-Fourier-mult-crossed-product}). We have $x = \lim_j M_{\varphi_j}(x)$ and $y = \lim_j M_{\varphi_j}(y)$ where the convergence holds in $\mathrm{C}^*_r(G)$ and in $\L^p(\VN(G))$. 

Thus also $xy = \lim_j M_{\varphi_j}(x)M_{\varphi_j}(y)$ in $\mathrm{C}^*_r(G)$. Using the lower semicontinuity of $\norm{\cdot}_{\Gamma,p}$ in the first inequality, we conclude that
\begin{align*}
\MoveEqLeft
\norm{xy}_{\Gamma,p} 
\leq \liminf_j \norm{M_{\varphi_j}(x) M_{\varphi_j}(y)}_{\Gamma,p} \\
&\ov{\eqref{Leibniz-rule-groupes}}{\leq} \liminf_j \norm{M_{\varphi_j}(x)}_{\mathrm{C}^*_r(G)} \norm{M_{\varphi_j}(y)}_{\Gamma,p} + \norm{M_{\varphi_j}(y)}_{\mathrm{C}^*_r(G)} \norm{M_{\varphi_j}(x)}_{\Gamma,p} \\
&\ov{\eqref{Def-Gamma-group}}{=} \norm{x}_{\mathrm{C}^*_r(G)} \norm{y}_{\Gamma,p} + \norm{y}_{\mathrm{C}^*_r(G)} \norm{x}_{\Gamma,p}.
\end{align*}
Here, we have also used \eqref{equ-Fourier-mult-on-column-space} in the last line to see that $M_{\check{\varphi}_j}$ are equally uniformly completely bounded Fourier multipliers approximating the identity.
The theorem is proved.
\end{proof}

\begin{remark} \normalfont
The argument for the lower semicontinuity is correct with the weaker assumption ``$(x_n)$ converges weakly to $x$ in $\L^p(\VN(G))$'' instead of $\norm{x_n-x}_{\mathrm{C}^*_r(G)} \to 0$. In the part 2, the argument proves that $\dom A^{\frac12}_p \cap \mathrm{C}^*_r(G)$ is a $*$-subalgebra of $\mathrm{C}^*_r(G)$.
\end{remark}

\begin{remark}
\label{rem-Fourier-quantum-compact-metric-exponent-restriction}
\normalfont
Note that in Theorem \ref{thm-Fourier-quantum-compact-metric} above, in contrast to its counterpart for Schur multipliers, Theorem \ref{th-quantum-metric-space} below, we have a restriction of the exponent $p > n$, where $n$ denotes the dimension of the representing Hilbert space. This restriction is in fact necessary in general as the example in the following proof shows.
\end{remark}

\begin{proof}
Consider the abelian discrete group $G = \Z^2$, the canonical identification $b_\psi \co \Z^2 \to \R^2$ of $\Z^2$ into $\R^2$ and the trivial homomorphism $\alpha \co \Z^2 \to \O(\R^2),\:(n,m) \mapsto \Id_{\R^2}$. The associated length function is $\psi \co \Z^2 \to \R_+, \: (n,m) \mapsto \psi(n,m) = \norm{b_\psi(n,m)}_{\R^2}^2 = n^2 + m^2$. Now pick the critical exponent $p = 2 = \dim \R^2$. Note that the operator space $\L^2(\VN(\Z^2)) = \L^2(\T^2)$ is an operator Hilbert space by \cite[Proposition 2.1 (iii)]{Pis4} and consequently has $\CCAP$\footnote{\thefootnote. More generally, if $1<p<\infty$ the operator space $\L^p(\T^2)$ has $\CCAP$ by \cite[Proposition 3.5]{JR1}
since an abelian group is weakly amenable with Cowling-Haagerup constant equal to 1.}. We also calculate 
\begin{align*}
\MoveEqLeft
\G_\psi 
\ov{\eqref{equ-Delta-A-prime-sgrp}}{=}\inf_{b_\psi(n,m) \neq b_\psi(n',m')} \norm{b_\psi(n,m) - b_\psi(n',m')}_{\R^2}^2 \\
&=\inf_{(n,m) \neq (n',m')} \norm{(n,m) - (n',m')}_{\R^2}^2
=1 
> 0.            
\end{align*}
However, $(\mathrm{C}^*_r(\Z^2),\norm{\cdot}_{\Gamma,2})$ is not a quantum compact metric space by \cite[Theorem 2.4.3]{Lat2}. Indeed, the set
$$
\big\{ x \in \mathrm{C}^*_r(\Z^2)_0 \cap \dom A^{\frac12}_2 : \: \norm{x}_{\Gamma,2} \leq 1 \big\}
$$
is not bounded, so in particular, not relatively compact. Indeed, consider the double sequence $(\alpha_{n,m})_{n,m \in \Z}$ defined by
\begin{equation}
\label{Def-alpha-m-n}
\alpha_{n,m} 
\ov{\mathrm{def}}{=} \frac{1}{(1+n^2 + m^2) \log(2 + n^2 + m^2)} \delta_{(n,m) \neq (0,0)}.
\end{equation}
Fix a number $N \in \N$ and consider the selfadjoint element $x_N \ov{\mathrm{def}}{=} \sum_{n^2 + m^2 \leq N} \alpha_{n,m} \e^{2\pi \i \langle (n,m),\cdot \rangle}$ of $\mathrm{C}^*_r(\Z^2)=\mathrm{C}(\T^2)$. We claim that $(x_N)$ is an unbounded sequence in $\mathrm{C}^*_r(\Z^2)_0=\mathrm{C}(\T^2)_0$ but bounded in $\norm{\cdot}_{\Gamma,2}$-seminorm. Indeed, observe that $x_N$ is a trigonometric polynomial without constant term, so belonging indeed to $\mathrm{C}^*_r(\Z^2)_0=\mathrm{C}(\T^2)_0$. Moreover, $\norm{x_N}_{\mathrm{C}(\T^2)} \geq |x_N(0)| = \sum_{n^2 + m^2 \leq N} \alpha_{n,m}$. But the series $\sum_{(n,m) \in \Z^2} \alpha_{n,m}$ diverges since using an integral test \cite[Proposition 7.57]{GhL1} and a change of variables to polar coordinates we have
\begin{align*}
\sum_{(n,m) \in \Z^2} \alpha_{n,m} 
&\cong \int_{\R^2} \frac{1}{(1+x^2 +y^2)\log(2+x^2 +y^2)} \d x \d y \\
&= C\int_0^\infty \frac{1}{(1+r^2)\log(2+r^2)} r \d r 
=\infty.
\end{align*}
On the other hand, using Plancherel theorem in the second equality and again an integral test and a change of variables to polar coordinates, we obtain 

\begin{align*}
\MoveEqLeft
\norm{x_N}_{\Gamma,2}^2 
\ov{\eqref{max-Gamma-2-group}}{\cong} \bnorm{A^{\frac12}(x_N)}_{\L^2(\T^2)}^2 
=\norm{ \left(\delta_{n^2+m^2 \leq N} \psi(n,m)^{\frac12} \alpha_{n,m}\right)_{(n,m) \in \Z^2} }_{\ell^2_{\Z^2}}^2 \\
&=\sum_{n^2 + m^2 \leq N} (n^2 + m^2) \alpha_{n,m}^2 
\ov{\eqref{Def-alpha-m-n}}{\leq} \sum_{(n,m) \in \Z^2} \frac{1}{(1 + n^2 + m^2) \left[\log(2+n^2 + m^2)\right]^2} \\
&\cong \int_{\R^2} \frac{1}{(1+x^2 +y^2)\left[\log(2+x^2 +y^2)\right]^2} \d x \d y 
=C \int_0^\infty \frac{1}{(1+r^2)\left[\log(2+r^2)\right]^2} r \d r < \infty.
\end{align*}
\end{proof}


We refer to the survey \cite{Cha1} and to \cite[pp.~628-629]{JM1} for more information. Recall that a finitely generated discrete group $G$ has rapid decay of order $r \geq 0$ if we have an estimate
$$
\norm{x}_{\VN(G)}
\lesssim k^r \norm{x}_{\L^2(\VN(G))}
$$
for any $x=\sum_{|s| \leq k} a_s \lambda_s$. Similarly using \cite[Lemma 1.3.1]{JM1}, we can obtain the following result.

\begin{thm}
\label{Th-quantum-Rdecay}
Let $G$ be a finitely generated discrete group with rapid $r$-decay. Suppose that there exists $\beta >0$ such that $2\frac{2r+1}{\beta} \leq p < \infty$ and that $\inf_{|s|=k} \psi(s) \gtrsim k^\beta$.
\begin{enumerate}
	\item If $\L^p(\VN(G))$ has $\CCAP$ and $\VN(G)$ has $\QWEP$ then $(\mathrm{C}^*_r(G),\norm{\cdot}_{\Gamma,p})$ is a quantum compact metric space.
	\item If $G$ is in addition weakly amenable then $(\mathrm{C}^*_r(G),\norm{\cdot}_{\Gamma,p})$ is a Leibniz quantum compact metric space.
\end{enumerate} 
\end{thm}

\begin{proof}
Using the notation $n \ov{\mathrm{def}}{=} 2\frac{2r+1}{\beta}$, we have by \cite[Lemma 1.3.1]{JM1} the estimate
$$
\norm{T_t}_{\L^2_0(\VN(G))\to \VN(G)}
\lesssim \frac{1}{t^{\frac{2r+1}{2\beta}}}
=\frac{1}{t^{\frac{n}{4}}}, \quad t>0.
$$
By \cite[Lemma 1.1.2]{JM1}, we deduce that $\norm{T_t}_{\L^1_0(\VN(G))\to \VN(G)} \lesssim \frac{1}{t^{\frac{n}{2}}}$. Moreover, by the argument of \cite[p.~628]{JM1}, the operator $A^{-1}$ is compact on $\L^2_0(\VN(G))$. Since $2\frac{2r+1}{\beta} \leq p$, we can use \cite[Theorem 1.1.7]{JM1} with $z=\frac{1}{2}$ and $q=\infty$. We deduce that $A^{-\frac{1}{2}} \co \L^p_0(\VN(G)) \to \VN(G)_0$ is compact. Hence the set $
\big\{ x \in \mathrm{C}^*_r(G)_{0} \cap \dom A_p^{\frac12} : \bnorm{A_p^{\frac12}(x)}_{\L^p(\VN(G))} \leq 1 \big\}$ is relatively compact in $\mathrm{C}^*_r(G)$. The end of the proof is similar to the proof of Theorem \ref{thm-Fourier-quantum-compact-metric} (note that the condition of rapid decay implies that we have a Lip-norm).
\end{proof}

\subsection{Gaps and estimates of norms of Schur multipliers}
\label{Sec-some-estimates}

Consider a semigroup $(T_t)_{t \geq 0}$ of selfadjoint unital completely positive Schur multipliers on $\B(\ell^2_I)$ as in \eqref{Semigroup-Schur}.
In order to find quantum compact metric spaces associated with $(T_t)_{t \geq 0}$ as defined in Subsection \ref{subsec-prelims-quantum}, we need some supplementary information on the semigroup.
One important of them is the notion of the gap, which we study in this subsection.
In all the subsection, we suppose $\dim H <\infty$.
We define the gap of $\alpha$ by
\begin{equation}
\label{equ-Delta-A-prime}
\G_\alpha \ov{\textrm{def}}{=} \inf_{\alpha_i - \alpha_j \neq \alpha_k - \alpha_l} \norm{(\alpha_i - \alpha_j) - (\alpha_k - \alpha_l)}_{H}^2.
\end{equation}
Note that by the proof of \cite[Proposition 5.4]{Arh1} and \cite[Theorem C.2.3]{BHV} then $\G_\alpha$ is independent of $\alpha$, that is, if $\alpha \co I \to H_\alpha,\:\beta \co I \to H_\beta$ are two families such that \eqref{equ-Schur-markovian-alpha} and \eqref{Semigroup-Schur} hold for both, then $\G_\alpha = \G_\beta$.

\begin{lemma} 
\label{CountingLemma}
If $\dim H = n$ and $\G_\alpha > 0$, for any integer $k \geq 1$, we have
$$
\card \Big\{ \alpha_i - \alpha_j : k^2 \G_\alpha
\leq \norm{\alpha_i - \alpha_j}_H^2 
\leq (k+1)^2 \G_\alpha \Big\} 
\leq (5^n-1) k^{n-1}.
$$
\end{lemma}

\begin{proof} 
If $\mathrm{B}_n$ denotes the open Euclidean unit ball in $H$ and if $\xi_1,\xi_2$ are distinct and can be written $\xi_1=\alpha_i - \alpha_j$ and $\xi_2=\alpha_k - \alpha_l$ for some $i,j,k,l \in I$, we have\footnote{\thefootnote. If $\eta \in \big( \xi_1 + \frac{\sqrt{\G_\alpha}}{2} \mathrm{B}_n \big) \cap \big( \xi_2 + \frac{\sqrt{\G_\alpha}}{2} \mathrm{B}_n \big)$ then $|\xi_1-\xi_2| \leq |\xi_1-\eta|+|\eta-\xi_2| < \frac{\sqrt{\G_\alpha}}{2}+\frac{\sqrt{\G_\alpha}}{2}=\sqrt{\G_\alpha}$ which is impossible.} 
$$
\Big( \xi_1 + \frac{\sqrt{\G_\alpha}}{2} \mathrm{B}_n \big) \cap \big( \xi_2 + \frac{\sqrt{\G_\alpha}}{2} \mathrm{B}_n \big) 
=\emptyset.
$$ 
Counting the maximum number of disjoint balls of radius $\frac{\sqrt{\G_\alpha}}{2}$ in the annulus $\Big( (k+1) \sqrt{\G_\alpha} + \frac{\sqrt{\G_\alpha}}{2} \Big) \mathrm{B}_n-\Big( k \sqrt{\G_\alpha} - \frac{\sqrt{\G_\alpha}}{2} \Big) \mathrm{B}_n$ combined with the binomial theorem, we obtain
\begin{align*}
\MoveEqLeft
\card  \Big\{ \alpha_i -\alpha_j : k^2 \G_\alpha \leq \norm{\alpha_i - \alpha_j}_H^2
\leq (k+1)^2 \G_\alpha \Big\} \\ 
&\leq \frac{\vol\Big(\Big( (k+1) \sqrt{\G_\alpha} + \frac{\sqrt{\G_\alpha}}{2} \Big) \mathrm{B}_n\Big) - \vol\Big( \Big( k \sqrt{\G_\alpha} - \frac{\sqrt{\G_\alpha}}{2} \Big) \mathrm{B}_n \Big)}{\vol\Big(\frac{\sqrt{\G_\alpha}}{2} \mathrm{B}_n\Big)}  \\ 
&=(2k+3)^n - (2k-1)^n 
=\sum_{j=0}^n {{n}\choose{j}} (2k)^{n-j} \big( 3^j - (-1)^j \big) \\
& \leq k^{n-1}\sum_{j=1}^n {{n}\choose{j}} 2^{n-j}\big(3^j - (-1)^j \big) \\ 
&\leq k^{n-1}\sum_{j=0}^n {{n}\choose{j}} 2^{n-j}3^j - k^{n-1}\sum_{j=0}^n 2^{n-j}(-1)^j 
= (5^n-1) k^{n-1}.  
\end{align*}
\end{proof}

We need the following lemma which tells that Schur multipliers with $0-1$ entries of diagonal block rectangular shape are completely contractive. 

\begin{lemma}
\label{lem-block-Schur}
Let $I$ be a non-empty index set.
Let $\{I_1,\ldots,I_N\}$ and $\{J_1,\ldots,J_N\}$ be two subpartitions of $I$, i.e. $I_k \subseteq I$ and $I_k \cap I_l = \emptyset$ for $k \neq l$ (and similarly for $J_1,\ldots,J_N$). The matrix $B = [b_{ij}]$, where
$$
b_{ij} 
=\begin{cases} 
1 & \text{if } i \in I_k,\:j \in J_k \text{ for the same }k \\
0 &  \text{otherwise} 
\end{cases},
$$
induces a completely contractive Schur multiplier $M_B \co S^\infty_I \to S^\infty_I$ and is also completely contractive on $S^p_I$ for any $1 \leq p \leq \infty$.
\end{lemma}

\begin{proof}
For any $i \in \{ 1 , \ldots , N \}$ consider the orthogonal projections $P_i \co \ell^2_I \to \ell^2_I$, $(\xi_k)_{k \in I} \to (\xi_k 1_{k \in I_i})_{k \in I}$ and $Q_i \co \ell^2_I \to \ell^2_I$, $(\xi_k)_{k \in I} \to (\xi_k 1_{k \in J_i})_{k \in I}$. If $i \neq j$, the ranges $P_i(\ell^2_I)$ and $P_j(\ell^2_I)$ (resp. $Q_i(\ell^2_I)$ and $Q_j(\ell^2_I)$) are orthogonal. Moreover, if $C \in S^\infty_I$ it is obvious that $M_B(C) = \sum_{i = 1}^N P_i C Q_i$. Then for any $\xi,\eta \in \ell^2_I$, using the Cauchy-Schwarz Inequality in the second inequality, we obtain
\begin{align*}
\MoveEqLeft
\big| \big\langle M_B(C)\xi, \eta \big\rangle_{\ell^2_I}\big| 
=\left| \bigg\langle \sum_{i=1}^N P_i C Q_i \xi, \eta \bigg\rangle_{\ell^2_I} \right|
=\left| \sum_{i=1}^N \langle C Q_i \xi, P_i \eta \rangle_{\ell^2_I} \right| 
\leq  \sum_{i=1}^N \big|\langle C Q_i \xi, P_i \eta \rangle_{\ell^2_I} \big| \\
&\leq \left( \sum_{i=1}^N \|CQ_i\xi\|_{\ell^2_I}^2 \right)^{\frac12} \left( \sum_{i = 1}^N \|P_i \eta\|_{\ell^2_I}^2 \right)^{\frac12} 
\leq \norm{C}_{S^\infty_I} \left( \sum_{i = 1}^N \norm{Q_i\xi}_{\ell^2_I}^2 \right)^{\frac12} \left( \sum_{i = 1}^N \norm{P_i \eta}_{\ell^2_I}^2 \right)^{\frac12} \\
&=\norm{C}_{S^\infty_I}   \norm{\sum_{i=1}^NQ_i\xi}_{\ell^2_I}   \norm{\sum_{i=1}^N P_i \eta}_{\ell^2_I}
\leq \norm{C}_{S^\infty_I} \norm{\xi}_{\ell^2_I} \norm{\eta}_{\ell^2_I}.
\end{align*}
We infer that $\norm{M_B(C)}_{S^\infty_I} \leq \norm{C}_{S^\infty_I}$. Thus $M_B \co S^\infty_I \to S^\infty_I$ is a contraction, and since it is a Schur multiplier, it is even a complete contraction \cite[Corollary 8.8]{Pau}.
Moreover, since $B$ has only real entries, $M_B$ is symmetric, so also (completely) contractive on $S^1_I$, and by interpolation also on $S^p_I$ for $1 \leq p \leq \infty$.
\end{proof}

Recall that $\SinftyI$ satisfies \eqref{equa-Sp0}. The following result is related to the Hardy-Littlewood-Sobolev theory \cite[Chapter II]{VSCC1}.

\begin{lemma} 
\label{LemWS}
If $\dim H = n$, $ \G_\alpha> 0$ and if the function $\widetilde{m} \co H \to \C$ satisfies
$$
\big| \widetilde{m}(\xi) \big|  
\leq c_n \norm{\xi}_H^{-(n+\epsi)} \quad \text{for some } \epsi > 0,
$$ 
then the Schur multiplier $M_{[\widetilde{m}(\alpha_i-\alpha_j)]} \co S^\infty_I \to S^\infty_I$ is completely bounded. Moreover, for any $t>0$, we have 
$$
\norm{T_t}_{\cb,\SinftyI \to \SinftyI} 
\leq \frac{c(n)}{(\G_\alpha\, t)^{\frac{n}{2}}}.
$$
\end{lemma}

\begin{proof}
Let $x = \sum_{i,j \in I'} x_{ij} \ot e_{ij} \in S^\infty_\ell(S^\infty_{I'}) \subseteq S^\infty_\ell(S^\infty_I)$, where $I' \subseteq I$ is finite and $\ell$ is an integer.
Let $\{I_1,\ldots,I_N\}$ be the partition of $I'$ corresponding to the equivalence relation $i \cong j \ov{\textrm{def}}{\Longleftrightarrow} \alpha_i = \alpha_j$, and let $J_1 = I_1, \ldots, J_N = I_N$. Let $B_0$ be the Schur multiplier symbol from Lemma \ref{lem-block-Schur} associated with these partitions. 

Moreover for $\xi \in \alpha(I') - \alpha(I')$, we let $\{I_1^\xi,\ldots,I_N^\xi\}$ and $\{J_1^\xi,\ldots,J_N^\xi\}$ be the (possibly empty) subpartitions of $I'$ such that $\alpha_i - \alpha_j = \xi \Leftrightarrow i \in I_k^\xi$ and $j \in J_k^\xi$ for the same $k$. We let $B_\xi$ be the Schur multiplier symbol from Lemma \ref{lem-block-Schur} associated with these partitions. Then, using Lemma \ref{CountingLemma} and our growth assumption on $\widetilde{m}$ in the third inequality, we obtain
\begin{align}
\MoveEqLeft
\label{Equa-divers-678}
\norm{\sum_{i,j \in I'} \widetilde{m}(\alpha_i - \alpha_j) x_{ij} \ot e_{ij}}_{S^\infty_\ell(S^{\infty}_I)} \nonumber \\ \nonumber
&=\Bgnorm{\sum_{\substack{i,j \in I' \\ \alpha_i = \alpha_j}} \widetilde{m}(\alpha_i - \alpha_j) x_{ij} \ot e_{ij}+\sum_{k \geq 1} \sum_{\substack{\xi \in \alpha(I') - \alpha(I') \\\nonumber
k^2 \G_\alpha \leq \norm{\xi}^2 < (k+1)^2 \G_\alpha}}\widetilde{m}(\alpha_i - \alpha_j) x_{ij} \ot e_{ij} }_{S^\infty_\ell(S^{\infty}_I)} \\ \nonumber
&\leq |\widetilde{m}(0)| \norm{ \sum_{\substack{i,j \in I' \\ \alpha_i = \alpha_j}} x_{ij} \ot e_{ij} }_{S^\infty_\ell(S^{\infty}_I)} \\
& + \sum_{k \geq 1} \sum_{\substack{\xi \in \alpha(I') - \alpha(I') \\ \nonumber
k^2 \G_\alpha \leq \norm{\xi}^2 < (k+1)^2 \G_\alpha}}
|\widetilde{m}(\xi)| \, \norm{ \sum_{\substack{i,j \in I' \\ \alpha_i - \alpha_j = \xi}} x_{ij} \ot e_{ij} }_{S^\infty_\ell(S^{\infty}_I)} \\ \nonumber 
&= |\widetilde{m}(0)| \bnorm{(\Id_{S^\infty_\ell} \ot M_{B_0})(x) }_{S^\infty_\ell(S^{\infty}_I)} \\ 
&+ \sum_{k \geq 1} \sum_{\begin{subarray}{c} \xi \in \alpha(I') - \alpha(I') \\ 
k^2 \G_\alpha \leq \|\xi\|^2 < (k+1)^2 \G_\alpha \end{subarray}}
|\widetilde{m}(\xi)| \, \bnorm{(\Id_{S^\infty_\ell} \ot M_{B_\xi})(x) }_{S^\infty_\ell(S^{\infty}_I)} \\ \nonumber
& \leq \bigg( |\widetilde{m}(0)| + \sum_{k \geq 1} \sum_{\begin{subarray}{c} \xi \in \alpha(I') - \alpha(I') \\ 
k^2 \G_\alpha \leq \norm{\xi}^2 < (k+1)^2 \G_\alpha \end{subarray}} |\widetilde{m}(\xi)| \bigg) \, \norm{x}_{S^\infty_\ell(S^\infty_I)} \\ 
& \leq 5^n \bigg( |\widetilde{m}(0)| + \sum_{k \geq 1} k^{n-1} \left(k \sqrt{\G_\alpha}\right)^{-(n+\epsi)} \bigg) \norm{x}_{S^\infty_\ell(S^\infty_I)} \\\nonumber
&= 5^n \bigg( |\widetilde{m}(0)| + \sum_{k \geq 1} k^{-1-\epsi} \left(\G_\alpha\right)^{-\frac{n-\epsi}{2}} \bigg) \norm{x}_{S^\infty_\ell(S^\infty_I)}
= c_{n,\epsi}(\G_\alpha) \norm{x}_{S^\infty_\ell(S^\infty_I)}. \nonumber
\end{align}  
Note that the Riemann series $\sum_{k \geq 1} k^{-1-\epsi}$ converges. 

For the second assertion, we use the function $\widetilde{m}=\e^{-t\norm{\cdot}_H}$. Since $x \in \SinftyI$ we may ignore the term $|\widetilde{m}(0)|$ above and from \eqref{Equa-divers-678}, it suffices to use the inequality 
$$
\dsp \sum_{k = 1}^\infty k^{n-1} \e^{-t \G_\alpha k^2} \leq  c(n) (\G_\alpha t)^{-\frac{n}{2}}.
$$
\end{proof}

\subsection{Seminorms associated to semigroups of Schur multipliers}
\label{Section-Leibniz-Schur}

In this subsection, we consider again a markovian semigroup $(T_t)_{t \geq 0}$ of Schur multipliers on $\B(\ell^2_I)$ satisfying Proposition \ref{def-Schur-markovian}. Recall from Subsection \ref{Sec-infos-on-Gamma} that we have a generator $A_p$ and also a carr\'e du champ $\Gamma$. We shall introduce a family of seminorms that will give rise in Subsection \ref{Sec-I-finite} below to some quantum compact metric spaces in the sense of Subsection \ref{subsec-prelims-quantum}. First suppose $2 \leq p<\infty$. For any $x \in \dom A^{\frac12}_p$ we let
\begin{equation}
\label{def-gamma-alpha-p>-2}
\norm{x}_{\Gamma,p}
\ov{\mathrm{def}}{=}
\max \Big\{ \bnorm{\Gamma(x,x)^\frac12}_{S^p_I}, \bnorm{ \Gamma(x^*,x^*)^\frac12}_{S^p_I} \Big\}.
\end{equation}
Hereby, we recall the definition of $\Gamma$ from \eqref{Def-Gamma} extended to $\dom A^{\frac12}_p$ in Lemma \ref{lem-Gamma-closure}.

For the proof of Proposition \ref{Prop-Leibniz-Schur-infinite}, we shall need the following short lemma.
Note that $\dom \partial_{\alpha,q,p}$ was defined in Proposition \ref{Prop-derivation-closable}.
\begin{lemma}
\label{lem-partial-closed-rad}
Let $2 \leq p < \infty$ and $-1 \leq q \leq 1$.
The mapping $\partial_{\alpha,q,p} \co \dom \partial_{\alpha,q,p} \subseteq S^p_I \to S^p_I(\L^2(\Gamma_q(H))_{\rad,p})$ is closed.
\end{lemma}

\begin{proof}
Let $(x_n)$ be a sequence in $\dom \partial_{\alpha,q,p}$ such that $x_n \to x$ for some element $x \in S^p_I$ and $\partial_{\alpha,q,p}(x_n) \to z$ where $z \in S^p_I(\L^2(\Gamma_q(H))_{\rad,p})$. We have $\bnorm{\partial_{\alpha,q,p}(y)}_{S^p_I(\L^2(\Gamma_q(H))_{\rad,p})} \ov{\eqref{Last-equality} \eqref{relation-partial-star}}{\cong} \norm{y}_{\Gamma,p} \ov{\eqref{equivalence-Ap-Schur-dom}}{\cong} \bnorm{A^{\frac12}_p(y)}_{S^p_I}$. Thus the sequence $(A^{\frac12}_p(x_n))$ converges in $S^p_I$. Since $A^{\frac12}_p$ is closed, $x$ belongs to $\dom A^{\frac12}_p = \dom \partial_{\alpha,q,p}$ and $A^{\frac12}_p(x) = \lim_n A^{\frac12}_p(x_n)$. According to \eqref{Equivalence-square-root-domaine-Schur}, the sequence $(\partial_{\alpha,q,p}(x_n))$ is then convergent in $\L^p(\Gamma_q(H) \otvn \B(\ell^2_I))$. Then as $\partial_{\alpha,q,p}$ is closed, seen as an unbounded operator $\dom \partial_{\alpha,q,p} \subseteq S^p_I \to \L^p(\Gamma_q(H) \otvn \B(\ell^2_I))$ by Proposition \ref{Prop-derivation-closable}, $x$ belongs to $\dom \partial_{\alpha,q,p}$ and $\partial_{\alpha,q,p}(x) = \lim_n \partial_{\alpha,q,p}(x_n)$, limit in $\L^p(\Gamma_q(H) \otvn \B(\ell^2_I))$. Since $p \geq 2$, by \cite[Proposition 2.5 (2)]{JMX} the latter space embeds into $S^p_I(\L^2(\Gamma_q(H))_{\rad,p})$, and the above limit also holds in $S^p_I(\L^2(\Gamma_q(H))_{\rad,p})$.
\end{proof}

\begin{prop}
\label{Prop-Leibniz-Schur-infinite}
Suppose $2 \leq p<\infty$. 
\begin{enumerate}
\item $\norm{\cdot}_{\Gamma,p}$ is a seminorm on $\dom A^{\frac12}_p$.
\item $\norm{\cdot}_{\Gamma,p}$ is lower semicontinuous on $\dom A^{\frac12}_p$ equipped with the topology induced by the weak topology of $S^p_I$.
\item For any $x,y \in \dom A^{\frac12}_p$, we have
\begin{equation}
\label{equ-1-Prop-Leibniz-Schur-infinite}
\norm{xy}_{\Gamma,p} 
\leq \norm{x}_{S^\infty_I} \norm{y}_{\Gamma,p}  + \norm{y}_{S^\infty_I} \norm{x}_{\Gamma,p}.
\end{equation}
\end{enumerate}

\end{prop}

\begin{proof}
1. If $x$ belongs to $\dom A^{\frac12}_p$ and $k \in \C$, we have
\begin{align*}
\MoveEqLeft
\norm{k x}_{\Gamma,p}
\ov{\eqref{def-gamma-alpha-p>-2}}{=} \max \Big\{ \bnorm{\Gamma(k x,k x)^\frac12}_{S^p_I}, \bnorm{ \Gamma((k x)^*,(k x)^*)^\frac12}_{S^p_I} \Big\} \\  
&=\max \Big\{ \bnorm{\Gamma(k x,k x)^\frac12}_{S^p_I}, \bnorm{ \Gamma(\ovl{k} x^*,\ovl{k} x^*)^\frac12}_{S^p_I}\Big\}
=\max \Big\{ |k| \bnorm{\Gamma(x,x)^\frac12}_{S^p_I}, |k|\bnorm{ \Gamma(x^*,x^*)^\frac12}_{S^p_I} \Big\}\\
&=|k|\max \Big\{ \bnorm{\Gamma(x,x)^\frac12}_{S^p_I}, \bnorm{ \Gamma(x^*,x^*)^\frac12}_{S^p_I} \Big\}
\ov{\eqref{def-gamma-alpha-p>-2}}{=}|k|\norm{x}_{\Gamma,p}.
\end{align*} 
Let us turn to the triangular inequality. Assume that $x,y \in \M_{I,\fin}$. Note that $S^{\frac{p}{2}}_I$ is a normed space since $p \geq 2$. According to the part 4 of Lemma \ref{Lemma-gamma-infos} and \eqref{Inegalite-triangulaire-1} applied with $Z=S^{\frac{p}{2}}_I$, we have
$$
\norm{\Gamma(x+y,x+y)}_{S^{\frac{p}{2}}_I}^{\frac{1}{2}}           
\ov{\eqref{Inegalite-triangulaire-1}}{\leq} \norm{\Gamma(x,x)}_{S^{\frac{p}{2}}_I}^{\frac{1}{2}}+\norm{\Gamma(y,y)}_{S^{\frac{p}{2}}_I}^{\frac{1}{2}}.
$$
We can rewrite this inequality under the form
\begin{equation}
\label{ine-65437}
\bnorm{\Gamma(x+y,x+y)^{\frac12}}_{S^p_I} 
\leq \bnorm{\Gamma(x,x)^{\frac12}}_{S^p_I} + \bnorm{\Gamma(y,y)^{\frac12}}_{S^p_I}.
\end{equation}
Then, for the general case where $x,y \in \dom A^{\frac12}_p$ consider some sequences $(x_n)$ and $(y_n)$ of $\M_{I,\fin}$ such that $x_n \xra[]{a} x$ and $y_n \xra[]{a} y$ with Lemma \ref{lem-Gamma-closure}. By Remark \ref{Remark-xra}, we have $x_n+y_n \xra[]{a} x+y$ and thus
\begin{align*}
\MoveEqLeft
\bnorm{\Gamma(x+y,x+y)^{\frac12}}_{S^p_I} 
=\bnorm{\Gamma(x+y,x+y)}_{S^{\frac{p}{2}}_I}^{\frac12} 
\ov{\eqref{Def-closure}}{=} \lim_n \bnorm{\Gamma(x_n+y_n,x_n+y_n)}_{S^{\frac{p}{2}}_I}^{\frac12} \\
&=\lim_n \bnorm{\Gamma(x_n+y_n,x_n+y_n)^{\frac12}}_{S^p_I} 
\ov{\eqref{ine-65437}}{\leq} \lim_j \bnorm{\Gamma(x_n,x_n)^{\frac12}}_{S^p_I} + \lim_n \bnorm{\Gamma(y_n,y_n)^{\frac12}}_{S^p_I} \\
&\ov{\eqref{Def-closure}}{=}\bnorm{\Gamma(x,x)^{\frac12}}_{S^p_I} + \bnorm{\Gamma(y,y)^{\frac12}}_{S^p_I}.
\end{align*}

A similar reasoning for $x,y$ replaced by $x^*,y^*$ gives 
\begin{equation}
\label{Equa-divers-1000001}
\bnorm{\Gamma(x^*+y^*,x^* + y^*)^{\frac12}}_{S^p_I} \leq \bnorm{\Gamma(x^*,x^*)^{\frac12}}_{S^p_I} + \bnorm{\Gamma(y^*,y^*)^{\frac12}}_{S^p_I}.
\end{equation}
Now we have
\begin{align*}
\MoveEqLeft
\norm{x+y}_{\Gamma,p} 
\ov{\eqref{def-gamma-alpha-p>-2}}{=} \max\left\{ \bnorm{\Gamma(x+y,x+y)^{\frac12}}_{S^p_I}, \bnorm{\Gamma(x^*+y^*,x^* + y^*)^{\frac12}}_{S^p_I} \right\} \\
&\ov{\eqref{ine-65437}\eqref{Equa-divers-1000001}}{\leq} \max \left\{ \bnorm{\Gamma(x,x)^{\frac12}}_{S^p_I} + \bnorm{\Gamma(y,y)^{\frac12}}_{S^p_I},\right. \nonumber\\
&\qquad \qquad \qquad \qquad \qquad \qquad \left. \bnorm{\Gamma(x^*,x^*)^{\frac12}}_{S^p_I} + \bnorm{\Gamma(y^*,y^*)^{\frac12}}_{S^p_I} \right\} \nonumber\\
&\leq \max \left\{ \bnorm{\Gamma(x,x)^{\frac12}}_{S^p_I}, \bnorm{\Gamma(x^*,x^*)^{\frac12}}_{S^p_I} \right\} \nonumber \\
&\qquad \quad+ \max \left\{ \bnorm{\Gamma(y,y)^{\frac12}}_{S^p_I}, \bnorm{\Gamma(y^*,y^*)^{\frac12}}_{S^p_I} \right\}
\ov{\eqref{def-gamma-alpha-p>-2}}{=} \norm{x}_{\Gamma,p} + \norm{y}_{\Gamma,p}. \nonumber
\end{align*}

2. Let $x \in S^p_I$ and $(x_\beta)$ be a net of elements of $\dom A^{\frac12}_p $ such that $(x_\beta)$ converges to $x$ for the weak topology of $S^p_I$  and, recalling the part 2 of Lemma \ref{Lemma-Hilbert-module2}, 
\[ 
\bnorm{x_\beta}_{\Gamma,p} 
= \bnorm{\partial_{\alpha,1,p}(x_{\beta})}_{S^p_I(\L^2(\Omega)_{\rad,p})} \leq 1.
\]
We have to show that $x$ belongs to $\dom A^{\frac12}_p$ and that $\bnorm{\partial_{\alpha,1,p}(x)}_{S^p_I(\L^2(\Omega)_{\rad,p})} \leq 1$. Note that the net $(\partial_{\alpha,1,p}(x_\beta))$ is bounded in $S^p_I(\L^2(\Omega)_{\rad,p})$ so admits a weakly convergent subnet $\big(\partial_{\alpha,1,p}(x_{\beta_i})\big)$. Then $\big((x_{\beta_i},\partial_{\alpha,1,p}(x_{\beta_i})\big)$ is a weakly convergent net in the graph of $\partial_{\alpha,1,p}$. Note that this graph is closed according to Lemma \ref{lem-partial-closed-rad} and convex, hence weakly closed. Thus the limit of $\big((x_{\beta_i},\partial_{\alpha,1,p}(x_{\beta_i})\big)$ belongs again to the graph and is of the form $\big(z,\partial_{\alpha,1,p}(z)\big)$ for some $z \in \dom \partial_{\alpha,1,p}$. In particular, $(x_{\beta_i})$ converges weakly to $z$ and $\partial_{\alpha,1,p}(x_{\beta_i})$ converges weakly to $\partial_{\alpha,1,p}(z)$. Then necessarily $x=z$. We conclude that $x$ belongs to $\dom \partial_{\alpha,1,p} \ov{}{=} \dom A_p^{\frac{1}{2}}$. Moreover, with \cite[Theorem 2.5.21]{Meg1}, we have
\begin{align*}
\MoveEqLeft
\bnorm{\partial_{\alpha,1,p}(x)}_{S^p_I(\L^2(\Omega)_{\rad,p})} 
\leq \liminf_{i} \bnorm{\partial_{\alpha,1,p}(x_{\beta_i})}_{S^p_I(\L^2(\Omega)_{\rad,p})}
\leq 1.
\end{align*}

3. Suppose first that $x,y \in \M_{I,\fin}$. Then using the structure of bimodule of $S^{p}_I(\L^2(\Omega)_{c,p})$ in the second inequality, we see that
\begin{align*}
  \bnorm{\Gamma(xy,xy)^{\frac{1}{2}}}_{S^p_I}
	&\ov{\eqref{nabla-norm-Lp-gradient}}{=}\bnorm{\partial_\alpha(xy)}_{S^{p}_I(\L^2(\Omega)_{c,p})}
	\ov{\eqref{Leibniz-Schur-gradient-mieux}}{=} \bnorm{x\partial_\alpha(y)+\partial_\alpha(x)y}_{S^{p}_I(\L^2(\Omega)_{c,p})}\\
&\leq \bnorm{x\partial_\alpha(y)}_{S^{p}_I(\L^2(\Omega)_{c,p})}+\bnorm{\partial_\alpha(x)y}_{S^{p}_I(\L^2(\Omega)_{c,p})}\\
&\leq \norm{x}_{\B(\ell^2_I)} \bnorm{\partial_\alpha(y)}_{S^{p}_I(\L^2(\Omega)_{c,p})}+\bnorm{\partial_\alpha(x)}_{S^{p}_I(\L^2(\Omega)_{c,p})} \norm{y}_{\B(\ell^2_I)}\\
&\ov{\eqref{nabla-norm-Lp-gradient}}{=} \norm{x}_{\B(\ell^2_I)} \bnorm{\Gamma(y,y)^{\frac{1}{2}}}_{S^p_I}+\bnorm{\Gamma(x,x)^{\frac{1}{2}}}_{S^p_I} \norm{y}_{\B(\ell^2_I)}.
\end{align*}
Similarly, we have
$$
\bnorm{\Gamma(y^*x^*,y^*x^*)^{\frac{1}{2}}}_{S^p_I}
\leq \norm{x^*}_{\B(\ell^2_I)} \bnorm{\Gamma(y^*,y^*)^{\frac{1}{2}}}_{S^p_I} + \bnorm{\Gamma(x^*,x^*)^{\frac{1}{2}}}_{S^p_I} \norm{y^*}_{\B(\ell^2_I)}.
$$ 
Therefore, $\norm{\cdot}_{\Gamma,p}$ satisfies the Leibniz condition \eqref{equ-1-Prop-Leibniz-Schur-infinite} for any $x,y \in \M_{I,\fin}$.

We will extend the Leibniz property \eqref{equ-1-Prop-Leibniz-Schur-infinite} to elements $x,y$ of $\dom A^{\frac12}_p$. We have $x = \lim_J \Tron_J(x)$ and $y = \lim_J \Tron_J(y)$ where the convergence holds in $S^\infty_I$ and in $S^p_I$. Thus also $xy = \lim_J \Tron_J(x)\Tron_J(y)$ in $S^p_I$ since $S^p$ is a Banach algebra for the usual product of operators by \cite[p.~225]{BLM}. Using the lower semicontinuity of $\norm{\cdot}_{\Gamma,p}$ in the first inequality, we conclude that $xy \in \dom A^{\frac12}_p$ and that
\begin{align*}
\norm{xy}_{\Gamma,p} 
& \leq \liminf_J \norm{\Tron_J(x) \Tron_J(y)}_{\Gamma,p} 
\ov{\eqref{equ-1-Prop-Leibniz-Schur-infinite}}{\leq} \liminf_J \norm{\Tron_J(x)}_{S^\infty_I} \norm{\Tron_J(y)}_{\Gamma,p} + \norm{\Tron_J(y)}_{S^\infty_I} \norm{\Tron_J(x)}_{\Gamma,p} \\
&\ov{\eqref{Def-Gamma-schur}}{=} \norm{x}_{S^\infty_I} \norm{y}_{\Gamma,p} + \norm{y}_{S^\infty_I} \norm{x}_{\Gamma,p}.
\end{align*}
The proposition is proved.
\end{proof}

\begin{remark} \normalfont
The argument proves that $\dom A^{\frac12}_p$ is a $*$-subalgebra of $S^\infty_I$.
\end{remark}

Suppose $1 <p \leq 2$. For any element $x$ of $\M_{I,\fin}$, we let
\begin{equation}
\label{def-gamma-alpha-p-2}
\norm{x}_{\Gamma,p}
\ov{\textrm{def}}{=}
\inf_{x=y+z} \Big\{ \bnorm{\Gamma(y,y)^{\frac12}}_{S^p_I}^{\frac{p}{2}} +\bnorm{\Gamma(z^*,z^*)^{\frac12}}_{S^p_I}^{\frac{p}{2}} \Big\}^{\frac{2}{p}}
\end{equation}
where the infimum is taken over all $y, z \in \M_{I,\fin}$ such that $x=y+z$.

\begin{prop}
For $1 < p \leq 2$, $\norm{\cdot}_{\Gamma,p}$ is a $\frac{p}{2}$-seminorm on $\M_{I,\fin}$, that is, the triangle inequality holds under the form $\norm{x+y}_{\Gamma,p}^{\frac{p}{2}} \leq \norm{x}_{\Gamma,p}^{\frac{p}{2}} + \norm{y}_{\Gamma,p}^{\frac{p}{2}}$.
\end{prop}

\begin{proof}
 The homogeneity of $\norm{\cdot}_{\Gamma,p}$ can be shown as in the case $p \geq 2$ and is left to the reader. Let us turn to the triangle inequality. For any $y,y' \in \M_{I,\fin}$, the part 4 of Lemma \ref{Lemma-gamma-infos} says that
\begin{equation}
\label{Gamma-CS-norm-prime}
\bnorm{\Gamma(y,y')}_{S^{\frac{p}{2}}_I}^{\frac{p}{2}}
\leq \bnorm{\Gamma(y,y)^{\frac{1}{2}}}_{S^p_I}^{\frac{p}{2}} \bnorm{\Gamma(y',y')^{\frac{1}{2}}}^{\frac{p}{2}}_{S^{p}_I}.	
\end{equation}
Recall that $S^{\frac{p}{2}}_I$ is a $\frac{p}{2}$-normed space for $p <  2$. With this inequality, we can now estimate for arbitrary $y,y' \in \M_{I,\fin}$ 
\begin{align*}
\MoveEqLeft
\bnorm{\Gamma(y+y',y+y')^{\frac12}}_{S^p_I}^p
= \bnorm{\Gamma(y+y',y+y')}_{S^{\frac{p}{2}}_I}^{\frac{p}{2}} 
=\bnorm{\Gamma(y,y) + \Gamma(y',y') + \Gamma(y,y') + \Gamma(y',y)}_{S^{\frac{p}{2}}_I}^{\frac{p}{2}} \nonumber \\
&\leq \bnorm{\Gamma(y,y)}_{S^{\frac{p}{2}}_I}^{\frac{p}{2}} + \bnorm{\Gamma(y',y')}_{S^{\frac{p}{2}}_I}^{\frac{p}{2}} + \bnorm{\Gamma(y,y')}_{S^{\frac{p}{2}}_I}^{\frac{p}{2}} + \bnorm{\Gamma(y',y)}_{S^{\frac{p}{2}}_I}^{\frac{p}{2}}\nonumber\\
& \ov{\eqref{Gamma-CS-norm-prime}}{\leq} \bnorm{\Gamma(y,y)^{\frac12}}_{S^p_I}^p + \bnorm{\Gamma(y',y')^{\frac12}}_{S^p_I}^p + 2 \bnorm{\Gamma(y,y)^{\frac12}}_{S^p_I}^{\frac{p}{2}} \bnorm{\Gamma(y',y')^{\frac12}}_{S^p_I}^{\frac{p}{2}}\nonumber\\
&= \left( \bnorm{\Gamma(y,y)^{\frac12}}_{S^p_I}^{\frac{p}{2}} + \bnorm{\Gamma(y',y')^{\frac12}}_{S^p_I}^{\frac{p}{2}} \right)^2. \nonumber
\end{align*}
Taking the square roots, we obtain
\begin{equation}
\label{norm-Gamma-triangular-inequality-prime}
\bnorm{\Gamma(y+y',y+y')^{\frac12}}_{S^p_I} ^{\frac{p}{2}}
\leq \bnorm{\Gamma(y,y)^{\frac12}}_{S^p_I}^{\frac{p}{2}} + \bnorm{\Gamma(y',y')^{\frac12}}_{S^p_I}^{\frac{p}{2}}.
\end{equation}
Now consider some elements $x$ and $x'$ of $\M_{I,\fin}$ and some decompositions $x=y+z$ and $x'=y'+z'$ where $y,z,y',z' \in \M_{I,\fin}$. We have $x+x'=y+y'+z+z'$. We obtain
\begin{align*}
\MoveEqLeft
 \norm{x+x'}_{\Gamma,p}^{\frac{p}{2}}        
		\ov{\eqref{def-gamma-alpha-p-2}}{\leq} \bnorm{\Gamma(y+y',y+y')^{\frac12}}_{S^p_I}^{\frac{p}{2}}+\bnorm{\Gamma((z+z')^*,(z+z')^*)^{\frac12}}_{S^p_I}^{\frac{p}{2}}\\
		&\ov{\eqref{norm-Gamma-triangular-inequality-prime}}{\leq} \bnorm{\Gamma(y,y)^{\frac12}}_{S^p_I}^{\frac{p}{2}} + \bnorm{\Gamma(y',y')^{\frac12}}_{S^p_I}^{\frac{p}{2}}+\bnorm{\Gamma(z^*,z^*)^{\frac12}}_{S^p_I}^{\frac{p}{2}} + \bnorm{\Gamma(z'^*,z'^*)^{\frac12}}_{S^p_I}^{\frac{p}{2}}.
\end{align*}
Passing to the infimum, we conclude that
$$
\norm{x+x'}_{\Gamma,p}^{\frac{p}{2}}
\leq \norm{x}_{\Gamma,p}^{\frac{p}{2}} +\norm{x'}_{\Gamma,p}^{\frac{p}{2}}.
$$
\end{proof}

Suppose $1 < p \leq \infty$. Let $I$ be an index set. Recall that $\SpI$ satisfies \eqref{equa-Sp0}. Suppose that the map $\alpha \co I \to H$ is injective. In this case, $\SpI$ is the space of elements of $S^p_{I}$ with null diagonal. We will also use the subspace $\MIfinzero$ of matrices with a finite number of non null entries with null diagonal, and $\big(\dom A^{\frac12}_p \big)_0 = \dom A^{\frac12}_p \cap \SpI$ the space of those matrices in $\dom A^{\frac12}_p$ that have null diagonal.

Note that  we can see $S^\infty_I \oplus \C\Id_{\ell^2_I}$ as the standard unitization of the $\mathrm{C}^*$-algebra $S^\infty_I$. We consider the subspace $(\MIfinzero) \oplus \C\Id_{\ell^2_I}$ of $S^\infty_I \oplus \C\Id_{\ell^2_I}$. 

Suppose $1 < p \leq \infty$. We define the (resp. $\frac{p}{2}$-) seminorm $\norm{\cdot}_{\Gamma,\alpha,p}\ov{\textrm{def}}{=}\norm{\cdot}_{\Gamma,p} \oplus 0$ on $\big(\dom A^{\frac12}_p\big)_0 \oplus \C\Id_{\ell^2_I}$ (resp. $\MIfinzero \oplus \C\Id_{\ell^2_I}$). That is, for any element $x=x_0+ \lambda \Id_{\ell^2_I}$ of $\big(\dom A^{\frac12}_p\big)_0 \oplus \C\Id_{\ell^2_I}$ (resp. $\MIfinzero \oplus \C\Id_{\ell^2_I}$), we have
\begin{equation}
\label{def-gamma-alpha-p}
\norm{x}_{\Gamma,\alpha,p}
\ov{\textrm{def}}{=}
\begin{cases} 
\inf_{x_0=y+z} \Big\{ \bnorm{\Gamma(y,y)^{\frac12}}_{S^p_I}^{\frac{p}{2}} +\bnorm{\Gamma(z^*,z^*)^{\frac12}}_{S^p_I}^{\frac{p}{2}} \Big\}^{\frac{2}{p}}
&\text{ if } 1 < p \leq 2\\
\max \Big\{ \bnorm{\Gamma(x_0,x_0)^\frac12}_{S^p_I}, \bnorm{ \Gamma(x_0^*,x_0^*)^{\frac12}}_{S^p_I} \Big\}	
& \text{ if } p \geq 2
\end{cases},
\end{equation}
where the infimum is taken over all $y, z \in \M_{I,\fin}$ such that $x_0=y+z$.
Note that when $p = \infty$, we suppose in addition that $I$ is finite, in order to have a well-defined domain of the seminorm above.

\begin{prop}
\label{Prop-C1-Schur}
Suppose $2 \leq p \leq \infty$.
If $p=\infty$, we suppose that $I$ is finite. We have 
\begin{equation}
\label{Lip-norm-Schur-2}
\Big\{x \in \big(\dom A^{\frac12}_p\big)_0 \oplus \C\Id_{\ell^2_I}: \norm{x}_{\Gamma,\alpha,p} = 0\Big\} 
=\C \Id_{\ell^2_I}.
\end{equation}

\end{prop}

\begin{proof}
\textit{Case $2 \leq p<\infty$.} Let $x=x_0+\lambda \Id_{\ell^2_I}$ be an element of $\big(\dom A^{\frac12}_p\big)_0 \oplus \C\Id_{\ell^2_I}$. According to Theorem \ref{cor-Riesz-equivalence-Schur}, we have $\bnorm{A_p^{\frac12}(x_0)}_{S^p_I} \ov{\eqref{max-Gamma-2}}{\lesssim} \max \left\{ \bnorm{\Gamma(x_0,x_0)^{\frac12}}_{S^p_I}, \bnorm{\Gamma(x_0^*,x_0^*)^{\frac12}}_{S^p_I} \right\} \ov{\eqref{def-gamma-alpha-p-2}}{=} \norm{x_0}_{\Gamma,p}$. 
This implies that when $\norm{x}_{\Gamma,\alpha,p} = 0$, that is $\norm{x_0}_{\Gamma,p} = 0$, we have $\bnorm{A_p^{\frac12}(x_0)}_{S^p_I} = 0$ and finally $A_p^{\frac12}(x_0)=0$. 
We recall from \cite[(10.5) p.~364]{HvNVW2} that since $x_0$ belongs to $\ovl{\Ran A_p} = \ovl{\Ran A_p^{\frac12}}$, we have \[x_0 = \lim_{t \to 0+} A_p^{\frac12} (t+A_p^{\frac12})^{-1}(x_0) = \lim_{t \to 0+}(t+A_p^{\frac12})^{-1} A_p^{\frac12}(x_0) = 0 .\]
The reverse inclusion is true by \eqref{def-gamma-alpha-p}.

\textit{Case $p=\infty$ and $I$ finite.} Fix some $2<p_0<\infty$. Let $x=x_0+\lambda \Id_{\ell^2_I}$ be an element of $\big(\dom A^{\frac12}_p\big)_0 \oplus \C\Id_{\ell^2_I}$. Using Theorem \ref{cor-Riesz-equivalence-Schur} in the first inequality, we can write
\begin{align}
\MoveEqLeft
\label{Inegalites-Case-p=infty}
\bnorm{A^{\frac12}(x_0)}_{S^{p_0}_I} 
\ov{\eqref{max-Gamma-2}}{\lesssim_{p_0}}\max \Big\{ \bnorm{\Gamma(x_0,x_0)^{\frac12} }_{S^{p_0}_I}, \bnorm{\Gamma(x_0^*,x_0^*)^\frac12}_{S^{p_0}_I} \Big\} \\
&\lesssim_I \max \Big\{ \bnorm{\Gamma(x_0,x_0)^\frac12}_{S^\infty_I}, \bnorm{ \Gamma(x_0^*,x_0^*)^\frac12}_{S^\infty_I} \Big\}
\ov{\eqref{def-gamma-alpha-p}}{=}\norm{x_0}_{\Gamma,\infty}. \nonumber           
\end{align}
The end of the proof is similar to the case $2 \leq p<\infty$.
\end{proof}

\subsection{Quantum metric spaces associated to semigroups of Schur multipliers}
\label{Sec-I-finite}

By \cite[p.~3]{Rie4}, the linear space of all selfadjoint elements of a unital $\mathrm{C}^*$-algebra is an order-unit space. So if $I$ is an index set (finite or infinite) then $(S^\infty_I)_\sa \oplus \R\Id_{\ell^2_I}$ is an order-unit space. Moreover, by \cite[Proposition 2.3]{Rie3}, a subspace of an order-unit space containing the order-unit is itself also an order-unit space. We conclude that $(\SinftyI)_{\sa} \oplus \R\Id_{\ell^2_I}$ is an order-unit space. Now, we will prove in this subsection that $(\SinftyI)_{\sa} \oplus \R\Id_{\ell^2_I}$ equipped with the restriction of $\norm{\cdot}_{\Gamma,\alpha,p}$ (defined by \eqref{def-gamma-alpha-p}) to this space (also denoted by $\norm{\cdot}_{\Gamma,\alpha,p}$) is a quantum \textit{compact} metric space \cite{Rie3}. 


Moreover, since Lipschitz pairs in the sense of Definition \ref{Def-Lipschitz-pair} are only defined for the moment for seminorms and not for quasi-seminorms, in view of Proposition \ref{Prop-C1-Schur} above, we restrict our family $\norm{\cdot}_{\Gamma,\alpha,p}$ in the next theorem to the case $2 \leq p \leq \infty$.

Recall that if $A$ is a sectorial operator on a Banach space $X$ with $\omega(A) <\frac{\pi}{2}$ (see \eqref{equ-sectorial-operator}) and if $(T_t)_{t \geq 0}$ is the associated semigroup, then by \cite[Corollary 15.20]{KW}, we have the following representation for the (unbounded) operator $A^{-\frac12}$
\begin{equation}
\label{equ-3-lem-A12-bounded}
A^{-\frac12}(x) 
=\frac{1}{\Gamma\left(\frac12\right)} \int_0^\infty t^{-\frac12} T_t(x) \d t, \quad x \in \Ran A.
\end{equation}

\begin{lemma}
\label{lem-A12-bounded}
Let $A$ be a sectorial operator on a Banach space $X$ with $\omega(A) <\frac{\pi}{2}$ and $(T_t)_{t \geq 0}$ be the associated semigroup. Assume that $\norm{T_t}_{\Ran(A) \to X} \leq C \frac{1}{t^{d}}$ for some $d > \frac{1}{2}$, where $\Ran(A)$ is normed as a subspace of $X$. Then $A^{-\frac12}$, initially defined on $\Ran A$, extends to a bounded operator on $\ovl{\Ran A}$.
\end{lemma}

\begin{proof}
If $x \in \Ran A$, we have
\begin{align*}
\MoveEqLeft
\bnorm{A^{-\frac12}(x)}_X             
=\norm{\frac{1}{\Gamma\left(\frac12\right)} \int_0^\infty t^{-\frac12} T_t(x) \d t}_{X}\leq \frac{1}{\Gamma\left(\frac12\right)} \int_0^\infty t^{-\frac12} \bnorm{T_t(x)}_{X} \d t\\
&\lesssim \int_0^1 t^{-\frac12} \norm{T_t(x)}_X \d t + \int_1^\infty t^{-\frac12} \norm{T_t(x)}_X \d t \nonumber 
\lesssim \int_0^1 t^{-\frac12} \norm{x}_X \d t +\int_1^\infty t^{-\frac12-d} \norm{x}_X \d t \nonumber \\
& \leq \left( \int_0^1 t^{-\frac12} \d t+\int_1^\infty t^{-(d+\frac12)} \d t \right) \norm{x}_X.
\end{align*}
Then $A^{-\frac12}$ extends to a bounded operator on $\ovl{\Ran A}$.
\end{proof}

Note that if in Lemma \ref{lem-A12-bounded}, $(T_t)_{t \geq 0}$ is a markovian semigroup of Schur multipliers,
then $A_\infty$, the generator of $T_t$ on $S^\infty_I$ is only $\frac{\pi}{2}$ sectorial in general.
However, since $\ovl{\Ran A_\infty} = \ovl{A(\M_{I,\fin})}$ in this case, it is easy to check that the conclusion of Lemma \ref{lem-A12-bounded} holds under the decay condition of $\norm{T_t}_{\Ran A_\infty \to S^\infty_I}$.


In the proof of the following result, we consider $\SinftyI$ living in the non-unital $\mathrm{C}^*$-algebra $S^\infty_I$.

\begin{lemma}
\label{lem-Schur-compact}
Let $1 < p < \infty$.
Assume that the Hilbert space $H$ is of finite dimension $n \in \N$ and that $\G_\alpha > 0$.

\begin{enumerate}

\item The operator $A_p^{-\frac12} \co \SpI \to \SpI \subseteq \SinftyI$ is bounded.

\item Suppose $-1 \leq q \leq 1$. Let $B_p = (\Id_{\L^p(\Gamma_q(H))} \ot A_p)|_{\ovl{\Ran \partial_{\alpha,q,p}}} \co \dom B_p \subseteq\ovl{\Ran \partial_{\alpha,q,p}} \to \ovl{\Ran \partial_{\alpha,q,p}}$. Then the operator $B_p^{-\frac12}$ is bounded. 
\end{enumerate}
\end{lemma}

\begin{proof}
1. We begin by showing that $A_2^{-\frac12} \co \StwoI \to \StwoI$ is bounded. Note that by the gap condition, we have for any $i,j \in I$ such that $\alpha_i - \alpha_j \neq 0$, $\norm{\alpha_i - \alpha_j}_H = \norm{\alpha_i - \alpha_j - 0}\geq \G_\alpha$.
Since $a_{ij} = \norm{\alpha_i - \alpha_j}_H^2$, we deduce that the diagonal operator $A_2^{-\frac12}$ is bounded on the Hilbert space $\StwoI$. 

Next we show that $A^{-\frac12} \co \SinftyI \to \SinftyI$ is bounded. Indeed, note that $(T_t)_{t \geq 0}$ extends to a bounded $C_0$-semigroup on $\SinftyI$. Moreover, Lemma \ref{LemWS} (where in case that $H$ is $1$-dimensional, we inject $H$ beforehand into a $2$-dimensional Hilbert space so that this lemma yields an estimate $\norm{T_t}_{\infty \to \infty} \leq C t^{-\frac{n}2}$ with $\frac{n}2 > \frac12$) together with Lemma \ref{lem-A12-bounded} yield that $A^{-\frac12}$ is bounded on $\SinftyI$. 

Now, if $2 < p < \infty$, it suffices to interpolate the operator $A^{-\frac12}$ between levels $2$ and $\infty$ and to note that we have $\SpI = (\StwoI,\SinftyI)_\theta$ for the right $\theta \in [0,1)$.
Indeed, the interpolation identity follows from the fact that the $\SpI$ are complemented subspaces of $S^p_I$ by the spectral projections \cite[(10.1) p.~361]{HvNVW2} (note that $(T_t)_{t \geq 0}$ is a bounded semigroup on $S^p_I$ for $2 \leq p \leq \infty$), and that these spectral projections are compatible for different values of $2 \leq p \leq \infty$.
Then if $1 < p < 2$, we use that $\ovl{\Ran A_p}$ is the dual space of $\ovl{\Ran A_{p^*}}$ via the usual duality bracket, and that $A^{-\frac12}_p$ on the first space is the adjoint of $A^{-\frac12}_{p^*}$ on the second space.



2. We want to show that $\Id_{\L^p(\Gamma_q(H))} \ot A_p^{-\frac12} \co \ovl{\Ran \partial_{\alpha,q,p}} \to \L^p(\Gamma_q(H) \otvn \B(\ell^2_I))$ is bounded.
Since $\M_{I,\fin}$ is a core for $\partial_{\alpha,q,p}$ according to Proposition \ref{Prop-derivation-closable}, it suffices to bound on $\partial_{\alpha,q,p}(\M_{I,\fin})$, or even on the bigger subset $\L^p(\Gamma_q(H)) \ot \Ran A_p$ of $\L^p(\Gamma_q(H) \otvn \B(\ell^2_I))$. Note already that according to 1., $\Id_{\L^p}\ot A_p^{-\frac12}$ is well-defined on $\L^p \ot \Ran A_p$. Let $P \co S^p_I \to \SpI$ be the spectral projection from \cite[(10.1) p.~361]{HvNVW2}. Then for $x \in \L^p \ot \Ran A_p$, we have $(\Id_{\L^p}\ot A_p^{-\frac12})(x) = (\Id_{\L^p} \ot A_p^{-\frac12})P(x)$. As $\Gamma_q(H)$ has $\QWEP$ according to \cite{Nou1}, with \cite[p.~984]{Jun} we are reduced to show that $A_p^{-\frac12} P \co S^p_I \to S^p_I$ is completely bounded. To this end, let $n \in \N$ and consider the amplified operators $(T_t \ot \Id_{\B(\ell^2_n)})_{t \geq 0}$, which are again a markovian semigroup of Schur multipliers, associated with the mapping $\alpha_n \co I \times \{1 ,\ldots,n\},\:(i,k) \mapsto \alpha_i$. It is easy to check that $\G_{\alpha_n} = \G_\alpha > 0$. We have $\Id_{\S^p_n}\ot ( A_p^{-\frac12} P ) = \tilde{A}_p^{-\frac12} \tilde{P}$, where the operators on the right hand side are those associated with the amplified markovian semigroup. Since $\norm{\tilde{P}} \leq 1 + \sup_{t \geq 0} \norm{T_t \ot \Id_{\B(\ell^2_n)}}$ is bounded independent of $n$, and the bound on $\tilde{A}_p^{-\frac12}$, which depends on $\sup_{t > 0}\norm{T_t \ot \Id_{\B(\ell^2_n)}}$ and on $\G_{\alpha_n}$, is equally independent of $n$, we conclude that $\Id_{\L^p} \ot A_p^{-\frac12} : \ovl{\Ran \partial_{\alpha,q,p}} \to \L^p(\Gamma_q(H) \otvn \B(\ell^2_I))$ is bounded.
Note that since $\Id_{\L^p} \ot A_p^{-\frac12}$ leaves $\partial_{\alpha,q,p}(\M_{I,\fin})$ invariant, it is easy to check that it actually maps to $\ovl{\Ran \partial_{\alpha,q,p}}$.

\end{proof}

In the following theorem, we consider the abelian $\mathrm{C}^*$-subalgebra $\mathfrak{M}$ of $S^\infty_I$ consisting of its diagonal operators.
We recall that we restrict to $p \geq 2$ to have a seminorm $\norm{\cdot}_{\Gamma,\alpha,p}$ in the usual sense.

\begin{thm} 
\label{th-quantum-metric-space}
Assume that the Hilbert space $H$ is of finite dimension and that $\G_\alpha > 0$.
\begin{enumerate}
\item Suppose $2 \leq p \leq \infty$ and that $I$ is finite. Then $\big(\big(\SinftyI\big)_{\sa} \oplus \R\Id_{\ell^2_I},\norm{\cdot}_{\Gamma,\alpha,p}\big)$ is a Leibniz quantum compact metric space (Definition \ref{Def-Leibniz-pair}).

\item Suppose $2 \leq p < \infty$ and that $I$ is infinite. We consider the abelian $\mathrm{C}^*$-subalgebra $\mathfrak{M}$ of diagonal operators of $S^\infty_I$. Then $\big(\big(\SinftyI\big)_{\sa} \oplus \R \Id_{\ell^2_I},\norm{\cdot}_{\Gamma,\alpha,p},\mathfrak{M}\big)$ is a bounded Leibniz quantum locally compact metric space (see Definitions \ref{Defi-locally-order-unit}, \ref{def-Leibniz-local} and \ref{Defi-bounded-locally}). 
\end{enumerate}
\end{thm}

\begin{proof} 
1. Note that $\SinftyI$ is a subspace of the space of matrices of $S^\infty_I$ with null diagonal and that we have a contractive projection from $S^\infty_I$ onto the diagonal of $S^\infty_I$. So it is easy to check that $\tr \oplus \Id_\R \co \SinftyI \oplus \R\Id_{\ell^2_I} \to \R$, $x+\lambda \mapsto \lambda$ is a state of the order-unit space $(\SinftyI)_\sa \oplus \R\Id_{\ell^2_I}$. By Proposition \ref{Prop-norm-precompact}, it now suffices to show that
\begin{equation}
\label{equ-proof-prop-Schur-quantum-compact-metric-2}
\left\{ x  \in \dom \norm{\cdot}_{\Gamma,\alpha,p} : \norm{x}_{\Gamma,\alpha,p} \leq 1 ,\: (\tr \oplus \Id_\C)(x) = 0 \right\} \text{ is relatively compact in }S^\infty_{I} \oplus \C\Id_{\ell^2_I}.
\end{equation}

\textit{Case $2 \leq p<\infty$.} According to Lemma \ref{lem-Schur-compact}, the operator $A_p^{-\frac12} \co \SpI \to \SinftyI$ is bounded. Applying this operator to the ball of $\SpI$ of radius $\bnorm{A_p^{-\frac12}}^{-1}$, we obtain that the set\footnote{\thefootnote. Since $I$ is finite, we have $\dom A_p^{\frac12}=S^p_I$.}
$$
\Big\{ x \in \ovl{\Ran A_p^{\frac{1}{2}}} : \bnorm{A_p^{\frac12}(x)}_{S^p_I} \leq 1 \Big\} 
$$ 
is bounded in $S^\infty_I$. Note that for any $x \in S^p_I$ we have
\begin{align}
\label{EQ-divers-13346}
\MoveEqLeft
\bnorm{A_p^{\frac12}(x)}_{S^p_I} 
\ov{\eqref{equivalence-Ap-Schur-dom}}{\lesssim_{p}} \max \Big\{ \bnorm{\Gamma(x,x)^{\frac12} }_{S^p_I}, \bnorm{\Gamma(x^*,x^*)^\frac12}_{S^p_I} \Big\} 
\ov{\eqref{def-gamma-alpha-p-2}}{=} \norm{x}_{\Gamma,\alpha,p}.           
\end{align} 
Since $\ovl{\Ran A_p^{\frac{1}{2}}}=\SinftyI$, we deduce from this inequality that 
\begin{equation}
\label{eq-divers-643}
\big\{x \in S^p_I \cap \SinftyI : \norm{x}_{\Gamma,\alpha,p} \leq 1  \big\}
\end{equation}
is bounded, hence by finite dimensionality, relatively compact in $S^\infty_I$. We have
$$
\big\{x \in \dom A^{\frac12}_p \cap \SinftyI: \norm{x}_{\Gamma,\alpha,p} \leq 1  \big\} 
= \big\{ x \in \dom \norm{\cdot}_{\Gamma,\alpha,p} : \norm{x}_{\Gamma,\alpha,p} \leq 1, \: (\tr \oplus \Id_\C)(x) = 0 \big\}.
$$ 
We deduce \eqref{equ-proof-prop-Schur-quantum-compact-metric-2}.

\textit{Case $p=\infty$.} Fix some $2<p_0<\infty$. Again according to Lemma \ref{lem-Schur-compact}, the operator $A_{p_0}^{-\frac12} \co \ovl{\Ran A_{p_0}} \to \SinftyI$ is bounded. This implies that the set\footnote{\thefootnote. Since $I$ is finite, we have $\dom A_{p_0}^{\frac12}=S^{p_0}_I$.}
$$
\big\{ x \in \ovl{\Ran A_{p_0}^{\frac{1}{2}}}  : \bnorm{A_{p_0}^{\frac12}(x)}_{S^{p_0}_I} \leq 1 \big\} 
$$ 
is bounded in $S^\infty_I$. The inequality \eqref{Inegalites-Case-p=infty} says that the set 
$$
\big\{x \in \SinftyI \oplus \C\Id_{\ell^2_I} : \norm{x}_{\Gamma,\alpha,\infty} \leq 1 \text{ and } (\tr\oplus \Id_\C)(x) = 0 \big\}
=\big\{x \in \SinftyI : \norm{x}_{\Gamma,\alpha,\infty} \leq 1\big\}  
$$ 
is bounded, i.e. relatively compact in $S^\infty_I$. The desired result follows from Proposition \ref{Prop-norm-precompact}.


For the Leibniz rule, see the end of the proof.

2. Suppose that $J$ is a finite subset of $I$. We denote by $\mu_J \co S^\infty_I \to \C$ the \textit{normalized} ``partial trace'' where $\mu_J(x)$ is the sum of the diagonal entries of index belonging to $J$. This state is local since $\mu_J(\mathrm{diag}(0,\ldots,0,1,\ldots,1,0,\ldots,0))=1$ and since $\mathrm{diag}(0,\ldots,0,1,\ldots,1,0,\ldots,0)$ belongs to $\mathfrak{M}$. This state extends on the unitization by $\mu=\mu_J \oplus \Id_\C$.  We claim that if $a,b \in \mathfrak{M}$ are compactly supported then the set
\begin{equation}
\label{equ-proof-prop-Schur-quantum-compact-metric-3}
a \big\{x \in \SinftyI \oplus \C \Id_{\ell^2_I} : \norm{x}_{\Gamma,\alpha,p} \leq 1, \mu(x)=0 \big\} b
\end{equation}
is relatively compact. Recall that the Gelfand spectrum of $\mathfrak{M}$ is $I$ (we can identify $\mathfrak{M}$ with $c_0(I)$). Thus, compactly supported elements are those lying in $\Span\{e_{ii}:\:i \in I\}$. So we can write $a = \sum_{i \in J_a} a_i e_{ii}$ and $b = \sum_{j \in J_b} b_j e_{jj}$ where $J_a$ and $J_b$ are finite subsets of $I$. Now, we see that the above set is contained in $\Span\{e_{ij}:\: i \in J_a,\: j \in J_b \}$, which is a finite-dimensional space since $J_a$ and $J_b$ are finite. Thus it suffices to show boundedness of the above set \eqref{equ-proof-prop-Schur-quantum-compact-metric-3}. Using the fact that matrices in $\SinftyI$ have null diagonal, note that \eqref{equ-proof-prop-Schur-quantum-compact-metric-3} is equal to 
\begin{equation}
\label{Divers-3531}
a \big\{x \in \ovl{\Ran A_\infty} \cap \M_{J_a,J_b} : \norm{x}_{\Gamma,\alpha,p} \leq 1 \big\} b
\end{equation}
where $\M_{J_a,J_b} $ is the space of $J_a \times J_b$ matrices. Essentially by \eqref{eq-divers-643} (note that at that point, finiteness of $I$ was not used), the subset$ \big\{x \in \ovl{\Ran A_\infty} \cap \M_{J_a,J_b} : \norm{x}_{\Gamma,\alpha,p} \leq 1 \big\}$ is bounded. Hence the subset \eqref{Divers-3531} is also bounded.

Similarly, according to Lemma \ref{lem-Schur-compact}, the set
$$
\Big\{ x \in \dom A^{\frac12}_p \cap \SinftyI  : \bnorm{A_p^{\frac12}(x)}_{S^p_I} \leq 1 \Big\} 
$$ 
is bounded in $S^\infty_I$. By \eqref{EQ-divers-13346}, we deduce that 
$$
\big\{x \in \dom A^{\frac12}_p \cap \SinftyI : \norm{x}_{\Gamma,\alpha,p} \leq 1  \big\}
$$
is bounded in $S^\infty_I$. So the quantum locally compact metric space is bounded by Definition \ref{Defi-bounded-locally}.

We check that the seminorm $\norm{\cdot}_{\Gamma,\alpha,p}$ equally satisfies the Leibniz inequality
\begin{equation}
\label{equ-rem-Schur-Leibniz}
\norm{xy}_{\Gamma,\alpha,p} 
\leq \norm{x}_{\B(\ell^2_I)} \norm{y}_{\Gamma,\alpha,p} + \norm{x}_{\Gamma,\alpha,p} \norm{y}_{\B(\ell^2_I)}, 
\quad x,y \in \big(\dom A_p^{\frac12} \cap \SinftyI \big) \oplus \C \Id_{\ell^2_I}.
\end{equation}
Indeed, if we write $x = x_0 + \lambda \Id_{\ell^2_I}$ and $y = y_0 + \mu \Id_{\ell^2_I}$, then we have $xy = x_0y_0 + \lambda y_0 + \mu x_0 + \lambda \mu \Id_{\ell^2_I}$.
Thus, we have with limits in $S^p_I$ norm, noting that $\Tron_J(\Id_{\ell^2_I})\Tron_J(z) = \Tron_J(z) \Tron_J(\Id_{\ell^2_I}) = \Tron_J(z)$ and the product $S^p_I \times S^p_I \to S^p_I$ is continuous \cite[p.~225]{BLM},
\begin{align*}
\MoveEqLeft
xy - \lambda \mu \Id_{\ell^2_I} = x_0y_0 + \lambda y_0 + \mu x_0 = \lim_J \left[ \Tron_J(x_0)\Tron_J(y_0) + \lambda \Tron_J(y_0) + \mu \Tron_J(x_0) \right] \\
& = \lim_J \left[ \Tron_J(x_0+\lambda \Id_{\ell^2_I})\Tron_J(y_0) + \mu \Tron_J(x_0)\right]\\
& = \lim_J \left[ \Tron_J(x_0+\lambda \Id_{\ell^2_I})\Tron_J(y_0+\mu \Id_{\ell^2_I}) - \lambda \mu \Tron_J(\Id_{\ell^2_I})\Tron_J(\Id_{\ell^2_I}) \right].
\end{align*}
Therefore, we have $xy - \lambda\mu \Id_{\ell^2_I} = \lim_J \Tron_J(x)\Tron_J(y) - \Tron_J(\lambda \mu \Id_{\ell^2_I})$, limit in $S^p_I$.
Then by the lower semicontinuity of the seminorm $\norm{\cdot}_{\Gamma,p}$ from Proposition \ref{Prop-Leibniz-Schur-infinite} and then the Leibniz property \eqref{equ-1-Prop-Leibniz-Schur-infinite}, we have $xy \in \dom A^{\frac12}_p$ and 
\begin{align*}
\MoveEqLeft
\norm{xy}_{\Gamma,\alpha,p} 
=\bnorm{xy - \lambda\mu \Id_{\ell^2_I}}_{\Gamma,p} 
\leq \liminf_J \bnorm{\Tron_J(x)\Tron_J(y) - \Tron_J(\lambda \mu \Id_{\ell^2_I})}_{\Gamma,p} \\
& \leq \liminf_J \big(\norm{\Tron_J(x)\Tron_J(y)}_{\Gamma,p} + \bnorm{\Tron_J(\lambda \mu \Id_{\ell^2_I})}_{\Gamma,p} \big)
\overset{\eqref{Def-Gamma-schur}}{=} \liminf_J \norm{\Tron_J(x)\Tron_J(y)}_{\Gamma,p} + 0 \\
& \leq \liminf_J \norm{\Tron_J(x)}_{S^\infty_I} \norm{\Tron_J(y)}_{\Gamma,p} + \norm{\Tron_J(y)}_{S^\infty_I} \norm{\Tron_J(x)}_{\Gamma,p} \\
& = \liminf_J \norm{\Tron_J(x)}_{S^\infty_I} \norm{\Tron_J(y_0)}_{\Gamma,p} + \norm{\Tron_J(y)}_{S^\infty_I} \norm{\Tron_J(x_0)}_{\Gamma,p} \\
& \overset{\eqref{Def-Gamma-schur}}{=} \norm{x}_{\B(\ell^2_I)} \norm{y_0}_{\Gamma,p} + \norm{y}_{\B(\ell^2_I)} \norm{x_0}_{\Gamma,p} 
= \norm{x}_{\B(\ell^2_I)} \norm{y}_{\Gamma,\alpha,p} + \norm{y}_{\B(\ell^2_I)} \norm{x}_{\Gamma,\alpha,p}.
\end{align*}

Finally, we check both properties of the point b of Definition \ref{Def-Leibniz-pair} as in the proof of \cite[Proposition 2.17]{Lat7}.
Thus, $(\SinftyI \oplus \R \Id_{\ell^2_I},\norm{\cdot})$ is a unital Leibniz pair.

We check the semicontinuity property.
To this end, let $(x_n)_n$ be a sequence in $\dom A_p^{\frac12} \cap \SinftyI \oplus \C \Id_{\ell^2_I} \subseteq \B(\ell^2_I)$ such that $x_n \to x$ in $\B(\ell^2_I)$ and $\norm{x_n}_{\Gamma,\alpha,p} \leq 1$.
We can write $x_n = x_{n,0} + \lambda_n \Id_{\ell^2_I}$ with $x_{n,0} \in \SinftyI$ and $\lambda_n \in \C$.
Since $\SinftyI \oplus \C \Id_{\ell^2_I}$ is a closed subspace of $\B(\ell^2_I)$, we have $x = x_0 + \lambda \Id_{\ell^2_I}$ for some $x_0 \in \SinftyI$ and $\lambda \in \C$.
As the assignment $\SinftyI \oplus \C \Id_{\ell^2_I} \to \C, \: x_0 + \lambda \Id_{\ell^2_I} \mapsto \lambda$ is continuous, we have $\lambda_n \to \lambda$, and therefore also $x_{n,0} \to x_0$ in $S^\infty_I$.
Appealing to Proposition \ref{Prop-Leibniz-Schur-infinite}, it suffices to show that $x_{n,0}$ converges weakly to $x_0$ in $S^p_I$. Since the convergence already holds in $S^\infty_I$ norm, it suffices to show that $x_{n,0}$ is bounded in $S^p_I$. But this follows from $\bnorm{A_p^{\frac12}(x_{n,0})}_{S^p_I} \lesssim \norm{x_{n,0}}_{\Gamma,\alpha,p} = \norm{x_n}_{\Gamma,\alpha,p} \leq 1$ and the fact that $A_p^{-\frac12} \co \SpI \to \SinftyI$ is bounded according to Lemma \ref{lem-Schur-compact}.

1. (Leibniz rule) Note that if $I$ is finite, then the above proof of Leibniz quantum locally compact metric space shows that $\big(\big(\SinftyI\big)_{\sa} \oplus \R \Id_{\ell^2_I}, \norm{\cdot}_{\Gamma,\alpha,p}\big)$ is a Leibniz quantum compact metric space, in case $2 \leq p < \infty$. We indicate how the same proof also works for $p = \infty$. Note first that when $I$ is finite, we have $S^\infty_I = \M_{I,\fin}$, so that the domain of $\norm{\cdot}_{\Gamma,\alpha,\infty}$ is the full space $\SinftyI \oplus \C \Id_{\ell^2_I}$: as in the proof of Proposition \ref{Prop-Leibniz-Schur-infinite}, one can show that for any $x,y \in \M_{I,\fin}$, one has
\[ 
\norm{xy}_{\Gamma,\infty} 
\leq \norm{x}_{S^\infty_I} \norm{y}_{\Gamma,\alpha,\infty} + \norm{y}_{S^\infty_I} \norm{x}_{\Gamma,\alpha,\infty}. 
\]
Thus, in particular the same holds if $x,y \in \dom \norm{\cdot}_{\Gamma,\alpha,\infty} = \SinftyI \oplus \C \Id_{\ell^2_I}$. For the lower semi-continuity, if suffices to note that the seminorm satisfies $\norm{x}_{\Gamma,\alpha,\infty} \lesssim \norm{x}_{S^\infty_I}$, and by the reversed triangle inequality, we also have $| \norm{x}_{\Gamma,\alpha,\infty} - \norm{y}_{\Gamma,\alpha,\infty} | \leq \norm{x-y}_{\Gamma,\alpha,\infty}$.
\end{proof}

\subsection{Gaps of some markovian semigroups of Schur and Fourier multipliers}
\label{subsec-calculs-gap}

In this subsection, we study some typical examples of markovian semigroups of Schur and Fourier multipliers. We equally calculate their gaps \eqref{equ-Delta-A-prime-sgrp} and \eqref{equ-Delta-A-prime} and examine the injectivity of their Hilbert space representation. This information is important for applications to compact quantum metric spaces in Subsections \ref{subsec-new-compact-quantum-metric-spaces-Fourier} and \ref{Sec-I-finite}.

\paragraph{Heat Schur and Poisson Schur semigroup}
In the following, $I$ is equal to $\{1,\ldots,n\}$, $\N$ or $\Z$. We consider the heat Schur semigroup $(T_t)_{t \geq 0}$ acting on $\B(\ell^2_I)$ defined by
\begin{equation}
\label{Heat-Schur}
T_t \co [x_{ij}] \mapsto \big[\e^{-|i-j|^2 t} x_{ij}\big].
\end{equation}
Moreover, we also consider the Poisson Schur semigroup $(T_t)_{t \geq 0}$ acting on $\B(\ell^2_I)$ defined by
\begin{equation}
\label{Poisson-Schur}
T_t \co [x_{ij}] \mapsto \big[\e^{-|i-j|t} x_{ij}\big].	
\end{equation}
These two semigroups are examples of noncommutative diffusion semigroups consisting of Schur multipliers. 

Indeed, for the first one, we can take the real Hilbert space $H = \R$ and put $\alpha_i = i$ for $\alpha \in I$. So we have $\norm{\alpha_i - \alpha_j}_\R^2 = |i-j|^2$. For the second one, we can consider the real Hilbert space $H = \ell^2_I$ and put $\alpha_i = \sum_{k=0}^i e_k$ for $i \in I$ with $i \geq 0$ and $\alpha_i = \sum_{k= 0}^{|i|} e_{-k}$ for $i \in I$ with $i < 0$ (if $I$ contains negative elements). Then for $i > j \geq 0$, using the orthogonality of the $e_l$'s in the third equality, we have
$$
\norm{\alpha_i - \alpha_j}_{\ell^2_I}^2 
=\norm{\sum_{k=0}^i e_k -\sum_{k=0}^j e_k}^2_{\ell^2_I} 
=\norm{ \sum_{k=j+1}^i e_k}^2_{\ell^2_I} 
=\sum_{k=j+1}^i \norm{e_k}^2_{\ell^2_I} 
=\sum_{k=j+1}^i 1 
=i-j 
=|i-j|.
$$ 
In a similar way, we obtain for $i,j$ in general position that $\norm{\alpha_i - \alpha_j}_{\ell^2_I}^2 = |i-j|$. Note that both mappings $\alpha$ are injective.

\begin{lemma}
\label{lem-calcul-Gap-heat-Schur}
Consider the above heat Schur semigroup from \eqref{Heat-Schur}. We have $\G_{\alpha} = 1$. 
\end{lemma}

\begin{proof}
We have 
$$
\G_\alpha \ov{\eqref{equ-Delta-A-prime}}{=}  \inf_{\alpha_i - \alpha_j \neq \alpha_k - \alpha_l} \norm{(\alpha_i - \alpha_j) - (\alpha_k - \alpha_l)}_\R^2=\inf_{i - j \neq k - l} \norm{(i-j) - (k-l)}_\R^2
$$ 
which is clearly equal to 1.
%
\end{proof}

%


We will now calculate $\G_\alpha$ for the Poisson Schur semigroup above.
First we have

\begin{lemma}
\label{lem-Poisson-Schur-orthogonal-relations}
Consider the above Poisson Schur semigroup from \eqref{Poisson-Schur} with $I=\Z$. Then for any $i,j,k,l \in \Z$, $i \geq j \geq k \geq l$, we have
$$
\langle \alpha_i - \alpha_j , \alpha_k - \alpha_l \rangle_H 
=0.
$$
\end{lemma}

\begin{proof}
Indeed, we have
$$
\langle \alpha_i - \alpha_j , \alpha_k - \alpha_l \rangle_H 
=\bigg\langle \sum_{r=0}^i e_r -\sum_{r=0}^j e_r, \sum_{r=0}^k e_r -\sum_{r=0}^l e_r\bigg\rangle_H
=\bigg\langle \sum_{r=j+1}^i e_r, \sum_{r=l+1}^k e_r\bigg\rangle_H
=0
$$
\end{proof}

\begin{lemma}
\label{lem-calcul-Gap-Poisson-Schur}
Consider again the above Poisson Schur semigroup from \eqref{Poisson-Schur}.
We have $\G_\alpha = 1$.
\end{lemma}

\begin{proof}
It is clear that it suffices to examine the case $I=\Z$. 

For the inequality $\G_\alpha \leq 1$, it suffices to take $i = 1$, $j = k = l = 0$, in which case we have $\norm{(\alpha_i - \alpha_j) - (\alpha_k - \alpha_l)}_H^2 = \norm{\alpha_i - \alpha_j}_H^2 = |i-j| = 1$.

For the reverse inequality $\G_\alpha \geq 1$, consider any $i,j,k,l$ such that $\alpha_i - \alpha_j \neq \alpha_k - \alpha_l$. We want to estimate $\norm{(\alpha_i - \alpha_j) - (\alpha_k - \alpha_l)}^2$ from below. Using that $\norm{x} = \norm{-x}$ and exchanging names of indices, we can assume without loss of generality that $\max(i,j,k,l) = i$.

\textit{First case:} We have $l = \min(i,j,k,l)$.
Then exchanging the names of indices $j$ and $k$ if necessary, we have $i \geq j \geq k \geq l$.
Thus, according to Lemma \ref{lem-Poisson-Schur-orthogonal-relations}, we have
$ \langle \alpha_i - \alpha_j , \alpha_k - \alpha_l \rangle = 0$.
Consequently,
\begin{align*}
\MoveEqLeft
\norm{(\alpha_i - \alpha_j) - (\alpha_k - \alpha_l)}^2 = \norm{\alpha_i - \alpha_j}^2 + \norm{\alpha_k - \alpha_l}^2 - 2 \langle \alpha_i - \alpha_j , \alpha_k - \alpha_l \rangle \\
& = |i-j| + |k-l| - 0.
\end{align*}
Clearly, if this expression is $0$, then $i = j$ and $k = l$, which is excluded by $\alpha_i - \alpha_j \neq \alpha_k - \alpha_l$, or $\alpha_i - \alpha_k \neq \alpha_j - \alpha_l$ in case that we had exchanged names of indices.
In any other case, this expression is $\geq 1$ since $i,j,k,l$ take entire values.

\textit{Second case:} We have $\min(i,j,k,l) \in \{j,k\}$ and $l$ is the second smallest value among $i,j,k,l$.
Then exchanging the names of indices $j$ and $k$ if necessary, we can suppose $\min(i,j,k,l) = k$.
So we have $i \geq j \geq l \geq k$, and thus by Lemma \ref{lem-Poisson-Schur-orthogonal-relations}
$ \langle \alpha_i - \alpha_j , \alpha_l - \alpha_k \rangle = 0$.
We calculate
\begin{align*}
\MoveEqLeft
\norm{(\alpha_i - \alpha_j) - (\alpha_k - \alpha_l)}^2 = \norm{\alpha_i - \alpha_j}^2 + \norm{\alpha_k - \alpha_l}^2 + 2 \langle \alpha_i - \alpha_j, \alpha_l - \alpha_k \rangle \\
& = |i-j| + |k-l| + 0.
\end{align*}
We argue as before to see that this quantity is $\geq 1$.

\textit{Third case: }We have $\min(i,j,k,l) \in \{j,k\}$ and $l$ is the second biggest value among $i,j,k,l$. Then exchanging the names of indices $j$ and $k$ if necessary, we can suppose $\min(i,j,k,l) = k$. So we have $i \geq l \geq j \geq k$, and thus by Lemma \ref{lem-Poisson-Schur-orthogonal-relations}
\begin{align*}
\MoveEqLeft
\langle \alpha_i - \alpha_j , \alpha_k - \alpha_l \rangle = \langle \alpha_i - \alpha_l, \alpha_k - \alpha_l \rangle + \langle \alpha_l - \alpha_j, \alpha_k - \alpha_l \rangle \\
& = - 0 + \langle \alpha_l - \alpha_j, \alpha_k - \alpha_j \rangle + \langle \alpha_l - \alpha_j, \alpha_j - \alpha_l \rangle 
= - 0 - 0 - \norm{\alpha_l - \alpha_j}^2 = - |l-j|.
\end{align*}
Then we calculate
\begin{align*}
\MoveEqLeft
\norm{(\alpha_i - \alpha_j) - (\alpha_k - \alpha_l)}^2 = \norm{\alpha_i - \alpha_j}^2 + \norm{\alpha_k - \alpha_l}^2 - 2 \langle \alpha_i - \alpha_j, \alpha_k - \alpha_l \rangle \\
& = |i-j| + |k-l| +2 |l-j| 
 = i-j + l - k + 2 l - 2j = i-k + 3(l-j).
\end{align*}
Again we argue as before to see that this quantity is $\geq 1$.
\end{proof}

\paragraph{Markovian semigroups of Herz-Schur multipliers vs. markovian semigroups of Fourier multipliers.}

Let $G$ be a discrete group and $\psi \co G \to \R$ be a function. Suppose that $\psi(e)=0$. Recall that by \cite[Corollary C.4.19]{BHV}, the function $\psi$ is conditionally negative definite if and only if for any $t \geq 0$, the function $\e^{-t\psi}$ is of positive type. On the one hand, by \cite[Proposition 4.2]{DCH} that exactly means that $\psi$ induces a completely positive Fourier multiplier $T_t \ov{\mathrm{def}}{=} M_{\exp(-t\psi)} \co \VN(G) \to \VN(G)$ for any $t \geq 0$. On the other hand, by \cite[Definition C.4.1]{BHV} and \cite[Theorem D.3]{BrO} that is equivalent to say that $\psi$ induces a completely positive Herz-Schur multiplier $T_t^\HS\co \B(\ell^2_G) \to \B(\ell^2_G)$ for any $t \geq 0$ (whose symbol is $\big[\exp(-t\phi(s,r))\big]_{s,r \in G}$ where $\phi(s,r) \ov{\mathrm{def}}{=} \psi(s^{-1}r)$ for any $s,r \in G$). In this case, by Proposition \ref{prop-Schoenberg}, we obtain a markovian semigroup $(T_t)_{t \geq 0}$ of Fourier multipliers and it is easy to check and well-known that we obtain a markovian semigroup $(T_t^\HS)_{t \geq 0}$ of Herz-Schur multipliers. Hence there is a bijective correspondence between markovian semigroups of Fourier multipliers on $\VN(G)$ and markovian semigroups of Herz-Schur multipliers on $\B(\ell^2_G)$. So any triple $(b,\pi,H)$ associated to a markovian semigroup $(T_t)_{t \geq 0}$ of Fourier multipliers by Proposition \ref{prop-Schoenberg} gives a couple $(\alpha,H)$ associated to $(T_t^\HS)_{t \geq 0}$ by \cite[Proposition 5.4]{Arh1} and conversely. More precisely, if $\psi(s)=\norm{b_\psi(s)}_H^2$ then for any $s,r \in G$
\begin{align*}
\MoveEqLeft
\phi(s,r)
=\psi(s^{-1}r)
=\norm{b_\psi(s^{-1}r)}_H^2 
\ov{\eqref{Cocycle-law}}{=} \norm{b_\psi(s)+\pi_s(b_\psi(t))}^2_H \\
&\ov{\eqref{cocycle-b--1}}{=} \norm{-\pi_{s^{-1}}(b_\psi(s))+\pi_{s^{-1}}(b_\psi(r))}_H^2
=\norm{b_\psi(r)-b_\psi(s)}_H^2.           
\end{align*}
So we can consider the couple $(b_\psi,H)$ for the semigroup $(T_t^\HS)_{t \geq 0}$.

Next we compare the gaps of Herz-Schur and Fourier markovian semigroups as we encountered them in Subsections \ref{Sec-I-finite} and \ref{subsec-new-compact-quantum-metric-spaces-Fourier}. We will see in Proposition \ref{prop-donut-rational-rational} that a strict inequality may occur in the following result.

\begin{prop}
\label{prop-lien-Gap-Fourier-Herz-Schur}
Let $(T_t)_{t \geq 0}$ and $(T_t^\HS)_{t \geq 0}$ be markovian semigroups of Fourier multipliers and Herz-Schur multipliers as above. Consider a triple $(b,\pi,H)$ for the first semigroup and the couple $(\alpha,H)$ such that $\alpha \ov{\mathrm{def}}{=} b$ for the second semigroup. We have 
\begin{equation}
\label{Inequality-Gaps}
\G_\alpha 
\leq \G_\psi
\end{equation}
where $\G_\alpha$ and $\G_\psi$ are defined in \eqref{equ-Delta-A-prime} and \eqref{equ-Delta-A-prime-sgrp}.
\end{prop}

\begin{proof}
For any $i,j,l,k \in G$, we have $\alpha_i - \alpha_j =b(i)-b(j)\ov{\eqref{Cocycle-law}}{=} \pi_j(b(j^{-1}i))=\pi_j(\alpha_{j^{-1}i})$ and $\alpha_k - \alpha_l = \pi_l(\alpha_{l^{-1}k})$.
Thus, we have 
\begin{align*}
\MoveEqLeft
\G_\alpha 
\ov{\eqref{equ-Delta-A-prime}}{=} \inf_{\alpha_i - \alpha_j \neq \alpha_k - \alpha_l} \norm{(\alpha_i - \alpha_j) - (\alpha_k - \alpha_l)}_{H}^2
=\inf_{\alpha_i - \alpha_j \neq \alpha_k - \alpha_l} \norm{\pi_j(\alpha_{j^{-1}i}) - \pi_l(\alpha_{l^{-1}k}) }_H^2 \\ 
&=\inf_{i,j,k,l} \norm{\pi_{l^{-1}j}(\alpha_{j^{-1}i}) - \alpha_{l^{-1}k} }_H^2 
=\inf_{j,s,r} \norm{\pi_j(\alpha_s) - \alpha_r }_H^2,
\end{align*}
where the infimum is taken over those $j,s,r \in G$ such that the considered norm is $\neq 0$. On the other hand, we have, since $b(s) = \alpha_s$ and by considering $j = e$ and $\pi_e = \Id_H$ in the inequality

$$
\G_\psi 
\ov{\eqref{equ-Delta-A-prime-sgrp}}{=} \inf_{b(s) \neq b(t)} \norm{b(s) - b(t)}_H^2 
=\inf_{\alpha_s \neq \alpha_r} \norm{\alpha_s - \alpha_r}_H^2 \geq \inf_{j,s,r} \norm{\pi_j(\alpha_s) - \alpha_r}_H^2 = \G_\alpha ,
$$
(again infimum over non-zero quantities). 
\end{proof}

\paragraph{Finite-dimensional Hilbert spaces and coboundary cocycles over infinite groups: a dichotomy between non-injective cocycles and the condition $\G_\psi = 0$.}

Now, we consider particular markovian semigroups of Fourier multipliers and their corresponding semigroups of Herz-Schur multipliers  and ask whether the representations $b =\alpha$ are injective and calculate the gaps $\G_\psi$ and $\G_\alpha$ defined in \eqref{equ-Delta-A-prime} and \eqref{equ-Delta-A-prime-sgrp}. We start with a general observation. Consider an infinite discrete group $G$ and a finite-dimensional orthogonal representation $\pi \co G \to \mathrm{O}(H)$.
Consider a 1-cocycle $b \co G \to H$ with respect to $\pi$. So we have a markovian semigroup of Fourier multipliers on $\VN(G)$ and a semigroup of Herz-Schur multipliers on $\B(\ell^2_G)$.

For 1-coboundaries, the situation is not as nice as in Subsection
\ref{subsec-new-compact-quantum-metric-spaces-Fourier}. Indeed, we have the following proposition.

\begin{prop}
\label{prop-dichotomy-injective-gap}
Let $G$ be an infinite discrete group and $b \co G \to H$ be a 1-coboundary with respect to some finite-dimensional orthogonal representation $\pi \co G \to \mathrm{O}(H)$. If $b$ is injective, then $\G_\psi = \G_\alpha = 0$. So if $\G_\psi > 0$ or $\G_\alpha > 0$ then $b$ is non-injective.
\end{prop}

\begin{proof}
By definition \cite[Definition 2.2.3]{BHV}, there exists $\xi \in H$ such that $b(s) = \pi_s(\xi) - \xi$ for any $s \in G$. Clearly, the formula $b(s) = \pi_s(\xi) - \xi$ implies that $\xi \neq 0$ in the case where $b$ is injective. For simplicity, we assume that $\norm{\xi}_H = 1$. Since $b$ is injective, $O = \{\pi_s(\xi) : \: s \in G \}$ is an infinite subset of the sphere $S$ of $H$. Since $H$ is finite-dimensional, $S$ is compact. So there exists some accumulation point $\eta$ in the sphere of the orbit $O$. We have thus a convergent sequence $(\pi_{s_n}(\xi))$ consisting by injectivity of $b$ of different points. We infer that 
$$
0 
\leq \G_\psi 
\ov{\eqref{equ-Delta-A-prime-sgrp}}{\leq} \lim_{n \to \infty} \norm{b(s_n)-b(s_{n+1})}^2 
=\lim_{n \to \infty} \norm{\pi_{s_{n+1}}(\xi) - \pi_{s_n}(\xi)}^2 
=0.
$$
Thus $0 \leq \G_\alpha \ov{\eqref{Inequality-Gaps}}{\leq} \G_\psi = 0$.
\end{proof}

\begin{remark} \normalfont
In the paper \cite[p.~1967]{JMP1}, the authors were able to find the original paper which contains the famous Bieberbach Theorem. Unfortunately, this theorem is badly written in \cite[p.~1967]{JMP1} (since a crucial assumption is missing; to compare with a textbook, e.g. \cite[Th 7.2.4 p.~306]{Rat1}). So the proof of \cite[Theorem 6.4]{JMP1} is doubtful since this ``more general version'' of Bieberbach Theorem is used in the proof. Nevertheless, we think that an additional ``proper'' assumption could lead to a correct statement. 
\end{remark}

Recall that a topological group has property (FH) if every affine isometric action of $G$ on a real Hilbert space has a fixed point, see \cite[Definition 2.1.4]{BHV}. Note that by \cite[Theorem 2.12.4]{BHV}, if a discrete group $G$ has (T) then $G$ has (FH) and the converse is true if $G$ is countable. A finite group has (FH). The groups $\Z^k$ and free groups $\F_k$ do not have (FH) if $k \geq 1$. Using \cite[Proposition 2.2.10]{BHV}, we deduce the following result. In \cite[p.~1968]{JMP1}, it is written that ``infinite groups satisfying Kazhdan property (T) do not admit finite-dimensional standard cocycles'' (i.e. injective, finite-dimensional with $\G_\psi >0$). But from our point of view, the proof is missing. Proposition \ref{prop-dichotomy-injective-gap} allows us to give a proof and to obtain a slightly more general version.

\begin{cor}
\label{Cor-Kaz}
Let $G$ be an infinite discrete group with property (FH) and $b \co G \to H$ be a 1-cocycle with respect to some finite-dimensional orthogonal representation $\pi \co G \to \mathrm{O}(H)$. If $b$ is injective, then $\G_\psi = \G_\alpha = 0$. 
\end{cor}

\paragraph{Semigroups on finite groups}
Let $G$ be a finite group.  By \cite[pp.~1970-1971]{JMP1}, there always exists an orthogonal representation $\pi \co G \to H$ on some finite-dimensional real Hilbert space $H$ and an injective 1-cocycle $b \co G \to H$ with respect to $\pi$. In this case, since $G$ is finite,  $b$ only takes a finite number of values. This implies that $\G_\psi \ov{\eqref{equ-Delta-A-prime-sgrp}}{=} \inf_{b(s) \neq b(t)} \norm{b(s) - b(t)}_{H}^2 > 0$. For example, we can consider the left regular representation $\pi \co G \to \B(\ell^2_G)$ defined by $\pi_s(e_t)=e_{st}$ for any $s,t \in G$ and the cocycle $b \co G \to \ell^2_G$, $s \mapsto \pi_s(\xi)-\xi$ where $\xi$ is some vector of $\ell^2_G$ satisfying $\pi_s(\xi) \not=\xi$ for any element $s$ of $G-\{e_G\}$. We refer to \cite[a) and b) p.~1971]{JMP1} for other interesting examples of 1-cocycles. Note in addition that in the context of Schur multipliers, we also have $\G_\alpha>0$ if $\alpha=b$.

\paragraph{Heat semigroup on $\T^n$}
Here $G=\Z^n$. We consider the Heat semigroup $(T_t)_{t \geq 0}$ on $\L^\infty(\T^n)=\VN(\Z^n)$ defined by $T_t \co \L^\infty(\T^n) \to \L^\infty(\T^n)$, $\e^{\i k \cdot } \mapsto \e^{-t|k|^2} \e^{\i k \cdot}$. The associated finite-dimensional injective cocycle is given by the canonical inclusion $b \co \Z^n \to \R^n$ equipped with the trivial action $\pi_k = \Id_{\R^n}$ for all $k \in \Z^n$. 

\begin{lemma}
\label{lem-calcul-Gap-heat-torus}
Consider the above heat semigroup. We have $\G_{\psi} = 1$.
\end{lemma}

\begin{proof}
Indeed, we have $\G_\psi \ov{\eqref{equ-Delta-A-prime-sgrp}}{=} \inf_{b_\psi(k) \neq b_\psi(l)} \norm{b_\psi(k) - b_\psi(l)}_{\R^n}^2=\inf_{k \neq l, k,l \in \Z^n} |k-l|^2=1$.
\end{proof}

\paragraph{The donut type markovian semigroup of Fourier and Herz-Schur multipliers.}

Consider now the example of donut type Fourier multipliers in the spirit of \cite[Section 5.3]{JMP1}. That is, we consider the group $G = \Z$ with cocycle
\begin{equation}
\label{Cocycle-donut}
b(n) 
\ov{\mathrm{def}}{=} \big(\e^{2 \pi \i \alpha n},\e^{2 \pi \i \beta n}\big) - (1,1) 
\in \C^2 
= \R^4,
\end{equation}
where $\alpha, \beta \in \R$ and we consider $\C^2$ as the real Hilbert space $\R^4$. The associated cocycle orthogonal representation is\footnote{\thefootnote. For any $n,m \in \Z$, we have the cocycle law
\begin{align*}
\MoveEqLeft
b(n)+\pi_n(b_m) 
=\big(\e^{2 \pi \i \alpha n},\e^{2 \pi \i \beta n}\big) - (1,1) +\pi_n\big[\big(\e^{2 \pi \i \alpha m},\e^{2 \pi \i \beta m}\big) - (1,1) \big] \\          
&=\big(\e^{2 \pi \i \alpha n},\e^{2 \pi \i \beta n}\big) - (1,1) + \big(\e^{2 \pi \i \alpha (n+m)},\e^{2 \pi \i \beta (n+m)}\big) -\big(\e^{2 \pi \i \alpha n},\e^{2 \pi \i \beta n}\big) \\
&=\big(\e^{2 \pi \i \alpha (n+m)},\e^{2 \pi \i \beta (n+m)}\big) - (1,1) 
=b(n+m).
\end{align*}}
\begin{equation}
\label{Donut-cocycle-action}
\pi_n(x,y) 
\ov{\mathrm{def}}{=}\big(\e^{2 \pi \i \alpha n}x,\e^{2\pi \i \beta n}y\big) 
\qquad n \in \Z, \: x,y \in \C.
\end{equation}

\begin{prop}
\label{prop-donut-rational-rational}
Consider the case that both $\alpha$ and $\beta$ take rational values. Then $b \co \Z \to \C^2$ is not injective, but $\G_\psi, \G_\alpha > 0$. Moreover, the strict inequality $\G_\alpha < \G_\psi$ may happen.
\end{prop}

\begin{proof}
Consider $p,q \in \Z$ and $N \in \N^*$ such that $\alpha = \frac{p}{N}$ and $\beta = \frac{q}{N}$. 
We have 
$$
b(n+N) 
=\big(\e^{2 \pi \i (p + \alpha n)},\e^{2 \pi \i (q + \beta n)}\big) - (1,1) 
=\big(\e^{2 \pi \i \alpha n}, \e^{2 \pi \i \beta n}\big) - (1,1) 
=b(n),
$$ 
so that $b$ is $N$-periodic, and hence only takes a finite number of values. In particular, the function $b \co \Z \to \C^2$ is not injective and the set $\{b(n) - b(m) : \: n,m \in \Z\}$ is finite, which readily implies that $\G_\psi = \inf_{b(n) \neq b(m)} \norm{b(n) - b(m)}_{\R^4}^2 > 0$. In the same manner, we obtain $\G_\alpha = \inf_{b(i) - b(j) \neq b(k) - b(l)} \norm{b(i) - b(j) - (b(k) - b(l))}_{\R^4}^2 > 0$.

We turn to the statement of strict inequality. To this end, we take $\alpha = \beta = \frac18$.
Then 
\begin{align*}
\MoveEqLeft
\G_\psi 
= 2 \inf \norm{(\e^{2 \pi \i \alpha n} - 1) - ( \e^{2 \pi \i \alpha m} -1 )}_\C^2 = 2 \inf \norm{\e^{2 \pi \i \alpha (n-m)} - 1}_\C^2 \\
& = 2 \left|\e^{2 \pi \i \frac18} - 1\right|^2 
= 2 \cdot \left| \frac{1 + \i}{\sqrt{2}} - 1 \right|^2 
= 2 \cdot \left( \left(\frac{1}{\sqrt{2}}-1\right)^2 + \frac{1}{(\sqrt{2})^2}\right) 
= 4 \left(1 - \frac{1}{\sqrt{2}} \right).
\end{align*}
On the other hand, according to the proof of Proposition \ref{prop-lien-Gap-Fourier-Herz-Schur}, we have
\begin{align*}
\MoveEqLeft
\G_\alpha 
= 2 \inf \norm{\e^{2 \pi \i \alpha n}(\e^{2 \pi \i \alpha r} - 1) - (\e^{2 \pi \i \alpha s} - 1) }_{\C}^2 
 \leq 2 \norm{\e^{2 \pi \i \alpha \cdot 1}(\e^{2 \pi \i \alpha \cdot 2} - 1) - (\e^{2 \pi \i \alpha \cdot 4} - 1)}_\C^2 \\
&=2 \norm{\frac{1 + \i}{\sqrt{2}}(\i-1) - (-1 - 1)}_\C^2 
=2 \left|\frac{1}{\sqrt{2}}(\i-1-1-\i) +2 \right|^2 
=2 \left|2 \left( 1 - \frac{1}{\sqrt{2}} \right) \right|^2 \\
&=8 \left( 1 - \frac{1}{\sqrt{2}} \right)^2 
<4 \left(1 - \frac{1}{\sqrt{2}} \right) 
=\G_\psi,
\end{align*}
since $1 - \frac{1}{\sqrt{2}} \in \left(0,\frac12\right)$.
\end{proof}

\begin{lemma}
\label{lem-donut-irrational-rational}
Suppose that at least one of the numbers $\alpha$ and $\beta$ is irrational. Then the cocycle $b$ is injective, and $\G_\alpha = \G_\psi = 0$.
\end{lemma}

\begin{proof}
If, say $\alpha$, is irrational, then $\e^{2\pi \i \alpha n} = 1$ for some $n \in \Z$ implies that $\alpha n \in \Z$, so $n = 0$. Then $b(n) = b(m)$ for some $n, m \in \Z$ implies that $\e^{2 \pi \i \alpha n} -1 = \e^{2 \pi \i \alpha m} - 1$, so $\e^{2 \pi \i \alpha (n-m)} = 1$, and consequently, $n - m = 0$. We infer that $b$ is injective. Note that $G = \Z$ is infinite and $H = \R^4$ is finite-dimensional. Furthermore, since $b(n) \ov{\eqref{Cocycle-donut}}{=} \big(\e^{2 \pi \i \alpha n},\e^{2 \pi \i \beta n}\big) - (1,1) \ov{\eqref{Donut-cocycle-action}}{=} \pi_n(1,1) - (1,1)$, the cocycle $b$ is a coboundary. Thus the assumptions of Proposition \ref{prop-dichotomy-injective-gap} are fulfilled. This proposition implies that $\G_\alpha = \G_\psi = 0$.
\end{proof}

\paragraph{The free group $\F_2$.}

Consider now the example of the action of the free group $\F_2$ with two generators $a_1$ and $a_2$ on $\R^3$ from \cite[Section 5.5]{JMP1}. That is, we consider the representation $\pi \co \F_2 \to \mathrm{O}(\R^3)$ from the Banach-Tarski paradox. Take an angle $\theta \in \R \backslash 2 \pi \Q$ and define $\pi$ uniquely by putting 
$$ 
\pi_{a_1} 
= \begin{bmatrix} \cos \theta & - \sin \theta & 0 \\ \sin \theta & \cos \theta & 0 \\ 0 & 0 & 1 \end{bmatrix} 
\quad \text{and} \quad 
\pi_{a_2} 
= \begin{bmatrix} 1 & 0 & 0 \\ 0 & \cos \theta & - \sin \theta \\ 0 & \sin \theta & \cos \theta \end{bmatrix}.
$$
Then it is known that $\pi \co \F_2 \to \mathrm{O}(\R^3)$ is injective. Further it is shown in \cite[Section 5.5]{JMP1} that there exists some $\xi \in \R^3 \backslash \{ 0 \}$ such that $b \co \F_2 \to \R^3$, $b(s) = \pi_s(\xi) - \xi$ is an injective coboundary cocycle. We appeal again to Proposition \ref{prop-dichotomy-injective-gap} to deduce that in this case $\G_\alpha = \G_\psi = 0$.

\subsection{Banach spectral triples}
\label{Sec-Banach-spectral-triples-first}


In the remainder of Section \ref{new-quantum-metric-spaces}, we will establish that our Hodge-Dirac operators defined in \eqref{equ-full-Hodge-Dirac-operator-Fourier} and in \eqref{equ-full-Hodge-Dirac-operator} associated with markovian semigroups of Fourier multipliers or Schur multipliers give rise to  (locally) compact Banach spectral triples, see Proposition \ref{First-spectral-triple}, Proposition \ref{prop-locally-compact-spectral-triple-Schur-I} and Theorem \ref{First-spectral-triple-2}. We will also establish such Banach spectral triples for two other Hodge-Dirac operators, see Theorem \ref{thm-spectral-triple} and Theorem \ref{thm-spectral-triple-group}. We refer to \cite{CGIS1}, \cite{CGRS1}, \cite{Con3}, \cite{GVF1} and \cite{Var1} for more information on spectral triples. Let us recall this notion. A (possibly kernel-degenerate, compact) spectral triple $(A,H,D)$ consists of a unital $\mathrm{C}^*$-algebra $A$, a Hilbert space $H$, a (densely defined, unbounded) selfadjoint operator $D$ and a representation $\pi \co A \to \B(H)$ which satisfy the following properties.
\begin{enumerate}
	\item $D^{-1}$ is compact on $\ovl{\Ran D}$.
	\item The set 
\begin{align*}
\MoveEqLeft
\Lip_D(A) 
\ov{\mathrm{def}}{=}\big\{ a \in A : \pi(a) \cdot \dom D \subseteq \dom D 
\text{ and the unbounded operator } \\
&\qquad \qquad \qquad [D,\pi(a)] \co \dom D \subseteq H \to H \text{ extends to an element of } \B(H)\big\} \nonumber
\end{align*} 
is dense in $A$.
\end{enumerate}

In the next subsections, we give new examples of spectral triples. Our examples can be generalized to the context of $\L^p$-spaces instead of Hilbert spaces. So, it is natural to state the following definition.

\begin{defi}
\label{Def-Banach-spectral-triple}
A (compact) Banach spectral triple $(A,X,D)$ consists of the following data: a reflexive\footnote{\thefootnote. It may perhaps be possible to replace the reflexivity by an assumption of weakly compactness, see \cite[p. 361]{HvNVW2}.} Banach space $X$, a bisectorial operator $D$ on $X$ with dense domain $\dom D \subseteq X$, a Banach algebra $A$ and a homomorphism $\pi \co A \to \B(X)$ such that for all $a \in A$ we have:
\begin{enumerate}
\item{} $D$ admits a bounded $\HI$ functional calculus on a bisector $\Sigma^\pm_\omega$.
\item{} $|D|^{-1}$ is a compact operator on $\ovl{\Ran D}$.
\item{} The set 
\begin{align}
\label{Lipschitz-algebra-def}
\MoveEqLeft
\Lip_D(A) 
\ov{\mathrm{def}}{=} \big\{a \in A : \pi(a) \cdot \dom D \subseteq \dom D 
\text{ and the unbounded operator } \\
&\qquad \qquad \qquad [D,\pi(a)] \co \dom D \subseteq X \to X \text{ extends to an element of } \B(X)\big\} \nonumber
\end{align} 
is dense in $A$.
\end{enumerate}
\end{defi}

In these conditions, $D$ is closed since any operator with non-empty resolvent set is closed. Moreover, by \cite[p.~448]{HvNVW2}, we have a direct sum decomposition $X=\ovl{\Ran D} \oplus \ker D$.
The operator $|D|^{-1}$ is defined through bisectorial functional calculus as $|D|^{-1} = m(D)$, where $m(\lambda) = \frac{1}{\sqrt{\lambda^2}}$ is a holomorphic function on any bisector of angle $< \frac{\pi}{2}$.
Finally, note that a spectral triple is a Banach spectral triple. Note that we can define $\Lip_D(A)$ in the case of an unbounded operator $D$ with dense domain. 

In the sequel, we will use sometimes the notation $a$ for $\pi(a)$. If $a \in \Lip_D(A)$, we let
\begin{equation}
\label{Lip-norm-Dirac}
\norm{a}_{D}
\ov{\mathrm{def}}{=} \bnorm{[D,a]}_{X \to X}.
\end{equation}

In the following, at several places we will use ideas of the proofs of \cite[Proposition 3.7]{Rie4} and \cite[Proposition 1.6]{Lat8}. We refer to \cite{Wea1} for related things. The point 1 says that the map $\partial_D \co \Lip_D(A) \to \B(X)$, $a \mapsto [D,a]$ is a derivation. 

\begin{prop}
\label{Prop-Lip-algebra}
Let $X$ be a reflexive Banach space. Consider a closed linear operator $D$ on $X$ with dense domain $\dom D \subseteq X$, a Banach algebra $A$ and a homomorphism $\pi \co A \to \B(X)$.
\begin{enumerate}
	\item The space $\Lip_D(A)$ is a subalgebra of $A$. Moreover, if $a,b \in \Lip_D(A)$, we have
\[
[D,ab]
=a[D,b] + [D,a]b
\quad \text{and} \quad
\bnorm{ab}_{D} 
\leq \norm{a}_{D}  \norm{\pi(b)}_{\B(X)} + \norm{\pi(a)}_{\B(X)}\norm{b}_{D}.
\]
\item $\norm{\cdot}_{D}$ is a seminorm on $\Lip_D(A)$. 

\item For any $a \in \Lip_D(A)$ we have $\pi(a)^* \cdot \dom D^* \subseteq \dom D^*$ and the linear operator 
$$
[D^*,\pi(a)^*] \co \dom D^* \subseteq X^* \to X^*
$$ 
extends to a bounded operator on $X^*$ denoted with the same notation and we have
\begin{equation}
\label{Adjoint-commutateur}
[D,\pi(a)]^*
=-[D^*,\pi(a)^*].
\end{equation} 

\item Suppose that $\pi \co A \to \B(X)$ is continuous when $A$ is equipped with the weak topology and when $\B(X)$ is equipped with the weak operator topology. Then $\norm{\cdot}_{D}$ is lower semicontinuous on $\Lip_D(A)$ when $\Lip_D(A)$ is equipped with the induced topology by the weak topology of $A$.
\end{enumerate}
\end{prop}

\begin{proof}
1. If $a,b \in \Lip_D(A)$ then $ab\cdot\dom D=a\cdot(b \cdot \dom D) \subseteq a\cdot\dom D \subseteq \dom D$. Moreover, if $\xi \in \dom{D}$ then:
\begin{align*}
Dab\xi - abD \xi 
= D ab \xi - aDb\xi+aDb\xi - abD\xi  
= [D,a] b\xi+ a[D,b]\xi.
\end{align*}
Thus, as operators on $\dom D$, we conclude that $[D,ab]=a[D,b] + [D,a]b$. So $ab \in \Lip_D(A)$. Moreover, we have
\begin{align*}
\norm{ab}_{D} 
&\ov{\eqref{Lip-norm-Dirac}}{=}\bnorm{[D,ab]}_{X \to X} 
=\bnorm{[D,a]b + a[D,b]}_{X \to X} \\
&\leq \bnorm{[D,a]}_{X \to X} \norm{\pi(b)}_{\B(X)} + \norm{\pi(a)}_{\B(X)}\bnorm{[D,b]}_{X \to X} \\
&=\norm{a}_{D} \norm{\pi(b)}_{\B(X)} + \norm{\pi(a)}_{\B(X)}\norm{b}_{D}.
\end{align*}

2. If $\lambda \in \C$ and $a \in \Lip_D(A)$, we have $\norm{\lambda a}_{D} \ov{\eqref{Lip-norm-Dirac}}{=} \bnorm{[D,\lambda a]}_{X \to X}=\bnorm{\lambda[D, a]}_{X \to X}=|\lambda|\bnorm{[D, a]}_{X \to X}=|\lambda|\bnorm{a}_{D}$. If $a,b \in \Lip_D(A)$, we have
$\norm{a+b}_{D} \ov{\eqref{Lip-norm-Dirac}}{=} \bnorm{[D,a+b]}_{X \to X}=\bnorm{[D,a]+[D,b]}_{X \to X} \leq \bnorm{[D,a]}_{X \to X}+\bnorm{[D,b]}_{X \to X}
\ov{\eqref{Lip-norm-Dirac}}{=}\norm{a}_{D}+\norm{b}_{D}$. 

3. Let $a \in \Lip_D(A)$. If $\xi \in \dom D$ and $\zeta \in \dom D^*$ then:
\begin{align*}
\MoveEqLeft
\big\langle \pi(a)^*\zeta,D(\xi) \big\rangle_{X^*,X}
=\big\langle \zeta,\pi(a)D(\xi) \big\rangle_{X^*,X} 
=\big\langle \zeta,D\pi(a)\xi \big\rangle_{X^*,X} - \big\langle \zeta,[D,\pi(a)]\xi \big\rangle_{X^*,X}\\
&=\big\langle D^*(\zeta),\pi(a)\xi \big\rangle_{X^*,X} - \big\langle \zeta,[D,\pi(a)]\xi \big\rangle_{X^*,X}.
\end{align*} 
Now, the linear map $\dom D \to \C,\: \xi \mapsto \big\langle D^*(\zeta), \pi(a)\xi \big\rangle_{X^*,X}$ is continuous. Since $[D,\pi(a)]$ is bounded, the linear map $\dom D \to \C$, $\xi \mapsto \big\langle \zeta, [D,\pi(a)]\xi \big\rangle_{X^*,X}$ is also continuous. Hence $\dom D \to \C$, $\xi \mapsto \big\langle \pi(a)^*\zeta, D(\xi) \big\rangle_{X^*,X}$ is continuous. Hence $\pi(a)^* \zeta$ belongs to $\dom D^{*}$. 

For any $a \in \Lip_D(A)$, any $\xi \in \dom D$ and any $\zeta \in \dom D^* $, we have
\begin{align*}
\MoveEqLeft
\big\langle [D,\pi(a)]\xi,\zeta \big\rangle_{X,X^*}            
=\big\langle (D\pi(a)-\pi(a)D)\xi, \zeta\big\rangle_{X,X^*} 
=\big\langle D\pi(a)\xi, \zeta\big\rangle_{X,X^*}-\big\langle \pi(a) D\xi, \zeta\big\rangle_{X,X^*} \\
&=\big\langle \pi(a)\xi, D^*\zeta\big\rangle_{X,X^*}-\big\langle D\xi, \pi(a)^*\zeta \big\rangle_{X,X^*}
=\big\langle \xi, \pi(a)^*D^*\zeta \big\rangle_{X,X^*}-\big\langle \xi, D^*\pi(a)^*\zeta\big\rangle_{X,X^*} \\
&=\big\langle \xi,-[D^*,\pi(a)^*]\zeta \big\rangle_{X,X^*}.
\end{align*}
Hence the operators $[D,\pi(a)]$ and $-[D^*,\pi(a)^*]$ are formal adjoint to each other in the sense of \cite[p.~167]{Kat1}. We infer that $-[D^*,\pi(a)^*] \subseteq [D,\pi(a)]^*$. Since this latter operator is bounded and since the domain $\dom D^*$ of $-[D^*,\pi(a)^*]$ is dense by \cite[Theorem 5.29 p.~168]{Kat1}, we obtain the result.  


4. Suppose that $\pi \co A \to \B(X)$ is continuous when $A$ is equipped with the weak topology and when $\B(X)$ is equipped with the weak operator topology.
We shall show that $\norm{\cdot}_{D}$ is lower semicontinuous on $\Lip_D(A)$ when $\Lip_D(A)$ is equipped with the induced topology by the weak topology of $A$.
Let $a \in A$ and $(a_j)$ be a net of elements of $\Lip_D(A)$ such that $(a_j)$ converges weakly to $a$ and $\norm{a_j}_{D} \leq 1$. We have to show that $a$ belongs to $\Lip_D(A)$ and that $\norm{a}_{D} \leq 1$. Note that 
\begin{equation}
\label{Divers-0978}
\bnorm{[D^*,\pi(a_j)^*]}_{X^* \to X^*}
\ov{\eqref{Adjoint-commutateur}}{=}\bnorm{[D,\pi(a_j)]}_{X \to X} \leq 1
\end{equation}
and that $\pi(a_j) \to \pi(a)$ converges for the weak operator topology. Let $\xi \in \dom D$ and let $\zeta \in \dom D^*$. Note that $\pi(a_j)^*(\zeta)$ belongs to $\dom D^*$. Then for any integer $n$
\begin{align*}
\MoveEqLeft
\big\langle \pi(a_j)\xi, D^*(\zeta) \big\rangle_{X,X^*}
=\big\langle \xi,\pi(a_j)^*D^*(\zeta) \big\rangle_{X,X^*} 
=\big\langle \xi,D^* \pi(a_j)^*(\zeta) \big\rangle_{X,X^*} - \big\langle \xi,[D^*, \pi(a_j)^*]\zeta \big\rangle_{X,X^*} \\
&=\big\langle D(\xi),\pi(a_j)^*\zeta \big\rangle_{X,X^*} - \big\langle \xi,[D^*,\pi(a_j)^*]\zeta \big\rangle_{X,X^*}.
\end{align*}
Passing to the limit and using the point 3, we deduce that
\begin{align*}
\MoveEqLeft
\left|\big\langle \pi(a) \xi,D^*(\zeta) \big\rangle_{X,X^*}\right| 
=\lim_{j} \left|\big\langle \pi(a_j)\xi, D^*(\zeta) \big\rangle_{X,X^*}\right| \\
&\leq \liminf_{j} \big[ \left| \langle D\xi,\pi(a_j)^*\zeta \rangle_{X,X^*}\right| + \left|\langle \xi, [D^*,\pi(a_j)^*]\zeta \rangle_{X,X^*} \right|\big]\\
&\ov{\eqref{Divers-0978}}{\leq} \left| \langle D\xi,\pi(a)^*\zeta \rangle_{X,X^*} \right| + \norm{\xi}_X \norm{\zeta}_{X^*} 
\leq \norm{\zeta}_{X^*} \big( \norm{D(\xi)}_X \norm{\pi(a)}_A + \norm{\xi}_{X} \big).
\end{align*} 
So the function $\dom D^* \to \mathbb{C}$, $\zeta \mapsto \big\langle \pi(a) \xi, D^*(\zeta) \big\rangle_{X,X^*}$ is continuous, and thus $\pi(a) \xi \in \dom D^{**}=\dom D$ by \cite[Theorem 5.29]{Kat1}. We conclude that $\pi(a) \cdot \dom D \subseteq \dom D$. If $\xi \in \dom D$ and $\zeta \in \dom D^*$ with $\norm{\xi}_X \leq 1$ and $\norm{\zeta}_{X^*} \leq 1$ then
\begin{align*}
\MoveEqLeft
\left| \big\langle \pi(a_j)\xi, D^*(\zeta) \big\rangle_{X,X^*} - \big\langle D(\xi),\pi(a_j)^*\zeta \big\rangle_{X,X^*} \right|
= \left|\big\langle D\pi(a_j)(\xi),\zeta \big\rangle_{X,X^*}-\langle \pi(a_j)D(\xi),\zeta \rangle_{X,X^*}\right| \\
&= \left| \big\langle (D \pi(a_j)-\pi(a_j)D)(\xi),\zeta \big\rangle_{X,X^*} \right|
= \left| \big\langle[D,\pi(a_j)](\xi),\zeta \big\rangle_{X,X^*} \right| \\
&\leq \bnorm{[D,\pi(a_j)]}_{X \to X} \norm{\xi}_X\norm{\zeta}_{X^*} 
\leq 1.
\end{align*}
Passing to the limit, we obtain
\begin{align*}
\MoveEqLeft
\left| \big\langle [D,\pi(a)]\xi,\zeta \big\rangle_{X,X^*} \right|
= \left| \big\langle (D\pi(a)-\pi(a)D)\xi,\zeta \big\rangle_{X,X^*} \right| \\
&= \left| \langle \pi(a)\xi,D^*\zeta \rangle_{X,X^*}- \langle D\xi,\pi(a)^*\zeta \rangle_{X,X^*} \right|
\leq 1.
\end{align*}
Hence $[D,\pi(a)]$ is bounded on $\dom{D}$ and thus extends to $X$ to a bounded operator and $\norm{[D,\pi(a)]} \leq 1$. 
\end{proof}

\begin{remark} \normalfont
\label{Natural-remark}
In the particular case of a spectral triple, we have $D=D^*$ and $\pi(a)^*=\pi(a^*)$. So the above says that $\Lip_D(A)$ is a $*$-subalgebra of $A$ (which is well-known folklore).
\end{remark}

\begin{remark} \normalfont
\label{Remark-123}
In the point 4, we can replace the weak topology on $A$ by a topology such that $\pi \co A \to \B(X)$ remains continuous when $\B(X)$ is equipped with the weak operator topology.
\end{remark}

The following is a Banach space generalization of \cite[Proposition 2.1]{FMR1} proved with a similar method. 

\begin{prop}
\label{Prop-magic-core}
Let $X$ be a Banach space and $D \co \dom D \subseteq X \to X$ be a closed operator. 
\begin{enumerate}
	\item Suppose that $C$ is a core of $D$ and that $a \in \B(X)$ satisfies
\begin{itemize}
\item[(a)] $a \cdot C \subseteq \dom D$
\item[(b)] $[D,a]|C \co C \to X$ is bounded on $C$.
\end{itemize}
Then $a \cdot \dom D \subseteq \dom D$ and the operator $[D,a] \co \dom D \to X$ is well-defined. 

\item If in addition $D$ is densely defined and if there exists an adjoint\footnote{\thefootnote. In the sense of \cite[p.~167]{Kat1}.} $T \co \dom T \subseteq X^* \to X^*$ of $D$, and a subspace $Y$ of $\dom T$ which is dense in $X^*$ such that $a^*\cdot Y \subseteq \dom T$, then $[D,a] \co \dom D \to X$ extends to a bounded operator on $X$.
\end{enumerate}
\end{prop}

\begin{proof}
1. Let $x \in \dom D$. By \eqref{Def-core} there exists a sequence $(x_n)$ of elements of $C$ such that $x_n \to x$ and $D(x_n) \to D (x)$ in $X$. Since $a \in \B(X)$, we have $a(x_n) \to a(x)$ in $X$. Moreover, for any integer $n,m$, we have
\begin{align*}
\bnorm{Da(x_n)-Da(x_m)}_X
&=\bnorm{aD(x_n)- a D(x_m)+Da(x_n)-aD(x_n)-Da(x_m)+a D(x_m)}_X\\
&=\bnorm{a D(x_n) - a D(x_m)+ [D,a](x_n)-[D,a](x_m)}_X \\
&\leq \norm{a}_{X \to X} \bnorm{D(x_n)-D(x_m)}_X+\bnorm{[D,a]}_{C \to X} \norm{x_n-x_m}_X.
\end{align*}
We infer that $(Da(x_n))$ is a Cauchy sequence in $X$ hence converges. Since $D$ is closed, we conclude that $a(x)$ belongs to $\dom D$.

2. Let $x \in \dom D$ and $y \in Y$. Then
\begin{align*}
\MoveEqLeft
\big\langle [D,a](x),y \big\rangle_{X,X^*}
=\big\langle (Da-aD)(x),y \big\rangle_{X,X^*}
=\big\langle Da(x),y \big\rangle_{X,X^*}-\big\langle aD(x),y \big\rangle_{X,X^*} \\
&=\big\langle a(x),T(y) \big\rangle_{X,X^*}-\big\langle D(x),a^*(y)\big\rangle_{X,X^*} 
=\big\langle x,a^*T(y) \big\rangle-\big\langle x,T a^*(y) \big\rangle_{X,X^*} \\
&=\big\langle x,-[T,a^*](y)\big\rangle_{X,X^*}.
\end{align*}
We infer that $y \in \dom [D,a]^*$. Hence $Y \subseteq \dom [D,a]^*$. Since $Y$ is dense in $X^*$, we deduce that $[D,a]^*$ is densely defined. By \cite[Theorem 5.28 p.~168]{Kat1}, this implies that $[D,a]$ is closable. From $[D,a]|_C \subseteq [D,a]$, we have $\ovl{[D,a]|_C} \subseteq \ovl{[D,a]}$ and the closure $\ovl{[D,a]|_C}$ is the bounded extension of $[D,a]|_C$ on $X$ (note that $C$ is dense in $X$ since $\dom D$ is dense). Hence the closed operator $\ovl{[D,a]}$ is defined on $X$ hence bounded by \cite[Theorem 5.20 p.~166]{Kat1}.
\end{proof}

Next we define the following notion of a locally compact spectral triple.
When $X$ is a Hilbert space, then compare to the one in \cite{SuZ1}, \cite[p.~588]{GGISV1}.

\begin{defi}
\label{def-locally-compact-spectral-triple}
Let $X$ be a reflexive Banach space, $D$ a densely defined closed operator on $X$ and $A$ a subalgebra of some Banach algebra.
Let $\pi \co A \to B(X)$ be a homomorphism.
We call $(A,X,D)$ a locally compact (Banach) spectral triple, provided that
\begin{enumerate}
\item $D$ is bisectorial on $X$ and has a bounded $\HI$ functional calculus on some bisector.
\item For any $a \in A$, we have $\pi(a) \cdot \dom D \subseteq \dom D$ and $[D,\pi(a)] \co \dom D \to X$ which is densely defined, extends to a bounded operator on $X$.
\item For any $a \in A$, $\pi(a) (\i \Id + D)^{-1}|\ovl{\Ran D}$ is a compact operator $\ovl{\Ran D} \to X$.
\end{enumerate}
\end{defi}

Recall that a spectral triple $(A,H,D)$ is even if there exists a selfadjoint unitary operator $\gamma \co H \to H$ such that $\gamma^2=\Id_H$, $\gamma D=-D\gamma$ and $\gamma \pi(a)=\pi(a)\gamma$ for all $a \in A$.

\begin{defi}
\label{Def-even}
We say that a (locally compact or compact) Banach spectral triple $(A,X,D)$ is even if there exists a surjective isometry $\gamma \co X \to X$ such that $\gamma^2=\Id_X$, $\gamma(\dom D) \subseteq \dom D$, $\gamma D=-D\gamma$ and $\gamma \pi(a)=\pi(a)\gamma$ for all $a \in A$.  
\end{defi}

\subsection{Spectral triples associated to semigroups of Fourier multipliers I}
\label{Sec-spectral-triples-Fourier-I}

In this subsection, we consider a markovian semigroup of Fourier multipliers as in  Proposition \ref{prop-Schoenberg}, together with the noncommutative gradient $\partial_{\alpha,q,p}$ and its adjoint. Now, we generalize the construction of \cite[p.~587-589]{JMP2} which corresponds to the case $q=1$ and $p=2$ below. So we obtain a scale of $\L^p$-Banach spectral triples associated to these semigroups. In particular, here we shall complete the picture. Indeed, we prove the compactness axiom (point 6 below) of Connes' spectral triples. Finally, we shall compare our result Theorem \ref{First-spectral-triple-2} with the result of \cite[p.~588]{JMP2} in Remark \ref{rem-Fourier-I-JMP2}. 



Suppose $1 < p<\infty$ and $-1 \leq q < 1$. Recall that the Hodge-Dirac operator is defined by
\begin{equation}
\label{Def-Dirac-operator-Fourier-2}
\D_{\psi,q,p} 
\ov{\mathrm{def}}{=}\begin{bmatrix} 
0 & (\partial_{\psi,q,p^*})^* \\ 
\partial_{\psi,q,p}& 0 
\end{bmatrix}
\end{equation}
on the subspace $\dom \D_{\psi,q,p} = \dom \partial_{\psi,q,p} \oplus \dom (\partial_{\psi,q,p^*})^* $ of $\L^p(\VN(G)) \oplus_p \L^p(\Gamma_q(H) \rtimes_{\alpha} G)$. If $a \in \mathrm{C}^*_r(G)$, we define the bounded operator $\pi(a) \co \L^p(\VN(G)) \oplus_p  \L^p(\Gamma_q(H) \rtimes_{\alpha} G) \to \L^p(\VN(G)) \oplus_p \L^p(\Gamma_q(H) \rtimes_{\alpha} G)$ by
\begin{equation}
\label{Def-pi-a-Fourier}
\pi(a)
\ov{\mathrm{def}}{=} \begin{bmatrix}
    \L_a & 0  \\
    0 & \tilde{\L}_a  \\
\end{bmatrix}, \quad a \in \mathrm{C}^*_r(G)
\end{equation}
where $\L_a \co \L^p(\VN(G)) \to \L^p(\VN(G))$, $x \mapsto ax$ is the left multiplication operator and where $\tilde{\L}_a \co \L^p(\Gamma_q(H) \rtimes_{\alpha} G) \to \L^p(\Gamma_q(H) \rtimes_{\alpha} G)$ is the left action of the bimodule. It is easy to check that $\pi \co \mathrm{C}^*_r(G)\to \B(\L^p(\VN(G)) \oplus_p  \L^p(\Gamma_q(H) \rtimes_{\alpha} G))$ is continuous when $\mathrm{C}^*_r(G)$ is equipped with the induced weak* topology of $\VN(G)$ and when $\B(\L^p(\VN(G)) \oplus_p  \L^p(\Gamma_q(H) \rtimes_{\alpha} G))$ is equipped with the weak operator topology. See also Lemma \ref{Lemma-continuous-Schur-I} below.

We will use the following proposition.

\begin{prop}
\label{prop-Fourier-vraiment-compact}
Let $1 < p < \infty$. Let $G$ be a discrete group. Consider the associated function $b_\psi \co G \to H$ from \eqref{liens-psi-bpsi}. Assume that the Hilbert space $H$ is finite-dimensional, that $b_\psi$ is injective and that $\G_\psi > 0$. Finally assume that the von Neumann crossed product $\Gamma_q(H) \rtimes_\alpha G$ has $\QWEP$.

\begin{enumerate}
\item
The operator $A_p^{-\frac12} \co \ovl{\Ran A_p} \to \ovl{\Ran A_p}$ is compact.

\item Suppose $-1 \leq q \leq 1$. Then the operator $B_p^{-\frac12}$ is compact where $B_p$ is defined in \eqref{Def-Bp-Fourier}.
\end{enumerate}
\end{prop}

\begin{proof}
1. As in the proof of Theorem \ref{thm-Fourier-quantum-compact-metric} the operator $A_2^{-\frac12} \co \ovl{\Ran A_2} \to \ovl{\Ran A_2}$ is compact. 

If we view $A_\infty$ as acting on $\mathrm{C}_r^*(G)$, we show that $A_\infty^{-\frac12} \co \ovl{\Ran A_\infty} \to \ovl{\Ran A_\infty}$ is bounded. Since $(T_{t,\infty})_{t \geq 0}$ is a semigroup with generator $A_\infty$, it suffices by Lemma \ref{lem-A12-bounded} to establish the bound $\norm{T_t}_{\Ran A_\infty \to \Ran A_\infty} \lesssim \frac{1}{t^d}$ for some $d > \frac12$. We conclude with \cite[Lemma 5.8]{JMP1} and the contractive inclusion $\VN(G) \subseteq \L^1(\VN(G))$. Thus, $A_{\infty}^{-\frac12} \co \ovl{\Ran A_\infty} \to \ovl{\Ran A_\infty}$ is bounded. 

Now, assume that $p > 2$. We will use complex interpolation. Since the resolvent of the operators $A_p$ are compatible for different values of $p$, the spectral projections \cite[p.~364]{HvNVW2} onto the spaces $\ovl{\Ran A_p}$ are compatible. Hence, the $\ovl{\Ran A_p}$'s form an interpolation scale. Observe that $\ovl{\Ran A_2}$ is a Hilbert space, hence a UMD space. Then we obtain the compactness of $A_p^{-\frac12} \co \ovl{\Ran A_p} \to \ovl{\Ran A_p}$ by means of complex interpolation between a compact and a bounded operator with Theorem \ref{Th-interpolation-Kalton}. 

If $p < 2$, we conclude by duality and Schauder's Theorem \cite[Theorem 3.4.15]{Meg1}, since $\ovl{\Ran A_p}$ is the dual space of $\ovl{\Ran A_{p^*}}$ and $A_p^{-\frac12}$ defined on the first space is the adjoint of $A_{p^*}^{-\frac12}$ defined on the second space.

2. 
According to the first part of the proposition, $A_2^{-\frac12} \co \StwoI \to \StwoI$ is compact. Since $A_2=(\partial_{\psi,q,2})^*\partial_{\psi,q,2}$ and $B_2=\partial_{\psi,q,2}(\partial_{\psi,q,2})^*|\ovl{\Ran\partial_{\psi,q,2}}$, \cite[Corollary 5.6]{Tha1} and Proposition \ref{Prop-liens-ranges-group} together with a functional calculus argument show that the operators $A_2^{-\frac12}|_{\ker( \partial_{\psi,q,2})^\perp}$ and $B_2^{-\frac12}$ are unitarily equivalent. Since $A_2^{-\frac12} \co \StwoI \to \StwoI$ is compact and $\ker (\partial_{\psi,q,2})^\perp = \StwoI$, $A_2^{-\frac12}|_{\ker( \partial_{\psi,q,2})^\perp}$ and finally $B_2^{-\frac12} \co \ovl{\Ran \partial_{\psi,q,2}} \to \ovl{\Ran \partial_{\psi,q,2}}$ is compact.

Recall that by \eqref{Def-Bp-Fourier}, $B_p$ is the generator of the (strongly continuous) semigroup $(\Id_{\L^p(\Gamma_q(H))} \rtimes T_t|\ovl{\Ran \partial_{\psi,q,p}})_{t \geq 0}$.  
Now, we use \cite[Lemma 5.8]{JMP1} and Proposition \ref{prop-Fourier-mult-crossed-product}  to obtain the bound $\norm{\Id_{} \rtimes T_t}_{\infty \to \infty} \lesssim \frac{1}{t^d}$ for some $d > \frac12$. According to Proposition \ref{Prop-liens-ranges-group}, we have $\ovl{\Ran B_p} = \ovl{\Ran \partial_{\psi,q,p}}$. By Lemma \ref{lem-A12-bounded}, we obtain that $B_p^{-\frac12} \co \ovl{\Ran B_p } \to \ovl{\Ran B_p}$ is bounded. 

Now, it suffices to interpolate with Theorem \ref{Th-interpolation-Kalton} the compactness of $B_2^{-\frac12} \co \ovl{\Ran \partial_{\psi,q,2}} \to \ovl{\Ran \partial_{\psi,q,2}}$ and the boundedness of $B_{p_0}^{-\frac12} \co \ovl{\Ran \partial_{\psi,q,p_0}} \to \ovl{\Ran \partial_{\psi,q,p_0}}$ to obtain compactness of $B_p^{-\frac12} \co \ovl{\Ran \partial_{\psi,q,p}} \to \ovl{\Ran \partial_{\psi,q,p}}$  for $2 \leq p < p_0 < \infty$  and the fact that on the $\L^2$-level, we have a Hilbert space hence a UMD space. Note that the spaces $\ovl{\Ran \partial_{\psi,q,p}} \subseteq \L^p(\Gamma_q(H) \otvn \B(\ell^2_I))$ interpolate by the complex interpolation method since we have bounded projections
$$
P_p \co \L^p(\Gamma_q(H) \rtimes G ) \to \ovl{\Ran \partial_{\psi,q,p}} \subseteq \L^p(\Gamma_q(H) \rtimes G)
$$
which are compatible for different values of $p$ according to \cite[p. 361 and p. 448]{HvNVW2}. We use duality if $p<2$.
\end{proof}

In the following theorem, recall the weak* closed operator $\partial_{\psi,q,\infty} \co \dom \partial_{\psi,q,\infty} \subseteq \Gamma_q(H) \rtimes_\alpha G \to \Gamma_q(H) \rtimes_\alpha G$ from Proposition \ref{Prop-closable-deriv-p-infty}.
The latter is valid if $G$ has $\mathrm{AP}$ and $q \neq 1$.

\begin{thm}
\label{First-spectral-triple-2}
Suppose $1<p<\infty$ and $-1 \leq q < 1$.
Consider the triple $(\mathrm{C}^*_r(G),\L^p(\VN(G)) \oplus_p \L^p(\Gamma_q(H) \rtimes_\alpha G), \D_{\psi,q,p})$.
It satisfies the following properties. In particular, it is a Banach spectral triple in the sense of Definition \ref{Def-Banach-spectral-triple}, in case that $b_\psi \co G \to H$ is injective, $\G_\psi > 0$, $H$ is finite dimensional and $\Gamma_q(H) \rtimes_\alpha G$ has $\QWEP$ (e.g. $G$ is amenable, or $G$ is a free group and $q = -1$).

\begin{enumerate}

\item We have $(\D_{\psi,q,p})^*=\D_{\psi,q,p^*}$. In particular, the operator $\D_{\psi,q,2}$ is selfadjoint.

\item We have
\begin{equation}
\label{Lip-algebra-I-Fourier}
\P_G 
\subseteq \Lip_{\D_{\psi,q,p}}(\mathrm{C}^*_r(G)).
\end{equation}

\item For any $a \in \P_G$, we have 
\begin{equation}
\label{Majo-commutateur-group}
\bnorm{\big[\D_{\psi,q,p},\pi(a)\big]}_{\L^p(\VN(G)) \oplus_p \L^p(\Gamma_q(H) \rtimes_{\alpha} G) \to \L^p(\VN(G)) \oplus_p \L^p(\Gamma_q(H) \rtimes_{\alpha} G)}
\leq \norm{\partial_{\psi,q}(a)}_{\Gamma_q(H) \rtimes_{\alpha} G}.
\end{equation}

\item Suppose that $G$ has $\mathrm{AP}$. We have
\begin{equation}
\label{Lip-algebra-I-Fourier}
\mathrm{C}^*_r(G) \cap \dom \partial_{\psi,q,\infty} 
\subseteq \Lip_{\D_{\psi,q,p} }(\mathrm{C}^*_r(G)).
\end{equation}

\item Suppose that $G$ has $\mathrm{AP}$. For any $a \in \mathrm{C}^*_r(G) \cap \dom \partial_{\psi,q,\infty}$, we have 
\begin{equation}
\label{norm-commutator-group}
\norm{\big[\D_{\psi,q,p},\pi(a)\big]}_{\L^p \oplus_p \L^p(\Gamma_q(H) \rtimes_{\alpha} G)  \to \L^p \oplus_p \L^p(\Gamma_q(H) \rtimes_{\alpha} G)} 
\leq \bnorm{\partial_{\psi,q,\infty}(a)}_{\Gamma_q(H)  \rtimes_{\alpha} G}.
\end{equation}

\item Assume that $\Gamma_q(H) \rtimes_\alpha G$ has $\QWEP$.
If $b_\psi \co G \to H$ is injective, $\G_\psi >0$ and if $H$ is finite-dimensional, then the operator 
$$
|\D_{\psi,q,p}|^{-1} \co \ovl{\Ran (\partial_{\psi,q,p^*})^*} \oplus \ovl{\Ran \partial_{\psi,q,p}} \to \ovl{\Ran (\partial_{\psi,q,p^*})^*} \oplus \ovl{\Ran \partial_{\psi,q,p}} 
$$ 
is compact.
\end{enumerate}
\end{thm}

\begin{proof}
1. An element $(z,t)$ of $\L^{p^*}(\VN(G)) \oplus_{p^*} \L^{p^*}(\Gamma_q(H) \rtimes_{\alpha} G)$ belongs to $\dom (\D_{\psi,q,p})^*$ if and only if there exists $(a,b) \in \L^{p^*}(\VN(G)) \oplus_{p^*} \L^{p^*}(\Gamma_q(H) \rtimes_{\alpha} G)$ such that for any $(x,y) \in \dom \partial_{\psi,q,p} \oplus \dom (\partial_{\psi,q,p^*})^*$ we have
\begin{align}
\label{Relation99-group}
\MoveEqLeft
\left\langle \begin{bmatrix} 
0 & (\partial_{\psi,q,p^*})^* \\ 
\partial_{\psi,q,p} & 0 
\end{bmatrix}\begin{bmatrix}
      x  \\
      y  \\
\end{bmatrix},\begin{bmatrix}
      z  \\
      t  \\
\end{bmatrix}\right\rangle            
=\left\langle 
\begin{bmatrix}
    x    \\
    y    \\
\end{bmatrix},\begin{bmatrix}
    a    \\
    b    \\
\end{bmatrix}\right\rangle \\
&\iff \big\langle (\partial_{\psi,q,p^*})^*(y), z\big\rangle
+\big\langle \partial_{\psi,q,p}(x),t \big\rangle
=\langle y,b\rangle+\langle x,a\rangle. \nonumber
\end{align} 
If $z \in \dom \partial_{\psi,q,p^*}$ and if $t \in \dom (\partial_{\psi,q,p})^*$ the latter holds with $b=\partial_{\psi,q,p^*}(z)$ and $a=(\partial_{\psi,q,p})^*(t)$. This proves that $\dom \partial_{\psi,q,p^*} \oplus \dom (\partial_{\psi,q,p})^* \subseteq \dom (\D_{\psi,q,p})^*$ and that 
$$
(\D_{\psi,q,p})^*(z,t)
=\big((\partial_{\psi,q,p})^*(t),\partial_{\psi,q,p^*}(z)\big)
= \begin{bmatrix} 
0 & (\partial_{\psi,q,p})^* \\ 
\partial_{\psi,q,p^*} & 0 
\end{bmatrix}
\begin{bmatrix}
    z    \\
    t    \\
\end{bmatrix}
\ov{\eqref{Def-Dirac-operator-Fourier-2}}{=} \D_{\psi,q,p^*}(z,t).
$$
Conversely, if $(z,t) \in \dom (\D_{\psi,q,p})^*$, choosing $y=0$ in \eqref{Relation99-group} we obtain $t \in \dom (\partial_{\psi,q,p})^*$ and taking $x=0$ we obtain $z \in \dom \partial_{\psi,q,p^*}$. 

2. We recall the dense subspaces $\P_G$ and $\P_{\rtimes,G}$ of $\L^p(\VN(G))$ and $\L^p(\Gamma_q(H) \rtimes_\alpha G)$ which are contained in the domains of $\partial_{\psi,q,p^*}$ and $(\partial_{\psi,q,p^*})^*$. So $\P_G \oplus \P_{\rtimes,G}$ is contained in $\dom \D_{\psi,q,p}$. For any $a \in \P_G$, we have $\L_a(\P_G) \subseteq \P_G$ and $\tilde{\L}_a(\P_{\rtimes,G}) \subseteq \P_{\rtimes,G}$. We infer that $\pi(a) \cdot (\P_G \oplus \P_{\rtimes,G}) \subseteq \dom \D_{\psi,q,p}$. Note also that $\pi(a)^* \cdot (\P_G \oplus \P_{\rtimes,G}) \subseteq \dom \D_{\psi,q,p^*}=\dom (\D_{\psi,q,p})^*$.

Let $a \in \P_G$. A simple computation shows that
\begin{align*}
\MoveEqLeft
\big[\D_{\psi,q,p},\pi(a)\big]
\ov{\eqref{Def-Dirac-operator-Fourier-2}\eqref{Def-pi-a-Fourier}}{=}\begin{bmatrix} 
0 & (\partial_{\psi,q,p^*})^* \\ 
\partial_{\psi,q,p}& 0
\end{bmatrix}\begin{bmatrix}
    \L_a & 0  \\
    0 & \tilde{\L}_a  \\
\end{bmatrix}
-\begin{bmatrix}
    \L_a & 0  \\
    0 & \tilde{\L}_a  \\
\end{bmatrix}
\begin{bmatrix} 
0 & (\partial_{\psi,q,p^*})^* \\ 
\partial_{\psi,q,p}& 0
\end{bmatrix}\\
&=\begin{bmatrix}
   0 &  (\partial_{\psi,q,p^*})^*\tilde{\L}_a \\
   \partial_{\psi,q,p}\L_a  &0  \\
\end{bmatrix}-\begin{bmatrix}
    0 &  \L_a (\partial_{\psi,q,p^*})^*\\
    \tilde{\L}_a\partial_{\psi,q,p} & 0  \\
\end{bmatrix} \\
&=\begin{bmatrix}
    0 &  (\partial_{\psi,q,p^*})^*\tilde{\L}_a-\L_a (\partial_{\psi,q,p^*})^* \\
   \partial_{\psi,q,p}\L_a-\tilde{\L}_a\partial_{\psi,q,p}  & 0  \\
\end{bmatrix}.
\end{align*}
We calculate the two non-zero components of the commutator. For the lower left corner, if $x \in \P_G$ and if we consider the canonical map $J \co \L^p(\VN(G)) \to \L^p(\Gamma_q(H) \rtimes_{\alpha} G)$, $x \mapsto 1 \rtimes x$, we have\footnote{\thefootnote. Recall that the term $\partial_{\psi,q,p}(a)x$ is by definition equal to $\partial_{\psi,q,p}(a)(1 \rtimes x)$.}
\begin{align}
\label{Bon-commutateur-group}
\MoveEqLeft
(\partial_{\psi,q,p}\L_a-\tilde{\L}_a\partial_{\psi,q,p})(x)           
=\partial_{\psi,q,p}\L_a(x)-\tilde{\L}_a\partial_{\psi,q,p}(x)
=\partial_{\psi,q,p}(ax)-a\partial_{\psi,q,p}(x) \\
&\ov{\eqref{Leibniz-Schur-gradient-mieux-sgrp}}{=} \partial_{\psi,q}(a)x 
=\L_{\partial_{\psi,q}(a)}J(x). \nonumber
\end{align}
For the upper right corner, if $\E$ is the conditional expectation associated to $J$, note that for any $y \in \P_{\rtimes,G}$ and any $x \in \P_G$, (we recall that we have the duality brackets $\langle f , g \rangle$ antilinear in the first variable)
\begin{align*}
\MoveEqLeft
\big\langle \big((\partial_{\psi,q,p^*})^*\tilde{\L}_a-\L_a (\partial_{\psi,q,p^*})^*\big)(y),x \big\rangle_{}            
=\big\langle (\partial_{\psi,q,p^*})^*\tilde{\L}_a(y),x \big\rangle -\big\langle \L_a (\partial_{\psi,q,p^*})^*(y),x \big\rangle_{}\\
&=\big\langle \tilde{\L}_a(y),\partial_{\psi,q,p^*}(x) \big\rangle -\big\langle  (\partial_{\psi,q,p^*})^*(y),\L_{a^*}(x) \big\rangle_{}
=\big\langle y,\tilde{\L}_{a^*}\partial_{\psi,q,p^*}(x) \big\rangle -\big\langle y,\partial_{\psi,q,p^*}\L_{a^*}(x) \big\rangle_{} \\
&=\big\langle y,\tilde{\L}_{a^*}\partial_{\psi,q,p^*}(x)-\partial_{\psi,q,p^*}\L_{a^*}(x) \big\rangle_{}
=\big\langle y,a^*\partial_{\psi,q,p^*}(x)-\partial_{\psi,q,p^*}(a^*x) \big\rangle_{} \\
&\ov{\eqref{Leibniz-Schur-gradient-mieux-sgrp}}{=}\big\langle y,-\partial_{\psi,q}(a^*)x \big\rangle_{}
=\big\langle y, -{\L}_{\partial_{\psi,q}(a^*)}(1 \rtimes x)\big\rangle_{} = \big\langle y , \L_{(\partial_{\psi,q}(a))^*}(1 \rtimes x) \big\rangle_{}\\
& =\big\langle {\L}_{\partial_{\psi,q}(a)}(y), 1 \rtimes x\big\rangle
=\big\langle \E{\L}_{\partial_{\psi,q}(a)}(y),x\big\rangle_{\L^p(\VN(G)),\L^{p^*}(\VN(G))}.
\end{align*} 
We conclude that
\begin{equation}
\label{Commutateur-etrange-group}
\big((\partial_{\psi,q,p^*})^*\tilde{\L}_a-\L_a (\partial_{\psi,q,p^*})^*\big)(y)
=\E{\L}_{\partial_{\psi,q}(a)}(y).
\end{equation}
The two non-zero components of the commutator are bounded linear operators on $\P_G$ and on $\P_{\rtimes,G}$. We deduce that $\big[\D_{\psi,q,p},\pi(a)\big]$ is bounded on $\P_G \oplus \P_{\rtimes,G}$.
Since $[\D_{\psi,q,p},\pi(a)]$ is closable\footnote{\thefootnote. Let $x_n \in X := \L^p(\VN(G)) \oplus \L^p(\Gamma_q(H) \rtimes_\alpha G)$ such that $x_n \to 0$, and $[\D_{\psi,q,p},\pi(a)]x_n \to y$.
If $z$ belongs to the dense subset $\P_G \oplus \P_{\rtimes,G}$ of $X^*$, we have on the one hand $\langle [\D_{\psi,q,p},\pi(a)]x_n, z \rangle \to \langle y , z \rangle$ and on the other hand, $\langle [\D_{\psi,q,p},\pi(a)]x_n, z \rangle = \langle x_n, -[\D_{\psi,q,p}^*,\pi(a)^*]z \rangle \to \langle 0 , z\rangle$. We infer that $y = 0$.}, densely defined and bounded on the dense subspace $\P_G \oplus \P_{\rtimes,G}$, it extends to a bounded operator on $\L^p(\VN(G)) \oplus_p \L^p(\Gamma_q(H) \rtimes_{\alpha} G)$. Hence $\P_G$ is a subset of $\Lip_{\D_{\psi,q,p} }(\mathrm{C}^*_r(G))$. 

3. If $(x,y) \in \dom \D_{\psi,q,p}$ and $a \in \P_G$, we have
\begin{align}
\label{Divers-5436-group}
\MoveEqLeft
\bnorm{\big[\D_{\psi,q,p},\pi(a)\big](x,y)}_{p}
=\norm{\big(\big((\partial_{\psi,q,p^*})^*\tilde{\L}_a - \L_a (\partial_{\psi,q,p^*})^* \big)y, \big(\partial_{\psi,q,p} \L_a - \tilde{\L}_a \partial_{\psi,q,p} \big)x \big)}_p \\          
&=\big(\norm{\big((\partial_{\psi,q,p^*})^*\tilde{\L}_a - \L_a (\partial_{\psi,q,p^*})^* \big)y}_{\L^p(\VN(G))}^p+\norm{\big(\partial_{\psi,q,p} \L_a - \tilde{\L}_a \partial_{\psi,q,p} \big)x}_p^p\big)^{\frac{1}{p}}\nonumber \\
&\ov{\eqref{Commutateur-etrange-group} \eqref{Bon-commutateur-group}}{=} \big(\norm{\E{\L}_{\partial_{\psi,q}(a)}(y)}_{\L^p(\VN(G))}^p+\norm{\partial_{\psi,q}(a) J(x)}_{\L^p(\Gamma_q(H) \rtimes_{\alpha} G)}^p \big)^{\frac1p} 
\nonumber\\
&\leq \norm{\partial_{\psi,q}(a)}_{\Gamma_q(H)  \rtimes_{\alpha} G} \norm{(x,y)}_p. \nonumber
\end{align}
So we obtain \eqref{Majo-commutateur-group}. 

4. Let $a \in \mathrm{C}^*_r(G) \cap \dom \partial_{\psi,q,\infty}$. Let $(a_j)$ be a net in $\P_G$ such that $a_j \to a$ and $\partial_{\psi,q,\infty}(a_j) \to \partial_{\psi,q,\infty}(a)$ both for the weak* topology. By Krein-Smulian Theorem, we can suppose that the nets $(a_j)$ and $(\partial_{\psi,q,\infty}(a_j))$ are bounded. By the point 4 of Proposition \ref{Prop-Lip-algebra} and Remark \ref{Remark-123}, we deduce that $a \in \Lip_{\D_{\psi,q,p} }(\mathrm{C}^*_r(G))$. By continuity of $\pi$, note that $\pi(a_j) \to \pi(a)$ for the weak operator topology. For any $\xi \in \dom \D_{\psi,q,p}$ and any $\zeta \in \dom (\D_{\psi,q,p})^*$, we have
\begin{align*}
\MoveEqLeft
\big\langle [\D_{\psi,q,p},\pi(a_j)]\xi,\zeta \big\rangle_{\L^p(\VN(G)) \oplus_p \L^p(\Gamma_q(H) \rtimes_{\alpha} G),\L^{p^*}(\VN(G)) \oplus_{p^*} \L^{p^*}(\Gamma_q(H) \rtimes_{\alpha} G)} \\           
&=\big\langle (\D_{\psi,q,p}\pi(a_j)-\pi(a_j)\D_{\psi,q,p})\xi, \zeta\big\rangle_{} 
=\big\langle \D_{\psi,q,p}\pi(a_j)\xi, \zeta\big\rangle-\big\langle \pi(a_j) \D_{\psi,q,p}\xi, \zeta\big\rangle \\
&=\big\langle \pi(a_j)\xi, (\D_{\psi,q,p})^*\zeta\big\rangle-\big\langle \pi(a_j)\D_{\psi,q,p}\xi, \zeta \big\rangle_{} 
\xra[j]{} \big\langle \pi(a)\xi, (\D_{\psi,q,p})^*\zeta\big\rangle-\big\langle \pi(a)\D_{\psi,q,p}\xi, \zeta \big\rangle_{} \\
&=\big\langle \D_{\psi,q,p}\pi(a)\xi, \zeta\big\rangle-\big\langle \pi(a)\D_{\psi,q,p}\xi, \zeta \big\rangle_{}=\big\langle [\D_{\psi,q,p},\pi(a)]\xi,\zeta \big\rangle_{}.
\end{align*}
Since the net $([\D_{\psi,q,p},\pi(a_j)])$ is bounded by \eqref{Majo-commutateur-group}, we deduce that $([\D_{\psi,q,p},\pi(a_j)])$ converges to $[\D_{\psi,q,p},\pi(a)]$ for the weak operator topology by a ``net version'' of \cite[Lemma 3.6 p.~151]{Kat1}. Furthermore, it is (really) easy to check that $\L_{\partial_{\psi,q,p}(a_j)}J \to \L_{\partial_{\psi,q,p}(a)}J$ and $-\E{\L}_{\partial_{\psi,q}(a_j)} \to -\E{\L}_{\partial_{\psi,q}(a)}$ both for the weak operator topology. By uniqueness of the limit, we deduce that the commutator is given by the same formula as that in the case of elements of $\P_G$. 

5. We obtain \eqref{norm-commutator-group} as in \eqref{Divers-5436-group}.

6. Recall that $\ovl{\Ran  A_p}=\ovl{\Ran (\partial_{\psi,q,p^*})^*}$ by Proposition \ref{Prop-liens-ranges-group} and that
$$
\D_{\psi,q,p}^2
\ov{\eqref{D-alpha-carre-egal-Fourier}}{=}
\begin{bmatrix} 
A_p & 0 \\ 
0 & (\Id_{\L^p(\Gamma_q(H))} \rtimes A_p)|\ovl{\Ran \partial_{\psi,q,p}}
\end{bmatrix}
\ov{\eqref{Def-Bp-Fourier}}{=}\begin{bmatrix} 
A_p & 0 \\ 
0 & B_p
\end{bmatrix}.
$$
So
\begin{align*}
\MoveEqLeft
|\D_{\psi,q,p}|^{-1}
=\begin{bmatrix} 
A_p|\ovl{\Ran  A_p} & 0 \\ 
0 & B_p
\end{bmatrix}^{-\frac{1}{2}}
=\begin{bmatrix} 
A_p^{-\frac{1}{2}}|\ovl{\Ran  A_p} & 0 \\ 
0 & B_p^{-\frac{1}{2}}
\end{bmatrix}.            
\end{align*} 
With Proposition \ref{prop-Fourier-vraiment-compact}, we conclude that $|\D_{\psi,q,p}|^{-1}$ is compact. 
\end{proof}


\begin{remark} \normalfont
\label{Remark-open-question-group}
Suppose that $G$ has $\mathrm{AP}$. We do not know if $\Lip_{\D_{\psi,q,p} }(\mathrm{C}^*_r(G))=\mathrm{C}^*_r(G) \cap \dom \partial_{\psi,q,\infty}$.  
\end{remark}

\begin{remark} \normalfont
\label{Rem-even-group-I}
Consider the assumptions of Theorem \ref{First-spectral-triple-2} be satisfied.
Note that the (Banach) spectral triple $(\mathrm{C}^*_r(G),\L^p(\VN(G)) \oplus_p \L^p(\Gamma_q(H) \rtimes_{\alpha} G),\D_{\psi,q,p})$ is even. Indeed, the Hodge-Dirac operator $\D_{\psi,q,p}$ anti-commutes with the involution 
\begin{align*}
\MoveEqLeft
\gamma_p
\ov{\mathrm{def}}{=}
\begin{bmatrix} 
-\Id_{\L^p(\VN(G))} & 0 \\ 
0& \Id_{\L^p(\Gamma_q(H) \rtimes_{\alpha} G)} 
\end{bmatrix} \co \\
& \L^p(\VN(G)) \oplus_p \L^p(\Gamma_q(H) \rtimes_{\alpha} G) \to \L^p(\VN(G)) \oplus_p \L^p(\Gamma_q(H) \rtimes_{\alpha} G). 
\end{align*}  
(which is selfadjoint if $p=2$), since
\begin{align*}
\MoveEqLeft
\D_{\psi,q,p}\gamma_p+\gamma_p\D_{\psi,q,p}\\
&\ov{\eqref{Def-Dirac-operator-Fourier-2}}{=} \begin{bmatrix} 
0 & (\partial_{\psi,q,p^*})^* \\ 
\partial_{\psi,q,p}& 0 
\end{bmatrix}
\begin{bmatrix} 
-\Id_{} & 0 \\ 
0& \Id_{}
\end{bmatrix}
+\begin{bmatrix} 
-\Id_{} & 0 \\ 
0& \Id_{}
\end{bmatrix}
\begin{bmatrix} 
0 & (\partial_{\psi,q,p^*})^* \\ 
\partial_{\psi,q,p}& 0 
\end{bmatrix}\\
&=\begin{bmatrix} 
0 & (\partial_{\psi,q,p^*})^* \\ 
-\partial_{\psi,q,p}&  0
\end{bmatrix}+\begin{bmatrix} 
0 & -(\partial_{\psi,q,p^*})^* \\ 
\partial_{\psi,q,p}&  0
\end{bmatrix}
=0.            
\end{align*} 
Moreover, for any $a \in \mathrm{C}^*_r(G)$, we have
\begin{align*}
\MoveEqLeft
\gamma_p\pi(a)
\ov{\eqref{Def-pi-a-Fourier}}{=} \begin{bmatrix} 
-\Id_{} & 0 \\ 
0& \Id_{}
\end{bmatrix} 
\begin{bmatrix}
    \L_a & 0  \\
    0 & \tilde{\L}_a  \\
\end{bmatrix}  
=\begin{bmatrix}
    -\L_a & 0  \\
    0 & \tilde{\L}_a  \\
\end{bmatrix}        
=\begin{bmatrix}
    \L_a & 0  \\
    0 & \tilde{\L}_a  \\
\end{bmatrix}\begin{bmatrix} 
-\Id_{} & 0 \\ 
0& \Id_{}
\end{bmatrix}
\ov{\eqref{Def-pi-a-Fourier}}{=} \pi(a)\gamma_p.
\end{align*} 
\end{remark}

\begin{remark} \normalfont
\label{rem-Fourier-I-not-optimal} 
The estimate \eqref{Majo-commutateur-group} is in general not optimal.
Indeed, already in the case $p = 2$ and $a = \lambda_s \in \P_G$ for some $s \in G$, we have according to \eqref{Bon-commutateur-group} and \eqref{Commutateur-etrange-group},
\begin{align}
\label{equ-1-rem-Fourier-I-not-optimal}
\MoveEqLeft
\norm{[\D_{\psi,q,2},\pi(a)]}_{\L^2(\VN(G))\oplus_2 \L^2(\Gamma_q(H)\rtimes_\alpha G) \to \L^2(\Gamma_q(H)) \oplus_2 \L^2(\Gamma_q(H) \rtimes_\alpha G)} \\
& \leq \max\left\{ \norm{\L_{\partial_{\psi,q}(a)}J}_{\L^2(\VN(G))\to \L^2(\Gamma_q(H) \rtimes_\alpha G)}, \norm{ \E \L_{\partial_{\psi,q}(a)} }_{\L^2(\Gamma_q(H) \rtimes_\alpha G) \to \L^2(\VN(G))} \right\} . \nonumber
\end{align}
Note that we have the Hilbert space adjoints $(\L_{\partial_{\psi,q}(a)}J)^* = J^* \L_{\partial_{\psi,q}(a)}^* = \E \L_{(\partial_{\psi,q}(a))^*} = - \E \L_{\partial_{\psi,q}(a^*)}$.
Thus, in the maximum of \eqref{equ-1-rem-Fourier-I-not-optimal}, it suffices to consider the second term. For any element $x = \sum_{t} x_{t} \rtimes \lambda_t$ of $\L^2(\Gamma_q(H) \rtimes_\alpha G)$, we have 
\begin{align*}
\MoveEqLeft
\E \L_{\partial_{\psi,q}(\lambda_s)}(x) 
= \E \bigg(\big(s_q(b_\psi(s)) \rtimes \lambda_s\big) \bigg(\sum_t x_t \rtimes \lambda_t\bigg) \bigg)
\ov{\eqref{Product-crossed-product}}{=} \E\bigg(\sum_t s_q(b_\psi(s)) \alpha_s(x_{t}) \rtimes \lambda_{st}\bigg) \\
& = \sum_t \tau(s_q(b_\psi(s))\alpha_s(x_{t})) \lambda_{st}.
\end{align*}
Thus we have
\begin{align*}
\MoveEqLeft
\norm{ \E \L_{\partial_{\psi,q}(\lambda_s)}(x) }_2^2 = \sum_t \left| \tau(s_q(b_\psi(s)) \alpha_s(x_{t})) \right|^2 
\leq \norm{s_q(b_\psi(s))}_{\L^2(\Gamma_q(H))}^2 \sum_t \norm{\alpha_s(x_{t})}_{\L^2(\Gamma_q(H))}^2  \\
& = \norm{s_q(b_\psi(s))}_{\L^2(\Gamma_q(H))}^2 \sum_t \norm{x_{t}}_{\L^2(\Gamma_q(H))}^2 \ov{\eqref{petit-Wick}}{=} \norm{b_\psi(s)}_{H}^2  \norm{x}_{\L^2(\Gamma_q(H) \otvn \B(\ell^2_I))}^2.
\end{align*}
We infer that
\begin{equation}
\label{equ-2-rem-Fourier-I-not-optimal}
\norm{[\D_{\psi,q,2},\pi(\lambda_s)]}_{\L^2 \oplus_2 \L^2 \to \L^2 \oplus_2 \L^2} 
\leq \norm{\E \L_{\partial_{\psi,q}(\lambda_s)} }_{2 \to 2} \leq \norm{b_\psi(s)}_H.
\end{equation}
In the case where $-1 < q < 1$ and $b_\psi(s) \neq 0$, this quantity is strictly less than
\begin{align*}
\MoveEqLeft
\norm{\partial_{\psi,q}(\lambda_s)}_{\Gamma_q(H) \rtimes_\alpha G} 
\ov{\eqref{def-partial-psi}}{=} \norm{ s_q(b_\psi(s)) \rtimes \lambda_s }_{\Gamma_q(H) \rtimes_\alpha G} = \norm{s_q(b_\psi(s))}_{\Gamma_q(H)} \\
&\overset{\cite[\textrm{Th. 1.10}]{BKS}}{=} \frac{2}{\sqrt{1-q}} \norm{b_\psi(s)}_H.
\end{align*}
\end{remark}

\begin{remark} \normalfont
\label{rem-Fourier-I-JMP2}
In \cite[p.~588]{JMP2}, the authors obtain a similar estimate to \eqref{Majo-commutateur-group} on the commutator for the case $q = 1$ and $p = 2$.
The differences are that the space on which the commutator is regarded, is the smaller $\L^2(\VN(G)) \oplus_2 \ovl{\Omega_\psi(G)}   \subsetneqq \L^2(\VN(G)) \oplus_2 \L^2(\Gamma_{1}(H) \rtimes_\alpha G)$ \cite[p.~587]{JMP2}, where $\Omega_\psi(G) = \vect\{ s_1(h) \rtimes \lambda_s : \: h \in H_\psi,\: s \in G \}$ and $H_\psi$ is the Hilbert space spanned by the $b_\psi(s),\: s \in G$ \cite[Lemma C.1]{JMP2}.
On the other hand, \cite[p.~588]{JMP2} obtains an exact estimate (without proof)
\begin{equation}
\label{equ-rem-Fourier-I-JMP2} \norm{[\D_{\psi,q=1,p=2},\pi(a)] }_{\L^2 \oplus_2 \ovl{\Omega_\psi(G)} \to \L^2 \oplus_2 \ovl{\Omega_\psi(G)}} 
= \max\big\{ \norm{\Gamma(a,a)}^{\frac12}_{\VN(G)}, \norm{\Gamma(a^*,a^*)}^{\frac12}_{\VN(G)} \big\}.
\end{equation}
Note that when $q = 1$ and $a \in \P_G$, then the right hand side of \eqref{equ-rem-Fourier-I-JMP2} is finite, whereas the right hand side of \eqref{Majo-commutateur-group}, that is, $\norm{\partial_{\psi,q=1,\infty}(a)}_{\Gamma_{1}(H) \rtimes_\alpha G}$ is infinite. We do not know whether \eqref{equ-rem-Fourier-I-JMP2} holds on our full space $\L^2(\VN(G)) \oplus_2 \L^2(\Gamma_q(H) \rtimes_\alpha G)$ or if some intermediate estimate out of \eqref{Majo-commutateur-group} and \eqref{equ-2-rem-Fourier-I-not-optimal} holds in that case.
\end{remark}

\subsection{Spectral triples associated to semigroups of Fourier multipliers II}
\label{Sec-spectral-triples-Fourier-II}

We consider in this subsection a discrete group $G$ and a cocycle $b_\psi \co G \to H$.
Now, we define another ``Hodge-Dirac operator'' by generalizing the construction of \cite[p.~588]{JMP2} which corresponds to the case $q=0$ and $p=2$ below.
Suppose $1 < p<\infty$ and $-1 \leq q < 1$. We let
\begin{equation}
\label{Def-Dirac-group-II}
\mathscr{D}_{\psi,q}(x \rtimes \lambda_s)
\ov{\mathrm{def}}{=} s_q(b_\psi(s))x \rtimes \lambda_s, \quad x \in \Gamma_q(H),s \in G. 
\end{equation}
We can see $\mathscr{D}_{\psi,q}$ as an unbounded operator acting on the subspace $\P_{\rtimes,G}$ of $\L^p(\Gamma_q(H) \rtimes_{\alpha} G)$.

Finally, note that \cite[p.~588]{JMP2} the authors refer to a real structure but we warn the reader that the antilinear isometry $J \co \L^2(\Gamma_q(H) \rtimes_{\alpha} G) \to \L^2(\Gamma_q(H) \rtimes_{\alpha} G)$, $x \mapsto x^*$ used in \cite{JMP2} does not\footnote{\thefootnote. For any $x \in \Gamma_q(H)$ and any $s \in G$, we have
\begin{align*}
\MoveEqLeft
J\mathscr{D}_{\psi,q}(x \rtimes \lambda_s)            
\ov{\eqref{Def-Dirac-free}}{=}J\big(s_q(b_\psi(s))x \rtimes \lambda_s\big)
=\big(s_q(b_\psi(s))x \rtimes \lambda_s\big)^*
\ov{\eqref{inverse-crossed-product}}{=} \alpha_{s^{-1}}(x^*s_q(b_\psi(s))) \rtimes \lambda_{s^{-1}}
\end{align*} 
and
\begin{align*}
\MoveEqLeft
\mathscr{D}_{\psi,q}J(x \rtimes \lambda_s)            
=\mathscr{D}_{\psi,q}\big(\alpha_{s^{-1}}(x^*) \rtimes \lambda_{s^{-1}}\big) 
\ov{\eqref{Def-Dirac-free}}{=} s_q(b_\psi(s^{-1}))\alpha_{s^{-1}}(x^*) \rtimes \lambda_{s^{-1}}
=-s_q(\pi_{s^{-1}}b_\psi(s))\alpha_{s^{-1}}(x^*) \rtimes \lambda_{s^{-1}}\\
&=-\alpha_{s^{-1}}(s_q(b_\psi(s)))\alpha_{s^{-1}}(x^*) \rtimes \lambda_{s^{-1}}
=-\alpha_{s^{-1}}(s_q(b_\psi(s))x^*) \rtimes \lambda_{s^{-1}}.
\end{align*}.} commute or anticommute with the Dirac operator $\mathscr{D}_{\psi,q}$. 

\begin{lemma}
\label{Lemma-closable-group-Dirac-II}
Suppose $1 < p<\infty$ and $-1 \leq q < 1$. 
\begin{enumerate}
	\item For any $a,b \in \P_{\rtimes,G}$ we have
\begin{equation}
\label{equ-1-Lemma-closable-group-Dirac-II}
\big\langle \mathscr{D}_{\psi,q}(a),b \big\rangle 
=\big\langle a , \mathscr{D}_{\psi,q}(b) \big\rangle
\end{equation}
where we use as usual the duality bracket $\langle x,y \rangle = \tau(x^* y)$.

\item The operator $\mathscr{D}_{\psi,q} \co \P_{\rtimes,G} \subseteq \L^p(\Gamma_q(H) \rtimes_{\alpha} G) \to \L^p(\Gamma_q(H) \rtimes_{\alpha} G)$ is closable.
\end{enumerate}
\end{lemma}

\begin{proof}
1. For any $s,t \in G$ and any $x,y \in \Gamma_q(H)$, we have
\begin{align*}
\MoveEqLeft
\big\langle \mathscr{D}_{\psi,q}(x \rtimes \lambda_s), y \rtimes \lambda_t \big\rangle 
\ov{\eqref{Def-Dirac-group-II}}{=} \big\langle s_q(b_\psi(s)) x \rtimes \lambda_s, y \rtimes \lambda_t \big\rangle 
=\tau_{\rtimes} \big((s_q(b_{\psi}(s))x \rtimes \lambda_s)^* (y \rtimes \lambda_t) \big) \\
&\ov{\eqref{inverse-crossed-product}}{=} \tau_{\rtimes} \big((\alpha_{s^{-1}}(x^*s_q(b_{\psi}(s))) \rtimes \lambda_{s^{-1}}) (y \rtimes \lambda_t)\big) 
\ov{\eqref{Product-crossed-product}}{=}\tau_{\rtimes} \big(\alpha_{s^{-1}}(x^*s_q(b_{\psi}(s))) \alpha_{s^{-1}}(y) \rtimes \lambda_{s^{-1}t} \big) \\
&=\tau \big(\alpha_{s^{-1}}(x^*s_q(b_{\psi}(s)) y)\big) \delta_{s=t}
\end{align*}
and
\begin{align*}
\MoveEqLeft
\big\langle x \rtimes \lambda_s, \mathscr{D}_{\psi,q}(y \rtimes \lambda_t) \big\rangle
\ov{\eqref{Def-Dirac-group-II}}{=} \big\langle x \rtimes \lambda_s, s_q(b_\psi(t))y \rtimes \lambda_t \big\rangle
=\tau_{\rtimes}\big((x \rtimes \lambda_s)^*(s_q(b_\psi(t))y \rtimes \lambda_t)\big)\\
&\ov{\eqref{inverse-crossed-product}}{=} \tau_{\rtimes}\big((\alpha_{s^{-1}}(x^*) \rtimes \lambda_{s^{-1}})(s_q(b_\psi(t))y \rtimes \lambda_t)\big)
\ov{\eqref{Product-crossed-product}}{=}\tau_{\rtimes}\big(\alpha_{s^{-1}}(x^*) \alpha_{s^{-1}}(s_q(b_\psi(t))y) \rtimes \lambda_{s^{-1}t}\big)\\
&=\tau_{}\big(\alpha_{s^{-1}}(x^*s_q(b_\psi(s))y)\big)\delta_{s=t}.
\end{align*}
Thus, \eqref{equ-1-Lemma-closable-group-Dirac-II} follows by linearity.

2. Since $\P_{\rtimes,G}$ is dense in $\L^{p^*}(\Gamma_q(H) \rtimes_\alpha G)$, this is a consequence of \cite[Theorem 5.28 p.~168]{Kat1}.
\end{proof}

We denote by $\mathscr{D}_{\psi,q,p} \co \dom \mathscr{D}_{\psi,q,p} \subseteq \L^p(\Gamma_q(H) \rtimes_{\alpha} G) \to \L^p(\Gamma_q(H) \rtimes_{\alpha} G)$ its closure. So $\P_{\rtimes,G}$ is a core of $\mathscr{D}_{\psi,q,p}$. We define the homomorphism $\pi \co \mathrm{C}^*_r(G) \to \B(\L^p(\Gamma_q(H) \rtimes_{\alpha} G))$ by
\begin{equation}
\label{Def-pi-free}
\pi(a)
\ov{\mathrm{def}}{=} 
\Id_{\Gamma_q(H)} \rtimes \L_a, \quad a \in \mathrm{C}^*_r(G)
\end{equation}
where $\L_a$ is the left multiplication operator by $a$. Note that $\pi(a)$ is equal to the map $\tilde{\L}_a$ of Subsection \ref{Sec-spectral-triples-Fourier-I}. 
Recall that using the idea of Lemma \ref{Lemma-continuous-Schur-I} below, it is not difficult to see that $\pi \co \mathrm{C}^*_r(G) \to \B(\L^p(\Gamma_q(H) \rtimes_{\alpha} G))$ is continuous when $\mathrm{C}^*_r(G)$ is equipped with the induced weak* topology of $\VN(G)$ and when $\B(\L^p(\Gamma_q(H) \rtimes_{\alpha} G))$ is equipped with the weak operator topology.

\begin{thm}
\label{thm-spectral-triple-group}
Suppose $1<p< \infty$ and $-1 \leq q < 1$. 
\begin{enumerate}
\item The operator $\mathscr{D}_{\psi,q,2}$ is selfadjoint.
\item Assume that $\L^{p^*}(\VN(G))$ has $\CCAP$ and that $\Gamma_q(H) \rtimes_{\alpha} G$ has $\QWEP$. For $1 < p < \infty$, the unbounded operators $\mathscr{D}_{\psi,q,p}$ and $\mathscr{D}_{\psi,q,p^*}$ are adjoint to each other (with respect to the duality bracket $\langle x,y \rangle = \tau(x^* y)$).

\item We have
\begin{equation}
\label{inclusion-Lip-algebra-II-Fourier}
\P_G
\subseteq \Lip_{\mathscr{D}_{\psi,q,p}}(\mathrm{C}^*_r(G)).
\end{equation}

\item For any $a \in \P_G$, we have
	\begin{equation}
\label{commutator-free-group-petit}
	[\mathscr{D}_{\psi,q,p},\pi(a)]
	=\L_{\partial_{\psi,q}(a)}
\end{equation}
and
\begin{equation}
\label{Equa-norm-commutator-Fourier-II-petit}
\bnorm{[\mathscr{D}_{\psi,q,p},\pi(a)]}_{\L^p(\Gamma_q(H) \rtimes_{\alpha} G) \to \L^p(\Gamma_q(H) \rtimes_{\alpha} G)}
=\bnorm{\partial_{\psi,q}(a)}_{\Gamma_q(H) \rtimes_{\alpha} G}. 
\end{equation} 

\item Suppose that $G$ has $\mathrm{AP}$. We have
\begin{equation}
\label{Lip-algebra-I-Fourier}
\mathrm{C}^*_r(G) \cap \dom \partial_{\psi,q,\infty}
\subseteq \Lip_{\mathscr{D}_{\psi,q,p}}(\mathrm{C}^*_r(G)).
\end{equation}

\item Suppose that $G$ is weakly amenable. We have
\begin{equation}
\label{Lip-algebra-I-Fourier}
\Lip_{\mathscr{D}_{\psi,q,p}}(\mathrm{C}^*_r(G))
=\mathrm{C}^*_r(G) \cap \dom \partial_{\psi,q,\infty}.
\end{equation}

	\item Suppose that $G$ has $\mathrm{AP}$. For any $a \in \mathrm{C}^*_r(G) \cap \dom \partial_{\psi,q,\infty}$, we have
	\begin{equation}
\label{commutator-free-group}
	[\mathscr{D}_{\psi,q,p},\pi(a)]
	=\L_{\partial_{\psi,q,\infty}(a)}.
\end{equation}
\item Suppose that $G$  has $\mathrm{AP}$. For any $a \in \mathrm{C}^*_r(G) \cap \dom \partial_{\psi,q,\infty}$, we have
\begin{equation}
\label{Equa-norm-commutator-Fourier-II}
\bnorm{[\mathscr{D}_{\psi,q,p},\pi(a)]}_{\L^p(\Gamma_q(H) \rtimes_{\alpha} G) \to \L^p(\Gamma_q(H) \rtimes_{\alpha} G)}
=\bnorm{\partial_{\psi,q,\infty}(a)}_{\Gamma_q(H) \rtimes_{\alpha} G}. 
\end{equation} 

\item Suppose $p = 2$ and $q = - 1$. We have $(\mathscr{D}_{\psi,-1,2})^2=\Id \rtimes A_2$.

\item 
If $\Gamma_{-1}(H) \rtimes_{\alpha} G$ has $\QWEP$, $b_\psi \co G \to H$ is injective, $\G_\psi >0$ and if $H$ is finite-dimensional then the operator $|\mathscr{D}_{\psi,-1,2}|^{-1} \co \ovl{\Ran \mathscr{D}_{\psi,-1,2}} \to \ovl{\Ran \mathscr{D}_{\psi,-1,2}}$ is compact.
\end{enumerate}
\end{thm}

\begin{proof}
1. It suffices to show that $\mathscr{D}_{\psi,q} \co \P_{\rtimes,G} \subseteq \L^2(\Gamma_q(H) \rtimes_\alpha G) \to \L^2(\Gamma_q(H) \rtimes_{\alpha} G)$ is essentially selfadjoint. By \eqref{equ-1-Lemma-closable-group-Dirac-II}, we infer that $\mathscr{D}_{\psi,q}$ is symmetric. By \cite[Corollary p.~257]{ReS1}, it suffices to prove that $\Ran(\mathscr{D}_{\psi,q}\pm \i\Id)$ is dense in $\L^2(\Gamma_q(H) \rtimes_\alpha G)$.

Let $z= \sum_s z_s \rtimes \lambda_s \in \L^2(\Gamma_q(H) \rtimes_\alpha G)$ be a vector which is orthogonal to $\Ran(\mathscr{D}_{\psi,q}+ \i\Id)$ in $\L^2(\Gamma_q(H) \rtimes_\alpha G)$. For any $s \in G$ and any $x \in \L^2(\Gamma_q(H))$, we have
\begin{align*}
\MoveEqLeft
\big\langle s_q(b_\psi(s)) x \rtimes \lambda_s, z\big\rangle_{\L^2(\Gamma_q(H) \rtimes_\alpha G)}-\i\tau(x^*z_s) \\
&\ov{\eqref{Def-Dirac-group-II}}{=} \big\langle \mathscr{D}_{\psi,q}(x \rtimes \lambda_s),z \big\rangle_{\L^2(\Gamma_q(H) \rtimes_\alpha G)}-\i\langle x \rtimes \lambda_s,z \rangle_{\L^2(\Gamma_q(H) \rtimes_\alpha G)} \\
&=\big\langle (\mathscr{D}_{\psi,q}- \i\Id)(x \rtimes \lambda_s),z \big\rangle_{\L^2(\Gamma_q(H) \rtimes_\alpha G)} 
=0.            
\end{align*} 

We infer that
$$
\big\langle x,s_q(b_\psi(s))z_s \big\rangle_{\L^2(\Gamma_q(H))}
=\tau( x^* s_q(b_\psi(s))z_s)
=\i \tau(x^*z_s)
=\langle x,\i z_s \rangle_{\L^2(\Gamma_q(H))}.
$$
We deduce that $s_q(b_\psi(s))z_s=\i z_s$. Recall that by \cite[p.~96]{JMX}, the map $\Gamma_q(H) \to \mathcal{F}_q(H)$, $y \mapsto y(\Omega)$ extends to an isometry $\Delta \co \L^2(\Gamma_q(H)) \to \mathcal{F}_q(H)$. For any $x \in \Gamma_q(H)$ and any $y \in \L^2(\Gamma_q(H))$, it is easy to check that $\Delta(xy)=x\Delta(y)$. We infer that 
$$
s_q(b_\psi(s)) \Delta(z_s)
=\Delta(s_q(b_\psi(s)) z_s)
=\Delta(\i z_s)
=\i\Delta(z_s).
$$ 
Thus the vector $\Delta(z_s)$ of $\mathcal{F}_q(H)$ is zero or an eigenvector of $s_q(b_\psi(s))$. Since $s_q(b_\psi(s))$ is selfadjoint, $\i$ is not an eigenvalue. So $\Delta(z_s) = 0$. Since $\Delta$ is injective, we infer that $z_s = 0$. It follows that $z = 0$. The case with $-\i$ instead of $\i$ is similar.

2. By \eqref{symmetry-Schur-Dirac-II} and \cite[Problem 5.24 p.~168]{Kat1}, $\mathscr{D}_{\psi,q,p}$ and $\mathscr{D}_{\psi,q,p^*}$ are formal adjoints. Hence $\mathscr{D}_{\psi,q,p^*} \subseteq (\mathscr{D}_{\psi,q,p})^*$ by \cite[p.~167]{Kat1}. Now, we will show the reverse inclusion. Let $z \in \dom (\mathscr{D}_{\psi,q,p})^*$. For any $y \in \dom \mathscr{D}_{\alpha,q,p}$, we have
\begin{equation}
\label{Divers-897}
\big\langle \mathscr{D}_{\psi,q,p}(y),z \big\rangle_{}
=\big\langle y,(\mathscr{D}_{\psi,q,p})^*(z) \big\rangle_{}.
\end{equation}
Note that $(\Id_{} \rtimes M_{\varphi_j})(z) \to z$ in $\L^{p^*}(\Gamma_q(H) \rtimes_\alpha G)$. Now, for any $y \in \P_{\rtimes,G}$, we have
\begin{align*}
\MoveEqLeft
\big\langle y,\mathscr{D}_{\psi,q,p^*}(\Id_{} \rtimes M_{\varphi_j})(z) \big\rangle_{}
=\big\langle \mathscr{D}_{\psi,q,p}(y),(\Id_{} \rtimes M_{\varphi_j})(z) \big\rangle_{} 
=\big\langle (\Id_{} \rtimes M_{\ovl{\varphi_j}})\mathscr{D}_{\psi,q,p}(y),z \big\rangle_{} \\
&=\big\langle \mathscr{D}_{\psi,q,p}(\Id_{} \rtimes M_{\ovl{\varphi_j}})(y),z \big\rangle_{}\ov{\eqref{Divers-897}}{=} \big\langle (\Id_{} \rtimes M_{\ovl{\varphi_j}})(y),(\mathscr{D}_{\psi,q,p})^*(z) \big\rangle_{} \\
& =\big\langle y,(\Id_{} \rtimes M_{\varphi_j})(\mathscr{D}_{\psi,q,p})^*(z) \big\rangle_{}.
\end{align*}
Hence
$$
\mathscr{D}_{\psi,q,p^*}(\Id_{} \rtimes M_{\varphi_j})(z)
=(\Id_{} \rtimes M_{\varphi_j})(\mathscr{D}_{\psi,q,p})^*(z)
\xra[j]{} (\mathscr{D}_{\psi,q,p})^*(z)
$$
By \eqref{Def-operateur-ferme}, we deduce that $z \in \dom \mathscr{D}_{\psi,q,p^*}$ and that
$$
\mathscr{D}_{\psi,q,p^*}(z)
=(\mathscr{D}_{\psi,q,p})^*(z).
$$

3. and 4. For any $s,t \in G$ and any $x \in \Gamma_q(H)$, using the relation $\alpha_s = \Gamma^\infty_q(\pi_s)$ in the sixth equality, we have
\begin{align*}
\MoveEqLeft
\big[\mathscr{D}_{\psi,q,p},\pi(\lambda_s)\big](x \rtimes \lambda_t)            
=\mathscr{D}_{\psi,q,p}\pi(\lambda_s)(x \rtimes \lambda_t)-\pi(\lambda_s)\mathscr{D}_{\psi,q,p}(x \rtimes \lambda_t) \\
&\ov{\eqref{Def-pi-free} \eqref{Def-Dirac-group-II}}{=} \mathscr{D}_{\psi,q,p}(\alpha_s(x) \rtimes \lambda_{st})-\pi(\lambda_s)(s_q(b_\psi(t)) x \rtimes \lambda_t)\\
&\ov{\eqref{Def-pi-free}}{=} \mathscr{D}_{\psi,q,p}(\alpha_s(x) \rtimes \lambda_{st})-\alpha_s(s_q(b_\psi(t)) x) \rtimes \lambda_{st}\\
&\ov{\eqref{Def-Dirac-group-II}}{=} s_q(b_\psi(st))\alpha_s(x) \rtimes \lambda_{st} -\alpha_s(s_q(b_\psi(t)) x) \rtimes \lambda_{st} \\
&\ov{\eqref{Cocycle-law}}{=} s_q(b_\psi(s))\alpha_s(x) \rtimes \lambda_{st}+s_q(\pi_s(b_\psi(t)))\alpha_s(x) \rtimes \lambda_{st} -\alpha_s(s_q(b_\psi(t)) x) \rtimes \lambda_{st} \\
&=s_q(b_\psi(s))\alpha_s(x) \rtimes \lambda_{st} 
\ov{\eqref{Product-crossed-product}}{=} \big(s_q(b_\psi(s)) \rtimes \lambda_s \big)(x \rtimes \lambda_t) 
\ov{\eqref{def-partial-psi}}{=} \partial_{\psi,q}(\lambda_s)(x \rtimes \lambda_t) \\ 
&=\L_{\partial_{\psi,q}(\lambda_s)} (x \rtimes \lambda_t).
\end{align*}
By linearity and density, we obtain $\big[\mathscr{D}_{\psi,q,p},\pi(\lambda_s)\big]=\L_{\partial_{\psi,q}(\lambda_s)}$. It suffices to use Proposition \ref{Prop-magic-core}. Finally by linearity, for any $a \in \P_G$ we obtain \eqref{commutator-free-group-petit}. In addition, we have
\begin{align}
\label{Equa-divers-12432-Fourier}
\bnorm{[\mathscr{D}_{\psi,q,p},\pi(a)]}_{\L^p(\Gamma_q(H) \rtimes_\alpha G) \to \L^p(\Gamma_q(H) \rtimes_\alpha G)}         
&\ov{\eqref{commutator-free-group-petit}}{=}\bnorm{\L_{\partial_{\psi,q}(a)}}_{\L^p(\Gamma_q(H) \rtimes_\alpha G) \to \L^p(\Gamma_q(H) \rtimes_\alpha G)}\\
&=\bnorm{\partial_{\psi,q}(a)}_{\Gamma_q(H) \rtimes_\alpha G}. \nonumber
\end{align} 

5. We next claim that
\begin{equation}
\label{equ-1-proof-thm-spectral-triple-Fourier-II}
\mathrm{C}^*_r(G) \cap \dom \partial_{\psi,q,\infty} 
\subseteq \Lip_{\mathscr{D}_{\psi,q,p}}(\mathrm{C}^*_r(G)).
\end{equation}
Let $a \in \mathrm{C}^*_r(G) \cap \dom \partial_{\psi,q,\infty}$ and consider a net $(a_j)$ of elements of $\P_G$ such that $a_j \to a$ and $\partial_{\psi,q}(a_j) \to \partial_{\psi,q,\infty}(a)$ both for the weak* topology. By Krein-Smulian Theorem (see Lemma \ref{lem-Krein-Smulian-2}), we can suppose that the nets $(a_j)$ and $(\partial_{\psi,q}(a_j))$ are bounded. So by \eqref{Equa-norm-commutator-Fourier-II-petit}, the net $(\big[\mathscr{D}_{\psi,q,p},\pi(a_j)\big])$ is also bounded. Now by the part 4 of Proposition \ref{Prop-Lip-algebra} and Remark \ref{Remark-123}, $a$ belongs to $\Lip_{\mathscr{D}_{\psi,q,p}}(\mathrm{C}^*_r(G))$. We have shown \eqref{equ-1-proof-thm-spectral-triple-Fourier-II}. 

6. Next we pass to the reverse inclusion of \eqref{equ-1-proof-thm-spectral-triple-Fourier-II} and claim 
$\Lip_{\mathscr{D}_{\psi,q,p}}(\mathrm{C}^*_r(G))
\subseteq \mathrm{C}^*_r(G) \cap \dom \partial_{\psi,q,\infty}$.
To this end, let $a \in \Lip_{\mathscr{D}_{\psi,q,p}}(\mathrm{C}^*_r(G))$ and denote $\hat{a}_r \ov{\mathrm{def}}{=} \langle a, \lambda_r \rangle = \tau(a \lambda_{r^{-1}})$ 
its Fourier coefficients for $r \in G$.
We define a linear form $T_a \co \P_{\rtimes,G} \to \C$ by 
$$
T_a(x \rtimes \lambda_s) 
\ov{\mathrm{def}}{=} \sum_{r \in G} \tau_\rtimes\big(( \hat{a}_r s_q(b_\psi(r)) \rtimes \lambda_r) \cdot (x \rtimes \lambda_s)\big).
$$
Since the trace vanishes for $r \neq s^{-1}$, the sum over $r$ is finite. We will show that it extends to a bounded linear form on $\L^p(\Gamma_q(H) \rtimes_\alpha G)$. To this end consider the bounded net $(a_j)$ in $\P_{\rtimes,G}$ defined by $a_j \ov{\mathrm{def}}{=} M_{\varphi_j}(a)$, where $(\varphi_j)$  is the approximating net guaranteed by the fact that $G$ has AP.
By Proposition \ref{Prop-weak-amenable-crossed}, we have $a_j \to a$ in the weak* topology.
Then for any $\xi,\eta \in \P_{\rtimes,G}$, using the point 2 in the first equality, we obtain
\begin{align*}
\MoveEqLeft
\big\langle [\mathscr{D}_{\psi,q,p},\pi(a)]\xi,\eta \big\rangle 
=\big\langle \pi(a) \xi, \mathscr{D}_{\psi,q,p^*} \eta \big\rangle 
-\big\langle \mathscr{D}_{\psi,q,p} \xi, \pi(a)^* \eta \big\rangle \\
&=\lim_j \big\langle \pi(a_j) \xi, \mathscr{D}_{\psi,q,p^*} \eta \big\rangle - \big\langle \mathscr{D}_{\psi,q,p} \xi, \pi(a_j)^* \eta \big\rangle 
=\lim_j \big\langle [\mathscr{D}_{\psi,q,p}, \pi(a_j)] \xi, \eta \big\rangle \\
&\ov{\eqref{commutator-free-group-petit}}{=} \lim_j \big\langle \partial_{\psi,q}(a_j) \xi, \eta \big\rangle 
=\tau_\rtimes((\partial_{\psi,q}(a_j) \xi)^* \eta) 
=\lim_j \tau_\rtimes(\xi^* \partial_{\psi,q}(a_j)^* \eta) \\
&=\lim_j \ovl{\tau_\rtimes(\eta^* \partial_{\psi,q}(a_j) \xi)} 
=\lim_j \ovl{\tau_\rtimes(\partial_{\psi,q}(a_j) \xi \eta^*)} 
= \ovl{T_a (\xi\eta^*)},
\end{align*}
where in the last equality, we used $\varphi_j(s) \to 1$ for each fixed $s \in G$ and the  definition of $T_a$. Since $a \in \Lip_{\mathscr{D}_{\psi,q,p}}(\mathrm{C}^*_r(G))$, for any $\xi,\eta \in \P_{\rtimes}$, we deduce the estimate $|T_a(\xi\eta^*)| \leq \norm{[\mathscr{D}_{\psi,q,p},\pi(a)]}_{p \to p} \norm{\xi}_p \norm{\eta}_{p^*}$. Choosing the element $\eta=1$ of $\P_{\rtimes,G}$, we get $|T_a(\xi)| \lesssim \norm{\xi}_p$, so that $T_a$ induces an element of $(\L^{p}(\Gamma_q(H)\rtimes_\alpha G))^*$. We infer that the Fourier coefficient sequence $(\hat{a}_r s_q(b_\psi(r)))_{r \in G}$ belongs to an element $b \in \L^{p^*}(\Gamma_q(H)\rtimes_\alpha G)$. But then the above calculation shows that for any $\xi,\eta \in \P_{\rtimes,G}$ 
\begin{equation}
\label{equ-4-proof-thm-spectral-triple-Fourier-II}
\big\langle [\mathscr{D}_{\psi,q,p},\pi(a)]\xi,\eta \big\rangle 
=\lim_j \langle \partial_{\psi,q}(a_j) \xi,\eta \rangle = \langle b \xi, \eta \rangle.
\end{equation}
It follows from $a \in \Lip_{\mathscr{D}_{\psi,q,p}}(\mathrm{C}^*_r(G))$ that $|\langle b \xi, \eta \rangle| \leq \norm{ [\mathscr{D}_{\psi,q,p},\pi(a)]}_{p \to p} \norm{\xi}_p \norm{\eta}_{p^*}$, so that $\norm{b \xi}_p \leq  \norm{ [\mathscr{D}_{\psi,q,p},\pi(a)]}_{p \to p} \norm{\xi}_p$. Thus, the element $b \in \L^{p^*}(\Gamma_q(H)\rtimes_\alpha G)$ is a pointwise multiplier $\L^p \to \L^p$. It follows \footnote{\thefootnote. For $t > 0$, let $\mu_t(b)$ denote the generalised singular number of $b$, so that $\tau_\rtimes(\phi(b)) = \int_0^\infty \phi(\mu_t) dt$ for any continuous function $\R_+ \to \R$ of bounded variation \cite[p.~30]{Xu}.
If $b$ were unbounded, then $\mu_t(b) \to \infty$ as $t \to 0$.
Taking $\xi = \phi(b) \neq 0$ with $\phi$ a smoothed indicator function of an interval $[x,y]$ with large $x$, it is not difficult to see that $\norm{b\xi}_p \geq x \norm{\xi}_p$, which is the desired contradiction.} that $b \in \Gamma_q(H) \rtimes_\alpha G$.

It remains to show that $a \in \dom \partial_{\psi,q,\infty}$. To this end, thanks to the AP property of $G$, we consider again the bounded net $(a_j)$ defined by $a_j \ov{\mathrm{def}}{=} M_{\varphi_j}(a)$. By Proposition \ref{Prop-weak-amenable-crossed}, we have $a_j \to a$ in the weak* topology. Moreover, we have
\begin{align*}
\MoveEqLeft
\norm{\partial_{\psi,q}(a_j)}_{\Gamma_q(H) \rtimes_\alpha G}
=\norm{ \partial_{\psi,q}\bigg(\sum_{r} \hat{a}_r\varphi_j(r)\lambda_r\bigg)}
=\norm{\sum_{r} \hat{a}_r \varphi_j(r)s_q(b_\psi(r)) \rtimes \lambda_r } \\
&=\norm{\sum_{r} \hat{a}_r (\Id \rtimes M_{\varphi_j})(s_q(b_\psi(r)) \rtimes \lambda_r) }
= \norm{(\Id \rtimes M_{\varphi_j})(b)}_{\Gamma_q(H) \rtimes_\alpha G} \\
&\leq \norm{\Id \rtimes M_{\varphi_j}}_{\Gamma_q(H) \rtimes_\alpha G \to \Gamma_q(H) \rtimes_\alpha G} \norm{b}_{\Gamma_q(H) \rtimes_\alpha G}.
\end{align*}
Using Proposition \ref{Prop-weak-amenable-crossed} in the last estimate (since $G$ is weakly amenable), we deduce that $(\partial_{\psi,q}(a_j))$ is a bounded net. Now it follows that the net $(\partial_{\psi,q}(a_j))$ converges for the weak* topology. Indeed, by uniform boundedness, it suffices\footnote{\thefootnote. Recall that if $D$ is a total subset of a Banach space $X$, then a bounded net $(y_j)$ converges to some $y$ of $X^*$ for the weak* topology if and only if $\langle y_j,x\rangle_{X^*,X} \to \langle y, x \rangle_{X^*,X}$ for all $x \in D$.} to test convergence against $x \rtimes \lambda_s$ and we have 
\begin{align*}
\MoveEqLeft
\big\langle \partial_{\psi,q}(a_j), x \rtimes \lambda_s \big\rangle 
=\big\langle a_j, (\partial_{\psi,q})_* (x \rtimes \lambda_s) \big\rangle 
\ov{\eqref{Adjoint-partial-Fourier}}{=} \big\langle s_q(b_\psi(s)), x \big\rangle\langle a_j, \lambda_s \rangle  
\xra[j]{}  
 \big\langle s_q(b_\psi(s)), x \big\rangle \langle a, \lambda_s\rangle.
\end{align*}
Thus, $\partial_{\psi,q}(a_j)$ converges in the weak* topology. 

7. Let $a \in \mathrm{C}^*_r(G) \cap \dom \partial_{\psi,q,\infty}$. Consider a net $(a_j)$ of elements of $\P_G$ such that $a_j \to a$ and $\partial_{\psi,q}(a_j) \to \partial_{\psi,q,\infty}(a)$ both for the weak* topology. By the consequence of the Krein-Smulian Theorem, Lemma \ref{lem-Krein-Smulian-2}, we can suppose that the nets $(a_j)$ and $(\partial_{\psi,q}(a_j))$ are bounded. By continuity of $\pi$, $\pi(a_j) \to \pi(a)$ for the weak operator topology of $\B(\L^p(\Gamma_q(H) \rtimes_\alpha G))$. For any $\xi \in \dom \mathscr{D}_{\psi,q,p}$ and any $\zeta \in \dom \mathscr{D}_{\psi,q,p^*}$, we have
\begin{align*}
\MoveEqLeft
\big\langle [\mathscr{D}_{\psi,q,p},\pi(a_j)]\xi,\zeta \big\rangle_{\L^p(\Gamma_q(H) \rtimes_\alpha G),\L^{p^*}(\Gamma_q(H) \rtimes_\alpha G)}            
=\big\langle (\mathscr{D}_{\psi,q,p}\pi(a_j)-\pi(a_j)\mathscr{D}_{\psi,q,p})\xi, \zeta\big\rangle_{} \\
&=\big\langle \mathscr{D}_{\psi,q,p}\pi(a_j)\xi, \zeta\big\rangle_{}-\big\langle \pi(a_j) \mathscr{D}_{\psi,q,p}\xi, \zeta\big\rangle_{} 
=\big\langle \pi(a_j)\xi, \mathscr{D}_{\psi,q,p^*} \zeta\big\rangle-\big\langle \pi(a_j)\mathscr{D}_{\psi,q,p}\xi, \zeta \big\rangle_{} \\
&\xra[j]{} \big\langle \pi(a)\xi, \mathscr{D}_{\psi,q,p^*}\zeta\big\rangle-\big\langle \pi(a)\mathscr{D}_{\psi,q,p}\xi, \zeta \big\rangle_{} 
=\big\langle [\mathscr{D}_{\psi,q,p},\pi(a)]\xi,\zeta \big\rangle_{}.
\end{align*}
Since the net $([\mathscr{D}_{\psi,q,p},\pi(a_j)])$ is bounded by \eqref{Equa-norm-commutator-Fourier-II-petit}, we deduce that $([\mathscr{D}_{\psi,q,p},\pi(a_j)])$ converges to $[\mathscr{D}_{\psi,q,p},\pi(a)]$ for the weak operator topology by a ``net version'' of \cite[Lemma 3.6 p.~151]{Kat1}. Furthermore, it is (really) easy to check that $\L_{\partial_{\psi,q,\infty}(a_j)} \to \L_{\partial_{\psi,q,\infty}(a)}$ for the weak operator topology of $\B(\L^p(\Gamma_q(H) \rtimes_\alpha G))$. By uniqueness of the limit, we deduce that the commutator is given by the same formula as that in the case of elements of $\P_G$.

8. We obtain \eqref{Equa-norm-commutator-Fourier-II} as in \eqref{Equa-divers-12432-Fourier}.

9. If $x \in \L^2(\Gamma_{-1}(H))$ and $s \in G$, note that 
\begin{align*}
\MoveEqLeft
(\mathscr{D}_{\psi,-1,2})^2(x \rtimes \lambda_s)
\ov{\eqref{Def-Dirac-group-II}}{=} \mathscr{D}_{\psi,-1,2}\big(s_{-1}(b_\psi(s))x \rtimes \lambda_s \big) 
\ov{\eqref{Def-Dirac-group-II}}{=}  s_{-1}(b_\psi(s))^2 x \rtimes \lambda_s \\
&\ov{\eqref{fermion-carre}}{=} x \rtimes \norm{b_\psi(s)}^2_H \lambda_s = x \rtimes \psi(s) \lambda_s
=(\Id_{\L^2(\Gamma_{-1}(H))} \rtimes A_2)(x \rtimes \lambda_s).
\end{align*} 
By \cite[Proposition G.2.4]{HvNVW2}, note that $\P_{\rtimes,G}$ is a core of $\Id_{\L^2(\Gamma_{-1}(H))} \rtimes A_2$. Then since $\mathscr{D}_{\psi,-1,2}^2$ is again selfadjoint, thus closed, we have $\Id_{\L^2(\Gamma_{-1}(H))} \rtimes A_2 \subseteq \mathscr{D}_{\psi,-1,2}^2$. Now it follows from the fact that both $T_1 \ov{\mathrm{def}}{=} \Id_{\L^2(\Gamma_{-1}(H))} \rtimes A_2$ and $T_2 \ov{\mathrm{def}}{=} \mathscr{D}_{\psi,-1,2}^2$ are selfadjoint that both operators are in fact equal. Indeed, consider some $\lambda \in \rho(T_1) \cap \rho(T_2)$ (e.g. $\lambda = \i$). Then $T_1-\lambda \subseteq T_2-\lambda$, the operator $T_1-\lambda$ is surjective and $T_2-\lambda$ is injective. Now, it suffices to use the result \cite[p.~5]{Sch1} to conclude that $T_1-\lambda=T_2-\lambda$ and thus $T_1 = T_2$.

10. By Proposition \ref{prop-Fourier-vraiment-compact}, the operator $A_2^{-\frac12} \co \ovl{\Ran A_2} \to \ovl{\Ran A_2}$ is compact. So there exists a sequence $(T_n)$ of finite rank bounded operators $T_n \co \ovl{\Ran A_2} \to \ovl{\Ran A_2}$ which approximate $A_2^{-\frac12}$ in norm. Similarly to the proof of \cite[Theorem 4.4]{JR1}, we may assume, without loss of generality, that the range of each $T_n$ is contained in $\P_G$ (consider a small perturbation if necessary). 
Composing on the right with the orthogonal projection from $\L^2(\VN(G))$ onto $\ovl{\Ran A_2}$, we can see $A_2^{-\frac12}$ and each $T_n$ as operators on $\L^2(\VN(G))$ which is an operator Hilbert space by \cite[p. 139]{Pis7}.
Moreover, $A_2^{-\frac12}$ is then a Fourier multiplier (with symbol $\psi(s)^{-\frac12} \delta_{s \neq e}$).
Since the operator space $\mathrm{OH}$ is homogeneous by \cite[Proposition 7.2 (iii)]{Pis7}, each $T_n$ and $A_2^{-\frac12}$ are completely bounded. Proposition \ref{prop-Fourier-mult-crossed-product} says that $\Id \rtimes A_2^{-1/2}$ is bounded. Using the projection $P_G$ of \cite[Corollary 4.6]{ArK1}, we obtain for any $n$ a completely bounded Fourier multiplier $P_G(T_n) \co \L^2(\VN(G)) \to \L^2(\VN(G))$ which has finite rank by the proof of \cite[Theorem 4.4]{JR1}. Moreover, by the contractivity of $P_G$ we have for any $n$
\begin{align*}
\MoveEqLeft
\bnorm{A_2^{-1/2}-P_G(T_n)}_{\cb,\L^2(\VN(G)) \to \L^2(\VN(G))}
=\bnorm{P_G(A_2^{-1/2}-T_n)}_{\cb,\L^2(\VN(G)) \to \L^2(\VN(G))} \\
&\leq \bnorm{A_2^{-1/2}-T_n}_{\cb,\L^2(\VN(G)) \to \L^2(\VN(G))}
=\bnorm{A_2^{-1/2}-T_n}_{\L^2(\VN(G)) \to \L^2(\VN(G))}  
\xra[n \to +\infty]{} 0.            
\end{align*}  
Now using the transference of Proposition \ref{prop-Fourier-mult-crossed-product}, we obtain
\begin{align*}
\MoveEqLeft
\bnorm{\Id \rtimes A_2^{-1/2} - \Id \rtimes  P_G(T_n)}_{\L^2(\Gamma_q(H) \rtimes_{\alpha} G) \to \L^2(\Gamma_q(H) \rtimes_{\alpha} G)}
=\bnorm{\Id \rtimes \big(A_2^{-1/2} -  P_G(T_n)\big)} \\
&\leq \bnorm{A_2^{-1/2} -  P_G(T_n)}_{\cb,\L^2(\VN(G)) \to \L^2(\VN(G))} 
\xra[n \to +\infty]{} 0.            
\end{align*}  
Since $H$ is finite-dimensional, note that the fermionic space $\mathcal{F}_{-1}(H)$ is finite-dimensional. It is clear that $\Id \rtimes P_G(T_n)$ is again a finite rank operator\footnote{\thefootnote. Note that if $(e_k)_{k = 1,\ldots, K}$ is an orthonormal basis of $\L^2(\Gamma_{-1}(H))$ and $(f_l)_{l=1,\ldots,L}$ is an orthonormal basis of $\Ran (P_G(T_n))$, then $(e_k \rtimes f_l)_{k = 1,\ldots,K; l = 1,\ldots,L}$ is an orthonormal basis of $\Ran (\Id \rtimes P_G(T_n))$.}. Hence $\Id \rtimes A_2^{-1/2}$ is a norm-limit of finite rank bounded operators, hence compact. 
Thus, $\Id \rtimes A_2^{-\frac12} \co \ovl{\Ran(\Id \rtimes A_2)} \to \ovl{\Ran(\Id_{\L^2(\Gamma_{-1}(H))} \rtimes A_2)}$ is also compact. Finally, we note that $\ovl{\Ran(\Id \rtimes A_2)} = \ovl{\Ran (\mathscr{D}_{\psi,-1,2})^2} =  \ovl{\Ran \mathscr{D}_{\psi,-1,2}}$, where the last equality follows from the fact that $\mathscr{D}_{\psi,-1,2}$ is selfadjoint, see the first point.
\end{proof}

\subsection{Spectral triples associated to semigroups of Schur multipliers I}
\label{Sec-new-spectral-triples-first}

In this subsection, we consider a markovian semigroup of Schur multipliers on $\B(\ell^2_I)$ with associated gradient $\partial_{\alpha,q,p}$.
Suppose $1<p<\infty$ and $-1 \leq q < 1$. Recall that the (full) Hodge-Dirac operator $\D_{\alpha,q,p}$ with domain $\dom \D_{\alpha,q,p} = \dom \partial_{\alpha,q,p} \oplus \dom (\partial_{\alpha,q,p^*})^*$ is defined in \eqref{equ-full-Hodge-Dirac-operator} by the formula
\begin{equation}
\label{Def-gros-Dirac}
\D_{\alpha,q,p} 
\ov{\mathrm{def}}{=} 
\begin{bmatrix} 
0 & (\partial_{\alpha,q,p^*})^* \\ 
\partial_{\alpha,q,p} & 0 
\end{bmatrix}.
\end{equation}
We will see in the main results of this section (Proposition \ref{First-spectral-triple}, Proposition \ref{prop-locally-compact-spectral-triple-Schur-I}) how this Hodge-Dirac operator gives rise to a Banach spectral triple. Note that the compactness criterion needs particular attention and supplementary assumptions. The Banach spectral triple of this subsection will be locally compact.

Let us turn to the description of the homomorphism $\pi$. For any $a \in \B(\ell^2_I)$, we denote by $\tilde{\L}_a \ov{\mathrm{def}}{=} \Id_{\L^p(\Gamma_q(H)} \ot \L_a \co \L^p(\Gamma_q(H) \otvn \B(\ell^2_I)) \to \L^p(\Gamma_q(H) \otvn \B(\ell^2_I)), \: f \ot e_{ij} \mapsto f \ot ae_{ij}$ the left action.
If $a \in S^\infty_I$, we define the bounded operator $\pi(a) \co S^p_I \oplus_p \L^p(\Gamma_q(H) \otvn \B(\ell^2_I)) \to S^p_I \oplus_p \L^p(\Gamma_q(H) \otvn \B(\ell^2_I))$ by
\begin{equation}
\label{Def-pi-a}
\pi(a)
\ov{\mathrm{def}}{=} 
\begin{bmatrix}
    \L_a & 0  \\
    0 & \tilde{\L}_a  \\
\end{bmatrix}, \quad a \in S^\infty_I
\end{equation}
where $\L_a \co S^p_I \to S^p_I$, $x \mapsto ax $ is the left multiplication operator. 
It is then easy to check that $\pi(a)^* = \pi(a^*)$ in case that $p = 2$.

\begin{lemma}
\label{Lemma-continuous-Schur-I}
Let $1 < p < \infty$ and $-1\leq q \leq 1$.
The map $\pi$ is continuous when $S^\infty_I$ is equipped with the weak topology and when $\B(S^p_I \oplus_p \L^p(\Gamma_q(H) \otvn \B(\ell^2_I)))$ is equipped with the weak* topology.
\end{lemma}

\begin{proof}
Indeed, we show that the map $\tilde{\pi} \co \B(\ell^2_I) \to \B(S^p_I \oplus_p \L^p(\Gamma_q(H) \otvn \B(\ell^2_I)))$, $a \mapsto \begin{bmatrix}
    \L_a & 0  \\
    0 & \tilde{\L}_a  \\
\end{bmatrix}$ is weak* continuous. Let $(a_j)$ be a \textit{bounded} net of $\B(\ell^2_I)$ converging in the weak* topology to $a$. It is obvious that the nets $(\L_{a_j})$ and $(\tilde{\L}_{a_j})$ are bounded. If $x \in S^p_I$ and if $y \in S^{p^*}_I$, we have $\langle \L_{a_j}(x),y \rangle_{S^p_I,S^{p^*}_I}=\tr((a_jx)^*y) = \tr(yx^*a_j^*)\xra[j ]{} \tr(yx^*
a^*)=\langle \L_a(x),y \rangle_{S^p_I,S^{p^*}_I}$ since $yx^* \in S^1_I$. So $(\L_{a_j})$ converges to $\L_{a}$ in the weak operator topology. Since $S^p_I$ is reflexive, the weak operator topology and the weak* topology of $\B(S^p_I)$ coincide on bounded sets. We conclude that $(\L_{a_j})$ converges to $(\L_{a})$ in the weak* topology. If $\sum_k z_k \ot t_k$ and if $\sum_l c_l \ot b_l$ are elements of $\Gamma_q(H) \ot \M_{I,\fin}$ we have 
\begin{align*}
\MoveEqLeft
\bigg\langle \tilde{\L}_{a_j}\bigg(\sum_k z_k \ot t_k\bigg),\sum_l c_l \ot b_l \bigg\rangle_{\L^p(\Gamma_q(H) \otvn \B(\ell^2_I)),\L^{p^*}(\Gamma_q(H) \otvn \B(\ell^2_I))}
=\sum_{k,l}  \tau(z_k^*c_l)\tr((a_jt_k)^*b_l) \\
&\xra[j ]{} \sum_{k,l}  \tau(z_k^*c_l)\tr((a t_k)^*b_l) 
=\bigg\langle \tilde{\L}_{a}\bigg(\sum_k z_k \ot t_k\bigg),\sum_l c_l \ot b_l \bigg\rangle_{\L^p(\Gamma_q(H) \otvn \B(\ell^2_I)),\L^{p^*}(\Gamma_q(H) \otvn \B(\ell^2_I))}.
\end{align*} 
By a ``net version'' of \cite[Lemma 3.6 p.~151]{Kat1}, the \textit{bounded} net $(\tilde{\L}_{a_j})$ converges to $\tilde{\L}_{a}$ in the weak operator topology. Once again, since $\L^p(\Gamma_q(H) \otvn \B(\ell^2_I))$ is reflexive, the weak operator topology and the weak* topology of $\B(\L^p(\Gamma_q(H) \otvn \B(\ell^2_I)))$ coincide on bounded sets. We infer that $(\tilde{\L}_{a_j})$ converges to $(\tilde{\L}_{a})$ in the weak* topology. By \cite[Theorem A.2.5 (2)]{BLM}, we conclude that $\tilde{\pi}$ is weak* continuous.    
\end{proof}

Recall in the following proposition the weak* closed operator 
\[\partial_{\alpha,q,\infty} \co \dom ( \partial_{\alpha,q,\infty} ) \subseteq \Gamma_q(H) \otvn \B(\ell^2_I) \to \Gamma_q(H) \otvn \B(\ell^2_I) \]
from Proposition \ref{Prop-derivation-closable}  point 6.
Note that the latter proposition is applicable since we suppose that $q \neq 1$ in the following.

\begin{prop}
\label{First-spectral-triple}
Let $1<p<\infty$ and $-1 \leq q < 1$. Then the triple 
\[ (S^\infty_I,S^p_I \oplus_p \L^p(\Gamma_q(H) \otvn \B(\ell^2_I)), \D_{\alpha,q,p})\]
satisfies the following properties.
\begin{enumerate}
\item We have $(\D_{\alpha,q,p})^*=\D_{\alpha,q,p^*}$. In particular, the operator $\D_{\alpha,q,2}$ is selfadjoint.

\item We have
\begin{equation}
\label{Lip-algebra-I-Schur}
S^\infty_I \cap \dom \partial_{\alpha,q,\infty}
\subseteq \Lip_{\D_{\alpha,q,p} }(S^\infty_I).
\end{equation}

\item For any $a \in S^\infty_I \cap \dom \partial_{\alpha,q,\infty}$, we have 
\begin{equation}
\label{norm-commutator}
\norm{\big[\D_{\alpha,q,p},\pi(a)\big]}_{S^p_I \oplus_p \L^p(\Gamma_q(H) \otvn \B(\ell^2_I)) \to S^p_I \oplus_p \L^p(\Gamma_q(H) \otvn \B(\ell^2_I))}
\leq \bnorm{\partial_{\alpha,q,\infty}(a)}_{\Gamma_q(H) \otvn \B(\ell^2_I)}.
\end{equation}

\end{enumerate}
\end{prop}

\begin{proof}
1. An element $(z,t)$ of $S^{p^*}_I \oplus_{p^*} \L^{p^*}(\Gamma_q(H) \otvn \B(\ell^2_I))$ belongs to $\dom (\D_{\alpha,q,p})^*$ if and only if there exists $(a,b) \in S^{p^*}_I \oplus_{p^*} \L^{p^*}(\Gamma_q(H) \otvn \B(\ell^2_I))$ such that for any $(x,y) \in \dom \partial_{\alpha,q,p} \oplus \dom (\partial_{\alpha,q,p^*})^*$ we have
\begin{align}
\label{Relation99}
\MoveEqLeft
\left\langle \begin{bmatrix} 
0 & (\partial_{\alpha,q,p^*})^* \\ 
\partial_{\alpha,q,p} & 0 
\end{bmatrix}\begin{bmatrix}
      x  \\
      y  \\
\end{bmatrix},\begin{bmatrix}
      z  \\
      t  \\
\end{bmatrix}\right\rangle            
=\left\langle 
\begin{bmatrix}
    x    \\
    y    \\
\end{bmatrix},\begin{bmatrix}
    a    \\
    b    \\
\end{bmatrix}\right\rangle \\
&\iff \big\langle (\partial_{\alpha,q,p^*})^*(y), z\big\rangle
+\big\langle \partial_{\alpha,q,p}(x),t \big\rangle
=\langle y,b\rangle+\langle x,a\rangle. \nonumber
\end{align} 
If $z \in \dom \partial_{\alpha,q,p^*}$ and if $t \in \dom (\partial_{\alpha,q,p})^*$ the latter holds with $b=\partial_{\alpha,q,p^*}(z)$ and $a=(\partial_{\alpha,q,p})^*(t)$. This proves that $\dom \partial_{\alpha,q,p^*} \oplus \dom (\partial_{\alpha,q,p})^* \subseteq \dom (\D_{\alpha,q,p})^*$ and that 
$$
(\D_{\alpha,q,p})^*(z,t)
=\big((\partial_{\alpha,q,p})^*(t),\partial_{\alpha,q,p^*}(z)\big)
= \begin{bmatrix} 
0 & (\partial_{\alpha,q,p})^* \\ 
\partial_{\alpha,q,p^*} & 0 
\end{bmatrix}
\begin{bmatrix}
    z    \\
    t    \\
\end{bmatrix}
\ov{\eqref{Def-gros-Dirac}}{=} \D_{\alpha,q,p^*}(z,t).
$$
Conversely, if $(z,t) \in \dom (\D_{\alpha,q,p})^*$, choosing $y=0$ in \eqref{Relation99} we obtain $t \in \dom (\partial_{\alpha,q,p})^*$ and taking $x=0$ we obtain $z \in \dom \partial_{\alpha,q,p^*}$. 

2. and 3. By Proposition \ref{Prop-derivation-closable} and Proposition \ref{Prop-core-1}, $\M_{I,\fin}$ and $\L^p(\Gamma_q(H)) \ot \M_{I,\fin}$ are cores of $\partial_{\alpha,q,p^*}$ and $(\partial_{\alpha,q,p^*})^*$. So $\M_{I,\fin} \oplus (\L^p(\Gamma_q(H)) \ot \M_{I,\fin})$ is a core of $\D_{\alpha,q,p}$. For any $a \in \M_{I,\fin}$, we have $\L_a(\M_{I,\fin}) \subseteq \M_{I,\fin}$ and $\tilde{\L}_a(\L^p(\Gamma_q(H)) \ot \M_{I,\fin}) \subseteq \L^p(\Gamma_q(H)) \ot \M_{I,\fin}$. We infer that $\pi(a) \cdot (\M_{I,\fin} \oplus (\L^p(\Gamma_q(H)) \ot \M_{I,\fin})) \subseteq \dom \D_{\alpha,q,p}$. So the condition (a) of the first point of Proposition \ref{Prop-magic-core} is satisfied. 
Note also that $\pi(a)^* \cdot (\M_{I,\fin} \oplus (\L^p(\Gamma_q(H)) \ot \M_{I,\fin})) \subseteq \dom \D_{\alpha,q,p^*}=\dom (\D_{\alpha,q,p})^*$ (condition 2 of Proposition \ref{Prop-magic-core}). 

Let $a \in \M_{I,\fin}$. A simple calculation shows that
\begin{align*}
\MoveEqLeft
\big[\D_{\alpha,q,p},\pi(a)\big]
\ov{\eqref{Def-gros-Dirac}\eqref{Def-pi-a}}{=}\begin{bmatrix} 
0 & (\partial_{\alpha,q,p^*})^* \\ 
\partial_{\alpha,q,p}& 0
\end{bmatrix}\begin{bmatrix}
    \L_a & 0  \\
    0 & \tilde{\L}_a  \\
\end{bmatrix}
-\begin{bmatrix}
    \L_a & 0  \\
    0 & \tilde{\L}_a  \\
\end{bmatrix}
\begin{bmatrix} 
0 & (\partial_{\alpha,q,p^*})^* \\ 
\partial_{\alpha,q,p}& 0
\end{bmatrix}\\
&=\begin{bmatrix}
   0 &  (\partial_{\alpha,q,p^*})^*\tilde{\L}_a \\
   \partial_{\alpha,q,p}\L_a  &0  \\
\end{bmatrix}-\begin{bmatrix}
    0 &  \L_a (\partial_{\alpha,q,p^*})^*\\
    \tilde{\L}_a\partial_{\alpha,q,p} & 0  \\
\end{bmatrix} \\
&=\begin{bmatrix}
    0 &  (\partial_{\alpha,q,p^*})^*\tilde{\L}_a-\L_a (\partial_{\alpha,q,p^*})^* \\
   \partial_{\alpha,q,p}\L_a-\tilde{\L}_a\partial_{\alpha,q,p}  & 0  \\
\end{bmatrix}.
\end{align*}
We calculate the two non-zero components of the commutator. For the lower left corner, if $x \in \M_{I,\fin}$ we have\footnote{\thefootnote. Recall that the term $\partial_{\alpha,q,p}(a)x$ is by definition equal to $\partial_{\alpha,q}(a)(1 \ot x)$.}
\begin{align}
\label{Bon-commutateur}
\MoveEqLeft
(\partial_{\alpha,q,p}\L_a-\tilde{\L}_a\partial_{\alpha,q,p})(x)           
=\partial_{\alpha,q,p}\L_a(x)-\tilde{\L}_a\partial_{\alpha,q,p}(x)
=\partial_{\alpha,q,p}(ax)-a\partial_{\alpha,q,p}(x) \\
&\ov{\eqref{Leibniz-Schur-gradient-mieux}}{=} \partial_{\alpha,q}(a)x 
=\L_{\partial_{\alpha,q}(a)}J(x) \nonumber
\end{align}
where $J \co S^p_I \to \L^p(\Gamma_q(H) \otvn \B(\ell^2_I))$, $x \mapsto 1 \ot x$. For the upper right corner, note that for any $y \in \L^p(\Gamma_q(H)) \ot \M_{I,\fin}$ and any $x \in \M_{I,\fin}$, (we recall that we have the duality brackets $\langle f , g \rangle$ antilinear in the first variable)
\begin{align*}
\MoveEqLeft
\big\langle \big((\partial_{\alpha,q,p^*})^*\tilde{\L}_a-\L_a (\partial_{\alpha,q,p^*})^*\big)(y),x \big\rangle_{}            
=\big\langle (\partial_{\alpha,q,p^*})^*\tilde{\L}_a(y),x \big\rangle -\big\langle \L_a (\partial_{\alpha,q,p^*})^*(y),x \big\rangle_{}\\
&=\big\langle \tilde{\L}_a(y),\partial_{\alpha,q,p^*}(x) \big\rangle -\big\langle  (\partial_{\alpha,q,p^*})^*(y),L_{a^*}(x) \big\rangle_{}
=\big\langle y,\tilde{\L}_{a^*}\partial_{\alpha,q,p^*}(x) \big\rangle -\big\langle y,\partial_{\alpha,q,p^*}\L_{a^*}(x) \big\rangle_{} \\
&=\big\langle y,\tilde{\L}_{a^*}\partial_{\alpha,q,p^*}(x)-\partial_{\alpha,q,p^*}\L_{a^*}(x) \big\rangle_{}
=\big\langle y,a^*\partial_{\alpha,q,p^*}(x)-\partial_{\alpha,q,p^*}(a ^*x) \big\rangle_{} \\
&\ov{\eqref{Leibniz-Schur-gradient-mieux}}{=} - \big\langle y,\partial_{\alpha,q}(a^*)x \big\rangle_{}
=\big\langle y, -{\L}_{\partial_{\alpha,q}(a^*)}(1 \ot x)\big\rangle_{}
= \big\langle y, \L_{(\partial_{\alpha,q}(a))^*}(1 \ot x) \big\rangle_{} \\
&=\big\langle {\L}_{\partial_{\alpha,q}(a)}(y), 1 \ot x\big\rangle
=\big\langle \E{\L}_{\partial_{\alpha,q}(a)}(y),x\big\rangle_{S^p_I,S^{p^*}_I}.
\end{align*} 
Here, $\E \co \L^p(\Gamma_q(H) \otvn \B(\ell^2_I)) \to S^p_I,\: t \ot z \mapsto \tau(t) z$ denotes the canonical normal conditional expectation.
We conclude that
\begin{equation}
\label{Commutateur-etrange}
\big((\partial_{\alpha,q,p^*})^*\tilde{\L}_a-\L_a (\partial_{\alpha,q,p^*})^*\big)(y)
=\E{\L}_{\partial_{\alpha,q}(a)}(y).
\end{equation}
The two non-zero components of the commutator are bounded linear operators on $\M_{I,\fin}$ and on $\L^p(\Gamma_q(H)) \ot \M_{I,\fin}$. We deduce that $\big[\D_{\alpha,q,p},\pi(a)\big]$ is bounded on the core $\M_{I,\fin} \oplus (\L^p(\Gamma_q(H)) \ot \M_{I,\fin})$ of $\D_{\alpha,q,p}$. By Proposition \ref{Prop-magic-core}, this operator extends to a bounded operator on $S^p_I \oplus_p \L^p(\Gamma_q(H) \otvn \B(\ell^2_I))$. Hence $\M_{I,\fin}$ is a subset of $\Lip_{\D_{\alpha,q,p} }(S^\infty_I)$. If $(x,y) \in \dom \D_{\alpha,q,p}$ and $a \in \M_{I,\fin}$, we have in addition
\begin{align}
\label{Divers-5436}
\MoveEqLeft
\bnorm{\big[\D_{\alpha,q,p},\pi(a)\big](x,y)}_{p}
=\norm{\big(\big((\partial_{\alpha,q,p^*})^*\tilde{\L}_a - \L_a (\partial_{\alpha,q,p^*})^* \big)y, \big(\partial_{\alpha,q,p} \L_a - \tilde{\L}_a \partial_{\alpha,q,p} \big)x \big)}_p \\          
&=\big(\norm{\big((\partial_{\alpha,q,p^*})^*\tilde{\L}_a - \L_a (\partial_{\alpha,q,p^*})^* \big)y}_{S^p_I}^p+\norm{\big(\partial_{\alpha,q,p} \L_a - \tilde{\L}_a \partial_{\alpha,q,p} \big)x}_p^p\big)^{\frac{1}{p}}\nonumber \\
&\ov{\eqref{Commutateur-etrange} \eqref{Bon-commutateur}}{=} \big(\norm{\E{\L}_{\partial_{\alpha,q}(a)}(y)}_{S_I^{p}}^p+\norm{\partial_{\alpha,q}(a) J(x)}_{\L^p(\Gamma_q(H) \otvn \B(\ell^2_I))}^p \big)^{\frac1p} 
\nonumber\\
&\leq \norm{\partial_{\alpha,q}(a)}_{\Gamma_q(H) \otvn \B(\ell^2_I)} \norm{(x,y)}_p. \nonumber
\end{align}
We conclude that
\begin{equation}
\label{Majo-commutateur}
\bnorm{\big[\D_{\alpha,q,p},\pi(a)\big]}_{_{S^p_I \oplus_p \L^p(\Gamma_q(H) \otvn \B(\ell^2_I)) \to S^p_I \oplus_p \L^p(\Gamma_q(H) \otvn \B(\ell^2_I))}}
\leq \norm{\partial_{\alpha,q}(a)}_{\Gamma_q(H) \otvn \B(\ell^2_I)}.
\end{equation}

Let $a \in S^\infty_I \cap \dom \partial_{\alpha,q,\infty}$. Let $(a_j)$ be a net in $\M_{I,\fin}$ such that $a_j \to a$ and $\partial_{\alpha,q,\infty}(a_j) \to \partial_{\alpha,q,\infty}(a)$ both for the weak* topology. By Krein-Smulian Theorem, we can suppose that the nets $(a_j)$ and $(\partial_{\alpha,q,\infty}(a_j))$ are bounded. Note that $a_j \to a$ for the weak topology of $S^\infty_I$. By the point 4 of Proposition \ref{Prop-Lip-algebra}, we deduce that $a \in \Lip_{\D_{\alpha,q,p} }(S^\infty_I)$. By continuity of $\pi$, note that $\pi(a_j) \to \pi(a)$ for the weak operator topology. For any $\xi \in \dom \D_{\alpha,q,p}$ and any $\zeta \in \dom (\D_{\alpha,q,p})^*$, we have
\begin{align*}
\MoveEqLeft
\big\langle [\D_{\alpha,q,p},\pi(a_j)]\xi,\zeta \big\rangle_{S^p_I \oplus_p \L^p(\Gamma_q(H)\otvn \B(\ell^2_I)),S^{p^*}_I \oplus_{p^*} \L^{p^*}}            
=\big\langle (\D_{\alpha,q,p}\pi(a_j)-\pi(a_j)\D_{\alpha,q,p})\xi, \zeta\big\rangle_{} \\
&=\big\langle \D_{\alpha,q,p}\pi(a_j)\xi, \zeta\big\rangle_{}-\big\langle \pi(a_j) \D_{\alpha,q,p}\xi, \zeta\big\rangle_{} 
=\big\langle \pi(a_j)\xi, (\D_{\alpha,q,p})^*\zeta\big\rangle-\big\langle \pi(a_j)\D_{\alpha,q,p}\xi, \zeta \big\rangle_{} \\
&\xra[j]{} \big\langle \pi(a)\xi, (\D_{\alpha,q,p})^*\zeta\big\rangle-\big\langle \pi(a)\D_{\alpha,q,p}\xi, \zeta \big\rangle_{} 
=\big\langle [\D_{\alpha,q,p},\pi(a)]\xi,\zeta \big\rangle_{}.
\end{align*}
Since the net $([\D_{\alpha,q,p},\pi(a_j)])$ is bounded by \eqref{Majo-commutateur}, we deduce that $([\D_{\alpha,q,p},\pi(a_j)])$ converges to $[\D_{\alpha,q,p},\pi(a)]$ for the weak operator topology by a ``net version'' of \cite[Lemma 3.6 p.~151]{Kat1}. Furthermore, it is (really) easy to check that $\L_{\partial_{\alpha,q,\infty}(a_j)}J \to \L_{\partial_{\alpha,q,\infty}(a)}J$ and $-\E{\L}_{\partial_{\alpha,q,\infty}(a_j)} \to -\E{\L}_{\partial_{\alpha,q,\infty}(a)}$ both for the weak operator topology. By uniqueness of the limit, we deduce that the commutator is given by the same formula that in the case of elements of $\M_{I,\fin}$. From here, we obtain \eqref{norm-commutator} as in \eqref{Divers-5436}.

\end{proof}

\begin{remark} \normalfont
\label{Remark-open-question}
We do not know if $\Lip_{\D_{\alpha,q,p} }(S^\infty_I)=S^\infty_I \cap \dom \partial_{\alpha,q,\infty}$.
\end{remark}

\begin{remark} \normalfont
\label{rem-Schur-I-not-optimal}
The estimate \eqref{norm-commutator} is in general not optimal.
Indeed, already in the case $p = 2$ and $a = e_{ij} \in \M_{I,\fin} \subseteq S^\infty_I \cap \dom \partial_{\alpha,q,\infty}$ for some $i,j \in I$, we have according to \eqref{Bon-commutateur} and \eqref{Commutateur-etrange},
\begin{align}
\label{equ-1-rem-Schur-I-not-optimal}
\MoveEqLeft
\norm{[\D_{\alpha,q,2},\pi(a)]}_{S^2_I \oplus_2 \L^2(\Gamma_q(H) \otvn \B(\ell^2_I)) \to S^2_I \oplus_2 \L^2(\Gamma_q(H) \otvn \B(\ell^2_I))} \\
& \leq \max\left\{ \norm{\L_{\partial_{\alpha,q}(a)}J}_{S^2_I \to \L^2(\Gamma_q(H) \otvn \B(\ell^2_I))}, \norm{ \E \L_{\partial_{\alpha,q}(a)} }_{\L^2(\Gamma_q(H) \otvn \B(\ell^2_I)) \to S^2_I} \right\} . \nonumber
\end{align}
Note that we have the Hilbert space adjoints $(\L_{\partial_{\alpha,q}(a)}J)^* = J^* \L_{\partial_{\alpha,q}(a)}^* = \E \L_{(\partial_{\alpha,q}(a))^*} = - \E \L_{\partial_{\alpha,q}(a^*)}$.
Thus, in the maximum of \eqref{equ-1-rem-Schur-I-not-optimal}, it suffices to consider the second term.
We have for $x = \sum_{k,l} x_{kl} \ot e_{kl} \in \L^2(\Gamma_q(H) \otvn \B(\ell^2_I))$
\begin{align*}
\MoveEqLeft
\E \L_{\partial_{\alpha,q}(e_{ij})}(x) 
= \E\bigg(\big(s_q(\alpha_i - \alpha_j) \ot e_{ij}\big)\bigg( \sum_{k,l} x_{kl} \ot e_{kl}\bigg)\bigg) 
= \E \bigg(\sum_l s_q(\alpha_i - \alpha_j) x_{jl} \ot e_{il}\bigg) \\
&= \sum_l \tau(s_q(\alpha_i-\alpha_j)x_{jl}) e_{il}.
\end{align*}
Thus,
\begin{align*}
\MoveEqLeft
\norm{ \E \L_{\partial_{\alpha,q}(e_{ij})}(x) }_2^2 = \sum_l \left| \tau(s_q(\alpha_i - \alpha_j)x_{jl}) \right|^2 
\leq \norm{s_q(\alpha_i-\alpha_j)}_{\L^2(\Gamma_q(H))}^2  \sum_l \norm{x_{jl}}_{\L^2(\Gamma_q(H))}^2  \\
& \leq \norm{s_q(\alpha_i-\alpha_j)}_{\L^2(\Gamma_q(H))}^2  \sum_{k,l} \norm{x_{kl}}_{\L^2(\Gamma_q(H))}^2 
= \norm{\alpha_i-\alpha_j}_{H}^2  \norm{x}_{\L^2(\Gamma_q(H) \otvn \B(\ell^2_I))}^2.
\end{align*}
We infer that
$$
\norm{[\D_{\alpha,q,2},\pi(a)]}_{S^2_I \oplus_2 \L^2(\Gamma_q(H) \otvn \B(\ell^2_I)) \to S^2_I \oplus_2 \L^2(\Gamma_q(H) \otvn \B(\ell^2_I))} 
= \norm{\E \L_{\partial_{\alpha,q}(e_{ij})} }_{2 \to 2} \leq \norm{\alpha_i - \alpha_j}_H.
$$
In the case where $-1 < q < 1$ and $\alpha_i - \alpha_j \neq 0$, this quantity is strictly less than
\begin{align*}
\MoveEqLeft
\norm{\partial_{\alpha,q}(e_{ij})}_{\Gamma_q(H) \otvn \B(\ell^2_I)} 
\ov{\eqref{def-delta-alpha}}{=} \norm{ s_q(\alpha_i - \alpha_j) \ot e_{ij} }_{\Gamma_q(H) \otvn \B(\ell^2_I)} 
= \norm{s_q(\alpha_i-\alpha_j)}_{\Gamma_q(H)} \norm{e_{ij}}_{\B(\ell^2_I)} \\
& = \norm{s_q(\alpha_i - \alpha_j)}_{\Gamma_q(H)} \overset{\cite[\textrm{Th. 1.10}]{BKS}}{=} \frac{2}{\sqrt{1-q}} \norm{\alpha_i - \alpha_j}_H.
\end{align*}
\end{remark}

Under additional assumptions on the family $(\alpha_i)_{i \in I}$ of our semigroup, we are able to prove that the triple from Proposition \ref{First-spectral-triple} is a locally compact spectral triple (see Definition \ref{def-locally-compact-spectral-triple}).
 In the following, we shall need the restriction $A_{R_p} = A_p|_{R_p}$ of $A_p$ to a row
$$ 
R_p 
\ov{\mathrm{def}}{=} \ovl{\Span\{e_{0j}:\:j \in I\}}, 
$$
where $0$ is some fixed element of $I$. Note that $A_p$ leaves clearly $R_p$ invariant.

\begin{prop}
\label{prop-Schur-vraiment-compact}
Let $1 < p < \infty$. Assume that $\alpha \co I \to H$ is injective, where $H$ is a Hilbert space of dimension $n \in \N$ and satisfies $\G_\alpha > 0$ where $\G_\alpha$ is defined in \eqref{equ-Delta-A-prime}.

\begin{enumerate}
\item
The operator $A_{R_p}^{-\frac12} \co \ovl{\Ran A_{R_p}} \to \ovl{\Ran A_{R_p}}$ is compact.
\item
Suppose $-1 \leq q \leq 1$. Let $B_{R_p} = \partial_{\alpha,q,p} \partial_{\alpha,q,p^*}^*|_{\ovl{\Ran \partial_{\alpha,q,p}|R_p}} \co \dom B_p \subseteq\ovl{\Ran \partial_{\alpha,q,p}|R_p} \to \ovl{\Ran \partial_{\alpha,q,p}|R_p}$.
Then $B_{R_p}$ is sectorial and injective, and the operator $B_{R_p}^{-\frac12}$ is compact. 
\end{enumerate}
\end{prop}

\begin{proof}
1. We begin by showing that $A_{R_2}^{-\frac12} \co \ovl{\Ran A_{R_2}} \to \ovl{\Ran A_{R_2}}$ is compact. Note that that it is obvious that $A_{R_2}$ is selfadjoint and that the $e_{0i}$'s where $i \in I$, form an orthonormal basis of $R_2$ consisting of eigenvectors of $A_{R_2}$. Thus the $e_{0i}$'s where $i \in I \backslash \{0\}$ form an orthonormal basis of $\ovl{\Ran {A_{R_2}}}$ consisting of eigenvectors of $A_{R_2}^{-\frac12} \co \ovl{\Ran A_{R_2}} \to \ovl{\Ran A_{R_2}}$ associated to the $\norm{\alpha_i - \alpha_0}_H^{-1}$. It suffices to show that  The condition $\G_\alpha > 0$, the injectivity of $\alpha$, together with the finite-dimensionality of $H$ imply that any bounded subset of $H$ meets the $\norm{\alpha_i - \alpha_0}_H$ only for a finite number of $i \in I$. Hence theses eigenvalues $\norm{\alpha_i - \alpha_0}_H^{-1}$vanish at infinity\footnote{\thefootnote. Recall that a family $(x_i)_{i \in I}$ vanishes at infinity means that for any $\epsi >0$, there exists a finite subset $J$ of $I$ such that for any $i \in I-J$ we have $|x_i| \leq \epsi$.}. We have proved the compactness. 

Next we show that $A_{R_\infty}^{-\frac12} \co \ovl{\Ran A_{R_\infty}} \to \ovl{\Ran A_{R_\infty}}$ is bounded where $R_\infty \subseteq S^\infty_I$. Since $(T_t|_{R_\infty})_{t \geq 0}$ is a semigroup with generator $A_{R_\infty}$, it suffices by Lemma \ref{lem-A12-bounded} to establish the bound $\norm{T_t|_{R_\infty}}_{\Ran A_{R_\infty} \to R_\infty} \lesssim \frac{1}{t^d}$ for some $d \geq \frac12$. We conclude with Lemma \ref{LemWS} and restriction. Thus, $A_{R_\infty}^{-\frac12} \co \ovl{\Ran A_{R_\infty}} \to \ovl{\Ran A_{R_\infty}}$ is bounded.

Now, assume that $p > 2$. We will use complex interpolation. Since the resolvent of the operators $A_{R_p}$ are compatible for different values of $p$ (in fact, they are equal since $R_{p_0} = R_{p_1}$ for $1 \leq p_0,p_1 \leq \infty$), the spectral projections \cite[p.~364]{HvNVW2} onto the spaces $\ovl{\Ran A_{R_p}}$ are compatible. Hence, the $\ovl{\Ran A_{R_p}}$'s form an interpolation scale. Observe that $\ovl{\Ran A_{R_2}}$ is a Hilbert space, hence a UMD space. Then we obtain the compactness of $A_{R_p}^{-\frac12} \co \ovl{\Ran A_{R_p}} \to \ovl{\Ran A_{R_p}}$ by means of complex interpolation between a compact and a bounded operator with Theorem \ref{Th-interpolation-Kalton}. 

If $p < 2$, we conclude by duality and Schauder's Theorem \cite[Theorem 3.4.15]{Meg1}, since $\ovl{\Ran A_{R_p}}$ is the dual space of $\ovl{\Ran A_{R_{p^*}}}$ and $A_{R_p}^{-\frac12}$ defined on the first space is the adjoint of $A_{R_{p^*}}^{-\frac12}$ defined on the second space.

2. We will use the shorthand notation $\partial_{p} \ov{\mathrm{def}}{=} \partial_{\alpha,q,p}|R_p$ and $\partial_p^* \ov{\mathrm{def}}{=} (\partial_{\alpha,q,p^*})^*|\ovl{\L^p(\Gamma_q(H)) \ot R_p}^p$. Note that $\partial_p$ and $\partial_p^*$ are again closed and densely defined.

We begin with the case $p = 2$. The operators $\partial_2$ and $\partial_2^*$ are indeed adjoints to each other (with the chosen domain). According to \cite[Corollary 5.6]{Tha1}, the unbounded operators $\partial_{2}^* \partial_{2}|(\ker \partial_{2})^\perp$ and $\partial_{2} \partial_{2}^*|(\ker \partial_{2}^*)^\perp$ are unitarily equivalent.

\vspace{0.2cm}

Note that we have $A_{R_2} = A_2|R_2 = \partial_{2}^* \partial_{2}$. Moreover, according to \cite[Probleme 5.27 p.~168]{Kat1}, $(\ker \partial_{2})^\perp = \ovl{\Ran \partial_{2}^*}$, which in turn equals $\ovl{\Ran (\partial_{\alpha,q,2})^*} \cap R_2  \ov{\text{Prop. }\ref{Prop--liens-ranges}}{=}\ovl{\Ran A_2} \cap R_2 = \ovl{\Ran A_{R_2}}$. By the first part of the proof, $A_{R_2}^{-\frac12} \co \ovl{\Ran A_{R_2}} \to \ovl{\Ran A_{R_2}}$ is compact, hence $\partial_{2}^* \partial_{2}|(\ker \partial_{2})^\perp$ is invertible, and by functional calculus, the above unitary equivalence also holds between $A_{R_2}^{-\frac12}$ on $\ovl{\Ran A_{R_2}}$ and $(\partial_{2} \partial_{2}^*)^{-\frac12}$ on $(\ker \partial_{2}^*)^\perp \ov{\cite[\textrm{p.}~168]{Kat1}}{=} \ovl{\Ran \partial_{2}}$. We infer that $B_{R_2}^{-\frac12}$ is compact.

We turn to the case of general $p$ and start by showing that $B_{R_p}^{-\frac12}$ is bounded. Note that $B_{R_p} = \partial_{\alpha,q,p} \partial_{\alpha,q,p^*}^*|\ovl{\Ran \partial_{p}} \ov{\text{Proposition }\ref{Prop-fundamental}}{=} B_p|\ovl{\Ran \partial_{p}}$. Since resolvents of $B_p$ leave $\ovl{\Ran \partial_{p}}$ invariant and $B_p$ is sectorial by Lemma \ref{lem-Bp-is-sectorial}, by \cite[Proposition 3.2.15]{Ege2}, $B_{R_p}$ is also sectorial. Since $B_p$ is injective according to Lemma \ref{lem-Bp-is-sectorial}, also $B_{R_p}$ is injective by restriction. Then again by \cite[Proposition 3.2.15]{Ege2}, $B_{R_p}^{-\frac12} = B_p^{-\frac12}|\ovl{\Ran \partial_{p}}$. By Lemma \ref{lem-Schur-compact} and restriction, we infer that $B_{R_p}^{-\frac12}$ is bounded. To improve boundedness to compactness, we pick some $2 < p < \infty$ and fix some auxiliary $p < p_0 < \infty$. Finally, it suffices interpolate compactness of $B_{R_2}^{-\frac12}$ and boundedness of $B_{R_{p_0}}^{-\frac12}$ to conclude with Theorem \ref{Th-interpolation-Kalton} that $B_{R_p}^{-\frac12}$ is compact. Details are left to the reader. For the case $p<2$, we use duality.
\end{proof}

\begin{prop}
\label{prop-locally-compact-spectral-triple-Schur-I}
Let $1 <  p < \infty$ and $-1 \leq q < 1$.
Assume that $H$ is finite dimensional, that $\alpha \co I \to H$ is injective and that $\G_\alpha > 0$.
Then $(S^\infty_I \cap \dom \partial_{\alpha,q,\infty},S^p_I \oplus \L^p(\Gamma_q(H) \otvn \B(\ell^2_I)),\mathcal{D}_{\alpha,q,p})$ is a locally compact Banach spectral triple.
In other words, we have the following properties.
\begin{enumerate}
\item $\D_{\alpha,q,p}$ is densely defined and has a bisectorial $\HI$ functional calculus.
\item For any $a \in S^\infty_I \cap \dom \partial_{\alpha,q,\infty}$, we have $a \in \Lip_{\D_{\alpha,q,p}}(S^\infty_I)$.
\item For any $a \in S^\infty_I \cap \dom \partial_{\alpha,q,\infty}$, $\pi(a) (\i \Id + \mathcal{D}_{\alpha,q,p})^{-1}$ and $\pi(a) |\mathcal{D}_{\alpha,q,p}|^{-1}$ are compact operators between the spaces $\ovl{\Ran \D_{\alpha,q,p}} \to S^p_I \oplus \L^p(\Gamma_q(H) \otvn \B(\ell^2_I))$.
Moreover, the same compactness holds for $a \in S^\infty_I$ a generic element.
\end{enumerate}
\end{prop}

\begin{proof}
Points 1. resp. 2. are already contained in Theorem \ref{Thm-full-operator-bisectorial} resp. Proposition \ref{First-spectral-triple}.
We turn to 3.
Note that by bisectorial $\HI$-calculus, we have a bounded operator $f(\D_{\alpha,q,p}) \co \ovl{\Ran \D_{\alpha,q,p}} \to \ovl{\Ran \D_{\alpha,q,p}}$, where $f(\lambda) = \sqrt{\lambda^2} (\i + \lambda)^{-1}$ belongs to $\HI(\Sigma_\omega^\pm)$.
Thus, recalling that $|\D_{\alpha,q,p}|^{-1}$ is the functional calculus of the function $n(\lambda) = \frac{1}{\sqrt{\lambda^2}}$, we obtain
$$
\pi(a) (\i \Id + \D_{\alpha,q,p})^{-1} 
= \pi(a)|\D_{\alpha,q,p}|^{-1} f(\D_{\alpha,q,p})
$$
as operators $\ovl{\Ran \D_{\alpha,q,p}} \to S^p_I \oplus \L^p(\Gamma_q(H) \otvn \B(\ell^2_I))$. By composition, if $\pi(a) |\D_{\alpha,q,p}|^{-1}$ is compact on $\ovl{\Ran \D_{\alpha,q,p}}$, then so is $\pi(a) (\i \Id + \D_{\alpha,q,p})^{-1}$. It thus suffices to consider the former in the sequel.

Recall that $\ovl{\Ran A_p}=\ovl{\Ran (\partial_{\alpha,q,p^*})^*}$ by Proposition \ref{Prop--liens-ranges} and that as operators on $\ovl{\Ran \D_{\alpha,q,p}}$,
$$
\D_{\alpha,q,p}^2
\ov{\eqref{D-alpha-carre-egal}}{=}
\begin{bmatrix} 
A_p|\ovl{\Ran A_p}& 0 \\ 
0 & (\Id_{\L^p(\Gamma_q(H))} \ot A_p)|\ovl{\Ran \partial_{\alpha,q,p}}
\end{bmatrix}.
$$
So as operators on $\ovl{\Ran \D_{\alpha,q,p}}$,
\begin{align*}
\MoveEqLeft
|\D_{\alpha,q,p}|^{-1}
=\begin{bmatrix} 
A_p|\ovl{\Ran  A_p} & 0 \\ 
0 & (\Id_{\L^p(\Gamma_q(H))} \ot A_p)|\ovl{\Ran \partial_{\alpha,q,p}}
\end{bmatrix}^{-\frac{1}{2}}
=\begin{bmatrix} 
A_p^{-\frac{1}{2}} & 0 \\ 
0 & B_p^{-\frac{1}{2}}
\end{bmatrix}.            
\end{align*} 
Hence we have
$$
\pi(a)|\D_{\alpha,q,p}|^{-1}
\ov{\eqref{Def-pi-a}}{=} \begin{bmatrix}
    \L_a & 0  \\
    0 & \tilde{\L}_a  \\
\end{bmatrix}\begin{bmatrix} 
A_p^{-\frac{1}{2}} & 0 \\ 
0 & B_p^{-\frac{1}{2}}
\end{bmatrix}
=\begin{bmatrix} 
\L_a A_p^{-\frac{1}{2}} & 0 \\ 
0 & \tilde{\L}_a B_p^{-\frac{1}{2}}
\end{bmatrix}.
$$
It suffices now to show that $\L_a A_p^{-\frac12} \co \ovl{\Ran A_p} \to S^p_I$ and $\tilde{\L}_a B_p^{-\frac12} \co \ovl{\Ran \partial_{\alpha,q,p}} \to \L^p(\Gamma_q(H) \otvn \B(\ell^2_I))$ are compact. We start with the first operator. Recall that $A_p(\M_{I,\fin})$ is a dense subspace of $\ovl{\Ran A_p}$. If $P_j \co S^p_I \to S^p_I$ is the Schur multiplier projecting onto the $j$th line associated with the matrix $[\delta_{j=k}]_{kl}$\footnote{\thefootnote. The entries are 1 on the $j$-row and zero anywhere else.}, if we choose first $a = e_{ij}$, and if $e_{kl} \in A_p(\M_{I,\fin})$ we have
\begin{align*}
\MoveEqLeft
\L_a A_p^{-\frac12}(e_{kl}) 
=\L_{e_{ij}} A_p^{-\frac12}(e_{kl}) 
= \norm{\alpha_k - \alpha_l}^{-1} e_{ij} e_{kl} \\
&= \delta_{j=k} \norm{\alpha_k - \alpha_l}^{-1}e_{il} 
=\delta_{j=k}\L_{e_{ij}} A_p^{-\frac12}(e_{kl}) 
=\L_{a} A_p^{-\frac12}P_j(e_{kl}).
\end{align*}
Then by span density of such $e_{kl}$ in $\ovl{\Ran A_p}$, we infer that $\L_a A_p^{-\frac12} = \L_a A_p^{-\frac12} P_j$ as bounded operators $\ovl{\Ran A_p} \to S^p_I$.
From Lemma \ref{lem-block-Schur}, we infer that $P_j$ is a completely contractive projection. It now suffices to show that $\L_a A_p^{-\frac12}P_j \co P_j(\ovl{\Ran A_p}) \to S^p_I$ is compact. Recall the space $R_p$ from Proposition \ref{prop-Schur-vraiment-compact}, where we choose $0 = j$ there. Since $A_p^{-\frac12}P_j = A_p^{-\frac12}|\ovl{\Ran A_{R_p}}$, we infer by Proposition \ref{prop-Schur-vraiment-compact} that $A_p^{-\frac12}P_j$ is indeed compact on $P_j(\ovl{\Ran A_p}) = \ovl{\Ran A_{R_p}}$. We infer by composition that $\L_a A_p^{-\frac12} P_j = \L_a A_p^{-\frac12}$ is compact as an operator on $\ovl{\Ran A_p}$, in case $a = e_{ij}$.
By linearity in $a$, $\L_a A_p^{-\frac12}$ is also compact for $a \in \M_{I,\fin}$.
Finally, if $a \in S^\infty_I$ is a generic element, then $a = \lim_{n \to \infty} a_n$ in $S^\infty_I$ for some sequence $a_n \in \M_{I,\fin}$, whence $\L_a = \lim_{n \to \infty}\L_{a_n}$ in $\B(S^p_I)$ and $\L_a A_p^{-\frac12} = \lim_{n \to \infty} \L_{a_n} A_p^{-\frac12}$ in $\B(\ovl{\Ran A_p},S^p_I)$.
Thus, $\L_a A_p^{-\frac12}$ is the operator norm limit of compact operators, and hence itself compact.

We turn to the second operator $\tilde{\L}_a B_p^{-\frac12}$. Let $x \ot e_{kl} \in \partial_{\alpha,q,p}(\M_{I,\fin})$ which is a dense subspace of $\ovl{\Ran \partial_{\alpha,q,p}}$ according to Proposition \ref{Prop-derivation-closable}.
We first fix $a = e_{ij}$ and turn to general $a \in S^\infty_I$ at the end.
\begin{align*}
\MoveEqLeft
\tilde{\L}_a B_p^{-\frac12}(x \ot e_{kl}) 
=\norm{\alpha_k - \alpha_l}^{-1} (1 \ot e_{ij}) (x \ot e_{kl}) 
=\delta_{j=k} \norm{\alpha_k - \alpha_l}^{-1} x \ot e_{il} \\
&=\delta_{j=k}(1 \ot e_{ij}) B_p^{-\frac12} (x \ot e_{kl})
=\tilde{\L}_a B_p^{-\frac12} (\Id_{\L^p} \ot P_j)(x \ot e_{kl}).
\end{align*}
By the same arguments as above, we infer that $\tilde{\L}_a B_p^{-\frac12} = \tilde{\L}_a B_p^{-\frac12} (\Id_{\L^p} \ot P_j)$ as bounded operators $\ovl{\Ran \partial_{\alpha,q,p}} \to \L^p(\Gamma_q(H) \otvn \B(\ell^2_I))$.
We infer by Proposition \ref{prop-Schur-vraiment-compact} that $B_p^{-\frac12}(\Id_{\L^p} \ot P_j)$ is compact, thus by composition, also $\tilde{\L}_a B_p^{-\frac12} = \tilde{\L}_a B_p^{-\frac12} (\Id_{\L^p} \ot P_j)$ is compact on $\ovl{\Ran \partial_{\alpha,q,p}}$.
Then by linearity (resp. operator norm limit), $\tilde{\L}_a B_p^{-\frac12}$ is also compact for any $a \in \M_{I,\fin}$ (resp. for any $a \in S^\infty_I$).
We finally deduce that $\pi(a) |\D_{\alpha,q,p}|^{-1}$ is compact on $\ovl{\Ran \D_{\alpha,q,p}}$.
\end{proof}

\begin{remark} \normalfont
\label{Rem-even-Schur-I}
Note that the locally compact (Banach) spectral triple 
\[(\M_{I,\fin},S^p_I \oplus_p \L^p(\Gamma_q(H) \otvn \B(\ell^2_I)),\D_{\alpha,q,p})\]
is even. Indeed, the Hodge-Dirac operator $\D_{\alpha,q,p}$ anti-commutes with the involution $
\gamma_p
\ov{\mathrm{def}}{=}
\begin{bmatrix} 
-\Id_{S^p_I} & 0 \\ 
0& \Id_{\L^p} 
\end{bmatrix}
\co S^p_I \oplus_p  \L^p(\Gamma_q(H) \otvn \B(\ell^2_I)) \to S^p_I \oplus_p  \L^p(\Gamma_q(H) \otvn \B(\ell^2_I)) $ (which is selfadjoint if $p=2$) since
\begin{align*}
\MoveEqLeft
\D_{\alpha,q,p}\gamma_p+\gamma_p\D_{\alpha,q,p}\\
&\ov{\eqref{Def-gros-Dirac}}{=} \begin{bmatrix} 
0 & (\partial_{\alpha,q,p^*})^* \\ 
\partial_{\alpha,q,p}& 0 
\end{bmatrix}
\begin{bmatrix} 
-\Id_{S^p_I} & 0 \\ 
0& \Id_{\L^p}
\end{bmatrix}
+\begin{bmatrix} 
-\Id_{S^p_I} & 0 \\ 
0& \Id_{\L^p}
\end{bmatrix}
\begin{bmatrix} 
0 & (\partial_{\alpha,q,p^*})^* \\ 
\partial_{\alpha,q,p}& 0 
\end{bmatrix}\\
&=\begin{bmatrix} 
0 & (\partial_{\alpha,q,p^*})^* \\ 
-\partial_{\alpha,q,p}&  0
\end{bmatrix}+\begin{bmatrix} 
0 & -(\partial_{\alpha,q,p^*})^* \\ 
\partial_{\alpha,q,p}&  0
\end{bmatrix}
=0.            
\end{align*} 
Moreover, for any $a \in S^\infty_I$, we have
\begin{align*}
\MoveEqLeft
\gamma_p\pi(a)
\ov{\eqref{Def-pi-a}}{=} \begin{bmatrix} 
-\Id_{S^p_I} & 0 \\ 
0& \Id_{\L^p}
\end{bmatrix} 
\begin{bmatrix}
    \L_a & 0  \\
    0 & \tilde{\L}_a  \\
\end{bmatrix}  
=\begin{bmatrix}
    -\L_a & 0  \\
    0 & \tilde{\L}_a  \\
\end{bmatrix}        
=\begin{bmatrix}
    \L_a & 0  \\
    0 & \tilde{\L}_a  \\
\end{bmatrix}\begin{bmatrix} 
-\Id_{S^p_I} & 0 \\ 
0& \Id_{\L^p}
\end{bmatrix}
\ov{\eqref{Def-pi-a}}{=} \pi(a)\gamma_p.
\end{align*} 
\end{remark}

\subsection{Spectral triples associated to semigroups of Schur multipliers II}
\label{subsec-new-spectral-triples-Schur-II}

In this subsection, we shall investigate another triple, defined on the second component of the reflexive Banach space of the preceding subsection. Suppose $1 < p<\infty$ and $-1 \leq q < 1$. We define another ``Hodge-Dirac operator'' by letting
\begin{equation}
\label{Def-Dirac-free}
\mathscr{D}_{\alpha,q}(x \ot e_{ij})
\ov{\mathrm{def}}{=} s_q(\alpha_i-\alpha_j)x \ot e_{ij}, \quad x \in \L^p(\Gamma_q(H)),i,j \in I. 
\end{equation}
We can see $\mathscr{D}_{\alpha,q}$ as an unbounded operator acting on the subspace $\L^p(\Gamma_q(H)) \ot \M_{I,\fin}$ of $\L^p(\Gamma_q(H) \otvn \B(\ell^2_I))$.

\begin{lemma}
\label{lem-Spectral-triple-Schur-II-is-closable}
Suppose $1 < p<\infty$ and $-1 \leq q < 1$. 
\begin{enumerate}
	\item For any $a \in \L^p(\Gamma_q(H)) \ot \M_{I,\fin}$ and $b \in \L^{p^*}(\Gamma_q(H)) \ot \M_{I,\fin}$, we have
\begin{equation}
\label{symmetry-Schur-Dirac-II}
\big\langle \mathscr{D}_{\alpha,q}(a),b \big\rangle 
=\big\langle a , \mathscr{D}_{\alpha,q}(b) \big\rangle
\end{equation}
where we use as usual the duality bracket $\langle x,y \rangle = \tau(x^* y)$.

	\item The operator $\mathscr{D}_{\alpha,q} \co \L^p(\Gamma_q(H)) \ot \M_{I,\fin} \subseteq \L^p(\Gamma_q(H) \otvn \B(\ell^2_I)) \to \L^p(\Gamma_q(H) \otvn \B(\ell^2_I))$ is closable.
\end{enumerate}
\end{lemma}

\begin{proof}
1. For any $x,y \in \L^p(\Gamma_q(H))$ and any $i,j,k,l \in I$, we have 
\begin{align*}
\label{D2-symmetric}
\MoveEqLeft
\big\langle \mathscr{D}_{\alpha,q}(x \ot e_{ij}), y \ot e_{kl}\big\rangle_{}            
\ov{\eqref{Def-Dirac-free}}{=} \big\langle s_q(\alpha_i-\alpha_j)x \ot e_{ij}, y \ot e_{kl}\big\rangle_{} \\
&=\tau\big(x^*s_q(\alpha_i-\alpha_j)y\big) \tr(e_{ij}^* e_{kl}) 
=\delta_{ik}\delta_{jl}\tau\big(x^*s_q(\alpha_k-\alpha_l)y\big)
=\tr(e_{ij}^* e_{kl}) \tau\big(x^* s_q(\alpha_k-\alpha_l)y\big) \nonumber \\
&=\big\langle x \ot e_{ij},s_q(\alpha_k-\alpha_l)y \ot e_{kl} \big\rangle_{}
\ov{\eqref{Def-Dirac-free}}{=} \big\langle x \ot e_{ij},\mathscr{D}_{\alpha,q}(y \ot e_{kl}) \big\rangle_{}. \nonumber
\end{align*}  
Thus \eqref{symmetry-Schur-Dirac-II} follows by linearity.

2. Since $\L^{p^*}(\Gamma_q(H)) \ot \M_{I,\fin}$ is dense in $\L^{p^*}(\Gamma_q(H) \otvn \B(\ell^2_I))$, this is a consequence of \cite[Theorem 5.28 p.~168]{Kat1}.
\end{proof}

We denote by $\mathscr{D}_{\alpha,q,p} \co \dom \mathscr{D}_{\alpha,q,p} \subseteq \L^p(\Gamma_q(H) \otvn \B(\ell^2_I)) \to \L^p(\Gamma_q(H) \otvn \B(\ell^2_I))$ its closure. By definition, the subspace $\L^p(\Gamma_q(H)) \ot \M_{I,\fin}$ is a core of $\mathscr{D}_{\alpha,q,p}$. We define the homomorphism $\pi \co S^\infty_I \to \B(\L^p(\Gamma_q(H) \otvn \B(\ell^2_I)))$ by
\begin{equation}
\label{Def-pi-Schur-II}
\pi(a)
\ov{\mathrm{def}}{=} 
\Id_{\L^p(\Gamma_q(H))} \ot \L_a, \quad a \in S^\infty_I
\end{equation}
where $\L_a$ is the left multiplication by $a$ on $S^p_I$. Note that $\pi(a)$ is equal to the map $\tilde{\L}_a$ of Subsection \ref{Sec-new-spectral-triples-first}. It is not difficult to see that $\pi$ is continuous when $S^\infty_I$ is equipped with the weak topology and when $\B(\L^p(\Gamma_q(H) \otvn \B(\ell^2_I)))$ is equipped with the weak operator topology, see also Lemma \ref{Lemma-continuous-Schur-I}.

\begin{thm}
\label{thm-spectral-triple}
Suppose $1<p< \infty$ and $-1 \leq q < 1$. 
\begin{enumerate}
\item We have $(\mathscr{D}_{\alpha,q,p})^*=\mathscr{D}_{\alpha,q,p^*}$ with respect to the duality bracket $\langle x,y \rangle = \tau(x^* y)$. In particular, the operator $\mathscr{D}_{\alpha,q,2}$ is selfadjoint.
\item We have
\begin{equation}
\label{Lip-algebra-II-Schur}
\Lip_{\mathscr{D}_{\alpha,q,p}}(S^\infty_I)
=S^\infty_I \cap \dom \partial_{\alpha,q,\infty}.
\end{equation}

	\item For any $a \in S^\infty_I \cap \dom \partial_{\alpha,q,\infty}$, we have
	\begin{equation}
\label{commutator-free}
	[\mathscr{D}_{\alpha,q,p},\pi(a)]
	=\L_{\partial_{\alpha,q,\infty}(a)}.
\end{equation}
\item For any $a \in S^\infty_I \cap \dom \partial_{\alpha,q,\infty}$, we have
\begin{equation}
\label{Norm-commutator-Schur-II}
\bnorm{[\mathscr{D}_{\alpha,q,p},\pi(a)]}_{\L^p(\Gamma_q(H) \otvn \B(\ell^2_I)) \to \L^p(\Gamma_q(H) \otvn \B(\ell^2_I))}
=\bnorm{\partial_{\alpha,q,\infty}(a)}_{\Gamma_q(H) \otvn \B(\ell^2_I)}. 
\end{equation}

\item Suppose $p=2$ and $q = -1$. We have $(\mathscr{D}_{\alpha,-1,2})^2=\Id_{\L^2(\Gamma_{-1}(H))} \ot A_2$.

\item 
If $I$ is finite and if $H$ is finite-dimensional
then the operator $|\mathscr{D}_{\alpha,-1,2}|^{-1} \co \ovl{\Ran \mathscr{D}_{\alpha,-1,2}} \to \ovl{\Ran \mathscr{D}_{\alpha,-1,2}}$ is compact\footnote{\thefootnote. Note that since $\Gamma_{-1}(H) \otvn \B(\ell^2_I)$ is then finite-dimensional, the spaces $\L^p(\Gamma_{-1}(H) \otvn \B(\ell^2_I))$ are all isomorphic for different values of $p$, $\ovl{\Ran \mathscr{D}_{\alpha,-1,2}} = \ovl{\Ran \mathscr{D}_{\alpha,-1,p}}$ for all $1 < p < \infty$ and thus, $|\mathscr{D}_{\alpha,-1,2}|^{-1}$ extends to a compact operator on $\ovl{\Ran \mathscr{D}_{\alpha,-1,p}}$ too.}.

\end{enumerate}
\end{thm}

\begin{proof} 
1. By \eqref{symmetry-Schur-Dirac-II} and \cite[Problem 5.24 p.~168]{Kat1}, $\mathscr{D}_{\alpha,q,p}$ and $\mathscr{D}_{\alpha,q,p^*}$ are formal adjoints. Hence $\mathscr{D}_{\alpha,q,p^*} \subseteq (\mathscr{D}_{\alpha,q,p})^*$ by \cite[p.~167]{Kat1}. Now, we will show the reverse inclusion. Let $z \in \dom (\mathscr{D}_{\alpha,q,p})^*$. For any $y \in \dom \mathscr{D}_{\alpha,q,p}$, we have
\begin{equation}
\label{Divers-897-Schur}
\big\langle \mathscr{D}_{\alpha,q,p}(y),z \big\rangle_{}
=\big\langle y,(\mathscr{D}_{\alpha,q,p})^*(z) \big\rangle_{}.
\end{equation}
Note that $(\Id_{\L^p(\Gamma_q(H))} \ot \Tron_J)(z) \to z$. Now, for any $y \in \L^p(\Gamma_q(H)) \ot \M_{I,\fin}$, we have
\begin{align*}
\MoveEqLeft
\big\langle y,\mathscr{D}_{\alpha,q,p^*}(\Id_{\L^p(\Gamma_q(H))} \ot \Tron_J)(z) \big\rangle_{}
=\big\langle \mathscr{D}_{\alpha,q,p}(y),(\Id_{\L^p(\Gamma_q(H))} \ot \Tron_J)(z) \big\rangle_{} \\
&=\big\langle (\Id_{\L^p(\Gamma_q(H))} \ot \Tron_J)\mathscr{D}_{\alpha,q,p}(y),z \big\rangle_{} 
=\big\langle \mathscr{D}_{\alpha,q,p}(\Id_{\L^p(\Gamma_q(H))} \ot \Tron_J)(y),z \big\rangle_{} \\
&\ov{\eqref{Divers-897-Schur}}{=} \big\langle (\Id_{\L^p(\Gamma_q(H))} \ot \Tron_J)(y),(\mathscr{D}_{\alpha,q,p})^*(z) \big\rangle_{} \\           
&=\big\langle y,(\Id_{\L^p(\Gamma_q(H))} \ot \Tron_J)(\mathscr{D}_{\alpha,q,p})^*(z) \big\rangle_{}.
\end{align*}
Hence
$$
\mathscr{D}_{\alpha,q,p^*}(\Id_{\L^p(\Gamma_q(H))} \ot \Tron_J)(z)
=(\Id_{\L^p(\Gamma_q(H))} \ot \Tron_J)(\mathscr{D}_{\alpha,q,p})^*(z)
\to (\mathscr{D}_{\alpha,q,p})^*(z)
$$
By \eqref{Def-operateur-ferme}, we deduce that $z \in \mathscr{D}_{\alpha,q,p^*}$ and that
$$
\mathscr{D}_{\alpha,q,p^*}(z)
=(\mathscr{D}_{\alpha,q,p})^*(z).
$$

2. and 3. For any $i,j,k,l \in I$ and any $x \in \L^p(\Gamma_q(H))$, we have
\begin{align}
\label{Equa-6543}
\MoveEqLeft
\big[\mathscr{D}_{\alpha,q,p},\pi(e_{ij})\big](x \ot e_{kl})            
=\mathscr{D}_{\alpha,q,p}\pi(e_{ij})(x \ot e_{kl})-\pi(e_{ij})\mathscr{D}_{\alpha,q,p}(x \ot e_{kl}) \\
&\ov{\eqref{Def-pi-Schur-II} \eqref{Def-Dirac-free}}{=}\mathscr{D}_{\alpha,q,p}(x \ot e_{ij}e_{kl})-\pi(e_{ij})(s_q(\alpha_k-\alpha_l) x \ot e_{kl})\\\nonumber
&=\delta_{j=k}\mathscr{D}_{\alpha,q,p}(x \ot e_{il})-s_q(\alpha_k-\alpha_l) x \ot e_{ij}e_{kl}\\\nonumber
&\ov{\eqref{Def-Dirac-free}}{=} \delta_{j=k}\big(s_q(\alpha_i-\alpha_l)x \ot e_{il}-s_q(\alpha_k-\alpha_l) x \ot e_{il} \big) \nonumber
=\delta_{j=k}s_q(\alpha_i-\alpha_k)x \ot e_{il}\\
&=s_q(\alpha_i-\alpha_j)x \ot e_{ij}e_{kl}
=\big(s_q(\alpha_i-\alpha_j) \ot e_{ij}\big)(x \ot e_{kl})
\ov{\eqref{def-delta-alpha}}{=} (\partial_{\alpha,q}(e_{ij}))(x \ot e_{kl}) \\
&=\L_{\partial_{\alpha,q}(e_{ij})}(x \ot e_{kl}). \nonumber
\end{align}
By linearity and density, we obtain $\big[\mathscr{D}_{\alpha,q,p},\pi(e_{ij})\big]=\L_{\partial_{\alpha,q}(e_{ij})}$. It suffices to use Proposition \ref{Prop-magic-core}. Finally by linearity, for any $a \in \M_{I,\fin}$ we obtain \eqref{commutator-free} and in addition
\begin{align}
\label{Equa-divers-12432}
\bnorm{[\mathscr{D}_{\alpha,q,p},\pi(a)]}_{\L^p(\Gamma_q(H) \otvn \B(\ell^2_I)) \to \L^p(\Gamma_q(H) \otvn \B(\ell^2_I))}         
&\ov{\eqref{commutator-free}}{=}\bnorm{\L_{\partial_{\alpha,q,\infty}(a)}}_{\L^p(\Gamma_q(H) \otvn \B(\ell^2_I)) \to \L^p(\Gamma_q(H) \otvn \B(\ell^2_I))}\\
&=\bnorm{\partial_{\alpha,q,\infty}(a)}_{\Gamma_q(H) \otvn \B(\ell^2_I)}. \nonumber
\end{align} 
We next claim that
\begin{equation}
\label{equ-1-proof-thm-spectral-triple-Schur-II}
S^\infty_I \cap \dom \partial_{\alpha,q,\infty} 
\subseteq \Lip_{\mathscr{D}_{\alpha,q,p}}(S^\infty_I).
\end{equation}
Let $a \in \dom \partial_{\alpha,q,\infty}$ and consider a net $(a_j)$ of elements of $\M_{I,\fin}$ such that $a_j \to a$ and $\partial_{\alpha,q}(a_j) \to \partial_{\alpha,q,\infty}(a)$ both for the weak* topology. By Krein-Smulian Theorem, we can suppose that the nets $(a_j)$ and $(\partial_{\alpha,q}(a_j))$ are bounded. So by \eqref{Equa-divers-12432}, the net $(\big[\mathscr{D}_{\alpha,q,p},\pi(a_j)\big])$ is also bounded. Now by the part 4 of Proposition \ref{Prop-Lip-algebra}, $a$ belongs to $\Lip_{\mathscr{D}_{\alpha,q,p}}(S^\infty_I)$. We have shown \eqref{equ-1-proof-thm-spectral-triple-Schur-II}. 

Next we claim that for $a \in S^\infty_I \cap \dom \partial_{\alpha,q,\infty}$, we have
\begin{equation}
\label{equ-2-proof-thm-spectral-triple-Schur-II}
\big[\mathscr{D}_{\alpha,q,p},\pi(a)\big] 
= \L_{\partial_{\alpha,q,\infty}(a)}.
\end{equation}
Note that $a_j \to a$ for the weak topology of $S^\infty_I$. By continuity of $\pi$, $\pi(a_j) \to \pi(a)$ for the weak operator topology of $\B(\L^p(\Gamma_q(H) \otvn \B(\ell^2_I)))$. For any $\xi \in \dom \mathscr{D}_{\alpha,q,p}$ and any $\zeta \in \dom \mathscr{D}_{\alpha,q,p^*}$, we have
\begin{align*}
\MoveEqLeft
\big\langle [\mathscr{D}_{\alpha,q,p},\pi(a_j)]\xi,\zeta \big\rangle_{\L^p(\Gamma_q(H) \otvn \B(\ell^2_I)),\L^{p^*}(\Gamma_q(H) \otvn \B(\ell^2_I))}            
=\big\langle (\mathscr{D}_{\alpha,q,p}\pi(a_j)-\pi(a_j)\mathscr{D}_{\alpha,q,p})\xi, \zeta\big\rangle_{} \\
&=\big\langle \mathscr{D}_{\alpha,q,p}\pi(a_j)\xi, \zeta\big\rangle_{}-\big\langle \pi(a_j) \mathscr{D}_{\alpha,q,p}\xi, \zeta\big\rangle_{} 
=\big\langle \pi(a_j)\xi, \mathscr{D}_{\alpha,q,p^*} \zeta\big\rangle-\big\langle \pi(a_j)\mathscr{D}_{\alpha,q,p}\xi, \zeta \big\rangle_{} \\
&\xra[j]{} \big\langle \pi(a)\xi, \mathscr{D}_{\alpha,q,p^*}\zeta\big\rangle-\big\langle \pi(a)\mathscr{D}_{\alpha,q,p}\xi, \zeta \big\rangle_{} 
=\big\langle [\mathscr{D}_{\alpha,q,p},\pi(a)]\xi,\zeta \big\rangle_{}.
\end{align*}
Since the net $([\mathscr{D}_{\alpha,q,p},\pi(a_j)])$ is bounded by \eqref{Equa-divers-12432}, we deduce that $([\mathscr{D}_{\alpha,q,p},\pi(a_j)])$ converges to $[\mathscr{D}_{\alpha,q,p},\pi(a)]$ for the weak operator topology by a ``net version'' of \cite[Lemma 3.6 p.~151]{Kat1}. Furthermore, it is (really) easy to check that $\L_{\partial_{\alpha,q,\infty}(a_j)} \to \L_{\partial_{\alpha,q,\infty}(a)}$ for the weak operator topology of $\B(\L^p(\Gamma_q(H) \otvn \B(\ell^2_I)))$. By uniqueness of the limit, we deduce that the commutator is given by the same formula as that in the case of elements of $\M_{I,\fin}$. From here, we obtain \eqref{Norm-commutator-Schur-II} as in \eqref{Equa-divers-12432}.
 
Next we pass to the reverse inclusion of \eqref{equ-1-proof-thm-spectral-triple-Schur-II} and claim
\begin{equation}
\label{equ-3-proof-thm-spectral-triple-Schur-II}
 \Lip_{\mathscr{D}_{\alpha,q,p}}(S^\infty_I) 
\subseteq S^\infty_I \cap \dom \partial_{\alpha,q,\infty}.
\end{equation}
To this end, let $a \in \Lip_{\mathscr{D}_{\alpha,q,p}}(S^\infty_I)$. Consider the bounded net $(a_J)$ defined by $a_J \ov{\mathrm{def}}{=} \Tron_J(a)$. We have $a_J(a) \to a$ in the weak* topology of $\B(\ell^2_I)$. The essential point will be to prove that
\begin{equation}
\label{equ-4-proof-thm-spectral-triple-Schur-II}
(\partial_{\alpha,q}(a_J)) \text{ is a bounded net.}
\end{equation}
Indeed, once \eqref{equ-4-proof-thm-spectral-triple-Schur-II} is shown, it follows that the net $(\partial_{\alpha,q}(a_J))$ converges for the weak* topology. Indeed, by uniform boundedness, it suffices\footnote{\thefootnote. Recall that if $D$ is a total subset of a Banach space $X$, then a bounded net $(y_J)$ converges to some $y$ of $X^*$ for the weak* topology if and only if $\langle y_J,x\rangle_{X^*,X} \to \langle y, x \rangle_{X^*,X}$ for all $x \in D$.} to test convergence against $f \ot e_{ij}$ and we have 
\begin{align*}
\MoveEqLeft
\big\langle \partial_{\alpha,q}(a_J), f \ot e_{ij} \big\rangle 
\ov{\eqref{Adjoint-partial-Schur}\eqref{Calcul-100000}}{=}  \big\langle s_q(\alpha_i - \alpha_j), f \big\rangle \langle a_J, e_{ij} \rangle
\xra[J]{}  
\big\langle s_q(\alpha_i - \alpha_j), f \big\rangle \langle a, e_{ij}\rangle.
\end{align*}
Thus $a \in \dom \partial_{\alpha,q,\infty}$ and \eqref{equ-3-proof-thm-spectral-triple-Schur-II} follows. 
To show \eqref{equ-4-proof-thm-spectral-triple-Schur-II}, we note that $\partial_{\alpha,q} (a_J)$ belongs to $\Gamma_q(H) \otvn \B(\ell^2_J)$. So
\begin{align*}
\MoveEqLeft
\bnorm{\partial_{\alpha,q} (a_J) }_{\Gamma_q(H) \otvn \B(\ell^2_I)}            
=\bnorm{\partial_{\alpha,q} (a_J) }_{\Gamma_q(H) \otvn \B(\ell^2_J)}
=\norm{\L_{\partial_{\alpha,q} (a_J)}}_{\L^p(\Gamma_q(H) \otvn \B(\ell^2_J)) \to \L^p(\Gamma_q(H) \otvn \B(\ell^2_J))}\\
&=\sup_{\norm{\xi}_{\L^p(\Gamma_q(H) \otvn \B(\ell^2_J))} \leq 1, \: \norm{\eta}_{\L^{p^*}(\Gamma_q(H) \otvn \B(\ell^2_J))} \leq 1} \left| \big\langle \partial_{\alpha,q}(a_J) \xi, \eta \big\rangle \right|.
\end{align*} 
Now, for any $i,j,k,l \in J$, any $x \in \L^p(\Gamma_q(H))$ and any $y \in \L^{p^*}(\Gamma_q(H))$, we have
\begin{align*}
\MoveEqLeft
\big\langle \big[\mathscr{D}_{\alpha,q,p}, \pi(a) \big] (x \ot e_{ij}), y \ot e_{kl} \big\rangle\\
&=\big\langle \pi(a)(x \ot e_{ij}), \mathscr{D}_{\alpha,q,p^*}(y \ot e_{kl})\big\rangle-\big\langle \pi(a)\mathscr{D}_{\alpha,q,p}(x \ot e_{ij}), y \ot e_{kl} \big\rangle_{} \\
&=\big\langle x \ot ae_{ij},s_q(\alpha_k-\alpha_l)y \ot e_{kl}\big\rangle-\big\langle \pi(a)(s_q(\alpha_i-\alpha_j)x \ot e_{ij}), y \ot e_{kl} \big\rangle_{}\\
&=\big\langle x \ot ae_{ij},s_q(\alpha_k-\alpha_l)y \ot e_{kl}\big\rangle-\big\langle s_q(\alpha_i-\alpha_j)x \ot ae_{ij}, y \ot e_{kl} \big\rangle_{}\\
&=\tau\big(x^*s_q(\alpha_k-\alpha_l)y)\tr(e_{ji}a^*e_{kl}\big)-\tau\big(x^*s_q(\alpha_i-\alpha_j) y)\tr(e_{ji}a^*e_{kl}\big)\\
&=\tau\big(x^*s_q(\alpha_k-\alpha_l+\alpha_j-\alpha_i)y)\tr(e_{ji}a^*e_{kl}\big)
=\delta_{j=l} \ovl{a_{ki}}\tau\big(x^*s_q(\alpha_k - \alpha_i)y\big)
\end{align*}
and
\begin{align*}
\MoveEqLeft
\big\langle \big[\mathscr{D}_{\alpha,q,p}, \pi(a_j) \big] (x \ot e_{ij}), y \ot e_{kl} \big\rangle  
=\bigg\langle \bigg(\sum_{r,s \in J} a_{rs} s_q(\alpha_r - \alpha_s) \ot e_{rs}\bigg)(x \ot e_{ij}), y \ot e_{kl} \bigg\rangle \\            
&=\sum_{r,s \in J} \ovl{a_{rs}}\big\langle \big(  s_q(\alpha_r - \alpha_s) \ot e_{rs}\big)(x \ot e_{ij}), y \ot e_{kl} \big\rangle \\  
&=\sum_{r,s \in J} \ovl{a_{rs}} \tau\big(x^*s_q(\alpha_r - \alpha_s)y\big)\tr (e_{ji}e_{sr}e_{kl})\\
&=\sum_{r,s \in J} \delta_{i=s}\delta_{r=k}\delta_{j=l} \ovl{a_{rs}}\tau\big(x^*s_q(\alpha_r - \alpha_s)y\big)
=\delta_{j=l}\ovl{a_{ki}} \tau(x^*s_q(\alpha_k-\alpha_i)y).
\end{align*}   
Fixing $\xi \in \L^p(\Gamma_q(H)) \ot \M_{J}$ and $\eta \in \L^{p^*}(\Gamma_q(H)) \ot \M_{J}$, we deduce that
\begin{align*}
\MoveEqLeft
\big\langle \partial_{\alpha,q}(a_J) \xi,\eta \big\rangle 
\ov{\eqref{commutator-free}}{=} \big\langle \big[ \mathscr{D}_{\alpha,q,p}, \pi(a_J) \big]\xi,\eta \big\rangle 
=\big\langle \big[\mathscr{D}_{\alpha,q,p}, \pi(a) \big]\xi,\eta \big\rangle.
\end{align*}
We conclude by estimating the absolute value of the last term against $\norm{ \big[\mathscr{D}_{\alpha,q,p}, \pi(a) \big] } \norm{\xi}_p \norm{\eta}_{p^*}$.

4. For any $a \in S^\infty_I \cap \dom \partial_{\alpha,q,\infty}$, we have
\begin{align*}
\bnorm{[\mathscr{D}_{\alpha,q,p},\pi(a)]}_{\L^p(\Gamma_q(H) \otvn \B(\ell^2_I)) \to \L^p(\Gamma_q(H) \otvn \B(\ell^2_I))}         
&\ov{\eqref{commutator-free}}{=}\bnorm{\L_{\partial_{\alpha,q,\infty}(a)}}_{\L^p(\Gamma_q(H) \otvn \B(\ell^2_I)) \to \L^p(\Gamma_q(H) \otvn \B(\ell^2_I))}\\
&=\bnorm{\partial_{\alpha,q,\infty}(a)}_{\Gamma_q(H) \otvn \B(\ell^2_I)}.
\end{align*} 

5. If $x \in \L^2(\Gamma_{-1}(H))$ and $i,j \in I$, note that 
\begin{align*}
\MoveEqLeft
(\mathscr{D}_{\alpha,-1,2})^2(x \ot e_{ij})
\ov{\eqref{Def-Dirac-free}}{=} \mathscr{D}_{\alpha,-1,2}\big(s_{-1}(\alpha_i-\alpha_j)x \ot e_{ij}\big) 
\ov{\eqref{Def-Dirac-free}}{=}  s_{-1}(\alpha_i-\alpha_j)^2 x \ot e_{ij} \\
&\ov{\eqref{fermion-carre}}{=} x \ot\norm{\alpha_i-\alpha_j}^2_H e_{ij}
=(\Id_{\L^2(\Gamma_{-1}(H))} \ot A_2)(x \ot e_{ij}).
\end{align*} 
By \cite[Proposition G.2.4]{HvNVW2}, note that $\L^2(\Gamma_{-1}(H)) \ot \M_{I,\fin}$ is a core of $\Id_{\L^2(\Gamma_{-1}(H))} \ot A_2$. Then since $\mathscr{D}_{\alpha,-1,2}^2$ is again selfadjoint, thus closed, we have $\Id_{\L^2(\Gamma_{-1}(H))} \ot A_2 \subseteq \mathscr{D}_{\alpha,-1,2}^2$. Now it follows from the fact that both $T_1 \ov{\mathrm{def}}{=} \Id_{\L^2(\Gamma_{-1}(H))} \ot A_2$ and $T_2 \ov{\mathrm{def}}{=} \mathscr{D}_{\alpha,-1,2}^2$ are selfadjoint that both operators are in fact equal. Indeed, consider some $\lambda \in \rho(T_1) \cap \rho(T_2)$ (e.g. $\lambda = \i$). Then $T_1-\lambda \subseteq T_2-\lambda$, the operator $T_1-\lambda$ is surjective and $T_2-\lambda$ is injective. Now, it suffices to use the result \cite[p.~5]{Sch1} to conclude that $T_1-\lambda=T_2-\lambda$ and thus $T_1 = T_2$.

6.
Since $\mathscr{D}_{\alpha,-1,2}$ is selfadjoint according to the first part, also $\mathscr{D}_{\alpha,-1,2}^2$ is selfadjoint.
As $H$ is finite-dimensional, $\Gamma_{-1}(H)$ is finite-dimensional.
As also $I$ is finite, $\L^2(\Gamma_{-1}(H) \otvn \B(\ell^2_I))$ is finite-dimensional.
Thus, $\ovl{\Ran \mathscr{D}_{\alpha,-1,2}} = \Ran \mathscr{D}_{\alpha,-1,2}$ and $\mathscr{D}_{\alpha,-1,2}^2|\ovl{\Ran \mathscr{D}_{\alpha,-1,2}} : \ovl{\Ran \mathscr{D}_{\alpha,-1,2}} \to \ovl{\Ran \mathscr{D}_{\alpha,-1,2}}$ is a bijective bounded operator.
We infer that its inverse $\mathscr{D}_{\alpha,-1,2}^{-2}$ and the square root of it $|\mathscr{D}_{\alpha,-1,2}|^{-1}$ are also bounded operators on the same space.
As $\L^2(\Gamma_{-1}(H) \otvn \B(\ell^2_I))$ is finite dimensional, $|\mathscr{D}_{\alpha,-1,2}|^{-1} : \ovl{\Ran \mathscr{D}_{\alpha,-1,2}} \to \ovl{\Ran \mathscr{D}_{\alpha,-1,2}}$ is compact.

\end{proof}

\begin{remark} \normalfont
In the case $p=2$, we have an alternative proof of the point 1 of the above Theorem \ref{thm-spectral-triple}. It suffices to show that $\mathscr{D}_{\alpha,q} \co \L^2(\Gamma_q(H)) \ot \M_{I,\fin} \subseteq \L^2(\Gamma_q(H) \otvn \B(\ell^2_I)) \to \L^2(\Gamma_q(H) \otvn \B(\ell^2_I))$ is essentially selfadjoint. 

\begin{proof}
By linearity, \eqref{symmetry-Schur-Dirac-II} says that $\mathscr{D}_{\alpha,q}$ is symmetric. By \cite[Corollary p.~257]{ReS1}, it suffices to prove that $\Ran(\mathscr{D}_{\alpha,q}\pm \i\Id)$ is dense in $\L^2(\Gamma_q(H) \otvn \B(\ell^2_I))$.

Let $z=[z_{ij}] \in \L^2(\Gamma_q(H) \otvn \B(\ell^2_I))$ be a vector which is orthogonal to $\Ran(\mathscr{D}_{\alpha,q}+ \i\Id)$ in $\L^2(\Gamma_q(H) \otvn \B(\ell^2_I))$. For any $i,j \in I$ and any $x \in \L^2(\Gamma_q(H))$, we have
\begin{align*}
\MoveEqLeft
\big\langle s_q(\alpha_i-\alpha_j)x \ot e_{ij}, z\big\rangle_{\L^2(\Gamma_q(H) \otvn \B(\ell^2_I))}-\i\tau(x^*z_{ij}) \\
&\ov{\eqref{Def-Dirac-free}}{=} \big\langle \mathscr{D}_{\alpha,q}(x \ot e_{ij}),z \big\rangle_{\L^2(\Gamma_q(H) \otvn \B(\ell^2_I))}-\i\langle x \ot e_{ij},z \rangle_{\L^2(\Gamma_q(H) \otvn \B(\ell^2_I))} \\
&=\big\langle (\mathscr{D}_{\alpha,q}- \i\Id)(x \ot e_{ij}),z \big\rangle_{\L^2(\Gamma_q(H) \otvn \B(\ell^2_I))} 
=0.            
\end{align*} 
We infer that
$$
\big\langle x,s_q(\alpha_i-\alpha_j)z_{ij} \big\rangle_{\L^2(\Gamma_q(H))}
=\tau( x^* s_q(\alpha_i-\alpha_j)z_{ij})
=\i \tau(x^*z_{ij})
=\big\langle x, \i z_{ij} \big\rangle_{\L^2(\Gamma_q(H))}.
$$
We deduce that $s_q(\alpha_i-\alpha_j)z_{ij}=\i z_{ij}$. Recall that by \cite[p.~96]{JMX}, the map $\Gamma_q(H) \to \mathcal{F}_q(H)$, $y \mapsto y(\Omega)$ extends to an isometry $\Delta \co \L^2(\Gamma_q(H)) \to \mathcal{F}_q(H)$. For any $x \in \Gamma_q(H)$ and any $y \in \L^2(\Gamma_q(H))$, it is easy to check that $\Delta(xy)=x\Delta(y)$.  We infer that 
$$
s_q(\alpha_i - \alpha_j) \Delta(z_{ij})
=\Delta(s_q(\alpha_i - \alpha_j) z_{ij})
=\Delta(\i z_{ij})
=\i\Delta(z_{ij}).
$$ 
Thus the vector $\Delta(z_{ij})$ of $\mathcal{F}_q(H)$ is zero or an eigenvector of $s_q(\alpha_i-\alpha_j)$. Since $s_q(\alpha_i-\alpha_j)$ is selfadjoint, $\i$ is not an eigenvalue. So $\Delta(z_{ij}) = 0$. Since $\Delta$ is injective, we infer that $z_{ij} = 0$. It follows that $z = 0$. The case with $-\i$ instead of $\i$ is similar.
\end{proof}
\end{remark}

In the case where $-1<q<1$, we have the following equivalence with the norm of the commutator. 

\begin{lemma}
\label{Norm-commutator}
Suppose $-1<q<1$. For any $x \in \M_{I,\fin}$, we have
\begin{equation}
\label{Equi-deriv-free}
\bnorm{\partial_{\alpha,q}(x)}_{\Gamma_{q}(H) \otvn \B(\ell^2_I)} 
\approx \max \Big\{\bnorm{ \Gamma(x,x)}_{\B(\ell^2_I)}^{\frac12}, \bnorm{ \Gamma(x^*,x^*)}_{\B(\ell^2_I)}^{\frac12} \Big\}.
\end{equation}
\end{lemma}

\begin{proof}
Here, we use an orthonormal basis $(e_k)_{k \in K}$ of $H$. Using the $\mathrm{C}^*$-identity, first note that the noncommutative Khintchine inequalities of \cite[Theorem 4.1]{BoS} can be rewritten under the following form for any element $f=\sum_{k \in K} s_q(e_k) \ot x_k$ where $(x_k)_{k \in K}$ is a finitely supported family of elements of $\B(\ell^2_I)$ 
\begin{align}
\MoveEqLeft
\label{Khint-q-infty}
\norm{f}_{\Gamma_{q}(H) \otvn \B(\ell^2_I)}            
\approx \max \Bigg\{\Bnorm{\sum_{k \in K} x_k^*x_k}_{\B(\ell^2_I)}^{\frac{1}{2}},\Bnorm{\sum_{k \in K} x_kx_k^*}_{\B(\ell^2_I)}^{\frac{1}{2}} \Bigg\}\\
&\ov{\eqref{petit-Wick}}{=} \max \Bigg\{\Bnorm{\sum_{k,j \in  K} \big\langle s_q(e_k),s_q(e_j) \big\rangle_{\L^2(\Gamma_q)} x_k^*x_j}_{\infty}^{\frac{1}{2}},\Bnorm{\sum_{k,j \in K} \big\langle s_q(e_k),s_q(e_j) \big\rangle_{\L^2(\Gamma_q)} x_kx_j^*}_{\infty}^{\frac{1}{2}} \Bigg\} \nonumber \\
&\ov{\eqref{JMX-equa-2.9}}{=}\max\Big\{\bnorm{f}_{\L^\infty(\B(\ell^2_I),\L^2(\Gamma_q)_{c})},\bnorm{f}_{\L^\infty(\B(\ell^2_I),\L^2(\Gamma_q)_{r})}\Big\} \nonumber\\
&\ov{\eqref{Egalite-fantastique}\eqref{Norm-conditional-general}}{=}\max\Big\{\bnorm{\E(f^* f)}_{\B(\ell^2_I)}^{\frac12},\bnorm{\E(f f^*)}_{\B(\ell^2_I)}^{\frac12}\Big\}.
\nonumber
\end{align} 
By weak* density, it is true for any element of $\Gauss_{q,\infty}(\B(\ell^2_I))$. Replacing $f$ by $\partial_{\alpha,q}(x)$ and using \eqref{Equa-Schur-grad-Gamma}, we obtain
\begin{align*}
\bnorm{\partial_{\alpha,q}(x)}_{\Gamma_q(H) \otvn \B(\ell^2_I)}             
&\ov{\eqref{Khint-q-infty}}{\approx} \max\Big\{\bnorm{\E\big((\partial_{\alpha,q}(x))^* \partial_{\alpha,q}(x)\big)}_{\B(\ell^2_I)}^{\frac12},\bnorm{\E\big(\partial_{\alpha,q}(x) (\partial_{\alpha,q}(x))^*\big)}_{\B(\ell^2_I)}^{\frac12}\Big\} \\
&\ov{\eqref{Equa-Schur-grad-Gamma}}{=} \max \Big\{\bnorm{ \Gamma(x,x)}_{\B(\ell^2_I)}^{\frac12}, \bnorm{ \Gamma(x^*,x^*)}_{\B(\ell^2_I)}^{\frac12} \Big\}.
\end{align*} 
\end{proof}

\begin{remark} \normalfont
Suppose $2 \leq p < \infty$. If $x \in \dom A_p^{\frac12}$, we have seen in the proof of Theorem \ref{cor-Riesz-equivalence-Schur} that
\begin{equation}
\label{}
\bnorm{\partial_{\alpha,q,p}(x)}_{\L^p(\Gamma_q(H) \otvn \B(\ell^2_I))}
\approx \max \left\{ \bnorm{\Gamma(x,x)^{\frac12}}_{S^p_I}, \bnorm{\Gamma(x^*,x^*)^{\frac12}}_{S^p_I} \right\}.
\end{equation}
\end{remark}

\subsection{Bisectoriality and functional calculus of the Dirac operator II}
\label{Sec-functional-calculus}

In the four preceding subsections, we have investigated four Hodge-Dirac operators -- one of which defines a compact Banach spectral triple and one defines a locally compact Banach spectral triple in the sense of Definition \ref{Def-Banach-spectral-triple} and Definition \ref{def-locally-compact-spectral-triple} under suitable assumptions. One of the remaining other triples was defined by the Hodge-Dirac operator of type Fourier II and the properties of bisectoriality and functional calculus was missing in order to have a spectral triple. In this subsection, we consider the Hodge-Dirac operators of type Fourier II and Schur II, that is, $\mathscr{D}_{\alpha,1,p}$ and $\mathscr{D}_{\psi,1,p}$ with parameter $q = 1$ and show that they are indeed bisectorial and do admit a bounded $\HI$ functional calculus. We equally obtain partial results for the case $-1 \leq q < 1$ in Remark \ref{rem-Schur-II-finite-I}.

\paragraph{Bisectoriality of Fourier II on $\L^p(\L^\infty(\Omega) \rtimes_\alpha G)$}

Consider the Dirac operator $\mathscr{D}_{\psi,1}(f \rtimes \lambda_s) = \W(b_\psi(s)) f \rtimes \lambda_s$ for $f \in \L^\infty(\Omega)$ and $s \in G$. 
Moreover, consider the subspace $\L^\infty_\mathrm{s}(\Omega) \ov{\mathrm{def}}{=} \{ f \in \L^\infty(\Omega) :\: \W(h)f \in \L^\infty(\Omega) \text{ for any }h \in H\}$ of $\L^\infty(\Omega)$ stable by ``Gaussian multiplication''.

%

\begin{prop}
\label{prop-Fourier-II-group}
\begin{enumerate}
	\item There is a weak* continuous group $(U_t)_{t \in \R}$ of trace preserving $*$-automorphisms $U_t \co \L^\infty(\Omega) \rtimes_\alpha G \to \L^\infty(\Omega) \rtimes_\alpha G$, such that its weak* generator $\i \mathscr{D}_{\psi,1,\infty}$ is an extension of the restriction $\i \mathscr{D}_{\psi,1}$ on $\Span\{ x \rtimes \lambda_s :\: x \in \L^\infty_\mathrm{s}(\Omega), \: s \in G \}$.

\item Let $1 \leq p < \infty$. Suppose that $\L^p(\VN(G))$ has $\CCAP$ and that $\L^\infty(\Omega) \rtimes_\alpha G$ has $\QWEP$. The operator $\i \mathscr{D}_{\psi,1} \co \P_{\rtimes,G} \subseteq \L^p(\L^\infty(\Omega) \rtimes_\alpha G) \to \L^p(\L^\infty(\Omega) \rtimes_\alpha G)$ is closable and its closure $\i \mathscr{D}_{\psi,1,p}$ generates the strongly continuous group $(U_{t,p})_{t \in \R}$ of isometries on $\L^p(\L^\infty(\Omega) \rtimes_\alpha G)$. 
\end{enumerate}
\end{prop}

\begin{proof}
We consider the weak* continuous group $(U_t)_{t \in \R}$ of trace preserving $*$-automorphisms from \cite[Section 3]{Arh7} given by (an additional factor $\sqrt{2}$ has notational reasons there) $U_t(x \rtimes \lambda_s) \ov{\mathrm{def}}{=}  \e^{\i t\W(b_\psi(s))} x \rtimes \lambda_s$. 
Thus according to \cite[Propositions 1.1.1 and 1.1.2 and Theorem 1.2.3]{Neer3}, we have a weak* closed and weak* densely defined generator $\i \mathscr{D}_{\psi,1,\infty}$ of $(U_t)_{t \in \R}$ and a generator $\i \mathscr{D}_{\psi,1,p}$ for $(U_{t,p})_{t \in \R}$. 

Now, let us show that $\i \mathscr{D}_{\psi,1,\infty}$ contains $\i \mathscr{D}_{\psi,1}$ restricted to $\Span\{ x \rtimes \lambda_s :\: x \in \L^\infty_\mathrm{s}(\Omega), \: s \in G \}$. For that, let $f \in \L^\infty_\mathrm{s}(\Omega)$, $s \in G$ and $y$ be an element of $\L^1(\L^\infty(\Omega) \rtimes_\alpha G)$. Using the differentiation under the integral sign by domination using $f \in \L^\infty_\mathrm{s}(\Omega)$, we obtain
\begin{align*}
\MoveEqLeft
\frac1t \big\langle y, (U_t - \Id)f \rtimes \lambda_s \big\rangle 
=\frac1t \tau\left( y^* (\e^{\i t \W(b_\psi(s))}-1) f \rtimes \lambda_s \right)  \\
&=\frac1t \int_\Omega (\e^{\i t \W(b_\psi(s))(\omega)} - 1)f(\omega) \ovl{y_s(\omega)} \d \mu(\omega) \\
&\xra[t \to 0]{}  \int_\Omega \i \W(b_\psi(s))(\omega) f(\omega) \ovl{y_s(\omega)} \d \mu(\omega) 
=\big\langle y, \i \W(b_\psi(s))f \rtimes \lambda_s \big\rangle. 
\end{align*}
We deduce that $f \rtimes \lambda_s$ belongs to $\dom \i \mathscr{D}_{\psi,1,\infty}$ and that $\i \mathscr{D}_{\psi,1,\infty}(x \rtimes \lambda_s)=\i \W(b_\psi(s))f \rtimes \lambda_s$.

Next we show that the generator $\i \mathscr{D}_{\psi,1,p}$ is the closure of $\i \mathscr{D}_{\psi,1} \co \P_{\rtimes,G} \subseteq \L^p(\L^\infty(\Omega) \rtimes_\alpha G) \to \L^p(\L^\infty(\Omega) \rtimes_\alpha G)$. Let $f \rtimes \lambda_s$ be an element of $\P_{\rtimes,G}$. Using \cite[Proposition 4.11]{EnN1}, we have
\begin{align*}
\MoveEqLeft
\frac1t (U_{t,p} - \Id)(f \rtimes \lambda_s)  
=\frac1t\big(\e^{\i t \W(b_\psi(s))}-1\big)f \rtimes \lambda_s
\ov{\eqref{semigroupe-vers-generateur}}{\xra[t \to 0]{}} 
\i \W(b_\psi(s))f \rtimes \lambda_s. 
\end{align*}
We infer that $f \rtimes \lambda_s$ belongs to $\dom \mathscr{D}_{\psi,1,p}$ and that we have $\i \mathscr{D}_{\psi,1,p}(f \rtimes \lambda_s) = \i \W(b_\psi(s)) f \rtimes \lambda_s$. So it suffices to prove that $\P_{\rtimes,G}$ is a core for $\i \mathscr{D}_{\psi,1,p}$. 

According to Lemma \ref{Lemma-carac-CBAP-Lp-spaces} and Proposition \ref{prop-Fourier-mult-crossed-product}, there exists a net $(\varphi_j)$ of finitely supported functions $G \to \C$ converging pointwise to $1$ such that $\Id_{\L^p(\Omega)} \rtimes M_{\varphi_j}$ converges to $\Id_{\L^p(\L^\infty(\Omega) \rtimes_\alpha G)}$ in the point norm topology and $\norm{\Id_{\L^p(\Omega)} \rtimes M_{\varphi_j}}_{\cb, \L^p \to \L^p} \leq 1$. 

Note that $\Id_{\L^p(\Omega)} \rtimes M_{\varphi_j}$ commutes with $U_{t,p}$. This is easy to see on elements of the form $f \rtimes \lambda_s$ and extends by a density argument to all of $\L^p(\L^\infty(\Omega) \rtimes_\alpha G)$. Hence, for $x \in \dom \mathscr{D}_{\psi,1,p}$, we have 
$$
\frac1t (U_{t,p} - \Id)(\Id_{\L^p(\Omega)} \rtimes M_{\varphi_j})(x) 
= \frac1t (\Id_{\L^p(\Omega)} \rtimes M_{\varphi_j}) (U_{t,p} - \Id)(x) \to \i (\Id_{\L^p(\Omega)} \rtimes M_{\varphi_j}) \mathscr{D}_{\psi,1,p}(x).
$$ 
We infer that $(\Id_{\L^p(\Omega)} \rtimes M_{\varphi_j})(x)$ belongs again to $\dom \mathscr{D}_{\psi,1,p}$ and 
$$
\mathscr{D}_{\psi,1,p}(\Id_{\L^p(\Omega)} \rtimes M_{\varphi_j})(x) 
= (\Id_{\L^p(\Omega)} \rtimes M_{\varphi_j}) \mathscr{D}_{\psi,1,p}(x)
\xra[j ]{} \mathscr{D}_{\psi,1,p} (x). 
$$
Moreover, $(\Id_{\L^p(\Omega)} \rtimes M_{\varphi_j})(x) \to x$ as $j \to \infty$. Therefore, to show the core property, it suffices to approximate an element $f \rtimes \lambda_s$ where $f \in \L^p(\Omega)$ in the graph norm of $\mathscr{D}_{\psi,1,p}$ by elements in $\P_{\rtimes,G}$. To this end, it suffices to take $1_{|f| \leq n} f \rtimes \lambda_s$ and to use the same reasoning, noting that 
$$
\big(1_{|f|\leq n} \rtimes 1_{\VN(G)}\big) U_{t,p}(f \rtimes \lambda_s) 
= U_{t,p}\big( 1_{|f | \leq n} f \rtimes \lambda_s\big).
$$ 
\end{proof}

\begin{cor}
\label{cor-Fourier-II-bisectoriel-et-calcul-fonctionnel}
Let $1 < p < \infty$. Then the operator $\mathscr{D}_{\psi,1,p}$ defined above is bisectorial and has a bounded $\HI(\Sigma_\omega^\pm)$ functional calculus to any angle $\omega > 0$.
\end{cor}

\begin{proof}
This follows at once from Proposition \ref{prop-Fourier-II-group} together with the fact that a noncommutative $\L^p$-space is UMD by \cite[Corollary 7.7]{PiX} together with the Hieber Pr\"uss Theorem \cite[Theorem 10.7.10]{HvNVW2} which says that the generator $iA$ of a bounded strongly continuous group on a $\UMD$-space has a bounded bisectorial $\H^\infty$ functional calculus of angle 0.
\end{proof}

\paragraph{Bisectoriality of Schur II on $\L^p(\L^\infty(\Omega) \otvn \B(\ell^2_I))$}
In the following, we consider the case $q = 1$ which corresponds to classical Gaussian space. Consider the ``Dirac operator'' $\mathscr{D}_{\alpha,1}$ defined on $\L^\infty(\Omega) \ot \M_{I,\fin}$ with values in $\L^0(\Omega) \ot \M_{I,\fin}$ by
$$
\mathscr{D}_{\alpha,1}(f \ot e_{ij}) 
\ov{\mathrm{def}}{=}  \W(\alpha_i - \alpha_j) f \ot e_{ij}, \quad x \in \L^\infty(\Omega), i,j \in I.
$$
Moreover, we consider the subspace $\L^\infty_\mathrm{s}(\Omega)$ as in the last paragraph.
See \cite{Arh4} for a generalization of part 1.


\begin{prop}
\label{prop-Schur-II-group}
\begin{enumerate}
	\item There is a weak* continuous group $(U_t)_{t \in \R}$ of trace preserving $*$-automorphisms $U_t \co \L^\infty(\Omega) \otvn \B(\ell^2_I) \to \L^\infty(\Omega) \otvn \B(\ell^2_I)$, such that its weak* generator $\i \mathscr{D}_{\alpha,1,\infty}$ is an extension of the restriction $\i \mathscr{D}_{\alpha,1}|_{\L^\infty_s(\Omega) \ot \M_{I,\fin}}$.

\item Let $1 \leq p < \infty$. The operator $\i \mathscr{D}_{\alpha,1} \co \L^\infty(\Omega) \ot \M_{I,\fin} \subseteq \L^p(\Omega,S^p_I) \to \L^p(\Omega,S^p_I)$ is closable and its closure $\i \mathscr{D}_{\alpha,1,p}$ generates the strongly continuous group $(U_{t,p})_{t \in \R}$ of isometries on $\L^p(\Omega,S^p_I)$. 
\end{enumerate}
\end{prop}

\begin{proof}
For $t \in \R$, we define the block diagonal operator $V_t \ov{\mathrm{def}}{=} \sum_{k \in I} \e^{\i t \W(\alpha_k)} \ot e_{kk} \in \L^\infty(\Omega) \otvn \B(\ell^2_I)$. Note that $(V_t)_{t \in \R}$ is a group of unitaries. Moreover, $t \mapsto V_t$ is weakly continuous. Indeed, for $\xi, \eta \in \L^2(\Omega) \ot_2 \ell^2_I$ and $t_0 \in \R$, we have using the dominated convergence theorem
$$ 
\langle V_t(\xi),\eta \rangle_{\L^2(\Omega) \ot_2 \ell^2_I} 
= \sum_{k \in I} \big\langle \e^{\i t \W(\alpha_k)} \xi_k , \eta_k \big\rangle_{\L^2(\Omega)} 
\xra[t \to t_0]{} \sum_{k \in I} \langle \e^{\i t_0 \W(\alpha_k)} \xi_k , \eta_k \rangle.
$$
Then by \cite[Lemma 13.4]{Str} or \cite[p.~239]{Tak2}, $t \mapsto V_t$ is even strongly continuous. By \cite[p.~238]{Tak2}, we conclude that $U_t \co \L^\infty(\Omega) \otvn \B(\ell^2_I) \to \L^\infty(\Omega) \otvn \B(\ell^2_I)$, $t \mapsto V_t xV_t^*$ define a weak* continuous group $(U_t)_{t \in \R}$ of $*$-automorphisms. For any $t \in \R$, it is easy to check that $U_t$ is trace preserving. 
So it induces for $1 \leq p < \infty$ a (uniquely) strongly continuous group $(U_{t,p})_{t \in \R}$ of complete isometries on $\L^p(\Omega,S^p_I)$. Thus according to \cite[Propositions 1.1.1 and 1.1.2 and Theorem 1.2.3]{Neer3}, we have a weak* closed and weak* densely defined generator $\i \mathscr{D}_{\alpha,1,\infty}$ of $(U_t)_{t \in \R}$ and a generator $\i \mathscr{D}_{\alpha,1,p}$ for $(U_{t,p})_{t \in \R}$. 

Now, let us show that $\i \mathscr{D}_{\alpha,1,\infty}$ contains $\i \mathscr{D}_{\alpha,1}|_{\L^\infty_\mathrm{s}(\Omega) \ot \M_{I,\fin}}$. To this end, let $f \in \L^\infty_\mathrm{s}(\Omega)$, $i,j \in I$ and $y$ be an element of $\L^1(\Omega,S^1_I)$. Using differentiation under the integral sign by domination (note that $x \in \L^\infty_s(\Omega)$), we obtain
\begin{align}
\MoveEqLeft
\frac1t \big\langle y, (U_t - \Id)(f \ot e_{ij}) \big\rangle_{\L^1(\Omega,S^1_I),\L^\infty(\Omega) \otvn \B(\ell^2_I)} 
=\frac1t \tau\left( y^* ((\e^{\i t \W(\alpha_i-\alpha_j)}-1)f \ot e_{ij}) \right) \nonumber \\
&=\frac1t \int_\Omega \big(\e^{\i t \W(\alpha_i-\alpha_j)(\omega)} - 1\big)f(\omega) \ovl{y_{ij}(\omega)} \d \mu(\omega) \nonumber \\
&\xra[t \to 0]{} \int_\Omega \i \W(\alpha_i - \alpha_j)(\omega) f(\omega) \ovl{y_{ij}(\omega)} \d \mu(\omega)
=\big\langle y , \i \W(\alpha_i - \alpha_j)f \ot e_{ij} \big\rangle. \label{equ-1-proof-prop-Schur-II-group}
\end{align}
We deduce that $f \ot e_{ij}$ belongs to $\dom \i \mathscr{D}_{\alpha,1,\infty}$ and that $\i \mathscr{D}_{\alpha,1,\infty}(x \ot e_{ij})=\i \W(\alpha_i - \alpha_j)f \ot e_{ij}$.

Next we show that the generator $\i \mathscr{D}_{\alpha,1,p}$ is the closure of 
$$
\i \mathscr{D}_{\alpha,1} \co \L^\infty(\Omega) \ot \M_{I,\fin} 
\subseteq \L^p(\Omega,S^p_I) \to \L^p(\Omega,S^p_I).
$$
For any $f \in \L^\infty(\Omega)$ and $i,j \in I$, using \cite[Proposition 4.11]{EnN1}, we have
\begin{align*}
\MoveEqLeft
\frac1t (U_{t,p} - \Id)(f \ot e_{ij})  
=\frac1t\big(\e^{\i t \W(\alpha_i-\alpha_j)}-1\big)f \ot e_{ij}
\ov{\eqref{semigroupe-vers-generateur}}{\xra[t \to 0]{}} \i \W(\alpha_i - \alpha_j)f \ot e_{ij}. 
\end{align*}
We infer that $f \ot e_{ij}$ belongs to $\dom \i \mathscr{D}_{\alpha,1,p}$ and
$\i \mathscr{D}_{\alpha,1,p}(f \ot e_{ij}) = \i \W(\alpha_i - \alpha_j) f \ot e_{ij}$. So it suffices to prove that $\L^\infty(\Omega) \ot \M_{I,\fin}$ is a core for $\i \mathscr{D}_{\alpha,1,p}$. Note that $\Id_{\L^p(\Omega)} \ot \Tron_J$ commutes with $U_{t,p}$. This is easy to check on elements of the form $f \ot e_{ij}$ and extends by a density argument to all of $\L^p(\Omega,S^p_I)$. Thus, for $x \in \dom \mathscr{D}_{\alpha,1,p}$, we have 
$$
\frac1t (U_{t,p} - \Id)(\Id_{\L^p(\Omega)} \ot \Tron_J)(x) 
= (\Id_{\L^p(\Omega)} \ot \Tron_J) \frac1t (U_{t,p} - \Id)(x) 
\to \i (\Id_{\L^p(\Omega)} \ot \Tron_J) \mathscr{D}_{\alpha,1,p}(x).
$$ 
We infer that $(\Id_{\L^p(\Omega)} \ot \Tron_J)(x)$ belongs to $\dom \mathscr{D}_{\alpha,1,p}$ and $\mathscr{D}_{\alpha,1,p}(\Id_{\L^p(\Omega)} \ot \Tron_J)(x) = (\Id_{\L^p(\Omega)} \ot \Tron_J) \mathscr{D}_{\alpha,1,p}(x)\xra[J \to I]{} \mathscr{D}_{\alpha,1,p} (x)$. Moreover, $(\Id_{\L^p(\Omega)} \ot \Tron_J)(x) \to x$ as $J \to I$. Therefore, to show the core property, it suffices to approximate an element $f \ot e_{ij}$ where $f \in \L^p(\Omega)$ in the graph norm of $\mathscr{D}_{\alpha,1,p}$ by elements in $\L^\infty(\Omega) \ot \M_{I,\fin}$. To this end, it suffices to take $1_{|f|\leq n} f \ot e_{ij}$ and to argue similar as beforehand, noting that 
$$
\big(1_{|f|\leq n} \ot 1_{\B(\ell^2_I)}\big) U_{t,p}(f \ot e_{ij}) 
= U_{t,p}\big( 1_{|f| \leq n} f \ot e_{ij}\big).
$$ 
\end{proof}

Similarly to Corollary \ref{cor-Fourier-II-bisectoriel-et-calcul-fonctionnel}, we obtain the following result.

\begin{cor}
\label{cor-Schur-II-bisectoriel-et-calcul-fonctionnel}
Let $1 < p < \infty$. Then the operator $\mathscr{D}_{\alpha,1,p}$ defined above is bisectorial and has a bounded $\HI(\Sigma_\omega^\pm)$ functional calculus to any angle $\omega > 0$.
\end{cor}

\begin{remark}
\label{rem-Schur-II-finite-I}
\normalfont
Now we study the properties of the operator $\mathscr{D}_{\alpha,q}$ from \eqref{Def-Dirac-free} in the case where $I$ is finite and $-1 \leq q<1$. For any $i,j \in I$, we will use the maps $J_{ij} \co \L^p(\Gamma_q(H)) \to \L^p(\Gamma_q(H) \otvn \B(\ell^2_I))$, $x \mapsto x \ot e_{ij}$ and $Q_{ij} \co \L^p(\Gamma_q(H) \otvn \B(\ell^2_I)) \to \L^p(\Gamma_q(H))$, $x \mapsto x_{ij}$. For any $k,l \in I$, we introduce the operator $\scr{L}_{kl} \co \L^p(\Gamma_q(H) \otvn \B(\ell^2_I)) \to \L^p(\Gamma_q(H) \otvn \B(\ell^2_I))$, $x \ot e_{ij} \mapsto \delta_{i=k,j=l}s_q(\alpha_i-\alpha_j)x \ot e_{ij}$. For any $i,j,j,l,k',l' \in I$ ad any $x \in \L^p(\Gamma_q(H))$, we have
\begin{align*}
\MoveEqLeft
\scr{L}_{kl}\scr{L}_{k'l'}(x \ot e_{ij})            
=\delta_{i=k',j=l'}\scr{L}_{kl}\big(s_q(\alpha_i-\alpha_j)x \ot e_{ij}\big) 
=\delta_{i=k',j=l'}\delta_{i=k,j=l}s_q(\alpha_i-\alpha_j)^2x \ot e_{ij}\\
&=\delta_{i=k,j=l}\scr{L}_{k'l'}\big(s_q(\alpha_i-\alpha_j)x \ot e_{ij}\big) 
=\scr{L}_{k'l'}\scr{L}_{kl}(x \ot e_{ij}).
\end{align*}  
So the operators $\scr{L}_{kl}$ commute. Moreover, we have 
\begin{equation}
\label{Dirac-as-sum}
\mathscr{D}_{\alpha,q}
=\sum_{k,l\in I} \scr{L}_{kl}.
\end{equation}
For any $k,l \in I$, note that
\begin{equation}
\label{Compression-L}
\scr{L}_{kl} 
=J_{kl}\L_{s_q(\alpha_k-\alpha_l)}Q_{kl}.
\end{equation}
Note that $s_q(\alpha_k-\alpha_l)$ is selfadjoint. So by \cite[Remark 2.25]{McM1} it is bisectorial of type 0 and admits a bounded bisectorial functional calculus for all $\theta \in (0, \frac{\pi}{2})$ (with $K_\theta = 1$). By adapting \cite[Proposition 8.4]{JMX}, the operator $\L_{s_q(\alpha_k-\alpha_l)}$ is bisectorial of type 0 on $\L^p(\Gamma_q(H))$. By \eqref{Compression-L}, the operator $\scr{L}_{kl}$ is also bisectorial of type 0 on $\L^p(\Gamma_q(H) \otvn \B(\ell^2_I))$ and admits a bounded functional calculus.
It is easy to check\footnote{\thefootnote.
We have
\begin{align*}
\MoveEqLeft
\sum_{k,l} J_{kl} R(\lambda,\L_{s_q(\alpha_k -\alpha_l)})Q_{kl} \left( \lambda - \sum_{i,j} J_{ij} \L_{s_q(\alpha_i - \alpha_j)} Q_{ij} \right) \\
& = \sum_{k,l}J_{kl} \lambda R(\lambda,\L_{s_q(\alpha_k - \alpha_l)})Q_{kl} - \sum_{k,l,i,j} J_{kl}R(\lambda,\L_{s_q(\alpha_k - \alpha_l)})Q_{kl}J_{ij} \L_{s_q(\alpha_i - \alpha_j)}Q_{ij} \\
& = \sum_{k,l} J_{kl}\lambda R(\lambda,\L_{s_q(\alpha_k - \alpha_l)})Q_{kl}- \sum_{k,l,i,j} \delta_{k=i} \delta_{l = j} J_{kl}R(\lambda,\L_{s_q(\alpha_k - \alpha_l)}) \L_{s_q(\alpha_i - \alpha_j)} Q_{ij} \\
& = \sum_{k,l}J_{kl} R(\lambda,\L_{s_q(\alpha_k - \alpha_l)})(\lambda - \L_{s_q(\alpha_k-\alpha_l)}) Q_{kl} \\
& = \sum_{k,l}J_{kl}\Id_{\L^p(\Gamma_q(H))} Q_{kl} \\
& = \Id_{\L^p(\Gamma_q(H) \otvn \B(\ell^2_I))},
\end{align*}
and similarly the other way around.
}
 that for any $\lambda \in \C \backslash \R$, we have $R(\lambda,\sum_{k,l \in I} \scr{L}_{kl}) = \sum_{k,l \in I} J_{kl} R(\lambda,\L_{s_q(\alpha_k-\alpha_l)})Q_{kl}$. Then the Cauchy integral formula yields for any $\theta \in (0,\frac{\pi}{2})$ and any $f \in \HI_0(\Sigma_\theta^\pm)$ that $f(\sum_{k,l \in I} \scr{L}_{kl}) = \sum_{k,l \in I} J_{kl} f(\L_{s_q(\alpha_k - \alpha_l)}) Q_{kl}$. Thus by \eqref{Dirac-as-sum}, we deduce that $\mathscr{D}_{\alpha,q}$ is bisectorial and admits a bounded functional calculus. The study of the operator \eqref{Def-Dirac-group-II} is similar and the verification is left to the reader. Here the assumption is that $G$ is a finite group.
\end{remark}

\section{Appendix : L\'evy measures and 1-cohomology}
\label{Sec-Levy-and-cocycles}

The following observation describes the link between L\'evy-Khinchin decompositions and 1-cocycles. Note that examples of explicit L\'evy measures are given in \cite[p.~184]{BeF1}. 

Let $G$ be a locally compact \textit{abelian} group. Recall that a 1-cocycle is defined by a strongly continuous unitary or orthogonal representation $\pi \co G \to \B(H)$ on a complex or real Hilbert space $H$ and a continuous map $b \co G \to H$ such that $b(s+t)=b(s)+\pi_s(b(t))$ for any $s,t \in G$. We also recall the following fundamental connection between continuous functions $\psi \co G \to \R$ which are conditionally of negative type in the sense of \cite[7.1 Definition]{BeF1} \cite[1.8 Definition p. 89]{BCR} and quadratic forms $q \co G \to \R_+$ together with L\'evy measures $\mu$ from \cite[18.20 Corollary p.~184]{BeF1}. If $\psi$ is such a function satisfying $\psi(0) = 0$, there exist a unique continuous quadratic form $q \co G \to \R_+$ and a unique positive symmetric measure $\mu$ on $\widehat{G} -\{0\}$ such that
\begin{equation}
\label{equ-Levy-representation}
\psi(s) = q(s) + \int_{\widehat{G}- \{ 0 \}} 1 - \Re \langle s , \chi \rangle \d \mu(\chi), \quad s \in G.
\end{equation}
Conversely, any continuous function of this type is conditionally of negative type and $\psi(0) = 0$. Note that a quadratic form is a map satisfying $q(s+t) + q(s-t) = 2q(s) + 2q(t)$ for any $s,t \in G$ \cite[7.18 Definition]{BeF1}. The quadratic form can be computed by $q(s) = \lim_{n \to \infty} \frac{\psi(ns)}{n^2}$, and $\mu$ is the L\'evy measure associated with $\psi$.

In the following result, note that $\pi_s$ is the multiplication operator by the function $\langle s, \cdot \rangle_{G,\hat{G}}$.

\begin{prop}
\label{Prop-Levy-cocycle}
Let $G$ be a locally compact abelian group. Let $\psi \co G \to \R$ be a continuous conditionally of negative type function such that $\psi(0)=0$ and $\lim_{n \to\infty} \frac{\psi(ns)}{n^2}=0$ for any $s \in G$. If $\mu$ is the L\'evy measure of $\psi$ on $\widehat{G}-\{0\}$ then $H=\L^2(\widehat{G}-\{0\},\frac{1}{2}\mu)$, $b \co G \to H$, $s \mapsto (\chi \mapsto 1-\langle s ,\chi\rangle)$ and $\pi \co G \to \B(H)$, $s \mapsto (f \mapsto \langle s, \cdot \rangle_{G,\hat{G}} f)$ define a 1-cocycle on $G$ such that 
$$
\psi(s)
=\norm{b(s)}_H^2 ,\quad s \in G.
$$
\end{prop}

\begin{proof}
For any $s \in G$, note that
\begin{equation}
\label{Petit-calcul-carac}
|1-\langle s ,\chi\rangle|^2
=(1-\langle s ,\chi\rangle)\ovl{1-\langle s ,\chi\rangle}
=1-\langle s ,\chi\rangle-\ovl{\langle s ,\chi\rangle}+1
=2-2\Re \langle s ,\chi\rangle.
\end{equation}
For any $s \in G$, using \eqref{equ-Levy-representation} in the last equality, we deduce that
\begin{align*}
\MoveEqLeft
\norm{b(s)}_{\L^2(\hat{G}-\{0\},\frac{1}{2}\mu)}^2            
=\frac{1}{2} \int_{\hat{G}-\{0\}} |(b(s))(\chi)|^2 \d \mu(\chi)
=\frac{1}{2}\int_{\hat{G}-\{0\}} |1-\langle s ,\chi\rangle|^2 \d \mu(\chi)\\
&\ov{\eqref{Petit-calcul-carac}}{=} \int_{\hat{G}-\{0\}} \big(1-\Re \langle s ,\chi\rangle\big) \d\mu(\chi) 
=\psi(s).
\end{align*} 
Moreover, for any $s,t \in G$ and any $\chi \in \widehat{G}-\{0\}$, we have
\begin{align*}
\MoveEqLeft
(b(s+t))(\chi)            
=1-\langle s+t ,\chi\rangle
=1-\langle s ,\chi\rangle\langle t ,\chi\rangle
=1-\langle s ,\chi\rangle+\langle s ,\chi\rangle(1-\langle t ,\chi\rangle)\\
&=b(s)(\chi)+\pi_s(b(t))(\chi)
=\big(b(s)+\pi_s(b(t))\big)(\chi).
\end{align*} 
Hence, we conclude that $b(s+t)=b(s)+\pi_s(b(t))$.
\end{proof}

In particular, by restriction of scalars, we can consider the real Hilbert space $H_{|\R}=\L^2(\widehat{G}-\{0\},\frac{1}{2}\mu)_{|\R}$, which comes with real scalar product $\langle f,g \rangle_{H_{|\R}} = \langle \Re f, \Re g \rangle_H + \langle \Im f, \Im g \rangle_H$. Then we take $b \co G \to H_{|\R}$, $s \mapsto (\chi \mapsto 1-\langle s ,\chi\rangle)$ and $\pi \co G \to \B(H_{|\R})$, $s \mapsto (f \mapsto \langle s, \cdot \rangle_{G,\hat{G}} f)$. Thus if $\lim_{n \to \infty} \frac{\psi(ns)}{n^2} = 0$ for all $s \in G$, then we obtain a 1-cocycle on $G$ as we considered in Proposition \ref{prop-Schoenberg} suitable for our markovian semigroups of Fourier multipliers.

In the case of an arbitrary continuous function $\psi \co G \to \R_+$ is with $\psi(0) = 0$ of conditionally negative type, we can use \cite[Lemma 3.1]{JMP1} and \cite[Lemma B1]{JMP2} on the quadratic part $q$ of $\psi$ with a direct sum argument for obtaining a  fairly concrete Hilbert space and an associated 1-cocycle.

\vspace{0.2cm}

\textbf{Acknowledgment}.
The authors are supported by the grant of the French National Research Agency ANR-18-CE40-0021 (project HASCON). The first author would like to thank Fran\c coise Lust-Piquard and Cyril L\'evy for some discussions. We are grateful to Marek Bo\.{z}ejko for the reference \cite{Boz2}, Quanhua Xu for the reference \cite{HaP1} and to Frederic Latremoliere for his confirmation of Proposition \ref{Prop-norm-precompact}. The second author acknowledges support by the grant ANR-17-CE40-0021 (project Front).


\vspace{0.2cm}

%
%
%
%
%
%
%
%
%
%
%
%
%
%
%

\footnotesize{
\noindent C\'edric Arhancet\\ 
\noindent13 rue Didier Daurat, 81000 Albi, France\\
URL: \href{http://sites.google.com/site/cedricarhancet}{https://sites.google.com/site/cedricarhancet}\\
cedric.arhancet@protonmail.com\\

\noindent Christoph Kriegler\\
Universit\'e Clermont Auvergne\\
CNRS\\
LMBP\\
F-63000 CLERMONT-FERRAND\\
FRANCE\\
URL: \href{https://lmbp.uca.fr/~kriegler/indexenglish.html}{https://lmbp.uca.fr/{\raise.17ex\hbox{$\scriptstyle\sim$}}\hspace{-0.1cm} kriegler/indexenglish.html}\\
christoph.kriegler@uca.fr

\end{document}